\documentclass[a4paper,leqno]{amsart}
\usepackage{amssymb}
\usepackage[hyphens]{url} 
\usepackage[utf8]{inputenc}
\usepackage{amsmath}
\usepackage{amsfonts}
\usepackage{amssymb}
\usepackage{enumerate}
\usepackage{amsthm}
\usepackage{multicol}
\usepackage{tikz, xcolor, fancyhdr, geometry}
\usepackage{graphicx}

\usepackage[colorlinks]{hyperref}
\hypersetup{
  urlcolor     = blue, 
  linkcolor    = blue, 
  citecolor   = blue   
} 
\DeclareMathOperator{\supp}{supp}

\DeclareMathOperator{\diam}{diam}
\DeclareMathOperator{\body}{body}
\DeclareMathOperator{\dist}{dist}
\def\l2sf{{\Big|\Big|\square _{I}^{\sigma,\mathbf{b}}f\Big|\Big|_{L^2(\sigma)}}}
\def\l2sg{{\Big|\Big|\square _{J}^{\omega,\mathbf{b}^{\ast }}g\Big|\Big|_{L^2(\omega)}}}
\def\A2{{\mathcal{A}_2^\alpha}}
\def\a_I{{E_{I^{\prime }}^{\sigma}\left( \widehat{\square }_{I}^{\sigma ,\flat ,\mathbf{b}}f\right)}}
\def\c_J{{\ E_{J^{\prime}}^{\omega }\left( \widehat{\square }_{J}^{\omega ,\flat ,\mathbf{b}^{\ast}}g\right)}}
\def\R{{\mathbb R}}
\def\Dl{{\frac{\frac{1}{|K^\ell_{in}|_\omega}\int_{K^\ell_{in}}b^\ast_Bd\omega}{\frac{1}{|K^\ell_{in}|_\omega}\int_{K^\ell_{in}}b^{\ast,orig}_Bd\omega}}}
\newtheorem{thm}{Theorem}[section]
\newtheorem{prop}[thm]{Proposition}
\newtheorem{lem}[thm]{Lemma}
\newtheorem{cor}[thm]{Corollary}

\newtheorem{dfn}[thm]{Definition}
\newtheorem{rem}[thm]{Remark}
\newtheorem{notation}[thm]{Notation}

\numberwithin{equation}{section}

\makeatletter
\renewcommand{\tocsection}[3]{%
  \indentlabel{\@ifnotempty{#2}{\bfseries\ignorespaces#1 #2\quad}}\bfseries#3}
\renewcommand{\tocsubsection}[3]{%
  \indentlabel{\@ifnotempty{#2}{\ignorespaces#1 #2\quad}}#3}
\renewcommand{\tocsubsubsection}[3]{%
  \indentlabel{\@ifnotempty{#2}{\hspace{1.5cm}\ignorespaces#1 #2\quad}}#3}

\newcommand\@dotsep{4.5}
\def\@tocline#1#2#3#4#5#6#7{\relax
  \ifnum #1>\c@tocdepth 
  \else
    \par \addpenalty\@secpenalty\addvspace{#2}%
    \begingroup \hyphenpenalty\@M
    \@ifempty{#4}{%
      \@tempdima\csname r@tocindent\number#1\endcsname\relax
    }{%
      \@tempdima#4\relax
    }%
    \parindent\z@ \leftskip#3\relax \advance\leftskip\@tempdima\relax
    \rightskip\@pnumwidth plus1em \parfillskip-\@pnumwidth
    #5\leavevmode\hskip-\@tempdima{#6}\nobreak
    \leaders\hbox{$\m@th\mkern \@dotsep mu\hbox{.}\mkern \@dotsep mu$}\hfill
    \nobreak
    \hbox to\@pnumwidth{\@tocpagenum{\ifnum#1=1\bfseries\fi#7}}\par
    \nobreak
    \endgroup
  \fi}
\AtBeginDocument{%
\expandafter\renewcommand\csname r@tocindent0\endcsname{0pt}
}
\def\l@subsection{\@tocline{2}{0pt}{2.5pc}{5pc}{}}
\makeatother

\begin{document}

\newgeometry{bindingoffset=0.2in, left=1in,right=1in,top=1in,bottom=1in, footskip=.25in}

\title[$n$-dimensional local $Tb$ Theorem]{A two weight local $Tb$ theorem for $n$-dimensional fractional singular integrals}

\author[C.~ Grigoriadis]{Christos Grigoriadis}
\address{Department of Mathematics \\
Michigan State University \\
East Lansing MI}
\email{grigori4@msu.edu}
\author[M.~ Paparizos]{Michail Paparizos}
\address{Department of Mathematics \\
Michigan State University \\East Lansing MI}
\email{paparizo@msu.edu}
\author[E.T. Sawyer]{Eric T. Sawyer}
\address{Department of Mathematics \& Statistics, McMaster University, 1280 Main Street West, Hamilton,\newline\noindent
\phantom{y}\hspace{2ex}Ontario, Canada L8S 4K1}
\email{sawyer@mcmaster.ca}
\thanks{E.T. Sawyer's research is supported in part by NSERC}
\author[C.-Y.~ Shen]{Chun-Yen Shen}
\address{Department of Mathematics,\\
National Taiwan University, 10617, Taipei}
\email{cyshen@math.ntu.edu.tw}
\thanks{C.-Y. Shen supported in part by MOST, through grant 104-2628-M-002
-015 -MY4}
\author[I.~ Uriarte-Tuero]{Ignacio Uriarte-Tuero}
\address{Department of Mathematics \\
University of Toronto \\
40 St. George St., Toronto, Ontario, Canada and \newline\noindent
\phantom{x}\hspace{1.2ex} Department of Mathematics, Michigan State University\\ East Lansing, Michigan}
\email{iuriarte@math.toronto.edu, ignacio@math.msu.edu}
\thanks{ I. Uriarte-Tuero has been partially supported by grant MTM2015-65792-P (MINECO, Spain), and by a Simons \newline\noindent \phantom{x}\hspace{1.2ex} Foundation Collaboration Grant for Mathematicians, Award Number: 637221.}

\maketitle
\begin{abstract}
We obtain a local two weight $Tb$ theorem with an energy side condition for
higher dimensional fractional Calder\'{o}n-Zygmund operators. The proof
follows the general outline of the proof for the corresponding
one-dimensional $Tb$ theorem in \cite{SaShUr12}, but encountering a number
of new challenges, including several arising from the failure in higher
dimensions of T. Hyt\"{o}nen's one-dimensional two weight $A_{2}$ inequality 
\cite{Hyt}.
\end{abstract}

\tableofcontents

\section{Introduction}

Boundedness properties of Calder\'{o}n-Zygmund singular integrals arise in
the most critical cases of the study of virtually all partial differential
equations, from Schr\"{o}dinger operators in quantum mechanics to
Navier-Stokes equations in fluid flow, as well as in the investigation of
a number of topics in geometry and analysis. In particular, the study of
boundedness of these operators from one weighted space $L^{2}\left( \mathbb{R%
}^{n};\sigma \right) $ to another $L^{2}\left( \mathbb{R}^{n};\omega \right) 
$, not only extends the scope of application in many cases, but reveals the
important properties of the kernels associated with the individual operators
under consideration, often hidden without such investigation into two weight
norm inequalities. The purpose of this monograph is to prove a general
characterization regarding boundedness of Calder\'{o}n-Zygmund singular
integrals from $L^{2}\left( \mathbb{R}^{n};\sigma \right) $ to $L^{2}\left( 
\mathbb{R}^{n};\omega \right) $, for locally finite positive Borel measures $%
\sigma $ and $\omega $, subject to some natural buffer conditions. This
result, a so-called local two weight $Tb$ theorem in $\mathbb{R}^{n}$,
includes much, if not most, of the known theory on two weight $L^{2}$%
-boundedness of singular integrals. We now digress to a brief history of
that part of this theory that is relevant to our purpose here.

Given a Calder\'{o}n-Zygmund kernel $K\left( x,y\right) $ in Euclidean space 
$\mathbb{R}^{n}$, a classical problem for some time was to identify optimal
cancellation conditions on $K$ so that there would exist an associated
singular integral operator $Tf\left( x\right) \sim \int K\left( x,y\right)
f\left( y\right) dy$ bounded on $L^{2}\left( \mathbb{R}^{n}\right) $. After
a long history, involving contributions by many authors\footnote{%
see e.g. \cite[page 53]{Ste} for references to the earlier work in this
direction}, this effort culminated in the decisive $T1$ theorem of David and
Journ\'{e} \cite{DaJo}, in which boundedness of an operator $T$ on $%
L^{2}\left( \mathbb{R}^{n}\right) $ associated to $K$, was characterized by 
\begin{equation*}
T\mathbf{1},T^{\ast }\mathbf{1}\in BMO,
\end{equation*}%
together with a weak boundedness property for some $\eta >0$,%
\begin{eqnarray}
&&\left\vert \int_{Q}T\varphi \left( x\right) \ \psi \left( x\right)
dx\right\vert \lesssim \sqrt{\left\Vert \varphi \right\Vert _{\infty
}\left\vert Q\right\vert +\left\Vert \varphi \right\Vert _{{Lip}\eta
}\left\vert Q\right\vert ^{1+\frac{\eta }{n}}}\sqrt{\left\Vert \psi
\right\Vert _{\infty }\left\vert Q\right\vert +\left\Vert \psi \right\Vert _{%
{Lip}\eta }\left\vert Q\right\vert ^{1+\frac{\eta }{n}}},  \label{WBPDJ} \\
&&\text{for all }\varphi ,\psi \in {Lip}\eta \text{ with }{supp}\varphi ,{%
supp}\psi \subset Q,\text{ and all cubes }Q\subset \mathbb{R}^{n};  \notag
\end{eqnarray}%
\emph{equivalently} by two testing conditions taken uniformly over
indicators of cubes,%
\begin{equation*}
\int_{Q}\left\vert T\mathbf{1}_{Q}\left( x\right) \right\vert ^{2}dx\lesssim
\left\vert Q\right\vert \text{ and }\int_{Q}\left\vert T^{\ast }\mathbf{1}%
_{Q}\left( x\right) \right\vert ^{2}dx\lesssim \left\vert Q\right\vert ,\ \
\ \ \ \text{all cubes }Q\subset \mathbb{R}^{n}.
\end{equation*}%
The optimal cancellation conditions, which in the words of Stein were `a
rather direct consequence of' the $T1$ theorem, were given in \cite[Theorem
4, page 306]{Ste}, involving integrals of the kernel over shells:%
\begin{eqnarray}
&&\int_{\left\vert x-x_{0}\right\vert <N}\left\vert \int_{\varepsilon
<\left\vert x-y\right\vert <N}K^{\alpha }\left( x,y\right) dy\right\vert
^{2}dx\leq \mathfrak{A}_{K^{\alpha }}\ \int_{\left\vert x_{0}-y\right\vert
<N}dy,  \label{can cond} \\
&&\ \ \ \ \ \ \ \ \ \ \ \ \ \ \ \ \ \ \ \ \ \ \ \ \ \text{for all }%
0<\varepsilon <N\text{ and }x_{0}\in \mathbb{R}^{n},  \notag
\end{eqnarray}%
together with a dual inequality.

We now come to a point of departure for two separate threads of further
research on cancellation conditions. The first thread treats extensions of
these testing conditions to the boundedness of Calder\'{o}n-Zygmund
operators\ on more general weighted spaces $L^{2}\left( w\right) \rightarrow
L^{2}\left( w\right) $, and even from one weighted space to another, $%
L^{2}\left( \sigma \right) \rightarrow L^{2}\left( \omega \right) $. The
second thread replaces the family of testing functions $\left\{ \mathbf{1}%
_{Q}\right\} _{Q\in \mathcal{D}}$ with families $\left\{ b_{Q}\right\}
_{Q\in \mathcal{D}}$ more amenable to the boundedness of the operator at
hand, subject of course to some sort of nondegeneracy conditions. Finally
the two threads recombine in the theorem of this paper. See the diagram.

\begin{center}
\includegraphics[height=8 cm, width=7 cm]{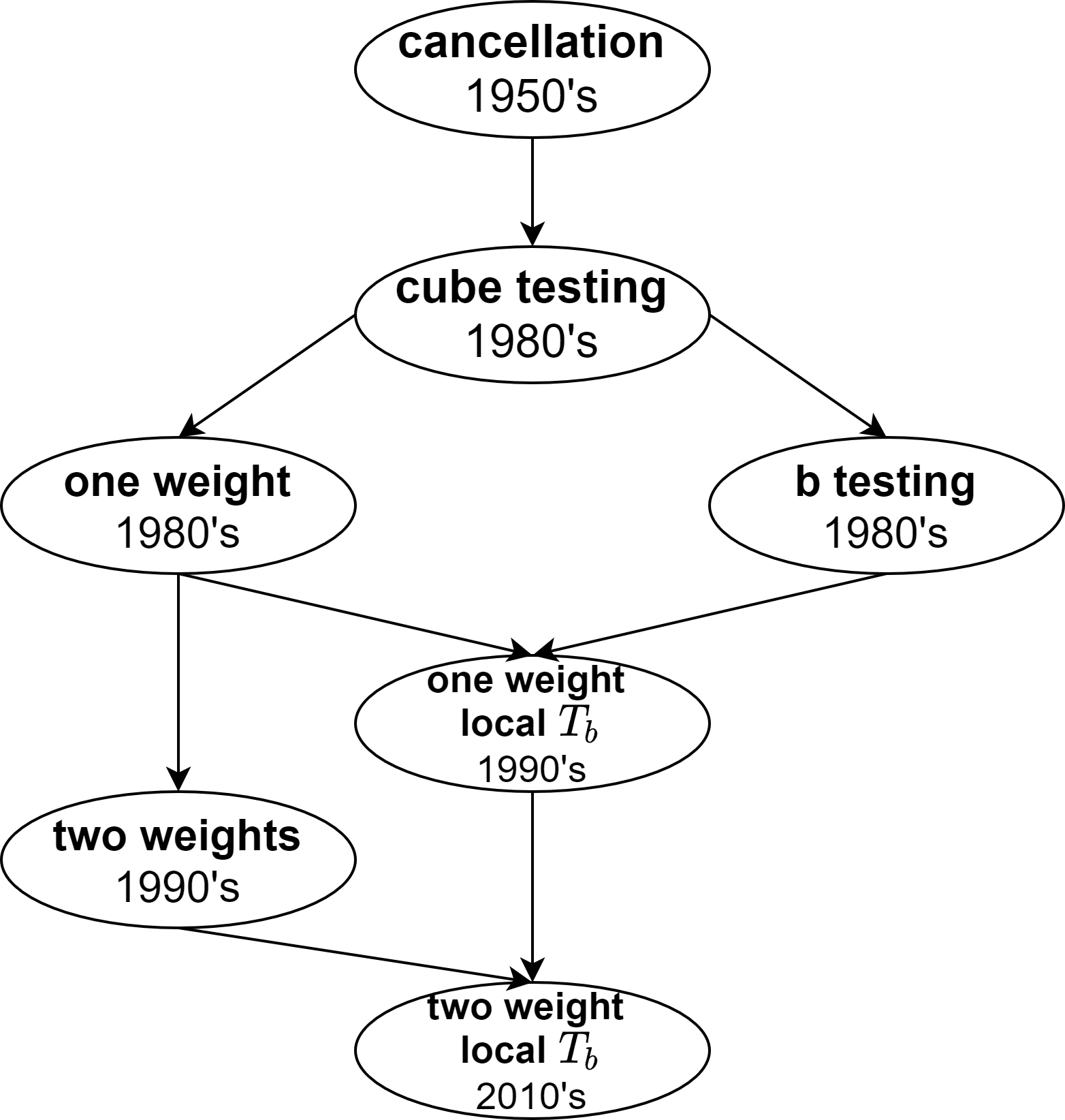}
\end{center}

\subsection{Weighted spaces}

An obvious next step was to replace Lebesgue measure with a fixed $A_{2}$
weight $w$,%
\begin{equation*}
\sup_{\text{cubes }Q\subset \mathbb{R}^{n}}\left( \frac{1}{\left\vert
Q\right\vert }\int_{Q}w\left( x\right) dx\right) \left( \frac{1}{\left\vert
Q\right\vert }\int_{Q}\frac{1}{w\left( x\right) }dx\right) \lesssim 1\ ,
\end{equation*}%
and ask when $T$ is bounded on $L^{2}\left( w\right) $, i.e. satisfies the
one weight norm inequality. For elliptic Calder\'{o}n-Zygmund operators $T$,
this question is reduced to the David Journ\'{e} theorem using two results
from decades ago, namely the 1956 Stein-Weiss interpolation with change of
measures theorem \cite{StWe}, and the 1974 Coifman and Fefferman extension 
\cite{CoFe} of the one weight Hilbert transform inequality of Hunt,
Muckenhoupt and Wheeden \cite{HuMuWh}, to a large class of general Calder%
\'{o}n-Zygmund operators $T$\footnote{%
Indeed, if $T$ is bounded on $L^{2}\left( w\right) $, then by duality it is
also bounded on $L^{2}\left( \frac{1}{w}\right) $, and the Stein-Weiss
interpolation theorem with change of measure shows that $T$ is bounded on
unweighted $L^{2}\left( \mathbb{R}^{n}\right) $. Conversely, if $T$ is
bounded on unweighted $L^{2}\left( \mathbb{R}^{n}\right) $, the proof in 
\cite{CoFe} shows that $T$ is bounded on $L^{2}\left( w\right) $ using $w\in
A_{2}$.}. A motivating example, for the case of the conjugate function $H$
on the unit circle, arose in the Helson-Szeg\"{o} theorem that characterized
the boundedness of $H$ on $L^{2}\left( w\right) $ by the existence of
bounded functions $u$ and $v$ on the circle with $\left\Vert v\right\Vert
_{\infty }<\frac{\pi }{2}$ and $w=e^{u+Hv}$. The equivalence with the $A_{2}$
condition on $w$ follows from the results just mentioned, and the question
of a direct argument linking the Helson-Szeg\"{o} condition to the $A_{2}$
condition has remained a tantalizing puzzle for decades since. See \cite[%
pages 222-227]{Ste} for this and other applications of one weight theory,
such as to the Dirichlet problem for elliptic divergence form operators with
bounded measurable coefficients.

However, for a pair of \emph{different} measures $\left( \sigma ,\omega
\right) $, the question is wide open in general, and we now focus our
discussion on the main problem\ considered in this monograph, that of
characterizing boundedness of a \emph{general} Calder\'{o}n-Zygmund operator 
$T$ from one $L^{2}\left( \sigma \right) $ space to another $L^{2}\left(
\omega \right) $ space, subject to natural buffer conditions on the weight
pair $\left( \sigma ,\omega \right) $. First we note that for the primordial
singular integral, namely the Hilbert transform $H$ in dimension one, the
two weight inequality was completely solved by establishing the NTV conjecture in the two part paper \cite%
{LaSaShUr2};\cite{Lac}, see also \cite{Hyt} for the general case permitting
common point masses, where it was shown that $H$ is bounded from $%
L^{2}\left( \sigma \right) $ to $L^{2}\left( \omega \right) $ if and only if
the testing and one-tailed Muckenhoupt conditions hold, i.e.%
\begin{eqnarray*}
&&\int_{I}\left\vert H\left( \mathbf{1}_{I}\sigma \right) \right\vert
^{2}d\omega \lesssim \int_{I}d\sigma \text{ and }\int_{I}\left\vert H\left( 
\mathbf{1}_{I}\omega \right) \right\vert ^{2}d\sigma \lesssim
\int_{I}d\omega ,\ \ \ \ \ \text{uniformly over all intervals }I\subset 
\mathbb{R}^{n}, \\
&&\left( \int_{\mathbb{R}}\frac{\left\vert I\right\vert }{\left\vert
I\right\vert ^{2}+\left\vert x-c_{I}\right\vert ^{2}}d\sigma \left( x\right)
\right) \left( \frac{1}{\left\vert I\right\vert }\int_{I}d\omega \right)
\lesssim 1,\ \text{and its dual},\ \ \ \ \ \text{uniformly over all
intervals }I\subset \mathbb{R}^{n}.
\end{eqnarray*}%
For $\alpha $-fractional Riesz transforms in higher dimensions $n\geq 2$, it
is known (except when $\alpha =n-1$) that the two weight norm inequality 
\emph{with doubling measures} is equivalent to the fractional one-tailed
Muckenhoupt and $T1$ cube testing conditions, see \cite[Theorem 1.4]{LaWi}
and \cite[Theorem 2.11]{SaShUr9}. Here a positive measure $\mu $ is doubling
if%
\begin{equation*}
\int_{2Q}d\mu \lesssim \int_{Q}d\mu ,\ \ \ \ \ \text{all cubes }Q\subset 
\mathbb{R}^{n}.
\end{equation*}%
However, these results rely on certain `positivity' properties of the
gradient of the kernel (which for the Hilbert transform kernel $\frac{1}{y-x}
$ is simply $\frac{d}{dx}\frac{1}{y-x}>0$ for $x\neq y$), something that is
not available for general elliptic, or even strongly elliptic, fractional
Calder\'{o}n-Zygmund operators.

Then in [Saw] this $T1$ theorem was extended to arbitrary \emph{smooth}
Calder\'{o}n-Zygmund operators and $\mathcal{A}_{2}$ measure pairs $\left(
\sigma ,\omega \right) $ with doubling comparable measures, where a pair of
doubling measures $\sigma $ and $\omega $ are \emph{comparable} in the sense
of Coifman and Fefferman \cite{CoFe}, if the measures are mutually
absolutely continuous, uniformly at all scales - i.e. there exist $0<\beta
,\gamma <1$ such that%
\begin{equation*}
\frac{\left\vert E\right\vert _{\sigma }}{\left\vert Q\right\vert _{\sigma }}%
<\beta \Longrightarrow \frac{\left\vert E\right\vert _{\omega }}{\left\vert
Q\right\vert _{\omega }}<\gamma \text{ for all Borel subsets }E\text{ of a
cube }Q.
\end{equation*}%
Subsequently, in \cite{Gr}, it was shown that the pivotal conditions of NTV
are implied by the two weight $\mathcal{A}_{2}$ condition if the weights are 
$A_{\infty }$, and pointed out that this then extends the $T1$ theorem to
pairs of $A_{\infty }$ weights for rougher Calder\'{o}n-Zygmund operators
upon applying the $T1$ theorem of \cite{SaShUr7}.

\subsection{$Tb$ theorems}

The original $T1$ theorem of David and Journ\'{e} \cite{DaJo}, which
characterized boundedness of a singular integral operator by testing over
indicators $\mathbf{1}_{Q\text{ }}$ of cubes $Q$, was quickly extended to a $%
Tb$ theorem by David, Journ\'{e} and Semmes \cite{DaJoSe}, in which the
indicators $\mathbf{1}_{Q\text{ }}$ were replaced by testing functions $b%
\mathbf{1}_{Q\text{ }}$ for an accretive function $b$, i.e. $0<c\leq {Re}%
b\leq \left\vert b\right\vert \leq C<\infty $. Here the accretive function $%
b $ could be chosen to adapt well to the operator at hand, resulting in
almost immediate verification of the $b$-testing conditions, despite
difficulty in verifying the $1$-testing conditions. One motivating example
of this phenomenon is the boundedness of the Cauchy integral on Lipschitz
curves, easily obtained from the above $Tb$ theorem\footnote{%
The problem reduces to boundedness on $L^{2}\left( \mathbb{R}\right) $ of
the singular integral operator $C_{A}$ with kernel $K_{A}\left( x,y\right)
\equiv \frac{1}{x-y+i\left( A\left( x\right) -A\left( y\right) \right) }$,
where the curve has graph $\left\{ x+iA\left( x\right) :x\in \mathbb{R}%
\right\} $. Now $b\left( x\right) \equiv 1+iA^{\prime }\left( x\right) $ is
accretive and the $b$-testing condition $\int_{I}\left\vert C_{A}\left( 
\mathbf{1}_{I}b\right) \left( x\right) \right\vert ^{2}dx\leq \mathfrak{T}%
_{H}^{b}\left\vert I\right\vert $ follows from $\left\vert C_{A}\left( 
\mathbf{1}_{I}b\right) \left( x\right) \right\vert ^{2}\approx \ln \frac{%
x-\alpha }{\beta -x}$, for $x\in I=\left[ \alpha ,\beta \right] $. In the
case of a $C^{1,\delta }$ curve, the kernel $K_{A}$ is $C^{1,\delta }$ and a 
$Tb$ theorem applies with $T=C_{A}$ and $\sigma =\omega =dx$, to show that $%
C_{A}$ is bounded on $L^{2}\left( \mathbb{R}\right) $.}. See e.g. \cite[%
pages 310-316]{Ste}.

Subsequently, M. Christ \cite{Chr} obtained a far more robust \emph{local} $%
Tb$ theorem in the setting of homogeneous spaces, in which the testing
functions could be further specialized to $b_{Q}\mathbf{1}_{Q\text{ }}$,
where now the accretive functions $b_{Q}$ can be \emph{chosen by the reader
to differ} for \emph{each} cube $Q$. Applications of the local $Tb$ theorem
included boundedness of layer potentials, see e.g. \cite{AAAHK} and
references there; and the Kato problem, see \cite{HoMc}, \cite{HoLaMc} and 
\cite{AuHoLaMcTc}: and many authors, including G. David \cite{Dav1};
Nazarov, Treil and Volberg \cite{NTV3}, \cite{NTV2}; Auscher, Hofmann,
Muscalu, Tao and Thiele \cite{AuHoMuTaTh}, Hyt\"{o}nen and Martikainen \cite%
{HyMa}, and more recently Lacey and Martikainen \cite{LaMa}, set about
proving extensions of the local $Tb$ theorem, for example to include a
single upper doubling weight together with weaker upper bounds on the
function $b$. But these extensions were modelled on the `nondoubling'
methods that arose in connection with upper doubling measures in the
analytic capacity problem, see Mattila, Melnikov and Verdera \cite{MaMeVe},
G. David \cite{Dav1}, \cite{Dav2}, X. Tolsa \cite{Tol}, and alsoVolberg \cite%
{Vol}, and were thus constrained to a single weight - a setting in which
both the Muckenhoupt and energy conditions follow from the upper doubling
condition.

More recently, in a precursor to the present paper, \cite{SaShUr12} obtained
a general two weight $Tb$ theorem for the Hilbert transform on the real
line. In this paper, we extend this precursor to higher dimensions. As in 
\cite{SaShUr12}, we adapt methods from the theory of two weight $T1$
theorems, which arose from \cite{NTV4}, \cite{Vol}, \cite{LaSaShUr2}, \cite%
{Lac}, \cite{SaShUr7} and \cite{SaShUr9}, and were used in \cite{HyMa} as
well, to prove a two weight local $Tb$ theorem. These methods involve the
`testing' perspective toward characterizing two weight norm inequalities for
an operator $T$. As suggested by work originating in \cite{DaJo} and \cite%
{Saw3}, it is plausible to conjecture that a given operator $T$ is bounded
from one weighted space to another if and only if both it and its dual are
bounded when tested over a suitable family of functions related
geometrically to $T$, e.g. testing over indicators of intervals for
fractional integrals $T$ as in \cite{Saw3}.

\subsection{Challenges in higher dimensional two weight $Tb$ theory}

A number of difficulties arise in {generalizing to higher dimensions the
work that was done in \cite{SaShUr12} for dimension $n=1$. The main
difficulty lies in the strictly one-dimensional nature of a fundamental
inequality of Hyt\"{o}nen, namely that local testing, i.e. testing the
integral of }$\left\vert T_{\sigma }\mathbf{1}_{Q}\right\vert ^{2}$ over the
cube $Q$,{\ together with the $A_{2}$ condition, implies full testing,
meaning that }$\left\vert T_{\sigma }\mathbf{1}_{Q}\right\vert ^{2}${\ is
integrated\ over the entire space }$\mathbb{R}^{n}${. For the proof\ of full
testing, Hyt\"{o}nen uses an inequality for the Hardy operator that is true
only in dimension }$n=1${\ - in fact it was recently proved in \cite{GP}
that this property of the Hardy operator is not available in higher
dimensions. Then with full testing in hand, we obtain a number of properties
that\ greatly simplify matters. Here are the main challenges encountered in
passing from the one-dimensional setting to the higher dimensional analog.}

\begin{enumerate}
\item \textbf{The nearby form}. The main difficulty in proving the $Tb$
theorem in dimensions $n>1$ arises in treating the nearby form in Chapter 5.
Full testing is used repeatedly everywhere in this chapter, and a demanding
technical approach involving random surgery and averaging, is needed to
circumvent full testing throughout this chapter. In particular, to obtain
estimates over adjacent cubes, we decompose one of the cubes into a smaller
rectangle that is separated from the other cube by a halo. The separated
part is estimated by Muckenhoupt's $A_{2}$ condition, while the halo part is
estimated by applying probability over grids. An illustrative example is the
following. Let $I$ be a cube in the grid associated to the function $f$ and $%
J$ a cube in the grid associated to the function $g$. Let also $%
b_{I},b_{J}^{\ast }$ be the testing functions used in the theorem for these
cubes.

We would like to estimate 
\begin{equation*}
\int T_{\sigma }^{\alpha }\left( b_{I}\mathbf{1}_{I\backslash J}\right)
b_{J}^{\ast }\mathbf{1}_{J}d{\omega .}
\end{equation*}%
The domains of integration inside the operator and inside the integral are
adjacent. In dimension $n=1$ we could use Hyt\"{o}nen's result. Now we
instead argue by splitting the integral as follows: 
\begin{equation*}
\left\vert \int T_{\sigma }^{\alpha }\left( b_{I}\mathbf{1}_{I\backslash
J}\right) b_{J}^{\ast }\mathbf{1}_{J}d{\omega }\right\vert \leq \left\vert
\int T_{\sigma }^{\alpha }\left( b_{I}\mathbf{1}_{I\backslash (1+\delta
)J}\right) b_{J}^{\ast }\mathbf{1}_{J}d{\omega }\right\vert +\left\vert \int
T_{\sigma }^{\alpha }\left( b_{I}\mathbf{1}_{(I\backslash J)\cap (1+\delta
)J}\right) b_{J}^{\ast }\mathbf{1}_{J}d{\omega }\right\vert .
\end{equation*}

The first term on the right hand side, where the domains inside the operator
and the integral are disjoint with positive distance, is bounded by a
constant multiple, depending on $\delta $ and $n$, times the $A_{2}$
constant. Using averaging over grids, the second term on the right hand side
is bounded by $\delta \mathfrak{N}_{T^{\alpha }}$ where the small $\delta $
gain comes from the fact that $|(I\backslash J)\cap (1+\delta )J|^{\frac{1}{n%
}}\approx \delta |I|$ where $|\cdot |$ denotes the Lebesque measure of the
cube.

\item \textbf{Splitting forms}. Here we begin with a pair of smooth
compactly supported functions $(f,g)$ and we would like to\ decompose the
functions into their Haar expansions. However, when we select a grid $%
\mathcal{G}$ for $f$, the support of $f$ may not be contained in any of the
dyadic cubes in the grid $\mathcal{G}$, with a similar problem when
selecting a grid $\mathcal{H}$ for $g$. To deal with this, we follow NTV by
adding and subtracting certain averages for these terms, resulting in four
integrals to be controlled by our hypotheses. In the one dimensional
setting, full testing was used to eliminate three out of the four such
integrals that appear after decomposing the functions in sums of martingale
differences. Here in this paper, the argument must be adjusted to avoid
using full testing by averaging over the two grids $\mathcal{G}$ and $%
\mathcal{H}$ associated with $f$ and $g$.

\item \textbf{Pointwise Lower Bound Property (PLBP)}. In \cite{SaShUr12} for 
$n=1$, the $PLBP$ was used to control terms involving certain `modified dual
martingale differences' in which a factor $b_{Q}$ had been removed.
Moreover, it was proved there that, without loss of generality, the $p$%
-weakly accretive families of testing functions $b_{Q}$ and $b_{Q}^{\ast }$
for $Q\in \mathcal{P}$ \ could be assumed to satisfy the \emph{pointwise
lower bound property}, written $PLBP$:%
\begin{equation*}
\left\vert b_{Q}\left( x\right) \right\vert \geq c_{1}>0\ \ \ \ \ \text{for}%
\ Q\in \mathcal{P}\text{ and }\sigma \text{-a.e. }x\in \mathbb{R},
\end{equation*}%
for some positive constant $c_{1}$. However, this reduction to assuming $%
PLBP $ depended heavily on Hyt\"{o}nen's $\mathcal{A}_{2}$ characterization
for supports on disjoint intervals, something that is unavailable in higher
dimensions \cite{GP}. To circumvent this difficulty we used an observation
(that goes back to Hyt\"{o}nen and Martikainen) that under the additional
assumption that the \emph{breaking cubes} $Q$, those for which there is a
dyadic child $Q^{\prime }$ of $Q$ with $b_{Q^{\prime }}\neq \mathbf{1}%
_{Q^{\prime }}b_{Q}$, satisfy an appropriate Carleson measure condition.

\item \textbf{Indented corona}. In chapter 8 (dealing with the stopping
form) we construct an `indented corona'. In dimension $n=1$ this
construction simply reduces to consideration of the `left and right ends' of
the intervals. In the absence of `right and left ends' in higher dimensions,
this simple construction is replaced by a more intricate tree of Carleson
cubes.
\end{enumerate}

\subsection{A higher dimensional two weight local $Tb$ theorem}

We begin with a discussion of the buffer conditions we will assume on the
pair $\left( \sigma ,\omega \right) $ of locally finite positive Borel
measures arising in the $Tb$ theorem.

\textbf{Muckenhoupt conditions}: Even for the simplest singular integral,
the Hilbert transform, testing over indicators of intervals no longer
suffices for boundedness \footnote{%
consider e.g. $d\omega \left( x\right) =\delta _{0}\left( x\right) $ and $%
d\sigma \left( x\right) =\left\vert x\right\vert dx$.}, and an additional
`side condition' on the weight pair is required - namely the Muckenhoupt $%
\mathfrak{A}_{2}$ condition, a simpler form of which was shown by Hunt,
Muckenhoupt and Wheeden \cite{HuMuWh} to characterize the one weight
inequality for the Hilbert transform. This side condition is a size
condition on the weight pair that is typically shown to be necessary by
testing over so-called tails of indicators of intervals, and indeed is known
to be necessary for boundedness of a broad class of fractional singular
integrals that are `strongly elliptic' \cite{SaShUr7}. Using this side condition of
Muckenhoupt, the solution of the NTV conjecture, due to three of the authors
and M. Lacey in the two part paper \cite{LaSaShUr2}-\cite{Lac}, shows that
the Hilbert transform $H$ is bounded between weighted $L^{2}$ spaces if and
only if the Muckenhoupt condition and the two testing conditions over
indicators of intervals all hold. However, the testing conditions for
singular integrals, unlike those for positive operators such as fractional
integrals, are extremely unstable and in principle difficult to check \cite%
{LaSaUr2}. On the other hand, given a weight pair, it may be possible to
produce a family of testing functions adapted to intervals on which the
boundedness of the operator is evident. In such a case, one would like to
conclude that finding an appropriately \emph{nondegenerate} family of such
testing functions, for which the corresponding testing conditions hold, is
enough to guarantee boundedness of the operator - bringing us back to a
local $Tb$ theorem. In any event, one would in general like to understand
the weakest testing conditions that are sufficient for two weight
boundedness of a given operator.

\textbf{Energy conditions}: Our $Tb$ theorem lies in this direction, but the
method of proof requires in addition a second `side condition', namely the
so-called energy condition, introduced in \cite{LaSaUr2}. The energy
condition is necessary for the boundedness of the Hilbert transform, and
actually follows there from testing over indicators of intervals and,
through the Muckenhoupt condition, testing over tails of indicators of intervals as well. More generally, it is known that the energy condition is
necessary for boundedness of \emph{gradient elliptic} fractional singular
integrals on the real line , but fails to be necessary for certain elliptic
singular integrals on the line and for even the nicest of singular operators
in higher dimensions \cite{SaShUr11}.

\textbf{Failure of sufficiency of Muckenhoupt and Energy conditions}:
However, the weight pair $\left( \omega ,\ddot{\sigma}\right) $ constructed
in \cite{LaSaUr2} satisfies the Muckenhoupt and energy conditions, yet is easily seen to fail to satisfy the norm inequality for the Hilbert transform. This shows that, even
assuming the necessary conditions of Muckenhoupt and energy, we still need
some sort of testing conditions, and our $Tb$ theorem essentially leaves the
choice of testing conditions at our disposal - subject only to nondegeneracy
and size conditions. 

$ $\\
\textbf{The main two weight local }$Tb$\textbf{\ theorem}: Here is a brief
statement of our main theorem.

\begin{thm}[local $Tb$ in higher dimensions]
\label{thm brief} Let $T^{\alpha }$ denote a Calder\'{o}n-Zygmund operator
on $\mathbb{R}^{n}$, and let $\sigma $ and $\omega $ be locally finite
positive Borel measures on $\mathbb{R}^{n}$ that satisfy the energy and
Muckenhoupt buffer conditions. Then $T_{\sigma }^{\alpha }$, where $%
T_{\sigma }^{\alpha}f\equiv T^{\alpha }\left( f\sigma \right) $, is bounded from
\thinspace $L^{2}\left( \sigma \right) $ to $L^{2}\left( \omega \right) $ 
\emph{if and only if} the $\mathbf{b}$-testing and $\mathbf{b}^{\ast }$%
-testing conditions%
\begin{equation}
\int_{I}\left\vert T_{\sigma }^{\alpha }b_{I}\right\vert ^{2}d\omega \leq
\left( \mathfrak{T}_{T^{\alpha }}^{\mathbf{b}}\right) ^{2}\left\vert
I\right\vert _{\sigma }\text{ and }\int_{J}\left\vert T_{\omega }^{\alpha
,\ast }b_{J}^{\ast }\right\vert ^{2}d\sigma \leq \left( \mathfrak{T}%
_{T^{\alpha }}^{\mathbf{b}^{\ast },\ast }\right) ^{2}\left\vert J\right\vert
_{\omega }\ ,  \label{b and b*  test}
\end{equation}%
taken over two families of test functions $\left\{ b_{I}\right\} _{I\in 
\mathcal{P}}$ and $\left\{ b_{J}^{\ast }\right\} _{J\in \mathcal{P}}$, where 
$b_{I}$ and $b_{J}^{\ast }$ are only required to be nondegenerate in an
average sense, and to be just slightly better than $L^{2}$ functions
themselves, namely $L^{p}$ for some $p>2$.
\end{thm}

The families of test functions $\left\{ b_{I}\right\} _{I\in \mathcal{P}}$
and $\left\{ b_{J}^{\ast }\right\} _{J\in \mathcal{P}}$ in the $Tb$ theorem
above are nondegenerate and slightly better than $L^{2}$ functions, but
otherwise remain at the disposal of the reader. It is this flexibility in
choosing families of test functions that distinguishes this characterization
as compared to the corresponding $T1$ theorem. The $Tb$ theorem here
generalizes many of the one-weight $Tb$ theorems, since in
the upper doubling case, the Muckenhoupt $\mathfrak{A}_{2}$ condition and
the energy condition easily follow from the upper doubling condition. Recall
that in the one-weight case with doubling and upper doubling measures $\mu $%
, there has been a long and sustained effort to relax the integrability
conditions of the testing functions: see e.g. S. Hofmann \cite{Hof} and
Alfonseca, Auscher, Axelsson, Hofmann and Kim \cite{AAAHK}. Subsequently, Hyt%
\"{o}nen- Martikainen \cite{HyMa} assumed $Tb$ in $L^{s}\left( \mu \right) $
for some $s>2$, and the one weight theorem with testing functions $b$ in $%
L^{2}\left( \mu \right) $ was attained by Lacey-Martikainen \cite{LaMa}, but
their argument strongly uses methods not immediately available in the two
weight setting.

\subsection{Application: a function theory characterization\ without buffer conditions}

Here we consider the $\alpha $-fractional vector Riesz transform $\mathbf{R}%
^{n,\alpha }=\left( R_{1}^{n,\alpha },,,,R_{n}^{n,\alpha }\right) $ on $%
\mathbb{R}^{n}$, whose components $R_{i}^{n,\alpha }$ have convolution
kernels $c_{n,\alpha }\frac{x_{i}}{\left\vert x\right\vert ^{n-\alpha +1}}$.
It is shown in \cite{SaShUr7} that the \emph{Muckenhoupt} conditions are
necessary for boundedness of $\mathbf{R}^{n,\alpha }:L^{2}\left( \sigma
\right) \rightarrow L^{2}\left( \omega \right) $. Moreover, there are in the
literature a number of geometric constraints on the measure pair $\left(
\sigma ,\omega \right) $ under which the \emph{energy} conditions are
necessary for boundedness of $\mathbf{R}^{n,\alpha }:L^{2}\left( \sigma
\right) \rightarrow L^{2}\left( \omega \right) $. For example, this is true
if

\begin{enumerate}
\item at least one of the two measures $\sigma ,\omega $ is compactly
supported on a $C^{1,\delta }$ curve in $\mathbb{R}^{n}$, see \cite{SaShUr8}, or

\item each measure $\sigma ,\omega $ satisfies a $k$-dispersed condition for
certain $k$ depending only $n$\ and $\alpha $, see \cite{SaShUr9}, or

\item each measure $\sigma ,\omega $ satisfies a uniformly full dimension
condition, see \cite{LaWi}.
\end{enumerate}

Under any of the above three geometric constraints on the measure pair $%
\left( \sigma ,\omega \right) $, the restriction of our $Tb$ Theorem \ref%
{thm brief} to the Riesz transform $\mathbf{R}^{\alpha }$ is thus improved
by eliminating the assumption of buffer conditions.

\begin{thm}
Let $\mathbf{R}^{n,\alpha }$ be the $\alpha $-fractional vector Riesz
transform on $\mathbb{R}^{n}$, and let $\sigma $ and $\omega $ be locally
finite positive Borel measures on $\mathbb{R}^{n}$ that satisfy at least one
of the three geometric constraints listed above. Then $\mathbf{R}^{n,\alpha }_\sigma$, where $\mathbf{R}^{n,\alpha }_\sigma f\equiv \mathbf{R}^{n,\alpha}\left( f\sigma \right) $, is bounded from \thinspace $L^{2}\left( \sigma \right) $ to $L^{2}\left( \omega \right) $ \emph{if and
only if} the energy and Muckenhoupt conditions hold, as well as the $\mathbf{%
b}$-testing and $\mathbf{b}^{\ast }$-testing conditions (\ref{b and b*  test}%
) for $T^{\alpha }=\mathbf{R}^{\alpha }$.
\end{thm}

An application of this theorem for $\mathbf{R}^{\alpha }$ in the case $n=2$
and $\alpha =1$ arises in characterizing Carleson measures for model spaces $%
K_{\theta }$, where $\theta $ is an inner function on the unit disk $\mathbb{%
D}$. See \cite{LaSaShUrWi} for terminology and a discussion of this problem.
The following theorem, but without condition (4), was proved in Lacey,
Sawyer, Shen, Uriarte-Tuero and Wick \cite[Theorem 1.15]{LaSaShUrWi}%
\footnote{%
The $T1$ testing in (3) is taken over only Carleson squares, whereas the $Tb$ testing in (4) is taken over all squares.}.\ Note that the Cauchy
transform $\mathcal{C}f$ is given by $R_{1}^{2,1}+iR_{2}^{2,1}$, and so its
boundedness is equivalent to that of the vector Riesz transform $\mathbf{R}%
^{2,1}$.

\begin{thm}
Let $\mu $ be a positive Borel measure on $\overline{\mathbb{D}}$ and let $%
\theta $ be an inner function with Clark measure $\sigma $. Set $\nu _{\mu
,\theta }\equiv \left\vert 1-\theta \right\vert ^{2}\mu $. Then the
following four conditions are equivalent (the equivalence of the first three
conditions is in \cite[Theorem 1.15]{LaSaShUrWi}):

\begin{enumerate}
\item $\mu $ is a Carleson measure for $K_{\theta }$,

\item The Cauchy transform $\mathcal{C}:L^{2}\left( \sigma \right)
\rightarrow L^{2}\left( \nu _{\mu ,\theta }\right) $ is bounded,

\item The Muckenhoupt and $T1$ testing conditions in \cite[(1.9), (1.10) and
(1.11)]{LaSaShUrWi} hold,

\item The Muckenhoupt, energy and $Tb$ testing conditions in (\ref{b and b* 
test}) hold.
\end{enumerate}
\end{thm}

An application of the Carleson measure property for $K_{\theta }$ was also
pointed out in \cite[Theorem 1.16]{LaSaShUrWi}, namely that if $\tau $ is a
positive Borel measure on $\mathbb{D}$ and $\varphi :\mathbb{D}\rightarrow 
\mathbb{D}$ is holomorphic, then the composition map $C_{\varphi }f\equiv
f\circ \varphi $ is bounded from the model space $K_{\theta }$ to the
weighted Hardy space $H_{\tau }^{2}$ if and only if the pushforward measure $%
\varphi _{\ast }\tau $ is a Carleson measure for $K_{\theta }$.

\subsection{History diagram and open problems}

Here is a list of open questions.

\begin{enumerate}
\item The most difficult and important problem in the theory of $T1$ and $Tb$
arises from the fact that, while the Muckenhoupt buffer conditions are
necessary for boundedness of a wide range of singular integrals, the \emph{%
energy} buffer conditions are only necessary for boundedness of the Hilbert
transform and some perturbations in dimension $n=1$, see \cite{Saw3}, \cite%
{SaShUr11}. What is a reasonable substitute for the energy buffer conditions
in a $T1$ or $Tb$ theorem?

\item Does Theorem \ref{thm brief} remain true in the case $p=2$, i.e. when $%
\mathbf{b}=\left\{ b_{Q}\right\} _{Q\in \mathcal{P}}$ is a $2$-weakly $%
\sigma $-accretive family of functions, and $\mathbf{b}^{\ast }=\left\{
b_{Q}^{\ast }\right\} _{Q\in \mathcal{P}}$ is a $2$-weakly $\omega $%
-accretive family of functions?

\item In the special case of the Hilbert transform in dimension $n=1$, are
the energy conditions in Theorem \ref{thm brief} already implied by the
Muckenhoupt, $\mathbf{b}$-testing and dual $\mathbf{b}^{\ast }$-testing
conditions for a pair of $p$-weakly accretive families, $p>2$?
\end{enumerate}

We end the introduction with a diagram detailing the relevant history of two weight theory for
this paper. Many important contributions are omitted, such as those dealing
with $L^{p},L^{q}$ assumptions in the case of Lebesgue measure, see for
example \cite{Hof1} and references there, and results for dyadic operators,
see for example \cite{AuHoMuTaTh} and references there. As is evident from
the diagram, the result of this paper (and its precursor for $n=1$) is the
first local $Tb$ theorem for two weights.

\begin{center}
\includegraphics[width=\textwidth]{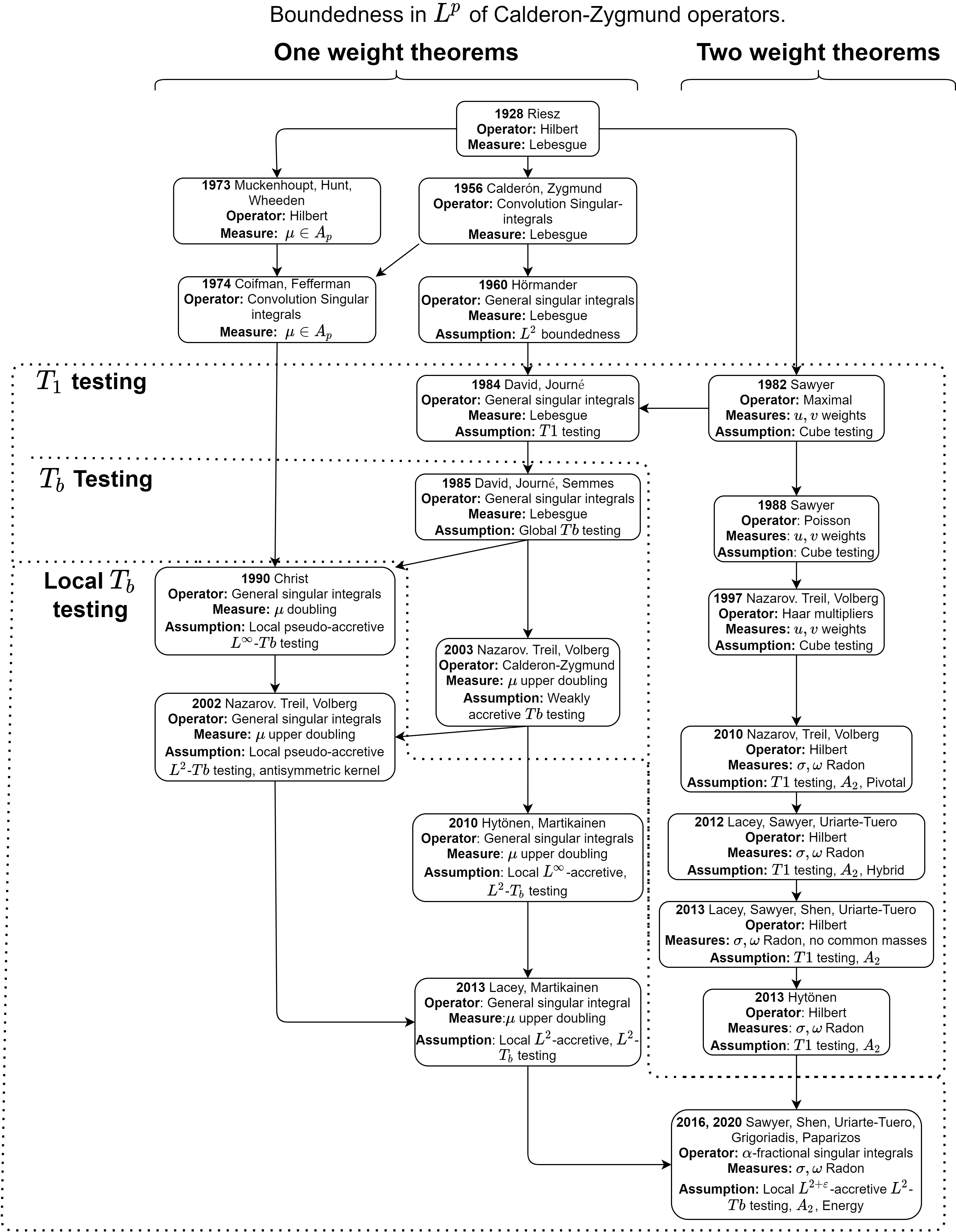}
\end{center}

\section{The local $Tb$ theorem and proof preliminaries}

\subsection{Standard fractional singular integrals}

Let $0\leq \alpha <n$. We define a standard $\alpha $-fractional CZ kernel $%
K^{\alpha }(x,y)$ to be a real-valued function defined on $\mathbb{R}%
^{n}\times \mathbb{R}^{n}$ satisfying the following fractional size and
smoothness conditions of order $1+\delta $ for some $\delta >0$: For $x\neq
y $,%
\begin{eqnarray}
\left\vert K^{\alpha }\left( x,y\right) \right\vert 
&\leq &\label{sizeandsmoothness'} 
C_{CZ}\left\vert
x-y\right\vert ^{\alpha -n}\\
\left\vert \nabla K^{\alpha }\left(
x,y\right) \right\vert 
&\leq&
C_{CZ}\left\vert x-y\right\vert ^{\alpha -n-1} \notag
\\
\left\vert \nabla K^{\alpha }\left( x,y\right) -\nabla K^{\alpha }\left(
x^{\prime },y\right) \right\vert &\leq &C_{CZ}\left( \frac{\left\vert
x-x^{\prime }\right\vert }{\left\vert x-y\right\vert }\right) ^{\delta
}\left\vert x-y\right\vert ^{\alpha -n-1},\ \ \ \ \ \frac{\left\vert
x-x^{\prime }\right\vert }{\left\vert x-y\right\vert }\leq \frac{1}{2}, 
\notag
\end{eqnarray}%
and the last inequality also holds for the adjoint kernel in which $x$ and $%
y $ are interchanged. We note that a more general definition of kernel has
only order of smoothness $\delta >0$, rather than $1+\delta $, but the use
of the Monotonicity and Energy Lemmas in arguments below involves first
order Taylor approximations to the kernel functions $K^{\alpha }\left( \cdot
,y\right) $.

\subsubsection{Defining the norm inequality}

We now turn to a precise definition of the weighted norm inequality%
\begin{equation}
\left\Vert T_{\sigma }^{\alpha }f\right\Vert _{L^{2}\left( \omega \right)
}\leq \mathfrak{N}_{T_{\sigma }^{\alpha }}\left\Vert f\right\Vert
_{L^{2}\left( \sigma \right) },\ \ \ \ \ f\in L^{2}\left( \sigma \right) .
\label{two weight'}
\end{equation}%
For this we introduce a family $\left\{ \eta _{\delta ,R}^{\alpha }\right\}
_{0<\delta <R<\infty }$ of nonnegative functions on $\left[ 0,\infty \right)$ so that the truncated kernels $K_{\delta ,R}^{\alpha }\left( x,y\right)
=\eta _{\delta ,R}^{\alpha }\left( \left\vert x-y\right\vert \right)
K^{\alpha }\left( x,y\right) $ are bounded with compact support for fixed $x$
or $y$. Then the truncated operators 
\begin{equation}
T_{\sigma ,\delta ,R}^{\alpha }f\left( x\right) \equiv \int_{\mathbb{R}^n
}K_{\delta ,R}^{\alpha }\left( x,y\right) f\left( y\right) d\sigma \left(
y\right) ,\ \ \ \ \ x\in \mathbb{R}^n, \label{def truncation}
\end{equation}
are pointwise well-defined, and we will refer to the pair $\left( K^{\alpha
},\left\{ \eta _{\delta ,R}^{\alpha }\right\} _{0<\delta <R<\infty }\right) $
as an $\alpha $-fractional singular integral operator, which we typically
denote by $T^{\alpha }$, suppressing the dependence on the truncations.

\begin{dfn}
We say that an $\alpha $-fractional singular integral operator $T^{\alpha} $ satisfies the norm inequality (\ref{two weight'})
provided%
\begin{equation*}
\left\Vert T_{\sigma ,\delta ,R}^{\alpha }f\right\Vert _{L^{2}\left( \omega
\right) }\leq \mathfrak{N}_{T_{\sigma }^{\alpha }}\left\Vert f\right\Vert
_{L^{2}\left( \sigma \right) },\ \ \ \ \ f\in L^{2}\left( \sigma \right)
,0<\delta <R<\infty .
\end{equation*}
\end{dfn}

It turns out that, in the presence of the Muckenhoupt conditions (\ref{def
A2}) below, the norm inequality (\ref{two weight'}) is essentially
independent of the choice of truncations used, and this is explained in some
detail in \cite{SaShUr10}. Thus, as in \cite{SaShUr10}, we are free to use
the tangent line truncations described there throughout the proofs of our
results.

\subsection{Weakly accretive functions}

Denote by $\mathcal{P}$ the collection of cubes in $\mathbb{R}^{n}$.
Note that we include an $L^{p}$ upper bound in our definition of `$p$-weakly
accretive family' of functions.

\begin{dfn}
Let $p\geq 2$ and let $\mu $ be a locally finite positive Borel measure on $%
\mathbb{R}^{n}$. We say that a family $\mathbf{b}=\left\{ b_{Q}\right\}
_{Q\in \mathcal{P}}$ of functions indexed by $\mathcal{P}$ is a $p$\emph{-weakly }$\mu $\emph{-accretive} family of functions on $\mathbb{R}^{n}$ if for $Q\in\mathcal{P}$,
\begin{eqnarray}
&&
\supp b_{Q}\subset Q 
 \notag\\
\label{acc} 0 < c_{\mathbf{b}} &\leq&  \frac{1}{\left\vert Q\right\vert _{\mu }}%
\int_{Q}b_{Q}d\mu  \leq \left( \frac{1}{\left\vert Q\right\vert
_{\mu }}\int_{Q}\left\vert b_{Q}\right\vert ^{p}d\mu \right) ^{\frac{1}{p}%
}\leq C_{\mathbf{b}}<\infty.
\end{eqnarray}
\end{dfn}

\subsection{b-testing conditions}

Suppose $\sigma $ and $\omega $ are locally finite positive Borel measures
on $\mathbb{R}^n$. The $\mathbf{b}$-testing conditions for $T^{\alpha }$ and $%
\mathbf{b}^{\ast }$-testing conditions for the dual $T^{\alpha ,\ast }$ are
given by%
\begin{eqnarray}
\int_{Q}\left\vert T_{\sigma }^{\alpha }b_{Q}\right\vert ^{2}d\omega &\leq
&\left( \mathfrak{T}_{T^{\alpha }}^{\mathbf{b}}\right) ^{2}\left\vert
Q\right\vert _{\sigma }\ ,\ \ \ \ \ \text{for all cubes }Q,
\label{b testing cond} \\
\int_{Q}\left\vert T_{\omega }^{\alpha ,\ast }b_{Q}^{\ast }\right\vert
^{2}d\sigma &\leq &\left( \mathfrak{T}_{T^{\alpha,\ast  }}^{\mathbf{b}^{\ast
}}\right) ^{2}\left\vert Q\right\vert _{\omega }\ ,\ \ \ \ \ \text{for
all cubes }Q.  \notag
\end{eqnarray}%

\subsection{Poisson integrals and the Muckenhoupt conditions 
}\label{def Poisson}
Let $\mu$ be a locally finite positive Borel measure on $\mathbb{R}^{n}$,
and suppose $Q$ is a cube in $\mathbb{R}^{n}$. Recall that $\left\vert
Q\right\vert^{\frac{1}{n}} =\ell \left( Q\right) $ for a cube $Q$. The two $\alpha $%
-fractional Poisson integrals of $\mu $ on a cube $Q$ are given by the
following expressions:%
\begin{eqnarray*}
\mathrm{P}^{\alpha }\left( Q,\mu \right) &\equiv &\int_{\mathbb{R}^n}\frac{%
\left\vert Q\right\vert ^{\frac{1}{n}}}{\left( \left\vert Q\right\vert ^{%
\frac{1}{n}}+\left\vert x-x_{Q}\right\vert \right) ^{n+1-\alpha }}d\mu
\left( x\right) , \\
\mathcal{P}^{\alpha }\left( Q,\mu \right) &\equiv &\int_{\mathbb{R}^n}\left( 
\frac{\left\vert Q\right\vert ^{\frac{1}{n}}}{\left( \left\vert Q\right\vert
^{\frac{1}{n}}+\left\vert x-x_{Q}\right\vert \right) ^{2}}\right) ^{n-\alpha
}d\mu \left( x\right) ,
\end{eqnarray*}%
where $\left\vert x-x_{Q}\right\vert $ denotes distance between $x$ and the center $x_{Q}$ of $Q$ and $\left\vert Q\right\vert $ denotes the Lebesgue measure of the
cube $Q$. We refer to $\mathrm{P}^{\alpha }$ as the \emph{standard} Poisson
integral and to $\mathcal{P}^{\alpha }$ as the \emph{reproducing} Poisson
integral. Note that these two kernels satisfy for all cubes $Q$ and positive measures $\mu$,
\begin{eqnarray*}
0 &\leq &\mathrm{P}^{\alpha }\left( Q,\mu \right) \leq C\mathcal{P}^{\alpha
}\left( Q,\mu \right) ,\ \ \ \ n-1\leq \alpha <n, \\
0 &\leq &\mathcal{P}^{\alpha }\left( Q,\mu \right) \leq C\mathrm{P}^{\alpha
}\left( Q,\mu \right) ,\ \ \ \ 0\leq \alpha <n-1.
\end{eqnarray*}

We now define the \emph{one-tailed} \emph{constant with holes }$\mathcal{A}%
_{2}^{\alpha }$ using the reproducing Poisson kernel $\mathcal{P}^{\alpha }$%
. On the other hand, the standard Poisson integral $\mathrm{P}^{\alpha }$
arises naturally throughout the proof of the $Tb$ theorem in estimating
oscillation of the fractional singular integral $T^{\alpha }$, and in the
definition of the energy conditions below.

\begin{dfn}\label{def call A2}
Suppose $\sigma $ and $\omega $ are locally finite positive Borel measures
on $\mathbb{R}^{n}$. The one-tailed constants $\mathcal{A}_{2}^{\alpha }$
and $\mathcal{A}_{2}^{\alpha ,\ast }$ with holes for the weight pair $\left(
\sigma ,\omega \right) $ are given by%
\begin{eqnarray*}
\mathcal{A}_{2}^{\alpha } &\equiv &\sup_{Q\in \mathcal{P}}\mathcal{P}%
^{\alpha }\left( Q,\mathbf{1}_{Q^{c}}\sigma \right) \frac{\left\vert
Q\right\vert _{\omega }}{\left\vert Q\right\vert ^{1-\frac{\alpha }{n}}}%
<\infty , \\
\mathcal{A}_{2}^{\alpha ,\ast } &\equiv &\sup_{Q\in \mathcal{P}}\mathcal{P}%
^{\alpha }\left( Q,\mathbf{1}_{Q^{c}}\omega \right) \frac{\left\vert
Q\right\vert _{\sigma }}{\left\vert Q\right\vert ^{1-\frac{\alpha }{n}}}%
<\infty .
\end{eqnarray*}
\end{dfn}

Note that these definitions are the conditions with `holes' introduced by Hyt%
\"{o}nen \cite{Hyt} - the supports of the measures $\mathbf{1}_{Q^{c}}\sigma$ and $\mathbf{1}_{Q^c}\omega $ in the definition of $\mathcal{A}_{2}^{\alpha
} $ are disjoint, and so any common point masses of $\sigma $ and $\omega $
do not appear simultaneously in the factors of any of the products $\mathcal{%
P}^{\alpha }\left( Q,\mathbf{1}_{Q^{c}}\sigma \right) \frac{\left\vert
Q\right\vert _{\omega }}{\left\vert Q\right\vert ^{1-\frac{\alpha}{n} }}$. Recall, the definition of the classical Muckenhoupt condition
$$
A_2^\alpha= \sup_{Q\in \mathcal{P}
}\frac{\left\vert Q\right\vert _{\omega }}{\left\vert Q\right\vert ^{1-\frac{\alpha }{n}}}\frac{\left\vert Q\right\vert _{\sigma }}{\left\vert
Q\right\vert ^{1-\frac{\alpha }{n}}}
$$
but it will find no use in the two weight setting with common point masses permitted.

Initially, these definitions of Muckenhoupt type were given in the following
`one weight' case, $d\omega \left( x\right) =w\left( x\right) dx$ and $%
d\sigma \left( x\right) =\frac{1}{w\left( x\right) }dx$, where $\mathcal{A}%
_{2}^{\alpha }\left( \lambda w,\left( \lambda w\right) ^{-1}\right) =%
\mathcal{A}_{2}^{\alpha }\left( w,w^{-1}\right) $ is homogeneous of degree $%
0 $. Of course the two weight version is homogeneous of degree $2$ in the
weight pair, $\mathcal{A}_{2}^{\alpha }\left( \lambda \sigma ,\lambda \omega
\right) =\lambda ^{2}\mathcal{A}_{2}^{\alpha }\left( \sigma ,\omega \right) $%
, while all of the other conditions we consider in connection with two
weight norm inequalities, including the operator norm $\mathfrak{N}%
_{T^{\alpha }}\left( \sigma ,\omega \right) $ itself, are homogeneous of
degree $1$ in the weight pair. This awkwardness regarding the homogeneity of
Muckenhoupt conditions could be rectified by simply taking the square root
of $\mathcal{A}_{2}^{\alpha }$ and renaming it, but the current definition
is so entrenched in the literature, in particular in connection with the $%
A_{2}$ conjecture, that we will leave it as is.

\subsubsection{Punctured $A_{2}^{\protect\alpha }$ conditions} \label{def punct}

The \emph{classical} $A_{2}^{\alpha }$ characteristic  fails to be finite when the measures $%
\sigma $ and $\omega $ have a common point mass - simply let $Q$ in the $%
\sup $ above shrink to a common mass point. But there is a substitute that
is quite similar in character that is motivated by the fact that for large
cubes $Q$, the $\sup $ above is problematic only if just \emph{one} of the
measures is \emph{mostly} a point mass when restricted to $Q$.

Given an at most countable set $\mathfrak{P}=\left\{ p_{k}\right\}
_{k=1}^{\infty }$ in $\mathbb{R}^n$, a cube $Q\in \mathcal{P}$, and a
positive locally finite Borel measure $\mu $, define 
\begin{equation}\label{puncture}
\mu \left( Q,\mathfrak{P}\right) \equiv \left\vert Q\right\vert _{\mu }-\sup
\left\{ \mu \left( p_{k}\right) :p_{k}\in Q\cap \mathfrak{P}\right\} ,
\end{equation}%
where the supremum is actually achieved since $\sum_{p_{k}\in Q\cap 
\mathfrak{P}}\mu \left( p_{k}\right) <\infty $ as $\mu $ is locally finite.
The quantity $\mu \left( Q,\mathfrak{P}\right) $ is simply the $\widetilde{%
\mu }$ measure of $Q$ where $\widetilde{\mu }$ is the measure $\mu $ with
its largest point mass from $\mathfrak{P}$ in $Q$ removed. Given a locally
finite positive measure pair $\left( \sigma ,\omega \right) $, let $%
\mathfrak{P}_{\left( \sigma ,\omega \right) }=\left\{ p_{k}\right\}
_{k=1}^{\infty }$ be the at most countable set of common point masses of $%
\sigma $ and $\omega $. Then the weighted norm inequality (\ref{two weight'}%
) typically implies finiteness of the following \emph{punctured} Muckenhoupt
conditions:%
\begin{eqnarray*}
A_{2}^{\alpha ,{punct}}\left( \sigma ,\omega \right) &\equiv
&\sup_{Q\in \mathcal{P}}\frac{\omega \left( Q,\mathfrak{P}_{\left( \sigma
,\omega \right) }\right) }{\left\vert Q\right\vert ^{1-\frac{\alpha }{n}}}%
\frac{\left\vert Q\right\vert _{\sigma }}{\left\vert Q\right\vert ^{1-\frac{%
\alpha }{n}}}, \\
A_{2}^{\alpha ,\ast ,{punct}}\left( \sigma ,\omega \right) &\equiv
&\sup_{Q\in \mathcal{P}}\frac{\left\vert Q\right\vert _{\omega }}{\left\vert
Q\right\vert ^{1-\frac{\alpha }{n}}}\frac{\sigma \left( Q,\mathfrak{P}%
_{\left( \sigma ,\omega \right) }\right) }{\left\vert Q\right\vert ^{1-\frac{%
\alpha }{n}}}.
\end{eqnarray*}%
In particular, all of the above Muckenhoupt conditions $\mathcal{A}%
_{2}^{\alpha }$, $\mathcal{A}_{2}^{\alpha ,\ast }$, $A_{2}^{\alpha ,{%
punct}}$ and $A_{2}^{\alpha ,\ast ,{punct}}$ are necessary for
boundedness of an elliptic $\alpha $-fractional singular integral $T_{\sigma
}^{\alpha }$ from\thinspace $L^{2}\left( \sigma \right) $ to $L^{2}\left(
\omega \right) $. It is convenient to define%
\begin{equation}
\mathfrak{A}_{2}^{\alpha }\equiv \mathcal{A}_{2}^{\alpha }+\mathcal{A}%
_{2}^{\alpha ,\ast }+A_{2}^{\alpha ,{punct}}+A_{2}^{\alpha ,\ast ,%
{punct}}\ .  \label{def A2}
\end{equation}

\subsection{Energy Conditions}

Here is the definition of the strong energy conditions, which we sometimes
refer to simply as the energy conditions. Let
$$
m_I^\mu\equiv\frac{1}{|I|_\mu}\int x d\mu(x)=\left\langle \frac{1}{|I|_\mu}\int x_1d\mu(x),...,\frac{1}{|I|_\mu}\int x_nd\mu(x)\right\rangle
$$
be the  average of $x$ with respect to the measure $\mu$, which we often abbreviate to $m_I$ when the measure $\mu$ is understood.

\begin{dfn}
\label{def strong quasienergy}Let $0\leq \alpha <n$. Suppose $\sigma $ and $%
\omega $ are locally finite positive Borel measures on $\mathbb{R}^{n}$.
Then the \emph{strong} energy constant $\mathcal{E}_{2}^{\alpha }$ is
defined by 
\begin{equation}
\left( \mathcal{E}_{2}^{\alpha }\right) ^{2}\equiv \sup_{I=\dot{\cup}I_{r}}%
\frac{1}{\left\vert I\right\vert _{\sigma }}\sum_{r=1}^{\infty }\left( \frac{%
\mathrm{P}^{\alpha }\left( I_{r},\mathbf{1}_{I}\sigma \right) }{\left\vert
I_{r}\right\vert ^{\frac{1}{n}}}\right) ^{2}\left\Vert
x-m^\omega_{I_{r}}\right\Vert _{L^{2}\left( \mathbf{1}_{I_{r}}\omega \right) }^{2}\
,  \label{strong b* energy}
\end{equation}%
where the supremum is taken over arbitrary decompositions of a cube $I$
using a pairwise disjoint union of subcubes $I_{r}$. Similarly, we define
the dual \emph{strong} energy constant $\mathcal{E}_{2 }^{\alpha,\ast }$ by switching the roles of $\sigma $ and $\omega $:%
\begin{equation}
\left( \mathcal{E}_{2}^{\alpha ,\ast }\right) ^{2}\equiv \sup_{I=\dot{\cup}%
I_{r}}\frac{1}{\left\vert I\right\vert _{\omega }}\sum_{r=1}^{\infty }\left( 
\frac{\mathrm{P}^{\alpha }\left( I_{r},\mathbf{1}_{I}\omega \right) }{%
\left\vert I_{r}\right\vert ^{\frac{1}{n}}}\right) ^{2}\left\Vert
x-m^\sigma_{I_{r}}\right\Vert _{L^{2}\left( \mathbf{1}_{I_{r}}\sigma \right) }^{2}\
.  \label{strong b energy}
\end{equation}
\end{dfn}

These energy conditions are necessary for
boundedness of  elliptic and gradient elliptic operators, including the Hilbert transform (but not for for certain elliptic singular operators that fail to be gradient elliptic) - see \cite{SaShUr11} and \cite{SaShUr12}. It is convenient to define
\begin{equation*}
\mathfrak{E}_{2}^{\alpha }\equiv \mathcal{E}_{2}^{\alpha }+\mathcal{E}_{2}^{\alpha ,\ast }
\end{equation*}%
as well as
\begin{equation}
\mathcal{NTV}_{\alpha }\equiv \mathfrak{T}_{T^{\alpha }}^{\mathbf{b}}+%
\mathfrak{T}_{T^{\alpha,\ast }}^{\mathbf{b}^{\ast } }+\sqrt{\mathfrak{A}
_{2}^{\alpha }}+\mathfrak{E}_{2}^{\alpha }\ .  \label{def NTV}
\end{equation}

\subsection{The two weight local $Tb$ Theorem}
Here we derive a higher dimensional local $Tb$ theorem based in part on the 
\emph{proof} of the one-dimensional analogue in \cite{SaShUr12}, \ which was
in turn based in part on the \emph{proof} of the $T1$ theorem in \cite%
{SaShUr7}, and in part on the \emph{proof} of a one weight $Tb$ theorem in
Hyt\"{o}nen and Martikainen \cite{HyMa}.
\begin{thm}
\label{dim high}Suppose that $\sigma $ and $\omega $ are locally finite
positive Borel measures on Euclidean space $\mathbb{R}^{n}$. Suppose that $%
T^{\alpha }$ is a standard $\alpha $-fractional singular integral operator
on $\mathbb{R}^{n}$, and set $T_{\sigma }^{\alpha }f=T^{\alpha }\left(
f\sigma \right) $ for any smooth truncation of $T_{\sigma }^{\alpha }$, so
that $T_{\sigma }^{\alpha }$ is \emph{apriori} bounded from $L^{2}\left(
\sigma \right) $ to $L^{2}\left( \omega \right) $. Assume the Muckenhoupt
and energy conditions hold, i.e. $\mathcal{A}_{2}^{\alpha },\mathcal{A}%
_{2}^{\alpha ,\ast },A_{2}^{\alpha ,{punct}},A_{2}^{\alpha ,\ast ,%
{punct}},\mathcal{E}_{2}^{\alpha },\mathcal{E}_{2}^{\alpha ,\ast
}<\infty $. Finally, let $p>2$ and let $\mathbf{b}=\left\{ b_{Q}\right\}
_{Q\in \mathcal{P}}$ be a $p$-weakly $\sigma$-accretive family of functions
on $\mathbb{R}^n$, and let $\mathbf{b}^{\ast }=\left\{ b_{Q}^{\ast }\right\}
_{Q\in \mathcal{P}}$ be a $p$-weakly $\omega $-accretive family of functions
on $\mathbb{R}^n$. Then for $0\leq \alpha <n$, the operator $T_{\sigma
}^{\alpha }$ is bounded from $L^{2}\left( \sigma \right) $ to $L^{2}\left(
\omega \right) $ with operator norm $\mathfrak{N}_{T_{\sigma }^{\alpha }}$,
i.e. 
\begin{equation*}
\left\Vert T_{\sigma }^{\alpha }f\right\Vert _{L^{2}\left( \omega \right)
}\leq \mathfrak{N}_{T_{\sigma }^{\alpha }}\left\Vert f\right\Vert
_{L^{2}\left( \sigma \right) },\ \ \ \ \ f\in L^{2}\left( \sigma \right) ,
\end{equation*}%
\textbf{uniformly} in smooth truncations of $T^{\alpha }$ \emph{if and only
if} the $\mathbf{b}$-testing conditions for $T^{\alpha }$ and the $\mathbf{b}%
^{\ast }$-testing conditions for the dual $T^{\alpha ,\ast }$ both hold.
Moreover, we have
\begin{equation*}
\mathfrak{N}_{T^{\alpha }}\lesssim \mathfrak{T}_{T^{\alpha }}^{\mathbf{b}}+%
\mathfrak{T}_{T^{\alpha }}^{\mathbf{b}^{\ast }}+\sqrt{\mathfrak{A}%
_{2}^{\alpha }}+\mathfrak{E}_{2}^{\alpha }\ .
\end{equation*}
\end{thm}

\begin{rem}
\label{special}In the special case that $\sigma =\omega =\mu $, the
classical Muckenhoupt $A_{2}^{\alpha }$ condition is%
\begin{equation*}
\sup_{Q\in \mathcal{P}}\frac{\left\vert Q\right\vert _{\mu }}{\left\vert
Q\right\vert ^{1-\frac{\alpha}{n} }}\frac{\left\vert Q\right\vert _{\mu }}{\left\vert
Q\right\vert ^{1-\frac{\alpha}{n} }}<\infty ,
\end{equation*}%
which is the upper doubling measure\ condition with exponent $n-\alpha $,
i.e. 
\begin{equation*}
\left\vert Q\right\vert _{\mu }\leq C\ell \left( Q\right) ^{n-\alpha },\ \ \ \text{for all cubes }Q,
\end{equation*}%
which of course prohibits point masses in $\mu $. Both Poisson integrals are
then bounded, 
\begin{eqnarray*}
\mathrm{P}^{\alpha }\left( Q,\mu \right) 
\!\!\!\!\!&\lesssim &\!\!\!\!\!
\sum_{k=0}^{\infty }
\frac{\left\vert Q\right\vert^\frac{1}{n} }{\left( 2^{k}\left\vert Q\right\vert^\frac{1}{n} \right)
^{n+1-\alpha }}\left\vert 2^{k}Q\right\vert _{\mu }
\lesssim
\sum_{k=0}^{\infty}\frac{\left\vert Q\right\vert^\frac{1}{n} }{\left( 2^{k}\left\vert Q\right\vert^\frac{1}{n} \right)
^{n+1-\alpha }}\left( 2^{k}\ell(Q) \right) ^{n-\alpha }
\!\!=2 \\
\ \ \mathcal{P}^{\alpha }\left( Q,\mu \right) 
\!\!\!\!\!&\lesssim &\!\!\!\!\!
\sum_{k=0}^{\infty}\left( \!\frac{\left\vert Q\right\vert^\frac{1}{n} }{\left( 2^{k}\left\vert Q\right\vert^\frac{1}{n}
\right) ^{2}}\!\right) ^{\!\!\! n-\alpha }\!\!\!\!\left\vert 2^{k}Q\right\vert _{\mu}
\lesssim
\sum_{k=0}^{\infty}\left( \!\frac{\left\vert Q\right\vert^\frac{1}{n} }{\left( 2^{k}\left\vert Q\right\vert^\frac{1}{n}
\right) ^{2}}\!\right) ^{\!\!\! n-\alpha }\!\!\!\!\left(
2^{k} \ell(Q) \right) ^{n-\alpha }\!\!=C_{\alpha }
\end{eqnarray*}
and it follows easily that the equal weight pair $\left( \mu ,\mu \right) $
satisfies not only the Muckenhoupt $\mathfrak{A}_{2}^{\alpha }$ condition,
but also the strong energy condition $\mathfrak{E}_{2}^{\alpha }$:%
\begin{eqnarray*}
\sum_{r=1}^{\infty }\left( \frac{\mathrm{P}^{\alpha }\left( I_{r},\mathbf{1}_{I}\sigma \right) }{\left\vert I_{r}\right\vert }\right) ^{2}\left\Vert x-m_{I_{r}}^{\omega }\right\Vert _{L^{2}\left( \omega \right)
}^{2}
&\leq&
C\sum_{r=1}^{\infty }\left\Vert \frac{x-m_{I_{r}}^{\omega
}}{\left\vert I_{r}\right\vert }\right\Vert _{L^{2}\left( \omega \right)}^{2}\leq
C\sum_{r=1}^{\infty }\left\vert I_{r}\right\vert _{\omega }\leq C\left\vert
I\right\vert _{\omega }=C\left\vert I\right\vert _{\sigma }\ ,
\end{eqnarray*}%
since $\omega =\sigma $. Thus Theorem \ref{dim high}, when restricted to a
single weight $\sigma =\omega $, recovers a slightly weaker,  due to our assumption that $p>2$, version of the
one weight theorem of Lacey and Martikainen \cite[Theorem 1.1]{LaMa} for
dimension $n=1$. On the other hand,
the possibility of a two weight theorem for a $2$-weakly $\mu $-accretive
family is highly problematic, as one of the key proof strategies used in 
\cite{LaMa} in the one weight case is a reduction to testing over $f$ and $g$
with controlled $L^{\infty }$ norm, a strategy that appears to be
unavailable in the two weight setting.
\end{rem}

In order to prove Theorem \ref{dim high}, it is
convenient to establish some improved properties for our $p$-weakly $\mu $%
-accretive family, and also necessary to establish some improved energy
conditions related to the families of testing functions $\mathbf{b}$ and $%
\mathbf{b}^{\ast }$. We turn to these matters in the next two subsections.

\begin{rem}
We alert the reader to the fact that a large portion of the argument
presented below originated in \cite{SaShUr12} \ in the case $n=1$, but that
significant differences arise in various places throughout, especially as
outlined in the introduction. As a consequence we repeat the arguments from 
\cite{SaShUr12} without further mention when needed.
\end{rem}

\subsection{Reduction to real bounded accretive families\label{str acc}}

We begin by noting that if $b_{Q}$ satisfies (\ref{acc}) with $\mu =\sigma $%
, and satisfies a given $\mathbf{b}$-testing condition for a weight pair $%
\left( \sigma ,\omega \right) $, then ${Re}b_{Q}$ satisfies 
$$
\left( 
\frac{1}{\left\vert Q\right\vert _{\mu }}\int_{Q}\left\vert {Re}%
b_{Q}\right\vert ^{p}d\mu \right) ^{\frac{1}{p}}\leq C_{\mathbf{b}}\left( p\right) 
$$
and the given $\mathbf{b}$-testing condition for $\left( \sigma
,\omega \right) $ with ${Re}b_{Q}$ in place of $b_{Q}$.

Thus we may assume throughout the proof of Theorem \ref{dim high} that our $p$%
-weakly\emph{\ }$\mu $-accretive families $\mathbf{b}\equiv \left\{
b_{Q}\right\} _{Q\in \mathcal{D}}$ and $\mathbf{b}^{\ast }\equiv \left\{
b_{Q}^{\ast }\right\} _{Q\in \mathcal{G}}$ consist of \textbf{real-valued} functions.

Next we show that the assumption of testing conditions for a fractional singular integral $T^\alpha$ and $p$-weakly $\mu $-accretive testing functions $%
\mathbf{b}=\left\{ b_{Q}\right\} _{Q\in \mathcal{P}}$ and $\mathbf{b}^{\ast
}=\left\{ b_{Q}^{\ast }\right\} _{Q\in \mathcal{P}}$ with $p>2$ can always
be replaced with real-valued $\infty $-weakly $\mu $-accretive testing
functions, thus reducing the $Tb$ theorem for the case $p>2$ to the case
when $p=\infty $. We now proceed to develop
a precise statement. We extend (\ref{acc}) to $2< p\leq \infty $ by%
\begin{eqnarray}
&&
\supp b_{Q}\subset Q\ ,\ \ \ \ \ Q\in \mathcal{P},
\label{acc infinity} \\
1 &\leq & 
\frac{1}{\left\vert Q\right\vert _{\mu }}
\int_{Q}b_{Q}d\mu  \leq \left\{ 
\begin{array}{cc}
\left( \frac{1}{\left\vert Q\right\vert _{\mu }}\int_{Q}\left\vert
b_{Q}\right\vert ^{p}d\mu \right) ^{\frac{1}{p}}\leq C_{\mathbf{b}}\left(
p\right) <\infty        & \text{for }2< p<\infty \\ 
\left\Vert b_{Q}\right\Vert _{L^{\infty }\left( \mu \right) }\leq C_{\mathbf{b}}\left( \infty \right) <\infty                 & \text{for }p=\infty
\end{array}
\right.   \notag
\end{eqnarray}

\begin{prop}
\label{conditional}Let $0\leq \alpha <1$, and let $\sigma $ and $\omega $ be
locally finite positive Borel measures on $\mathbb{R}^n$, and
let $T^{\alpha }$ be a standard $\alpha $-fractional elliptic and gradient
elliptic singular integral operator on $\mathbb{R}^n$. Set $T_{\sigma
}^{\alpha }f=T^{\alpha }\left( f\sigma \right) $ for any smooth truncation
of $T_{\sigma }^{\alpha }$, so that $T_{\sigma }^{\alpha }$ is \emph{apriori}
bounded from $L^{2}\left( \sigma \right) $ to $L^{2}\left( \omega \right) $.
Finally, define the sequence of positive extended real numbers%
\begin{equation*}
\left\{ p_{m}\right\} _{m=0}^{\infty }=\left\{ \frac{2}{1-\left( \frac{2}{3}%
\right) ^{m}}\right\} _{m=0}^{\infty }=\left\{ \infty ,6,\frac{18}{5},\frac{%
162}{65},...\right\} .
\end{equation*}%
Suppose that the following statement is true:\medskip

\begin{description}
\item[$\left( \mathcal{S}_{\infty }\right) $] If $\mathbf{b}=\left\{
b_{Q}\right\} _{Q\in \mathcal{P}}$ is an $\infty $-weakly $\sigma $%
-accretive family of functions on $\mathbb{R}^n$ and if $\mathbf{b}^{\ast
}=\left\{ b_{Q}^{\ast }\right\} _{Q\in \mathcal{P}}$ is an $\infty $%
-weakly $\omega $-accretive family of functions on $\mathbb{R}^n$, then the
operator norm $\mathfrak{N}_{T_{\sigma }^{\alpha }}$ of $T_{\sigma }^{\alpha
}$ from $L^{2}\left( \sigma \right) $ to $L^{2}\left( \omega \right) $, i.e.
the best constant in%
\begin{equation*}
\left\Vert T_{\sigma }^{\alpha }f\right\Vert _{L^{2}\left( \omega \right)
}\leq \mathfrak{N}_{T_{\sigma }^{\alpha }}\left\Vert f\right\Vert
_{L^{2}\left( \sigma \right) },\ \ \ \ \ f\in L^{2}\left( \sigma \right) ,
\end{equation*}%
\textbf{uniformly} in smooth truncations of $T^{\alpha }$, satisfies%
\begin{equation*}
 \mathfrak{N}_{T^{\alpha }}\lesssim \left( C_{\mathbf{b}}\left(
\infty \right) +C_{\mathbf{b}^{\ast }}\left( \infty \right) \right) \left( 
\mathfrak{T}_{T^{\alpha }}^{\mathbf{b}}+\mathfrak{T}_{T^{\alpha }}^{\mathbf{b%
}^{\ast }}+\sqrt{\mathfrak{A}_{2}^{\alpha }}+\mathfrak{E}_{2}^{\alpha
}\right) \ ,
\end{equation*}%
where $C_{\mathbf{b}}\left( \infty \right) ,C_{\mathbf{b}^{\ast }}\left(
\infty \right) $ are the accretivity constants in (\ref{acc infinity}), and
the constants implied by $\lesssim $ depend on $\alpha $ and the constant $C_{CZ}$ in (\ref{sizeandsmoothness'}).\\
Then for each $m\geq 0$, the following statements hold:\medskip

\item[$\left( \mathcal{S}_{m}\right) $] Let $p\in \left( p_{m+1},p_{m}\right]
$. If $\mathbf{b}=\left\{ b_{Q}\right\} _{Q\in \mathcal{P}}$ is a $p$%
-weakly $\sigma $-accretive family of functions on $\mathbb{R}^n$, and if $\mathbf{b}^{\ast }=\left\{ b_{Q}^{\ast }\right\} _{Q\in \mathcal{P}}$ is a $%
p $-weakly $\omega $-accretive family of functions on $\mathbb{R}^n$, then
the operator norm $\mathfrak{N}_{T_{\sigma }^{\alpha }}$ of $T_{\sigma
}^{\alpha }$ from $L^{2}\left( \sigma \right) $ to $L^{2}\left( \omega\right)$, \textbf{uniformly} in smooth truncations of $T^{\alpha }$, satisfies
\begin{equation*}
 \mathfrak{N}_{T^{\alpha }}\lesssim \left( C_{\mathbf{b}}\left(
p\right) +C_{\mathbf{b}^{\ast }}\left( p\right) \right) ^{3^{m+1}}\left( 
\mathfrak{T}_{T^{\alpha }}^{\mathbf{b}}+\mathfrak{T}_{T^{\alpha }}^{\mathbf{b%
}^{\ast }}+\sqrt{\mathfrak{A}_{2}^{\alpha }}+\mathfrak{E}_{2}^{\alpha
}\right) \ ,
\end{equation*}%
where $C_{\mathbf{b}}\left( p\right) ,C_{\mathbf{b}^{\ast }}\left( p\right) $
are the accretivity constants in (\ref{acc}), and the constants implied by $%
\lesssim $ depend on $p$, $\alpha $, and the constant $C_{CZ}$ in (\ref{sizeandsmoothness'}).
\end{description}
\end{prop}

\begin{proof}[Proof of Proposition \protect\ref{conditional}]
We will prove it by  induction. We first prove $\left( \mathcal{S}_{0}\right) $. So fix $p\in \left(
p_{1},p_{0}\right) =\left( 6,\infty \right) $, and let $\mathbf{b}=\left\{
b_{Q}\right\} _{Q\in \mathcal{P}}$ be a $p$-weakly $\sigma $-accretive
family of functions on $\mathbb{R}^n$, and let $\mathbf{b}^{\ast }=\left\{
b_{Q}^{\ast }\right\} _{Q\in \mathcal{P}}$ be a $p$-weakly $\omega $%
-accretive family of functions on $\mathbb{R}^n$. Let $0<\varepsilon <1$ (to
be chosen differently at various points in the argument below) and define%
\begin{equation}
\lambda =\lambda \left( \varepsilon \right) =\left( \frac{p}{p-2}C_{\mathbf{b%
}}\left( p\right) ^{p}\frac{1}{\varepsilon }\right) ^{\frac{1}{p-2}}
\label{lambda choice}
\end{equation}%
and a new collection of test functions,%
\begin{equation}
\widehat{b}_{Q}\equiv 2b_{Q}\left( \mathbf{1}_{\left\{ \left\vert
b_{Q}\right\vert \leq \lambda \right\} }+\frac{\lambda }{\left\vert
b_{Q}\right\vert }\mathbf{1}_{\left\{ \left\vert b_{Q}\right\vert >\lambda
\right\} }\right) ,\ \ \ \ \ Q\in \mathcal{P},  \label{new}
\end{equation}%
We compute
\begin{eqnarray*}
\int_{\left\{ \left\vert b_{Q}\right\vert >\lambda \right\} }\!\!\!\!\left\vert
b_{Q}\right\vert ^{2}d\sigma
\!\!\!\!\!&=&\!\!\!\!\!
\int_{\left\{ \left\vert b_{Q}\right\vert
>\lambda \right\} }\left[ \int_{0}^{\left\vert b_{Q}\right\vert }2tdt\right]
d\sigma \\
&=&\!\!\!\!\!
\int \int_{\left\{ \left( x,t\right) \in \mathbb{R}^n\times \left( 0,\infty
\right) :\max \left\{ t,\lambda \right\} <\left\vert b_{Q}\left( x\right)
\right\vert \right\} }2tdtd\sigma \left( x\right) \\
&=&\!\!\!\!\!
\int_{0}^{\lambda}\int_{\left\{ x\in \mathbb{R}^n:\lambda
<\left\vert b_{Q}\left( x\right) \right\vert \right\} }\!\!\!\!\!d\sigma \left(x\right) 2tdt
+
\int_{\lambda}^\infty \int_{\left\{ x\in
\mathbb{R}^n:t<\left\vert b_{Q}\left( x\right) \right\vert \right\} }\!\!\!\!\! d\sigma
\left( x\right) 2tdt \\
&=&
\!\!\!\!\!\lambda ^{2}\left\vert \left\{ \left\vert b_{Q}\right\vert >\lambda
\right\} \right\vert _{\sigma }+\int_{\lambda }^{\infty }\left\vert \left\{
\left\vert b_{Q}\right\vert >t\right\} \right\vert _{\sigma }2tdt,
\end{eqnarray*}%
and hence%
\begin{eqnarray}
\int_{\left\{ \left\vert b_{Q}\right\vert >\lambda \right\} }\left\vert
b_{Q}\right\vert ^{2}d\sigma &\leq &\lambda ^{2}\frac{1}{\lambda ^{p}}\left(
\int \left\vert b_{Q}\right\vert ^{p}d\sigma \right) +\int_{\lambda
}^{\infty }\frac{1}{t^{p}}\left( \int \left\vert b_{Q}\right\vert
^{p}d\sigma \right) 2tdt  \label{hence} \\
&\leq&
\left\{ \lambda ^{2-p}+\int_{\lambda }^{\infty }2t^{1-p}dt\right\} C_{\mathbf{b}}\left( p\right) ^{p}\left\vert Q\right\vert _{\sigma }  \notag \\
&=&
\frac{p}{p-2}\lambda ^{2-p}C_{\mathbf{b}}\left( p\right) ^{p}\left\vert
Q\right\vert _{\sigma }=\varepsilon \left\vert Q\right\vert _{\sigma }\ , 
\notag
\end{eqnarray}%
by (\ref{lambda choice}). Thus we have the lower bound,
\begin{eqnarray}
&& \label{low}\\ \notag
\left\vert \frac{1}{\left\vert Q\right\vert _{\sigma }}\int_{Q}\widehat{b}%
_{Q}d\sigma \right\vert 
&=&
2\left\vert \frac{1}{\left\vert Q\right\vert
_{\sigma }}\int_{Q}b_{Q}d\sigma -\frac{1}{\left\vert Q\right\vert _{\sigma }}
\int_{Q}b_{Q}\left( 1-\frac{\lambda }{\left\vert b_{Q}\right\vert }\right) 
\mathbf{1}_{\left\{ \left\vert b_{Q}\right\vert >\lambda \right\} }d\sigma
\right\vert  \\
&\geq &
2\left\vert \frac{1}{\left\vert Q\right\vert _{\sigma }}
\int_{Q}b_{Q}d\sigma \right\vert -2\left( \frac{1}{\left\vert Q\right\vert
_{\sigma }}\int_{Q}\left\vert b_{Q}\right\vert ^{2}\mathbf{1}_{\left\{\left\vert b_{Q}\right\vert >\lambda \right\} }d\sigma \right) ^{\frac{1}{2}}
\notag \\
&\geq &
2-2\left( \frac{1}{\left\vert Q\right\vert _{\sigma }}\varepsilon
\left\vert Q\right\vert _{\sigma }\right) ^{\frac{1}{2}}=2-2\sqrt{\varepsilon }
\geq 1>0,\ \ \ \ \ Q\in \mathcal{P}.  \notag
\end{eqnarray}
For an upper bound we have 
\begin{equation*}
\left\Vert \widehat{b}_{Q}\right\Vert _{L^{\infty }\left( \sigma \right)
}\leq 2\lambda =2\lambda \left( \varepsilon \right) =2\left( \frac{p}{p-2}C_{%
\mathbf{b}}\left( p\right) ^{p}\frac{1}{\varepsilon }\right) ^{\frac{1}{p-2}%
},
\end{equation*}%
which altogether shows that%
\begin{eqnarray}\label{cinfty}
C_{\widehat{\mathbf{b}}}\left( \infty \right) 
&\leq &
2\left( \frac{p}{p-2}C_{\mathbf{b}}\left( p\right) ^{p}\frac{1}{\varepsilon }\right) ^{\frac{1}{p-2}
}=2\left( \frac{p}{p-2}\right) ^{\frac{1}{p-2}}C_{\mathbf{b}}\left( p\right)
^{\frac{p}{p-2}}\varepsilon ^{-\frac{1}{p-2}}
\end{eqnarray}%
if we choose $0<\varepsilon \leq \frac{1}{4}$. Similarly we have
\begin{eqnarray*}
C_{\widehat{\mathbf{b}}^{\ast }}\left( \infty \right) &\leq &
2\left( \frac{p}{p-2}C_{\mathbf{b}^{\ast }}\left( p\right) ^{p}\frac{1}{\varepsilon ^{\ast }%
}\right) ^{\frac{1}{p-2}}=2\left( \frac{p}{p-2}\right) ^{\frac{1}{p-2}}C_{\mathbf{b}^{\ast }}\left( p\right) ^{\frac{p}{p-2}}\left( \varepsilon ^{\ast
}\right) ^{-\frac{1}{p-2}}
\end{eqnarray*}
for $0<\varepsilon ^{\ast}\leq \frac{1}{4}$. 
Moreover, we also have, using (\ref{hence}),
\begin{eqnarray*}
\sqrt{\int_{Q}\left\vert T_{\sigma }^{\alpha }\widehat{b}_{Q}\right\vert
^{2}d\omega }
\!\!\!&\leq &\!\!\!
2\sqrt{\int_{Q}\left\vert T_{\sigma }^{\alpha
}b_{Q}\right\vert ^{2}d\omega }+2\sqrt{\int_{Q}\left\vert T_{\sigma
}^{\alpha }\mathbf{1}_{\left\{ \left\vert b_{Q}\right\vert >\lambda \right\}
}\left( \frac{\lambda }{\left\vert b_{Q}\right\vert }-1\right)
b_{Q}\right\vert ^{2}d\omega } \\
&\leq &
2\mathfrak{T}_{T^{\alpha }}^{\mathbf{b}}\sqrt{\left\vert Q\right\vert
_{\sigma }}+2\mathfrak{N}_{T^{\alpha }}\sqrt{\int_{\left\{ \left\vert
b_{Q}\right\vert >\lambda \right\} }\left\vert b_{Q}\right\vert ^{2}d\sigma }\\
&\leq&
2\left\{ \mathfrak{T}_{T^{\alpha }}^{\mathbf{b}}+\sqrt{\varepsilon }%
\mathfrak{N}_{T^{\alpha }}\right\} \sqrt{\left\vert Q\right\vert _{\sigma }}%
\ ,\ \ \ \ \ \text{for all cubes }Q,
\end{eqnarray*}
which shows that
\begin{equation}
\mathfrak{T}_{T^{\alpha }}^{\widehat{\mathbf{b}}}\leq 2\mathfrak{T}%
_{T^{\alpha }}^{\mathbf{b}}+2\sqrt{\varepsilon }\mathfrak{N}_{T^{\alpha }}\ .
\label{test}
\end{equation}
Now we apply the fact that $\left( \mathcal{S}_{\infty }\right) $ holds to obtain
\begin{eqnarray*}
\ \ \ \mathfrak{N}_{T^{\alpha }}
\!\!\!&\lesssim &\!\!\! 
\left( C_{\widehat{\mathbf{b}}}\left(
\infty \right) +C_{\widehat{\mathbf{b}}^{\ast }}\left( \infty \right)
\right) \left\{ \mathfrak{T}_{T^{\alpha }}^{\widehat{\mathbf{b}}}+\mathfrak{T%
}_{T^{\alpha ,\ast }}^{\widehat{\mathbf{b}}^{\ast }}+\sqrt{\mathfrak{A}%
_{2}^{\alpha }}+\mathfrak{E}_{2}^{\alpha }\right\}
\end{eqnarray*}%
and take $\varepsilon =\varepsilon ^{\ast }$ to conclude, using (\ref{cinfty}) and (\ref{test}), that
\begin{eqnarray}\label{norm after cimplied}
&&\ \ \ \  \mathfrak{N}_{T^\alpha}\lesssim 
C_{implied}\left( C_{\mathbf{b}}\left( p\right) +C_{\mathbf{b}^{\ast
}}\left( p\right) \right) ^{\frac{p}{p-2}}\varepsilon ^{-\frac{1}{p-2}%
}\left\{ \mathfrak{T}_{T^{\alpha }}^{\mathbf{b}}+\mathfrak{T}_{T^{\alpha
,\ast }}^{\mathbf{b}^{\ast }}+\sqrt{\mathfrak{A}_{2}^{\alpha }}+\mathfrak{E}
_{2}^{\alpha }\right\} \\
&&\notag
\hspace{4cm}+ C_{implied}
\left( C_{\mathbf{b}}\left(
p\right) +C_{\mathbf{b}^{\ast }}\left( p\right) \right) ^{\frac{p}{p-2}%
}\varepsilon ^{\frac{1}{2}-\frac{1}{p-2}}\mathfrak{N}_{T^{\alpha }}
\end{eqnarray}
Now we choose 
\begin{equation*}
\varepsilon =\frac{1}{\Gamma }\left( C_{\mathbf{b}}\left( p\right) +C_{%
\mathbf{b}^{\ast }}\left( p\right) \right) ^{-\frac{\frac{p}{p-2}}{\frac{1}{2%
}-\frac{1}{p-2}}}
\end{equation*}
with $\Gamma =\left( 2C_{{implied}}\right) ^{4} $, which satisfies $\Gamma\geq 1$, so that the final term on the right satisfies
\begin{eqnarray*}
C_{{implied}}\left( C_{\mathbf{b}}\left( p\right) +C_{\mathbf{b}
^{\ast }}\left( p\right) \right) ^{\frac{p}{p-2}}\ \varepsilon ^{\frac{1}{2}-%
\frac{1}{p-2}}\ \mathfrak{N}_{T^{\alpha }}
\leq 
C_{{implied}}\ \left( \frac{1}{\Gamma }\right) ^{\frac{1}{4}%
}\ \mathfrak{N}_{T^{\alpha }}=\frac{1}{2}\mathfrak{N}_{T^{\alpha }}
\end{eqnarray*}%
where we have used $\frac{1}{2}-\frac{1}{p-2}\geq \frac{1}{4}$ for $p>6$.
This term can then be absorbed into the left hand side of (\ref{norm after cimplied}) to obtain
\begin{eqnarray*}
\mathfrak{N}_{T^{\alpha }} 
&\lesssim &
\left( C_{\mathbf{b}}\left( p\right) +C_{\mathbf{b}^{\ast
}}\left( p\right) \right) ^{\frac{p}{p-2}\left\{ 1+\frac{\frac{1}{p-2}}{%
\frac{1}{2}-\frac{1}{p-2}}\right\} }\left\{ \mathfrak{T}_{T^{\alpha }}^{%
\mathbf{b}}+\mathfrak{T}_{T^{\alpha ,\ast }}^{\mathbf{b}^{\ast }}+\sqrt{%
\mathfrak{A}_{2}^{\alpha }}+\mathfrak{E}_{2}^{\alpha }\right\}
\end{eqnarray*}%
Since 
\begin{equation*}
\frac{p}{p-2}\left\{ 1+\frac{\frac{1}{p-2}}{\frac{1}{2}-\frac{1}{p-2}}%
\right\} =\left( 1+\frac{2}{p-2}\right) \left( 1+\frac{2}{p-4}\right) \leq 3%
\text{ for }p>6,
\end{equation*}%
we get%
\begin{equation*}
\mathfrak{N}_{T^{\alpha }}\lesssim \left( C_{\mathbf{b}}\left( p\right) +C_{%
\mathbf{b}^{\ast }}\left( p\right) \right) ^{3}\left\{ \mathfrak{T}%
_{T^{\alpha }}^{\mathbf{b}}+\mathfrak{T}_{T^{\alpha ,\ast }}^{\mathbf{b}%
^{\ast }}+\sqrt{\mathfrak{A}_{2}^{\alpha }}+\mathfrak{E}_{2}^{\alpha
}\right\} ,
\end{equation*}%
which completes the proof of $\left( \mathcal{S}_{0}\right) $.

We now show that $\left( \mathcal{S}_{p}\right) $ holds for all $p\!\in\! \left(
p_{m+1},p_{m}\right] $. So fix $m\geq 1$, $p\!\in\! \left( p_{m+1},p_{m}\right] $,
and suppose that $\mathbf{b}=\left\{ b_{Q}\right\} _{Q\in \mathcal{P}}$ is a 
$p$-weakly $\sigma $-accretive family of functions on $\mathbb{R}^n$ and that 
$\mathbf{b}^{\ast }=\left\{ b_{Q}^{\ast }\right\} _{Q\in \mathcal{P}}$ is a $p$-weakly $\omega $-accretive family of functions on $\mathbb{R}^n$. Note that
the sequence $\left\{ p_{m}\right\} _{m=0}^{\infty }=\left\{ \frac{2}{1-\left( \frac{2}{3}\right) ^{m}}\right\} _{m=0}^{\infty }$ satisfies the
recursion relation
\begin{equation*}
p_{m+1}=\frac{6}{1+\frac{4}{p_{m}}},\text{ equivalently, }p_{m}=\frac{4}{%
\frac{6}{p_{m+1}}-1},\text{\ \ }m\geq 0.
\end{equation*}%
Choose $q\in \left( p_{m},p_{m-1}\right] $ so that%
\begin{equation}
p>\frac{6}{1+\frac{4}{q}}=\frac{6q}{q+4},\text{ i.e. }q<\frac{4}{\frac{6}{p}-1}=\frac{4p}{6-p},
\label{restrict}
\end{equation}%
which can be done since $p>p_{m+1}=\frac{2}{1-\left( \frac{2}{3}\right)
^{m+1}}$ is equivalent to $p_{m}=\frac{2}{1-\left( \frac{2}{3}\right) ^{m}}\!<\!\frac{4}{\frac{6}{p}-1}$, which leaves room to choose $q$ satisfying $%
p_{m}<q<\frac{4}{\frac{6}{p}-1}$.

Now let $0<\varepsilon <1$ (to be fixed later), define $\lambda =\lambda
\left( \varepsilon \right) $ as in (\ref{lambda choice}), and define $%
\widehat{b}_{Q}$ as in (\ref{new}). Recall from (\ref{hence}) and (\ref{low}) that we then have
\begin{equation*}
\int_{\left\{ \left\vert b_{Q}\right\vert >\lambda \right\} }\left\vert
b_{Q}\right\vert ^{2}d\sigma \leq \varepsilon \left\vert Q\right\vert
_{\sigma }
\text{ and }
\left\vert \frac{1}{\left\vert Q\right\vert _{\sigma }}\int_{Q}\widehat{b}%
_{Q}d\sigma \right\vert \geq 1,\ \ \ \ \ Q\in \mathcal{P}\ ,
\end{equation*}%
if we choose $0<\varepsilon \leq \frac{1}{4}$. We of course have the
previous upper bound 
\begin{equation*}
\left\Vert \widehat{b}_{Q}\right\Vert _{L^{\infty }\left( \sigma \right)
}\leq 2\lambda =2\lambda \left( \varepsilon \right) =2\left( \frac{p}{p-2}C_{%
\mathbf{b}}\left( p\right) ^{p}\frac{1}{\varepsilon }\right) ^{\frac{1}{p-2}}
\end{equation*}%
and while this turned out to be sufficient in the case $m=0$, we must do
better than $O\left( \frac{1}{\varepsilon }\right) ^{\frac{1}{p-2}}$ in the
case $m\geq 1$. In fact we compute the $L^{q}$ norm instead, recalling that $q>p$ and using Chebysev's inequality,

\begin{eqnarray*}
\left( \frac{1}{\left\vert Q\right\vert _{\mu }}\int_{Q}\left\vert 
\widehat{b}_{Q}\right\vert ^{q}d\mu \right) ^{\!\!\frac{1}{q}}
\!\!\!\!\!\!&=&\!\!\!\!
2\left( \frac{1}{%
\left\vert Q\right\vert _{\mu }}\int_{Q}\left\vert b_{Q}\left( \mathbf{1}%
_{\left\{ \left\vert b_{Q}\right\vert \leq \lambda \right\} }+\frac{\lambda 
}{\left\vert b_{Q}\right\vert }\mathbf{1}_{\left\{ \left\vert
b_{Q}\right\vert >\lambda \right\} }\right) \right\vert ^{q}d\mu \right) ^{%
\frac{1}{q}} \\
&=&\!\!\!\!
2\left( \frac{1}{\left\vert Q\right\vert _{\mu }}\int_{\left\{ \left\vert
b_{Q}\right\vert \leq \lambda \right\} }\left[ \int_{0}^{\left\vert
b_{Q}\right\vert }qt^{q-1}dt\right] d\sigma +\frac{\lambda ^{q}\left\vert
\left\{ \left\vert b_{Q}\right\vert >\lambda \right\} \right\vert _{\mu }}{%
\left\vert Q\right\vert _{\mu }}\right) ^{\!\!\frac{1}{q}} \\
&\leq &\!\!\!\!
2\left( \frac{1}{\left\vert Q\right\vert _{\mu }}\int_{0}^{\lambda }%
\left[ \int_{\left\{ t<\left\vert b_{Q}\right\vert \leq \lambda \right\}
}d\sigma \right] qt^{q-1}dt+C_{\mathbf{b}}\left( p\right) ^{p}\lambda
^{q-p}\right) ^{\frac{1}{q}} \\
&\leq &\!\!\!\!
2\left( \frac{1}{\left\vert Q\right\vert _{\mu }}\int_{0}^{\lambda }%
\left[ \frac{1}{t^{p}}\int \left\vert b_{Q}\right\vert ^{p}d\sigma \right]
qt^{q-1}dt+C_{\mathbf{b}}\left( p\right) ^{p}\lambda ^{q-p}\right) ^{\frac{1%
}{q}} \\
&\leq &\!\!\!\!
2C_{\mathbf{b}}\left( p\right) ^{\frac{p}{q}}\left( \int_{0}^{\lambda
}qt^{q-p-1}dt+\lambda ^{q-p}\right) ^{\frac{1}{q}}\\
&=&\!\!\!\!
2C_{\mathbf{b}}\left( p\right) ^{\frac{p}{q}}\left( \frac{2q-p}{q-p}\lambda ^{q-p}\right) ^{\frac{1}{q}}
\end{eqnarray*}%
which shows that $C_{\widehat{\mathbf{b}}}\left( q\right) $ satisfies the
estimate%
\begin{eqnarray*}
C_{\widehat{\mathbf{b}}}\left( q\right) 
&\leq &
2C_{\mathbf{b}}\left(
p\right) ^{\frac{p}{q}}\left( \frac{2q-p}{q-p}\right) ^{\frac{1}{q}}\left[
\left( \frac{p}{p-2}C_{\mathbf{b}}\left( p\right) ^{p}\frac{1}{\varepsilon }%
\right) ^{\frac{1}{p-2}}\right] ^{1-\frac{p}{q}} \\
&\lesssim &
C_{\mathbf{b}}\left( p\right) ^{\frac{p}{q}\left( \frac{q-2}{p-2}\right) }\varepsilon ^{-\frac{1-\frac{p}{q}}{p-2}}\lesssim C_{\mathbf{b}%
}\left( p\right) ^{\frac{3}{2}}\varepsilon ^{-\frac{1-\frac{p}{q}}{p-2}},
\end{eqnarray*}%
a significant improvement over the bound $O\left( \varepsilon ^{-\frac{1}{p-2%
}}\right) $. Here we have used that if $p>\frac{6q}{%
q+4}$, then 
$$
\frac{p}{q}\left(\frac{q-2}{p-2}\right)<\frac{\frac{6q}{q-4}}{\frac{6q}{q-4}-2}\frac{q-2}{q} <\frac{3}{2}
$$
as the function $x\mapsto\frac{x}{x-2}$ is decreasing when $x>2$. Moreover, from (\ref{test}) we also have%
\begin{equation*}
\mathfrak{T}_{T^{\alpha }}^{\widehat{\mathbf{b}}}\leq 2\mathfrak{T}%
_{T^{\alpha }}^{\mathbf{b}}+2\sqrt{\varepsilon }\mathfrak{N}_{T^{\alpha }}\ .
\end{equation*}%

We can do the same for the dual testing functions $\mathbf{b}^{\ast
}=\left\{ b_{Q}^{\ast }\right\} _{Q\in \mathcal{P}}$ and then altogether, provided $0<\varepsilon\leq\frac{1}{4}$, we have both
\begin{eqnarray*}
1 &\leq &\left\vert \frac{1}{\left\vert Q\right\vert _{\sigma }}\int_{Q}%
\widehat{b}_{Q}d\sigma \right\vert \leq \left\Vert \widehat{b}%
_{Q}\right\Vert _{L^{q}\left( \sigma \right) }\leq C_{\mathbf{b}}\left(
p\right) ^{\frac{3}{2}}\varepsilon ^{-\frac{1-\frac{p}{q}}{p-2}},\ \ \ \ \
Q\in \mathcal{P}\ , \\
&&\ \ \ \ \ \ \ \ \ \ \mathfrak{T}_{T^{\alpha }}^{\widehat{\mathbf{b}}}\leq 2%
\mathfrak{T}_{T^{\alpha }}^{\mathbf{b}}+2\sqrt{\varepsilon }\mathfrak{N}%
_{T^{\alpha }}\ ,
\end{eqnarray*}
as well as
\begin{eqnarray*}
1 &\leq &\left\vert \frac{1}{\left\vert Q\right\vert _{\omega }}\int_{Q}
\widehat{b^{\ast }}_{Q}d\omega \right\vert \leq \left\Vert \widehat{b^{\ast }
}_{Q}\right\Vert _{L^{q}\left( \omega \right) }\leq C_{\mathbf{b}^{\ast
}}\left( p\right) ^{\frac{3}{2}}\varepsilon ^{-\frac{1-\frac{p}{q}}{p-2}},\
\ \ \ \ Q\in \mathcal{P}\ , \\
&&\ \ \ \ \ \ \ \ \ \ \mathfrak{T}_{T^{\alpha }}^{\widehat{\mathbf{b}^{\ast }
}}\leq 2\mathfrak{T}_{T^{\alpha }}^{\mathbf{b}^{\ast }}+2\sqrt{\varepsilon}\mathfrak{N}_{T^{\alpha }}
\end{eqnarray*}
We now use these estimates, together with the fact that $\left( \mathcal{S}%
_{m-1}\right) $ holds, to obtain%
\begin{eqnarray*}
\mathfrak{N}_{T^{\alpha }} 
\!\!\!\!\!\!&\lesssim & \!\!\!\!\!
\left( C_{\widehat{\mathbf{b}}}\left(
q\right) \!+\! C_{\widehat{\mathbf{b}}^{\ast }}\left( q\right) \right)
^{3^{n}}\left\{ \mathfrak{T}_{T^{\alpha }}^{\widehat{\mathbf{b}}}+\mathfrak{T%
}_{T^{\alpha ,\ast }}^{\widehat{\mathbf{b}}^{\ast }}+\sqrt{\mathfrak{A}%
_{2}^{\alpha }}+\mathfrak{E}_{2}^{\alpha }\right\} \\
&\lesssim &\!\!\!\!\!
\left( C_{\mathbf{b}}\left( p\right) \!+\! C_{\mathbf{b}^{\ast
}}\left( p\right) \right) ^{\frac{3}{2}3^{n}}\!\!\varepsilon ^{-\frac{1-\frac{p}{q}}{p-2}} 
\!\left\{\! \left[ \mathfrak{T}
_{T^{\alpha }}^{\mathbf{b}}\!+\!\sqrt{\varepsilon }\mathfrak{N}_{T^{\alpha }}%
\!\right]\! +\! \left[ \mathfrak{T}_{T^{\alpha ,\ast }}^{\mathbf{b}^{\ast }}\!+\!\sqrt{\varepsilon }\mathfrak{N}_{T^{\alpha }}\!\right]\! +\! \sqrt{\mathfrak{A}%
_{2}^{\alpha }}\!+\!\mathfrak{E}_{2}^{\alpha }\!\right\} \\
&\lesssim &\!\!\!\!\!
\left( C_{\mathbf{b}}\left( p\right) +C_{\mathbf{b}^{\ast
}}\left( p\right) \right) ^{\frac{3}{2}3^{n}}
\!\!\left(\!\varepsilon ^{-\frac{1-\frac{p}{q}}{p-2}}\!\left\{ \mathfrak{T}_{T^{\alpha }}^{\mathbf{b}}\!+\!\mathfrak{T}%
_{T^{\alpha ,\ast }}^{\mathbf{b}^{\ast }}\!+\!\sqrt{\mathfrak{A}_{2}^{\alpha }}\!+\!
\mathfrak{E}_{2}^{\alpha }\right\}\! +\!
\sqrt{\varepsilon }\varepsilon ^{-\frac{1-\frac{p}{q}}{p-2}}\mathfrak{N}
_{T^{\alpha }}   \right)
\end{eqnarray*}
We can absorb the term 
$
\left( C_{\mathbf{b}}\left( p\right) +C_{\mathbf{b}^{\ast
}}\left( p\right) \right) ^{\frac{3}{2}3^{n}}\sqrt{\varepsilon }\varepsilon ^{-\frac{1-\frac{p}{q}}{p-2}}\mathfrak{N}
_{T^{\alpha }}
$
into the left hand
side as before, by choosing
\begin{equation*}
\varepsilon =\frac{1}{\Gamma }\left( C_{\mathbf{b}}\left( p\right) +C_{\mathbf{b}^{\ast }}\left( p\right) \right) ^{
\left( \frac{\frac{3}{2}3^{n}}{\frac{1-\frac{p}{q}}{p-2}-\frac{1}{2}}\right) }
\end{equation*}
with $\Gamma $ sufficiently large, depending only on the implied constant, since (\ref{restrict}) gives $%
\frac{\frac{6}{p}-1}{2}<\frac{2}{q}$, and hence

\begin{equation}
\frac{1}{2}-\frac{1-\frac{p}{q}}{p-2}=\frac{p\left( 1+\frac{2}{q}\right) -4}{%
2p-4}>\frac{p\left( 1+\frac{\frac{6}{p}-1}{2}\right) -4}{2p-4}=\frac{1}{4}.
\label{quarter}
\end{equation}%
Thus,
\begin{eqnarray*}
\mathfrak{N}_{T^{\alpha }} 
&\lesssim &
\left( C_{\mathbf{b}}\left( p\right) +C_{\mathbf{b}^{\ast
}}\left( p\right) \right) ^{ \frac{3}{2}3^{n} \left( 1+1\right)
}\left\{ \mathfrak{T}_{T^{\alpha }}^{\mathbf{b}}+\mathfrak{T}_{T^{\alpha
,\ast }}^{\mathbf{b}^{\ast }}+\sqrt{\mathfrak{A}_{2}^{\alpha }}+\mathfrak{E}%
_{2}^{\alpha }\right\} .
\end{eqnarray*}%
Here we have used that (\ref{quarter}) implies
$\displaystyle
\frac{\frac{1-\frac{p}{q}}{p-2}}{\frac{1}{2}-\frac{1-\frac{p}{q}}{p-2}}<4%
\frac{1-\frac{p}{q}}{p-2}\leq 1.
$
So we finally have
\begin{equation*}
\mathfrak{N}_{T^{\alpha }}\lesssim \left( C_{\mathbf{b}}\left( p\right) +C_{%
\mathbf{b}^{\ast }}\left( p\right) \right) ^{3^{n+1}}\left\{ \mathfrak{T}%
_{T^{\alpha }}^{\mathbf{b}}+\mathfrak{T}_{T^{\alpha ,\ast }}^{\mathbf{b}%
^{\ast }}+\sqrt{\mathfrak{A}_{2}^{\alpha }}+\mathfrak{E}_{2}^{\alpha
}\right\} ,
\end{equation*}%
which completes the proof of Proposition \ref{conditional}.
\end{proof}

Thus we may assume for the proof of Theorem \ref{dim high}
given below that $p=\infty $ and that the testing functions are real-valued
and satisfy%
\begin{eqnarray}
&& {\supp}b_{Q}\subset Q\ ,\ \ \ \ \ Q\in \mathcal{P},
\label{acc infinity'} \\
1 &\leq &\frac{1}{\left\vert Q\right\vert _{\mu }}\int_{Q}b_{Q}d\mu \leq
\left\Vert b_{Q}\right\Vert _{L^{\infty }\left( \mu \right) }\leq C_{\mathbf{%
b}}\left( \infty \right) <\infty ,\ \ \ \ \ Q\in \mathcal{P}\ .  \notag
\end{eqnarray}

\subsection{Reverse H\"{o}lder control of children}

Here we begin to further reduce the proof of Theorem \ref{dim high} to the
case of bounded real testing functions $\mathbf{b}=\left\{ b_{Q}\right\}
_{Q\in \mathcal{P}}$ having reverse H\"{o}lder control 
\begin{equation}
\left\vert \frac{1}{\left\vert Q^{\prime }\right\vert _{\sigma }}%
\int_{Q^{\prime }}b_{Q}d\sigma \right\vert \geq c\left\Vert \mathbf{1}%
_{Q^{\prime }}b_{Q}\right\Vert _{L^{\infty }\left( \sigma \right) }>0,
\label{rev Hol con}
\end{equation}%
for all children $Q^{\prime }\in \mathfrak{C}\left( Q\right) $ with $%
\left\vert Q^{\prime }\right\vert _{\sigma }>0$ and $Q\in \mathcal{P}$.

\subsubsection{Control of averages over children}

\begin{lem}
\label{further red}Suppose that $\sigma $ and $\omega $ are locally finite
positive Borel measures on $\mathbb{R}^n$. Assume that $%
T^{\alpha }$ is a standard $\alpha $-fractional elliptic and gradient
elliptic singular integral operator on $\mathbb{R}^n$, and set $T_{\sigma
}^{\alpha }f=T^{\alpha }\left( f\sigma \right) $ for any smooth truncation
of $T_{\sigma }^{\alpha }$, so that $T_{\sigma }^{\alpha }$ is \emph{apriori}
bounded from $L^{2}\left( \sigma \right) $ to $L^{2}\left( \omega \right) $.
Let $Q\in \mathcal{P}$ and let $\mathfrak{N}_{T^{\alpha }}\left( Q\right) $
be the best constant in the local inequality%
\begin{equation*}
\sqrt{\int_{Q^{\prime }}\left\vert T_{\sigma }^{\alpha }\left( \mathbf{1}%
_{Q}f\right) \right\vert ^{2}d\omega }\leq \mathfrak{N}_{T^{\alpha }}\left(
Q\right) \sqrt{\int_{Q}\left\vert f\right\vert ^{2}d\sigma }\ ,\ \ \ \ \
f\in L^{2}\left( \mathbf{1}_{Q}\sigma \right) .
\end{equation*}%
Suppose that $b_{Q}$ is a real-valued function supported in $Q$ such that 
\begin{eqnarray*}
&&
1\leq \frac{1}{\left\vert Q\right\vert _{\sigma }}\int_{Q}b_{Q}d\sigma
\leq \left\Vert \mathbf{1}_{Q}b_{Q}\right\Vert _{L^{\infty }\left( \sigma
\right) }\leq C_{\mathbf{b}}\ , \\
&&
\sqrt{\int_{Q}\left\vert T_{\sigma }^{\alpha }b_{Q}\right\vert ^{2}d\omega 
}\leq \mathfrak{T}_{T^{\alpha }}^{b_{Q}}\left( Q\right) \sqrt{\left\vert
Q\right\vert _{\sigma }}\ .
\end{eqnarray*}%
Then for every $0<\delta <\frac{1}{2^{n+1}C_{\mathbf{b}}^{3}}$, there exists a
real-valued function $\widetilde{b}_{Q}$ supported in $Q$ such that 
\begin{enumerate}[(1).]
\item
$\displaystyle 1\leq \frac{1}{\left\vert Q\right\vert _{\sigma }}\int_{Q}\widetilde{b}
_{Q}d\sigma \leq \left\Vert \mathbf{1}_{Q}\widetilde{b}_{Q}\right\Vert
_{L^{\infty }\left( \sigma \right) }\leq 2\left( 1+\sqrt{C_{\mathbf{b}}}\right)
C_{\mathbf{b}}$\ , 
\item
$\displaystyle\sqrt{\int_{Q}\left\vert T_{\sigma }^{\alpha }\widetilde{b}_{Q}\right\vert
^{2}d\omega }\leq \left[ \mathfrak{T}_{T^{\alpha }}^{b_{Q}}\left( Q\right)+2
C_{\mathbf{b}}^{\frac{3}{4}}\delta ^{\frac{1}{4}}\mathfrak{N}_{T^{\alpha
}}\left( Q\right) \right] \sqrt{\left\vert Q\right\vert _{\sigma }}$\ , 
\item
$\displaystyle 0<\left\Vert \mathbf{1}_{Q_{i}}\widetilde{b}_{Q}\right\Vert _{L^{\infty
}\left( \sigma \right) }\leq \frac{16C_{\mathbf{b}}}{\delta }\left\vert \frac{1}{%
\left\vert Q_{i}\right\vert _{\sigma }}\int_{Q_{i}}\widetilde{b}_{Q}d\sigma
\right\vert \ ,\ \ \ \ \ Q_{i}\in \mathfrak{C}\left( Q\right)$ .
\end{enumerate}

\end{lem}

\begin{proof}
Let $0<\delta <1$ and fix $Q\in \mathcal{P}$. By assumption we have%
\begin{equation*}
1\leq \frac{1}{\left\vert Q\right\vert _{\sigma }}\int_{Q}b_{Q}d\sigma \leq
\left\Vert \mathbf{1}_{Q}b_{Q}\right\Vert _{L^{\infty }\left( \sigma \right)
}\leq C_{\mathbf{b}}.
\end{equation*}%
Let $Q_i$ be the children of $Q$.
We now define $\widetilde{b}_{Q}$. First we note that the inequality%
\begin{equation}
\left\vert \frac{1}{\left\vert Q_i\right\vert _{\sigma }}%
\int_{Q_i}b_{Q}d\sigma \right\vert <\frac{\delta }{C_{\mathbf{b}}}
\left\Vert \mathbf{1}_{Q_i}b_{Q}\right\Vert _{L^{\infty }\left(
\sigma \right) }  \label{Q' big}
\end{equation}%
cannot hold for \textit{all} $Q_i$, since otherwise we obtain the contradiction 
\begin{eqnarray*}
\left\vert \int_{Q}b_{Q}d\sigma \right\vert 
&\leq &
\sum_{i=1}^{2^n}\left\vert \int_{Q_i}b_{Q}d\sigma \right\vert 
<
\frac{\delta }{C_{\mathbf{b}}}\sum_{i=1}^{2^n} \left\vert Q_{i}\right\vert _{\sigma }\left\Vert \mathbf{1}_{Q_{i}}b_{Q}\right\Vert _{L^{\infty }\left( \sigma \right) }  \\
&\leq &
\frac{\delta }{C_{\mathbf{b}}}\left\vert Q\right\vert _{\sigma }\left\Vert 
\mathbf{1}_{Q}b_{Q}\right\Vert _{L^{\infty }\left( \sigma \right) }\leq
\delta \left\vert \int_{Q}b_{Q}d\sigma \right\vert <\left\vert
\int_{Q}b_{Q}d\sigma \right\vert .
\end{eqnarray*}%
If (\ref{Q' big}) holds for none of the $Q_{i}$, then we simply define $\widetilde{b}_{Q}=b_{Q}$, and trivially all the conclusions of the Lemma \ref{further red} hold. If (\ref{Q' big})
holds for at least one of the children, say $Q_{i_0}$, then we define 
$\widetilde{b}_{Q}$ differently according to how large the $L^{1}\left(
\sigma \right) $-average $\frac{1}{\left\vert Q_{i_0}\right\vert
_{\sigma }}\int_{Q_{i_0}}\left\vert b_{Q}\right\vert d\sigma $ is. In this case, define $\tilde{G}$ to be the set of indices for which (\ref{Q' big}) holds and $G$ the set of indices for which (\ref{Q' big}) fails. We define
\begin{eqnarray*}
\widetilde{b}_{Q} 
&\equiv &
\sum_{i\in G}b_Q\mathbf{1}_Q 
+
\sum_{i\in G_{0}}\delta \mathbf{1}%
_{Q_{i}}
+
\sum_{i\in G_{+}}\left( \frac{1}{\left\vert Q_{i}\right\vert
_{\sigma }}\int_{Q_{i}}\left\vert b_{Q}\right\vert d\sigma \right) \mathbf{1}%
_{Q_{i}} \\
&&
+\sum_{i\in B_{-}}\left(  p_{i}-n_{i}\left( 1+\sqrt{C_{\mathbf{b}}\delta }\right) 
 \right) \mathbf{1}_{Q_{i}} 
+\sum_{i\in B_{+}}\left(  \left( 1+\sqrt{C_{\mathbf{b}}\delta }\right) p_{i}-n_{i}
\right) \mathbf{1}_{Q_{i}} 
\end{eqnarray*}%
where
\begin{eqnarray*}
G_{0} 
&\equiv &
\left\{ i\in\tilde{G}:\frac{1}{\left\vert Q_{i}\right\vert _{\sigma }}%
\int_{Q_{i}}\left\vert b_{Q}\right\vert d\sigma =0\right\} \\
G_{+} 
&\equiv &
\left\{ i\in\tilde{G}:0<\frac{1}{\left\vert Q_{i}\right\vert _{\sigma }}%
\int_{Q_{i}}\left\vert b_{Q}\right\vert d\sigma \leq \sqrt{C_{\mathbf{b}%
}\delta }\right\} , \\
B_{-} 
&\equiv &
\left\{ i\in\tilde{G}:\frac{1}{\left\vert Q_{i}\right\vert _{\sigma }}%
\int_{Q_{i}}\left\vert b_{Q}\right\vert d\sigma >\sqrt{C_{\mathbf{b}}\delta }%
\text{ and }\int_{Q_{i}}n_id\sigma >\int_{Q_{i}}p_i d\sigma \right\} , \\
B_{+} 
&\equiv &
\left\{ i\in\tilde{G}:\frac{1}{\left\vert Q_{i}\right\vert _{\sigma }}%
\int_{Q_{i}}\left\vert b_{Q}\right\vert d\sigma >\sqrt{C_{\mathbf{b}}\delta }%
\text{ and }\int_{Q_{i}}p_i d\sigma \geq \int_{Q_{i}}n_i d\sigma \right\} .
\end{eqnarray*}
and $p_i, n_i$ are the positive and negative parts of $b_Q$ respectively on $Q_i$, i.e.
\begin{eqnarray*}
\mathbf{1}_{Q_{i}}\left( x\right) b_{Q}\left( x\right)
&=&
p_i\left( x\right) -n_i\left( x\right) , \\
\mathbf{1}_{Q_{i}}\left( x\right) \left\vert b_{Q}\left(
x\right) \right\vert 
&=&
p_i\left( x\right) +n_i\left( x\right) ,
\end{eqnarray*}%
Now let us check the conclusions of the Lemma \ref{further red}. For \textit{(1)} we have 
\begin{eqnarray*}
 1&\leq& \frac{1}{\left\vert Q\right\vert _{\sigma }}\int_{Q}b_{Q}d\sigma \\
 &\leq&
 \frac{1}{\left\vert Q\right\vert _{\sigma }}\int_{Q}\widetilde{b}
_{Q}d\sigma
+
\frac{1}{|Q|_\sigma}\sum_{i \in B_-}\int_{Q_i}n_i\sqrt{C_b\delta}d\sigma
-
\frac{1}{|Q|_\sigma}\sum_{i \in B_+}\int_{Q_i}p_i\sqrt{C_b\delta}d\sigma\\
&\leq &
\frac{1}{\left\vert Q\right\vert _{\sigma }}\int_{Q}\widetilde{b}
_{Q}d\sigma+\sqrt{C_b\delta}C_b \frac{1}{|Q|_\sigma}\sum_{i \in B_-}|Q_i|_\sigma
\leq
\frac{1}{\left\vert Q\right\vert _{\sigma }}\int_{Q}\widetilde{b}
_{Q}d\sigma+C_b^\frac{3}{2}\sqrt{\delta}
\end{eqnarray*}
and choosing $\delta$ small enough we get
$$
\frac{1}{2}\leq \frac{1}{\left\vert Q\right\vert _{\sigma }}\int_{Q}\widetilde{b}
_{Q}d\sigma \leq\left\Vert\mathbf{1}_Q\widetilde{b}
_{Q} \right\Vert_{L^{\infty}(\sigma)},
$$
which in turn is bounded by
$$
\sup_{Q_i\in \mathfrak{C}(Q)} \left\Vert \mathbf{1}_{Q_i}\widetilde{b}_{Q_i}\right\Vert
_{L^{\infty }\left( \sigma \right) }
\leq 2\left( 1+\sqrt{C_{\mathbf{b}}}\right)
C_{\mathbf{b}}
$$
by taking the different cases on $Q_i$:
\begin{enumerate}[(a)]

    \item For $i \in G_0$, $\left\Vert \mathbf{1}_{Q_i}\widetilde{b}_{Q_i}\right\Vert_{L^{\infty }}\leq \delta$,
    \item For $i \in G_+$, $\left\Vert \mathbf{1}_{Q_i}\widetilde{b}_{Q_i}\right\Vert
_{L^{\infty }}\leq C_\mathbf{b}$,
    \item For $i \in B_-\cup B_+$, $\left\Vert \mathbf{1}_{Q_i}\widetilde{b}_{Q_i}\right\Vert
_{L^{\infty }}\leq 2(1+\sqrt{C_\mathbf{b}})C_\mathbf{b}$.
\end{enumerate}
This completes the proof for \textit{(1)}.

For \textit{(2)}, we have from Minkowski's inequality
\begin{eqnarray*}
\sqrt{\frac{1}{\left\vert Q\right\vert _{\sigma }}\int_{Q}\left\vert
T_{\sigma }^{\alpha }\widetilde{b}_{Q}\right\vert ^{2}d\omega }
\!\!\!\!&\leq &\!\!\!\!
\sqrt{\frac{1}{\left\vert Q\right\vert _{\sigma }}\int_{Q}\left\vert T_{\sigma
}^{\alpha }b_{Q}\right\vert ^{2}\!d\omega } 
+
\sqrt{\frac{1}{\left\vert
Q\right\vert _{\sigma }}\int_{Q}\left\vert T_{\sigma }^{\alpha }\left( 
\widetilde{b}_{Q}-b_{Q}\right) \right\vert ^{2}\!\!d\omega } \\
&\leq &
\mathfrak{T}_{T^{\alpha }}^{b_{Q}}\left( Q\right) +\mathfrak{N}%
_{T^{\alpha }}\left( Q\right) \sqrt{\frac{1}{\left\vert Q\right\vert
_{\sigma }}\int_{Q}\left\vert \widetilde{b}_{Q}-b_{Q}\right\vert ^{2}\!d\sigma}\\
&=&
\mathfrak{T}_{T^{\alpha }}^{b_{Q}}\left( Q\right) +\mathfrak{N}%
_{T^{\alpha }}\left( Q\right)\sqrt{\frac{1}{\left\vert Q\right\vert
_{\sigma }}\sum_{Q_i \in \mathfrak{C}(Q)}\int_{Q_i}\left\vert \widetilde{b}_{Q}-b_{Q}\right\vert ^{2}\!d\sigma}
\end{eqnarray*}
and this last term is bounded by:
$$
\left(\sum_{i \in G}\ +\ \sum_{i \in G_0}\ +\ \sum_{i \in G_+}\ +\ \sum_{i \in B_-}\ +\ \sum_{i \in B_+}        \right)
\sqrt{\frac{1}{\left\vert Q\right\vert
_{\sigma }} \int_{Q_i}\left\vert \widetilde{b}_{Q}-b_{Q}\right\vert ^{2}\!d\sigma}
$$
and since we have:
\begin{enumerate}[(a)]
\item for $i\in G$,
$$
\frac{1}{\left\vert Q\right\vert
     _{\sigma }}\int_{Q_i}\left\vert   \widetilde{b}_{Q}-b_{Q}\right\vert ^{2}\!d\sigma=0
$$
\item for $i\in G_0$, 
  \begin{eqnarray*}
     \frac{1}{\left\vert Q\right\vert
     _{\sigma }}\int_{Q_i}\left\vert   \widetilde{b}_{Q}-b_{Q}\right\vert ^{2}\!d\sigma
    &\leq&
    \frac{1}{\left\vert Q\right\vert
    _{\sigma }}\left( \int_{Q_i}\delta^2d\sigma +   \int_{Q_i} |b_Q|^2d\sigma     \right)\\
    &\leq&
   \frac{1}{\left\vert Q\right\vert
    _{\sigma }}\left( \delta^2|Q_i|_\sigma + C_{\mathbf{b}}  \int_{Q_i} |b_Q|d\sigma     \right)
    =
    \delta^2 \frac{|Q_i|_\sigma}{|Q|_\sigma}
    \end{eqnarray*}
by the accretivity of $b_Q$ and the definition of $G_0$.
\item for $i\in G_+$,
\begin{eqnarray*}
\frac{1}{\left\vert Q\right\vert _{\sigma }}\int_{Q_i}\left\vert 
\widetilde{b}_{Q}-b_{Q}\right\vert ^{2}d\omega
&=&
\frac{1}{\left\vert Q\right\vert _{\sigma }}\int_{Q_{i}}\left\vert
\left( \frac{1}{\left\vert Q_{i}\right\vert _{\sigma }}\int_{Q_{i}}\left\vert b_{Q}\right\vert d\sigma \right)-b_{Q}\right\vert ^{2}d\sigma \\
&\leq &
\frac{1}{\left\vert Q\right\vert _{\sigma }} \left( \int_{Q_{i}}\left\vert \frac{1}{\left\vert Q_{i}\right\vert _{\sigma
}}\int_{Q_{i}}\left\vert b_{Q}\right\vert d\sigma \right\vert
^{2}\!\! d\sigma 
+
\int_{Q_{i}}\left\vert b_{Q}\right\vert ^{2}d\sigma \!\!\right)  \\
&\leq &
\frac{1}{\left\vert Q\right\vert _{\sigma }} \left( \int_{Q_{i}}C_{\mathbf{b}}\delta d\sigma 
+
C_{\mathbf{b}}\int_{Q_{i}}\left\vert
b_{Q}\right\vert d\sigma \right)\\
&\leq &
\left(C_{\mathbf{b}}\delta 
+
C_{\mathbf{b}}\sqrt{C_{\mathbf{b}}\delta }\right)
\frac{|Q_i|_\sigma}{|Q|_\sigma}
\leq
2C_{\mathbf{b}}^{\frac{3}{2}}\delta ^{\frac{1}{2}}\frac{|Q_i|_\sigma}{|Q|_\sigma}.
\end{eqnarray*}
\item for $i\in B_-$,
\begin{eqnarray*}
\frac{1}{\left\vert Q\right\vert _{\sigma }}\int_{Q_i}\left\vert 
\widetilde{b}_{Q}-b_{Q}\right\vert ^{2}d\sigma 
\!\!\!&=&\!\!\!
\frac{1}{\left\vert Q\right\vert _{\sigma }}\int_{Q_{i}}\left\vert 
C_{\mathbf{b}}\delta n_i\right\vert ^{2}d\sigma =
C_{\mathbf{b}}\delta
\frac{1}{\left\vert Q\right\vert _{\sigma }} \int_{Q_{i}}\left\vert n_i\right\vert ^{2}d\sigma  \\
&\leq &
C_{\mathbf{b}}^{3}\delta\frac{|Q_i|_\sigma}{|Q|_\sigma}.
\end{eqnarray*}%
\item and for $i\in B_+$, the same estimate as in the previous case,
\end{enumerate}
we obtain
\begin{equation*}
\sqrt{\frac{1}{\left\vert Q\right\vert _{\sigma }}\int_{Q}\left\vert
T_{\sigma }^{\alpha }\widetilde{b}_{Q}\right\vert ^{2}d\omega }\leq 
\mathfrak{T}_{T^{\alpha }}^{b_{Q}}\left( Q\right) +2\cdot 2^n C_{\mathbf{b}}^{\frac{3}{4}%
}\delta ^{\frac{1}{4}}\mathfrak{N}_{T^{\alpha }}\left( Q\right) \ .
\end{equation*}
where the dimensional constant comes from 
$$
\frac{1}{\sqrt{|Q|_\sigma}}\sum_{i=1}^{2^n}\sqrt{|Q_i|_\sigma}\leq 2^n .
$$
Now we are left with verifying \textit{(3)}. Note that 
\begin{enumerate}[(a)]
    \item for $i\in G$, the inequality (\ref{Q' big}) does not hold and as $\widetilde{b}_Q=b_Q$ there, immediately we obtain 
    $$
    \left\Vert \mathbf{1}_{Q_i}\widetilde{b}_{Q}\right\Vert _{L^{\infty }\left(\sigma \right) } \leq \left\vert \frac{C_{\mathbf{b}}}{\delta} \int_{Q_i}\widetilde{b}_{Q}d\sigma  \right\vert
    $$
     \item for $i\in G_0\cup G_+$, 
$$ 
\frac{\left\Vert \mathbf{1}_{Q_{i}}\widetilde{b}_{Q}\right\Vert
_{L^{\infty }\left( \sigma \right) }}{\left\vert \frac{1}{\left\vert Q_{{i}}\right\vert _{\sigma }}\int_{Q_{i}}\widetilde{b}_{Q}d\sigma \right\vert } 
=
1<\frac{C_{\mathbf{b}}}{\delta}
$$
\item for $i\in B_-$, 
\begin{eqnarray*}
\frac{\left\Vert \mathbf{1}_{Q_{i}}\widetilde{b}_{Q}\right\Vert
_{L^{\infty }\left( \sigma \right) }}{\left\vert \frac{1}{\left\vert Q_{i}\right\vert _{\sigma }}\int_{Q_{i}}\widetilde{b}%
_{Q}d\sigma \right\vert } 
&\leq &
\frac{\left( 1+\sqrt{C_{\mathbf{b}}\delta }%
\right) C_{\mathbf{b}}}{\left\vert \frac{1}{\left\vert Q_{i}\right\vert _{\sigma }}\int_{Q_{i}}\left[ p_i-n_i\left( 1+\sqrt{C_{\mathbf{b}}\delta }\right) \right] d\sigma \right\vert } \\
&\leq &
\frac{\left( 1+\sqrt{C_{\mathbf{b}}\delta }\right) C_{\mathbf{b}}}{\left\vert 
\sqrt{C_{\mathbf{b}}\delta }\frac{1}{\left\vert Q_{i}\right\vert
_{\sigma }}\int_{Q_{i}}n_i d\sigma \right\vert }\\
&\leq&
\frac{2( 1+\sqrt{C_{\mathbf{b}}\delta}) C_{\mathbf{b}}}{\sqrt{C_{\mathbf{b}}\delta }\frac{1}{\left\vert Q_{i%
}\right\vert _{\sigma }}\int_{Q_{i}}\left\vert b_{Q}\right\vert
d\sigma } \\
&\leq &
\frac{4C_{\mathbf{b}}}{C_{\mathbf{b}}\delta }=\frac{4}{\delta }\text{\ }, 
\end{eqnarray*}%
as, by taking $0<\delta<\frac{1}{4C_{\mathbf{b}}^{3}}$, we have $1+\sqrt{C_{\mathbf{b}}\delta}<2$.
\item and for $i\in B_+$ similarly as in the previous case.
\end{enumerate}

In order to obtain the inequalities for $\widetilde{b}_{Q}$ in the
conclusion of Lemma \ref{further red}, we simply multiply the above function 
$\widetilde{b}_{Q}$ by a factor of $2$.

Finally, if $\left\vert b_{Q}\right\vert \geq c_{1}>0$, we easily see that $\left\vert
\widetilde{b}_{Q}\right\vert \geq \left\vert b_{Q}\right\vert\geq c_1>0 $ as well. This
completes the proof of Lemma \ref{further red}.
\end{proof}

\subsubsection{Control of averages in coronas}

Let $\mathcal{D}_{Q}$ be the grid of dyadic subcubes of $Q$. In the
construction of the triple corona below, we will need to repeat the
construction in the previous subsubsection for a subdecomposition $\left\{
Q_{i}\right\} _{i=1}^{\infty }$ of dyadic subcubes $Q_{i}\in \mathcal{D}%
_{Q}$ of a cube $Q$. Define the corona corresponding to the
subdecomposition $\left\{ Q_{i}\right\} _{i=1}^{\infty }$ by 
\begin{equation*}
\mathcal{C}_{Q}\equiv \mathcal{D}_{Q}\backslash \bigcup_{i=1}^{\infty }%
\mathcal{D}_{Q_{i}}\ .
\end{equation*}

\begin{lem}
\label{prelim control of corona}Suppose that $\sigma $ and $\omega $ are
locally finite positive Borel measures on $\mathbb{R}^n$. Assume
that $T^{\alpha }$ is a standard $\alpha $-fractional elliptic and gradient
elliptic singular integral operator on $\mathbb{R}^n$, and set $T_{\sigma
}^{\alpha }f=T^{\alpha }\left( f\sigma \right) $ for any smooth truncation
of $T_{\sigma }^{\alpha }$, so that $T_{\sigma }^{\alpha }$ is \emph{apriori}
bounded from $L^{2}\left( \sigma \right) $ to $L^{2}\left( \omega \right) $.
Let $Q\in \mathcal{P}$ and let $\mathfrak{N}_{T^{\alpha }}\left( Q\right) $
be the best constant in the local inequality%
\begin{equation*}
\sqrt{\int_{Q}\left\vert T_{\sigma }^{\alpha }\left( \mathbf{1}_{Q}f\right)
\right\vert ^{2}d\omega }\leq \mathfrak{N}_{T^{\alpha }}\left( Q\right) 
\sqrt{\int_{Q}\left\vert f\right\vert ^{2}d\sigma }\ ,\ \ \ \ \ f\in
L^{2}\left( \mathbf{1}_{Q}\sigma \right) .
\end{equation*}%
Let $\left\{ Q_{i}\right\} _{i=1}^{\infty }\subset \mathcal{D}_{Q}$\ be a
collection of pairwise disjoint dyadic subcubes of $Q$. Suppose that $%
b_{Q}$ is a real-valued function supported in $Q$ such that 
\begin{eqnarray*}
&&1\leq \frac{1}{\left\vert Q^{\prime }\right\vert _{\sigma }}%
\int_{Q^{\prime }}b_{Q}d\sigma \leq \left\Vert \mathbf{1}_{Q^{\prime
}}b_{Q}\right\Vert _{L^{\infty }\left( \sigma \right) }\leq C_{\mathbf{b}}\
,\ \ \ \ \ Q^{\prime }\in \mathcal{C}_{Q}\ , \\
&&\sqrt{\int_{Q}\left\vert T_{\sigma }^{\alpha }b_{Q}\right\vert ^{2}d\omega 
}\leq \mathfrak{T}_{T^{\alpha }}^{b_{Q}}\left( Q\right) \sqrt{\left\vert
Q\right\vert _{\sigma }}\ .
\end{eqnarray*}%
Then for every $0<\delta <\frac{1}{4C_{\mathbf{b}}^{3}}$, there exists a
real-valued function $\widetilde{b}_{Q}$ supported in $Q$ such that%
\begin{eqnarray*}
&&
1\leq \frac{1}{\left\vert Q^{\prime }\right\vert _{\sigma }}%
\int_{Q^{\prime }}\widetilde{b}_{Q}d\sigma \leq \left\Vert \mathbf{1}%
_{Q^{\prime }}\widetilde{b}_{Q}\right\Vert _{L^{\infty }\left( \sigma
\right) }\leq 2\left( 1+\sqrt{C_{\mathbf{b}}}\right) C_{\mathbf{b}}\ ,\ \ \
\ \ Q^{\prime }\in \mathcal{C}_{Q}\ , \\
&&
\sqrt{\int_{Q}\left\vert T_{\sigma }^{\alpha }\widetilde{b}_{Q}\right\vert
^{2}d\omega }\leq \left[ 2\mathfrak{T}_{T^{\alpha }}^{b_{Q}}\left( Q\right)
+4C_{\mathbf{b}}^{\frac{3}{2}}\delta ^{\frac{1}{4}}\mathfrak{N}_{T^{\alpha
}}\left( Q\right) \right] \sqrt{\left\vert Q\right\vert _{\sigma }}\ , \\
&&
0<\left\Vert \mathbf{1}_{Q_{i}}\widetilde{b}_{Q}\right\Vert _{L^{\infty
}\left( \sigma \right) }\leq \frac{16C_{\mathbf{b}}}{\delta }\left\vert 
\frac{1}{\left\vert Q_{i}\right\vert _{\sigma }}\int_{Q_{i}}\widetilde{b}%
_{Q}d\sigma \right\vert \ ,\ \ \ \ \ 1\leq i<\infty .
\end{eqnarray*}%
Moreover, if $\left\vert b_{Q}\right\vert \geq c_{1}>0$, then we may take $%
\left\vert \widetilde{b}_{Q}\right\vert \geq c_{1}$ as well.
\end{lem}

The additional gain in the lemma is in the final line that controls the
degeneracy of $\widetilde{b}_{Q}$ at the `bottom' of the corona $\mathcal{C}%
_{Q}$ by establishing a reverse H\"{o}lder control. Note that if we combine
this control with the accretivity control in the corona $\mathcal{C}_{Q}$,
namely 
\begin{equation*}
\left\Vert \mathbf{1}_{Q^{\prime }}\widetilde{b}_{Q}\right\Vert _{L^{\infty
}\left( \sigma \right) }\leq 2\left( 1+\sqrt{C_{\mathbf{b}}}\right) C_{%
\mathbf{b}}\leq 2\left( 1+\sqrt{C_{\mathbf{b}}}\right) C_{\mathbf{b}}\frac{1%
}{\left\vert Q^{\prime }\right\vert _{\sigma }}\int_{Q^{\prime }}\widetilde{b%
}_{Q}d\sigma ,
\end{equation*}%
we obtain reverse H\"{o}lder control throughout the entire collection $%
\mathcal{C}_{Q}\cup \left\{ Q_{i}\right\} _{i=1}^{\infty }$: 
\begin{equation*}
\left\Vert \mathbf{1}_{I}\widetilde{b}_{Q^{\prime }}\right\Vert _{L^{\infty
}\left( \sigma \right) }\leq C_{\delta ,\mathbf{b}}\left\vert \frac{1}{%
\left\vert I\right\vert _{\sigma }}\int_{I}\widetilde{b}_{Q^{\prime
}}d\sigma \right\vert ,\ \ \ \ \ I\in \mathfrak{C}\left( Q^{\prime }\right)
,Q^{\prime }\in \mathcal{C}_{Q}\text{ }.
\end{equation*}%
This has the crucial consequence that the martingale and dual martingale
differences $\bigtriangleup _{Q^{\prime }}^{\sigma ,\mathbf{b}}$ and $\square _{Q^{\prime }}^{\sigma ,\mathbf{b}}$ associated with these functions
as defined in (\ref{def diff}), satisfy%
\begin{equation}
\left\vert \bigtriangleup _{Q^{\prime }}^{\sigma ,\mathbf{b}}h\right\vert
,\left\vert \square _{Q^{\prime }}^{\sigma ,\mathbf{b}}h\right\vert \leq
C_{\delta ,\mathbf{b}}\sum_{I\in \mathfrak{C}\left( Q^{\prime }\right)
}\left( \frac{1}{\left\vert I\right\vert _{\sigma }}\int_{I}\left\vert
h\right\vert d\sigma +\frac{1}{\left\vert Q^{\prime }\right\vert _{\sigma }}%
\int_{Q^{\prime }}\left\vert h\right\vert d\sigma \right) \mathbf{1}_{I}\ .
\label{cruc conseq}
\end{equation}%
However, the defect in this lemma is
that we lose the weak testing condition for $\widetilde{b}_{Q}$ in the
corona even if we had assumed it at the outset for $b_{Q}$.

\begin{proof}
The proof of Lemma \ref{prelim control of corona} is similar to that of the Lemma \ref{further red}. Indeed, we define
\begin{eqnarray*}
\widetilde{b}_{Q} 
&\equiv &
\sum_{i\in G_{0}}\delta \mathbf{1}%
_{Q_{i}}+\sum_{i\in G_{+}}\left( \frac{1}{\left\vert Q_{i}\right\vert
_{\sigma }}\int_{Q_{i}}\left\vert b_{Q}\right\vert d\sigma \right) \mathbf{1}%
_{Q_{i}} \\
&&
+\sum_{i\in B_{-}}\left( \frac{1}{\left\vert Q_{i}\right\vert _{\sigma }}%
\int_{Q_{i}}\left[ p_{i}-n_{i}\left( 1+\sqrt{C_{\mathbf{b}}\delta }\right) 
\right] d\sigma \right) \mathbf{1}_{Q_{i}} \\
&&
+\sum_{i\in B_{+}}\left( \frac{1}{\left\vert Q_{i}\right\vert _{\sigma }}
\int_{Q_{i}}\left[ \left( 1+\sqrt{C_{\mathbf{b}}\delta }\right) p_{i}-n_{i}\right] d\sigma \right) \mathbf{1}_{Q_{i}} \\
&&
+b_{Q}\mathbf{1}_{Q\backslash \cup _{i=1}^{\infty }Q_{i}}\ ,
\end{eqnarray*}
where
\begin{eqnarray*}
G_{0}
&\equiv &
\left\{ i:\frac{1}{\left\vert Q_{i}\right\vert _{\sigma }}%
\int_{Q_{i}}\left\vert b_{Q}\right\vert d\sigma =0\right\}, \\
G_{+} &\equiv &\left\{ i:0<\frac{1}{\left\vert Q_{i}\right\vert _{\sigma }}%
\int_{Q_{i}}\left\vert b_{Q}\right\vert d\sigma \leq \sqrt{C_{\mathbf{b}%
}\delta }\right\} , \\
B_{-} &\equiv &\left\{ i:\frac{1}{\left\vert Q_{i}\right\vert _{\sigma }}%
\int_{Q_{i}}\left\vert b_{Q}\right\vert d\sigma >\sqrt{C_{\mathbf{b}}\delta }%
\text{ and }\int_{Q_{i}}n_id\sigma >\int_{Q_{i}}p_id\sigma \right\} , \\
B_{+} &\equiv &\left\{ i:\frac{1}{\left\vert Q_{i}\right\vert _{\sigma }}%
\int_{Q_{i}}\left\vert b_{Q}\right\vert d\sigma >\sqrt{C_{\mathbf{b}}\delta }%
\text{ and }\int_{Q_{i}}p_id\sigma \geq \int_{Q_{i}}n_id\sigma \right\} .
\end{eqnarray*}%
and $p_i, n_i$ the positive and negative parts of $b_Q$ on each $Q_i$. The proof of Lemma \ref{further red} can be applied verbatim. We emphasise only that when estimating the testing condition, we need the bound
\begin{equation*}
\int_{Q}\left\vert \widetilde{b}_{Q}-b_{Q}\right\vert ^{2}d\sigma \leq
C\left( C_{\mathbf{b}}\right) \delta ^{\frac{1}{4}}\sum_{i=1}^{\infty
}\left\vert Q_{i}\right\vert _{\sigma }\leq C\left( C_{\mathbf{b}}\right)
\delta ^{\frac{1}{4}}\left\vert Q\right\vert _{\sigma }.
\end{equation*}
\end{proof}
\begin{rem}
The estimate $\int_{Q}\left\vert \widetilde{b}_{Q}-b_{Q}\right\vert
^{2}d\sigma \leq C\left( C_{\mathbf{b}}\right) \delta ^{\frac{1}{4}%
}\sum_{i=1}^{\infty }\left\vert Q_{i}\right\vert _{\sigma }$ in the last
line of the above proof is of course too large in general to be dominated by
a fixed multiple of $\left\vert Q^{\prime }\right\vert _{\sigma }$ for $%
Q^{\prime }\in \mathcal{C}_{Q}$, and this is the reason we have no control
of weak testing for $\widetilde{b}_{Q}$ in the rest of the corona even if we
assume weak testing for $b_{Q}$ in the corona $\mathcal{C}_{Q}$. This defect
is addressed in the next subsection below.
\end{rem}

\subsection{Three corona decompositions}

We will use multiple corona constructions, namely a Calder\'{o}n-Zygmund
decomposition, an accretive/testing decomposition, and an energy
decomposition, in order to reduce matters to the stopping form, which is
treated in Section \ref{Sec stop} by adapting the bottom/up stopping time and recursion
of M. Lacey in \cite{Lac}. We will then iterate these corona decompositions
into a single corona decomposition, which we refer to as the \emph{triple
corona}. More precisely, we iterate the first generation of common stopping times
with an infusion of the reverse H\"{o}lder condition on children, followed
by another iteration of the first generation of weak testing stopping times. Recall that we must show the bilinear inequality%
\begin{equation*}
\left\vert \int \left( T_{\sigma }^{\alpha }f\right) gd\omega \right\vert
\leq \mathfrak{N}_{T^{\alpha }}\left\Vert f\right\Vert _{L^{2}\left( \sigma
\right) }\left\Vert g\right\Vert _{L^{2}\left( \omega \right) },\ \ \ \ \
f\in L^{2}\left( \sigma \right) \text{ and }g\in L^{2}\left( \omega \right) .
\end{equation*}

\subsubsection{The Calder\'{o}n-Zygmund corona decomposition}

In this section, we introduce the Calder\'{o}n-Zygmund stopping times $\mathcal{F}$ for
a function $\phi \in L^{2}\left( \mu \right) $ relative to a cube $S_{0}$
and a positive constant $C_{0}\geq 4$. Let $\mathcal{F}=\left\{ F\right\}
_{F\in \mathcal{F}}$ be the collection of Calder\'{o}n-Zygmund stopping
cubes for $\phi $ defined so that $F\subset S_{0}$, $S_{0}\in \mathcal{F}$,
and for all $F\in \mathcal{F}$ with $F\subsetneqq S_{0}$ we have%
\begin{eqnarray*}
\frac{1}{\left\vert F\right\vert _{\mu }}\int_{F}\left\vert \phi \right\vert
d\mu &>&C_{0}\frac{1}{\left\vert \pi _{\mathcal{F}}F\right\vert _{\mu }}%
\int_{F}\left\vert \phi \right\vert d\mu ; \\
\frac{1}{\left\vert F^{\prime }\right\vert _{\mu }}\int_{F^{\prime
}}\left\vert \phi \right\vert d\mu &\leq &C_{0}\frac{1}{\left\vert \pi _{%
\mathcal{F}}F\right\vert _{\mu }}\int_{F}\left\vert \phi \right\vert d\mu \
\ \ \ \text{ for }F\subsetneqq F^{\prime }\subset \pi _{\mathcal{F}}F.
\end{eqnarray*}%
We denote by $\pi _{\mathcal{F}}F$ be the smallest member of $\mathcal{F}$ that \emph{strictly} contains $F$. For a cube $I\in \mathcal{D}$
let $\pi _{\mathcal{D}}I$ be the $\mathcal{D}$-parent of $I$ in the grid $\mathcal{D}$. For $F,F^{\prime }\in 
\mathcal{F}$, we say that $F^{\prime }$ is an $\mathcal{F}$-child of $F$ if $%
\pi _{\mathcal{F}}\left( F^{\prime }\right) =F$ (it could be that $F=\pi _{%
\mathcal{D}}F^{\prime }$), and we denote by $\mathfrak{C}_{\mathcal{F}%
}\left( F\right) $ the set of $\mathcal{F}$-children of $F$. We call $\pi _{%
\mathcal{F}}\left( F^{\prime }\right) $ the $\mathcal{F}$-parent of $%
F^{\prime }\in \mathcal{F}$. 

To achieve the construction above we use the following definition.
\begin{dfn}\label{CZ stopping times}
Let $C_0\geq 4$. Given a dyadic grid $\mathcal{D}$ and a cube $S_{0}\in\mathcal{D}$, define $\mathcal{S}\left(S_{0}\right) $ to be the \emph{maximal} $\mathcal{D}$-subcubes $I\subset S_{0}$ such that%
\begin{equation*}
\frac{1}{\left\vert I\right\vert _{\mu }}\int_{I}\left\vert \phi \right\vert
d\mu >C_{0}\frac{1}{\left\vert S_{0}\right\vert _{\mu }}\int_{S_{0}}\left%
\vert \phi \right\vert d\mu \ ,
\end{equation*}%
and then define the Calder\'{o}n-Zygmund stopping cubes of $S_{0}$ to be the collection 
\begin{equation*}
\mathcal{F}=\left\{ S_{0}\right\} \cup \bigcup\limits_{m=0}^{\infty }
\mathcal{S}_{m}
\end{equation*}
where $\mathcal{S}_{0}=\mathcal{S}\left( S_{0}\right) $ and $\mathcal{S}_{m+1}=\bigcup\limits_{S\in \mathcal{S}_{m}}\mathcal{S}\left( S\right) $
for $m\geq 0$.
\end{dfn}
Define the corona of $F$ by
\begin{equation*}
\mathcal{C}_{F}\equiv \left\{ F^{\prime }\in \mathcal{D}:F\supset F^{\prime
}\supsetneqq H\text{ for some }H\in \mathfrak{C}_{\mathcal{F}}\left(
F\right) \right\} .
\end{equation*}
The stopping cubes $\mathcal{F}$ above satisfy a Carleson condition:
\begin{equation*}
\sum_{F\in \mathcal{F}:\ F\subset \Omega }\left\vert F\right\vert _{\mu
}\leq C\left\vert \Omega \right\vert _{\mu }\ ,\ \ \ \ \ \text{for all open
sets }\Omega .
\end{equation*}%
Indeed, 
\begin{equation*}
\sum_{F^{\prime }\in \mathfrak{C}_{\mathcal{F}}\left( F\right) }\left\vert
F^{\prime }\right\vert _{\mu }\leq \sum_{F^{\prime }\in \mathfrak{C}_{%
\mathcal{F}}\left( F\right) }\frac{\int_{F^{\prime }}\left\vert \phi
\right\vert d\mu }{C_{0}\frac{1}{\left\vert F\right\vert _{\mu }}%
\int_{F}\left\vert \phi \right\vert d\mu }\leq \frac{1}{C_{0}}|F|,
\end{equation*}%
and standard arguments now complete the proof of the Carleson condition.

We emphasize that accretive functions $b$ play no role in the Calder\'{o}n-Zygmund corona decomposition.

\subsubsection{The accretive/testing corona decomposition}

We use a corona construction modelled after that of Hyt\"{o}nen and
Martikainen \cite{HyMa}, that delivers a \emph{weak corona testing condition}
that coincides with the testing condition itself \textbf{only} at the tops
of the coronas. This corona decomposition is developed to optimize the
choice of a new family of real valued testing functions $\left\{ \widehat{b}%
_{Q}\right\} _{Q\in \mathcal{D}}$ taken from the vector $\mathbf{b}\equiv
\left\{ b_{Q}\right\} _{Q\in \mathcal{D}}$ so that we have

\begin{enumerate}
\item the telescoping property at our disposal in each accretive corona,

\item a weak corona testing condition remains in force for the new testing
functions $\widehat{b}_{Q}$ that coincides with the testing condition at
the tops of the coronas,

\item the tops of the coronas, i.e. the stopping cubes, enjoy a Carleson
condition.
\end{enumerate}
We will henceforth refer to the old family as the \emph{original} family,
and denote it by $\left\{ b_{Q}^{{orig}}\right\} _{Q\in \mathcal{D}}$%
. The original family will reappear later in helping to estimate the nearby
form.

Let $\sigma $ and $\omega $ be locally finite Borel measures on $\mathbb{R}^n$
. We assume that the vector of `testing functions' $\mathbf{b}\equiv \left\{
b_{Q}\right\} _{Q\in \mathcal{D}}$ is a $\infty$-weakly $\sigma $-accretive
family, i.e. for $Q\in\mathcal{D}$
\begin{eqnarray*}
&&\supp b_{Q}\subset Q ,\\
0<c_{\mathbf{b}} \!\!\!\!\!&\leq& \!\!\! \frac{1}{\left\vert Q\right\vert _{\mu }}%
\int_{Q}b_{Q}d\sigma  \leq \left\vert\left\vert b_Q\right\vert\right\vert_{L^\infty(\sigma)}\leq C_{\mathbf{b}}<\infty 
\end{eqnarray*}%
and also that $\mathbf{b}^{\ast }\equiv \left\{ b_{Q}\right\} _{Q\in 
\mathcal{D}}$ is a $\infty$-weakly $\omega $-accretive family, and we assume in
addition the testing conditions%
\begin{eqnarray*}
\int_{Q}\left\vert T_{\sigma }^{\alpha }\left( \mathbf{1}_{Q}b_{Q}\right)
\right\vert ^{2}d\omega &\leq & \left(\mathfrak{T}_{T^{\alpha }}^{\mathbf{b}
}\right)^2 \left\vert Q\right\vert _{\sigma }\ ,\ \ \ \ \ \text{for all cubes }Q, \\
\int_{Q}\left\vert T_{\omega }^{\alpha ,\ast }\left( \mathbf{1}
_{Q}b_{Q}^{\ast }\right) \right\vert ^{2}d\sigma 
&\leq &
\left( \mathfrak{T}_{T^{\alpha,\ast}}^{\mathbf{b}^{\ast } }\right)^2\left\vert Q\right\vert _{\omega }\
,\ \ \ \ \ \text{for all cubes }Q.
\end{eqnarray*}

\begin{dfn}
\label{accretive stopping times gen} Given a cube $S_{0}$, define $\mathcal{S}\left( S_{0}\right) $ to be the \emph{maximal} subcubes $%
I\subset S_{0}$ such that satisfy one of the following
\begin{enumerate}[(a).]
\item 
$ \displaystyle
\left\vert \frac{1}{\left\vert I\right\vert _{\mu }}
\int_{I}b_{S_{0}}d\sigma \right\vert <\gamma 
$, or
\item 
$\displaystyle \int_{I}\left\vert T_{\sigma }^{\alpha }\left(
b_{S_{0}}\right) \right\vert ^{2}d\omega >\Gamma \left( \mathfrak{T}%
_{T^{\alpha }}^{\mathbf{b}}\right) ^{2}\left\vert I\right\vert _{\sigma }
$
\end{enumerate}
where the positive constants $\gamma ,\Gamma $ satisfy $0<\gamma <1<\Gamma
<\infty $. Then define the $\mathbf{b}$-accretive stopping cubes of $S_{0}$
to be the collection 
\begin{equation*}
\mathcal{F}=\left\{ S_{0}\right\} \cup \bigcup\limits_{m=0}^{\infty }%
\mathcal{S}_{m}
\end{equation*}%
where $\mathcal{S}_{0}=\mathcal{S}\left( S_{0}\right) $ and $\mathcal{S}%
_{m+1}=\bigcup\limits_{S\in \mathcal{S}_{m}}\mathcal{S}\left( S\right) $
for $m\geq 0$.
\end{dfn}

For $\varepsilon >0$ chosen small enough depending on $p>2$, the $\mathbf{b}$%
-accretive stopping cubes satisfy a $\sigma $-Carleson condition relative to
the measure $\sigma $, and the new testing functions $\left\{ \widetilde{b}%
_{Q}\right\} _{Q\in \mathcal{D}}$, defined by $\widetilde{b_{S}}=\mathbf{1}%
_{S}b_{S_{0}}$ for $S\in \mathcal{C}_{S_{0}}$, satisfy \emph{weak} testing
inequalities. The following lemma is essentially in \cite{HyMa}, but we include a proof for completeness.

\begin{lem}
\label{Car and Test gen} For $\gamma$ small enough and $\Gamma$ large enough, we have the following:

\begin{enumerate}[(1).]
\item For every open set $\Omega $ we have we have the inequality,
\begin{equation}
\sum_{S\in \mathcal{F}:\ S\subset \Omega }\left\vert S\right\vert _{\sigma
}\leq C\left\vert
\Omega \right\vert _{\sigma }\ .  \label{Car gen}
\end{equation}
\item For every cube $S\in \mathcal{C}_{S_{0}}$ we have the weak corona
testing inequality,
\begin{equation}
\int_{S}\left\vert T_{\sigma }^{\alpha }b_{S_{0}}\right\vert ^{2}d\omega
\leq C\left( \mathfrak{T}_{T^{\alpha }}^{\mathbf{b}}\right) ^{2}\left\vert
S\right\vert _{\sigma }\ .  \label{Test gen}
\end{equation}
\end{enumerate}
\end{lem}

\begin{proof}
Inequality (\ref{Test gen}) is immediate from the definition of $\mathcal{F}$ in the definition \ref{accretive stopping times gen}. We now address the Carleson condition (\ref{Car gen}). A standard argument
reduces matters to the case where $\Omega $ is a cube $Q\in \mathcal{F}$
with $\left\vert Q\right\vert _{\sigma }>0$. It suffices to consider each of
the two stopping criteria separately. We first address the stopping
condition $\left\vert \frac{1}{\left\vert I\right\vert _{\sigma }}%
\int_{I}b_{S_{0}}d\sigma \right\vert <\gamma$. Throughout
this proof we will denote the union of these children $\mathcal{S}\left(
Q\right) $ of $Q$ by $E\left( Q\right) \equiv \bigcup\limits_{S\in \mathcal{%
S}\left( Q\right) }S$. Then we have%
\begin{equation*}
\left\vert \int_{E\left( Q\right) }b_{Q}d\sigma \right\vert \leq \sum_{S\in 
\mathcal{S}\left( Q\right) }\left\vert \int_{S}b_{Q}d\sigma \right\vert
<\gamma \sum_{S\in \mathcal{S}\left( Q\right) }\left\vert
S\right\vert _{\sigma }\leq \gamma \left\vert Q\right\vert
_{\sigma }\ ,
\end{equation*}
which together with our hypotheses on $b_{Q}$ gives
\begin{eqnarray*}
\left\vert Q\right\vert _{\sigma } 
&\leq&
\left\vert\int_{Q}b_{Q}d\sigma \right\vert 
\leq
\left\vert \int_{E\left( Q\right)
}b_{Q}d\sigma \right\vert +\left\vert \int_{Q\backslash E\left( Q\right)
}b_{Q}d\sigma \right\vert \\
&\leq &
\gamma\left\vert Q\right\vert _{\sigma }+\sqrt{\int_{Q\backslash E\left( Q\right) }\left\vert b_{Q}\right\vert ^{2}d\sigma }
\sqrt{\left\vert Q\backslash E\left( Q\right) \right\vert _{\sigma }} \\
&\leq &
\gamma \left\vert Q\right\vert _{\sigma }+C_{\mathbf{b}}\sqrt{\left\vert Q\right\vert _{\sigma }}\sqrt{\left\vert
Q\backslash E\left( Q\right) \right\vert _{\sigma }}.
\end{eqnarray*}%
Rearranging the inequality yields 
\begin{eqnarray*}
\left( 1-\gamma \right) \left\vert Q\right\vert _{\sigma }
&\leq &
C_{\mathbf{b}}\sqrt{\left\vert Q\right\vert _{\sigma }}\sqrt{\left\vert Q\backslash E\left( Q\right) \right\vert _{\sigma }}
\end{eqnarray*}
or
$$
\frac{\left( 1-\gamma \right) ^{2}}{C_{\mathbf{b}}^{2}}\left\vert Q\right\vert _{\sigma } 
\leq 
\left\vert Q\backslash
E\left( Q\right) \right\vert _{\sigma }\ ,
$$
which in turn gives
\begin{eqnarray*}
\sum_{S\in \mathcal{S}\left( Q\right) }\left\vert S\right\vert _{\sigma }
&=&
|E(Q)|= \left\vert Q\right\vert _{\sigma }-\left\vert Q\backslash E\left( Q\right)
\right\vert _{\sigma } \\
&\leq &
\left\vert Q\right\vert _{\sigma }-\frac{\left( 1-\gamma \right)^{2}}{C_{\mathbf{b}}^{2}}\left\vert
Q\right\vert _{\sigma }
=\bigg( 1-\frac{\left( 1-\gamma \right) ^{2}}{C_{\mathbf{b}}^{2}}\bigg) \left\vert
Q\right\vert _{\sigma }\equiv \beta \left\vert Q\right\vert _{\sigma }\ ,
\end{eqnarray*}
where $0<\beta <1$ since $1\leq C_{\mathbf{b}}$. If we now
iterate this inequality, we obtain for each $k\geq 1$,
\begin{eqnarray*}
\sum_{\substack{ S\in \mathcal{F}:\ S\subset Q  \\ \pi _{\mathcal{F}
}^{\left( k\right) }\left( S\right) =Q}}\left\vert S\right\vert _{\sigma }
&=&
\sum _{\substack{ S\in \mathcal{F}:\ S\subset Q  \\ \pi _{\mathcal{F}
}^{\left( k-1\right) }\left( S\right) =Q}}\sum_{S^{\prime }\in \mathcal{S}
\left( S\right) }\left\vert S^{\prime }\right\vert _{\sigma }\leq \sum 
_{\substack{ S\in \mathcal{F}:\ S\subset Q  \\ \pi _{\mathcal{F}}^{\left(
k-1\right) }\left( S\right) =Q}}\beta \left\vert S\right\vert _{\sigma } \\
&&\vdots \\
&\leq &\sum_{\substack{ S\in \mathcal{F}:\ S\subset Q  \\ \pi _{\mathcal{F}%
}^{\left( 1\right) }\left( S\right) =Q}}\beta
^{k-1}\left\vert S\right\vert _{\sigma }\leq \beta ^{k}\left\vert
Q\right\vert _{\sigma }\ .
\end{eqnarray*}%
Finally then%
\begin{equation*}
\sum_{S\in \mathcal{F}:\ S\subset Q}\left\vert S\right\vert _{\sigma }\leq\sum_{k=0}^{\infty }\sum_{\substack{ S\in 
\mathcal{F}:\ S\subset Q  \\ \pi _{\mathcal{F}}^{\left( k\right) }\left(
S\right) =Q}}\left\vert S\right\vert _{\sigma }\leq \sum_{k=0}^{\infty
}\beta ^{k}\left\vert Q\right\vert _{\sigma }=\frac{1}{1-\beta }\left\vert
Q\right\vert _{\sigma }=\frac{C_{\mathbf{b}}^{2}}{\left( 1-\gamma\right) ^{2}}\left\vert Q\right\vert _{\sigma }.
\end{equation*}
Now we turn to the second stopping criterion 
$
\int_{I}\left\vert T_{\sigma }^{\alpha }\left( b_{S_{0}}\right) \right\vert
^{2}d\omega >\Gamma \left( \mathfrak{T}_{T^{\alpha }}^{\mathbf{b}}\right)
^{2}\left\vert I\right\vert _{\sigma }.
$
We have
\begin{eqnarray*}
\sum_{S\in \mathfrak{C}_{\mathcal{F}}\left( S_{0}\right) }\left\vert
S\right\vert _{\sigma } 
&\leq &
\frac{1}{\Gamma \left( \mathfrak{T}_{T^{\alpha }}^{\mathbf{b}}\right) ^{2}}\sum_{S\in \mathfrak{C}_{\mathcal{F}
}\left( S_{0}\right) }\int_{S}\left\vert T_{\sigma }^{\alpha }\left(
b_{S_{0}}\right) \right\vert ^{2}d\omega \\
&\leq &
\frac{1}{\Gamma \left( \mathfrak{T}_{T^{\alpha }}^{\mathbf{b}}\right)
^{2}}\int_{S_0}\left\vert T_{\sigma }^{\alpha }\left(b_{S_{0}}\right) \right\vert ^{2}d\omega \leq \frac{1}{\Gamma }\left\vert
S_{0}\right\vert _{\sigma }.
\end{eqnarray*}
Iterating this inequality gives
\begin{equation*}
\sum_{\substack{ S\in \mathcal{F}  \\ S\subset S_{0}}}\left\vert
S\right\vert _{\sigma }\leq \sum_{k=0}^{\infty }\frac{1}{\Gamma ^{k}}%
\left\vert S_{0}\right\vert _{\sigma }=\frac{\Gamma }{\Gamma -1}\left\vert
S_{0}\right\vert _{\sigma },
\end{equation*}
and then
\begin{equation*}
\sum_{\substack{ S\in \mathcal{F}  \\ S\subset \Omega }}\left\vert
S\right\vert _{\sigma }=\sum_{\substack{ \text{maximal }S_{0}\in \mathcal{F} 
\\ S_{0}\subset \Omega }}\sum_{\substack{ S\in \mathcal{F}  \\ S\subset
S_{0} }}\left\vert S\right\vert _{\sigma }\leq \frac{\Gamma }{\Gamma -1}\sum 
_{\substack{ \text{maximal }S_{0}\in \mathcal{F}  \\ S_{0}\subset \Omega }}%
\left\vert S_{0}\right\vert _{\sigma }=\frac{\Gamma }{\Gamma -1}\left\vert
\Omega \right\vert _{\sigma }\ .
\end{equation*}
This completes the proof of Lemma \ref{Car and Test gen}.
\end{proof}

\subsubsection{The energy corona decompositions}

Given a weight pair $\left( \sigma ,\omega \right) $, we  construct an energy
corona decomposition for $\sigma $ and an energy corona decomposition for $\omega $, that uniformize estimates (c.f. \cite{NTV3}, \cite{LaSaUr2}, \cite%
{SaShUr6} and \cite{SaShUr7}). In order to define these constructions, we
recall that the energy condition constant $\mathcal{E}_{2}^{\alpha}$ is given by 
\begin{equation*}
\left( \mathcal{E}_{2}^{\alpha }\right) ^{2}\equiv \sup_{\substack{ Q\in 
\mathcal{P}  \\ Q= \dot{\cup}J_{r}}}\frac{1}{\left\vert Q\right\vert
_{\sigma }}\sum_{r=1}^{\infty }\left( \frac{\mathrm{P}^{\alpha }\left( J_{r},
\mathbf{1}_{Q}\sigma \right) }{\left\vert J_{r}\right\vert^\frac{1}{n} }\right)
^{2}\left\Vert x-m_{J_{r}}\right\Vert _{L^{2}\left( \mathbf{1}_{J_{r}}\omega
\right) }^{2}\ ,
\end{equation*}%
where $\dot{\cup}J_{r}$ is an arbitrary subdecomposition of $Q$ into
cubes $J_{r}\in \mathcal{P}$ and interchanging the roles of $\sigma$ and $\omega$ we have the constant $\mathcal{E}_2^{\alpha,*}$. Also recall that $\mathfrak{E}_{2}^{\alpha }=\mathcal{E}_{2}^{\alpha}+\mathcal{E}_2^{\alpha,*}$. In the next definition we restrict the
cubes $Q$ to a dyadic grid $\mathcal{D}$, but keep the subcubes $J_{r}$ unrestricted.

\begin{dfn}
\label{def energy corona 3}Given a dyadic grid $\mathcal{D}$ and a cube 
$S_{0}\in \mathcal{D}$, define $\mathcal{S}\left( S_{0}\right) $ to be the 
\emph{maximal} $\mathcal{D}$-subcubes $I\subset S_{0}$ such that%
\begin{equation}
\sup_{I\supset \dot{\cup}J_{r}}\sum_{r=1}^{\infty }\left( \frac{\mathrm{P}
^{\alpha }\left( J_{r},\mathbf{1}_{I}\sigma \right) }{\left\vert
J_{r}\right\vert^\frac{1}{n} }\right) ^{2}\left\Vert x-m_{J_{r}}\right\Vert
_{L^{2}\left( \mathbf{1}_{J_{r}}\omega \right) }^{2}\geq C_{{en}}
\left[ \left( \mathfrak{E}_{2}^{\alpha }\right) ^{2}+\mathfrak{A}
_{2}^{\alpha }\right] \ \left\vert I\right\vert _{\sigma },
\label{def stop 3}
\end{equation}
where the cubes $J_{r}\in \mathcal{P}$ are pairwise disjoint in $I$, $\mathfrak{E}_{2}^{\alpha }$ is the energy condition constant, and $C_{{en}}$ is a sufficiently large positive constant depending only
on $\alpha $. Then define the $\sigma $-energy stopping cubes of $S_{0}$
to be the collection 
\begin{equation*}
\mathcal{F}=\left\{ S_{0}\right\} \cup \bigcup\limits_{m=0}^{\infty }
\mathcal{S}_{m}
\end{equation*}
where $\mathcal{S}_{0}=\mathcal{S}\left( S_{0}\right) $ and $\mathcal{S}
_{m+1}=\bigcup\limits_{S\in \mathcal{S}_{m}}\mathcal{S}\left( S\right) $
for $m\geq 0$.
\end{dfn}

We now claim that from the energy condition $\mathfrak{E}_{2}^{\alpha }<\infty$, we obtain the $\sigma $-Carleson estimate,%
\begin{equation}
\sum_{S\in \mathcal{S}:\ S\subset I}\left\vert S\right\vert _{\sigma }\leq
2\left\vert I\right\vert _{\sigma },\ \ \ \ \ I\in \mathcal{D}.
\label{sigma Carleson 3}
\end{equation}%
Indeed, for any $S_{1}\in \mathcal{F}$ we have%
\begin{eqnarray*}
\sum_{S\in \mathfrak{C}_{\mathcal{F}}\left( S_{1}\right) }\!\!\!\!\left\vert S\right\vert _{\sigma } 
\!\!\!\!\!\!&\leq &\!\!\!\!\!\!\!
\frac{1}{C_{{en}}\!\Big(\!\mathfrak{A}_2^\alpha\! +\! \!\left( \mathcal{E}_{2}^{\alpha }\right) ^{\!2}\!\!\Big)\!}\!\sum_{S\in \mathfrak{C}_{\mathcal{F}}\left(
S_{1}\right) }\!\sup_{S\supset \dot{\cup}J_{r}}\sum_{r=1}^{\infty }\!\left( 
\frac{\mathrm{P}^{\alpha }\!\left( J_{r},\mathbf{1}_{S}\sigma\! \right) }{\left\vert J_{r}\right\vert^\frac{1}{n} }\!\right)^{\!2}\!\!\!\left\Vert x\!-\!m_{J_{r}}\right\Vert^2
_{L^{2}\!(\mathbf{1}_{\!J_{r}}\omega)} \\
&\leq &\!\!\!\!\!\!
\frac{1}{C_{{en}}\left( \mathcal{E}_{2}^{\alpha }\right)
^{2}}\left( \mathcal{E}_{2}^{\alpha }\right) ^{2}\left\vert S_{1}\right\vert
_{\sigma }=\frac{1}{C_{{en}}}\left\vert S_{1}\right\vert
_{\sigma }\ ,
\end{eqnarray*}%
upon noting that the union of the subdecompositions $\dot{\cup}J_{r}\subset
S $ over $S\in \mathfrak{C}_{\mathcal{F}}\left( S_{1}\right) $ is a
subdecomposition of $S_{1}$, and the proof of the Carleson estimate is now
finished by iteration in the standard way.

Finally, we record the reason for introducing energy stopping times. If 
\begin{equation}
\mathbf{X}_{\alpha }\left( \mathcal{C}_{S}\right) ^{2}\equiv \sup_{I\in 
\mathcal{C}_{S}}\frac{1}{\left\vert I\right\vert _{\sigma }}\sup_{I\supset 
\dot{\cup}J_{r}}\sum_{r=1}^{\infty }\left( \frac{\mathrm{P}^{\alpha }\left(
J_{r},\mathbf{1}_{I}\sigma \right) }{\left\vert J_{r}\right\vert^\frac{1}{n} }\right)
^{2}\left\Vert x-m_{J_{r}}\right\Vert _{L^{2}\left( \mathbf{1}_{J_{r}}\omega
\right) }^{2}  \label{def stopping energy 3}
\end{equation}%
is (the square of) the $\alpha $\emph{-stopping energy} of the weight pair $\left( \sigma ,\omega \right) $ with respect to the corona $\mathcal{C}_{S}$
, then we have the \emph{stopping energy bounds}%
\begin{equation}
\mathbf{X}_{\alpha }\left( \mathcal{C}_{S}\right) \leq \sqrt{C_{{en}}}\sqrt{\left( \mathfrak{E}_{2}^{\alpha }\right) ^{2}+\mathfrak{A}_{2}^{\alpha }},\ \ \ \ \ S\in \mathcal{F},  \label{def stopping bounds 3}
\end{equation}
where $\mathfrak{A}_{2}^{\alpha }$ and the energy constant $\mathfrak{E}_{2}^{\alpha }$ are controlled by the assumptions in Theorem \ref{dim high}.

\subsection{Iterated coronas and general stopping data}

We will use a construction that permits \emph{iteration} of the above three
corona decompositions by combining Definitions \ref{CZ stopping times}, \ref%
{accretive stopping times gen} and \ref{def energy corona 3} into a single
stopping condition. However, there is one remaining difficulty with the
triple corona constructed in this way, namely if a stopping cube $I\in 
\mathcal{A}$ is a child of a cube $Q$ in the corona $\mathcal{C}_{A}$,
then the modulus of the average $\left\vert \frac{1}{\left\vert I\right\vert
_{\sigma }}\int_{I}b_{Q}d\sigma \right\vert $ of $b_{Q}$ on $I$ may be far
smaller than the sup norm of $\left\vert b_{Q}\right\vert $ on the child $I$%
, indeed it may be that $\frac{1}{\left\vert I\right\vert _{\sigma }}%
\int_{I}b_{Q}d\sigma =0$. This of course destroys any reasonable estimation
of the martingale and dual martingale differences $\bigtriangleup
_{Q}^{\sigma ,\mathbf{b}}f$ and $\square _{Q}^{\sigma ,\mathbf{b}}f$ used in
the proof of Theorem \ref{dim high}, and so we will use Lemma \ref{prelim
control of corona} on the function $b_{A}$ to obtain a new function $%
\widetilde{b}_{A}$ for which this problem is circumvented at the `bottom' of
the corona, i.e. for those $A^{\prime }\in \mathfrak{C}_{\mathcal{A}}\left( A\right) $.
We then refer to the stopping times $A^{\prime }\in \mathfrak{C}_{\mathcal{A}}\left(
A\right) $ as `shadow' stopping times since we have lost control of the weak
testing condition relative to the new function $\widetilde{b}_{A}$. Thus we
must redo the weak testing stopping times for the new function $\widetilde{b}%
_{A}$, but also stopping if we hit one of the shadow stopping times. Here
are the details.

\begin{dfn}
\label{def shadow}Let $C_{0}\geq 4$, $0<\gamma <1$ and $1<\Gamma <\infty $.
Suppose that $\mathbf{b}=\left\{ b_{Q}\right\} _{Q\in \mathcal{P}}$ is an $\infty $-weakly $\sigma $-accretive family on $\mathbb{R}^n$. Given a dyadic
grid $\mathcal{D}$ and a cube $Q\in \mathcal{D}$, define the collection
of `shadow' stopping times $\mathcal{S}_{{shadow}}\left( Q\right) $
to be the \emph{maximal} $\mathcal{D}$-subcubes $I\subset Q$ such that one of the following holds:
\begin{enumerate}[(a).]
    \item 
\begin{equation*}
\frac{1}{\left\vert I\right\vert _{\sigma }}\int_{I}\left\vert f\right\vert
d\sigma >C_{0}\frac{1}{\left\vert Q\right\vert _{\sigma }}\int_{Q}\left\vert
f\right\vert d\sigma \ ,
\end{equation*}%
\item
\begin{equation*}
\left\vert \frac{1}{\left\vert I\right\vert _{\mu }}\int_{I}b_{Q}d\sigma
\right\vert <\gamma \ \text{or }\int_{I}\left\vert T_{\sigma }^{\alpha
}\left( b_{Q}\right) \right\vert ^{2}d\omega >\Gamma \left( \mathfrak{T}%
_{T^{\alpha }}^{\mathbf{b}}\right) ^{2}\left\vert I\right\vert _{\sigma }%
\text{ },
\end{equation*}%
\item
\begin{equation*}
\sup_{I\supset \dot{\cup}J_{r}}\sum_{r=1}^{\infty }\left( \frac{\mathrm{P}%
^{\alpha }\left( J_{r},\sigma \right) }{\left\vert J_{r}\right\vert^\frac{1}{n} }\right)
^{2}\left\Vert x-m_{J_{r}}\right\Vert _{L^{2}\left( \mathbf{1}_{J_{r}}\omega
\right) }^{2}\geq C_{{en}}\left[ \left( \mathfrak{E}_{2}^{\alpha}\right) ^{2}+\mathfrak{A}_{2}^{\alpha }\right]
\ \left\vert I\right\vert _{\sigma }\ .
\end{equation*}
\end{enumerate}
\end{dfn}

Now we apply Lemma \ref{prelim control of corona}\ to the function $b_{Q}$
with $\mathcal{S}_{{shadow}}\left( Q\right)
\equiv \left\{ Q_{i}\right\} _{i=1}^{\infty }$ to obtain a new function $%
\widetilde{b}_{Q}$ satisfying the properties%
\begin{eqnarray}
&& 
\supp\widetilde{b}_{Q}\subset Q\ ,  \label{props} \\
&& 
1\leq \frac{1}{\left\vert Q^{\prime }\right\vert _{\sigma }}%
\int_{Q^{\prime }}\widetilde{b}_{Q}d\sigma \leq \left\Vert \mathbf{1}%
_{Q^{\prime }}\widetilde{b}_{Q}\right\Vert _{L^{\infty }\left( \sigma
\right) }\leq 2\left( 1+\sqrt{C_{\mathbf{b}}}\right) C_{\mathbf{b}}\ ,\ \ \
\ \ Q^{\prime }\in \mathcal{C}_{Q}\ ,  \notag \\
&&
\sqrt{\int_{Q}\left\vert T_{\sigma }^{\alpha }b_{Q}\right\vert ^{2}d\omega 
}\leq \left[ 2\mathfrak{T}_{T^{\alpha }}^{\mathbf{b}}\left( Q\right) +4C_{\mathbf{b}}^{\frac{3}{2}}\delta ^{\frac{1}{4}}\mathfrak{N}_{T^{\alpha
}}\left( Q\right) \right] \sqrt{\left\vert Q\right\vert _{\sigma }}\ , 
\notag \\
&&
\left\Vert \mathbf{1}_{Q_{i}}\widetilde{b}_{Q}\right\Vert _{L^{\infty
}\left( \sigma \right) }\leq \frac{16C_{\mathbf{b}}}{\delta }\left\vert 
\frac{1}{\left\vert Q_{i}\right\vert _{\sigma }}\int_{Q_{i}}\widetilde{b}%
_{Q}d\sigma \right\vert \ ,\ \ \ \ \ 1\leq i<\infty .  \notag
\end{eqnarray}%
Note that each of the functions $\widetilde{b}_{Q^{\prime }}\equiv \mathbf{1}%
_{Q^{\prime }}\widetilde{b}_{Q}$, for $Q^{\prime }\in \mathcal{C}_{Q}$, now
satisfies the crucial reverse H\"{o}lder property%
\begin{equation*}
\left\Vert \mathbf{1}_{I}\widetilde{b}_{Q^{\prime }}\right\Vert _{L^{\infty
}\left( \sigma \right) }\leq C_{\delta ,\mathbf{b}}\left\vert \frac{1}{%
\left\vert I\right\vert _{\sigma }}\int_{I}\widetilde{b}_{Q^{\prime
}}d\sigma \right\vert \ ,\ \ \ \ \ \text{for all }I\in \mathfrak{C}\left(
Q^{\prime }\right) ,\ Q^{\prime }\in \mathcal{C}_{Q}.
\end{equation*}%
Indeed, if $I$ equals one of the $Q_{i}$ then the reverse H\"{o}lder
condition in the last line of (\ref{props}) applies, while if $I\in \mathcal{%
C}_{Q}$ then the accretivity in the second line of (\ref{props}) applies.

Since we have lost the weak testing condition in the corona for this new
function $\widetilde{b}_{Q}$, the next step is to run again the weak testing
construction of stopping times, but this time starting with the new function 
$\widetilde{b}_{Q}$, and also stopping if we hit one of the `shadow'
stopping times $Q_{i}$. Here is the new stopping criterion.

\begin{dfn}
\label{def iterated}Let $C_{0}\geq 4$ and $1<\Gamma <\infty $. Let $\mathcal{%
S}_{{shadow}}\left( Q\right) \equiv \left\{ Q_{i}\right\}
_{i=1}^{\infty }$ be as in Definition \ref{def shadow}. Define $\mathcal{S}_{
{iterated}}\left( Q\right) $ to be the \emph{maximal} $\mathcal{D}$%
-subcubes $I\subset Q$ such that either%
\begin{equation*}
\int_{I}\left\vert T_{\sigma }^{\alpha }\left( \widetilde{b}_{Q}\right)
\right\vert ^{2}d\omega >\Gamma \left( \mathfrak{T}_{T^{\alpha }}^{\widetilde{\mathbf{b}}}\right) ^{2}\left\vert I\right\vert _{\sigma }\text{ },
\end{equation*}%
or
\begin{equation*}
I=Q_{i}\text{ for some }1\leq i<\infty .
\end{equation*}
\end{dfn}

Thus for each cube $Q$ we have now constructed \emph{iterated stopping
children} $\mathcal{S}_{{iterated}}\left( Q\right) $ by first
constructing shadow stopping times $\mathcal{S}_{{shadow}}\left(
Q\right) $ using one step of the triple corona construction, then modifying
the testing function to have reverse H\"{o}lder controlled children, and
finally running again the weak testing stopping time construction to get $\mathcal{S}_{{iterated}}\left( Q\right) $. These iterated stopping
times $\mathcal{S}_{{iterated}}\left( Q\right) $ have control of CZ
averages of $f$ and energy control of $\sigma $ and $\omega $, simply
because these controls were achieved in the shadow construction, and were
unaffected by either the application of Lemma \ref{prelim control of corona}
or the rerunning of the weak testing stopping criterion for $\widetilde{b}_{Q}$. And of course we now have weak testing within the corona determined
by $Q$ and $\mathcal{S}_{{iterated}}\left( Q\right) $, and we also
have the crucial reverse H\"{o}lder condition on all the children of
cubes in the corona. With all of this in hand, here then is the
definition of the construction of iterated coronas.

\begin{dfn}
\label{iterated stopping times}Let $C_{0}\geq 4$, $0<\gamma <1$ and $%
1<\Gamma <\infty $. Suppose that $\mathbf{b}=\left\{ b_{Q}\right\} _{Q\in 
\mathcal{P}}$ is an $\infty $-weakly $\sigma $-accretive family on $%
\mathbb{R}^n$. Given a dyadic grid $\mathcal{D}$ and a cube $S_{0}$ in $%
\mathcal{D}$, define the iterated stopping cubes of $S_{0}$ to be the
collection 
\begin{equation*}
\mathcal{F}=\left\{ S_{0}\right\} \cup \bigcup\limits_{m=0}^{\infty }
\mathcal{S}_{m}
\end{equation*}
where $\mathcal{S}_{0}=\mathcal{S}_{{iterated}}\left( S_{0}\right) $
and $\mathcal{S}_{m+1}=\bigcup\limits_{S\in \mathcal{S}_{m}}\mathcal{S}_{{iterated}}\left( S\right) $ for $m\geq 0$, and where $\mathcal{S}_{{iterated}}\left( Q\right) $ is defined in Definition \ref{def
iterated}.
\end{dfn}

It is useful to append to the notion of stopping times $\mathcal{S}$ in the
above $\sigma $-iterated corona decomposition a positive constant $A_{0}$
and an additional structure $\alpha _{\mathcal{S}}$ called stopping bounds
for a function $f$. We will refer to the resulting\ triple $\left( A_{0},\mathcal{F},\alpha _{\mathcal{F}}\right) $ as constituting stopping data for 
$f$. If $\mathcal{F}$ is a grid, we define $F^{\prime }\prec F$ if $F^{\prime }\subsetneqq F$ and $F^{\prime },F\in \mathcal{F}$. Recall that $\pi _{\mathcal{F}}F^{\prime }$ is the smallest $F\in \mathcal{F}$ such that $F^{\prime }\prec F$.

Suppose we are given a positive constant $A_{0}\geq 4$, a subset $\mathcal{F}$ of the dyadic grid $\mathcal{D}$
(called the stopping times), and a corresponding sequence $\alpha _{\mathcal{F}}\equiv \left\{ \alpha _{\mathcal{F}}\left( F\right) \right\} _{F\in 
\mathcal{F}}$ of nonnegative numbers $\alpha _{\mathcal{F}}\left( F\right)
\geq 0$ (called the stopping bounds). Let $\left( \mathcal{F},\prec ,\pi _{\mathcal{F}}\right) $ be the tree structure on $\mathcal{F}$ inherited from $\mathcal{D}$, and for each $F\in \mathcal{F}$ denote by $\mathcal{C}%
_{F}=\left\{ I\in \mathcal{D}:\pi _{\mathcal{F}}I=F\right\} $ the corona
associated with $F$: 
\begin{equation*}
\mathcal{C}_{F}=\left\{ I\in \mathcal{D}:I\subset F\text{ and }I\not\subset
F^{\prime }\text{ for any }F^{\prime }\prec F\right\} .
\end{equation*}%

\begin{dfn}
\label{general stopping data}
We say the triple $\left( A_{0},\mathcal{F},\alpha _{\mathcal{F}}\right) $
constitutes \emph{stopping data} for a function $f\in L_{loc}^{1}\left(
\sigma \right) $ if

\begin{enumerate}[(1).]
\item ${E}_{I}^{\sigma }\left\vert f\right\vert \leq \alpha _{\mathcal{F}}\left( F\right) $ for all $I\in \mathcal{C}_{F}$ and $F\in 
\mathcal{F}$,

\item $\sum_{F^{\prime }\preceq F}\left\vert F^{\prime }\right\vert _{\sigma
}\leq A_{0}\left\vert F\right\vert _{\sigma }$ for all $F\in \mathcal{F}$,

\item $\sum_{F\in \mathcal{F}}\alpha _{\mathcal{F}}\left( F\right)
^{2}\left\vert F\right\vert _{\sigma }\mathbf{\leq }A_{0}^{2}\left\Vert
f\right\Vert _{L^{2}\left( \sigma \right) }^{2}$,

\item $\alpha _{\mathcal{F}}\left( F\right) \leq \alpha _{\mathcal{F}}\left(
F^{\prime }\right) $ whenever $F^{\prime },F\in \mathcal{F}$ with $F^{\prime
}\subset F$.
\end{enumerate}
\end{dfn}

Property \textit{(1)} says that $\alpha _{\mathcal{F}}\left( F\right) $ bounds the
averages of $f$ in the corona $\mathcal{C}_{F}$, and property \textit{(2)} says that
the cubes at the tops of the coronas satisfy a Carleson condition relative
to the weight $\sigma $. Note that a standard `maximal cube' argument
extends the Carleson condition in property \textit{(2)} to the inequality%
\begin{equation}
\sum_{F^{\prime }\in \mathcal{F}:\ F^{\prime }\subset A}\left\vert F^{\prime
}\right\vert _{\sigma }\leq A_{0}\left\vert A\right\vert _{\sigma }\text{
for all open sets }A\subset \mathbb{R}^n.  \label{Car ext}
\end{equation}%
Property \textit{(3)} is the quasi-orthogonality condition that says the sequence of
functions $\left\{ \alpha _{\mathcal{F}}\left( F\right) \mathbf{1}%
_{F}\right\} _{F\in \mathcal{F}}$ is in the vector-valued space $L^{2}\left(
\ell ^{2};\sigma \right) $ with control and is often referred to as a Carleson embedding theorem, and property \textit{(4)} says that the control on
stopping data is nondecreasing on the stopping tree $\mathcal{F}$. We
emphasize that we are \emph{not} assuming in this definition the stronger
property that there is $C>1$ such that $\alpha _{\mathcal{F}}\left(
F^{\prime }\right) >C\alpha _{\mathcal{F}}\left( F\right) $ whenever $%
F^{\prime },F\in \mathcal{F}$ with $F^{\prime }\subsetneqq F$. Instead, the
properties \textit{(2)} and \textit{(3)} substitute for this lack. Of course the stronger
property \emph{does} hold for the familiar \emph{Calder\'{o}n-Zygmund}
stopping data determined by the following requirements for $C>1$,%
\begin{equation*}
E_{F^{\prime }}^{\sigma }\left\vert f\right\vert >C E_{F}^{\sigma }\left\vert f\right\vert \text{ whenever }F^{\prime },F\in 
\mathcal{F}\text{ with }F^{\prime }\subsetneqq F,
\end{equation*}
\begin{equation*}
E_{I}^{\sigma }\left\vert f\right\vert \leq C E_{F}^{\sigma
}\left\vert f\right\vert \text{ for }I\in \mathcal{C}_{F},
\end{equation*}%
which are themselves sufficiently strong to automatically force properties
\textit{(2)} and \textit{(3)} with $\alpha _{\mathcal{F}}\left( F\right) =\mathbb{E}%
_{F}^{\sigma }\left\vert f\right\vert $.

We have the following useful consequence of \textit{(2)} and \textit{(3)} that says the
sequence $\left\{ \alpha _{\mathcal{F}}\left( F\right) \mathbf{1}%
_{F}\right\} _{F\in \mathcal{F}}$ has a \emph{quasi-orthogonal} property
relative to $f$ with a constant $C_{0}^{\prime }$ depending only on $C_{0}$:%
\begin{equation}
\left\Vert \sum_{F\in \mathcal{F}}\alpha _{\mathcal{F}}\left( F\right) 
\mathbf{1}_{F}\right\Vert _{L^{2}\left( \sigma \right) }^{2}\leq
C_{0}^{\prime }\left\Vert f\right\Vert _{L^{2}\left( \sigma \right) }^{2}.
\label{q orth}
\end{equation}%
$\ $

\begin{prop}
\label{data}Let $f\in L^{2}\left( \sigma \right) $, let $\mathcal{F}$ be as in Definition \ref{iterated stopping times}, and define stopping data $\alpha _{\mathcal{F}}$
by $\alpha _{F}=\frac{1}{\left\vert F\right\vert _{\sigma }}
\int_{F}\left\vert f\right\vert d\sigma $. Then there is $A_{0}\geq 4$,
depending only on the constant $C_{0}$ in Definition \ref{CZ stopping times}
, such that the triple $\left( A_{0},\mathcal{F},\alpha _{\mathcal{F}%
}\right) $ constitutes \emph{stopping data} for the function $f$.
\end{prop}

\begin{proof}
This is an easy exercise using  (\ref{Car gen}) and (\ref{sigma Carleson 3}), and is left for the reader.
\end{proof}

\subsection{Reduction to good functions}

We begin with a specification of the various parameters that will arise
during the proof, as well as the extension of goodness introduced in \cite%
{HyMa}.

\begin{dfn}
\label{def parameters}The parameters $\mathbf{r}$, $\mathbf{\tau }$ and $%
\mathbf{\rho }$ will be fixed below to satisfy 
\begin{equation*}
\mathbf{\tau }>\mathbf{r}\text{ and }\mathbf{\rho }>\mathbf{r}+\mathbf{\tau }%
,
\end{equation*}%
where $\mathbf{r}$ is the goodness parameter fixed in (\ref{choice of r}).
\end{dfn}

Let $0<\varepsilon <1$ to be chosen later. Define $J$ to be $\varepsilon -%
{good}$ in a cube $K$ if 
\begin{equation*}
d\left( J,{skel}K\right) >2\left\vert J\right\vert ^{\varepsilon
}\left\vert K\right\vert ^{1-\varepsilon },
\end{equation*}%
where the skeleton ${skel}K\equiv \bigcup\limits_{K^{\prime }\in 
\mathfrak{C}\left( K\right) }\partial K^{\prime }$ of a cube $K$ consists of
the boundaries of all the children $K^{\prime }$ of $K$. Define $\mathcal{G}%
_{\left( k,\varepsilon \right) -{good}}^{\mathcal{D}}$ to consist of
those $J\in \mathcal{G}$ such that $J$ is good in every supercube $K\in 
\mathcal{D}$ that lies at least $k$ levels above $J$. We also define $J$ to
be $\varepsilon -{good}$ in a cube $K$ and beyond if $J\in \mathcal{%
G}_{\left( k,\varepsilon \right) -{good}}^{\mathcal{D}}$ where $%
k=\log _{2}\frac{\ell \left( K\right) }{\ell \left( J\right) }$. We can now
say that $J\in \mathcal{G}_{\left( k,\varepsilon \right) -{good}}^{%
\mathcal{D}}$ if and only if $J$ is $\varepsilon -{good}$ in $\pi
^{k}J$ and beyond. As the goodness parameter $\varepsilon $ will eventually
be fixed throughout the proof, we sometimes suppress it, and simply say "$J$
is ${good}$ in a cube $K$ and beyond" instead of "$J$ is $%
\varepsilon -{good}$ in a cube $K$ and beyond".

As pointed out on page 14 of \cite{HyMa} by Hyt\"{o}nen and Martikainen,
there are subtle difficulties associated in using dual martingale
decompositions of functions which depend on the entire dyadic grid, rather
than on just the local cube in the grid. We will proceed at first in the
spirit of \cite{HyMa}. The goodness that we will infuse below into the main
`below' form $\mathsf{B}_{\Subset _{\mathbf{\rho }}}\left( f,g\right) $ will
be the Hyt\"{o}nen-Martikainen `weak' goodness: every pair $\left(
I,J\right) \in \mathcal{D}\times \mathcal{G}$ that arises in the form $%
\mathsf{B}_{\Subset _{\mathbf{\rho }}}\left( f,g\right) $ will satisfy $J\in 
\mathcal{G}_{\left( k,\varepsilon \right) -{good}}^{\mathcal{D}}$
where $\ell \left( I\right) =2^{k}\ell \left( J\right) $.

It is important to use \emph{two} independent random grids, one for each
function $f$ and $g$ simultaneously, as this is necessary in order to apply
probabilistic methods to the dual martingale averages $\square _{I}^{\mu ,%
\mathbf{b}}$ that depend, not only on $I$, but also on the underlying grid
in which $I$ lives. The proof methods for functional energy from \cite%
{SaShUr7} and \cite{SaShUr6} relied heavily on the use of a single grid, and
this must now be modified to accomodate two independent grids.

\subsubsection{Parameterizations of dyadic grids}

It is important to use two independent grids, one for each function $f$ and $g$ simultaneously, as it is necessary in order to apply probabilistic methods to the dual martingale averages $\square_I^{
\mu,\mathbf{b}}$ that depend not only on $I$ but also on the underlying grid in which $I$ lives.

Now we recall the construction from the paper \cite{SaShUr10}. We momentarily fix a large
positive integer $M\in \mathbb{N}$, and consider the tiling of $\mathbb{R}^n$
by the family of cubes $\mathbb{D}_{M}\equiv \left\{ I_{\alpha
}^{M}\right\} _{\alpha \in \mathbb{Z}}$ having side length $2^{-M}$ and
given by $I_{\alpha }^{M}\equiv I_{0}^{M}+\alpha\cdot 2^{-M}$ where $I_{0}^{M}=%
\left[ 0,2^{-M}\right) $. A \emph{dyadic grid} $\mathcal{D}$ built on $%
\mathbb{D}_{M}$ is\ defined to be a family of cubes $\mathcal{D}$
satisfying:

\begin{enumerate}
\item Each $I\in \mathcal{D}$ has side length $2^{-\ell }$ for some $\ell
\in \mathbb{Z}$ with $\ell \leq M$, and $I$ is a union of $2^{n(M-\ell)}$
cubes from the tiling $\mathbb{D}_{M}$,

\item For $\ell \leq M$, the collection $\mathcal{D}_{\ell }$ of cubes
in $\mathcal{D}$ having side length $2^{-\ell }$ forms a pairwise disjoint
decomposition of the space $\mathbb{R}^n$,

\item Given $I\in \mathcal{D}_{i}$ and $J\in \mathcal{D}_{j}$ with $j\leq
i\leq M$, it is the case that either $I\cap J=\emptyset $ or $I\subset J$.
\end{enumerate}

We now momentarily fix a \emph{negative} integer $N\in -\mathbb{N}$, and
restrict the above grids to cubes of side length at most $2^{-N}$:%
\begin{equation*}
\mathcal{D}^{N}\equiv \left\{ I\in \mathcal{D}:\text{side length of }I\text{
is at most }2^{-N}\right\} \text{.}
\end{equation*}%
We refer to such grids $\mathcal{D}^{N}$ as a (truncated) dyadic grid $%
\mathcal{D}$ built on $\mathbb{D}_{M}$ of size $2^{-N}$. There are now two
traditional means of constructing probability measures on collections of
such dyadic grids, namely parameterization by choice of parent, and
parameterization by translation.

\textbf{Construction \#1}: Consider first the special case of dimension $n=1$. For any 
\begin{equation*}
\beta =\{\beta _{i}\}_{i\in \mathbb{Z}_{M}^{N}}\in \omega _{m}^{N}\equiv \left\{
0,1\right\} ^{\mathbb{Z}_{M}^{N}},
\end{equation*}%
where $\mathbb{Z}_{M}^{N}\equiv \left\{ \ell \in \mathbb{Z}:N\leq \ell \leq
M\right\} $, define the dyadic grid $\mathcal{D}_{\beta }$ built on $\mathbb{%
D}_{m}$ of size $2^{-N}$ by 
\begin{equation}
\mathcal{D}_{\beta }=\left\{ 2^{-\ell }\left( [0,1)+k+\sum_{i:\ \ell <i\leq
M}2^{-i+\ell }\beta _{i}\right) \right\} _{N\leq \ell \leq M,\,k\in {\mathbb{Z}}}  \label{def dyadic grid}
\end{equation}%
Place the uniform probability measure $\rho _{M}^{N}$ on the finite index
space $\omega _{M}^{N}=\left\{ 0,1\right\} ^{\mathbb{Z}_{M}^{N}}$, namely
that which charges each $\beta \in \omega _{M}^{N}$ equally. This
construction is then extended to Euclidean space $\mathbb{R}^n$ by taking
products in the usual way and using the product index space $\Omega
_{M}^{N}\equiv (\omega _{M}^{N})^n$ and the uniform product probability measure $%
\mu _{M}^{N}=\rho _{M}^{N}\times ...\times \rho _{M}^{N}$.

\textbf{Construction \#2}: Momentarily fix a (truncated) dyadic grid $%
\mathcal{D}$ built on $\mathbb{D}_{M}$ of size $2^{-N}$. For any 
\begin{equation*}
\gamma \in \Gamma _{M}^{N}\equiv \left\{ 2^{-M}\mathbb{Z}^n_{+}:\left\vert
\gamma _{i}\right\vert <2^{-N}\right\} ,
\end{equation*}%
where $\mathbb{Z}^{n}_+=(\mathbb{N}\cup\{0\})^n$, define the dyadic grid $\mathcal{D}^{\gamma }$ built on $\mathbb{D}_{m}$ of
size $2^{-N}$ by%
\begin{equation*}
\mathcal{D}^{\gamma }\equiv \mathcal{D}+\gamma .
\end{equation*}%
Place the uniform probability measure $\nu _{M}^{N}$ on the finite index set 
$\Gamma _{M}^{N}$, namely that which charges each multiindex $\gamma $ in $%
\Gamma _{M}^{N}$ equally.

The two probability spaces $\left( \left\{ \mathcal{D}_{\beta }\right\}
_{\beta \in \Omega _{M}^{N}},\mu _{M}^{N}\right) $ and $\left( \left\{ 
\mathcal{D}^{\gamma }\right\} _{\gamma \in \Gamma _{M}^{N}},\nu
_{M}^{N}\right) $ are isomorphic since both collections $\left\{ \mathcal{D}%
_{\beta }\right\} _{\beta \in \Omega _{M}^{N}}$ and $\left\{ \mathcal{D}%
^{\gamma }\right\} _{\gamma \in \Gamma _{M}^{N}}$ describe the set $%
\boldsymbol{A}_{M}^{N}$ of \textbf{all} (truncated) dyadic grids $\mathcal{D}%
^{\gamma }$ built on $\mathbb{D}_{m}$ of size $2^{-N}$, and since both
measures $\mu _{M}^{N}$ and $\nu _{M}^{N}$ are the uniform measure on this
space. The first construction may be thought of as being \emph{parameterized
by scales} - each component $\beta _{i}$ in $\beta =\{\beta _{i}\}_{i\in
\mathbb{Z}_{M}^{N}}\in \omega _{M}^{N}$ amounting to a choice of the two possible
tilings at level $i$ that respect the choice of tiling at the level below -
and since any grid in $\boldsymbol{A}_{M}^{N}$ is determined by a choice of
scales , we see that $\left\{ \mathcal{D}_{\beta }\right\} _{\beta \in
\Omega _{M}^{N}}=\boldsymbol{A}_{M}^{N}$. The second construction may be
thought of as being \emph{parameterized by translation} - each $\gamma \in
\Gamma _{M}^{N}$ amounting to a choice of translation of the grid $\mathcal{D%
}$ fixed in construction \#2\ - and since any grid in $\boldsymbol{A}_{M}^{N}
$ is determined by any of the cubes at the top level, i.e. with side
length $2^{-N}$, we see that $\left\{ \mathcal{D}^{\gamma }\right\} _{\gamma
\in \Gamma _{M}^{N}}=\boldsymbol{A}_{M}^{N}$ as well, since every cube
at the top level in $\boldsymbol{A}_{M}^{N}$ has the form $Q+\gamma $ for
some $\gamma \in \Gamma _{M}^{N}$ and $Q\in \mathcal{D}$ at the top level in 
$\boldsymbol{A}_{M}^{N}$ (i.e. every cube at the top level in $%
\boldsymbol{A}_{M}^{N}$ is a union of small cubes in $\mathbb{D}_{m}$,
and so must be a translate of some $Q\in \mathcal{D}$ by an amount $2^{-M}$
times an element of $\mathbb{Z}_{+}$). Note also that in all dimensions, $%
\#\Omega _{M}^{N}=\#\Gamma _{M}^{N}=2^{n(M-N)}$. We will use $\boldsymbol{E}_{\Omega _{M}^{N}}$ to denote expectation with respect to this common
probability measure on $\boldsymbol{A}_{M}^{N}$.

\begin{notation}
\label{suppress M and N}For purposes of notation and clarity, we now
suppress all reference to $M$ and $N$ in our families of grids, and in the
notations $\Omega $ and $\Gamma $ for the parameter sets, and we use $%
\boldsymbol{P}_{\Omega }$ and $\boldsymbol{E}_{\Omega }$ to denote
probability and expectation with respect to families of grids, and instead
proceed as if all grids considered are unrestricted. The careful reader can
supply the modifications necessary to handle the assumptions made above on
the grids $\mathcal{D}$ and the functions $f$ and $g$ regarding $M$ and $N$.
\end{notation}

\subsection{Formulas for martingale averages} We need the following formulas defined on Appendix A of \cite{SaShUr12}.
\begin{eqnarray}
\mathbb{E}_{Q}^{\mu ,\mathbf{b}}f\left( x\right) &\equiv &\mathbf{1}%
_{Q}\left( x\right) \frac{1}{\int_{Q}b_{Q}d\mu }\int_{Q}fb_{Q}d\mu ,\ \ \ \
\ Q\in \mathcal{P}\ ,\label{def expectation} \\
\mathbb{F}_{Q}^{\mu ,\mathbf{b}}f\left( x\right) &\equiv &\mathbf{1}%
_{Q}\left( x\right) b_{Q}\left( x\right) \frac{1}{\int_{Q}b_{Q}d\mu }%
\int_{Q}fd\mu ,\ \ \ \ \ Q\in \mathcal{P}\ ,\notag
\end{eqnarray}
\begin{equation}
\widehat{\mathbb{F}}_{Q}^{\mu ,\mathbf{b}}f\left( x\right) \equiv \mathbf{1}%
_{Q}\left( x\right) \frac{1}{\int_{Q}b_{Q}d\mu }\int_{Q}fd\mu ,\ \ \ \ \
Q\in \mathcal{P}\ .  \label{F hat}
\end{equation}%
and 
\begin{eqnarray}
\bigtriangleup _{Q}^{\mu ,\mathbf{b}}f\left( x\right) \!\!\!\!\!&\equiv &\!\!\!\!\!
\left(\sum_{Q^{\prime }\in \mathfrak{C}\left( Q\right) }\mathbb{E}_{Q^{\prime
}}^{\mu ,\mathbf{b}}f\left( x\right) \!\!\right) \!\!-\mathbb{E}_{Q}^{\mu ,\mathbf{b}%
}f\left( x\right)
= 
\!\!\!\sum_{Q^{\prime }\in \mathfrak{C}\left( Q\right) }\!\!\!\!\!\!\mathbf{1}_{Q^{\prime }}(x) \left( \mathbb{E}_{Q^{\prime }}^{\mu ,\mathbf{b}}f\!\left( x\right) \!-\!\mathbb{E}_{Q}^{\mu ,\mathbf{b}}f\!\left(
x\right) \right)  \label{def diff} \\
\square _{Q}^{\mu ,\mathbf{b}}f\left( x\right) 
\!\!\!\!\!&\equiv &\!\!\!\!\!
\left(
\sum_{Q^{\prime }\in \mathfrak{C}\left( Q\right) }\mathbb{F}_{Q^{\prime
}}^{\mu ,\mathbf{b}}f\left( x\right)\!\! \right) \!\!-\mathbb{F}_{Q}^{\mu ,\mathbf{b}}f\left( x\right) 
=
\!\!\!\sum_{Q^{\prime }\in \mathfrak{C}\left( Q\right) }\!\!\!\!\!\!\mathbf{1}_{Q^{\prime }}\left( x\right) \left( \mathbb{F}_{Q^{\prime }}^{\mu,\mathbf{b}}f\left( x\right) -\mathbb{F}_{Q}^{\mu ,\mathbf{b}}f\left(
x\right) \right)\notag
\end{eqnarray}
We also need
\begin{eqnarray}
\bigtriangledown _{Q}^{\mu }f 
&\equiv &
\sum_{Q^{\prime }\in \mathfrak{C}_{%
{brok}}\left( Q\right) }\left( \frac{1}{\left\vert Q^{\prime
}\right\vert _{\mu }}\int_{Q^{\prime }}\left\vert f\right\vert d\mu \right) 
\mathbf{1}_{Q^{\prime }},  \label{Carleson avg op} \\
\widehat{\bigtriangledown }_{Q}^{\mu }f 
&\equiv &
\sum_{Q^{\prime }\in 
\mathfrak{C}_{{brok}}\left( Q\right) }\left( \frac{1}{\left\vert
Q^{\prime }\right\vert _{\mu }}\int_{Q^{\prime }}\left\vert f\right\vert
d\mu +\frac{1}{\left\vert Q\right\vert _{\mu }}\int_{Q}\left\vert
f\right\vert d\mu \right) \mathbf{1}_{Q^{\prime }},  \notag
\end{eqnarray}
\begin{equation}
\sum_{Q\in \mathcal{D}}\left\Vert \widehat{\bigtriangledown }_{Q}^{\mu
}f\right\Vert _{L^{2}\left( \mu \right) }^{2}\lesssim \left\Vert
f\right\Vert _{L^{2}\left( \mu \right) }^{2}\ .  \label{Car embed}
\end{equation}%
and
\begin{eqnarray} 
\square _{Q}^{\mu ,\pi ,\mathbf{b}}f &=&\left[ \sum_{Q^{\prime }\in 
\mathfrak{C}\left( Q\right) }\mathbb{F}_{Q^{\prime }}^{\mu ,\pi ,\mathbf{b}}f%
\right] -\mathbb{F}_{Q}^{\mu ,\mathbf{b}}f=\sum_{Q^{\prime }\in \mathfrak{C}%
\left( Q\right) }\mathbb{F}_{Q^{\prime }}^{\mu ,b_{Q}}f-\mathbb{F}_{Q}^{\mu
,b_{Q}}f,  \label{def pi box} \\
\mathbb{F}_{Q}^{\mu ,\pi ,\mathbf{b}}f &=&\mathbf{1}_{Q}\frac{b_{\pi Q}}{%
\int_{Q}b_{\pi Q}d\mu }\int_{Q}fd\mu ,\label{def pi exp} \\
\square _{Q}^{\mu ,\mathbf{b}} &=&\square _{Q}^{\mu ,\pi ,\mathbf{b}}\square
_{Q}^{\mu ,\pi ,\mathbf{b}}+\square _{Q,{brok}}^{\mu ,\mathbf{b}}  \text{ and }%
\square _{Q}^{\mu ,\mathbf{b}}=\square _{Q}^{\mu ,\pi ,\mathbf{b}}+\square
_{Q,{brok}}^{\mu ,\pi ,\mathbf{b}}  \label{box pi equals} \\
\square _{Q,{brok}}^{\mu ,\mathbf{b}}f&=&\sum_{Q^{\prime
}\in \mathfrak{C}_{{brok}}\left( Q\right) }\mathbb{F}_{Q^{\prime
}}^{\mu ,b_{Q^{\prime }}}f-\mathbb{F}_{Q^{\prime }}^{\mu ,b_{Q}}f,\notag \\
\left\vert \square _{Q,{brok}}^{\mu ,\pi
,\mathbf{b}}f\right\vert &\lesssim &\left\vert \widehat{\bigtriangledown }%
_{Q}^{\mu }f\right\vert ,  \label{F est}
\end{eqnarray}%
with similar equalities and inequalities for $\bigtriangleup $ and $\mathbb{E%
}$. Here $\mathfrak{C}_{{brok}}\left( Q\right) $ denotes the set
of broken children, i.e. those $Q^{\prime }\in \mathfrak{C}\left( Q\right) $
for which $b_{Q^{\prime }}\neq \mathbf{1}_{Q^{\prime }}b_{Q}$, and more
generally and typically, $\mathfrak{C}_{{brok}}\left( Q\right) =%
\mathfrak{C}\left( Q\right) \cap \mathcal{A}$ where $\mathcal{A}$ is a
collection of stopping cubes that includes the broken children and
satisfies a $\sigma $-Carleson condition and $\pi Q$ is the dyadic father of $Q$.

Define another modified dual martingale difference by 
\begin{equation}
\square _{I}^{\sigma ,\flat ,\mathbf{b}}f\equiv \square _{I}^{\sigma ,%
\mathbf{b}}f-\sum_{I^{\prime }\in \mathfrak{C}_{{brok}}\left(
I\right) }\mathbb{F}_{I^{\prime }}^{\sigma ,\mathbf{b}}f=\left(
\sum_{I^{\prime }\in \mathfrak{C}_{{nat}}\left( I\right) }%
\mathbb{F}_{I^{\prime }}^{\sigma ,\mathbf{b}}f\right) -\mathbb{F}%
_{I}^{\sigma ,\mathbf{b}}f,  \label{flat box}
\end{equation}%
where we have removed the averages over broken children from $\square
_{I}^{\sigma ,\mathbf{b}}f$, but left the average over $I$ intact. On any
child $I^{\prime }$ of $I$, the function $\square _{I}^{\sigma ,\flat ,%
\mathbf{b}}f$ is thus a constant multiple of $b_{I}$, and so we have%
\begin{eqnarray}
\hspace{1cm}\square _{I}^{\sigma ,\flat ,\mathbf{b}}f
&=&
b_{I}\sum_{I^{\prime }\in 
\mathfrak{C}\left( I\right) }\mathbf{1}_{I^{\prime }}E_{I^{\prime }}^{\sigma
}\left( \frac{1}{b_{I}}\square _{I}^{\sigma ,\flat ,\mathbf{b}}f\right)
=b_{I}\ \sum_{I^{\prime }\in \mathfrak{C}\left( I\right) }\mathbf{1}%
_{I^{\prime }}E_{I^{\prime }}^{\sigma }\left( \widehat{\square }_{I}^{\sigma
,\flat ,\mathbf{b}}f\right) ;  \label{flat box hat} \\
\widehat{\square }_{I}^{\sigma ,\flat ,\mathbf{b}}f &\equiv &
\sum_{I^{\prime
}\in \mathfrak{C}\left( I\right) }\mathbf{1}_{I^{\prime }}\ E_{I^{\prime
}}^{\sigma }\left( \frac{1}{b_{I}}\square _{I}^{\sigma ,\flat ,\mathbf{b}%
}f\right),  \notag\\
&=&
\!\!\!\!
\sum_{I^{\prime }\in 
\mathfrak{C}_{{nat}}\left( I\right) }\!\!\!\!\mathbf{1}_{I^{\prime }}\!\!
\left[ \frac{1}{\int_{I^{\prime }}b_{I}d\mu }\int_{I^{\prime }}fd\mu \!-\!\frac{1%
}{\int_{I}b_{I}d\mu }\int_{I}fd\mu \right] 
-\!\!\!
\sum_{I^{\prime }\in \mathfrak{C}%
_{{brok}}\left( I\right) }\!\!\!\!\!\!\mathbf{1}_{I^{\prime }}\!\! \left[ \frac{1}{\int_{I}b_{I}d\mu }\int_{I}fd\mu \right]\notag
\end{eqnarray}%

Thus for $I\in \mathcal{C}_{A}$ we have 
\begin{equation}
\square _{I}^{\sigma ,\flat ,\mathbf{b}}f=b_{A}\sum_{I^{\prime }\in 
\mathfrak{C}\left( I\right) }\mathbf{1}_{I^{\prime }}E_{I^{\prime }}^{\sigma
}\left( \widehat{\square }_{I}^{\sigma ,\flat ,\mathbf{b}}f\right) =b_{A}%
\widehat{\square }_{I}^{\sigma ,\flat ,\mathbf{b}}f,  \label{factor b_A}
\end{equation}%
where the averages $E_{I^{\prime }}^{\sigma }\left( \widehat{\square }%
_{I}^{\sigma ,\flat ,\mathbf{b}}f\right) $ satisfy the following telescoping
property for all $K\in \left( \mathcal{C}_{A}\setminus \left\{ A\right\}
\right) \cup \left( \bigcup_{A^{\prime }\in \mathfrak{C}_{\mathcal{A}}\left(
A\right) }A^{\prime }\right) $ and $L\in \mathcal{C}_{A}$ with $K\subset L$:%
\begin{equation}
\sum_{I:\ \pi K\subset I\subset L}E_{I_{K}}^{\sigma }\left( \widehat{\square 
}_{I}^{\sigma ,\flat ,\mathbf{b}}f\right) =\left\{ 
\begin{array}{ccc}
-E_{L}^{\sigma }\widehat{\mathbb{F}}_{L}^{\sigma }f & \text{ if } & K\in 
\mathfrak{C}_{\mathcal{A}}\left( A\right) \\ 
E_{K}^{\sigma }\widehat{\mathbb{F}}_{K}^{\sigma }f-E_{L}^{\sigma }\widehat{%
\mathbb{F}}_{L}^{\sigma }f & \text{ if } & K\in \mathcal{C}_{A}%
\end{array}%
\right. ,  \label{telescoping}
\end{equation}%
where $\widehat{\mathbb{F}}_{K}^{\sigma }$ is defined in (\ref{F hat})
above.

Finally, in analogy with the broken differences $\bigtriangleup _{Q,{%
brok}}^{\mu ,\pi ,\mathbf{b}}$ and $\square _{Q,{brok}}^{\mu
,\pi ,\mathbf{b}}$ introduced above, we define%
\begin{equation}
\bigtriangleup _{I,{brok}}^{\mu ,\flat ,\mathbf{b}}f\equiv
\sum_{I^{\prime }\in \mathfrak{C}_{{brok}}\left( I\right) }\mathbb{%
E}_{I^{\prime }}^{\sigma ,\mathbf{b}}f\text{\ \ and \ \ }\square _{I,{brok}}^{\mu ,\flat ,\mathbf{b}}f\equiv \sum_{I^{\prime }\in \mathfrak{C}_{%
{brok}}\left( I\right) }\mathbb{F}_{I^{\prime }}^{\sigma ,\mathbf{b%
}}f\ ,  \label{def flat broken}
\end{equation}%
so that%
\begin{equation}
\bigtriangleup _{I}^{\mu ,\mathbf{b}}=\bigtriangleup _{I}^{\mu ,\flat ,%
\mathbf{b}}+\bigtriangleup _{I,{brok}}^{\mu ,\flat ,\mathbf{b}}%
\text{\ \ and\ \ }\square _{I}^{\mu ,\mathbf{b}}=\square _{I}^{\mu ,\flat ,\mathbf{%
b}}+\square _{I,{brok}}^{\mu ,\flat ,\mathbf{b}}\ .
\label{flat broken}
\end{equation}%
These modified differences and the identities (\ref{factor b_A}) and (\ref%
{telescoping}) play a useful role in the analysis of the nearby and
paraproduct forms.

\begin{lem}
\label{b proj}For dyadic cubes $R$ and $Q$ we have%
\begin{equation*}
\bigtriangleup _{R}^{\mu ,b}\bigtriangleup _{Q}^{\mu ,b}=\left\{ 
\begin{array}{ccc}
\bigtriangleup _{Q}^{\mu ,b} & \text{ if } & R=Q \\ 
0 & \text{ if } & R\not=Q%
\end{array}%
\right. .
\end{equation*}
\end{lem}

For the reader's convenience we now collect the various martingale and
probability estimates that will be used in the proof that follows. First we
summarize the martingale identities and estimates that we will use in our proof. Suppose $\mu $ is a positive locally finite Borel
measure, and that $\mathbf{b}$ is a $\infty$-weakly $\mu $-controlled accretive
family. Then,\\
\textbf{Martingale identities:} Both of the following identities hold pointwise 
$\mu $-almost everywhere, as well as in the sense of strong convergence in $%
L^{2}\left( \mu \right) $:%
\begin{eqnarray*}
f &=&\sum_{I\in \mathcal{D}:\ I\subset I_{\infty },\ \ell \left( I\right) \geq
2^{-N}}\!\!\!\!\!\!\!\!\square _{I}^{\sigma ,\mathbf{b}}f+\mathbb{F}_{I_{\infty }}^{\sigma ,%
\mathbf{b}}f, \\
f &=&\sum_{I\in \mathcal{D}:\ I\subset I_{\infty },\ \ell \left( I\right) \geq
2^{-N}}\!\!\!\!\!\!\!\!\bigtriangleup _{I}^{\sigma ,\mathbf{b}}f+\mathbb{E}_{I_{\infty }}^{\sigma ,%
\mathbf{b}}f.
\end{eqnarray*}
\textbf{Frame estimates:} Both of the following frame estimates hold:%
\begin{eqnarray}
\left\Vert f\right\Vert _{L^{2}\left( \mu \right) }^{2} 
&\approx &
\sum_{Q\in 
\mathcal{D}}\left\{ \left\Vert \square _{Q}^{\mu ,\mathbf{b}}f\right\Vert
_{L^{2}\left( \mu \right) }^{2}+\left\Vert \bigtriangledown _{Q}^{\mu ,\mathbf{b}%
}f\right\Vert _{L^{2}\left( \mu \right) }^{2}\right\}  \label{FRAME} \\
&\approx &
\sum_{Q\in \mathcal{D}}\left\{ \left\Vert \bigtriangleup _{Q}^{\mu
,\mathbf{b}}f\right\Vert _{L^{2}\left( \mu \right) }^{2}+\left\Vert
\bigtriangledown _{Q}^{\mu ,\mathbf{b}}f\right\Vert _{L^{2}\left( \mu
\right) }^{2}\right\} \ .  \notag
\end{eqnarray}
\textbf{Weak upper Riesz estimates:} Define the pseudoprojections, 
\begin{eqnarray}
\Psi _{\mathcal{B}}^{\mu ,\mathbf{b}}f &\equiv &\sum_{I\in \mathcal{B}%
}\square _{I}^{\mu ,\mathbf{b}}f, \label{Psi op}\\
\left( \Psi _{\mathcal{B}}^{\mu ,\mathbf{b}}\right) ^{\ast }f &\equiv
&\sum_{I\in \mathcal{B}}\left( \square _{I}^{\mu ,\mathbf{b}}\right) ^{\ast
}f=\sum_{I\in \mathcal{B}}\bigtriangleup _{I}^{\mu ,\mathbf{b}}f.\notag
\end{eqnarray}%
We have the `upper Riesz' inequalities for pseudoprojections $\Psi _{%
\mathcal{B}}^{\mu ,\mathbf{b}}$ and $\left( \Psi _{\mathcal{B}}^{\mu ,%
\mathbf{b}}\right) ^{\ast }$:%
\begin{eqnarray}
\left\Vert \Psi _{\mathcal{B}}^{\mu ,\mathbf{b}}f\right\Vert _{L^{2}\left(
\mu \right) }^{2} &\leq &C\sum_{I\in \mathcal{B}}\left\Vert \square
_{I}^{\mu ,\mathbf{b}}f\right\Vert _{L^{2}\left( \mu \right)
}^{2}+\sum_{I\in \mathcal{B}}\left\Vert \widehat{\bigtriangledown} _{I}^{\mu ,\mathbf{b}%
}f\right\Vert _{L^{2}\left( \mu \right) }^{2},  \label{UPPER RIESZ} \\
\left\Vert \left( \Psi _{\mathcal{B}}^{\mu ,\mathbf{b}}\right) ^{\ast
}f\right\Vert _{L^{2}\left( \mu \right) }^{2} &\leq &C\sum_{I\in \mathcal{B}%
}\left\Vert \bigtriangleup _{I}^{\mu ,\mathbf{b}}f\right\Vert _{L^{2}\left(
\mu \right) }^{2}+\sum_{I\in \mathcal{B}}\left\Vert \left( \widehat{\bigtriangledown} _{I}^{\mu ,%
\mathbf{b}}\right) ^{\ast }f\right\Vert _{L^{2}\left( \mu \right) }^{2}, 
\notag
\end{eqnarray}
for all $f\in L^{2}\left( \mu \right) $ and all subsets $\mathcal{B}$ of the
grid $\mathcal{D}$. Here the positive constant $C$ and depends
only on the accretivity constants, and is \emph{independent} of the subset $%
\mathcal{B}$ and the testing family $\mathbf{b}$. The Haar martingale
differences $\bigtriangleup _{Q}^{\mu ,\mathbf{b}}$ are independent of both
the testing families and the grid, while the Carleson averaging operators $%
\bigtriangledown _{Q}^{\mu }$  depend on the grid only through
the choice of broken children of $Q$.

\subsection{Monotonicity Lemma}
As in virtually all proofs of a two weight $T1$ theorem (see e.g. \cite{Lac}, 
\cite{LaSaShUr2} , \cite{SaShUr7} and/or \cite{SaShUr6}), the key to
starting an estimate for any of the forms we consider below, is the
Monotonicity Lemma and the Energy Lemma, to which we now turn. In dimension $%
n=1$ (\cite{LaSaShUr2}, \cite{Lac}) the Haar functions have opposite sign on
their children, and this was exploited in a simple but powerful monotonicity
argument. In higher dimensions, this simple argument no longer holds and
that Monotonicity Lemma is replaced with the Lacey-Wick formulation of the
Monotonicity Lemma (see \cite{LaWi}, and also \cite{SaShUr6}) involving the
smaller Poisson operator. As the martingale differences with test functions $%
b_{Q\,}$ here are no longer of one sign on children, we will adapt the
Lacey-Wick formulation of the Monotonicity Lemma to the operator $T^{\alpha
} $ and the dual martingale differences $\left\{ \square _{J}^{\omega ,%
\mathbf{b}^{\ast }}\right\} _{J\in \mathcal{G}}$, bearing in mind that the
operators $\square _{J}^{\omega ,\mathbf{b}^{\ast }}$ are no longer
projections, which results in only a one-sided estimate with additional
terms on the right hand side. It is here that we need the crucial property
that the Range of $\square _{J}^{\omega ,\mathbf{b}^{\ast }}$ is
orthogonal to constants, $\int \left( \square _{J}^{\omega ,\mathbf{b}^{\ast
}}\Psi \right) d\sigma =\int \left( \triangle _{J}^{\sigma ,\mathbf{b}^{\ast
}}1\right) \Psi d\omega =\int \left( 0\right) \Psi d\omega =0$.

We will also need the smaller Poisson integral used in the Lacey-Wick
formulation of the Monotonicity Lemma,
\begin{equation*}
\mathrm{P}_{1+\delta }^{\alpha }\left( J,\mu \right) \equiv \int \frac{%
\left\vert J\right\vert ^{\frac{1+\delta }{n}}}{\left( \left\vert
J\right\vert +\left\vert y-c_{J}\right\vert \right) ^{n+1+\delta -\alpha }}%
d\mu \left( y\right) ,
\end{equation*}%
which is discussed in more detail below.

\begin{lem}[Monotonicity Lemma]
\label{mono}Suppose that$\ I$ and $J$ are cubes in $\mathbb{R}^n$ such
that $J\subset \gamma J\subset I$ for some $\gamma >1$, and that $\mu $ is a
signed measure on $\mathbb{R}^n$ supported outside $I$. Let $0<\delta <1$ and
let $\Psi \in L^{2}\left( \omega \right) $. Finally suppose that $T^{\alpha
} $ is a standard fractional singular integral on $\mathbb{R}^n$ with $0\leq \alpha <1$, and
suppose that $\mathbf{b}^{\ast }$ is an $\infty $-weakly $\mu $-controlled
accretive family on $\mathbb{R}^n$. Then we have the estimate%
\begin{equation}
\left\vert \left\langle T^{\alpha }\mu ,\square _{J}^{\omega ,\mathbf{b}%
^{\ast }}\Psi \right\rangle _{\omega }\right\vert \lesssim C_{\mathbf{b}%
^{\ast }}C_{CZ}\ \Phi ^{\alpha }\left( J,\left\vert \mu \right\vert \right)
\ \left\Vert \square _{J}^{\omega ,\mathbf{b}^{\ast }}\Psi \right\Vert
_{L^{2}\left( \omega \right) }^{\bigstar },  \label{estimate}
\end{equation}%
where%
\begin{eqnarray*}
\Phi ^{\alpha }\left( J,\left\vert \mu \right\vert \right) &\equiv &\frac{%
\mathrm{P}^{\alpha }\left( J,\left\vert \mu \right\vert \right) }{\left\vert
J\right\vert }\left\Vert \bigtriangleup _{J}^{\omega ,\mathbf{b}^{\ast
}}x\right\Vert _{L^{2}\left( \omega \right) }^{\spadesuit }+\frac{\mathrm{P}%
_{1+\delta }^{\alpha }\left( J,\left\vert \mu \right\vert \right) }{%
\left\vert J\right\vert }\left\Vert x-m_{J}\right\Vert _{L^{2}\left( \mathbf{%
1}_{J}\omega \right) }, \\
\left\Vert \bigtriangleup _{J}^{\omega ,\mathbf{b}^{\ast }}x\right\Vert
_{L^{2}\left( \omega \right) }^{\spadesuit 2} &\equiv &\left\Vert
\bigtriangleup _{J}^{\omega ,\mathbf{b}^{\ast }}x\right\Vert _{L^{2}\left(
\omega \right) }^{2}+\inf_{z\in \mathbb{R}}\sum_{J^{\prime }\in \mathfrak{C}%
_{{brok}}\left( J\right) }\left\vert J^{\prime }\right\vert
_{\omega }\left( E_{J^{\prime }}^{\omega }\left\vert x-z\right\vert \right)
^{2}, \\
\left\Vert \square _{J}^{\omega ,\mathbf{b}^{\ast }}\Psi \right\Vert
_{L^{2}\left( \mu \right) }^{\bigstar 2} &\equiv &\left\Vert \square
_{J}^{\omega ,\mathbf{b}^{\ast }}\Psi \right\Vert _{L^{2}\left( \mu \right)
}^{2}+\sum_{J^{\prime }\in \mathfrak{C}_{{brok}}\left( J\right)
}\left\vert J^{\prime }\right\vert _{\omega }\left[ E_{J^{\prime }}^{\omega
}\left\vert \Psi \right\vert \right] ^{2}.
\end{eqnarray*}%
All of the implied constants above depend only on $\gamma >1$, $0<\delta <1$
and $0<\alpha <1$.
\end{lem}

Using $\displaystyle \bigtriangledown _{J}^{\omega }h=\sum_{J^{\prime }\in \mathfrak{C}_{{brok}}\left( J\right) }\left( E_{J^{\prime }}^{\omega
}\left\vert h\right\vert \right) \mathbf{1}_{J^{\prime }}$ defined in (\ref{Carleson avg op}), we can rewrite the expressions $%
\left\Vert \bigtriangleup _{J}^{\omega ,\mathbf{b}^{\ast }}x\right\Vert
_{L^{2}\left( \omega \right) }^{\spadesuit 2}$ and $\left\Vert \square
_{J}^{\omega ,\mathbf{b}^{\ast }}\Psi \right\Vert _{L^{2}\left( \mu \right)
}^{\bigstar 2}$ as%
\begin{eqnarray*}
\left\Vert \bigtriangleup _{J}^{\omega ,\mathbf{b}^{\ast }}x\right\Vert
_{L^{2}\left( \omega \right) }^{\spadesuit 2} &\equiv &\left\Vert
\bigtriangleup _{J}^{\omega ,\mathbf{b}^{\ast }}x\right\Vert _{L^{2}\left(
\omega \right) }^{2}+\inf_{z\in \mathbb{R}}\left\Vert \bigtriangledown
_{J}^{\omega }\left( x-z\right) \right\Vert _{L^{2}\left( \omega \right)
}^{2}, \\
\left\Vert \square _{J}^{\omega ,\mathbf{b}^{\ast }}\Psi \right\Vert
_{L^{2}\left( \mu \right) }^{\bigstar 2} &\equiv &\left\Vert \square
_{J}^{\omega ,\mathbf{b}^{\ast }}\Psi \right\Vert _{L^{2}\left( \mu \right)
}^{2}+\left\Vert \bigtriangledown _{J}^{\omega }\Psi \right\Vert
_{L^{2}\left( \omega \right) }^{2}.
\end{eqnarray*}

\begin{proof}

 Using $\square _{J}^{\omega ,%
\mathbf{b}^{\ast }}=\square _{J}^{\omega ,\pi ,\mathbf{b}^{\ast }}\square
_{J}^{\omega ,\pi ,\mathbf{b}^{\ast }}+\square _{J,{brok}}^{\omega
,\pi, \mathbf{b}^{\ast }}$, we write%
\begin{eqnarray*}
\left\vert \left\langle T^{\alpha }\mu ,\square _{J}^{\omega ,\mathbf{b}%
^{\ast }}\Psi \right\rangle _{\omega }\right\vert &=&\left\vert \left\langle
T^{\alpha }\mu ,\left( \square _{J}^{\omega ,\pi ,\mathbf{b}^{\ast }}\square
_{J}^{\omega ,\pi ,\mathbf{b}^{\ast }}+\square _{J,{brok}}^{\omega,\pi
,\mathbf{b}^{\ast }}\right) \Psi \right\rangle _{\omega }\right\vert \\
&\leq &
\left\vert \left\langle T^{\alpha }\mu ,\square _{J}^{\omega ,\pi ,%
\mathbf{b}^{\ast }}\square _{J}^{\omega ,\pi ,\mathbf{b}^{\ast }}\Psi
\right\rangle _{\omega }\right\vert +\left\vert \left\langle T^{\alpha }\mu
,\square _{J,{brok}}^{\omega,\pi,\mathbf{b}^{\ast }}\Psi
\right\rangle _{\omega }\right\vert \\ &\equiv& \mathrm{I+II}.
\end{eqnarray*}%
Since $\left\langle 1,\square _{J}^{\omega ,\pi ,\mathbf{b}^{\ast
}}h\right\rangle _{\omega }=0$, we use $\displaystyle m_{J}=\frac{1}{\left\vert
J\right\vert _{\omega }}\int_{J}xd\omega \left( x\right) $ to obtain
\begin{eqnarray*}
T^{\alpha }\mu \left( x\right) -T^{\alpha }\mu \left( m_{J}\right) 
&=&
\int \left[ \left( K^{\alpha }\right) \left( x,y\right) -\left( K^{\alpha
}\right) \left( m_{J},y\right) \right] d\mu \left( y\right) \\
&=&
\int \left[
\nabla( K^{\alpha })^T\left( \theta \left( x,m_{J}\right),y
\right)\boldsymbol{\cdot} \left( x-m_{J}\right) \right] d\mu \left( y\right)
\end{eqnarray*}
for some $\theta \left( x,m_{J}\right) \in J$ to obtain%
\begin{eqnarray*}
\mathrm{I} &=&\left\vert \int \left[ T^{\alpha }\mu \left( x\right) -T^{\alpha }\mu
\left( m_{J}\right) \right] \ \square _{J}^{\omega ,\pi ,\mathbf{b}^{\ast
}}\square _{J}^{\omega ,\pi ,\mathbf{b}^{\ast }}\Psi \left( x\right) d\omega
\left( x\right) \right\vert \\
&=&
\left\vert \int \left\{ \int  \nabla( K^{\alpha })^T\left( \theta \left( x,m_{J}\right) \right)  d\mu \left(
y\right) \right\} \boldsymbol{\cdot} \left( x-m_{J}\right) \ \square _{J}^{\omega ,\pi ,%
\mathbf{b}^{\ast }}\square _{J}^{\omega ,\pi ,\mathbf{b}^{\ast }}\Psi \left(
x\right) d\omega \left( x\right) \right\vert \\
&\leq &
\left\vert \int \left\{ \int  \nabla( K^{\alpha })^T\left( m_{J},y\right)  d\mu \left( y\right) \right\}\boldsymbol{\cdot} \ \left(
x-m_{J}\right) \ \square _{J}^{\omega ,\pi ,\mathbf{b}^{\ast }}\square
_{J}^{\omega ,\pi ,\mathbf{b}^{\ast }}\Psi \left( x\right) d\omega \left(
x\right) \right\vert \\
& & +
\Bigg| \int \!\!\left\{\! \int \!\left[ \nabla( K^{\alpha })^T\left( \theta \left( x,m_{J}\right),y \right) \!\!-\!\! \nabla( K^{\alpha })^T\left( m_{J}, y\right) \right] d\mu \left( y\right)\! \right\}\! \boldsymbol{\cdot} \!
\left( x-m_{J}\right)\  \square _{J}^{\omega ,\pi ,\mathbf{b}^{\ast
}}\square _{J}^{\omega ,\pi ,\mathbf{b}^{\ast }}\Psi \left( x\right) d\omega
\left( x\right)\! \Bigg| 
\\
&\equiv&
\mathrm{I}_1+\mathrm{I}_2
\end{eqnarray*}
Now we estimate%
\begin{eqnarray*}
\mathrm{I}_{1} &=&
\left\vert\left[ \int  \nabla( K^{\alpha })\left( m_{J},y\right) d\mu \left( y\right) \right]^T \ \boldsymbol{\cdot}
\int \left( x-m_{J}\right) \ \square _{J}^{\omega ,\pi ,\mathbf{b}^{\ast
}}\square _{J}^{\omega ,\pi ,\mathbf{b}^{\ast }}\Psi \left( x\right) d\omega
\left( x\right) \right\vert \\
&\leq&
n \int \int \left\vert  \nabla( K^{\alpha })\left(
m_{J},y\right)  \right\vert d|\mu| \left( y\right) \ \left\vert \bigtriangleup _{J}^{\omega ,\pi ,\mathbf{b}^{\ast }}x \right\vert \
\left\vert \square _{J}^{\omega ,\pi ,\mathbf{b}^{\ast }}\Psi \left( x\right)
\right\vert \ d\omega \left( x\right)  \\
&\lesssim &
n\cdot C_{CZ}\frac{\mathrm{P}^{\alpha }\left( J,\left\vert \mu
\right\vert \right) }{\left\vert J\right\vert^\frac{1}{n} }\ \left\Vert \bigtriangleup
_{J}^{\omega ,\pi ,\mathbf{b}^{\ast }}x\right\Vert _{L^{2}\left( \omega
\right) }\left\Vert \square _{J}^{\omega ,\pi ,\mathbf{b}^{\ast }}\Psi
\right\Vert _{L^{2}\left( \omega \right) }\ 
\end{eqnarray*}%
and%
\begin{eqnarray*}
\mathrm{I}_{2} 
&\lesssim &C_{CZ}\frac{\mathrm{P}_{1+\delta }^{\alpha }\left(
J,\left\vert \mu \right\vert \right) }{\left\vert J\right\vert^\frac{1}{n} }\int
\left\vert x-m_{J}\right\vert \left\vert \square _{J}^{\omega ,\pi ,\mathbf{b%
}^{\ast }}\square _{J}^{\omega ,\pi ,\mathbf{b}^{\ast }}\Psi \left( x\right)
\right\vert d\omega \left( x\right) \\
&\lesssim &C_{CZ}\frac{\mathrm{P}_{1+\delta }^{\alpha }\left( J,\left\vert
\mu \right\vert \right) }{\left\vert J\right\vert }\sqrt{\int_{J}\left\vert
x-m_{J}\right\vert ^{2}d\omega \left( x\right) }\left\Vert \square
_{J}^{\omega ,\pi ,\mathbf{b}^{\ast }}\square _{J}^{\omega ,\pi ,\mathbf{b}%
^{\ast }}\Psi \right\Vert _{L^{2}\left( \omega \right) } \\
&\lesssim &C_{CZ}\frac{\mathrm{P}_{1+\delta }^{\alpha }\left( J,\left\vert
\mu \right\vert \right) }{\left\vert J\right\vert }\left\Vert
x-m_{J}\right\Vert _{L^{2}\left( \mathbf{1}_{J}\omega \right) }\left\Vert
\square _{J}^{\omega ,\pi ,\mathbf{b}^{\ast }}\Psi \right\Vert _{L^{2}\left(
\omega \right) }\ .
\end{eqnarray*}%
For term $\mathrm{II}$ we fix $z\in \overline{J}$ for the moment. Then since 
$$
\left\langle 1,\square _{J,{brok}}^{\omega ,\mathbf{b}^{\ast
}}h\right\rangle _{\omega }=\left\langle 1,\square _{J}^{\omega ,\mathbf{b}%
^{\ast }}h-\square _{J}^{\omega ,\pi ,\mathbf{b}^{\ast }}h\right\rangle
_{\omega }=0
$$
we have

\begin{eqnarray*} 
\mathrm{II} 
&=& 
\left\vert \left\langle T^{\alpha }\mu ,\square _{J,{brok}%
}^{\omega ,\mathbf{b}^{\ast }}\Psi \right\rangle _{\omega }\right\vert   \\
&=&
\left\vert \int \left\{ \int  \nabla( K^{\alpha })^T\left( \theta \left( x,z\right),y \right)  d\mu \left( y\right)
\right\} \boldsymbol{\cdot}\left( x-z\right) \ \square _{J,{brok}}^{\omega ,\pi,\mathbf{b}^{\ast }}\Psi  \left( x\right) d\omega \left( x\right) \right\vert \\
&\leq&
C_{CZ}\frac{\mathrm{P}^{\alpha }\left( J,\left\vert \mu\right\vert \right) }{\left\vert J\right\vert^\frac{1}{n} }  \int |x-z|\cdot\Big|\square _{J,{brok}}^{\omega ,\pi,\mathbf{b}^{\ast }}\Psi  \left( x\right)\Big|d\omega(x)\\
&\leq&
C_{CZ}\frac{\mathrm{P}^{\alpha }\left( J,\left\vert \mu\right\vert \right) }{\left\vert J\right\vert^\frac{1}{n} }\sum_{J'\in\mathfrak{C}_{{brok}(J)}}\int_{J'} |x-z|\cdot\mathbf{1}_{J^{\prime }}E_{J^{\prime }}^{\omega }\left\vert \Psi \right\vert d\omega(x)
\end{eqnarray*}
having used the reverse H\"{o}lder control of children (\ref{rev Hol con}) to obtain
\begin{equation*}
\left\vert \square _{J,{brok}}^{\omega ,\mathbf{b}^{\ast }}\Psi
\right\vert 
=
\left\vert \sum_{J^{\prime }\in \mathfrak{C}_{{brok}%
}\left( JQ\right) }\left( \mathbb{F}_{J^{\prime }}^{\omega ,b_{J^{\prime }}}-%
\mathbb{F}_{J^{\prime }}^{\omega ,b_{J}}\right) \Psi \right\vert 
\lesssim
\sum_{J^{\prime }\in \mathfrak{C}_{{brok}}\left( J\right) }\mathbf{%
1}_{J^{\prime }}E_{J^{\prime }}^{\omega }\left\vert \Psi \right\vert ,
\end{equation*}%
and since 
$$
\int_{J'} |x-z|\cdot\mathbf{1}_{J^{\prime }}E_{J^{\prime }}^{\omega }\left\vert \Psi \right\vert d\omega(x) 
=
\int_{J'} \frac{|x-z|}{\sqrt{|J'|_\omega}} \frac{\mathbf{1}_{J'}\int_{J'}|\Psi|d\omega(x)}{\sqrt{|J'|_\omega}}d\omega(x)
$$
we get
\begin{equation*}
\mathrm{II} \leq 
C_{CZ}\frac{\mathrm{P}^{\alpha }\left( J,\left\vert \mu
\right\vert \right) }{\left\vert J\right\vert^\frac{1}{n} }\sqrt{\sum_{J^{\prime }\in 
\mathfrak{C}_{{brok}}\left( J\right) }\left\vert J^{\prime
}\right\vert _{\omega }\left( E_{J^{\prime }}^{\omega }\left\vert
x-z\right\vert \right) ^{2}}\sqrt{\sum_{J^{\prime }\in \mathfrak{C}_{%
{brok}}\left( J\right) }\left\vert J^{\prime }\right\vert _{\omega
}\left[ E_{J^{\prime }}^{\omega }\left\vert \Psi \right\vert \right] ^{2}}.
\end{equation*}

Combining the estimates for terms $\mathrm{I}$ and $\mathrm{II}$, we obtain%
\begin{eqnarray*}
&&
\left\vert \left\langle T^{\alpha }\mu ,\square _{J}^{\omega ,\mathbf{b}%
^{\ast }}\Psi \right\rangle _{\omega }\right\vert \\
&\lesssim &
C_{CZ}\frac{\mathrm{P}^{\alpha }\left( J,\left\vert \mu
\right\vert \right) }{\left\vert J\right\vert^\frac{1}{n} }\ \left\Vert \bigtriangleup
_{J}^{\omega ,\pi ,\mathbf{b}^{\ast }}x\right\Vert _{L^{2}\left( \omega
\right) }\left\Vert \square _{J}^{\omega ,\pi ,\mathbf{b}^{\ast }}\Psi
\right\Vert _{L^{2}\left( \omega \right) } \\
&+&
C_{CZ}\frac{\mathrm{P}_{1+\delta }^{\alpha }\left( J,\left\vert \mu
\right\vert \right) }{\left\vert J\right\vert^\frac{1}{n} }\left\Vert x-m_{J}\right\Vert
_{L^{2}\left( \mathbf{1}_{J}\omega \right) }\left\Vert \square _{J}^{\omega
,\pi ,\mathbf{b}^{\ast }}\Psi \right\Vert _{L^{2}\left( \omega \right) } \\
&+&C_{CZ}\frac{\mathrm{P}^{\alpha }\left( J,\left\vert \mu \right\vert
\right) }{\left\vert J\right\vert^\frac{1}{n} }\inf_{z\in \overline{J}}\sqrt{%
\sum_{J^{\prime }\in \mathfrak{C}_{{brok}}\left( J\right)
}\left\vert J^{\prime }\right\vert _{\omega }\left( E_{J^{\prime }}^{\omega
}\left\vert x-z_{2}\right\vert \right) ^{2}} 
\sqrt{\sum_{J^{\prime }\in 
\mathfrak{C}_{{brok}}\left( J\right) }\left\vert J^{\prime
}\right\vert _{\omega }\left[ E_{J^{\prime }}^{\omega }\left\vert \Psi
\right\vert \right] ^{2}}
\end{eqnarray*}%
and then noting that the infimum over $z\in \mathbb{R}$ is achieved for $%
z\in \overline{J}$, and using the triangle inequality on $\square
_{J}^{\omega ,\pi ,\mathbf{b}^{\ast }}=\square _{J}^{\omega ,\mathbf{b}%
^{\ast }}-\square _{J,{brok}}^{\omega,\pi ,\mathbf{b}^{\ast }}$ we get
(\ref{estimate}).
\end{proof}

The right hand side of (\ref{estimate}) in the Monotonicity Lemma will be
typically estimated in what follows using the frame inequalities for any cube $K$,%
\begin{eqnarray*}
\sum_{J\subset K}\left\Vert \square _{J}^{\omega ,\mathbf{b}^{\ast }}\Psi
\right\Vert _{L^{2}\left( \omega \right) }^{\bigstar 2} &\lesssim
&\left\Vert \Psi \right\Vert _{L^{2}\left( \omega \right) }^{2}\ , \\
\sum_{J\subset K}\left\Vert \bigtriangleup _{J}^{\omega ,\mathbf{b}^{\ast
}}x\right\Vert _{L^{2}\left( \omega \right) }^{\spadesuit 2} 
&\lesssim&
\int_{K}\left\vert x-m_{K}\right\vert ^{2}d\omega \left( x\right) \ ,
\end{eqnarray*}%
together with these inequalities for the square function expressions.
To see the last one, write $x=(x_1,\dots, x_n)$ and note that for $J\subset K$,
\begin{eqnarray*}
\left\Vert \bigtriangleup _{J}^{\omega ,\mathbf{b}^{\ast
}}x\right\Vert _{L^{2}\left( \omega \right) }^{\spadesuit 2}\!=\!\int_J \left\vert \bigtriangleup _{J}^{\omega ,\mathbf{b}^{\ast
}}x\right\vert^{2}d\omega
\!=\!
\int_J \sum_{i=1}^n\left\vert \bigtriangleup _{J}^{\omega ,\mathbf{b}^{\ast
}}x_i\right\vert^{2}d\omega
\!\!\!\!&\leq& \!\!\!\!
\sum_{i=1}^n\int_{K}\left\vert x_i-{m_{K}}_i\right\vert ^{2}d\omega
\!=\! 
||x-m_k||_{L^2(\textbf{1}_K\omega)}^2
\end{eqnarray*}
using the one-variable result from \cite{SaShUr12}.

\begin{lem}
For any cube $K$ we have%
\begin{eqnarray}
\sum_{J\subset K}\sum_{J^{\prime }\in \mathfrak{C}_{{brok}}\left(
J\right) }\left\vert J^{\prime }\right\vert _{\omega }\left[ E_{J^{\prime
}}^{\omega }\left\vert \Psi \right\vert \left( x\right) \right] ^{2}
&\lesssim &\int_{K}\left\vert \Psi \left( x\right) \right\vert ^{2}d\omega
\left( x\right) ,  \label{with both} \\
\text{and }\sum_{J\subset K}\inf_{z\in \mathbb{R}}\sum_{J^{\prime }\in 
\mathfrak{C}_{{brok}}\left( J\right) }\left\vert J^{\prime
}\right\vert _{\omega }\left( E_{J^{\prime }}^{\omega }\left\vert
x-z\right\vert \right) ^{2} &\lesssim &\int_{K}\left\vert x-m_{K}\right\vert
^{2}d\omega \left( x\right) .  \notag
\end{eqnarray}
\end{lem}

\begin{proof}
The first inequality in (\ref{with both}) is just the Carleson embedding
theorem since the cubes\\ $\left\{ J^{\prime }\in \mathfrak{C}_{{%
brok}}\left( J\right) :J\subset K\right\} $ satisfy an $\omega $-Carleson
condition, and the second inequality in (\ref{with both}) follows by
choosing $z=m_{K}$ to obtain%
\begin{equation*}
\inf_{z\in \mathbb{R}}\sum_{J^{\prime }\in \mathfrak{C}_{{brok}%
}\left( J\right) }\left\vert J^{\prime }\right\vert _{\omega }\left(
E_{J^{\prime }}^{\omega }\left\vert x-z\right\vert \right) ^{2}\leq
\sum_{J^{\prime }\in \mathfrak{C}_{{brok}}\left( J\right)
}\left\vert J^{\prime }\right\vert _{\omega }\left( E_{J^{\prime }}^{\omega
}\left\vert x-m_{K}\right\vert \right) ^{2},
\end{equation*}%
and then applying the Carleson embedding theorem again:%
\begin{equation*}
\sum_{J\subset K}\sum_{J^{\prime }\in \mathfrak{C}_{{brok}}\left(
J\right) }\left\vert J^{\prime }\right\vert _{\omega }\left( E_{J^{\prime
}}^{\omega }\left\vert x-m_{K}\right\vert \right) ^{2}\lesssim
\int_{K}\left\vert x-m_{K}\right\vert ^{2}d\omega \left( x\right) .
\end{equation*}
\end{proof}

\subsubsection{The smaller Poisson integral}

The expressions $$\inf_{z\in \mathbb{R}}\frac{\mathrm{P}_{1+\delta }^{\alpha
}\left( J,\left\vert \mu \right\vert \right) }{\left\vert J\right\vert }%
\left\Vert x-z\right\Vert _{L^{2}\left( \mathbf{1}_{J}\omega \right)
}\left\Vert \square _{J}^{\omega ,\mathbf{b}^{\ast }}\Psi \right\Vert
_{L^{2}\left( \omega \right) }^{\bigstar }$$ are typically easier to sum due
to the small Poisson operator $\mathrm{P}_{1+\delta }^{\alpha }\left(
J,\left\vert \mu \right\vert \right) $. To illlustrate, we show here one way
in which we can exploit the additional decay in the Poisson integral $%
\mathrm{P}_{1+\delta }^{\alpha }$. Suppose that $J$ is good in $I$ with $%
\ell \left( J\right) =2^{-s}\ell \left( I\right) $ (see Definition \ref{good
arb} below for `goodness'). We then compute%
\begin{eqnarray*}
\frac{\mathrm{P}_{1+\delta }^{\alpha }\left( J,\mathbf{1}_{A\backslash
I}\sigma \right) }{\left\vert J\right\vert ^{\frac{1}{n}}}
&\approx&
\int_{A\backslash I}\frac{\left\vert J\right\vert ^{\frac{\delta }{n}}}{%
\left\vert y-c_{J}\right\vert ^{n+1+\delta -\alpha }}d\sigma \left( y\right)
\\
&\leq &
\int_{A\backslash I}\left( \frac{\left\vert J\right\vert ^{\frac{1}{n}}%
}{{\dist}\left( c_{J},I^{c}\right) }\right) ^{\delta }\frac{1}{%
\left\vert y-c_{J}\right\vert ^{n+1-\alpha }}d\sigma \left( y\right) \\
&\lesssim &
\left( \frac{\left\vert J\right\vert ^{\frac{1}{n}}}{{%
\dist}\left( c_{J},I^{c}\right) }\right) ^{\delta }\frac{\mathrm{P}^{\alpha
}\left( J,\mathbf{1}_{A\backslash I}\sigma \right) }{\left\vert J\right\vert
^{\frac{1}{n}}},
\end{eqnarray*}%
and use the goodness inequality,%
\begin{equation*}
\dist\left( c_{J},I^{c}\right) \geq 2\ell \left( I\right)
^{1-\varepsilon }\ell \left( J\right) ^{\varepsilon }\geq 2\cdot
2^{s\left( 1-\varepsilon \right) }\ell \left( J\right) ,
\end{equation*}%
to conclude that%
\begin{equation}
\left( \frac{\mathrm{P}_{1+\delta }^{\alpha }\left( J,\mathbf{1}_{A\backslash
I}\sigma \right) }{\left\vert J\right\vert ^{\frac{1}{n}}}\right) \lesssim
2^{-s\delta \left( 1-\varepsilon \right) }\frac{\mathrm{P}^{\alpha }\left( J,%
\mathbf{1}_{A\backslash I}\sigma \right) }{\left\vert J\right\vert ^{\frac{1}{%
n}}}  
\end{equation}%
Now we can estimate%
\begin{eqnarray*}
&&
\sum_{J\subset K:\ J\text{ good}\text{ in }K}\inf_{z\in 
\mathbb{R}}\frac{\mathrm{P}_{1+\delta }^{\alpha }\left( J,\mathbf{1}%
_{K^{c}}\left\vert \mu \right\vert \right) }{\left\vert J\right\vert^\frac{1}{n} }
\left\Vert x-z\right\Vert _{L^{2}\left( \mathbf{1}_{J}\omega \right)
}\left\Vert \square _{J}^{\omega ,\mathbf{b}^{\ast }}\Psi \right\Vert
_{L^{2}\left( \omega \right) }^{\bigstar } \\
&\leq &\!\!\!\!\!
\sqrt{\sum_{\substack{J\subset K\\ J\text{ good}\text{ in }K}} \!\!\!\!\!\!\!\left( \frac{\mathrm{P}_{1+\delta }^{\alpha }\left( J,\mathbf{1}
_{K^{c}}\left\vert \mu \right\vert \right) }{\left\vert J\right\vert^\frac{1}{n} }%
\right)^{\!\!2}\inf_{z\in \mathbb{R}}\left\Vert x-z\right\Vert _{L^{2}\left( 
\mathbf{1}_{J}\omega \right) }^{2}} \sqrt{\sum_{\substack{J\subset K\\ J\text{ good}\text{ in }K}}\!\!\!\!\!\!\left\Vert \square _{J}^{\omega ,\mathbf{b}^{\ast
}}\Psi \right\Vert _{L^{2}\left( \omega \right) }^{\bigstar 2}}
\end{eqnarray*}%
where
\begin{eqnarray*}
&&
\sum_{J\subset K:\ J\text{ good}\text{ in }K}\left( \frac{%
\mathrm{P}_{1+\delta }^{\alpha }\left( J,\mathbf{1}_{K^{c}}\left\vert \mu
\right\vert \right) }{\left\vert J\right\vert }\right) ^{2}\inf_{z\in 
\mathbb{R}}\left\Vert x-z\right\Vert _{L^{2}\left( \mathbf{1}_{J}\omega
\right) }^{2} \\
&=&
\sum_{s=0}^{\infty }\sum_{\substack{ J\subset K:\ J\text{ good}\text{ in }K  \\ \ell \left( J\right) =2^{-s}\ell \left( I\right) }}%
\left( \frac{\mathrm{P}_{1+\delta }^{\alpha }\left( J,\mathbf{1}%
_{K^{c}}\left\vert \mu \right\vert \right) }{\left\vert J\right\vert }%
\right) ^{2}\inf_{z\in \mathbb{R}}\left\Vert x-z\right\Vert _{L^{2}\left( 
\mathbf{1}_{J}\omega \right) }^{2} \\
&\leq &
\sum_{s=0}^{\infty }\sum_{\substack{ J\subset K:\ J\text{ good}\text{ in }K  \\ \ell \left( J\right) =2^{-s}\ell \left(
I\right) }}\left( 2^{-s\delta \left( 1-\varepsilon \right) }\frac{\mathrm{P}%
^{\alpha }\left( J,\mathbf{1}_{K^{c}}\sigma \right) }{\left\vert
J\right\vert ^{\frac{1}{n}}}\right) ^{2}\inf_{z\in \mathbb{R}}\left\Vert
x-z\right\Vert _{L^{2}\left( \mathbf{1}_{J}\omega \right) }^{2} \\
&\leq &
\left( \frac{\mathrm{P}^{\alpha }\left( K,\mathbf{1}_{K^{c}}\sigma
\right) }{\left\vert K\right\vert ^{\frac{1}{n}}}\right)
^{2}\sum_{s=0}^{\infty }\sum_{\substack{ J\subset K:\ J\text{ good}\text{ in }K  \\ \ell \left( J\right) =2^{-s}\ell \left( I\right) }}%
2^{-2s\delta \left( 1-\varepsilon \right) }\inf_{z\in \mathbb{R}}\left\Vert
x-z\right\Vert _{L^{2}\left( \mathbf{1}_{K}\omega \right) }^{2}\\
&\lesssim &
\left( \frac{\mathrm{P}^{\alpha }\left( K,\mathbf{1}%
_{K^{c}}\sigma \right) }{\left\vert K\right\vert ^{\frac{1}{n}}}\right)
^{2}\inf_{z\in \mathbb{R}}\left\Vert x-z\right\Vert _{L^{2}\left( \mathbf{1}%
_{K}\omega \right) }^{2}\ ,
\end{eqnarray*}%
and where we have used (\ref{Poisson inequality}), which gives in particular%
\begin{equation*}
\mathrm{P}^{\alpha }(J,\mu \mathbf{1}_{I^{c}})\lesssim \left( \frac{\ell
\left( J\right) }{\ell \left( I\right) }\right) ^{1-\varepsilon \left(
n+1-\alpha \right) }\mathrm{P}^{\alpha }(I,\mu \mathbf{1}_{I^{c}}).
\end{equation*}%
for $J\subset I$ and $d\left( J,\partial I\right) >2\ell \left(
J\right) ^{\varepsilon }\ell \left( I\right) ^{1-\varepsilon }$. We will use
such arguments repeatedly in the sequel.

Armed with the Monotonicity Lemma and the lower frame inequality%
\begin{equation*}
\sum_{I\in \mathcal{D}}\left\Vert \square _{I}^{\omega ,\mathbf{b}^{\ast
}}g\right\Vert _{L^{2}\left( \mu \right) }^{\bigstar 2}\lesssim \left\Vert
g\right\Vert _{L^{2}\left( \omega \right) }^{2}\ ,
\end{equation*}%
we can obtain a $\mathbf{b}^{\ast }$-analogue of the Energy Lemma as in \cite%
{SaShUr7} and/or \cite{SaShUr6}.

\subsubsection{The Energy Lemma}

Suppose now we are given a subset $\mathcal{H}$ of the dyadic grid $\mathcal{%
G}$. \label{nonstandard norm}Due to the failure of both martingale and dual
martingale pseudoprojections $\mathsf{Q}_{\mathcal{H}}^{\omega ,\mathbf{b}%
^{\ast }}x$ and $\mathsf{P}_{\mathcal{H}}^{\omega ,\mathbf{b}^{\ast }}g$ (see below for definition) to satisfy inequalities of the
form $\left\Vert \mathsf{P}_{\mathcal{H}}^{\omega ,\mathbf{b}^{\ast
}}g\right\Vert _{L^{2}\left( \omega \right) }\lesssim \left\Vert
g\right\Vert _{L^{2}\left( \omega \right) }$ when the children `break', it
is convenient to define the `square function norms' $\left\Vert \mathsf{Q}_{%
\mathcal{H}}^{\omega ,\mathbf{b}^{\ast }}x\right\Vert _{L^{2}\left( \omega
\right) }^{\spadesuit }$ and $\left\Vert \mathsf{P}_{\mathcal{H}}^{\omega ,%
\mathbf{b}^{\ast }}g\right\Vert _{L^{2}\left( \omega \right) }^{\bigstar }$
of the pseudoprojections 
\begin{equation*}
\mathsf{Q}_{\mathcal{H}}^{\omega ,\mathbf{b}^{\ast }}x=\sum_{J\in \mathcal{H}%
}\bigtriangleup _{J}^{\omega ,\mathbf{b}^{\ast }}x\text{ and }\mathsf{P}_{%
\mathcal{H}}^{\omega ,\mathbf{b}^{\ast }}g=\sum_{J\in \mathcal{H}}\square
_{J}^{\omega ,\mathbf{b}^{\ast }}g\ ,
\end{equation*}%
by%
\begin{eqnarray*}
\left\Vert \mathsf{Q}_{\mathcal{H}}^{\omega ,\mathbf{b}^{\ast }}x\right\Vert
_{L^{2}\left( \omega \right) }^{\spadesuit 2} 
&\equiv &
\sum_{J\in \mathcal{H}%
}\left\Vert \bigtriangleup _{J}^{\omega ,\mathbf{b}^{\ast }}x\right\Vert
_{L^{2}\left( \omega \right) }^{\spadesuit 2}\\ 
&=&
\sum_{J\in \mathcal{H}%
}\left\Vert \bigtriangleup _{J}^{\omega ,\mathbf{b}^{\ast }}x\right\Vert
_{L^{2}\left( \omega \right) }^{2}+\sum_{J\in \mathcal{H}}\inf_{z\in \mathbb{%
R}}\sum_{J^{\prime }\in \mathfrak{C}_{{brok}}\left( J\right)
}\left\vert J^{\prime }\right\vert _{\omega }\left( E_{J^{\prime }}^{\omega
}\left\vert x-z\right\vert \right) ^{2}\\
\left\Vert \mathsf{P}_{\mathcal{H}}^{\omega ,\mathbf{b}^{\ast }}g\right\Vert
_{L^{2}\left( \omega \right) }^{\bigstar 2} 
&\equiv &
\sum_{J\in \mathcal{H}%
}\left\Vert \square _{J}^{\omega ,\mathbf{b}^{\ast }}g\right\Vert
_{L^{2}\left( \omega \right) }^{\bigstar 2}\\
&=&\sum_{J\in \mathcal{H}}\left\Vert \square _{J}^{\omega ,\mathbf{b}^{\ast }}g\right\Vert
_{L^{2}\left( \omega \right) }^{2}+\sum_{J\in \mathcal{H}}\sum_{J^{\prime
}\in \mathfrak{C}_{{brok}}\left( J\right) }\left\vert J^{\prime
}\right\vert _{\omega }\left[ E_{J^{\prime }}^{\omega }\left\vert
g\right\vert  \right] ^{2}
\end{eqnarray*}%
for any subset $\mathcal{H}\subset \mathcal{G}$. The average $E_{J}^{\omega
}\left\vert x-z\right\vert $ above is taken with respect to the variable $x$%
, i.e. $E_{J}^{\omega }\left\vert x-z\right\vert =\frac{1}{\left\vert
J\right\vert _{\omega }}\int \left\vert x-z\right\vert d\omega \left(
x\right) $, and it is important that the infimum $\inf_{z\in \mathbb{R}}$ is
taken \emph{inside} the sum $\sum_{J\in \mathcal{H}}$.

Note that we are defining here square function expressions related to
pseudoprojections, which depend not only on the functions $\mathsf{Q}_{%
\mathcal{H}}^{\omega ,\mathbf{b}^{\ast }}x$ and $\mathsf{P}_{\mathcal{H}%
}^{\omega ,\mathbf{b}^{\ast }}g$, but also on the particular representations 
$\sum_{J\in \mathcal{H}}\bigtriangleup _{J}^{\omega ,\mathbf{b}^{\ast }}x$
and $\sum_{J\in \mathcal{H}}\square _{J}^{\omega ,\mathbf{b}^{\ast }}g$.
This slight abuse of notation should not cause confusion, and it provides a
useful way of bookkeeping the sums of squares of norms of martingale and
dual martingale differences $\left\Vert \bigtriangleup _{J}^{\omega ,\mathbf{%
b}^{\ast }}x\right\Vert _{L^{2}\left( \omega \right) }^{2}$ and $\left\Vert
\square _{J}^{\omega ,\mathbf{b}^{\ast }}g\right\Vert _{L^{2}\left( \omega
\right) }^{2}$, along with the norms of the associated Carleson square
function expressions 
\begin{eqnarray*}
\sum_{J\in \mathcal{H}}\inf_{z\in \mathbb{R}}\left\Vert \nabla _{J}^{\omega
}\left( x-z\right) \right\Vert _{L^{2}\left( \omega \right) }^{2}
&=&\sum_{J\in \mathcal{H}}\inf_{z\in \mathbb{R}}\sum_{J^{\prime }\in 
\mathfrak{C}_{{brok}}\left( J\right) }\left\vert J^{\prime
}\right\vert _{\omega }\left( E_{J^{\prime }}^{\omega }\left\vert
x-z\right\vert \right) ^{2} \\
\sum_{J\in \mathcal{H}}\left\Vert \nabla _{J}^{\omega }\Psi \right\Vert
_{L^{2}\left( \omega \right) }^{2} &=&\sum_{J\in \mathcal{H}}\sum_{J^{\prime
}\in \mathfrak{C}_{{brok}}\left( J\right) }\left\vert J^{\prime
}\right\vert _{\omega }\left[ E_{J^{\prime }}^{\omega }\left\vert \Psi
\right\vert \right] ^{2}.
\end{eqnarray*}%
Note also that the upper weak Riesz inequalities  yield
the inequalities%
\begin{eqnarray*}
\left\Vert \mathsf{Q}_{\mathcal{H}}^{\omega ,\mathbf{b}^{\ast }}x\right\Vert
_{L^{2}\left( \omega \right) }^{2}
&\lesssim &
\sum_{J\in \mathcal{H}%
}\left\Vert \bigtriangleup _{J}^{\omega ,\mathbf{b}^{\ast }}x\right\Vert
_{L^{2}\left( \omega \right) }^{2}\leq \left\Vert \mathsf{Q}_{\mathcal{H}%
}^{\omega ,\mathbf{b}^{\ast }}x\right\Vert _{L^{2}\left( \omega \right)}^{\spadesuit 2} \\
\left\Vert \mathsf{P}_{\mathcal{H}}^{\omega ,\mathbf{b}^{\ast }}g\right\Vert
_{L^{2}\left( \omega \right) }^{2}
&\lesssim &
\sum_{J\in \mathcal{H}%
}\left\Vert \square _{J}^{\omega ,\mathbf{b}^{\ast }}g\right\Vert
_{L^{2}\left( \omega \right) }^{2}\leq \left\Vert \mathsf{P}_{\mathcal{H}%
}^{\omega ,\mathbf{b}^{\ast }}g\right\Vert _{L^{2}\left( \omega \right)
}^{\bigstar 2}
\end{eqnarray*}%
We will exclusively use $\left\Vert \mathsf{Q}_{\mathcal{H}}^{\omega ,%
\mathbf{b}^{\ast }}x\right\Vert _{L^{2}\left( \omega \right) }^{\spadesuit
2} $ in connection with energy terms, and use $\left\Vert \mathsf{P}_{%
\mathcal{H}}^{\sigma ,\mathbf{b}^{\ast }}f\right\Vert _{L^{2}\left( \sigma
\right) }^{\bigstar 2}$ and $\left\Vert \mathsf{P}_{\mathcal{H}}^{\omega ,%
\mathbf{b}^{\ast }}g\right\Vert _{L^{2}\left( \omega \right) }^{\bigstar 2}$%
\ in connection with functions $f\in L^{2}\left( \sigma \right) $ and $g\in
L^{2}\left( \omega \right) $. Finally, note that $\mathsf{Q}_{\mathcal{H}%
}^{\omega ,\mathbf{b}^{\ast }}x=\mathsf{Q}_{\mathcal{H}}^{\omega ,\mathbf{b}%
^{\ast }}\left( x-m\right) $ for any constant $m$. 

Recall that%
\begin{equation*}
\Phi ^{\alpha }\left( J,\nu \right) \equiv \frac{\mathrm{P}^{\alpha }\left(
J,\nu \right) }{\left\vert J\right\vert^\frac{1}{n} }\left\Vert \bigtriangleup
_{J}^{\omega ,\mathbf{b}^{\ast }}x\right\Vert _{L^{2}\left( \omega \right)
}^{\spadesuit }+\frac{\mathrm{P}_{1+\delta }^{\alpha }\left( J,\nu \right) }{%
\left\vert J\right\vert^\frac{1}{n} }\left\Vert x-m_{J}\right\Vert _{L^{2}\left( \mathbf{%
1}_{J}\omega \right) }\ .
\end{equation*}

\begin{lem}[\textbf{Energy Lemma}]
\label{ener}Let $J\ $be a cube in $\mathcal{G}$. Let $\Psi _{J}$ be an $%
L^{2}\left( \omega \right) $ function supported in $J$ with vanishing $%
\omega $-mean, and let $\mathcal{H}\subset \mathcal{G}$ be such that $%
J^{\prime }\subset J$ for every $J^{\prime }\in \mathcal{H}$. Let $\nu $ be
a positive measure supported in $\mathbb{R}\backslash \gamma J$ with $\gamma
>1$, and for each $J^{\prime }\in \mathcal{H}$, let $d\nu _{J^{\prime
}}=\varphi _{J^{\prime }}d\nu $ with $\left\vert \varphi _{J^{\prime
}}\right\vert \leq 1$. Suppose that $\mathbf{b}^{\ast }$ is an $\infty $%
-weakly $\mu $-controlled accretive family on $\mathbb{R}^n$. Let $T^{\alpha }$
be a standard $\alpha $-fractional singular integral operator with $0\leq
\alpha <1$. Then we have%
\begin{eqnarray*}
&&
\left\vert \sum_{J^{\prime }\in \mathcal{H}}\left\langle T^{\alpha }\left(
\nu _{J^{\prime }}\right) ,\square _{J^{\prime }}^{\omega ,\mathbf{b}^{\ast
}}\Psi _{J}\right\rangle _{\omega }\right\vert 
\lesssim
C_{\gamma}\sum_{J^{\prime }\in \mathcal{H}}\Phi ^{\alpha }\left( J^{\prime },\nu\right) \left\Vert \square _{J^{\prime }}^{\omega ,\mathbf{b}^{\ast }}\Psi
_{J}\right\Vert _{L^{2}\left( \mu \right) }^{\bigstar } \\
&\lesssim &
C_{\gamma }\sqrt{\sum_{J^{\prime }\in \mathcal{H}}\Phi ^{\alpha
}\left( J^{\prime },\nu \right) ^{2}}\sqrt{\sum_{J^{\prime }\in \mathcal{H}%
}\left\Vert \square _{J^{\prime }}^{\omega ,\mathbf{b}^{\ast }}\Psi
_{J}\right\Vert _{L^{2}\left( \mu \right) }^{\bigstar 2}} \\
&\lesssim &
\left( \frac{\mathrm{P}^{\alpha }\left( J,\nu \right) }{%
\left\vert J\right\vert }\left\Vert \mathsf{Q}_{\mathcal{H}}^{\omega ,%
\mathbf{b}^{\ast }}x\right\Vert _{L^{2}\left( \omega \right) }^{\spadesuit }+%
\frac{\mathrm{P}_{1+\delta }^{\alpha }\left( J,\nu \right) }{\left\vert J\right\vert^\frac{1}{n} }\left\Vert x-m_{J}\right\Vert _{L^{2}\left( \mathbf{1}%
_{J}\omega \right) }\right) \left\Vert \mathsf{P}_{\mathcal{H}}^{\omega ,%
\mathbf{b}^{\ast }}\Psi _{J}\right\Vert _{L^{2}\left( \mu \right)
}^{\bigstar }
\end{eqnarray*}%
and in particular the `energy' estimate%
\begin{eqnarray*}
&&\left\vert \left\langle T^{\alpha }\varphi \nu ,\Psi _{J}\right\rangle
_{\omega }\right\vert \\
&\!\!\!\leq&
\!\!\!\!\!C_{\gamma }\!\!\left(\! \frac{\mathrm{P}^{\alpha
}\left( J,\nu \right) }{\left\vert J\right\vert ^\frac{1}{n}}\left\Vert \mathsf{Q}%
_{J}^{\omega ,\mathbf{b}^{\ast }}x\right\Vert _{L^{2}\left( \omega \right)
}^{\spadesuit }\!\!\!+\!\frac{\mathrm{P}_{1+\delta }^{\alpha }\left( J,\nu \right) }{%
\left\vert J\right\vert ^\frac{1}{n}}\left\Vert x-m_{J}\right\Vert _{L^{2}\left( \mathbf{%
1}_{J}\omega \right) }\!\!\right)
\!\!\left\Vert \sum_{J^{\prime }\subset J}\!\!\square
_{J^{\prime }}^{\omega ,\mathbf{b}^{\ast }}\Psi _{J}\right\Vert
_{L^{2}\left( \mu\right) }^{\bigstar }\
\end{eqnarray*}%

where $\displaystyle \left\Vert \sum_{J^{\prime }\subset J}\square _{J^{\prime }}^{\omega ,%
\mathbf{b}^{\ast }}\Psi _{J}\right\Vert _{L^{2}\left( \mu \right)
}^{\bigstar }\lesssim \left\Vert \Psi _{J}\right\Vert _{L^{2}\left( \mu
\right) }$, and the `pivotal' bound%
\begin{equation*}
\left\vert \left\langle T^{\alpha }\left( \varphi \nu \right) ,\Psi
_{J}\right\rangle _{\omega }\right\vert \lesssim C_{\gamma }\mathrm{P}%
^{\alpha }\left( J,\left\vert \nu \right\vert \right) \sqrt{\left\vert
J\right\vert _{\omega }}\left\Vert \Psi _{J}\right\Vert _{L^{2}\left( \omega
\right) }\ ,
\end{equation*}%
for any function $\varphi $ with $\left\vert \varphi \right\vert \leq 1$.
\end{lem}

\begin{proof}
Using the Monotonicity Lemma \ref{mono}, followed by $\left\vert \nu _{J^{\prime}}\right\vert \leq \nu $, the Poisson equivalence 
\begin{equation}
\frac{\mathrm{P}^{\alpha }\left( J^{\prime },\nu \right) }{\left\vert
J^{\prime }\right\vert^\frac{1}{n} }\approx \frac{\mathrm{P}^{\alpha }\left( J,\nu
\right) }{\left\vert J\right\vert^\frac{1}{n} },\ \ \ \ \ J^{\prime }\subset J\subset
\gamma J,\ \ \ {\supp}\nu \cap \gamma J=\emptyset ,
\label{Poisson equiv}
\end{equation}%
and the weak frame inequalities for dual martingale differences, we have
\begin{eqnarray*}
&&
\left\vert \sum_{J^{\prime }\in \mathcal{H}}\left\langle T^{\alpha }\left(
\nu _{J^{\prime }}\right) ,\square _{J^{\prime }}^{\omega ,\mathbf{b}^{\ast
}}\Psi _{J}\right\rangle _{\omega }\right\vert 
\lesssim 
\sum_{J^{\prime }\in 
\mathcal{H}}\Phi ^{\alpha }\left( J^{\prime },\left\vert \mu \right\vert
\right) \left\Vert \square _{J^{\prime }}^{\omega ,\mathbf{b}^{\ast }}\Psi
_{J}\right\Vert _{L^{2}\left( \mu \right) }^{\bigstar } \\
&\lesssim &
\left( \sum_{J^{\prime }\in \mathcal{H}}\left( \frac{\mathrm{P}%
^{\alpha }\left( J^{\prime },\nu \right) }{\left\vert J^{\prime }\right\vert^\frac{1}{n} 
}\right) ^{2}\left\Vert \bigtriangleup _{J^{\prime }}^{\omega ,\mathbf{b}%
^{\ast }}x\right\Vert _{L^{2}\left( \omega \right)}^{\spadesuit 2}\right)
^{\frac{1}{2}}\left( \sum_{J^{\prime }\in \mathcal{H}}\left\Vert \square
_{J^{\prime }}^{\omega ,\mathbf{b}^{\ast }}\Psi _{J}\right\Vert
_{L^{2}\left( \omega \right) }^{\bigstar 2}\right) ^{\frac{1}{2}} \\
&&
+\left( \sum_{J^{\prime }\in \mathcal{H}}\left( \frac{\mathrm{P}_{1+\delta
}^{\alpha }\left( J^{\prime },\left\vert \mu \right\vert \right) }{%
\left\vert J^{\prime }\right\vert^\frac{1}{n} }\right) ^{2}\left\Vert x-m_{J^{\prime
}}\right\Vert _{L^{2}\left( \mathbf{1}_{J^{\prime }}\omega \right)
}^{2}\right) ^{\frac{1}{2}}\left( \sum_{J^{\prime }\in \mathcal{H}%
}\left\Vert \square _{J^{\prime }}^{\omega ,\mathbf{b}^{\ast }}\Psi
_{J}\right\Vert _{L^{2}\left( \omega \right) }^{\bigstar 2}\right) ^{\frac{1%
}{2}} \\
&\lesssim &
 \!\!\!\!\!\frac{\mathrm{P}^{\alpha }\left( J,\nu \right) }{%
\left\vert J\right\vert^\frac{1}{n} } \!\!\left\Vert \mathsf{Q}_{\mathcal{H}}^{\omega
,\mathbf{b}^{\ast }}x\right\Vert _{L^{2}\left( \omega \right) }^{\spadesuit
}\!\!\left\Vert \Psi _{J}\right\Vert _{L^{2}\left( \omega \right) }\!+\!\frac{1}{\gamma ^{\delta ^{\prime }}} \frac{\mathrm{P}_{1+\delta ^{\prime
}}^{\alpha }\left( J,\nu \right) }{\left\vert J\right\vert^\frac{1}{n} }
\!\!\left\Vert x-m_{J}\right\Vert _{L^{2}\left( \mathbf{1}_{J}\omega \right)
}\!\!\left\Vert \Psi _{J}\right\Vert _{L^{2}\left( \omega \right) }\ .
\end{eqnarray*}%
The last inequality follows from the following calculation using Haar
projections $\bigtriangleup _{K}^{\omega }$:

\begin{eqnarray}
&&\!\!\!\!\!\!\!\!\!\!\!\!\!\!\!\!\!\!\!\!\!\!\!\!\!
\sum_{J^{\prime }\in \mathcal{H}}\left( \frac{\mathrm{P}_{1+\delta
}^{\alpha }\left( J^{\prime },\nu \right) }{\left\vert J^{\prime
}\right\vert^\frac{1}{n} }\right) ^{2}\left\Vert x-m_{J^{\prime }}\right\Vert
_{L^{2}\left( \mathbf{1}_{J^{\prime }}\omega \right) }^{2}
\label{Haar trick} \\
&=&
\sum_{J^{\prime }\in \mathcal{H}}\left( \frac{\mathrm{P}_{1+\delta
}^{\alpha }\left( J^{\prime },\nu \right) }{\left\vert J^{\prime
}\right\vert^\frac{1}{n} }\right) ^{2}\sum_{J^{\prime \prime }\subset J^{\prime
}}\left\Vert \bigtriangleup _{J^{\prime \prime }}^{\omega }x\right\Vert
_{L^{2}\left( \omega \right) }^{2}\notag \\
&=&
\sum_{J^{\prime \prime }\subset J}\left\{
\sum_{J^{\prime }:\ J^{\prime \prime }\subset J^{\prime }\subset J}\left( \frac{\mathrm{P}_{1+\delta }^{\alpha }\left( J^{\prime },\nu \right) }{%
\left\vert J^{\prime }\right\vert^\frac{1}{n} }\right) ^{2}\right\} \left\Vert
\bigtriangleup _{J^{\prime \prime }}^{\omega }x\right\Vert _{L^{2}\left(
\omega \right) }^{2}  \notag \\
&\lesssim &
\frac{1}{\gamma ^{2\delta ^{\prime }}}\sum_{J^{\prime \prime
}\subset J}\left( \frac{\mathrm{P}_{1+\delta ^{\prime }}^{\alpha }\left(
J^{\prime \prime },\nu \right) }{\left\vert J^{\prime \prime }\right\vert^\frac{1}{n} }%
\right) ^{2}\left\Vert \bigtriangleup _{J^{\prime \prime }}^{\omega
}x\right\Vert _{L^{2}\left( \omega \right) }^{2}\notag \\
&\leq&
\frac{1}{\gamma
^{2\delta ^{\prime }}}\left( \frac{\mathrm{P}_{1+\delta ^{\prime }}^{\alpha
}\left( J,\nu \right) }{\left\vert J\right\vert^\frac{1}{n} }\right) ^{2}\sum_{J^{\prime
\prime }\subset J}\left\Vert \bigtriangleup _{J^{\prime \prime }}^{\omega
}x\right\Vert _{L^{2}\left( \omega \right) }^{2}\ ,  \notag
\end{eqnarray}%
which in turn follows from (recalling $\delta =2\delta ^{\prime }$ and  $\left\vert J^{\prime }\right\vert^\frac{1}{n} +\left\vert y-c_{J^{\prime }}\right\vert
\approx \left\vert J\right\vert ^\frac{1}{n}+\left\vert y-c_{J}\right\vert $ and $\frac{%
\left\vert J\right\vert }{\left\vert J\right\vert +\left\vert
y-c_{J}\right\vert }\leq \frac{1}{\gamma }$ for $y\in \mathbb{R}^n\backslash
\gamma J$)%
\begin{eqnarray*}
\sum_{J^{\prime }:\ J^{\prime \prime }\subset J^{\prime }\subset J}\left( 
\frac{\mathrm{P}_{1+\delta }^{\alpha }\left( J^{\prime },\nu \right) }{%
\left\vert J^{\prime }\right\vert^\frac{1}{n} }\right) ^{2}
&=&\!\!\!\!\!\!
\sum_{J^{\prime }:\
J^{\prime \prime }\subset J^{\prime }\subset J}\left\vert J^{\prime
}\right\vert ^\frac{2\delta}{n}\Bigg( \int_{\mathbb{R}^n\backslash \gamma J}\frac{1}{\left( \left\vert J^{\prime }\right\vert^\frac{1}{n} +\left\vert y-c_{J^{\prime
}}\right\vert \right) ^{n+1+\delta -\alpha }}d\nu \left( y\right) \Bigg) ^{2}
\\
&\lesssim &\!\!\!\!\!\!
\sum_{J^{\prime }:\ J^{\prime \prime }\subset J^{\prime }\subset
J}\frac{1}{\gamma ^{2\delta ^{\prime }}}\frac{\left\vert J^{\prime
}\right\vert ^\frac{2\delta}{n}}{\left\vert J\right\vert ^\frac{2\delta}{n}}\Bigg( \int_{\mathbb{R}^n\backslash \gamma J}\frac{\left\vert J\right\vert ^\frac{\delta ^{\prime
}}{n}}{\left( \left\vert J\right\vert^\frac{1}{n} +\left\vert y-c_{J}\right\vert \right)
^{n+1+\delta ^{\prime }-\alpha }}d\nu \left( y\right) \Bigg) ^{2} \\
&=&\!\!\!\!\!\!
\frac{1}{\gamma^{2\delta'} }\left( \sum_{J^{\prime }:\
J^{\prime \prime }\subset J^{\prime }\subset J}\frac{\left\vert J^{\prime}\right\vert ^\frac{2\delta}{n}}{\left\vert J\right\vert ^\frac{2\delta}{n}}\right) \left( 
\frac{\mathrm{P}_{1+\delta ^{\prime }}^{\alpha }\left( J,\nu \right) }{%
\left\vert J\right\vert^\frac{1}{n} }\right) ^{2}
\lesssim
\frac{1}{\gamma ^{2\delta
^{\prime }}}\left( \frac{\mathrm{P}_{1+\delta ^{\prime }}^{\alpha }\left(
J,\nu \right) }{\left\vert J\right\vert^\frac{1}{n} }\right) ^{2}.
\end{eqnarray*}%
Finally we obtain the `energy' estimate from the equality%
\begin{equation*}
\Psi _{J}=\sum_{J^{\prime }\subset J}\square _{J^{\prime }}^{\omega ,\mathbf{%
b}^{\ast }}\Psi _{J}\ ,\ \ \ \ \ (\text{since }\Psi _{J}\text{ has vanishing 
}\omega \text{-mean)},\text{ }
\end{equation*}%
and we obtain the `pivotal' bound from the inequality%
\begin{equation*}
\sum_{J^{\prime \prime }\subset J}\left\Vert \bigtriangleup _{J^{\prime
\prime }}^{\omega ,\mathbf{b}^{\ast }}x\right\Vert _{L^{2}\left( \omega
\right) }^{\spadesuit 2}\lesssim \left\Vert \left( x-m_{J}\right)
\right\Vert _{L^{2}\left( \mathbf{1}_{J}\omega \right) }^{2}\leq \left\vert
J\right\vert ^{2}\left\vert J\right\vert _{\omega }\ .
\end{equation*}
\end{proof}

\subsection{Organization of the proof}

We adapt the proof of the main theorem in \cite{SaShUr9}, but beginning
instead with the decomposition of Hyt\"{o}nen and Martikainen \cite{HyMa},
to obtain the norm inequality%
\begin{equation*}
\mathfrak{N}_{T^{\alpha }}\lesssim \mathfrak{T}_{T^{\alpha }}^{\mathbf{b}}+%
\mathfrak{T}_{T^{\alpha }}^{\mathbf{b}^{\ast }}+\sqrt{\mathfrak{A}%
_{2}^{\alpha }}+\mathfrak{E}_{2}^{\alpha }
\end{equation*}%
under the \emph{apriori} assumption $\mathfrak{N}_{T^{\alpha }}<\infty $,
which is achieved by considering one of the truncations $T_{\sigma ,\delta
,R}^{\alpha }$ defined in (\ref{def truncation}) above. This will be carried
out in the next four sections of this paper. In the next section we consider
the various form splittings and reduce matters to the \emph{disjoint} form,
the \emph{nearby} form and the \emph{main below} form. Then these latter
three forms are taken up in the subsequent three sections, using material
from the appendices.

A major source of difficulty will arise in the infusion of goodness for the
cubes $J$ into the below form where the sum is taken over all pairs $\left(
I,J\right) $ such that $\ell \left( J\right) \leq \ell \left( I\right) $. We
will infuse goodness in a weak way pioneered by Hyt\"{o}nen and Martikainen
in a one weight setting. This weak form of goodness is then exploited in all
subsequent constructions by typically replacing $J$ by $J^{\maltese }$ in
defining relations, where $J^{\maltese }$ is the smallest cube $K$ for which $J$ is good w.r.t. $K$
and beyond.

Another source of difficulty arises in the treatment of the nearby form in
the setting of two weights. The one weight proofs in \cite{HyMa} and \cite%
{LaMa} relied strongly on a property peculiar to the one weight setting -
namely the fact already pointed out in Remark \ref{special}\ above that both
of the Poisson integrals are bounded, namely $\mathrm{P}^{\alpha }\left(
Q,\mu \right) \lesssim 1$ and $\mathcal{P}^{\alpha }\left( Q,\mu \right)
\lesssim 1$. We will circumvent this difficulty by combining a recursive
energy argument with the full testing conditions assumed for the \emph{%
original} testing functions $b_{Q}^{{orig}}$, before these
conditions were suppressed by corona constructions that delivered only weak
testing conditions for the new testing functions $b_{Q}$.

Of particular importance will be a result proved in the Appendix, where we
show that the functional energy for an arbitrary pair of grids is controlled
by the Muckenhoupt and energy side conditions. The somewhat lengthy proof of
this latter assertion is similar to the corresponding proof in the $T1$
setting - see e.g. \cite{SaShUr9} - but requires a different decomposition
of the stopping cubes into `Whitney cubes' in order to accomodate the weaker
notion of goodness used here.

Here is
a brief schematic diagram of the splittings and decompositions we will
describe below, with associated bounds given in a box.
\begin{equation*}
\fbox{$%
\begin{array}{ccccccc}
\Theta \left( f,g\right) &  &  &  &  &  &  \\ 
\downarrow &  &  &  &  &  &  \\ 
\Theta _{2}^{{good}}=\mathsf{B}_{\Subset _{\mathbf{\rho }}}\left(
f,g\right) & + & \Theta _{1}=\mathsf{B}_{\cap }\left( f,g\right) & + & 
\Theta _{3}=\mathsf{B}_{\diagup }\left( f,g\right) & + & \Theta _{2}^{%
{bad}}\left( f,g\right) \\ 
\downarrow &  & \fbox{$\mathcal{NTV}_{\alpha }$} &  & \fbox{$\mathcal{NTV}%
_{\alpha }+\sqrt{\delta }\mathfrak{N}_{T^{\alpha }}$} &  & \fbox{$2^{-%
\mathbf{r}\varepsilon }\mathfrak{N}_{T^{\alpha }}$} \\ 
\downarrow &  &  &  &  &  &  \\ 
\mathsf{T}_{{diagonal}}\left( f,g\right) & + & \mathsf{T}_{{%
far}{below}}\left( f,g\right) & + & \mathsf{T}_{{far}%
{above}}\left( f,g\right) & + & \mathsf{T}_{{disjoint}%
}\left( f,g\right) \\ 
\downarrow &  & \downarrow &  & \fbox{$\emptyset $} &  & \fbox{$\emptyset $}
\\ 
\downarrow &  & \downarrow &  &  &  &  \\ 
\mathsf{B}_{\Subset _{\mathbf{\rho }}}^{A}\left( f,g\right) &  & \mathsf{T}_{%
{far}{below}}^{1}\left( f,g\right) & + & \mathsf{T}_{%
{far}{below}}^{2}\left( f,g\right) &  &  \\ 
\downarrow &  & \fbox{$\mathcal{NTV}_{\alpha }+\mathcal{E}^{\alpha }_{2}$} &  & \fbox{$\mathcal{NTV}_{\alpha }$} &  &  \\ 
\downarrow &  &  &  &  &  &  \\ 
\mathsf{B}_{stop}^{A}\left( f,g\right) & + & \mathsf{B}_{paraproduct}^{A}%
\left( f,g\right) & + & \mathsf{B}_{neighbour}^{A}\left( f,g\right) &  &  \\ 
\fbox{$\mathcal{E}^{\alpha }_{2}+\sqrt{\mathfrak{A}_{2}^{\alpha }}$} & 
& \fbox{$\mathfrak{T}_{T^{\alpha }}$} &  & \fbox{$\sqrt{\mathfrak{A}_2^{\alpha }}$} & 
& 
\end{array}%
$}
\end{equation*}
\section{Form splittings}

\begin{notation}
Fix grids $\mathcal{D}$ and $\mathcal{G}$. We will use $\mathcal{D}$ to
denote the grid associated with $f\in L^{2}\left( \sigma \right) $, and we
will use $\mathcal{G}$ to denote the grid associated with $g\in L^{2}\left(\omega \right) $.
\end{notation}

Now we turn to the probability estimates for martingale differences and
halos that we will use. Recall that given $\overrightarrow{\lambda}=(\lambda_1,...,\lambda_n)$, $0<\lambda_i <\frac{1}{2}$ for all $1\leq i \leq n$, the $%
\lambda $-halo of $J$ is defined to be 
\begin{equation*}
\partial _{\overrightarrow{\lambda} }J\equiv \left( 1+\overrightarrow{\lambda} \right) J\backslash \left(
1-\overrightarrow{\lambda} \right) J.
\end{equation*}%
Suppose $\mu $ is a positive locally finite Borel measure, and that $\mathbf{%
b}$ is a $p$-weakly $\mu $-controlled accretive family for some $p>2$. Then
the following probability estimate holds.\\
$\text{}$\\
\textbf{Bad cube probability estimates.} Suppose that $\mathcal{D}$ and $%
\mathcal{G}$ are independent random dyadic grids. With $\Psi _{\mathcal{G}%
_{k-{bad}}^{\mathcal{D}}}^{\mu ,\mathbf{b}^{\ast }}g\equiv
\sum_{J\in \mathcal{G}_{k-{bad}}^{\mathcal{D}}}\square _{J}^{\mu ,%
\mathbf{b}^{\ast }}g$ equal to the pseudoprojection of $g$ onto $k$-bad $%
\mathcal{G}$-cubes, we have%
\begin{eqnarray}
\boldsymbol{E}_{\Omega }^{\mathcal{D}}\left( \left\Vert \Psi _{\mathcal{G}%
_{k-{bad}}^{\mathcal{D}}}^{\mu ,\mathbf{b}^{\ast }}g\right\Vert
_{L^{2}\left( \mu \right) }^{2}\right) \notag
&\lesssim&
\boldsymbol{E}_{\Omega }^{%
\mathcal{D}}\left( \sum_{J\in \mathcal{G}_{k-{bad}}^{\mathcal{D}}}%
\left[ \left\Vert \square _{J,\mathcal{G}}^{\mu ,\mathbf{b}^{\ast
}}g\right\Vert _{L^{2}\left( \mu \right) }^{2}+\left\Vert \nabla _{J,%
\mathcal{G}}^{\mu }g\right\Vert _{L^{2}\left( \mu \right) }^{2}\right]
\right) \\&\leq& Ce^{-k\varepsilon }\left\Vert g\right\Vert _{L^{2}\left( \mu
\right) }^{2}\ ,
\end{eqnarray}%
where the first inequality is the `weak upper half Riesz' inequality from Appendix A of \cite{SaShUr12} for the
pseudoprojection $\Psi _{\mathcal{G}_{k-{bad}}^{\mathcal{D}}}^{\mu ,%
\mathbf{b}^{\ast }}$, and the second inequality is proved using the frame
inequality in (\ref{main bad prob}) below.

\textbf{Halo probability estimates.} Suppose that $\mathcal{D}$ and $\mathcal{G}
$ are independent random grids. Using the \emph{parameterization by
translations}\ of grids and taking the average over certain translates $\tau
+\mathcal{D}$ of the grid $\mathcal{D}$ we have%
\begin{eqnarray}
\boldsymbol{E}_{\Omega }^{\mathcal{D}}\sum_{I^{\prime }\in \mathcal{D}:\
\ell \left( I^{\prime }\right) \approx \ell \left( J^{\prime }\right)
}\int_{J^{\prime }\cap \partial _{\delta }I^{\prime }}d\omega &\lesssim &%
\mathbb{\delta }\int_{J^{\prime }}d\omega ,\ \ \ \ \ J^{\prime }\in 
\mathfrak{C}\left( J\right) ,J\in \mathcal{G},  \label{hand'} \\
\boldsymbol{E}_{\Omega }^{\mathcal{G}}\sum_{J^{\prime }\in \mathcal{G}:\
\ell \left( J^{\prime }\right) \approx \ell \left( I^{\prime }\right)
}\int_{I^{\prime }\cap \partial _{\delta }J^{\prime }}d\sigma &\lesssim &%
\mathbb{\delta }\int_{I^{\prime }}d\sigma ,\ \ \ \ \ I^{\prime }\in 
\mathfrak{C}\left( I\right) ,I\in \mathcal{D},  \notag
\end{eqnarray}%
and where the expectations $\boldsymbol{E}_{\Omega }^{\mathcal{D}}$ and $%
\boldsymbol{E}_{\Omega }^{\mathcal{G}}$ are taken over grids $\mathcal{D}$
and $\mathcal{G}$ respectively. Indeed, it is geometrically evident that for
any fixed pair of side lengths $\ell _{1}\approx \ell _{2}$, the average of
the measure $\left\vert J^{\prime }\cap \partial _{\delta }I^{\prime
}\right\vert _{\omega }$ of the set $J^{\prime }\cap \partial _{\delta
}I^{\prime }$, as a cube $I^{\prime }\in \mathcal{D}$ with side length $%
\ell \left( I^{\prime }\right) =\ell _{1}$ is translated across a cube $%
J^{\prime }\in \mathcal{G}$ of side length $\ell \left( J^{\prime }\right)
=\ell _{2}$, is at most $C\left\vert J^{\prime }\right\vert _{\omega }$.
Using this observation it is now easy to see that (\ref{hand'}) holds.

 In the $\sigma $-iterated corona construction we redefined
the family $\mathbf{b}=\left\{ b_{Q}\right\} _{Q\in \mathcal{D}}$ so that
the new functions $b_{Q}^{{new}}$ are given in terms of the original
functions $b_{Q}^{{orig}}$ by $b_{Q}^{{new}}=\mathbf{1}%
_{Q}b_{A}^{{orig}}$ for $Q\in \mathcal{C}_{A}^{\sigma }$, and of
course we then dropped the superscript ${new}$. We continue to refer
to the triple stopping cubes $A$ as `breaking' cubes even if $b_{A}$ happens
to equal $\mathbf{1}_{A}b_{\pi A}$. The results of Appendix A of \cite{SaShUr12} apply with
this more inclusive definition of `breaking' cubes, and the associated
definition of `broken' children, since only the Carleson condition\ on
stopping cubes is relevant here.

This and Proposition \ref{data} give us the \emph{triple corona decomposition} of $f\!\!=\!\!\!\sum\limits_{A\in \mathcal{A}}\!\mathsf{P}_{\mathcal{C}_{A}}^{\sigma }f$,  where the pseudoprojection $\mathsf{P}_{\mathcal{C}_{A}}^{\sigma }$ is defined as:
$$
\mathsf{P}_{\mathcal{C}_{A}}^{\sigma }f=\sum_{I\in\mathcal{C}_A}\square_I^{\mu,\mathbf{b}}f
$$ 
We now record the main facts proved above for the triple corona.

\begin{lem}
Let $f\in L^2(\sigma)$. We have
\begin{equation*}
f=\sum_{A\in \mathcal{A}}\mathsf{P}_{\mathcal{C}_{A}}^{\sigma }f
\end{equation*}
both in the sense of norm convergence in $L^{2}\left( \sigma \right)$ and pointwise $\sigma$-a.e.
The corona tops $\mathcal{A}$ and stopping bounds $\left\{ \alpha _{\mathcal{%
A}}\left( A\right) \right\} _{A\in \mathcal{A}}$ satisfy properties (1),
(2), (3) and (4) in Definition \ref{general stopping data}, hence constitute
stopping data for $f$. Moreover, $\mathbf{b}=\left\{ b_{I}\right\} _{I\in 
\mathcal{D}}$ is a $\infty$-weakly $\sigma $-controlled accretive family on $%
\mathcal{D}$ with corona tops $\mathcal{A\subset D}$, where $b_{I}=\mathbf{1}%
_{I}b_{A}$ for all $I\in \mathcal{C}_{A}$, and the weak corona forward
testing condition holds uniformly in coronas, i.e.%
\begin{equation*}
\frac{1}{\left\vert I\right\vert _{\sigma }}\int_{I}\left\vert T_{\sigma
}^{\alpha }b_{A}\right\vert ^{2}d\sigma \leq C,\ \ \ \ \ I\in \mathcal{C}%
_{A}^{\sigma }\ .
\end{equation*}
\end{lem}
Similar statements hold for $g\in L^2(\omega)$.

We have defined corona decompositions of $f$ and $g$ in the $\sigma $%
-iterated triple corona construction above, but in order to start these
corona decompositions for $f$ and $g$ respectively within the dyadic grids $%
\mathcal{D}$ and $\mathcal{G}$, we need to first restrict $f$ and $g$ to be
supported in a large common cube $Q_{\infty }$. Then we cover $Q_{\infty
}$ with $2^n$ pairwise disjoint cubes $I_{\infty }\in \mathcal{D}$ with $%
\ell \left( I_{\infty }\right) =\ell \left( Q_{\infty }\right) $, and
similarly cover $Q_{\infty }$ with $2^n$ pairwise disjoint cubes $%
J_{\infty }\in \mathcal{G}$ with $\ell \left( J_{\infty }\right) =\ell
\left( Q_{\infty }\right) $. We can now use the broken martingale decompositions, together with random surgery, to reduce matters to
consideration of the four forms%
\begin{equation*}
\sum_{I\in \mathcal{D}:\ I\subset I_{\infty }}\sum_{J\in \mathcal{G}:\
J\subset J_{\infty }}\int \left( T_{\sigma }^{\alpha }\square _{I}^{\sigma ,%
\mathbf{b}}f\right) \square _{J}^{\omega ,\mathbf{b}^{\ast }}gd\omega ,
\end{equation*}%
with $I_{\infty }$ and $J_{\infty }$ as above, and where we can then use the
cubes $I_{\infty }$ and $J_{\infty }$ as the starting cubes in our
corona constructions below. Indeed, the identities in \cite[Lemma 3.5]{HyMa}), give%
\begin{eqnarray*}
f &=&\sum_{I\in \mathcal{D}:\ I\subset I_{\infty },\ \ell \left( I\right) \geq
2^{-N}}\!\!\!\!\!\!\!\!\square _{I}^{\sigma ,\mathbf{b}}f+\mathbb{F}_{I_{\infty }}^{\sigma ,%
\mathbf{b}}f, \\
g 
&=&
\sum_{J\in \mathcal{G}:\ J\subset J_{\infty },\ \ell \left( J\right)
\geq 2^{-N}}\!\!\!\!\!\!\!\!\!\square _{J}^{\omega ,\mathbf{b}^{\ast }}g+\mathbb{F}_{J_{\infty
}}^{\omega ,\mathbf{b}^{\ast }}g,
\end{eqnarray*}%
which can then be used to write the bilinear form $\int \left( T_{\sigma
}f\right) gd\omega $ as a sum of the forms

\begin{eqnarray}
&&
\sum_{\substack{2^{n+1} \text{pairs}  \\ \left( I_{\infty },J_{\infty
}\right) }}\!\left\{ \sum_{\substack{I\in \mathcal{D}\\ I\subset I_{\infty }}}\sum_{\substack{J\in 
\mathcal{G}\\ J\subset J_{\infty }}}\int \left( T_{\sigma }^{\alpha }\square
_{I}^{\sigma ,\mathbf{b}}f\right) \square _{J}^{\omega ,\mathbf{b}^{\ast
}}gd\omega
+\sum_{\substack{I\in \mathcal{D}\\ I\subset I_{\infty }}}\int \left(
T_{\sigma }^{\alpha }\square _{I}^{\sigma ,\mathbf{b}}f\right) \mathbb{F}%
_{J_{\infty }}^{\omega ,\mathbf{b}^{\ast }}g d\omega \right.\notag \\
&& 
\left. \hspace{1.7cm} +\sum_{J\in \mathcal{G}%
:\ J\subset J_{\infty }}\int \left( T_{\sigma }^{\alpha }\mathbb{F}%
_{I_{\infty }}^{\sigma ,\mathbf{b}}f\right) \square _{J}^{\omega ,\mathbf{b}%
^{\ast }}gd\omega   
+\int \left( T_{\sigma }^{\alpha }\mathbb{F}_{I_{\infty
}}^{\sigma ,\mathbf{b}}f\right) \mathbb{F}_{J_{\infty }}^{\omega ,\mathbf{b}%
^{\ast }}gd\omega \right\}
\label{sum of forms} 
\end{eqnarray}
taken over the $2^{n+1}$ pairs of cubes $\left( I_{\infty },J_{\infty
}\right) $ above. The second, third and fourth sums in (\ref{sum of forms}) can be controlled
using testing and random surgery. For example, for the second sum we have

\begin{eqnarray*}
\left\vert \sum_{I\in \mathcal{D}:\ I\subset I_{\infty }}\!\!\int \!\!\left(
T_{\sigma }^{\alpha }\square _{I}^{\sigma ,\mathbf{b}}f\right) \mathbb{F}%
_{J_{\infty }}^{\omega ,\mathbf{b}^{\ast }}gd\omega \right\vert
&\leq &
\left\vert
\int_{I_\infty \cap J_\infty} \left( \sum_{I\in \mathcal{D}:\ I\subset I_{\infty }}\!\!\!\!\!\square
_{I}^{\sigma ,\mathbf{b}}f\right) T_{\omega }^{\alpha ,\ast }\left( \mathbb{F}_{J_{\infty }}^{\omega ,\mathbf{b}^{\ast }}g\right) d\sigma \right\vert\\
&&+
\left\vert   \int_{I_\infty \cap\left( (1+\delta)J_\infty\backslash J_\infty\right)} \left( \sum_{I\in \mathcal{D}:\ I\subset I_{\infty }}\!\!\!\!\!\square
_{I}^{\sigma ,\mathbf{b}}f\right) T_{\omega }^{\alpha ,\ast }\left( \mathbb{F}_{J_{\infty }}^{\omega ,\mathbf{b}^{\ast }}g\right) d\sigma    \right\vert
\\
&&+
\left\vert   \int_{I_\infty \backslash (1+\delta)J_\infty} \left( \sum_{I\in \mathcal{D}:\ I\subset I_{\infty }}\!\!\!\!\!\square
_{I}^{\sigma ,\mathbf{b}}f\right) T_{\omega }^{\alpha ,\ast }\left( \mathbb{F}_{J_{\infty }}^{\omega ,\mathbf{b}^{\ast }}g\right) d\sigma    \right\vert\\
&\equiv&
A_1+A_2+A_3
\end{eqnarray*}
So we are left with bounding $A_1,A_2,A_3$. We have
\begin{eqnarray*}
A_1\leq 
\left(\int_{I_\infty} \left\vert \sum_{I\in \mathcal{D}:\ I\subset I_{\infty }}\!\!\!\!\!\square
_{I}^{\sigma ,\mathbf{b}}f\right\vert^2d\sigma\right)^\frac{1}{2}\left( \int_{J_\infty}\left\vert T_{\omega }^{\alpha ,\ast }\left( \mathbb{F}_{J_{\infty }}^{\omega ,\mathbf{b}^{\ast }}g\right)\right\vert^2 d\sigma\right)^\frac{1}{2} 
\end{eqnarray*}
and since $\mathbb{F}_{J_{\infty }}^{\omega ,\mathbf{b}^{\ast }}g=b_{J_{\infty
}}^{\ast }\frac{E_{J_{\infty }}^{\omega }g}{E_{J_{\infty }}^{\omega
}b_{J_{\infty }}^{\ast }}$ is $b_{J_{\infty }}^{\ast }$ times an `accretive'
average of $g$ on $J_{\infty }$,
we get
\begin{eqnarray*}
A_1&\leq& \left\Vert \sum_{I\in \mathcal{D}:\ I\subset I_{\infty }}\!\!\!\!\!\square
_{I}^{\sigma ,\mathbf{b}}f\right\Vert _{L^{2}\left( \sigma \right)
}\left(\int_{J_\infty}
\left|T_{\omega}^{\alpha,*}(\mathbf{1}_{J_\infty}b^{*}_{J_\infty})\right|^2 d\sigma \right)^{\!\!\frac{1}{2}}
|E^\omega_{J_\infty}g|\cdot \frac{1}{c_{\mathbf{b}^{*}}|J_\infty|_\omega}\\
&\lesssim& \mathfrak{T}^\textbf{b*}_{T^{\!\alpha,\!*}}\left\Vert f\right\Vert _{L^{2}\left( \sigma \right)}\left\Vert g\right\Vert _{L^{2}\left( \omega \right)}
\end{eqnarray*}
where in the last inequality we used the frame estimates \eqref{FRAME} and the dual testing condition on $b_{J_\infty}^\ast$. 

For $A_2$ we use expectation on the grid $\mathcal{G}$.

\begin{eqnarray*}
\boldsymbol{E}^{\mathcal{G}}A_2 
\!\!\!\!&\leq &\!\!\!\!
\boldsymbol{E}^{\mathcal{G}%
}\int_{I_{\infty }\cap \left[ \left( 1\!+\!\delta \right) J_{\infty }\setminus
J_{\infty }\right] }\left\vert \sum_{I\in \mathcal{D}:\ I\subset I_{\infty
}}\square _{I}^{\sigma ,\mathbf{b}}f\right\vert \left\vert T_{\omega
}^{\alpha ,\ast }\left( \mathbb{F}_{J_{\infty }}^{\omega ,\mathbf{b}^{\ast
}}g\right) \right\vert d\sigma  \\
&\leq &\!\!\!\!
\boldsymbol{E}^{\mathcal{G}}\left( \int_{I_{\infty }\cap \left[
\left( 1\!+\!\delta \right) J_{\infty }\setminus J_{\infty }\right] }\left\vert
\sum_{I\in \mathcal{D}:\ I\subset I_{\infty }}\!\!\!\!\!\square _{I}^{\sigma ,\mathbf{b}}f\right\vert ^{2}\!\!d\sigma \right) ^{\!\!\frac{1}{2}}\!\!\left( \int_{I_{\infty
}\cap \left[ \left( 1\!+\!\delta \right) J_{\infty }\setminus J_{\infty }\right]
}\left\vert T_{\omega }^{\alpha ,\ast }\left( \mathbb{F}_{J_{\infty
}}^{\omega ,\mathbf{b}^{\ast }}g\right) \right\vert ^{2}\!d\sigma \right) ^{\!\!\!
\frac{1}{2}}\\
&\leq &\!\!\!\!
\left( \boldsymbol{E}^{\mathcal{G}}\int_{I_{\infty }\cap \left[
\left( 1+\delta \right) J_{\infty }\setminus J_{\infty }\right] }\left\vert
\sum_{I\in \mathcal{D}:\ I\subset I_{\infty }}\square _{I}^{\sigma ,\mathbf{b%
}}f\right\vert ^{2}d\sigma \right) ^{\frac{1}{2}}\left( \mathfrak{N}%
_{T^{\alpha }}\int \left\vert g\right\vert ^{2}d\omega \right) ^{\frac{1}{2}}\\
&\leq &\!\!\!\!
\left( C\delta \int_{I_{\infty }}\left\vert \sum_{I\in \mathcal{D}:\
I\subset I_{\infty }}\square _{I}^{\sigma ,\mathbf{b}}f\right\vert
^{2}d\sigma \right) ^{\frac{1}{2}}\left( \mathfrak{N}_{T^{\alpha }}\int
\left\vert g\right\vert ^{2}d\omega \right) ^{\frac{1}{2}}\\
&\leq&\!\!\!\!
\sqrt{C\delta 
\mathfrak{N}_{T^{\alpha }}}\left\Vert f\right\Vert _{L^{2}\left( \sigma \right)}\left\Vert g\right\Vert _{L^{2}\left( \omega \right)}
\end{eqnarray*}

Finally for $A_3$ we use lemma \ref{lemma1} since $\dist(I_\infty\backslash (1+\delta)J_\infty,J_\infty)\approx \delta \ell(J_\infty)$ to get 
$$
A_3\lesssim \sqrt{\mathfrak{A}_2^\alpha}\delta^{\alpha-n}\left\Vert f\right\Vert _{L^{2}\left( \sigma \right)}\left\Vert g\right\Vert _{L^{2}\left( \omega \right)}.
$$
Altogether we get

$$
\boldsymbol{E}_\Omega^\mathcal{G}\left\vert \sum_{\substack{I\in \mathcal{D}\\ I\subset I_{\infty }}}\int \left(
T_{\sigma }^{\alpha }\square _{I}^{\sigma ,\mathbf{b}}f\right) \mathbb{F}%
_{J_{\infty }}^{\omega ,\mathbf{b}^{\ast }}gd\omega \right\vert \!\!  \lesssim \!\! \left( \mathfrak{T}^\textbf{b}_{T^\alpha}\!+\!\sqrt{\mathfrak{A}_2^\alpha}\delta^{\alpha-n}\!+\!\delta\mathfrak{N}_{T^\alpha}\!\!\right)\left\Vert f\right\Vert_{L^2(\sigma)}\!\left\Vert g\right\Vert_{L^2(\omega)}
$$
Similarly we deal with the third and fourth sum of (\ref{sum of forms}). We are left to deal with the first sum in \eqref{sum of forms}.

\subsection{The Hyt\"{o}nen-Martikainen decomposition and weak goodness\label{Subsec HM}}

Now we turn to the various splittings of forms, beginning with the two
weight analogue of the decomposition of Hyt\"{o}nen and Martikainen \cite%
{HyMa}. Let $\mathbf{b}$ (respectively $\mathbf{b}^{\ast }$) be a $\infty$-weakly $\sigma 
$-controlled (respectively $\omega $-controlled) accretive family. Fix the stopping data $\mathcal{A}$ and $\left\{ \alpha _{\mathcal{A}%
}\left( A\right) \right\} _{A\in \mathcal{A}}$ and dual martingale
differences $\square _{I}^{\sigma ,\mathbf{b}}$ constructed above with the
triple iterated coronas, as well as the corresponding data for $g$. We are left with the estimation of the bilinear form $\int\left( T_{\sigma }f\right) gd\omega$
to that of the sum 
\begin{equation*}
\sum_{I\in \mathcal{D}}\sum_{J\in 
\mathcal{G}}\int \left( T^\alpha_{\sigma }\square _{I}^{\sigma ,\mathbf{b}}f\right)
\square _{J}^{\omega ,\mathbf{b}^{\ast }}gd\omega ,
\end{equation*}
We split the form $\left\langle T_{\sigma }^{\alpha }f,g\right\rangle _{\omega }$ into the sum
of two essentially symmetric forms by cube size,%
\begin{eqnarray}
\int \left( T_{\sigma }f\right) gd\omega 
&=&
\left\{ \sum_{\substack{ I\in 
\mathcal{D}:\ J\in \mathcal{G}  \\ \ell \left( J\right) \leq \ell \left(
I\right) }}+\sum_{\substack{ I\in \mathcal{D}:\ J\in \mathcal{G}  \\ \ell
\left( J\right) >\ell \left( I\right) }}\right\} \int \left( T^\alpha_{\sigma
}\square _{I}^{\sigma ,\mathbf{b}}f\right) \square _{J}^{\omega ,\mathbf{b}%
^{\ast }}gd\omega ,  \label{ess symm}\\
&\equiv&
\Theta(f,g)+\Theta^{*}(f,g) \notag
\end{eqnarray}
and focus on the first sum,%
\begin{equation*}
\Theta \left( f,g\right) =\sum_{I\in \mathcal{D}\text{ and }J\in \mathcal{G}%
:\ \ell \left( J\right) \leq \ell \left( I\right) }\left\langle T_{\sigma
}^{\alpha }\square _{I}^{\sigma ,\mathbf{b}}f,\square _{J}^{\omega ,\mathbf{b%
}^{\ast }}\right\rangle _{\omega },
\end{equation*}%
since the second sum is handled dually, but is easier due to the missing
diagonal. Before introducing goodness into the sum, we follow \cite{HyMa}
and split the form $\Theta \left( f,g\right) $ into 3 pieces:
\begin{eqnarray*}
&&\hspace{-0.9cm}
\sum_{I\in \mathcal{D}}\!\left\{ \!\sum_{\substack{ J\in \mathcal{G}:\ \ell
\left( J\right) \leq \ell \left( I\right)  \\ d\left( J,I\right) >2%
\ell \left( J\right) ^{\varepsilon }\ell \left( I\right) ^{1-\varepsilon }}}
\!\!\!+\!\!
\sum_{\substack{ J\in \mathcal{G}:\ \ell \left( J\right) \leq 2^{-\mathbf{%
\mathbf{r} }}\ell \left( I\right)  \\ d\left( J,I\right) \leq 2\ell \left(
J\right) ^{\varepsilon }\ell \left( I\right) ^{1-\varepsilon }}}
\!\!\!+\!\!
\sum_{\substack{ J\in \mathcal{G}:\ 2^{-\mathbf{\mathbf{r} }}\ell \left( I\right)
<\ell \left( J\right) \leq \ell \left( I\right)  \\ d\left( J,I\right) \leq 2\ell \left( J\right) ^{\varepsilon }\ell \left( I\right)
^{1-\varepsilon }}}\!\!\right\} \!\!
\left\langle\! T_{\sigma
}^{\alpha }\square _{I}^{\sigma ,\mathbf{b}}\!f,\square _{J}^{\omega ,\mathbf{b%
}^{\ast }}\!\right\rangle _{\!\omega }  \\
&\equiv &
\Theta _{1}\left( f,g\right) +\Theta _{2}\left( f,g\right) +\Theta
_{3}\left( f,g\right) \ ,
\end{eqnarray*}%
where $\varepsilon >0$ will be chosen to satisfy $0<\varepsilon <\frac{1}{%
n+1-\alpha }$ later. Now the disjoint form $\Theta _{1}\left( f,g\right) $
can be handled by `long-range' and `short-range' arguments which we give in a
section below, and the nearby form $\Theta _{3}\left( f,g\right) $ will
be handled using surgery methods and a new recursive
argument involving energy conditions and the `original' testing functions
discarded in the corona construction. The remaining form $\Theta_2(f,g)$ will be treated further in this section after introducing weak goodness.

\subsubsection{Good cubes with `body'}. We begin with the weaker extension of goodness introduced in \cite{HyMa},
except that we will make it a bit stronger by replacing the skeleton `$%
{skel}K$' of a cube $K$, as used in \cite{HyMa}, by a larger
collection of points `${body}K$', which we call the dyadic body of $%
K $. This modification will prove useful in establishing the Straddling
Lemma in the treatment of the stopping form in Section \ref{Sec stop} below.
Let $\mathcal{P}$ denote the collection of all cubes in $\mathbb{R}^n$.
The content of the next four definitions is inspired by, or sometimes
identical with, that already appearing in the work of Nazarov, Treil and
Volberg in \cite{NTV1} and \cite{NTV3}.

\begin{dfn}
Given a dyadic cube $K \in \R^n$, we define $W(K)$ to be the Whitney cubes in $K$. Namely, $S \in W(K)$ if:
\begin{itemize}
    \item $3S\subset K$.
    \item $S' \cap S\neq \emptyset$ and $3S' \subset K$ imply $S' \subset S$.
\end{itemize}
\end{dfn}

\begin{dfn}\label{body}
We define the \textit{dyadic body} `$body K$' of a dyadic cube $K\in \R^n$ by
$$
bodyK=\bigcup_{S \in W(K)} \partial S
$$
where $\partial S$ is the boundary of $S$.
\end{dfn}

\begin{dfn}\label{good arb}
Let $0<\epsilon<1$. For dyadic cubes $J,K \in \R^n$ with $\ell(J)\leq \ell(K)$ we define $J$ to be $\epsilon-$good in $K$ if
\begin{equation}\label{good}
\dist(J,bodyK)>2\ell(J)^\epsilon\ell(K)^{1-\epsilon}
\end{equation}
and we say it is $\epsilon-$bad in $K$ if \eqref{good} fails. 
\end{dfn}
\begin{dfn}
\label{good two grids}
Let $\mathcal{D}$ and $\mathcal{G}$ be two dyadic grids in $\R^n$. Define $\mathcal{G}^\mathcal{D}_{(k,\epsilon)-good}$ to consist of those cubes $J \in \mathcal{G}$ such that $J$ is $\epsilon-$good inside every cube $K\in \mathcal{D}$ with $K\cap J \neq \emptyset$ and $\ell(K)\geq 2^k\ell(J)$.
\end{dfn}

\subsubsection{Grid probability}

As pointed out on page 14 of \cite{HyMa} by Hyt\"{o}nen and Martikainen,
there are subtle difficulties associated in using dual martingale
decompositions of functions which depend on the entire dyadic grid, rather
than on just the local cube in the grid. We will proceed at first in the
spirit of \cite{HyMa}, and the goodness that we will infuse below into the
main `below' form $\mathsf{B}_{\Subset _{\mathbf{r}}}\left( f,g\right) $
will be the Hyt\"{o}nen-Martikainen `weak' version of NTV goodness, but
using the body `${body}I$' of a cube rather than\ its skeleton `%
${skel}I$': every pair $\left( I,J\right) \in \mathcal{D}\times 
\mathcal{G}$ that arises in the form $\mathsf{B}_{\Subset _{\mathbf{r}%
}}\left( f,g\right) $ will satisfy $J\in \mathcal{G}_{\left( k,\varepsilon
\right) -{good}}^{\mathcal{D}}$ where $\ell \left( I\right)
=2^{k}\ell \left( J\right) $.

Now we return to the martingale differences $\square _{I}^{\sigma ,\mathbf{b}%
}$ and $\square _{J}^{\omega ,\mathbf{b}^{\ast }}$ with controlled families $%
\mathbf{b}$ and $\mathbf{b}^{\ast }$ in $\mathbb{R}^n$. When we
want to emphasize that the grid in use is $\mathcal{D}$ or $\mathcal{G}$, we
will denote the martingale difference by $\square _{I,\mathcal{D}}^{\sigma ,%
\mathbf{b}}$, and similarly for $\square _{J,\mathcal{G}}^{\omega ,\mathbf{b}%
^{\ast }}$. Recall Definition \ref{good arb} for the meaning of when an
cube $J$ is $\varepsilon $-${bad}$ with respect to another
cube $K$.

\begin{dfn}
\label{bad in grid}We say that $J\in \mathcal{P}$ is $k$-${bad}$ in
a grid $\mathcal{D}$ if there is a cube $K\in \mathcal{D}$ with $\ell
\left( K\right) =2^{k}\ell \left( J\right) $ such that $J$ is $\varepsilon $-%
${bad}$ with respect to $K$ (context should eliminate any ambiguity
between the different use of $k$-${bad}$ when $k\in \mathbb{N}$ and $%
\varepsilon $-${bad}$ when $0<\varepsilon <\frac{1}{2}$).
\end{dfn}
Following \cite{SaShUr12} we know that in one dimension for an interval $J$ and grids $\mathcal{D}_0$
\begin{equation}
\boldsymbol{P}_{\Omega }^{\mathcal{D}_0}\left( \mathcal{D}_0:J\text{ is }k\text{-%
}{bad}\text{ in }\mathcal{D}_0\right) \equiv \int_{\Omega }\mathbf{1}%
_{\left\{ \mathcal{D}_0:\ J\text{ is }k\text{-}{bad}\text{ in }%
\mathcal{D}_0\right\} }d\mu _{\Omega }\left( \mathcal{D}_0\right) \leq
C\varepsilon k2^{-\varepsilon k}.  \label{key prob}
\end{equation}
Thus we conclude:
\begin{equation}
\boldsymbol{P}_{\Omega }^{\mathcal{D}_0}\left( \mathcal{D}_0:J\text{ is }k\text{-%
}{good}\text{ in }\mathcal{D}_0\right) \geq
1-C\varepsilon k2^{-\varepsilon k}.  \label{key prob good}
\end{equation}
Now for a cube $J$ to be good in our $n$-dimensional setting, it needs to be good in each side. So, we conclude that 
\begin{equation}
\boldsymbol{P}_{\Omega }^{\mathcal{D}}\left( \mathcal{D}:J\text{ is }k\text{-%
}{good}\text{ in }\mathcal{D}\right) \geq
(1-C\varepsilon k2^{-\varepsilon k})^n.  \label{key prob good n dimensions}
\end{equation}
and therefore a cube is bad with probability bounded by:
\begin{equation}
\boldsymbol{P}_{\Omega }^{\mathcal{D}}\left( \mathcal{D}:J\text{ is }k\text{-%
}{bad}\text{ in }\mathcal{D}\right) \leq
1-(1-C\varepsilon k2^{-\varepsilon k})^n.  \label{key prob bad n dimensions}
\end{equation}

Then we obtain from (\ref{key prob bad n dimensions}), using the lower frame inequality, the
expectation estimate%
\begin{eqnarray*}
&&\int_{\Omega }\sum_{J\in \mathcal{G}_{k-{bad}}^{\mathcal{D}}}\left[
\left\Vert \square _{J,\mathcal{G}}^{\omega ,\mathbf{b}^{\ast }}g\right\Vert
_{L^{2}\left( \omega \right) }^{2}+\left\Vert \bigtriangledown _{J,\mathcal{G}%
}^{\omega }g\right\Vert _{L^{2}\left( \omega \right) }^{2}\right] d\mu
_{\Omega }\left( \mathcal{D}\right) \\
&=&\sum_{J\in \mathcal{G}}\left[ \left\Vert \square _{J,\mathcal{G}}^{\omega
,\mathbf{b}^{\ast }}g\right\Vert _{L^{2}\left( \omega \right)
}^{2}+\left\Vert \bigtriangledown _{J,\mathcal{G}}^{\omega }g\right\Vert _{L^{2}\left(
\omega \right) }^{2}\right] \int_{\Omega }\mathbf{1}_{\left\{ \mathcal{D}:\ J%
\text{ is }k\text{-}{bad}\text{ in }\mathcal{D}\right\} }d\mu
_{\Omega }\left( \mathcal{D}\right) \\
&\leq &(1-(1-C\varepsilon k2^{-\varepsilon k})^n)\sum_{J\in \mathcal{G}}\left[
\left\Vert \square _{J,\mathcal{G}}^{\omega ,\mathbf{b}^{\ast }}g\right\Vert
_{L^{2}\left( \omega \right) }^{2}+\left\Vert \bigtriangledown _{J,\mathcal{G}%
}^{\omega }g\right\Vert _{L^{2}\left( \omega \right) }^{2}\right]\\ 
&\leq&
(1-(1-C\varepsilon k2^{-\varepsilon k})^n)\left\Vert g\right\Vert _{L^{2}\left( \omega
\right) }^{2}\ ,
\end{eqnarray*}%
where $\bigtriangledown _{J,\mathcal{G}}^{\omega }$ denotes the `broken' Carleson
averaging operator in (\ref{Carleson avg op}) that depends on the broken
children in the grid $\mathcal{G}$. Altogether then it follows easily that%
\begin{equation}
\boldsymbol{E}_{\Omega }^{\mathcal{D}}\left( \sum_{J\in \bigcup_{\ell
=k}^{\infty }\mathcal{G}_{\ell -{bad}}^{\mathcal{D}}}\!\!\!\left[
\left\Vert \square _{J,\mathcal{G}}^{\omega ,\mathbf{b}^{\ast }}g\right\Vert
_{L^{2}\left( \omega \right) }^{2}\!+\!\left\Vert \bigtriangledown _{J,\mathcal{G}%
}^{\omega }g\right\Vert _{L^{2}\left( \omega \right) }^{2}\right] \right)
\!\leq\! (1\!-\!(1\!-\!C\varepsilon k2^{-\varepsilon k})^n)\!\left\Vert g\right\Vert _{L^{2}\left(
\omega \right) }^{2} \label{main bad prob}
\end{equation}%
for some large positive constant $C$.

From such inequalities summed for $k\geq \mathbf{r}$, it can be concluded as
in \cite{NTV3} that there is an absolute choice of $\mathbf{r}$ depending on 
$0<\varepsilon <\frac{1}{2}$ so that the following holds. Let $%
T\;:\;L^{2}(\sigma )\rightarrow L^{2}(\omega )$ be a bounded linear
operator. We then have the following traditional inequality for two random
grids in the case that $\mathbf{b}$ is an $\infty $-weakly $\mu $%
-controlled accretive family: 
\begin{equation}
\left\Vert T\right\Vert  \leq
2\sup_{\left\Vert f\right\Vert _{L^{2}(\sigma )}=1}\sup_{\left\Vert
g\right\Vert _{L^{2}(\omega )}=1}\boldsymbol{E}_{\Omega }\boldsymbol{E}%
_{\Omega ^{\prime }}\left\vert \left\langle \sum_{I,J\in \mathcal{D}_{%
\mathbf{r}-{good}}^{\mathcal{G}}}T\left( \square _{I,\mathcal{D}%
}^{\sigma ,\mathbf{b}}f\right) f,\square _{J,\mathcal{D}}^{\omega ,\mathbf{b}%
^{\ast }}g\right\rangle _{\omega }\right\vert \,.  \label{e.Tgood'}
\end{equation}

However, this traditional method of introducing goodness is flawed here in
the general setting of dual martingale differences, since these differences
are no longer orthogonal projections, and as emphasized in \cite{HyMa}, we
cannot simply add back in bad cubes whenever we want telescoping
identities to hold - but these are needed in order to control the right hand
side of (\ref{e.Tgood'}). In fact, in the analysis of the form $\Theta
\left( f,g\right) $ above, it is necessary to have goodness for the
cubes $J$ and telescoping for the cubes $I$. On the other hand, in
the analysis of the form $\Theta ^{\ast }\left( f,g\right) $ above, it is
necessary to have just the opposite - namely goodness for the cubes $I$
and telescoping for the cubes $J$.

Thus, because in this unfortunate set of circumstances we can no longer `add
back in' bad cubes to achieve telescoping, we are prevented from introducing
goodness in the \emph{full} sum (\ref{ess symm})\ over all $I$ and $J$,
prior to splitting according to side lengths of $I$ and $J$. Thus the
infusion of goodness must come \emph{after} the splitting by side length,
but one must work much harder to introduce goodness directly into the form $%
\Theta \left( f,g\right) $ \emph{after} we have restricted the sum to
cubes $J$ that have smaller side length than $I$. This is accomplished
in the next subsubsection using the \emph{weaker form of NTV goodness}
introduced by Hyt\"{o}nen and Martikainen in \cite{HyMa} (that permits
certain additional pairs $\left( I,J\right) $ in the good forms where $\ell
\left( J\right) \leq 2^{-\mathbf{r}}\ell \left( I\right) $ and yet $J$ is $%
{bad}$ in the traditional sense), and that will prevail later in the
treatment of the far below forms $\mathsf{T}_{{far}{below}%
}^{1}\left( f,g\right) $, and of the local forms $\mathsf{B}_{\Subset _{%
\mathbf{r}}}^{A}\left( f,g\right) $ (see Subsection \ref{Sub wrapup}) where
the need for using the `body' of a cube will become apparent in dealing
with the stopping form, and also in the treatment of the functional energy
in Appendix .

\subsubsection{Weak goodness}

Let $\mathcal{D}$ and $\mathcal{G}$ be dyadic grids. It remains to estimate
the form $\Theta _{2}\left( f,g\right) $ which, following \cite{HyMa}, we
will split into a `bad' part and a `good' part. For this we introduce our
main definition associated with the above modification of the weak goodness
of Hyt\"{o}nen and Martikainen, namely the definition of the cube $%
R^{\maltese }$ in a grid $\mathcal{D}$, given an arbitrary cube $R\in 
\mathcal{P}$.

\begin{dfn}
\label{def sharp cross}Let $\mathcal{D}$ be a dyadic grid. Given $R\in 
\mathcal{P}$, let $R^{\maltese }$ be the smallest (if any such exist) $%
\mathcal{D}$-dyadic supercube $Q$ of $R$ such that $R$ is good inside 
\textbf{all} $\mathcal{D}$-dyadic supercubes $K$ of $Q$. Of course $%
R^{\maltese }$ will not exist if there is no $\mathcal{D}$-dyadic cube $%
Q $ containing $R$ in which $R$ is good. For cubes $R,Q\in \mathcal{P}$
let $\kappa \left( Q,R\right) =\log _{2}\frac{\ell \left( Q\right) }{\ell
\left( R\right) }$. For $R\in \mathcal{P}$ for which $R^{\maltese }$ exists,
let $\kappa \left( R\right) \equiv \kappa \left( R^{\maltese },R\right) $.
\end{dfn}

Note that we typically suppress the dependence of $R^{\maltese }$ on the
grid $\mathcal{D}$, since the grid is usually understood from context. If $%
R^{\maltese }$ exists, we thus have that $R$ is good inside all $\mathcal{D}$%
-dyadic supercubes $K$ of $R$ with $\ell \left( K\right) \geq \ell
\left( R^{\maltese }\right) $. Note in particular the monotonicity property
for $J^{\prime },J\in \mathcal{P}$:%
\begin{equation*}
J^{\prime }\subset J\Longrightarrow \left( J^{\prime }\right) ^{\maltese
}\subset J^{\maltese }.
\end{equation*}%
Here now is the decomposition:%
\begin{eqnarray*}
\Theta _{2}\left( f,g\right) &=&\sum_{I\in \mathcal{D}}\sum_{\substack{ J\in 
\mathcal{G}:\ J^{\maltese }\not\subsetneqq I\text{, }\ell \left( J\right)
\leq 2^{-\mathbf{r}}\ell \left( I\right)  \\ d\left( J,I\right) \leq 2\ell
\left( J\right) ^{\varepsilon }\ell \left( I\right) ^{1-\varepsilon }}}\int
\left( T_{\sigma }^{\alpha }\square _{I}^{\sigma ,\mathbf{b}}f\right)
\square _{J}^{\omega ,\mathbf{b}^{\ast }}gd\omega \\
&&+\sum_{I\in \mathcal{D}}\sum_{\substack{ J\in \mathcal{G}:\ J^{\maltese
}\subsetneqq I\text{, }\ell \left( J\right) \leq 2^{-\mathbf{r}}\ell \left(
I\right)  \\ d\left( J,I\right) \leq 2\ell \left( J\right) ^{\varepsilon
}\ell \left( I\right) ^{1-\varepsilon }}}\int \left( T_{\sigma }^{\alpha
}\square _{I}^{\sigma ,\mathbf{b}}f\right) \square _{J}^{\omega ,\mathbf{b}%
^{\ast }}gd\omega \\
&\equiv &\Theta _{2}^{{bad}}\left( f,g\right) +\Theta _{2}^{{%
good}}\left( f,g\right) \ ,
\end{eqnarray*}%
and where if $J^{\maltese }$ fails to exist, we assume by convention that $%
J^{\maltese }\not\subsetneqq I$, i.e. $J^{\maltese }$ is \emph{not} strictly
contained in $I$, so that the pair $\left( I,J\right) $ is then included in
the bad form $\Theta _{2}^{{bad}}\left( f,g\right) $. We will in
fact estimate a larger quantity corresponding to the bad form, namely%
\begin{equation}
\Theta _{2}^{{bad}\natural }\left( f,g\right) \equiv \sum_{I\in 
\mathcal{D}}\sum_{\substack{ J\in \mathcal{G}:\ J^{\maltese }\not\subsetneqq
I\text{, }\ell \left( J\right) \leq 2^{-\mathbf{r}}\ell \left( I\right)  \\ %
d\left( J,I\right) \leq 2\ell \left( J\right) ^{\varepsilon }\ell \left(
I\right) ^{1-\varepsilon }}}\left\vert \int \left( T_{\sigma }^{\alpha
}\square _{I}^{\sigma ,\mathbf{b}}f\right) \square _{J}^{\omega ,\mathbf{b}%
^{\ast }}gd\omega \right\vert  \label{Theta_2^bad sharp}
\end{equation}%
with absolute value signs \emph{inside} the sum.

\begin{rem}
We now make some general comments on where we now stand and where we are
going.

\begin{enumerate}
\item In the first sum $\Theta _{2}^{{bad}}\left( f,g\right) $
above, we are roughly keeping the pairs of cubes $\left( I,J\right) $
such that $J$ is ${bad}$ with respect to some `nearby' cube
having side length larger than that of $I$.

\item We have defined energy and dual energy conditions that are independent
of the testing families (because the definition of $\mathsf{E}\left(
J,\omega \right) =\mathbb{E}_{J}^{\omega ,x}\mathbb{E}_{J}^{\omega
,x^{\prime }}\left( \left\vert \frac{x-x^{\prime }}{\ell \left( J\right) }%
\right\vert ^{2}\right) $ does not involve pseudoprojections $\square _{J,%
\mathcal{D}}^{\omega ,\mathbf{b}^{\ast }}$), but the functional energy
condition defined below \emph{does} involve the dual martingale
pseudoprojections $\square _{J,\mathcal{D}}^{\omega ,\mathbf{b}^{\ast }}$.

\item Using the notion of weak goodness above, we will be able to eliminate
all pairs of cubes with $J$ bad in $I$, which then permits control of
the short range form in Section \ref{Sec disj form} and the neighbour form
in Section \ref{Sec Main below} provided $0<\varepsilon <\frac{1}{n+1-\alpha }$%
. Defining shifted coronas in terms of $J^{\maltese }$ will then allow
existing arguments to prove the Intertwining Proposition and obtain control
of the functional energy in Appendix , as well as permitting control of the
stopping form in Section \ref{Sec stop}, but all of this with some new
twists, for example the introduction of a top/down `indented corona' in the
analysis of the stopping form.

\item The nearby form $\Theta _{3}\left( f,g\right) $ is handled in Section %
\ref{Sec nearby} using the energy condition assumption along with the
original testing functions $b_{Q}^{{orig}}$ discarded during the
construction of the testing/accretive corona.
\end{enumerate}
\end{rem}

These remarks will become clear in this and the following sections. Recall
that we earlier defined in Definition \ref{good two grids}, the set $%
\mathcal{G}_{k-{good}}^{\mathcal{D}}=\mathcal{G}_{\left(
k,\varepsilon \right) -{good}}^{\mathcal{D}}$ to consist of those $%
J\in \mathcal{G}$ such that $J$ is $\varepsilon -{good}$ inside
every cube $K\in \mathcal{D}$ with $K\cap J\neq \emptyset $ that lies at
least $k$ levels `above' $J$, i.e. $\ell \left( K\right) \geq 2^{k}\ell
\left( J\right) $. We now define an analogous notion of $\mathcal{G}_{k-%
{bad}}^{\mathcal{D}}$.

\begin{dfn}
\label{def Gbad}Let $\varepsilon >0$. Define the set $\mathcal{G}_{k-%
{bad}}^{\mathcal{D}}=\mathcal{G}_{\left( k,\varepsilon \right) -%
{bad}}^{\mathcal{D}}$ to consist of all $J\in \mathcal{G}$ such that
there is a $\mathcal{D}$-cube $K$ with sidelength $\ell \left( K\right)
=2^{k}\ell \left( J\right) $ for which $J$ is $\varepsilon -{bad}$
with respect to $K$.
\end{dfn}

Note that for grids $\mathcal{D}$ and $\mathcal{G}$, the complement of $%
\mathcal{G}_{k-{good}}^{\mathcal{D}}$ is the union of $\mathcal{G}%
_{\ell -{bad}}^{\mathcal{D}}$ for $\ell \geq k$, i.e.%
\begin{equation*}
\mathcal{G\setminus G}_{k-{good}}^{\mathcal{D}}=\bigcup_{\ell \geq k}%
\mathcal{G}_{\ell -{bad}}^{\mathcal{D}}\ .
\end{equation*}%
Now assume $\varepsilon >0$. We then have the following important property,
namely for all cubes $R$, and all $k\geq \mathbf{r}$ (where the goodness
parameter $\mathbf{r}$ will be fixed given $\varepsilon >0$ in (\ref{choice
of r}) below): 
\begin{equation}
\#\left\{ Q:\kappa \left( Q,R\right) =k\text{ and }d\left( R,Q\right) \leq
2\ell \left( R\right) ^{\varepsilon }\ell \left( Q\right) ^{1-\varepsilon
}\right\} \lesssim 1.  \label{imp}
\end{equation}%
As in \cite{HyMa}, set%
\begin{equation*}
\mathcal{G}_{{bad},n}^{\mathcal{D}}\equiv \left\{ J\in \mathcal{G}:J%
\text{ is }\varepsilon -{bad}\ \text{with respect to some }K\in 
\mathcal{D}\text{ with }\ell \left( K\right) \geq n\right\} .
\end{equation*}%
We will now use the set equality%
\begin{eqnarray}
&&\left\{ J\in \mathcal{G}:\ J^{\maltese }\not\subset I,\text{ }\ell \left(
J\right) \leq 2^{-\mathbf{r}}\ell \left( I\right) ,\ d\left( J,I\right) \leq
2\ell \left( J\right) ^{\varepsilon }\ell \left( I\right) ^{1-\varepsilon
}\right\}  \label{set equ} \\
&=&\left\{ R\in \mathcal{G}_{{bad},\ell \left( Q\right) }^{\mathcal{D%
}}:\ \mathbf{r}\leq \kappa \left( Q,R\right) <\kappa \left( R\right) ,\
d\left( R,Q\right) \leq 2\ell \left( R\right) ^{\varepsilon }\ell \left(
Q\right) ^{1-\varepsilon }\right\} ,  \notag
\end{eqnarray}%
which the careful reader can prove by painstakingly verifying both
containments.

Assuming only that $\mathbf{b}$ is $2$-weakly $\mu $-controlled accretive, and following the
proof in \cite{HyMa}, we use (\ref{set equ}) to show that for any fixed
grids $\mathcal{D}$ and $\mathcal{G}$, and any bounded linear operator $%
T_{\sigma }^{\alpha }$ we have the following inequality for the form $\Theta
_{2}^{{bad}\natural ,{strict}}\left( f,g\right) $, defined
to be $\Theta _{2}^{{bad}\natural }\left( f,g\right) $ as in (\ref%
{Theta_2^bad sharp}) with the pairs $\left( I,J\right) $ removed when $%
J^{\maltese }=I$. We use $\varepsilon _{Q,R}=\pm 1$ to obtain 
\begin{eqnarray*}
&&\Theta _{2}^{{bad}\natural ,{strict}}\left( f,g\right)
=\sum_{Q\in \mathcal{D}}\sum_{\substack{ R\in \mathcal{G}_{{bad}%
,\ell \left( Q\right) }^{\mathcal{D}}:\ \mathbf{r}\leq \kappa \left(
Q,R\right) <\kappa \left( R\right)  \\ d\left( R,Q\right) \leq 2\ell \left(
R\right) ^{\varepsilon }\ell \left( Q\right) ^{1-\varepsilon }}}\left\vert
\left\langle T_{\sigma }^{\alpha }\left( \square _{Q,\mathcal{D}}^{\sigma ,%
\mathbf{b}}f\right) ,\square _{R,\mathcal{G}}^{\omega ,\mathbf{b}^{\ast
}}g\right\rangle \right\vert \\
&=&\sum_{Q\in \mathcal{D}}\sum_{\substack{ R\in \mathcal{G}_{{bad}%
,\ell \left( Q\right) }^{\mathcal{D}}:\ \mathbf{r}\leq \kappa \left(
Q,R\right) <\kappa \left( R\right)  \\ d\left( R,Q\right) \leq 2\ell \left(
R\right) ^{\varepsilon }\ell \left( Q\right) ^{1-\varepsilon }}}\varepsilon
_{Q,R}\left\langle T_{\sigma }^{\alpha }\left( \square _{Q,\mathcal{D}%
}^{\sigma ,\mathbf{b}}f\right) ,\square _{R,\mathcal{G}}^{\omega ,\mathbf{b}%
^{\ast }}g\right\rangle \\
&\leq &
\sum_{Q\in \mathcal{D}}\left\vert \left\langle T_{\sigma }^{\alpha
}\left( \square _{Q,\mathcal{D}}^{\sigma ,\mathbf{b}}f\right) ,\sum 
_{\substack{ R\in \mathcal{G}_{{bad},\ell \left( Q\right) }^{%
\mathcal{D}}:\ \mathbf{r}\leq \kappa \left( Q,R\right) <\kappa \left(
R\right)  \\ d\left( R,Q\right) \leq 2\ell \left( R\right) ^{\varepsilon
}\ell \left( Q\right) ^{1-\varepsilon }}}\varepsilon _{Q,R}\square _{R,%
\mathcal{G}}^{\omega ,\mathbf{b}^{\ast }}g\right\rangle \right\vert 
\end{eqnarray*}
\begin{eqnarray*}
&\leq &\mathfrak{N}_{T^{\alpha }}\sum_{Q\in \mathcal{D}}\left\Vert \square
_{Q,\mathcal{D}}^{\sigma ,\mathbf{b}}f\right\Vert _{L^{2}\left( \sigma
\right) }\left\Vert \sum_{\substack{ R\in \mathcal{G}_{{bad},\ell
\left( Q\right) }^{\mathcal{D}}:\ \mathbf{r}\leq \kappa \left( Q,R\right)
<\kappa \left( R\right)  \\ d\left( R,Q\right) \leq 2\ell \left( R\right)
^{\varepsilon }\ell \left( Q\right) ^{1-\varepsilon }}}\varepsilon
_{Q,R}\square _{R,\mathcal{G}}^{\omega ,\mathbf{b}^{\ast }}g\right\Vert
_{L^{2}\left( \omega \right) } \\
&\leq &
\mathfrak{N}_{T^{\alpha }}\sum_{Q\in \mathcal{D}}\left\Vert \square
_{Q,\mathcal{D}}^{\sigma ,\mathbf{b}}f\right\Vert _{L^{2}\left( \sigma
\right) }\sum_{k=\mathbf{r}}^{\infty }\left\Vert \sum_{\substack{ R\in 
\mathcal{G}_{{bad},\ell \left( Q\right) }^{\mathcal{D}}:k=\kappa
\left( Q,R\right) <\kappa \left( R\right)  \\ d\left( R,Q\right) \leq 2\ell
\left( R\right) ^{\varepsilon }\ell \left( Q\right) ^{1-\varepsilon }}}%
\varepsilon _{Q,R}\square _{R,\mathcal{G}}^{\omega ,\mathbf{b}^{\ast
}}g\right\Vert _{L^{2}\left( \omega \right) }\ ,
\end{eqnarray*}%
by Minkowski's inequality, and we continue with
\begin{eqnarray*}
&\leq &
2\mathfrak{N}_{T^{\alpha }}\sum_{k=\mathbf{r}}^{\infty }\left(
\sum_{Q\in \mathcal{D}}\left\Vert \square _{Q,\mathcal{D}}^{\sigma ,\mathbf{b%
}}f\right\Vert _{L^{2}\left( \sigma \right) }^{2}\right) ^{\frac{1}{2}%
} \cdot \\ 
&&
\left( \sum_{Q\in \mathcal{D}}\sum_{\substack{ R\in \mathcal{G}_{{%
bad},\ell \left( Q\right) }^{\mathcal{D}}:\ k=\kappa \left( Q,R\right)
<\kappa \left( R\right)  \\ d\left( R,Q\right) \leq 2\ell \left( R\right)
^{\varepsilon }\ell \left( Q\right) ^{1-\varepsilon }}}\left( \left\Vert
\square _{R,\mathcal{G}}^{\omega ,\mathbf{b}^{\ast }}g\right\Vert
_{L^{2}\left( \omega \right) }^{2}+\left\Vert \bigtriangledown _{R,\mathcal{G}%
}^{\omega }g\right\Vert _{L^{2}\left( \omega \right) }^{2}\right) \right) ^{%
\frac{1}{2}} \\
&\lesssim &\mathfrak{N}_{T^{\alpha }}\left\Vert f\right\Vert _{L^{2}\left(
\sigma \right) }\sum_{k=\mathbf{r}}^{\infty }\left( \sum_{R\in \mathcal{G}_{%
{bad},2^{k}\ell \left( R\right) }^{\mathcal{D}}}\left( \left\Vert
\square _{R,\mathcal{G}}^{\omega ,\mathbf{b}^{\ast }}g\right\Vert
_{L^{2}\left( \omega \right) }^{2}+\left\Vert \bigtriangledown _{R,\mathcal{G}%
}^{\omega }g\right\Vert _{L^{2}\left( \omega \right) }^{2}\right) \right) ^{%
\frac{1}{2}},
\end{eqnarray*}%
where $\bigtriangledown _{R,\mathcal{G}}^{\omega }$ denotes the `broken' Carleson
averaging operator in (\ref{Carleson avg op}) that depends on the grid $%
\mathcal{G}$, and

\begin{enumerate}
\item the penultimate inequality uses Cauchy-Schwarz in $Q$ and the weak
upper Riesz inequalities (\ref{UPPER RIESZ}) for $\displaystyle \sum_{\substack{ R\in 
\mathcal{G}_{{bad},\ell \left( Q\right) }^{\mathcal{D}}:\ k=\kappa
\left( Q,R\right) <\kappa \left( R\right)  \\ d\left( R,Q\right) \leq 2\ell
\left( R\right) ^{\varepsilon }\ell \left( Q\right) ^{1-\varepsilon }}}%
\varepsilon _{Q,R}\square _{R,\mathcal{G}}^{\omega ,\mathbf{b}^{\ast }}$,
once for the sum when $\varepsilon _{Q,R}=1$, and again for the sum when $%
\varepsilon _{Q,R}=-1$. However, we note that since the sum in $R$ is
pigeonholed by $k=\kappa \left( Q,R\right) $, the $R$'s are pairwise
disjoint cubes and the pseudoprojections $\square _{R,\mathcal{G}%
}^{\omega ,\mathbf{b}^{\ast }}g$ are pairwise orthogonal. Thus we could
instead apply Cauchy-Schwarz first in $R$, and then in $Q$ as was done in 
\cite{HyMa}, but we must still apply weak upper Riesz inequalities as above.

\item and the final inequality uses the frame inequality (\ref{FRAME})
together with (\ref{imp}), namely the fact that there are at most $C$
cubes $Q$ such that $\kappa \left( Q,R\right) \geq \mathbf{r}$ is fixed
and $d\left( R,Q\right) \leq 2\ell \left( R\right) ^{\varepsilon }\ell
\left( Q\right) ^{1-\varepsilon }$.
\end{enumerate}

Now it is easy to verify that we have the same inequality for the pairs $%
\left( J^{\maltese },J\right) $ that were removed, and then we take grid
expectations and use the probability estimate (\ref{main bad prob}) to
obtain for $\varepsilon ^{\prime }=\frac{1}{2}\varepsilon $ that $\boldsymbol{E}_{\Omega }^{\mathcal{D}}\left( \Theta _{2}^{{bad}%
\natural }\left( f,g\right) \right)$ is bounded by
\begin{eqnarray}
&&\label{HM bad}\\
&\leq &
\boldsymbol{E}_{\Omega }^{\mathcal{D}}\mathfrak{N}_{T^{\alpha
}}\left\Vert f\right\Vert _{L^{2}\left( \sigma \right) }\sum_{k=\mathbf{r}%
}^{\infty }\left( \sum_{R\in \mathcal{G}_{{bad},2^{k}\ell \left(
R\right) }^{\mathcal{D}}}\left( \left\Vert \square _{R,\mathcal{G}}^{\omega ,%
\mathbf{b}^{\ast }}g\right\Vert _{L^{2}\left( \omega \right)
}^{2}+\left\Vert \bigtriangledown _{R,\mathcal{G}}^{\omega }g\right\Vert _{L^{2}\left(
\omega \right) }^{2}\right) \right) ^{\frac{1}{2}}  \notag \\
&\leq &
\mathfrak{N}_{T^{\alpha }}\left\Vert f\right\Vert _{L^{2}\left(
\sigma \right) }\sum_{k=\mathbf{r}}^{\infty }\left( \boldsymbol{E}_{\Omega
}^{\mathcal{D}}\sum_{R\in \mathcal{G}_{{bad},2^{k}\ell \left(
R\right) }^{\mathcal{D}}}\left( \left\Vert \square _{R,\mathcal{G}}^{\omega ,%
\mathbf{b}^{\ast }}g\right\Vert _{L^{2}\left( \omega \right)
}^{2}+\left\Vert \bigtriangledown _{R,\mathcal{G}}^{\omega }\right\Vert _{L^{2}\left(
\omega \right) }^{2}\right) \right) ^{\frac{1}{2}}  \notag \\
&\lesssim &
2^{-\frac{1}{2}\varepsilon ^{\prime }\mathbf{r}}\mathfrak{N}%
_{T^{\alpha }}\left\Vert f\right\Vert _{L^{2}\left( \sigma \right) }\sum_{k=%
\mathbf{r}}^{\infty }\left( (1-(C_{1}2^{-\varepsilon k})^n)\left\Vert g\right\Vert
_{L^{2}\left( \omega \right) }^{2}\right) ^{\frac{1}{2}}\notag\\
&\leq&
C_{{good}} 2^{-\frac{1}{2}\varepsilon \mathbf{r}}\mathfrak{N}_{T^{\alpha }}\left\Vert
f\right\Vert _{L^{2}\left( \sigma \right) }\left\Vert g\right\Vert
_{L^{2}\left( \omega \right) }.  \notag
\end{eqnarray}%
Clearly we can now fix $\mathbf{r}$ sufficiently large depending on $%
\varepsilon >0$ so that%
\begin{equation}
C_{{good}}2^{-\frac{1}{2}\varepsilon \mathbf{r}}<\frac{1}{100},
\label{choice of r}
\end{equation}%
and then the final term above, namely $C_{{good}}2^{-\frac{1}{2}%
\varepsilon \mathbf{r}}\mathfrak{N}_{T^{\alpha }}\left\Vert f\right\Vert
_{L^{2}\left( \sigma \right) }\left\Vert g\right\Vert _{L^{2}\left( \omega
\right) }$, can be absorbed at the end of the proof in Subsection \ref{Sub
wrapup}. Note that (\ref{choice of r}) fixes our choice of the parameter $%
\mathbf{r}$ for any given $\varepsilon >0$. Later we will choose $%
0<\varepsilon <\frac{1}{2}\leq \frac{1}{n+1-\alpha }$. It is this type of weak
goodness that we will exploit in the local forms $\mathsf{B}_{\Subset _{%
\mathbf{r}}}^{A}\left( f,g\right) $ treated below in Section \ref{Sec Main
below}.

We are now left with the following `good' form to control:%
\begin{equation*}
\Theta _{2}^{{good}}\left( f,g\right) =\sum_{I\in \mathcal{D}}\sum 
_{\substack{ J^{\maltese }\subsetneqq I:\ \ell \left( J\right) \leq 2^{-%
\mathbf{r}}\ell \left( I\right)  \\ d\left( J,I\right) \leq 2\ell \left(
J\right) ^{\varepsilon }\ell \left( I\right) ^{1-\varepsilon }}}\int \left(
T_{\sigma }^{\alpha }\square _{I}^{\sigma ,\mathbf{b}}f\right) \square
_{J}^{\omega ,\mathbf{b}^{\ast }}gd\omega .
\end{equation*}%
The first thing we observe regarding this form is that the cubes $J$
which arise in the sum for $\Theta _{2}^{{good}}\left( f,g\right) $
must lie entirely inside $I$ since $J\subset J^{\maltese }\subsetneqq I$.
Then in the remainder of the paper, we proceed to analyze 
\begin{equation}
\Theta _{2}^{{good}}\left( f,g\right) =\sum_{I\in \mathcal{D}%
}\sum_{J^{\maltese }\subsetneqq I:\ \ell \left( J\right) \leq 2^{-\mathbf{r}%
}\ell \left( I\right) }\int \left( T_{\sigma }^{\alpha }\square _{I}^{\sigma
,\mathbf{b}}f\right) \square _{J}^{\omega ,\mathbf{b}^{\ast }}gd\omega ,
\label{def Theta 2 good}
\end{equation}%
in the same way we analyzed the below term $\mathsf{B}_{\Subset _{\mathbf{r}%
}}\left( f,g\right) $ in \cite{SaShUr6}; namely, by implementing the
canonical corona splitting and the decomposition into paraproduct, neighbour
and stopping forms, but now with an additional broken form. We have $\left(
\kappa ,\varepsilon \right) $-goodness available for all the cubes $J\in 
\mathcal{G}$ arising in the form $\Theta _{2}^{{good}}\left(
f,g\right) $, and moreover, the cubes $I\in \mathcal{D}$ arising in the
form $\Theta _{2}^{{good}}\left( f,g\right) $ for a fixed $J$ are
tree-connected, so that telescoping identities hold for these cubes $I$.
This will prove decisive in the following three sections of the paper.

The forms $\Theta _{1}\left( f,g\right) $ and $\Theta _{3}\left( f,g\right) $
are analogous to the disjoint and nearby forms $\mathsf{B}_{\cap }\left(
f,g\right) $ and $\mathsf{B}_{/}\left( f,g\right) $ in \cite{SaShUr6}
respectively. In the next\ two sections, we control the disjoint form $%
\Theta _{1}\left( f,g\right) $ in essentially the same way that the disjoint
form $\mathsf{B}_{\cap }\left( f,g\right) $ was treated in \cite{SaShUr6}
and in earlier papers of many authors beginning with Nazarov, Treil and
Volberg (see e.g. \cite{Vol}), and we control the nearby form $\Theta
_{3}\left( f,g\right) $ using the probabilistic surgery of Hyt\"{o}nen and
Martikainen building on that of NTV, together with a new deterministic
surgery involving the energy condition and the original testing functions.
But first we recall, in the following subsection, the characterization of
boundedness of one-dimensional forms supported on disjoint cubes \cite%
{Hyt2}.

\section{Disjoint form\label{Sec disj form}}

Here we control the disjoint form $\Theta _{1}\left( f,g\right) $ by further
decomposing it as follows:%
\begin{eqnarray*}
\Theta _{1}\left( f,g\right)  &=&\sum_{I\in \mathcal{D}}\sum_{\substack{ %
J\in \mathcal{G}:\ \ell \left( J\right) \leq \ell \left( I\right)  \\ %
d\left( J,I\right) >2\ell \left( J\right) ^{\varepsilon }\ell \left(
I\right) ^{1-\varepsilon }}}\int \left( T_{\sigma }\square _{I}^{\sigma ,%
\mathbf{b}}f\right) \square _{J}^{\omega ,\mathbf{b}^{\ast }}gd\omega  
\end{eqnarray*}
which can be rewritten as
\begin{eqnarray*}
&&
\!\!\!\!\!\!\!\!\!\!\!\!\!\!\!\sum_{I\in \mathcal{D}}\!\!\left\{ \sum_{\substack{ J\in \mathcal{G}:\ \ell
\left( J\right) \leq \ell \left( I\right)  \\ d\left( J,I\right) >\max (\ell
\left( I\right),2\ell(J)^\varepsilon\ell(I)^{1-\varepsilon}) }}
\!\!+\!\!
\sum_{\substack{ J\in \mathcal{G}:\ \ell \left( J\right) \leq \ell \left( I\right) \\ \ell \left( I\right) \geq d\left( J,I\right) >2\ell \left(
J\right) ^{\varepsilon }\ell \left( I\right) ^{1-\varepsilon }}}\!\!\right\}\!\! \int\!\!\! \left( T_{\sigma }\square _{I}^{\sigma ,%
\mathbf{b}}f\right) \square _{J}^{\omega ,\mathbf{b}^{\ast }}gd\omega  \\
&\equiv &
\Theta _{1}^{{long}}\left( f,g\right) +\Theta _{1}^{%
{short}}\left( f,g\right) ,
\end{eqnarray*}%
where $\Theta _{1}^{{long}}\left( f,g\right) $ is a `long range'
form in which $J$ is far from $I$, and where $\Theta _{1}^{{short}}\left( f,g\right) $ is a short range form. It should be noted that the goodness plays no role in treating the disjoint form.

\subsection{Long range form}

\begin{lem}
\label{delta long}We have%
\begin{equation*}
\sum_{I\in \mathcal{D}}\sum_{\substack{ J\in \mathcal{G}:\ \ell \left(
J\right) \leq \ell \left( I\right)  \\ d\left( J,I\right) >\ell \left(
I\right) }}\left\vert \int \left( T_{\sigma }\square _{I}^{\sigma ,\mathbf{b}%
}f\right) \square _{J}^{\omega ,\mathbf{b}^{\ast }}gd\omega \right\vert
\lesssim \sqrt{\mathfrak{A}_{2}^{\alpha }}\left\Vert f\right\Vert _{L^{2}\left( \sigma
\right) }\left\Vert g\right\Vert _{L^{2}\left( \omega \right) }
\end{equation*}
\end{lem}

\begin{proof}
Since $J$ and $I$ are separated by at least $\max \left\{ \ell \left(
J\right) ,\ell \left( I\right) \right\} $, we have the inequality%
\begin{eqnarray*}
\mathrm{P}^{\alpha }\left( J,\left\vert \square _{I}^{\sigma ,\mathbf{b}%
}f\right\vert \sigma \right) 
&\approx&
\int_{I}\frac{\ell \left( J\right) }{\left\vert y-c_{J}\right\vert ^{n+1-\alpha }}\left\vert \square _{I}^{\sigma
,\mathbf{b}}f\left( y\right) \right\vert d\sigma \left( y\right)\\ &\lesssim&
\left\Vert \square _{I}^{\sigma ,\mathbf{b}}f\right\Vert _{L^{2}\left(
\sigma \right) }\frac{\ell \left( J\right) \sqrt{\left\vert I\right\vert
_{\sigma }}}{d\left( I,J\right) ^{n+1-\alpha }},
\end{eqnarray*}
since $\displaystyle \int_{I}\left\vert \square _{I}^{\sigma ,\mathbf{b}}f\left( y\right)
\right\vert d\sigma \left( y\right) \leq \left\Vert \square _{I}^{\sigma ,\mathbf{b}}f\right\Vert _{L^{2}\left( \sigma \right) }\sqrt{\left\vert
I\right\vert _{\sigma }}$. Thus if $A\left( f,g\right) $ denotes the left
hand side of the conclusion of Lemma \ref{delta long}, we have using first the Energy Lemma,
\begin{eqnarray*}
&\!\!\!\!\!\!\!\!&\!\!\!\!\!\!\!\hspace{-2cm}\!\!\!\!\!\!\! A\left( f,g\right)  
\lesssim 
\sum_{I\in \mathcal{D}}\sum_{\substack{J\;:\;\ell
\left( J\right) \leq \ell \left( I\right) \\ d\left( I,J\right) \geq \ell
\left( I\right) }}\left\Vert \square _{I}^{\sigma ,\mathbf{b}}f\right\Vert
_{L^{2}\left( \sigma \right) }\left\Vert \square _{J}^{\omega ,\mathbf{b}%
^{\ast }}g\right\Vert _{L^{2}\left( \omega \right) }\frac{\ell \left( J\right) }{d\left(
I,J\right) ^{n+1-\alpha }}\sqrt{\left\vert I\right\vert _{\sigma }}\sqrt{%
\left\vert J\right\vert _{\omega }} \\
&\!\!\!\!\!\!\!\!\!\!\!\!\!\!\!\!\!\!\!\!\!\!\!\!\!\!\!\!\!\!\!\!\equiv &
\!\!\!\!\!\!\!\!\!\!\!\!\!\!\!\!\!\!\!\!\!\!\!\!\!
\sum_{\left( I,J\right) \in \mathcal{P}}\left\Vert \square
_{I}^{\sigma ,\mathbf{b}}f\right\Vert _{L^{2}\left( \sigma \right)
}\left\Vert \square _{J}^{\omega ,\mathbf{b}^{\ast }}g\right\Vert
_{L^{2}\left( \omega \right) }A\left( I,J\right) ; \\
\text{with }A\left( I,J\right)  &\equiv &\frac{\ell \left( J\right) }{%
d\left( I,J\right) ^{n+1-\alpha }}\sqrt{\left\vert I\right\vert _{\sigma }}%
\sqrt{\left\vert J\right\vert _{\omega }}; \\
\text{ and }\mathcal{P} &\equiv &\left\{ \left( I,J\right) \in \mathcal{D}%
\times \mathcal{G}:\ell \left( J\right) \leq \ell \left( I\right) \text{ and 
}d\left( I,J\right) \geq \ell \left( I\right) \right\} .
\end{eqnarray*}%
Now let $\mathcal{D}_{N}\equiv \left\{ K\in \mathcal{D}:\ell \left( K\right)
=2^{N}\right\} $ for each $N\in \mathbb{Z}$. For $N\in \mathbb{Z}$ and $s\in 
\mathbb{Z}_{+}$, we further decompose $A\left( f,g\right) $ by pigeonholing
the sidelengths of $I$ and $J$ by $2^{N}$ and $2^{N-s}$ respectively: 
\begin{eqnarray*}
A\left( f,g\right)  &=&\sum_{s=0}^{\infty }\sum_{N\in \mathbb{Z}%
}A_{N}^{s}\left( f,g\right) ; \\
A_{N}^{s}\left( f,g\right)  &\equiv &\sum_{\left( I,J\right) \in \mathcal{P}%
_{N}^{s}}\left\Vert \square _{I}^{\sigma ,\mathbf{b}}f\right\Vert
_{L^{2}\left( \sigma \right) }\left\Vert \square _{J}^{\omega ,\mathbf{b}%
^{\ast }}g\right\Vert _{L^{2}\left( \omega \right) }A\left( I,J\right)  \\
\text{where }\mathcal{P}_{N}^{s} &\equiv &\left\{ \left( I,J\right) \in 
\mathcal{D}_{N}\times \mathcal{G}_{N-s}:d\left( I,J\right) \geq \ell \left(
I\right) \right\} .
\end{eqnarray*}

Now let $\mathsf{P}_{M}^{\sigma }=\sum\limits_{K\in \mathcal{D}_{M}}\square
_{K}^{\sigma ,\mathbf{b}}$ denote the dual martingale pseudoprojection onto Span$\left\{ \square _{K}^{\sigma ,\mathbf{b}}\right\} _{K\in 
\mathcal{D}_{M}}$. Since the cubes $K$ in $\mathcal{D}_M$ are pairwise disjoint, the pseudoprojections $\square _{K}^{\sigma ,\mathbf{b}}$ are mutually orthogonal, which means that $\left\Vert\mathsf{P}_{M}^{\sigma }f\right\Vert^2_{L^2(\sigma)}=\sum\limits_{K\in \mathcal{D}_{M}}\left\Vert\square
_{K}^{\sigma ,\mathbf{b}}f\right\Vert^2_{L^2(\sigma)}$. We claim that%
\begin{equation}
\left\vert A_{N}^{s}\left( f,g\right) \right\vert \leq C2^{-s}\sqrt{%
\mathfrak{A}_{2}^{\alpha }}\left\Vert \mathsf{P}_{N}^{\sigma }f\right\Vert
_{L^{2}\left( \sigma \right) }^{\bigstar }\left\Vert \mathsf{P}%
_{N-s}^{\omega }g\right\Vert _{L^{2}\left( \omega \right) }^{\bigstar },\ \
\ \ \ \text{for }s\geq 0\text{ and }N\in \mathbb{Z}.  \label{AsN}
\end{equation}%
With this proved, we can then obtain%
\begin{eqnarray*}
A\left( f,g\right) &=&\sum_{s=0}^{\infty }\sum_{N\in \mathbb{Z}%
}A_{N}^{s}\left( f,g\right) =\sum_{s=0}^{\infty }\sum_{N\in \mathbb{Z}%
}A_{N}^{s}\left( f,g\right) \\
&\leq&
C\sqrt{\mathfrak{A}_{2}^{\alpha }}%
\sum_{s=0}^{\infty }2^{-s}\sum_{N\in \mathbb{Z}}\left\Vert \mathsf{P}%
_{N}^{\sigma }f\right\Vert^\bigstar _{L^{2}\left( \sigma \right) }\left\Vert \mathsf{P}_{N-s}^{\omega }g\right\Vert^\bigstar _{L^{2}\left( \omega \right) } \\
&\leq &
C\sqrt{\mathfrak{A}_{2}^{\alpha }}\sum_{s=0}^{\infty }2^{-s}\left(
\sum_{N\in \mathbb{Z}}\left\Vert \mathsf{P}_{N}^{\sigma }f\right\Vert
_{L^{2}\left( \sigma \right) }^{\bigstar 2}\right) ^{\frac{1}{2}}\left(
\sum_{N\in \mathbb{Z}}\left\Vert \mathsf{P}_{N-s}^{\omega }g\right\Vert
_{L^{2}\left( \omega \right) }^{\bigstar 2}\right) ^{\frac{1}{2}} \\
&\leq &C\sqrt{\mathfrak{A}_{2}^{\alpha }}\sum_{s=0}^{\infty
}2^{-s}\left\Vert f\right\Vert _{L^{2}\left( \sigma \right) }\left\Vert
g\right\Vert _{L^{2}\left( \omega \right) }=C\sqrt{\mathfrak{A}_{2}^{\alpha }}%
\left\Vert f\right\Vert _{L^{2}\left( \sigma \right) }\left\Vert
g\right\Vert _{L^{2}\left( \omega \right) }.
\end{eqnarray*}

To prove (\ref{AsN}), we pigeonhole the distance between $I$ and $J$:%
\begin{eqnarray*}
A_{N}^{s}\left( f,g\right) &=&\sum\limits_{\ell =0}^{\infty }A_{N,\ell
}^{s}\left( f,g\right) ; \\
A_{N,\ell }^{s}\left( f,g\right) &\equiv &\sum_{\left( I,J\right) \in 
\mathcal{P}_{N,\ell }^{s}}\left\Vert \square _{I}^{\sigma ,\mathbf{b}%
}f\right\Vert _{L^{2}\left( \sigma \right) }\left\Vert \square _{J}^{\omega ,%
\mathbf{b}^{\ast }}g\right\Vert _{L^{2}\left( \omega \right) }A\left(
I,J\right) \\
\text{where }\mathcal{P}_{N,\ell }^{s} &\equiv &\left\{ \left( I,J\right)
\in \mathcal{D}_{N}\times \mathcal{G}_{N-s}:d\left( I,J\right) \approx
2^{N+\ell }\right\} .
\end{eqnarray*}%
If we define $\mathcal{H}\left( A_{N,\ell }^{s}\right) $ to be the bilinear
form on $\ell ^{2}\times \ell ^{2}$ with matrix $\left[ A\left( I,J\right) %
\right] _{\left( I,J\right) \in \mathcal{P}_{N,\ell }^{s}}$, then it remains
to show that the norm $\left\Vert \mathcal{H}\left( A_{N,\ell }^{s}\right)
\right\Vert _{\ell ^{2}\rightarrow \ell ^{2}}$ of $\mathcal{H}\left(
A_{N,\ell }^{s}\right) $ on the sequence space $\ell ^{2}$ is bounded by $%
C2^{-s-\ell }\sqrt{\mathfrak{A}_{2}^{\alpha }}$. In turn, this is equivalent
to showing that the norm $\left\Vert \mathcal{H}\left( B_{N,\ell
}^{s}\right) \right\Vert _{\ell ^{2}\rightarrow \ell ^{2}}$ of the bilinear
form $\mathcal{H}\left( B_{N,\ell }^{s}\right) \equiv \mathcal{H}\left(
A_{N,\ell }^{s}\right) ^{{tr}}\mathcal{H}\left( A_{N,\ell
}^{s}\right) $ on the sequence space $\ell ^{2}$ is bounded by $%
C^{2}2^{-2s-2\ell }\mathfrak{A}_{2}^{\alpha }$. Here $\mathcal{H}\left(
B_{N,\ell }^{s}\right) $ is the quadratic form with matrix kernel $\left[
B_{N,\ell }^{s}\left( J,J^{\prime }\right) \right] _{J,J^{\prime }\in 
\mathcal{D}_{N-s}}$ having entries:%
\begin{equation*}
B_{N,\ell }^{s}\left( J,J^{\prime }\right) \equiv \sum_{I\in \mathcal{D}%
_{N}:\ d\left( I,J\right) \approx d\left( I,J^{\prime }\right) \approx
2^{N+\ell }}A\left( I,J\right) A\left( I,J^{\prime }\right) ,\ \ \ \ \ \text{%
for }J,J^{\prime }\in \mathcal{G}_{N-s}.
\end{equation*}

We are reduced to showing the bilinear form inequality,%
\begin{equation*}
\left\Vert \mathcal{H}\left( B_{N,\ell }^{s}\right) \right\Vert _{\ell
^{2}\rightarrow \ell ^{2}}\leq C2^{-2s-2\ell }\mathfrak{A}_{2}^{\alpha }\ \
\ \ \text{for }s\geq 0\text{, }\ell \geq 0\text{ and }N\in \mathbb{Z}.
\end{equation*}%
We begin by computing $B_{N,\ell }^{s}\left( J,J^{\prime }\right)$:
\begin{eqnarray*}
B_{N,\ell }^{s}\left( J,J^{\prime }\right) 
&=&\!\!\!\!\!\!\!\!\!\!\!
\sum_{\substack{ I\in 
\mathcal{D}_{N} \\ d\left( I,J\right) \approx d\left( I,J^{\prime }\right)\approx 2^{N+\ell }}}\!\!\!\!\!\!\!\!\!\!\!\!
\frac{\ell \left( J\right) }{d\left( I,J\right)
^{n+1-\alpha }}\sqrt{\left\vert I\right\vert _{\sigma }}\sqrt{\left\vert
J\right\vert _{\omega }}\frac{\ell \left( J^{\prime }\right) }{d\left(
I,J^{\prime }\right) ^{n+1-\alpha }}\sqrt{\left\vert I\right\vert _{\sigma }}%
\sqrt{\left\vert J^{\prime }\right\vert _{\omega }} \\
&=&
\!\!\!\!\!\!\!\!\!\!\!\!\!\sum_{\substack{ I\in \mathcal{D}_{N} \\ d\left( I,J\right)
\approx d\left( I,J^{\prime }\right) \approx 2^{N+\ell }}}\!\!\!\!\!\!\!\!\!\!\!\!\!\!\frac{|I|_\sigma}{d\left( I,J\right) ^{n+1-\alpha }d\left(
I,J^{\prime }\right) ^{n+1-\alpha }}\cdot \ell \left( J\right) \ell
\left( J^{\prime }\right) \sqrt{\left\vert J\right\vert _{\omega }}\sqrt{%
\left\vert J^{\prime }\right\vert _{\omega }}.
\end{eqnarray*}%
Now we show that%
\begin{equation}
\left\Vert B_{N,\ell }^{s}\right\Vert _{\ell ^{2}\rightarrow \ell
^{2}}\lesssim 2^{-2s-2\ell }\mathfrak{A}_{2}^{\alpha }\ ,  \label{Schur s}
\end{equation}%
by applying the proof of Schur's lemma. Fix $\ell \geq 0$ and $s\geq 0$.
Choose the Schur function $\beta \left( K\right) =\frac{1}{\sqrt{\left\vert
K\right\vert _{\omega }}}$. Fix $J\in \mathcal{D}_{N-s}$. We now group those 
$I\in \mathcal{D}_{N}$ with $d\left( I,J\right) \approx 2^{N+\ell }$ into
finitely many groups $G_{1},...G_{C}$ for which the union of the $I$ in each
group is contained in a cube of side length roughly $\frac{1}{100}2^{N+\ell
}$ , and we set $I_{k}^{\ast }\equiv \bigcup\limits_{I\in G_{k}}I$ for $%
1\leq k\leq C$ (note that $I_{k}^{\ast }$ is not a cube). We then have%
\begin{eqnarray*}
&&\sum_{J^{\prime }\in \mathcal{G}_{N-s}}\frac{\beta \left( J\right) }{\beta
\left( J^{\prime }\right) }B_{N,\ell }^{s}\left( J,J^{\prime }\right)  \\
&=&\sum_{\substack{ J^{\prime }\in \mathcal{G}_{N-s} \\ d\left( J^{\prime
},J\right) \leq \frac{1}{100}2^{N+\ell +2}}}\frac{\beta \left( J\right) }{%
\beta \left( J^{\prime }\right) }B_{N,\ell }^{s}\left( J,J^{\prime }\right)
+
\sum_{\substack{ J^{\prime }\in \mathcal{G}_{N-s} \\ d\left( J^{\prime
},J\right) >\frac{1}{100}2^{N+\ell +2}}}\frac{\beta \left( J\right) }{\beta
\left( J^{\prime }\right) }B_{N,\ell }^{s}\left( J,J^{\prime }\right)  \\
&\equiv& A+B,
\end{eqnarray*}%
where%
\begin{eqnarray*}
A
&\lesssim &
\sum_{\substack{ J^{\prime }\in \mathcal{G}_{N-s} \\ d\left(
J,J^{\prime }\right) \leq \frac{1}{100}2^{N+\ell +2}}}\left\{ \sum
_{\substack{ I\in \mathcal{D}_{N} \\ d\left( I,J\right) \approx 2^{N+\ell }}}%
\left\vert I\right\vert _{\sigma }\right\} \ \frac{2^{2\left( N-s\right) }}{%
2^{2\left( \ell +N\right) \left( n+1-\alpha \right) }}\left\vert J^{\prime
}\right\vert _{\omega } \\
&=&
\sum_{\substack{ J^{\prime }\in \mathcal{G}_{N-s} \\ d\left( J,J^{\prime
}\right) \leq \frac{1}{100}2^{N+\ell +2}}}\left\{ \sum_{k=1}^{C}\left\vert
I_{k}^{\ast }\right\vert _{\sigma }\right\} \ \frac{2^{2\left( N-s\right) }}{%
2^{2\left( \ell +N\right) \left( n+1-\alpha \right) }}\left\vert J^{\prime
}\right\vert _{\omega }\\
&=&
\frac{2^{2\left( N-s\right) }}{2^{2\left( \ell
+N\right) \left( n+1-\alpha \right) }}\sum_{k=1}^{C}\sum_{\substack{ %
J^{\prime }\in \mathcal{G}_{N-s} \\ d\left( J,J^{\prime }\right) \leq \frac{1%
}{100}2^{N+\ell +2}}}\left\vert I_{k}^{\ast }\right\vert _{\sigma }\
\left\vert J^{\prime }\right\vert _{\omega } \\
&\lesssim &
2^{-2s-2\ell }\sum_{k=1}^{C}\frac{\left\vert I_{k}^{\ast
}\right\vert _{\sigma }}{2^{\left( \ell +N\right) \left( n-\alpha \right) }}%
\frac{\left\vert \frac{1}{100}2^{N+\ell +2}J\right\vert _{\omega }}{%
2^{\left( \ell +N\right) \left( n-\alpha \right) }}\lesssim 2^{-2s-2\ell }%
\mathfrak{A}_{2}^{\alpha },
\end{eqnarray*}%
since $I_{k}^{\ast }$ is contained in a cube $\tilde{I_k^{\ast}}$ such that $|I_{k}^{\ast }|\approx|\tilde{I_k^\ast}|$, with an implied constant depending only on dimension, and $\tilde{I_{k}^{\ast}}$, $\frac{1}{100}2^{N+\ell +2}J$ are well
separated. If we let $Q_{k}$ be the smallest cube containing the set%
\begin{equation*}
E_{k}\equiv \bigcup\limits_{\substack{ J^{\prime }\in \mathcal{D}_{N-s}:\
d\left( I_{k}^{\ast },J^{\prime }\right) \approx 2^{N+\ell } \\ d\left(
J,J^{\prime }\right) >\frac{1}{100}2^{N+\ell +2}}}J^{\prime }
\end{equation*}%
we then have

\begin{eqnarray*}
B &\lesssim &\sum_{\substack{ J^{\prime }\in \mathcal{D}_{N-s} \\ d\left(
J,J^{\prime }\right) >\frac{1}{100}2^{N+\ell +2}}}\left\{ \sum_{\substack{ %
I\in \mathcal{D}_{N} \\ d\left( I,J^{\prime }\right) \approx d\left(
I,J\right) \approx 2^{N+\ell }}}\left\vert I\right\vert _{\sigma }\right\} \ 
\frac{2^{2\left( N-s\right) }}{2^{2\left( \ell +N\right) \left( n+1-\alpha
\right) }}\left\vert J^{\prime }\right\vert _{\omega } \\
&\lesssim &\sum_{\substack{ J^{\prime }\in \mathcal{D}_{N-s} \\ d\left(
J,J^{\prime }\right) >\frac{1}{100}2^{N+\ell +2}}}\left\{ \sum_{k:\ d\left(
I_{k}^{\ast },J^{\prime }\right) \approx 2^{N+\ell }}\left\vert I_{k}^{\ast
}\right\vert _{\sigma }\right\} \ \frac{2^{2\left( N-s\right) }}{2^{2\left(
\ell +N\right) \left( n+1-\alpha \right) }}\left\vert J^{\prime }\right\vert
_{\omega } \\
&\lesssim &\frac{2^{2\left( N-s\right) }}{2^{2\left( \ell +N\right) \left(
n+1-\alpha \right) }}\sum_{k=1}^{C}\left\vert I_{k}^{\ast }\right\vert
_{\sigma }\left\vert E_{k}\right\vert _{\omega } \\
&\lesssim &2^{-2s-2\ell }\sum_{k=1}^{C}\frac{\left\vert I_{k}^{\ast
}\right\vert _{\sigma }}{2^{\left( \ell +N\right) \left( n-\alpha \right) }}%
\frac{\left\vert Q_{k}\right\vert _{\omega }}{2^{\left( \ell +N\right)
\left( n-\alpha \right) }}\lesssim 2^{-2s-2\ell }\mathfrak{A}_{2}^{\alpha },
\end{eqnarray*}%
since $I_{k}^{\ast }$ is contained in a cube $\tilde{I_k^{\ast}}$ such that $|I_{k}^{\ast }|\approx|\tilde{I_k^\ast}|$, with an implied constant depending only on dimension, and $\tilde{I_{k}^{\ast}}$, $\frac{1}{100}2^{N+\ell +2}J$ are well
separated.
Thus we can now apply Schur's argument with $\displaystyle \sum_{J}\left( a_{J}\right)
^{2}=\sum_{J^{\prime }}\left( b_{J^{\prime }}\right) ^{2}=1$ to obtain%
\begin{eqnarray*}
\sum_{J,J^{\prime }\in \mathcal{G}_{N-s}}a_{J}b_{J^{\prime }}B_{N,\ell
}^{s}\left( J,J^{\prime }\right) 
\!\!\!\!&=&\!\!\!\!
\sum_{J,J^{\prime }\in \mathcal{G}_{N-s}}a_{J}\beta \left( J\right)
b_{J^{\prime }}\beta \left( J^{\prime }\right) \frac{B_{N,\ell }^{s}\left(
J,J^{\prime }\right) }{\beta \left( J\right) \beta \left( J^{\prime }\right)
} \\
&\leq &
\!\!\!\!
\sum_{J}\left( a_{J}\beta \left( J\right) \right) ^{2}\sum_{J^{\prime
}}\frac{B_{N,\ell }^{s}\left( J,J^{\prime }\right) }{\beta \left( J\right)
\beta \left( J^{\prime }\right) }+\sum_{J^{\prime }}\left( b_{J^{\prime
}}\beta \left( J^{\prime }\right) \right) ^{2}\sum_{J}\frac{B_{N,\ell
}^{s}\left( J,J^{\prime }\right) }{\beta \left( J\right) \beta \left(
J^{\prime }\right) } \\
&=&\!\!\!\!
\sum_{J}\left( a_{J}\right) ^{2}\left\{ \sum_{J^{\prime }}\frac{\beta
\left( J\right) }{\beta \left( J^{\prime }\right) }B_{N,\ell }^{s}\left(
J,J^{\prime }\right) \right\}\! +\!  \sum_{J^{\prime }}\left( b_{J^{\prime
}}\right) ^{2}\left\{ \sum_{J}\frac{\beta \left( J^{\prime }\right) }{\beta
\left( J\right) }B_{N,\ell }^{s}\left( J,J^{\prime }\right) \right\} \\
&\lesssim &
\!\!\!\!
2^{-2s-2\ell }A_{2}^{\alpha }\left( \sum_{J}\left( a_{J}\right)
^{2}+\sum_{J^{\prime }}\left( b_{J^{\prime }}\right) ^{2}\right)
=2^{1-2s-2\ell }\mathfrak{A}_{2}^{\alpha }.
\end{eqnarray*}%
This completes the proof of (\ref{Schur s}). We can now sum in $\ell $ to
get (\ref{AsN}) and we are done. This completes our proof of the long range
estimate%
\begin{equation*}
\mathcal{A}\left( f,g\right) \lesssim \sqrt{A_{2}^{\alpha }}\left\Vert
f\right\Vert _{L^{2}\left( \sigma \right) }\left\Vert g\right\Vert
_{L^{2}\left( \omega \right) }\ .
\end{equation*}
\end{proof}

\subsection{Short range form}

The form $\Theta _{1}^{{short}}\left( f,g\right) $ is
handled by the following lemma.
\begin{lem}
\label{delta short}We have%
\begin{equation*}
\sum_{I\in \mathcal{D}}\sum_{\substack{ J\in \mathcal{G}:\ \ell \left(
J\right) \leq 2^{-\mathbf{\rho }}\ell \left( I\right)  \\ \ell \left(
I\right) \geq d\left( J,I\right) >2\ell \left( J\right)
^{\varepsilon }\ell \left( I\right) ^{1-\varepsilon }}}\left\vert \int
\left( T_{\sigma }\square _{I}^{\sigma ,\mathbf{b}}f\right) \square
_{J}^{\omega ,\mathbf{b}^{\ast }}gd\omega \right\vert \lesssim \sqrt{\mathfrak{A}_{2}^{\alpha }}\left\Vert f\right\Vert _{L^{2}\left( \sigma \right)
}\left\Vert g\right\Vert _{L^{2}\left( \omega \right) }
\end{equation*}
\end{lem}

\begin{proof}
The pairs $(I,J)$ that occur in the sum above satisfy $J\subset 4I\backslash I$, so we consider
\begin{equation*}
\mathcal{P}\equiv \left\{ \left( I,J\right)\! \in\! \mathcal{D}\!\times\! \mathcal{G}\!:\!\ell\left( J\right) \!\leq\! 2^{-\mathbf{\rho }}\ell \left( I\right), \ell
\left( I\right)\!\geq\! d\left( J,I\right) \!>\!2\ell \left(J\right)
^{\varepsilon }\ell \left( I\right) ^{1-\varepsilon }, J\!\subset\!
4I\backslash I\!\right\} 
\end{equation*}%
For $\left( I,J\right) \in \mathcal{P}$, the `pivotal' estimate from the
Energy Lemma \ref{ener} gives
\begin{equation*}
\left\vert \left\langle T_{\sigma }^{\alpha }\left( \square _{I}^{\sigma ,%
\mathbf{b}}f\right) ,\square _{J}^{\omega ,\mathbf{b}^{\ast }}g\right\rangle
_{\omega }\right\vert \lesssim \left\Vert \square _{J}^{\omega ,\mathbf{b}%
^{\ast }}g\right\Vert _{L^{2}\left( \omega \right) }\mathrm{P}^{\alpha
}\left( J,\left\vert \bigtriangleup _{I}^{\sigma }f\right\vert \sigma
\right) \sqrt{\left\vert J\right\vert _{\omega }}\,.
\end{equation*}
Now we pigeonhole the lengths of $I$ and $J$ and the distance between them
by defining
\begin{equation*}
\mathcal{P}_{N,d}^{s}\equiv \left\{ \left( I,J\right) \in \mathcal{P}:\ell
\left( I\right) =2^{N},\ \ell \left( J\right) =2^{N-s},\ 2^{d-1}\leq d\left(
I,J\right) \leq 2^{d},\text{ }J\subset 4I\backslash I\right\}.
\end{equation*}
Note that the closest a cube $J$ can come to $I$ is determined by: 
\begin{eqnarray*}
&&
2^{d}\geq 2\ell \left( I\right) ^{1-\varepsilon }\ell \left(
J\right) ^{\varepsilon }=2^{1+N\left( 1-\varepsilon \right)
}2^{\left( N-s\right) \varepsilon }=2^{1+N-\varepsilon s}; \\
&&
\text{which implies }N-\varepsilon s+1\leq d\leq N.
\end{eqnarray*}%
Thus we have%
\begin{eqnarray*}
&&
\sum\limits_{\left( I,J\right) \in \mathcal{P}}\left\vert \left\langle
T_{\sigma }^{\alpha }\left( \square _{I}^{\sigma ,\mathbf{b}}f\right)
,\square _{J}^{\omega ,\mathbf{b}^{\ast }}g\right\rangle _{\omega
}\right\vert \\
&\lesssim&
\sum\limits_{\left( I,J\right) \in \mathcal{P}
}\left\Vert \square _{J}^{\omega ,\mathbf{b}^{\ast }}g\right\Vert
_{L^{2}\left( \omega \right) }\mathrm{P}^{\alpha }\left( J,\left\vert
\square _{I}^{\sigma ,\mathbf{b}}f\right\vert \sigma \right) \sqrt{%
\left\vert J\right\vert _{\omega }} \\
&=&
\sum\limits_{s=0}^{\infty }\ \sum_{N\in \mathbb{Z}%
}\ \sum_{d=N-\varepsilon s+1}^{N}\ \sum_{\left( I,J\right) \in \mathcal{P}%
_{N,d}^{s}}\ \left\Vert \square _{J}^{\omega ,\mathbf{b}^{\ast
}}g\right\Vert _{L^{2}\left( \omega \right) }\mathrm{P}^{\alpha }\left(
J,\left\vert \square _{I}^{\sigma ,\mathbf{b}}f\right\vert \sigma \right) 
\sqrt{\left\vert J\right\vert _{\omega }}.
\end{eqnarray*}%
Now we use%
\begin{eqnarray*}
\mathrm{P}^{\alpha }\left( J,\left\vert \square _{I}^{\sigma ,\mathbf{b}%
}f\right\vert \sigma \right)  &=&\int_{I}\frac{\ell \left( J\right) }{\left(
\ell \left( J\right) +\left\vert y-c_{J}\right\vert \right) ^{n+1-\alpha }}%
\left\vert \square _{I}^{\sigma ,\mathbf{b}}f\left( y\right) \right\vert
d\sigma \left( y\right)  \\
&\lesssim &\frac{2^{N-s}}{2^{d\left( n+1-\alpha \right) }}\left\Vert \square
_{I}^{\sigma ,\mathbf{b}}f\right\Vert _{L^{2}\left( \sigma \right) }\sqrt{%
\left\vert I\right\vert _{\sigma }}
\end{eqnarray*}%
and apply Cauchy-Schwarz in $J$ and use $J\subset 4I\backslash I$ to get

\begin{eqnarray*}
&&
\sum\limits_{\left( I,J\right) \in \mathcal{P}}\left\vert \left\langle
T_{\sigma }^{\alpha }\left( \square _{I}^{\sigma ,\mathbf{b}}f\right)
,\square _{J}^{\omega ,\mathbf{b}^{\ast }}g\right\rangle _{\omega
}\right\vert  \\
&\lesssim &
\sum\limits_{s=0}^{\infty }\ \sum_{N\in \mathbb{Z}%
}\ \sum_{d=N-\varepsilon s-1}^{N}\ \sum_{I\in \mathcal{D}_{N}}\frac{%
2^{N-s}2^{N\left( n-\alpha \right) }}{2^{d\left( n+1-\alpha \right) }}%
\left\Vert \square _{I}^{\sigma ,\mathbf{b}}f\right\Vert _{L^{2}\left(
\sigma \right) }\frac{\sqrt{\left\vert I\right\vert _{\sigma }}\sqrt{%
\left\vert 4I\backslash I\right\vert _{\omega }}}{2^{N\left( n-\alpha \right)
}}\cdot \\
&&
\hspace{6cm}\cdot\sqrt{\sum_{\substack{ J\in \mathcal{G}_{N-s} \\ J\subset 4I\backslash I\text{ and }
d\left( I,J\right) \approx 2^{d}}}\!\!\!\!\!\!\!\!\!\left\Vert \square _{J}^{\omega ,\mathbf{b}^{\ast }}g\right\Vert _{L^{2}\left( \omega \right) }^{2}} \\
&\lesssim &
\left( 1+\varepsilon s\right)\!\! \sum\limits_{s=0}^{\infty } \sum_{N\in \mathbb{Z}}\frac{2^{N-s}2^{N\left( n-\alpha \right) }%
}{2^{\left( N-\varepsilon s\right) \left( n+1-\alpha \right) }}\!\sqrt{\mathfrak{A}_{2}^{\alpha }}\!\!\sum_{I\in \mathcal{D}_{N}}\!\!\left\Vert \square _{I}^{\sigma ,\mathbf{b}}f\right\Vert _{L^{2}\left( \sigma \right) }\!\!\sqrt{\sum_{\substack{ J\in \mathcal{G}_{N-s} \\ J\subset 4I\backslash I}}\!\!\!\!\!\left\Vert \square
_{J}^{\omega ,\mathbf{b}^{\ast }}g\right\Vert _{L^{2}\left( \omega \right)
}^{2}} \\
&\lesssim &
\left( 1+\varepsilon s\right) \sum\limits_{s=\mathbf{\rho }%
}^{\infty }2^{-s\left[ 1-\varepsilon \left( n+1-\alpha \right) \right] }%
\sqrt{\mathfrak{A}_{2}^{\alpha }}\left\Vert f\right\Vert _{L^{2}\left( \sigma \right)
}\left\Vert g\right\Vert _{L^{2}\left( \omega \right)}
\lesssim 
\sqrt{\mathfrak{A}_{2}^{\alpha }}\left\Vert f\right\Vert _{L^{2}\left( \sigma \right)
}\left\Vert g\right\Vert _{L^{2}\left( \omega \right) }
\end{eqnarray*}%
where in the third line above we have used $\displaystyle \sum_{d=N-\varepsilon
s-1}^{N}1\lesssim 1+\varepsilon s$, and in the last line $\displaystyle \frac{2^{N-s}2^{N\left( n-\alpha \right) }}{2^{\left( N-\varepsilon s\right)
\left( n+1-\alpha \right) }}=2^{-s\left[ 1-\varepsilon \left( n+1-\alpha
\right) \right] }$ followed by Cauchy-Schwarz in $I$ and $N$, using that we
have bounded overlap, depending only on dimension and the goodness constant in the quadruples of $I$ for $I\in \mathcal{D}_{N}$. More
precisely, if we define $f_{k}\equiv \Psi _{\mathcal{D}_{k}}^{\sigma ,%
\mathbf{b}}f=\sum\limits_{I\in \mathcal{D}_{k}}\square _{I}^{\sigma ,\mathbf{b}}f$
and $g_{k}\equiv \Psi _{\mathcal{G}_{k}}^{\sigma ,\mathbf{b}^{\ast
}}g=\sum\limits_{J\in \mathcal{G}_{k}}\square _{J}^{\omega ,\mathbf{b}^{\ast }}g$,
then we have the quasi-orthogonality inequality 
\begin{eqnarray*}
\sum_{N\in \mathbb{Z}}\left\Vert f_{N}\right\Vert _{L^{2}\left( \sigma
\right) }\left\Vert g_{N-s}\right\Vert _{L^{2}\left( \omega \right) } 
\leq
\left( \sum_{N\in \mathbb{Z}}\left\Vert f_{N}\right\Vert _{L^{2}\left(
\sigma \right) }^{2}\right) ^{\frac{1}{2}}\left( \sum_{N\in \mathbb{Z}%
}\left\Vert g_{N-s}\right\Vert _{L^{2}\left( \omega \right) }^{2}\right) ^{%
\frac{1}{2}} 
\lesssim 
\left\Vert f\right\Vert _{L^{2}\left( \sigma \right) }\left\Vert
g\right\Vert _{L^{2}\left( \omega \right) }.
\end{eqnarray*}%
We have assumed that 
\begin{equation}
0<\varepsilon <\frac{1}{n+1-\alpha }  \label{short requirement}
\end{equation}%
in the calculations above, and this completes the proof of Lemma \ref{delta
short}.
\end{proof}

\section{Nearby form\label{Sec nearby}}

We dominate the nearby form $\Theta_3(f,g)$ by

\begin{equation*}
\left\vert \Theta _{3}\left( f,g\right) \right\vert \leq \sum_{I\in \mathcal{%
D}}\sum_{\substack{ J\in \mathcal{G}:\ 2^{-\mathbf{r}n}|I|
<|J| \leq|I|  \\ d\left( J,I\right) \leq
2\ell(J) ^{\varepsilon }\ell(I) ^{1-\varepsilon }}}%
\left\vert \int \left( T_{\sigma }^{\alpha }\square _{I}^{\sigma ,\mathbf{b}%
}f\right) \square _{J}^{\omega ,\mathbf{b}^{\ast }}gd\omega \right\vert ,
\end{equation*}
and prove the following proposition that controls the expectation, over two
independent grids, of the nearby form $\Theta _{3}\left( f,g\right) $. It
should be noted that weak goodness plays no role in treating the nearby form. Note also that in various steps we will use a small $\delta>0$. In all those different instances $\delta$ is free of any dependence. Our goal is the following proposition.

\begin{prop}\label{The prop}
Suppose $T^{\alpha }$ is a standard fractional singular
integral with $0\leq \alpha <n$. Let $\theta \in \left( 0,1\right) $ be
sufficiently small depending only on $\alpha, n$. Then there is a
constant $C_{\theta }$ such that for $f\in L^{2}\left( \sigma \right) $ and $%
g\in L^{2}\left( \omega \right) $, and dual martingale differences $\square
_{I}^{\sigma ,\mathbf{b}}$ and $\square _{J}^{\omega ,\mathbf{b}^{\ast }}$
with $\infty $-weakly accretive families of test functions $\mathbf{b}$
and $\mathbf{b}^{\ast }$, we have%
\begin{eqnarray}
&&\boldsymbol{E}_{\Omega }^{\mathcal{D}}\boldsymbol{E}_{\Omega }^{\mathcal{G}%
}\sum_{I\in \mathcal{D}}\sum_{\substack{ J\in \mathcal{G}:\ 2^{-\mathbf{r} n}|I| <|J| \leq |I|  \\ %
d\left( J,I\right) \leq 2\ell(J) ^{\varepsilon }\ell(I) ^{1-\varepsilon }}}\!\!\!\left\vert \left\langle T_{\sigma }^{\alpha
}\left( \square _{I}^{\sigma ,\mathbf{b}}f\right) ,\square _{J}^{\omega ,%
\mathbf{b}^{\ast }}g\right\rangle _{\omega }\right\vert  \label{delta near}
\\
&\lesssim &
\left( C_{\theta }\mathcal{NTV}_{\alpha }+\sqrt{\theta }\mathfrak{%
N}_{T^{\alpha }}\right) \left\Vert f\right\Vert _{L^{2}\left( \sigma \right)
}\left\Vert g\right\Vert _{L^{2}\left( \omega \right) }\ .  \notag
\end{eqnarray}
\end{prop}

 \newpage
 The following diagram is a sketch of the proof of proposition \eqref{The prop}.
\begin{center}
    \includegraphics[height=20 cm, width=18.5 cm]{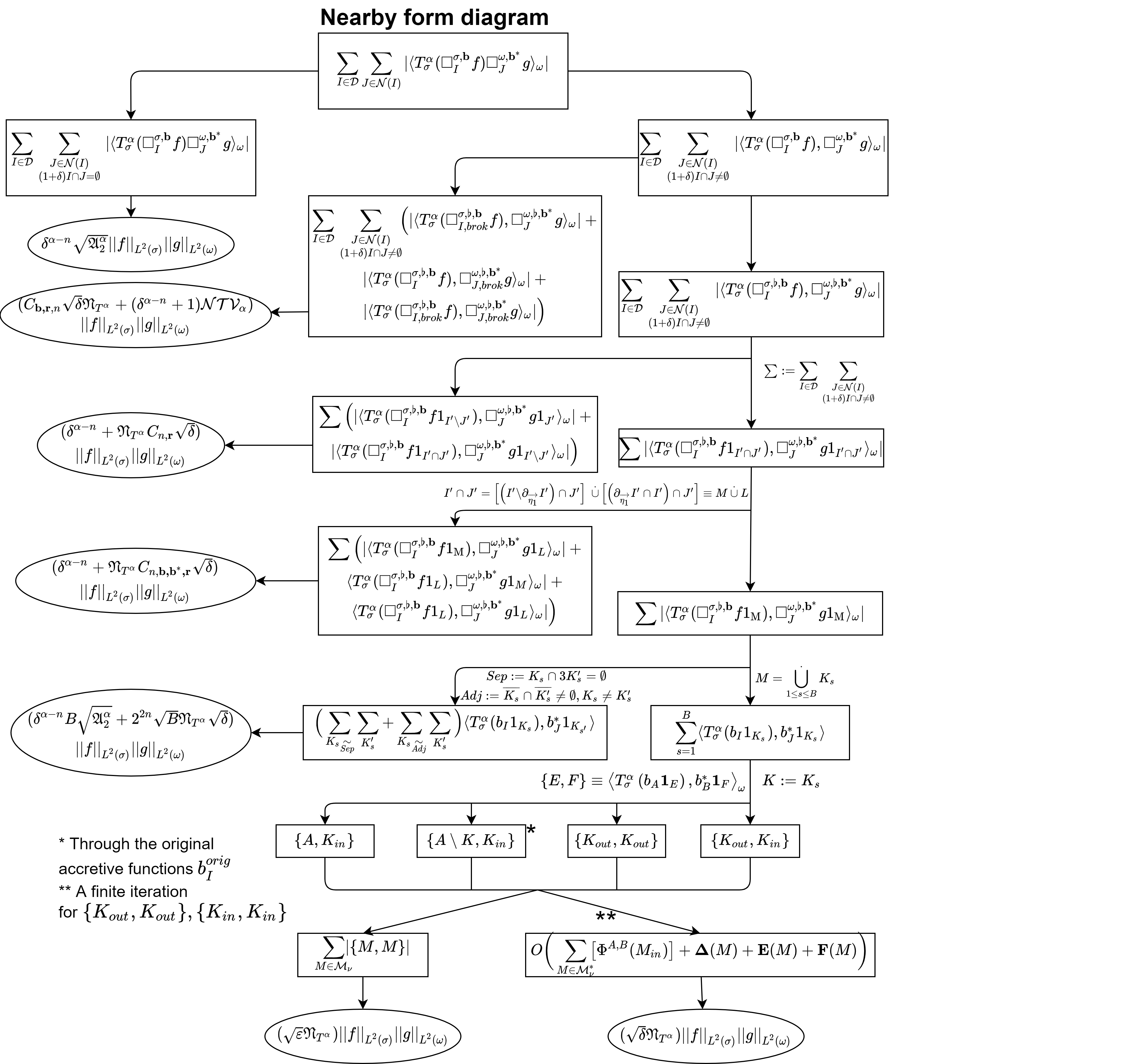}
\end{center}

\newpage

Before we proceed any further let us mention that
we will repeatedly use the inequality
\begin{equation}\label{PLBP removed}
\left\Vert \widehat{\square }%
_{I}^{\sigma ,\flat ,\mathbf{b}}f\right\Vert _{L^{2}\left( \sigma \right)
}\lesssim \left\Vert \square _{I}^{\sigma ,\mathbf{b}}f\right\Vert
_{L^{2}\left( \sigma \right) }^{\bigstar } 
\end{equation}

\begin{lem} \label{lemma star}
For $f \!\in\! L^2(\sigma)$ and $I\!\in\! \mathcal{C}_\mathcal{A}(A)$ we have
$\left\Vert \widehat{\square }%
_{I}^{\sigma ,\flat ,\mathbf{b}}f\right\Vert _{L^{2}\left( \sigma \right)
}\lesssim \left\Vert \square _{I}^{\sigma ,\mathbf{b}}f\right\Vert
_{L^{2}\left( \sigma \right) }^{\bigstar }$.
\end{lem}
\begin{proof}
Let $I^{\prime }\in \mathfrak{C}_{%
\mathcal{D}}\left( I\right) \cap \mathcal{C}_{\mathcal{A}}\left( A\right)$. Since $I^\prime\in \mathcal{C}_{\mathcal{A}}\left( A\right)$, from the corona construction we have
\begin{equation}
\label{accretivity2}
\left\vert \frac{1}{\left\vert I^{\prime }\right\vert _{\sigma }}%
\int_{I^{\prime }}b_{A}d\sigma \right\vert >\gamma.
\end{equation}
Now let $\{I^\prime_j\}_{j \in \mathbb{N}}$ be the collection of maximal subcubes $S$ of $I^\prime$ such that 
\begin{equation*}
\left\vert \frac{1}{\left\vert S\right\vert _{\sigma }}\int_{S}b_{A}d\sigma
\right\vert <\gamma ^{2}.
\end{equation*}
Let $\displaystyle E=\bigcup_jI^\prime_j$. We then have
\begin{equation*}
\left\vert \int_{E }b_{A}d\sigma \right\vert \leq \sum_{j}\left\vert \int_{I^\prime_j}b_{A}d\sigma \right\vert
<\gamma^2 \sum_{j}\left\vert I^\prime_j\right\vert
_{\sigma }\leq \gamma^2 \left\vert I^\prime\right\vert _{\sigma }\,
\end{equation*}
which together with \eqref{accretivity2} gives%
\begin{eqnarray*}
\gamma\left\vert I^\prime\right\vert _{\sigma } 
&<&
\left\vert \int_{I^\prime}b_{A}d\sigma
\right\vert =\left\vert \int_{E}b_{A}d\sigma \right\vert
+\left\vert \int_{I^\prime\setminus E}b_{A}d\sigma \right\vert \\
&\leq &
\gamma^2 \left\vert I^\prime\right\vert _{\sigma }+\sqrt{\int_{I^\prime\setminus
E }\left\vert b_{A}\right\vert ^{2}d\sigma }\sqrt{\left\vert
I^\prime\setminus E \right\vert _{\sigma }} \\
&\leq &
\gamma^2 \left\vert I^\prime\right\vert _{\sigma }+ C_{\mathbf{b}}\left\vert I^\prime\setminus E \right\vert _{\sigma },
\end{eqnarray*}%
where in the last inequality we used the $\infty$-accretivity of $b_A$. Rearranging the inequality yields successively%
\begin{eqnarray*}
\gamma\left( 1-\gamma \right) \left\vert I^\prime\right\vert _{\sigma } 
&\leq & 
C_{\mathbf{b}}\left\vert
I^\prime\setminus E \right\vert _{\sigma }; \\
\frac{\gamma\left( 1-\gamma \right)}{C_{\mathbf{b}}}\left\vert
I^\prime\right\vert _{\sigma }
&\leq &
\left\vert I^\prime\setminus E
\right\vert _{\sigma }\ ,
\end{eqnarray*}%
which in turn gives%
\begin{eqnarray}
\label{carleson2}
\ \ \ \ \ \sum_{j}\left\vert I^\prime_j\right\vert _{\sigma }\!\!\!
&=&\!\!
\left\vert I^\prime\right\vert _{\sigma }-\left\vert I^\prime\setminus E
\right\vert _{\sigma } \\\notag
&\leq &\!\!
\left\vert I^\prime\right\vert _{\sigma }-\frac{\gamma\left( 1-\gamma \right) 
}{C_{\mathbf{b}}}\left\vert I^\prime\right\vert _{\sigma }=\left( 1-
\frac{\gamma\left( 1-\gamma \right)}{C_{\mathbf{b}}}\right)
\left\vert I^\prime\right\vert _{\sigma }\equiv \beta \left\vert I^\prime\right\vert
_{\sigma } \notag
\end{eqnarray}%
where $0<\beta <1$ since $1\leq C_{\mathbf{b}}$. This implies
$$
\left\vert I^\prime\right\vert_\sigma \leq\frac{1}{1-\beta} \left\vert I^\prime\setminus E
\right\vert _{\sigma }
$$
Having that in hand and the fact that $\widehat{\square }_{I}^{\sigma ,\flat ,\mathbf{b}%
}f$ is constant on $I^\prime$, say  $\mathbf{1}_{I^\prime} \widehat{\square }_{I}^{\sigma ,\flat ,\mathbf{b}%
}f=c_{I^\prime}$ we can now calculate:
\begin{eqnarray*} 
\left\Vert 
\mathbf{1}_{I^{\prime }}\widehat{\square }_{I}^{\sigma ,\flat ,\mathbf{b}%
}f\right\Vert _{L^{2}\left( \sigma \right) }^{2}
&=&
\int_{I^{\prime }}\left\vert \widehat{\square }_{I}^{\sigma ,\flat ,\mathbf{b%
}}f\right\vert ^{2}d\sigma  =\left\vert I^{\prime }\right\vert _{\sigma
}\left\vert c_{I^{\prime }}\right\vert ^{2}\\
&\leq&
\frac{\displaystyle\frac{1}{\left\vert
I^{\prime }\setminus E\right\vert _{\sigma }}\int_{I^{\prime
}\setminus E}\left\vert b_{A}\right\vert ^{2}d\sigma }{\gamma^{4}}%
\left\vert I^{\prime }\right\vert _{\sigma }\left\vert c_{I^{\prime
}}\right\vert ^{2} \notag\\
&=&
\frac{1}{\gamma^4}\frac{\left\vert I^{\prime }\right\vert _{\sigma }}{%
\left\vert I^{\prime }\setminus E\right\vert _{\sigma }}%
\int_{I^{\prime }\setminus E}\left\vert b_{A}\right\vert
^{2}\left\vert c_{I^{\prime }}\right\vert ^{2}d\sigma \notag \\
&\leq &
\frac{1}{\gamma^4}\frac{\left\vert I^{\prime }\right\vert _{\sigma }}{
\left\vert I^{\prime }\setminus E\right\vert _{\sigma }}%
\int_{I^{\prime }}\left\vert b_{A}\widehat{\square }_{I}^{\sigma ,\flat ,
\mathbf{b}}f\right\vert ^{2}d\sigma\notag\\
&\leq&
\frac{1}{\gamma^4}\frac{1}{1-\beta}%
\int_{I^{\prime }}\left\vert b_{A}\widehat{\square }_{I}^{\sigma ,\flat ,
\mathbf{b}}f\right\vert ^{2}d\sigma,\notag
\end{eqnarray*}
and thus for $I'\in \mathcal{C}_{A}$ we obtain
\begin{equation*}
\int_{I^{\prime }}\left\vert \widehat{\square }_{I}^{\sigma ,\flat ,\mathbf{b%
}}f\right\vert ^{2}d\sigma \lesssim \int_{I^{\prime }}\left\vert b_{A}%
\widehat{\square }_{I}^{\sigma ,\flat ,\mathbf{b}}f\right\vert ^{2}d\sigma ,
\end{equation*}%
which in turn gives, after summing over all $I^{\prime }\in \mathfrak{C}_{%
\mathcal{D}}\left( I\right) \cap C_{\mathcal{A}}\left( A\right)$,
\begin{equation*}
\sum_{I^{\prime }\in \mathfrak{C}_{\mathcal{D}}\left( I\right) \cap C_{%
\mathcal{A}}\left( A\right) }\left\Vert \mathbf{1}_{I^{\prime }}\widehat{%
\square }_{I}^{\sigma ,\flat ,\mathbf{b}}f\right\Vert _{L^{2}\left( \sigma
\right) }^{2}\lesssim \left\Vert \mathbf{1}_{I}b_{A}\widehat{%
\square }_{I}^{\sigma ,\flat ,\mathbf{b}}f\right\Vert _{L^{2}\left( \sigma
\right) }^{2}\leq \left\Vert b_{A}\widehat{\square }_{I}^{\sigma ,\flat ,%
\mathbf{b}}f\right\Vert _{L^{2}\left( \sigma \right) }^{2}.
\end{equation*}%

Now if $I^\prime \in \mathfrak{C}_{\mathcal{D}}\left( I\right) \cap \mathcal{A}$, from the definition of $\widehat{\nabla }_{Q}^{\mu }f$ in (\ref{Carleson avg op}),
$$
\sum_{I^{\prime }\in \mathfrak{C}_{\mathcal{D}}\left( I\right) \cap \mathcal{%
A}}\left\Vert \mathbf{1}_{I^{\prime }}\widehat{\square }_{I}^{\sigma ,\flat ,%
\mathbf{b}}f\right\Vert _{L^{2}\left( \sigma \right) }^{2}\lesssim
\left\Vert \widehat{\nabla }_{I}^{\sigma }f\right\Vert _{L^{2}\left( \sigma
\right) }^{2}\ .
$$
Now we are ready to prove (\ref{PLBP removed}). As $b_A=b_I$ and
\begin{eqnarray*}
\left\Vert \widehat{\square }_{I}^{\sigma ,\flat ,\mathbf{b}}f\right\Vert
_{L^{2}\left( \sigma \right) }^{2}
&=& \!\!\!\!\!
\sum_{I^{\prime }\in \mathfrak{C}_{\mathcal{D}}\left( I\right) \cap C_{\mathcal{A}}\left( A\right) }\left\Vert 
\mathbf{1}_{I^{\prime }}\widehat{\square }_{I}^{\sigma ,\flat ,\mathbf{b}%
}f\right\Vert _{L^{2}\left( \sigma \right) }^{2}+\sum_{I^{\prime }\in 
\mathfrak{C}_{\mathcal{D}}\left( I\right) \cap \mathcal{A}}\left\Vert 
\mathbf{1}_{I^{\prime }}\widehat{\square }_{I}^{\sigma ,\flat ,\mathbf{b}%
}f\right\Vert _{L^{2}\left( \sigma \right) }^{2}\\
&\lesssim&
\left\Vert b_{I}\widehat{\square }
_{I}^{\sigma ,\flat ,\mathbf{b}}f\right\Vert _{L^{2}\left( \sigma \right)
}^2 
+
\left\Vert \widehat{\nabla }_{I}^{\sigma }f\right\Vert _{L^{2}\left(
\sigma \right) }^2
\end{eqnarray*}
we obtain

\begin{eqnarray*}
\left\Vert \widehat{\square }_{I}^{\sigma ,\flat ,\mathbf{b}}f\right\Vert
_{L^{2}\left( \sigma \right) } 
\!\!\!\!&\lesssim &\!\!\!\!
\left\Vert b_{I}\widehat{\square }%
_{I}^{\sigma ,\flat ,\mathbf{b}}f\right\Vert _{L^{2}\left( \sigma \right)
}+\left\Vert \widehat{\nabla }_{I}^{\sigma }f\right\Vert _{L^{2}\left(
\sigma \right) }=\left\Vert \square _{I}^{\sigma ,\flat ,\mathbf{b}%
}f\right\Vert _{L^{2}\left( \sigma \right) }+\left\Vert \widehat{\nabla }%
_{I}^{\sigma }f\right\Vert _{L^{2}\left( \sigma \right) } \\
&\leq &\!\!\!\!
\left\Vert \square _{I}^{\sigma ,\mathbf{b}}f\right\Vert
_{L^{2}\left( \sigma \right) }+\left\Vert \square _{I,{broken}%
}^{\sigma ,\flat ,\mathbf{b}}f\right\Vert _{L^{2}\left( \sigma \right)
}+\left\Vert \widehat{\nabla }_{I}^{\sigma }f\right\Vert _{L^{2}\left(
\sigma \right) }\lesssim \left\Vert \square _{I}^{\sigma ,\mathbf{b}%
}f\right\Vert _{L^{2}\left( \sigma \right) }^{\bigstar }.
\end{eqnarray*}

\end{proof}

Now from quasiorthogonality and (\ref{PLBP removed}) we get,
\begin{eqnarray*}
\sum_{J\in \mathcal{G}}\sum_{J^{\prime }\in \mathfrak{C}\left( J\right)
}\left\vert J^{\prime }\right\vert _{\omega }\left\vert E_{J^{\prime
}}^{\omega }\left( \widehat{\square }_{J}^{\omega ,\flat ,\mathbf{b}^{\ast
}}g\right) \right\vert ^{2}
&\lesssim &\!\!\!
\sum_{J\in \mathcal{G}}\left\Vert 
\widehat{\square }_{J}^{\omega ,\flat ,\mathbf{b}^{\ast }}g\right\Vert
_{L^{2}\left( \omega\right) }^{2}
\lesssim
\sum_{J\in \mathcal{G}}\left\Vert
\square _{J}^{\omega ,\flat ,\mathbf{b}^{\ast }}g\right\Vert _{L^{2}\left(\omega\right) }^{2} \\
&\lesssim &\!\!\!
\sum_{J\in \mathcal{G}}\left( \left\Vert \square _{J}^{\omega ,%
\mathbf{b}^{\ast }}g\right\Vert _{L^{2}\left( \omega\right) }^{2}+\left\Vert
\nabla _{J}^{\omega }g\right\Vert _{L^{2}\left( \omega\right) }^{2}\right)
\lesssim
\left\Vert g\right\Vert _{L^{2}\left( \omega \right) }^{2}\ .
\end{eqnarray*}

We also need the following lemma, that controls the above inner product for cubes of positive distance.
\begin{lem}\label{lemma1}
Given the $\infty $-weakly accretive families of test functions $\mathbf{b}$ and $\mathbf{b}^{\ast }$ and cubes $Q,R \subset \R^n$, we have
\begin{equation}
    |\langle T^\alpha_\sigma(b_Q\mathbf{1}_Q),b^*_R\mathbf{1}_{R\backslash (1+\delta)Q}\rangle_\omega|\lesssim \delta^{\alpha-n}\sqrt{\mathfrak{A}^{\alpha}_2}
\sqrt{|Q|_\sigma}\sqrt{|R|_\omega}\end{equation}
where the implied constant depends on the accretivity constants of the families $\mathbf{b, b^*}$ and the dimension $n$.
\end{lem}

\begin{proof}
We have that
$\left\vert \left\langle T_{\sigma }^{\alpha }\left( b_{Q}\mathbf{1}%
_{Q}\right) , 
b_{R}^{\ast }\mathbf{1}_{R\backslash \left( 1+\delta \right)
Q}\right\rangle _{\omega }\right\vert$
\begin{eqnarray*} 
&\leq &
\int_{R\backslash \left(
1+\delta \right) Q}\left\vert T_{\sigma }^{\alpha }\left( b_{Q}\mathbf{1}%
_{Q}\right) \right\vert \left\vert b_{R}^{\ast }\right\vert d\omega \\
&\leq &
\left( \int_{R\backslash \left( 1+\delta \right) Q}\left\vert
T_{\sigma }^{\alpha }\left( b_{Q}\mathbf{1}_{Q}\right) \right\vert
^{2}d\omega \right) ^{\frac{1}{2}}\left( \int_{R\backslash \left( 1+\delta
\right) Q}\left\vert b_{R}^{\ast }\right\vert ^{2}d\omega \right) ^{\frac{1}{%
2}} \\
&\lesssim &
\left( \int_{\mathbb{R}^{n}\backslash \left( 1+\delta \right)
Q}\left( \int_{Q}\left\vert x-y\right\vert ^{\alpha -n}\left\vert
b_{Q}\left( y\right) \right\vert d\sigma \left( y\right) \right) ^{2}d\omega
\left( x\right) \right) ^{\frac{1}{2}}\left( \int_{R}\left\vert b_{R}^{\ast
}\right\vert ^{2}d\omega \right) ^{\frac{1}{2}} 
\end{eqnarray*}
\begin{eqnarray*}
&\lesssim &
\left( \int_{\mathbb{R}^{n}\backslash \left( 1+\delta \right)
Q}\left( \int_{Q}\left( \delta \left\vert x-c_{Q}\right\vert \right)
^{\alpha -n}\left\vert b_{Q}\left( y\right) \right\vert d\sigma \left(
y\right) \right) ^{2}d\omega \left( x\right) \right) ^{\frac{1}{2}}\sqrt{%
\left\vert R\right\vert _{\omega }} \\
&\lesssim &
\delta ^{\alpha -n}\left( \int_{\mathbb{R}^{n}\backslash \left(
1+\delta \right) Q}\left\vert x-c_{Q}\right\vert ^{2\left( \alpha -n\right)
}d\omega \left( x\right) \right) ^{\frac{1}{2}}\left( \int_{Q}\left\vert
b_{Q}\left( y\right) \right\vert d\sigma \left( y\right) \right) \sqrt{%
\left\vert R\right\vert _{\omega }} \\
&\lesssim &
\delta ^{\alpha -n}\left( \int_{\mathbb{R}^{n}\backslash \left(
1+\delta \right) Q}\left\vert x-c_{Q}\right\vert ^{2\left( \alpha -n\right)
}d\omega \left( x\right) \right) ^{\frac{1}{2}}\left\vert Q\right\vert
_{\sigma }\sqrt{\left\vert R\right\vert _{\omega }}\\
&\leq& 
\delta ^{\alpha -n}%
\sqrt{\mathfrak{A}_{2}^{\alpha }}\sqrt{\left\vert Q\right\vert _{\sigma }}%
\sqrt{\left\vert R\right\vert _{\omega }}
\end{eqnarray*}%
since%
\begin{eqnarray*}
\left( \int_{\mathbb{R}^{n}\backslash \left( 1+\delta \right) Q}\!\!\!\!\!\!\!\!\!\!\!\!\!\!\!\!\!\!\!\!\left\vert
x-c_{Q}\right\vert ^{2\left( \alpha -n\right) }d\omega \left( x\right)
\right) \left\vert Q\right\vert _{\sigma } \!\!\!\!\!
&=& \!\!\!\!\!
\left( \int_{\mathbb{R}%
^{n}\backslash \left( 1+\delta \right) Q}\left( \frac{\left\vert Q\right\vert
^{\frac{1}{n}}}{\left\vert x-c_{Q}\right\vert ^{2}}\right) ^{n-\alpha
}d\omega \left( x\right) \right) \frac{\left\vert Q\right\vert _{\sigma }}{%
\left\vert Q\right\vert ^{1-\frac{\alpha }{n}}} \\
&\lesssim &
\mathcal{P}^{\alpha }\left( Q,\omega \right) \frac{\left\vert
Q\right\vert _{\sigma }}{\left\vert Q\right\vert ^{1-\frac{\alpha }{n}}}\leq 
\mathcal{A}_{2}^{\alpha ,\ast }.
\end{eqnarray*}
\end{proof}

 As usual, we continue to write the independent grids for $f$
and $g$ as $\mathcal{D}$ and $\mathcal{G}$ respectively. Write the dual
martingale averages $\square _{I}^{\sigma ,\mathbf{b}}f$ and $\square
_{J}^{\omega ,\mathbf{b}^{\ast }}g$ as linear combinations 
\begin{eqnarray*}
\square _{I}^{\sigma ,\mathbf{b}}f 
&=&
b_{I}\ \sum_{I^{\prime }\in \mathfrak{C%
}_{{nat}}\left( I\right) }\mathbf{1}_{I^{\prime }}\ E_{I^{\prime
}}^{\sigma }\left( \widehat{\square }_{I}^{\sigma ,\mathbf{b}}f\right)
+
\sum_{I^{\prime }\in \mathfrak{C}_{ {brok}}\left( I\right)
}b_{I^{\prime }}\ \mathbf{1}_{I^{\prime }}\widehat{\mathbb{F}}_{I^{\prime
}}^{\sigma ,b_{I^{\prime }}}f-b_{I}\ \sum_{I^{\prime }\in \mathfrak{C}_{
 {brok}}\left( I\right) }\mathbf{1}_{I^{\prime }}\widehat{\mathbb{F}%
}_{I}^{\sigma ,b_{I}}f, \\
\square _{J}^{\omega ,\mathbf{b}^{\ast }}g 
&=&
b_{J}^{\ast }\ \sum_{J^{\prime
}\in \mathfrak{C}_{ {nat}}\left( J\right) }\mathbf{1}_{J^{\prime
}}\ E_{J^{\prime }}^{\omega }\left( \widehat{\square }_{J}^{\omega ,\mathbf{b%
}^{\ast }}g\right) +
 \sum_{J^{\prime }\in \mathfrak{C}_{ {brok}%
}\left( J\right) }b_{J^{\prime }}^{\ast }\ \mathbf{1}_{J^{\prime }}\widehat{%
\mathbb{F}}_{J^{\prime }}^{\omega ,b_{J^{\prime }}^{\ast }}g-b_{J}^{\ast }\
\sum_{J^{\prime }\in \mathfrak{C}_{ {brok}}\left( J\right) }\mathbf{%
1}_{J^{\prime }}\widehat{\mathbb{F}}_{J}^{\omega ,b_{J}^{\ast }}g,
\end{eqnarray*}%
of the appropriate function $b$ times the indicators of their children,
denoted $I^{\prime }$ and $J^{\prime }$ respectively. We will regroup the
terms as needed below.

On the natural child $I^{\prime }$, the expression $\widehat{\square }%
_{I}^{\sigma ,\mathbf{b}}f=\frac{1}{b_{I}}\square _{I}^{\sigma ,\mathbf{b}}f$
simply denotes the dual martingale average with $b_{I}$ removed, so that we
need not assume $\left\vert b_{I}\right\vert $ is bounded below in order to
make sense of $\frac{1}{b_{I}}\square _{I}^{\sigma ,\mathbf{b}}f$. Similar
comments apply to the expressions $\widehat{\mathbb{F}}_{I^{\prime
}}^{\sigma ,b_{I^{\prime }}}f=\frac{1}{b_{I^{\prime }}}\mathbb{F}_{I^{\prime
}}^{\sigma ,b_{I^{\prime }}}f$ and $\widehat{\mathbb{F}}_{I}^{\sigma
,b_{I}}f=\frac{1}{b_{I}}\mathbb{F}_{I}^{\sigma ,b_{I}}f$. 
Now if we set 
$$
\mathcal{N}(I)=\{ J\in \mathcal{G}:\ 2^{-\mathbf{r}n}|I|
<|J| \leq |I|  , d\left( J,I\right) \leq
2\ell(J) ^{\varepsilon }\ell(I) ^{1-\varepsilon }\}
$$
for the cubes or similar size to $I$, the left hand side of (\ref{delta near}) is bounded by
\begin{eqnarray}\label{1st}
\mathrm{\mathbf{I+II}} 
&\equiv&
\sum_{I\in \mathcal{D}}\!\!\sum_{\substack{J\in\mathcal{N}(I)\\ (1+\delta)I\cap J=\emptyset}}\!\! \left\vert \left\langle T_{\sigma }^{\alpha
}\left( \square _{I}^{\sigma ,\mathbf{b}}f\right) ,\square _{J}^{\omega ,
\mathbf{b}^{\ast }}g\right\rangle _{\omega }\right\vert  \\
&&+
\sum_{I\in \mathcal{D}}\sum_{\substack{J\in\mathcal{N}(I)\\ (1+\delta)I\cap J\neq\emptyset}} \left\vert \left\langle T_{\sigma }^{\alpha
}\left( \square _{I}^{\sigma ,\mathbf{b}}f\right) ,\square _{J}^{\omega ,%
\mathbf{b}^{\ast }}g\right\rangle _{\omega }\right\vert \nonumber 
\end{eqnarray}

When working in higher dimensions, run the proof pretending you have Hyt\"{o}nen's estimate (which is of course not true due to \cite{GP}). Then wherever we were supposed to use Hyt\"{o}nen, we use the delta separation trick. The $\delta$-separated part is easily seen to be bounded by the Muckenhoupt conditions, and the $\delta$-close part will give a $\sqrt{\delta}$ estimate. But $\delta$, which can be chosen at the end, is independent of everything else (it is the Hyt\"{o}nen-delta, not related to anything else in the proof). So, provided the proof only deals with finite estimates and finitely many constructions (like the Cantor set construction, that only does finitely many iterations), those $\sqrt{\delta}$ terms will be absorbable at the end. Here are the details:

\subsection{The case of $\delta$-separated cubes.}
In this subsection we are estimating $\mathrm{\mathbf{I}}$ in \eqref{1st} by using Lemma \ref{lemma1}.
\begin{dfn}
We say that the cubes $J$ and $I$ are $\delta$-\emph{separated}, where $\delta>0$, if $J\cap(1+\delta)I=\emptyset $.
\end{dfn}
For the first sum in \eqref{1st} we have,  following the proof of Lemma \ref%
{lemma1}, the satisfactory estimate%
\begin{equation*}
\left\vert \left\langle T_{\sigma }^{\alpha }\left( \square _{I}^{\sigma ,%
\mathbf{b}}f\right) ,\square _{J}^{\omega ,\mathbf{b}^{\ast }}g\right\rangle
_{\omega }\right\vert 
\lesssim
\delta^{\alpha-n} \sqrt{
\mathfrak{A}^\alpha_{2}}\left\Vert \square _{I}^{\sigma ,\mathbf{b}%
}f\right\Vert _{L^{2}\left( \sigma \right) }\left\Vert \square _{J}^{\omega ,%
\mathbf{b}^{\ast }}g\right\Vert _{L^{2}\left( \omega \right) }.
\end{equation*}%
Indeed,
\begin{eqnarray*}
& &\left\vert \left\langle T_{\sigma }^{\alpha }\left( \square _{I}^{\sigma ,%
\mathbf{b}}f\right) , \square _{J}^{\omega ,\mathbf{b}^{\ast }}g\right\rangle
_{\omega }\right\vert \\
&\leq &
\int_{J\backslash \left(
1+\delta \right) I}\left\vert T_{\sigma }^{\alpha }\left( \square _{I}^{\sigma ,%
\mathbf{b}}f\right) \right\vert \left\vert\square _{J}^{\omega ,\mathbf{b}^{\ast }}g\right\vert d\omega \\
&\leq &
\left( \int_{J\backslash \left( 1+\delta \right) I}\left\vert
T_{\sigma }^{\alpha }\left( \square _{I}^{\sigma ,%
\mathbf{b}}f\right) \right\vert
^{2}d\omega \right) ^{\frac{1}{2}}\left( \int_{J}\left\vert\square _{J}^{\omega ,\mathbf{b}^{\ast }}g\right\vert ^{2}d\omega \right) ^{\frac{1}{%
2}} \\
&\lesssim &
\delta ^{\alpha -n}\left( \int_{\mathbb{R}^{n}\backslash \left(
1+\delta \right) I}\!\!\!\!\!\!\!\!\!\left\vert x-c_{I}\right\vert ^{2\left( \alpha -n\right)
}d\omega \left( x\right) \right) ^{\frac{1}{2}}\left( \int_{I}\left\vert
\square _{I}^{\sigma ,%
\mathbf{b}}f  \right\vert d\sigma \left( y\right) \right) 
\Big|\Big|\square _{J}^{\omega ,\mathbf{b}^{\ast }}g\Big|\Big|_{L^2(\omega)}\\
&\lesssim &
\delta ^{\alpha -n}\left( \int_{\mathbb{R}^{n}\backslash \left(
1+\delta \right) I}\!\!\!\!\!\!\left\vert x-c_{I}\right\vert ^{2\left( \alpha -n\right)
}d\omega \left( x\right) \right) ^{\frac{1}{2}}\sqrt{\left\vert I\right\vert
_{\sigma }}\Big|\Big|\square _{I}^{\sigma ,%
\mathbf{b}}f\Big|\Big|_{L^2(\sigma)}\Big|\Big|\square _{J}^{\omega ,\mathbf{b}^{\ast }}g\Big|\Big|_{L^2(\omega)}\\
&\leq& 
\delta ^{\alpha -n}%
\sqrt{\mathfrak{A}_{2}^{\alpha}}\Big|\Big|\square _{I}^{\sigma ,%
\mathbf{b}}f\Big|\Big|_{L^2(\sigma)}\Big|\Big|\square _{J}^{\omega ,\mathbf{b}^{\ast }}g\Big|\Big|_{L^2(\omega)}
\end{eqnarray*}%
So combining all the above we get for the $\delta$-separated cubes that
\begin{eqnarray}
&&\label{first bound}\\
\mathrm{\mathbf{I}}
&\leq&
\sum_{I\in \mathcal{D}}\sum_{\substack{J\in\mathcal{N}(I)\\ (1+\delta)I\cap J=\emptyset}} \delta^{\alpha-n} \sqrt{
\mathfrak{A}^\alpha_{2}}\left\Vert \square _{I}^{\sigma ,\mathbf{b}%
}f\right\Vert _{L^{2}\left( \sigma \right) }\left\Vert \square _{J}^{\omega,
\mathbf{b}^{\ast }}g\right\Vert _{L^{2}\left( \omega \right) }\notag\\
&\leq &\!\!\!\!\!
\delta^{\alpha-n} \sqrt{
\mathfrak{A}^\alpha_{2}} \left(\sum_{I\in \mathcal{D}}\sum_{\substack{J\in\mathcal{N}(I)\\ (1+\delta)I\cap J=\emptyset}} \left\Vert \square _{I}^{\sigma ,\mathbf{b}%
}f\right\Vert _{L^{2}\left( \sigma \right) }^2\right)^\frac{1}{2}\!\!
\left(\sum_{I\in \mathcal{D}}\sum_{\substack{J\in\mathcal{N}(I)\\ (1+\delta)I\cap J=\emptyset}}\left\Vert \square _{J}^{\omega,
\mathbf{b}^{\ast }}g\right\Vert _{L^{2}\left( \omega \right) }^2\right)^\frac{1}{2}\notag\\
&\lesssim&\!\!\!\!\!
\delta^{\alpha-n} \sqrt{
\mathfrak{A}^\alpha_{2}}||f||_{L^2(\sigma)}||g||_{L^2(\omega)}\notag
\end{eqnarray}
where the implied constant in the last line depends only on the goodness parameter $\mathbf{r}$ and the finite repetition of $I$ and $J$ in each sum respectively.

\subsection{The case of $\delta$-close cubes.}
Now we turn to the second sum in \eqref{1st} which we will bound by using random surgery and expectation.
\begin{dfn}
We say that the cubes $J$ and $I$ are $\delta$-close, if $J\cap (1+\delta)I\neq \emptyset.$
\end{dfn}
We have
\begin{eqnarray}\label{broken split}
\left\langle T_{\sigma }^{\alpha }\left( \square _{I}^{\sigma ,\mathbf{b}%
}f\right) ,\square _{J}^{\omega ,\mathbf{b}^{\ast }}g\right\rangle _{\omega
} \label{broken}
&=&
\left\langle T_{\sigma }^{\alpha }\left( \square _{I}^{\sigma ,\flat ,%
\mathbf{b}}f\right) ,\square _{J}^{\omega ,\flat ,\mathbf{b}^{\ast
}}g\right\rangle _{\omega }\\
&&+
\left\langle T_{\sigma }^{\alpha }\left( \square
_{I,{brok}}^{\sigma ,\flat ,\mathbf{b}}f\right) ,\square _{J,%
{brok}}^{\omega ,\flat ,\mathbf{b}^{\ast }}g\right\rangle _{\omega
}  \nonumber \\
&&+
\left\langle T_{\sigma }^{\alpha }\left( \square _{I}^{\sigma ,\flat ,%
\mathbf{b}}f\right) ,\square _{J,{brok}}^{\omega ,\flat ,\mathbf{b}%
^{\ast }}g\right\rangle _{\omega }\nonumber \\
&&+
\left\langle T_{\sigma }^{\alpha }\left(
\square _{I,{brok}}^{\sigma ,\flat ,\mathbf{b}}f\right) ,\square
_{J}^{\omega ,\flat ,\mathbf{b}^{\ast }}g\right\rangle _{\omega }\ .  \notag
\end{eqnarray}%
The estimation of the latter three inner products, i.e. those in which a
broken operator $\square _{I,{brok}}^{\sigma ,\flat ,\mathbf{b}}$
or $\square _{J,{brok}}^{\omega ,\flat ,\mathbf{b}^{\ast }}$
arises, is simpler,  but still requires 
the use of random surgery in order to avoid the full testing condition that was available in one dimension \cite{GP}.
Indeed, recall that 
\begin{eqnarray*}
\square _{I,{brok}}^{\sigma ,\flat ,\mathbf{b}}f
&=&
\sum_{I^{\prime }\in \mathfrak{C}_{{brok}}\left( I\right) }
\mathbb{F}_{I^{\prime }}^{\sigma ,\mathbf{b}}f=\sum_{I^{\prime }\in 
\mathfrak{C}_{{brok}}\left( I\right) }\left( E_{I^{\prime
}}^{\sigma }\widehat{\mathbb{F}}_{I^{\prime }}^{\sigma ,\mathbf{b}}f\right)
b_{I^{\prime }} \\
\square _{J,{brok}}^{\omega ,\flat ,\mathbf{b}^*}g
&=&
\sum_{J^{\prime }\in \mathfrak{C}_{{brok}}\left( J\right) }
\mathbb{F}_{J^{\prime }}^{\omega ,\mathbf{b}^{\ast }}g=\sum_{J^{\prime }\in 
\mathfrak{C}_{{brok}}\left( J\right) }\left( E_{J^{\prime
}}^{\omega }\widehat{\mathbb{F}}_{J^{\prime }}^{\omega ,\mathbf{b}^{\ast
}}g\right) b_{J^{\prime }}^{\ast }
\end{eqnarray*}%
so that if at least one broken difference appears in the inner product, as
is the case for the latter three inner products in \eqref{broken},
we need to use random surgery to get the necessary bound. For example, the fourth term satisfies
\begin{eqnarray*}
&&\left\vert \left\langle T_{\sigma }^{\alpha }\left( \square _{I,{%
brok}}^{\sigma ,\flat ,\mathbf{b}}f\right) ,\square _{J}^{\omega ,\flat ,%
\mathbf{b}^{\ast }}g\right\rangle _{\omega }\right\vert
=
\left\vert
\sum_{I^{\prime }\in \mathfrak{C}_{{brok}}\left( I\right) }\left(
E_{I^{\prime }}^{\sigma }\widehat{\mathbb{F}}_{I^{\prime }}^{\sigma ,\mathbf{%
b}}f\right) \left\langle T_{\sigma }^{\alpha }b_{I^{\prime }},\square
_{J}^{\omega ,\flat ,\mathbf{b}^{\ast }}g\right\rangle _{\omega }\right\vert
\end{eqnarray*}
and since
\begin{eqnarray*}
\left\langle T_{\sigma }^{\alpha }b_{I^{\prime }},\square
_{J}^{\omega ,\flat ,\mathbf{b}^{\ast }}g\right\rangle _{\omega }
&=&
\left\langle \mathbf{1}_{I'\cap J} T_{\sigma }^{\alpha }b_{I^{\prime }},\square
_{J}^{\omega ,\flat ,\mathbf{b}^{\ast }}g\right\rangle _{\omega }
+
\left\langle \mathbf{1}_{J\backslash(1+\delta)I'} T_{\sigma }^{\alpha }b_{I^{\prime }},\square
_{J}^{\omega ,\flat ,\mathbf{b}^{\ast }}g\right\rangle _{\omega }\\
&&+
\left\langle \mathbf{1}_{(J\backslash I^\prime)\cap(1+\delta)I'}T_{\sigma }^{\alpha }b_{I^{\prime }},\square
_{J}^{\omega ,\flat ,\mathbf{b}^{\ast }}g\right\rangle _{\omega }\\
&\equiv&
A(f,g)+B(f,g)+C(f,g)
\end{eqnarray*}
we have
\begin{eqnarray*}
&&\left\vert
\sum_{I^{\prime }\in \mathfrak{C}_{{brok}}\left( I\right) }\left(
E_{I^{\prime }}^{\sigma }\widehat{\mathbb{F}}_{I^{\prime }}^{\sigma ,\mathbf{%
b}}f\right) A(f,g)\right\vert\\
&\leq&
C_{\mathbf{b,b^*}}
\sum_{I^{\prime }\in \mathfrak{C}_{{brok}}\left(
I\right) }\left\vert E_{I^{\prime }}^{\sigma }\widehat{\mathbb{F}}%
_{I^{\prime }}^{\sigma ,\mathbf{b}}f\right\vert \mathfrak{T}_{T^{\alpha }}^{%
\mathbf{b}}\sqrt{\left\vert I^{\prime }\right\vert _{\sigma }}\left\Vert
\square _{J}^{\omega ,\flat ,\mathbf{b}^{\ast }}g\right\Vert _{L^{2}\left(
\omega \right) } \\
&\leq &
\mathfrak{T}_{T^{\alpha }}^{\mathbf{b}}\left\Vert \nabla
_{I}^{\sigma }f\right\Vert _{L^{2}\left( \sigma \right) }\!\!\Bigg(\!\!\sum_{I^{\prime
}\in \mathfrak{C}_{{brok}}\left( I\right) }\!\!\!\left( \left\Vert
\square _{J}^{\omega ,\mathbf{b}^{\ast }}g\right\Vert^2 _{L^{2}\left( \omega
\right) }+\left\Vert \square _{J,{brok}}^{\omega ,\flat ,\mathbf{b}%
^{\ast }}g\right\Vert^2 _{L^{2}\left( \omega \right) }\right)\!\!\Bigg)^\frac{1}{2} \\
&\lesssim &
\mathfrak{T}_{T^{\alpha }}^{\mathbf{b}}\left\Vert \square _{I}^{\sigma ,\mathbf{b}%
}f\right\Vert _{L^{2}\left( \sigma \right) }^{\bigstar }\left\Vert \square
_{J}^{\omega ,\mathbf{b}^{\ast }}g\right\Vert _{L^{2}\left( \omega \right)
}^{\bigstar }\
\end{eqnarray*}
Next by Lemma \ref{lemma1},
\begin{eqnarray*}
\left\vert
\sum_{ I^{\prime }\in \mathfrak{C}_{{brok}}\left( I\right) }\!\!\!\!\!\left(
E_{I^{\prime }}^{\sigma }\widehat{\mathbb{F}}_{I^{\prime }}^{\sigma ,\mathbf{%
b}}f\right) B(f,g)\right\vert
&\!\!\!\!\leq &\!\!\!\!\!\!\!
\sum_{I^{\prime }\in \mathfrak{C}_{{brok}}\left(
I\right) }\!\!\!\!\!\!\!\left\vert E_{I^{\prime }}^{\sigma }\widehat{\mathbb{F}}_{I^{\prime }}^{\sigma ,\mathbf{b}}f\right\vert \delta^{\alpha-n}\sqrt{\mathfrak{A}_2^\alpha}\sqrt{\left\vert I^{\prime }\right\vert _{\sigma }}\left\Vert
\square _{J}^{\omega ,\flat ,\mathbf{b}^{\ast }}g\right\Vert _{L^{2}\left(
\omega \right) }\\
&\!\!\!\!\leq&
\delta^{\alpha-n}\sqrt{\mathfrak{A}_2^{\alpha }}\left\Vert \square _{I}^{\sigma ,\mathbf{b}%
}f\right\Vert _{L^{2}\left( \sigma \right) }^{\bigstar }\left\Vert \square
_{J}^{\omega ,\mathbf{b}^{\ast }}g\right\Vert _{L^{2}\left( \omega \right)
}^{\bigstar }\
\end{eqnarray*}
Finally, using Cauchy-Schwarz, the norm inequality and accretivity we get 
\begin{eqnarray*}
&&\sum_{I\in \mathcal{D}}\sum_{\substack{J\in\mathcal{N}(I)\\ I\cap J\neq\emptyset}}\left\vert
\sum_{I^{\prime }\in \mathfrak{C}_{{brok}}\left( I\right) }\left(
E_{I^{\prime }}^{\sigma }\widehat{\mathbb{F}}_{I^{\prime }}^{\sigma ,\mathbf{%
b}}f\right) C(f,g)\right\vert \\
&\leq&
C_{\mathbf{b}}\mathfrak{N}_{T^\alpha}\sum_{I\in \mathcal{D}}\sum_{\substack{J\in\mathcal{N}(I)\\ I\cap J\neq\emptyset}}
\sum_{I^{\prime }\in \mathfrak{C}_{{brok}}(I)}\Big|E_{I^{\prime }}^{\sigma }\widehat{\mathbb{F}}_{I^{\prime }}^{\sigma ,\mathbf{%
b}}f\Big|\sqrt{|I'|_\sigma}\cdot\\
&&
\hspace{1cm}\cdot\bigg(\sum_{J'\in\mathfrak{C}(J)}\left[\c_J\right]^2\Big|\Big((J\backslash I^\prime)\cap(1+\delta)I'\Big)\cap J'  \Big|_\omega\bigg)^\frac{1}{2}\\
&\leq&
C_{\mathbf{b,r},n}\mathfrak{N}_{T^\alpha}||f||_{L^2(\sigma)}\cdot\\
&&
\cdot\bigg(\sum_{I\in \mathcal{D}}\sum_{\substack{J\in\mathcal{N}(I)\\ I\cap J\neq\emptyset}}\sum_{I^{\prime }\in \mathfrak{C}_{{brok}}(I)}\sum_{J'\in\mathfrak{C}(J)}
\!\!\!\!\!\left[\c_J\right]^2\Big|\Big((J\backslash I^\prime)\cap(1+\delta)I'\Big)\cap J'  \Big|_\omega\bigg)^\frac{1}{2}.
\end{eqnarray*}
Now, it is geometrically evident that for the Lebesque measure we have 
$$
\Big|\Big((J\backslash I^\prime)\cap(1+\delta)I'\Big)\cap J'  \Big|\lesssim \delta |J'|.
$$
Taking averages over the grid $\mathcal{D}$ we get the same inequality for the $\omega$ measure:
$$
\boldsymbol{E}_{\Omega }^{\mathcal{D}}\Big|\Big((J\backslash I^\prime)\cap(1+\delta)I'\Big)\cap J'  \Big|_\omega \lesssim \delta \left|J'\right|_\omega.
$$
Thus, if we fix $J'$, there are only finitely many $I'$ involved that contribute (are non-zero), and then the expectation in $\mathcal{D}$ can "go through" the sum in $I'$ to get the estimate
$$
\boldsymbol{E}_{\Omega }^{\mathcal{D}}\sum_{I\in \mathcal{D}}\sum_{\substack{J\in\mathcal{N}(I)\\ I\cap J\neq\emptyset}}\left\vert
\sum_{I^{\prime }\in \mathfrak{C}_{{brok}}\left( I\right) }\left(
E_{I^{\prime }}^{\sigma }\widehat{\mathbb{F}}_{I^{\prime }}^{\sigma ,\mathbf{%
b}}f\right) C(f,g)\right\vert
\leq
C_{\mathbf{b,r},n}\sqrt{\delta}\mathfrak{N}_{T^\alpha}||f||_{L^2(\sigma)}||g||_{L^2(\omega)}.
$$
The constant $C_{\mathbf{b,r},n}$ depends on the accretivity constant of the family $\mathbf{b}$, the dimension $n$ and the \textit{finite} repetition of the intervals $J'$ appearing in the sum.

The third term in (\ref{broken split}) is handled similarly if we change to $\left\langle \square _{I}^{\sigma
,\flat ,\mathbf{b}}f,T_{\omega }^{\alpha ,\ast }\left( \square _{J,{brok}}^{\omega ,\flat ,\mathbf{b}^{\ast }}g\right) \right\rangle _{\sigma
} $, the dual operator. For the second term in \eqref{broken split} the proof is somewhat different: it does not use probability, it is easier because the terms involving $g$ can be estimated as the terms involving $f$ in the proof just done for the fourth term, and then using Carleson estimates. 
So combining the above we get the following 
\begin{eqnarray}\label{flat reduction}
&&
\boldsymbol{E}_{\Omega }^{\mathcal{D}}\sum_{I\in \mathcal{D}}\sum_{\substack{J\in\mathcal{N}(I)\\ (1+\delta)I\cap J\neq\emptyset}} \left\vert \left\langle T_{\sigma }^{\alpha
}\left( \square _{I}^{\sigma ,\mathbf{b}}f\right) ,\square _{J}^{\omega,\mathbf{b}^{\ast }}g\right\rangle _{\omega }\right\vert \\
&\leq&  \notag
\sum_{I\in \mathcal{D}}\sum_{\substack{J\in\mathcal{N}(I)\\ (1+\delta)I\cap J\neq\emptyset}} \left\vert\left\langle T_{\sigma }^{\alpha }\left( \square _{I}^{\sigma ,\flat ,\mathbf{b}}f\right) ,\square _{J}^{\omega ,\flat ,\mathbf{b}^{\ast
}}g\right\rangle _{\omega }\right\vert\\
&&\notag
\hspace{1cm}+\left(C_{\mathbf{b,r},n}\sqrt{\delta}\mathfrak{N}_{T^\alpha}+(\delta^{\alpha-n}+1)\mathcal{NTV}_{\alpha }\right)||f||_{L^2(\sigma)}||g||_{L^2(\omega)}
\end{eqnarray}

Thus it remains to consider the first inner product $\left\langle T_{\sigma
}^{\alpha }\left( \square _{I}^{\sigma ,\flat ,\mathbf{b}}f\right) ,\square
_{J}^{\omega ,\flat ,\mathbf{b}^{\ast }}g\right\rangle _{\omega }$ on the
right hand side of (\ref{flat reduction}), which we call the problematic term,
and write it as%
\begin{eqnarray}\label{PIJ}
\nonumber P\left( I,J\right) &\equiv &\left\langle T_{\sigma }^{\alpha }\left( \square
_{I}^{\sigma ,\flat ,\mathbf{b}}f\right) ,\square _{J}^{\omega ,\flat ,%
\mathbf{b}^{\ast }}g\right\rangle _{\omega }   \\
\nonumber&=&\sum_{I^{\prime }\in \mathfrak{C}\left( I\right), J^{\prime
}\in \mathfrak{C}\left( J\right) }\left\langle T_{\sigma }^{\alpha }\left( 
\mathbf{1}_{I^{\prime }}\square _{I}^{\sigma ,\flat ,\mathbf{b}}f\right) ,%
\mathbf{1}_{J^{\prime }}\square _{J}^{\omega ,\flat ,\mathbf{b}^{\ast
}}g\right\rangle _{\omega }  \notag \\
&=&\sum_{I^{\prime }\in \mathfrak{C}\left( I\right), J^{\prime
}\in \mathfrak{C}\left( J\right) }E_{I^{\prime }}^{\sigma }\left( \widehat{%
\square }_{I}^{\sigma ,\flat ,\mathbf{b}}f\right) \ \left\langle T_{\sigma
}^{\alpha }\left( \mathbf{1}_{I^{\prime }}b_{I}\right) ,\mathbf{1}%
_{J^{\prime }}b_{J}^{\ast }\right\rangle _{\omega }\ E_{J^{\prime }}^{\omega
}\left( \widehat{\square }_{J}^{\omega ,\flat ,\mathbf{b}^{\ast }}g\right).
\end{eqnarray}

It now remains to show that
\begin{equation}
\boldsymbol{E}_{\Omega }^{\mathcal{D}}\boldsymbol{E}_{\Omega }^{\mathcal{G}%
}\sum_{I\in \mathcal{D}}\sum_{ J\in \mathcal{N}(I)}\left\vert P\left( I,J\right) \right\vert
\lesssim \left( C_{\theta }\mathcal{NTV}_{\alpha }+\sqrt{\theta }\mathfrak{N}%
_{T^{\alpha }}\right) \left\Vert f\right\Vert _{L^{2}\left( \sigma \right)
}\left\Vert g\right\Vert _{L^{2}\left( \omega \right) }.
\label{must show final}
\end{equation}%

Suppose now that $I\in \mathcal{C}_{A}$ for $A\in \mathcal{A}$, and that $%
J\in \mathcal{C}_{B}$ for $B\in \mathcal{B}$. Then the inner product in the
third line of (\ref{PIJ}) becomes%
\begin{equation*}
\left\langle T_{\sigma }^{\alpha }\left( b_{I}\mathbf{1}_{I^{\prime
}}\right) ,b_{J}^{\ast }\mathbf{1}_{J^{\prime }}\right\rangle _{\omega
}=\left\langle T_{\sigma }^{\alpha }\left( b_{A}\mathbf{1}_{I^{\prime
}}\right) ,b_{B}^{\ast }\mathbf{1}_{J^{\prime }}\right\rangle _{\omega }\ ,
\end{equation*}%
and we will write this inner product in either form, depending on context.
We also introduce the following notation:%
\begin{equation*}
P_{\left( I,J\right) }\left( E,F\right) \equiv \left\langle T_{\sigma
}^{\alpha }\left( b_{I}\mathbf{1}_{E}\right) ,b_{J}^{\ast }\mathbf{1}%
_{F}\right\rangle _{\omega },\ \ \ \ \ \text{for any sets }E\text{ and }F,
\end{equation*}%
so that%
\begin{equation*}
P\left( I,J\right) =\sum_{I^{\prime }\in \mathfrak{C}\left( I\right) \text{
and }J^{\prime }\in \mathfrak{C}\left( J\right) }E_{I^{\prime }}^{\sigma
}\left( \widehat{\square }_{I}^{\sigma ,\flat ,\mathbf{b}}f\right) \
P_{\left( I,J\right) }\left( I^{\prime },J^{\prime }\right) \ E_{J^{\prime
}}^{\omega }\left( \widehat{\square }_{J}^{\omega ,\flat ,\mathbf{b}^{\ast
}}g\right) .
\end{equation*}%
The first thing we do is reduce matters to showing inequality (\ref{must
show final}) in the case that $P_{\left( I,J\right) }\left( I^{\prime },J^{\prime }\right) $ is replaced by 
$$P_{\left( I,J\right) }\left( I^{\prime }\cap J^{\prime },I^{\prime }\cap
J^{\prime }\right) $$ in the terms $P\left( I,J\right) $ appearing in (\ref%
{must show final}). To see this,
write $\left\langle T_{\sigma }^{\alpha }\left( b_{I}\mathbf{1}%
_{I^{\prime }}\right) ,b_{J}^{\ast }\mathbf{1}_{J^{\prime }}\right\rangle
_{\omega }$ as
$$
\left\langle T_{\sigma }^{\alpha }\left( b_{I}\mathbf{1}%
_{I^{\prime }\backslash J^{\prime }}\right) ,b_{J}^{\ast }\mathbf{1}%
_{J^{\prime }}\right\rangle _{\omega }
+
\left\langle
T_{\sigma }^{\alpha }\left( b_{I}\mathbf{1}_{I^{\prime }\cap J^{\prime }}\right) ,b_{J}^{\ast }\mathbf{%
1}_{J^{\prime }\backslash I^{\prime }}\right\rangle _{\omega } 
+
\left\langle T_{\sigma }^{\alpha }\left( b_{I}\mathbf{1}%
_{I^{\prime }\cap J^{\prime }}\right) ,b_{J}^{\ast }\mathbf{1}_{I^{\prime }\cap J^{\prime }}\right\rangle_\omega 
$$
Set  $$\mathrm{I}=\left\langle T_{\sigma }^{\alpha }\left( b_{I}\mathbf{1}%
_{I^{\prime }\backslash J^{\prime }}\right) ,b_{J}^{\ast }\mathbf{1}%
_{J^{\prime }}\right\rangle _{\omega }$$
$$
\mathrm{II}= \left\langle
T_{\sigma }^{\alpha }\left( b_{I}\mathbf{1}_{I^{\prime }\cap J^{\prime }}\right) ,b_{J}^{\ast }\mathbf{%
1}_{J^{\prime }\backslash I^{\prime }}\right\rangle _{\omega }
\,\,\,\text{and}\,\,\, 
\mathrm{III}=\left\langle T_{\sigma }^{\alpha }\left( b_{I}\mathbf{1}%
_{I^{\prime }\cap J^{\prime }}\right) ,b_{J}^{\ast }\mathbf{1}_{I^{\prime }\cap J^{\prime }}\right\rangle_\omega 
$$
For the first one, we have
$$
\mathrm{I}\leq \left\vert \left\langle T_{\sigma }^{\alpha }\left( b_{I}\mathbf{1}%
_{I^{\prime }\backslash(1+\delta) J^{\prime }}\right) ,b_{J}^{\ast }\mathbf{1}%
_{J^{\prime }}\right\rangle _{\omega }\right\vert 
+
\left\vert \left\langle T_{\sigma }^{\alpha }\left( b_{I}\mathbf{1}%
_{(I^{\prime }\backslash J^{\prime })\cap(1+\delta) J^\prime}\right) ,b_{J}^{\ast }\mathbf{1}%
_{J^{\prime }}\right\rangle _{\omega }\right\vert 
\equiv
\mathrm{I_1+I_2}
$$
Using Lemma \ref{lemma1}, $\mathrm{I}_1\lesssim \delta^{\alpha-n}\sqrt{\mathfrak{A}_2^\alpha}\sqrt{|I'|_\sigma}\sqrt{|J'|_\omega}$ and for $\mathrm{I}_2$ we need to use random surgery. Summing all the terms for $\mathrm{I}_2$ and using Lemma \ref{lemma star}, we have
\begin{eqnarray}\label{reduction to intersection}
&&\ \ \ \ \ \ \boldsymbol{E}_{\Omega }^{\mathcal{G}}\sum_{I\in\mathcal{D}}\sum_{J\in\mathcal{N}(I)}\sum_{I'\in\mathfrak{C}(I)}\sum_{J'\in\mathfrak{C}(J)}\mathfrak{N}_{T^{\alpha}}\left|\a_I\right|\Big(\int_{(I^{\prime }\backslash J^{\prime })\cap(1+\delta) J^\prime}\!\!\!\!\!\!\!\!\!\!\!\!\!\!\!\!\!\!\!|b_I|^2d\sigma\Big)^{\frac{1}{2}}\cdot\\
&&\hspace{6cm}\cdot
\left|\c_J\right| \Big(\int_{J'}|b_J|^2d\omega\Big)^{\frac{1}{2}} \notag \\
&\lesssim&\!\!\!\!
\mathfrak{N}_{T^{\alpha}}\boldsymbol{E}_{\Omega }^{\mathcal{G}}\sum_{I\in\mathcal{D}}\sum_{J\in\mathcal{N}(I)}\sum_{I'\in\mathfrak{C}(I)}\sum_{J'\in\mathfrak{C}(J)}\Big|\a_I\Big|\Big|(I'\backslash J')\cap(1+\delta) J'\Big|_\sigma^{\frac{1}{2}}
\Big|\c_J\Big| \big|J'\big|_\omega^{\frac{1}{2}} \notag \\
&\leq&\!\!\!\!
\mathfrak{N}_{T^{\alpha}}\boldsymbol{E}_{\Omega }^{\mathcal{G}}\Bigg(\!\!\!\sum\!\! \left[\a_I\!\right]^{\!2} \Big|(I'\backslash J')\cap(1+\delta) J'\Big|_\sigma\!\! \Bigg)^{\!\!\frac{1}{2}}\bigg(\!\!\!\sum \left[\!\!\!\c_J\!\right]^{\!2}\!|J'|_\omega\!\!\bigg)^{\!\!\frac{1}{2}} \notag \\
&\leq&\!\!\!\!
\mathfrak{N}_{T^{\alpha}} C_{n,\mathbf{r}}||g||_{L^2(\omega)}\bigg(\sum_I\sum_{I'}\left[\a_I\right]^2\boldsymbol{E}_{\Omega }^{\mathcal{G}}\sum_J\sum_{J'}\Big|(I'\backslash J')\cap(1+\delta) J'\Big|_\sigma\bigg )^{\!\!\!\frac{1}{2}}\notag \\
&\leq&\!\!\!\!
\mathfrak{N}_{T^{\alpha}} C_{n,\mathbf{r}}||g||_{L^2(\omega)}\bigg(\sum_I\sum_{I'}\left[\a_I\right]^2 \delta |I'|_\sigma\bigg )^\frac{1}{2} \notag \\
&\leq&\!\!\!\!
\mathfrak{N}_{T^{\alpha}} C_{n,\mathbf{r}}\sqrt{\delta}||f||_{L^2(\sigma)}||g||_{L^2(\omega)}\notag
\end{eqnarray}
Similarly, we get the bound for $\mathrm{II}$. 

We are left then with $\mathrm{III}$ where we are integrating over $I' \cap J'$. We have to overcome two difficulties at this step. First, $I' \cap J'$ is not necessarily a cube, so we cannot apply any of the testing conditions available. Second, $I' \cap J'$, even if it is a cube, does not need to belong in either of 
the grids $\mathcal{D}$ or $\mathcal{G}$. We would like to split $I' \cap J'$ in smaller cubes of the grid $\mathcal{G}$. The problem is that the boundary of $I'\cap J'$ does not necessarily align with the grid $\mathcal{G}$. To deal with this, we cut a slice around $I' \cap J'$ so that what is left inside can be split in cubes of the grid $\mathcal{G}$. This small slice will be bounded using once again random surgery. While for the remaining cubes, we will use a more involved random surgery technique along with the $\mathcal{A}_2$ and testing condition.

Here are the details: Let $\eta_0=2^{-m}$ for $m$ large enough. For any cube $L$ we define the $\overrightarrow{\eta_1}$-halo for $\overrightarrow{\eta_1}=(\eta_1^1,\dots,\eta_1^n)$ by 
\begin{eqnarray*}
\partial_{\overrightarrow{\eta_1}}L &=& (1+\overrightarrow{\eta_1})L-(1-\overrightarrow{\eta_1})L
\end{eqnarray*}
where $(1+\overrightarrow{\eta_1})L$ means a dilation of each coordinate of $L$ according to the corresponding coordinates of $1+\overrightarrow {\eta_1}$.
Choose the coordinates of $\overrightarrow{\eta_1}$ such that 
$\frac{\eta_0}{2}\leq\eta_1^i<\eta_0$ for all $1\leq i\leq n$ and such that if
\begin{equation}\label{defIJ}
I'\cap J'=\bigg[\Big(I'\backslash\partial_{\overrightarrow{\eta_1}}I'\Big)\cap 
J'\bigg]\overset{\cdot}{\cup}\bigg[\left(\partial_{\overrightarrow{\eta_1}}I'\cap I'\right)\cap 
J'\bigg]\equiv M\overset{\cdot}{\cup} L
\end{equation}
then $M$ consists of $B\lesssim 2^{n\cdot m}$ cubes $K_s\in\mathcal{G}$ with 
$\ell(K_s)\geq 2^{-m-1}\ell(J')$. Note that either $M$ or $L$ might be empty 
depending on where $J'$ is located, but this is not a problem.
Thus 
\begin{eqnarray*}
\left\langle T_{\sigma }^{\alpha }\left( b_{I}\mathbf{1}%
_{I^{\prime }\cap J^{\prime }}\right) ,b_{J}^{\ast }\mathbf{1}_{I^{\prime }\cap J^{\prime }}\right\rangle_\omega 
=\!\!\!\!\!
 &\left\langle T_{\sigma }^{\alpha }\left( b_{I}\mathbf{1}_{M}\right) ,b_{J}^{\ast }\mathbf{1}_{L}\right\rangle_\omega& 
\!\!\!+
\left\langle T_{\sigma }^{\alpha }\left( b_{I}\mathbf{1}%
_{L}\right) ,b_{J}^{\ast }\mathbf{1}_{M}\right\rangle_\omega\\
&& 
\!\!\!\!\!\!\!\!\!\hspace{-1cm}\!\!\!\!\!\!\!\!\!\!\!\!\!
+\left\langle T_{\sigma }^{\alpha }\left( b_{I}\mathbf{1}%
_{L}\right) ,b_{J}^{\ast }\mathbf{1}_{L}\right\rangle_\omega
+
\left\langle T_{\sigma }^{\alpha }\left( b_{I}\mathbf{1}%
_{M}\right) ,b_{J}^{\ast }\mathbf{1}_{M}\right\rangle_\omega
\end{eqnarray*}
The first two can be estimated using Lemma \ref{lemma1} and a random surgery. It is important to mention here that the averages will be taken on the grid $\mathcal{D}$, so that we do not have common intersection among the different translations of the halo. Indeed,
\begin{eqnarray*}
 \left\langle T_{\sigma }^{\alpha }\left( b_{I}\mathbf{1}_{M}\right) ,b_{J}^{\ast }\mathbf{1}_{L}\right\rangle_\omega 
=
 \left\langle T_{\sigma }^{\alpha }\left( b_{I}\mathbf{1}_{M}\right) ,b_{J}^{\ast }\mathbf{1}_{L\backslash (1+\delta)M}\right\rangle_\omega
 +
 \left\langle T_{\sigma }^{\alpha }\left( b_{I}\mathbf{1}%
_{M}\right) ,b_{J}^{\ast }\mathbf{1}_{L\cap(1+\delta)M}\right\rangle_\omega
\equiv
\mathrm{A_1+A_2}
\end{eqnarray*}
and
\begin{eqnarray*}
\left\langle T_{\sigma }^{\alpha }\left( b_{I}\mathbf{1}%
_{L}\right) ,b_{J}^{\ast }\mathbf{1}_{M}\right\rangle_\omega
=
\left\langle T_{\sigma }^{\alpha }\left( b_{I}\mathbf{1}_{L}\right) ,b_{J}^{\ast }\mathbf{1}_{M\backslash (1+\delta)L}\right\rangle_\omega
+
\left\langle T_{\sigma }^{\alpha }\left( b_{I}\mathbf{1}%
_{L}\right) ,b_{J}^{\ast }\mathbf{1}_{M\cap(1+\delta)L}\right\rangle_\omega
\equiv
\mathrm{A_3+A_4}
\end{eqnarray*}
The first terms on the right hand side of both displays, $A_1$ and $A_3$, are bounded, by applying the proof of Lemma \ref{lemma1} for $M$ and $L$ and using the fact that $M$ consists of $B\lesssim 2^{nm}$ cubes. The bound is a constant multiple of $2^n\delta^{\alpha-n}\sqrt{\mathfrak{A}_2^\alpha}\sqrt{|I'|_\sigma}\sqrt{|J'|_\omega}$, which when plugged into the left hand side of (\ref{must show
final}) we get by using Cauchy-Schwarz that
\begin{eqnarray}\label{ML inequalities}
&&\sum_{I\in\mathcal{D}}\sum_{J\in\mathcal{N}(I)}\sum_{\substack{I'\in\mathfrak{C}(I)\\J'\in\mathfrak{C}(J)  }}\!\!\!\bigg|\a_I\bigg|(\mathrm{A_1+A_3})\bigg|\c_J\bigg|\\
&\lesssim&
\sum_{I\in\mathcal{D}}\sum_{J\in\mathcal{N}(I)}\sum_{\substack{I'\in\mathfrak{C}(I)\\J'\in\mathfrak{C}(J)     }}\!\!\!\Big|\a_I\Big| \delta^{\alpha-n}\sqrt{\mathfrak{A}_2^\alpha}\sqrt{|I'|_\sigma}\sqrt{|J'|_\omega}\bigg|\!\!\c_J\bigg|\notag\\
&\lesssim&
\delta^{\alpha-n}\sqrt{\mathfrak{A}_2^\alpha}||f||_{L^2(\sigma)}||g||_{L^2(\omega)}\notag
\end{eqnarray}
For $\mathrm{A_2}$ (and similarly for $\mathrm{A_4}$), we have
\begin{eqnarray}\label{LM surgery}
&&
\quad \boldsymbol{E}^\mathcal{D}_\Omega\sum_{I\in\mathcal{D}}\sum_{J\in\mathcal{N}(I)}\sum_{\substack{I'\in\mathfrak{C}(I), J'\in\mathfrak{C}(J) }}\bigg|\a_I\bigg|\cdot 
\bigg|\bigg|T_\sigma^\alpha(b_I\mathbf{1}_M)\bigg|\bigg|_{L^2(\omega)}\cdot \\
&&
\hspace{4 cm}\cdot \bigg|\bigg|b_J^*\mathbf{1}_{L\cap(1+\delta)M})\bigg|\bigg|_{L^2(\omega)}
\bigg|\c_J\bigg|\notag\\
&\leq&
\mathfrak{N}_{T^\alpha}C_{\mathbf{b}}\boldsymbol{E}^\mathcal{D}_\Omega  \!\!\!\!\!\!\!
\sum_{\substack{I'\in\mathfrak{C}(I)\& J'\in\mathfrak{C}(J)\\J\in \mathcal{N}(I)} }\!\!\!\!\!\!\!\!\Big|\a_I\Big|\Big|M\Big|_\sigma^\frac{1}{2}\Big|L\cap(1+\delta)M\Big|_\omega^\frac{1}{2} \bigg|\c_J\bigg|\notag\\
&\leq&
\mathfrak{N}_{T^\alpha}C_{\mathbf{b,b^*,r},n}\left(
\sum_{\substack{I'\in\mathfrak{C}(I) \& J'\in\mathfrak{C}(J)\\J\in \mathcal{N}(I)}} \bigg|\a_I\bigg|^2|M|_\sigma\right)^{\!\!\!\frac{1}{2}}\cdot\notag\\
&&\hspace{3cm}\cdot
\bigg(\boldsymbol{E}^\mathcal{D}_\Omega\!\!\!\!
\sum_{\substack{I'\in\mathfrak{C}(I)\& J'\in\mathfrak{C}(J)\\J\in \mathcal{N}(I)}}\!\! \bigg|\c_J\bigg|^2|L\cap(1+\delta)M|_\omega \bigg)^{\!\!\frac{1}{2}}\notag\\
&\leq&
\mathfrak{N}_{T^\alpha}C_{\mathbf{b,b^*,r},n}\sqrt{\delta}||f||_{L^2(\sigma)} ||g||_{L^2(\omega)}\notag
\end{eqnarray}
by noting that $(1+\delta)M \cap L$ is a halo of width $\delta$, much smaller than $\eta_0$ (so as to get the estimate by $\sqrt{\delta}$, not $\sqrt{\eta_0}$). Although an estimate of $\sqrt{\eta_0}$ is easy to obtain (as $L$ already has width $\eta_0$) and is sufficient for the purposes of this term, the estimate of $\sqrt{\delta}$ will be crucially used later in \eqref{KK surgery} to kill the $B$ term.
Note also that we can take the averages over all directions, so that we avoid common intersection along the different translations. Notice that $L,M$ are "moving" together. This is not a problem since by "moving" they cover different parts of the cube $J'$. 

Thus we only need to estimate  
$
\left\langle T_{\sigma }^{\alpha }\left( b_{I}\mathbf{1}%
_{L}\right) ,b_{J}^{\ast }\mathbf{1}_{L}\right\rangle_\omega
+
\left\langle T_{\sigma }^{\alpha }\left( b_{I}\mathbf{1}%
_{M}\right) ,b_{J}^{\ast }\mathbf{1}_{M}\right\rangle_\omega
$.
Applying one more time random surgery to the first term we get that 
\begin{eqnarray}\label{LL surgery}
&&
\boldsymbol{E}_{\Omega }^{\mathcal{D}}\boldsymbol{E}_{\Omega }^{\mathcal{G}%
}\sum_{I\in \mathcal{D}}\sum_{ J\in \mathcal{N}(I)}\sum_{\substack{ I^{\prime }\in \mathfrak{C}%
\left( I\right)  \\ J^{\prime }\in \mathfrak{C}\left( J\right) }} \left\vert
E_{I^{\prime }}^{\sigma }\left( \widehat{\square }_{I}^{\sigma ,\flat ,%
\mathbf{b}}f\right) \ \left\langle T_{\sigma }^{\alpha }\left( b_{I}\mathbf{1%
}_{L}\right) ,b_{J}^{\ast }\mathbf{1}_{L}\right\rangle _{\omega } E_{J^{\prime }}^{\omega }\left( \widehat{\square }_{J}^{\omega ,\flat ,%
\mathbf{b}^{\ast }}g\right) \right\vert \notag \\
&\lesssim &
\boldsymbol{E}_{\Omega }^{\mathcal{G}}\mathfrak{N}_{T^{\alpha
}}\left\Vert f\right\Vert _{L^{2}\left( \sigma \right) }\boldsymbol{E}%
_{\Omega }^{\mathcal{D}}\sqrt{\sum_{I\in \mathcal{D}}\sum_{J\in\mathcal{N}(I)
}\!\sum_{\substack{ %
I^{\prime }\in \mathfrak{C}\left( I\right)  \\ J^{\prime }\in \mathfrak{C}%
\left( J\right) }}\!\!\!\left( \int_{\partial _{\eta_{1}}I^{\prime }\cap
J^{\prime }}\!\!\!\left\vert b_{J}^{\ast }\right\vert ^{2}d\omega \right)
\!\!\left\vert E_{J^{\prime }}^{\omega }\left( \widehat{\square }_{J}^{\omega
,\flat ,\mathbf{b}^{\ast }}g\right) \right\vert ^{2}}\notag
\end{eqnarray}%
using (\ref{PLBP removed}) and the frame inequalities again.
Then using Cauchy-Schwarz on the expectation $\boldsymbol{E}_{\Omega }^{%
\mathcal{D}}$, this is dominated by 

\begin{eqnarray*}
& &\!\!\!
\boldsymbol{E}_{\Omega }^{\mathcal{G}}\mathfrak{N}_{T^{\alpha
}}\!\left\Vert f\right\Vert _{L^{2}\left( \sigma \right) } \!\!\!\!\sqrt{\sum_{J\in 
\mathcal{G}}\sum_{J^{\prime }\in \mathfrak{C}\left( J\right) }\left( 
\boldsymbol{E}_{\Omega }^{\mathcal{D}}\!\!\!\!\!\!\!\!\!\!\!\!\sum_{\substack{ I\in \mathcal{D}:\
2^{-\mathbf{r}n}|I| <|J| \leq |I|  \\ d\left( J,I\right) \leq 2\ell(J) ^{\varepsilon
}\ell(I) ^{1-\varepsilon }  \\ I^{\prime }\in \mathfrak{C}%
\left( I\right) }}\!\!\!\!\!\!\!\!\!\!\!\!
\left\vert \partial _{\overrightarrow{\eta_{1}}}I^{\prime }\cap J^{\prime
}\right\vert _{\omega }\!\!\right) \!\!\left\vert E_{J^{\prime }}^{\omega }\left( 
\widehat{\square }_{J}^{\omega ,\flat ,\mathbf{b}^{\ast }}g\right)
\right\vert ^{2}} \\
\end{eqnarray*}%

\begin{eqnarray*}
&\lesssim &
\!\!\!\boldsymbol{E}_{\Omega }^{\mathcal{G}}\mathfrak{N}_{T^{\alpha
}}\left\Vert f\right\Vert _{L^{2}\left( \sigma \right) } \!\!\sqrt{\sum_{J\in 
\mathcal{G}}\sum_{J^{\prime }\in \mathfrak{C}\left( J\right) }\!\!\!2^{\mathbf{r}}\!\!\left( \boldsymbol{E}_{\Omega }^{\mathcal{D}}\!\!\!\!\!\!\!\!\!\!\!\!\sum_{I^{\prime }\in \mathcal{D}:|J'| \leq |I'| \leq 2^{\mathbf{r}}|J'| }\!\!\!\!\!\!\!\!\!\!\!\!\left\vert \partial _{\eta
_{0}}I^{\prime }\cap J^{\prime }\right\vert _{\omega }\right)\!\! \left\vert
E_{J^{\prime }}^{\omega }\left( \widehat{\square }_{J}^{\omega ,\flat ,%
\mathbf{b}^{\ast }}g\right) \right\vert ^{2}} \\
&\lesssim &
\!\!\!\sqrt{\eta_{0}}\mathfrak{N}_{T^{\alpha }}\left\Vert f\right\Vert
_{L^{2}\left( \sigma \right) }\left\Vert g\right\Vert _{L^{2}\left( \omega
\right) }
\leq
\sqrt{\lambda}\mathfrak{N}_{T^{\alpha }}\left\Vert f\right\Vert
_{L^{2}\left( \sigma \right) }\left\Vert g\right\Vert _{L^{2}\left( \omega
\right) }
\end{eqnarray*}%
where in the last line we have used $\eta_{1}^i\leq \eta_{0}$, and then 
\begin{equation*}
\boldsymbol{E}_{\Omega }^{\mathcal{D}}\sum_{I^{\prime }\in \mathcal{D}:|J'|\leq |I'| \leq 2^{%
\mathbf{r}}|J'| }\left\vert \partial _{\eta_{0}}I^{\prime }\cap J^{\prime }\right\vert _{\omega }\lesssim \eta_{0}\left\vert J^{\prime }\right\vert _{\omega }
\end{equation*}%
as long as we choose $\eta_{0} \ll 2^{-\mathbf{r}%
}$. 

This leaves us to estimate the term $\left\langle T_{\sigma }^{\alpha }\left( b_{I}\mathbf{1}%
_{M}\right) ,b_{J}^{\ast }\mathbf{1}_{M}\right\rangle_\omega$.  It is at this point
that we will use the decomposition $M=\displaystyle\bigcup _{1\leq
s\leq B}^\cdot K_{s}$ constructed above. 
We have
\begin{eqnarray*} 
\left\langle T_{\sigma }^{\alpha }\left( b_{I}\mathbf{1}_{M}\right)
,b_{J}^{\ast }\mathbf{1}_{M}\right\rangle _{\omega
}  &=&  \sum\limits_{s,s^{\prime }=1}^{B}\left\langle T_{\sigma }^{\alpha }\left(
b_{I}\mathbf{1}_{K_{s}}\right) ,b_{J}^{\ast }\mathbf{1}_{K_{s^{\prime
}}}\right\rangle _{\omega }
\end{eqnarray*}
which can be rewritten as 

\begin{equation}\label{kk}
\sum\limits_{s=1}^{B}\left\langle T_{\sigma
}^{\alpha }\left( b_{I}\mathbf{1}_{K_{s}}\right), b_{J}^{\ast }\mathbf{1}
_{K_{s}}\right\rangle _{\omega }
+
\bigg(\sum\!\!\!\!\!\!\!\sum_{K_s\underset{Sep}{\sim} K_{s'}}
+
\sum\!\!\!\!\!\!\!\sum_{K_s\underset{Adj}{\sim} K_{s'}}\bigg)\left\langle T_{\sigma }^{\alpha }\left( b_{I}\mathbf{1}_{K_{s}}\right),b_{J}^{\ast }\mathbf{1}
_{K_{s^{\prime }}}\right\rangle _{\omega }
\end{equation}
where we call $K_s\underset{Sep}{\sim} K_{s'}$ the separated cubes, i.e. $3K_s\cap K_{s'}=\emptyset$, while  by $K_s\underset{Adj}{\sim} K_{s'}$ are the adjacent cubes, i.e. $K_s\cap K_{s'}=\emptyset$ and $\overline{K_s}\cap \overline{K_{s'}}\neq \emptyset$. 
The separated terms sum can be estimated directly by $\sqrt{\mathfrak{A}_2^\alpha}$. Indeed, as in the proof of Lemma \ref{lemma1},
\begin{eqnarray*}
\left\langle T_{\sigma }^{\alpha }\left( b_{I}\mathbf{1}_{K_{s}}\right),b_{J}^{\ast }\mathbf{1}
_{K_{s^{\prime }}}\right\rangle _{\omega }
\!\!\!\!&\lesssim&\!\!\!\!
\left( \int_{K_{s'}}\left( \int_{K_s}\left\vert x-y\right\vert ^{\alpha -n}\left\vert b_{I}\left( y\right) \right\vert d\sigma \left(
y\right) \right) ^{2}d\omega \left( x\right) \right) ^{\frac{1}{2}}\!\!\!\!\sqrt{\left\vert K_{s'}\right\vert _{\omega }} \\
&\lesssim&\!\!\!\!
\left( \int_{\mathbb{R}^n\backslash K_s}\left\vert x-x_{K_s}\right\vert ^{2\alpha -2n} d\omega (x) \right)^{\frac{1}{2}}|K_s|_\sigma \sqrt{%
\left\vert K_{s'}\right\vert _{\omega }}\\
&\lesssim&\!\!\!\!
\sqrt{\mathfrak{A}_2^\alpha}\sqrt{|K_s|_\sigma}\sqrt{|K_{s'}|_\omega}
\end{eqnarray*}
thus,
\begin{eqnarray}\label{KK disjoint}
\quad \quad \quad \sum\!\!\!\!\!\!\!\!\sum_{K_s\underset{Sep}{\sim} K_{s'}}\left\langle T_{\sigma }^{\alpha }\left( b_{I}\mathbf{1}_{K_{s}}\right),b_{J}^{\ast }\mathbf{1}
_{K_{s^{\prime }}}\right\rangle _{\omega }
\leq
C_{\mathbf{b}}\sum\!\!\!\!\!\!\!\!\sum_{K_s\underset{Sep}{\sim} K_{s'}}\!\!\!\!\!\!\!\ \sqrt{\mathfrak{A}_2^\alpha}\sqrt{|K_s|_\sigma}\sqrt{|K_{s'}|_\omega}
\end{eqnarray}
which plugged into \eqref{PIJ} appropriately, we get the bound $B\sqrt{\mathfrak{A}_2^\alpha}\sqrt{|I'|_\sigma}\sqrt{|J'|_\omega}$. 

To deal with the adjacent cubes term in \eqref{kk}, we write  
\begin{eqnarray*}
&&
\sum\!\!\!\!\!\!\!\sum_{K_s\underset{Adj}{\sim} K_{s'}}\left\langle T_{\sigma }^{\alpha }\left( b_{I}\mathbf{1}_{K_{s}}\right),b_{J}^{\ast }\mathbf{1}
_{K_{s^{\prime }}}\right\rangle _{\omega }=\sum\!\!\!\!\!\!\!\sum_{K_s\underset{Adj}{\sim} K_{s'}}\left\langle  b_{I}\mathbf{1}_{K_{s}},T_{\omega }^{\alpha,\ast }\left(b_{J}^{\ast }\mathbf{1}
_{K_{s^{\prime }}}\right)\right\rangle _{\sigma }\\
&\!\!\!\!\!\!\!\!\!\!\!\!\!\!\!\!\!\!\!\!\!\!\!\!\!\!\!\!\!\!\!\!\!=&\!\!\!\!\!\!\!\!\!\!\!\!\!\!\!\!
\sum\!\!\!\!\!\!\!\sum_{K_s\underset{Adj}{\sim} K_{s'}}\left\langle  b_{I}\mathbf{1}_{K_{s}\cap(1+\delta)K_{s'}},T_{\omega }^{\alpha,\ast }\left(b_{J}^{\ast }\mathbf{1}
_{K_{s'}}\right)\right\rangle _{\sigma }
+
\sum\!\!\!\!\!\!\!\sum_{K_s\underset{Adj}{\sim} K_{s'}}\left\langle  b_{I}\mathbf{1}_{K_{s}\backslash(1+\delta)K_{s'}},T_{\omega }^{\alpha,\ast }\left(b_{J}^{\ast }\mathbf{1}
_{K_{s'}}\right)\right\rangle _{\sigma }\\
&\!\!\!\!\!\!\!\!\!\!\!\!\!\!\!\!\!\!\!\!\!\!\!\!\!\!\!\!\!\!\!\!\!\equiv&\!\!\!\!\!\!\!\!\!\!\!\!\!\!\!\!\!\!
\mathrm{\overset{\sim}{I}+\overset{\sim}{II}}
\end{eqnarray*}
For $\mathrm{\overset{\sim}{II}}$ we use Lemma \ref{lemma1} to get
\begin{eqnarray}\label{KK adjacent}
\mathrm{\overset{\sim}{II}}   &\lesssim &
\delta^{\alpha-n}\sqrt{\mathfrak{A}_2^\alpha}\left(\sum\limits_{s=1}^{B}|K_s|_\sigma\right)^\frac{1}{2}\left(\sum\limits_{s=1}^{B}\left(\sum_{s'\geq s}\Big|K_{s'}\Big|_\omega^\frac{1}{2}\right)^2\right)^\frac{1}{2}\\
&\lesssim&
\delta^{\alpha-n} B \sqrt{\mathfrak{A}_2^\alpha} \sqrt{|I'|_\sigma}\sqrt{|J'|_\omega}\notag
\end{eqnarray}
while summing $\mathrm{\overset{\sim}{I}}$ over 
$$
\mathcal{T}=\{I\in\mathcal{D}, J\in\mathcal{N}(I), I'\in\mathfrak{C}_{nat}(I), J'\in\mathfrak{C}_{nat}(J)\}
$$
and using Cauchy-Schwarz, accretivity, taking averages and using Jensen, we get
\begin{eqnarray}\label{KK surgery}
\\
&&\!\!\!\!\!\!
\boldsymbol{E}^\mathcal{G}_\Omega \sum_{\mathcal{T}}\!\bigg|\a_I\! \c_J\bigg|\! \sum\!\!\!\!\!\!\!\sum_{K_s\underset{Adj}{\sim} K_{s'}}\!\!\!\!\!\left\langle  b_{I}\mathbf{1}_{K_{s}\cap(1+\delta)K_{s'}},T_{\omega }^{\alpha,\ast }\left(b_{J}^{\ast }\mathbf{1}
_{K_{s'}}\right)\right\rangle _{\sigma }\notag\\
&\lesssim&\!\!\!\!
\boldsymbol{E}^\mathcal{G}_\Omega \sum_{\mathcal{T}}\bigg|\a_I \c_J\bigg| \sum\!\!\!\!\!\!\!\sum_{K_s\underset{Adj}{\sim} K_{s'}}\!\!\!\!\!\!
\mathfrak{N}_{T^\alpha}\sqrt{|K_{s}\cap(1+\delta)K_{s'}|_\sigma}\sqrt{|K_{s'}|_\omega}\notag\\
&\lesssim&\!\!\!\!
\mathfrak{N}_{T^\alpha}\boldsymbol{E}^\mathcal{G}_\Omega \sum_{\mathcal{T}}\bigg|\a_I \c_J\bigg| \Big(\sum_{s=1}^B|K_{s'}|_\omega\Big)^\frac{1}{2}
\bigg(\sum_{s=1}^B \Big(\sum_{s\leq s'}
\sqrt{|K_{s}\cap(1+\delta)K_{s'}|_\sigma}\Big)^2\bigg)^\frac{1}{2}\notag\\
&\lesssim&\!\!\!\!
\mathfrak{N}_{T^\alpha}\boldsymbol{E}^\mathcal{G}_\Omega \sum_{\mathcal{T}}\bigg|\a_I \c_J\bigg| \sqrt{|J'|_\omega}
\bigg(\sum_{s=1}^B \Big(\sum_{s\leq s'}
|K_{s}\cap(1+\delta)K_{s'}|_\sigma\cdot \sum_{s\leq s'} 1\Big)\bigg)^\frac{1}{2}\notag\\
&\lesssim&\!\!\!\!
\mathfrak{N}_{T^\alpha}\sqrt{B}\left\Vert g\right\Vert _{L^{2}\left( \omega
\right) }\!\! \left(\!\sum_{\underset{  I'\in\mathfrak{C}_{nat}(I)}{I\in \mathcal{D}}}\!\!\!\!\!\!\bigg|\a_I\bigg|^2 \boldsymbol{E}^\mathcal{G}_\Omega\!\!\!\!\! \sum_{\underset{J'\in\mathfrak{C}_{nat}(J)}{J\in\mathcal{N}(I)}}\!\! \sum_{s=1}^B \sum_{s\leq s'}
|K_{s}\!\cap\!(1\!+\!\delta)K_{s'}|_\sigma\!\!\right)^{\!\!\!\frac{1}{2}}\notag\\
&\lesssim&\!\!\!\!
\mathfrak{N}_{T^\alpha}\sqrt{B}\left\Vert g\right\Vert _{L^{2}\left( \omega
\right) } \left(\sum_{\mathcal{T}}\bigg|\a_I\bigg|^2 2^{n}\delta|I'|_\sigma\right)^\frac{1}{2}\notag\\
&\lesssim&
\mathfrak{N}_{T^\alpha}2^{2n}\sqrt{B}\sqrt{\delta}\left\Vert f\right\Vert
_{L^{2}\left( \sigma \right) }\left\Vert g\right\Vert _{L^{2}\left( \omega
\right) } \notag
\end{eqnarray}
because there are up to $2^n$ adjacent cubes $K_{s'}$ for a given $K_s$. The implied constant depends on $\mathbf{r}$ of the nearby form. Note that $\delta$ is independent of $B$ or $\mathbf{r}$ and will later be chosen small enough so that the terms containing the norm inequality constant will be absorbed.

Thus now we are left only with the first term of (\ref{kk}), i.e. we need to estimate 
$$
\sum\limits_{s=1}^{B}\left\langle T_{\sigma
}^{\alpha }\left( b_{I}\mathbf{1}_{K_{s}}\right), b_{J}^{\ast }\mathbf{1}%
_{K_{s}}\right\rangle _{\omega }
$$

Before proceeding further it will prove convenient to introduce some
additional notation, namely we will write the energy estimate in the second
display of the Energy Lemma as%
\begin{equation}
\left\vert \left\langle T^{\alpha }\nu ,\Psi _{J}\right\rangle _{\omega
}\right\vert \lesssim C_{\gamma,\delta }\ \mathrm{P}_{\delta }^{\alpha }\mathsf{Q}%
^{\omega }\left( J,\upsilon \right) \ \left\Vert \Psi _{J}\right\Vert
_{L^{2}\left( \mu \right) }\ \ \ \text{if }\int \Psi
_{J}d\omega =0\text{ and }\gamma J\cap {\supp}\nu =\emptyset
\label{star}
\end{equation}
where
\begin{equation*}
\mathrm{P}_{\delta }^{\alpha }\mathsf{Q}^{\omega }\left( J,\upsilon \right)
\equiv \frac{\mathrm{P}^{\alpha }\left( J,\nu \right) }{\left\vert
J\right\vert }\left\Vert \mathsf{Q}_{J}^{\omega ,\mathbf{b}^{\ast
}}x\right\Vert _{L^{2}\left( \omega \right) }^{\spadesuit }+\frac{\mathrm{P}%
_{1+\delta }^{\alpha }\left( J,\nu \right) }{\left\vert J\right\vert }%
\left\Vert x-m_{J}\right\Vert _{L^{2}\left( \mathbf{1}_{J}\omega \right) }\ .
\label{def compact}
\end{equation*}%
The use of the compact notation $\mathrm{P}_{\delta }^{\alpha }\mathsf{Q}%
^{\omega }\left( J,\upsilon \right) $ to denote the complicated expression
on the right hand side will considerably reduce the size of many subsequent
displays.

We now consider the inner product $\left\langle T_{\sigma }^{\alpha }\left(
b_{A}\mathbf{1}_{K}\right) ,b_{B}^{\ast }\mathbf{1}_{K}\right\rangle
_{\omega }$ and estimate the case when%
\begin{equation*}
K\in \mathcal{G},\ K\subset I^{\prime }\cap J^{\prime },\ I^{\prime }\in 
\mathfrak{C}\left( I\right) ,\ J^{\prime }\in \mathfrak{C}\left( J\right) ,\
I\in \mathcal{C}_{A}^{\mathcal{A}},\ J\in \mathcal{C}_{B}^{\mathcal{B}},\
\ell(K) =2^{-m-1}\ell(J') \mathfrak{.}
\end{equation*}
 For subsets $E,F\subset A\cap B$ and cubes $%
K\subset A\cap B$ we define%
\begin{equation}
\left\{ E,F\right\} \equiv \left\langle T_{\sigma }^{\alpha }\left( b_{A}\mathbf{1}_{E}\right) ,b_{B}^{\ast }\mathbf{1}_{F}\right\rangle _{\omega }\ ,
\label{def E,F}
\end{equation}
and $K_{in}$ the $2^n$ grandchildren of $K$ that do not intersect the boundary of $K$ while $K_{out}$ the rest $4^n-2^n$ grandchildren of $K$ that intersect its boundary i.e.
\begin{equation*}
K_{{in}}= \left\{
K^{\prime \prime }\in \mathfrak{C}^{\left( 2\right) }\left( K\right)
:\partial K^{\prime \prime }\cap \partial K=\emptyset \right\}
\end{equation*}
\begin{equation*}
K_{{out}}= \left\{ K^{\prime \prime }\in \mathfrak{C}%
^{\left( 2\right) }\left( K\right) :\partial K^{\prime \prime }\cap \partial
K\neq \emptyset \right\}
\end{equation*}
We can write 
\begin{equation}
\left\{ K,K\right\} =\left\{ A,K_{in}\right\} -\left\{ A\backslash
K,K_{in}\right\} +\left\{ K_{out},K_{out}\right\} +\left\{ K_{in},K_{out}\right\} .  \label{K,K}
\end{equation}
Note that the first two terms on the right hand side of (\ref{K,K})
decompose the inner product $\left\{ K,K_{in}\right\} $, which
`includes' one of the difficult symmetric inner product $\left\{ K_{
in},K_{in}\right\} $, and where the other difficult symmetric
inner products are contained in $\left\{ K_{out},K_{out}\right\} $, which can be handled recursively. Thus the difficult symmetric
inner products are ultimately controlled by testing on the cube $A$ to
handle the `paraproduct' term$\ \left\{ A,K_{in}\right\} $, and by
using the energy condition and a trick that resurrects the original testing
functions $\left\{ b_{J}^{\ast ,orig}\right\} _{J\in \mathcal{G}}$%
, discarded in the corona constructions above, to handle the `stopping' term 
$\left\{ A\backslash K,K_{in}\right\} $. More precisely, these
original testing functions $b_{J}^{\ast ,orig}$ are the testing
functions obtained after reducing matters to the case of bounded testing
functions.

The first term on the right side of (\ref{K,K}) satisfies%
\begin{eqnarray}
\,\,\,\,\,\,\,\left\vert \left\{ A,K_{in}\right\} \right\vert 
&=&
\left\vert
\int_{K_{in}}\left( T_{\sigma }^{\alpha }b_{A}\right) b_{B}^{\ast
}d\omega \right\vert \leq \left\Vert \mathbf{1}_{K_{in}}T_{\sigma
}^{\alpha }b_{A}\right\Vert _{L^{2}\left( \omega \right) }\left\Vert \mathbf{%
1}_{K_{in}}b_{B}^{\ast }\right\Vert _{L^{2}\left( \omega \right)}\\
&\leq&
 \left\Vert b_{B}^{\ast }\right\Vert _{\infty }\left\Vert \mathbf{1}_{K_{in}}T_{\sigma }^{\alpha }b_{A}\right\Vert _{L^{2}\left(\omega \right) }\sqrt{\left\vert K_{in}\right\vert _{\omega }}\ \nonumber.
\label{AKin}
\end{eqnarray}

We now turn to the term $\left\{ A\backslash K,K_{in}\right\}$. Decompose $\mathbf{1}_{K_{in}}b_{B}^{\ast }$ as%
\begin{equation*}
\mathbf{1}_{K_{in}}b_{B}^{\ast }=\sum_{\ell =1}^{2^n}\mathbf{1}%
_{K^\ell_{in}}\left( b_{B}^{\ast }-\frac{1}{\left\vert
K^\ell_{in}\right\vert _{\omega }}
\int_{K^\ell_{in}}b_{B}^{\ast }d\omega \right) +\sum_{\ell =1}^{2^n}\mathbf{1}_{K^\ell_{in}}\frac{1}{\left\vert K^\ell_{in}\right\vert
_{\omega }}\int_{K^\ell_{in}}b_{B}^{\ast }d\omega ,
\end{equation*}%
and then apply the Energy Lemma to the function%
\begin{equation*}
k_{K_{in}}^{\ast }\equiv \sum_{\ell =1}^{2^n}\mathbf{1}_{K^\ell_{in}}\left( b_{B}^{\ast }-\frac{1}{\left\vert K^\ell_{in}\right\vert _{\omega }}\int_{K^\ell_{in}}b_{B}^{\ast }d\omega \right) \equiv \sum_{j=1}^{2^n} k_{K_{in}}^{\ast ,j}
\end{equation*}%
which does indeed satisfy $\square _{K^{\prime }}^{\omega ,\mathbf{b}^{\ast
}}k_{K_{in}}^{\ast }=0$ unless $K^{\prime }$ is a dyadic
subcube of $K$ that is contained in $K_{in}$. (Furthermore, we
could even replace grandchildren by $m$-grandchildren in this argument in
order that $\square _{K^{\prime }}^{\omega ,\mathbf{b}^{\ast }}k_{K_{in}}^{\ast }=0$ unless $K^{\prime }$ is a dyadic $m$-grandchild of 
$K$ that is contained in $K_{in}$, but we will not need this.) We obtain
\begin{eqnarray}
\left\langle T_{\sigma }^{\alpha }\left( b_{A}\mathbf{1}_{A\backslash
K}\right) ,\mathbf{1}_{K_{in}}b_{B}^{\ast }\right\rangle _{\omega} 
\!\!\!\!&=&
\left\langle T_{\sigma }^{\alpha }\left( b_{A}\mathbf{1}_{A\backslash
K}\right) ,k_{K_{in}}^{\ast }\right\rangle _{\omega }
\label{reach'''} \notag \\
\!\!\!\!& &
+\left\langle T_{\sigma }^{\alpha }\left( b_{A}\mathbf{1}_{A\backslash
K}\right) ,\sum_{\ell =1}^{2^n}\mathbf{1}_{K^\ell_{in}}\left( 
\frac{1}{\left\vert K^\ell_{in}\right\vert _{\omega }}%
\int_{K^\ell_{in}}b_{B}^{\ast }d\omega \right) \right\rangle
_{\omega }  
\end{eqnarray}
and%
\begin{eqnarray}\label{after mon'''}
\left\vert \left\langle T_{\sigma }^{\alpha }\left( b_{A}\mathbf{1}%
_{A\backslash K}\right) ,k_{K_{in}}^{\ast }\right\rangle _{\omega
}\right\vert
&\leq&
\sum_{\ell=1}^{2^n}\left\vert \left\langle T_{\sigma }^{\alpha }\left( b_{A}\mathbf{1}_{A\backslash K}\right) ,k_{K_{in}}^{\ast
,\ell}\right\rangle _{\omega }\right\vert \notag \\
&\leq&
 C_{\eta_0,n }\left[ \sum_{\ell=1}^{2^n} \mathrm{P}_{\delta }^{\alpha }\mathsf{Q}^{\omega }\left(
K^\ell_{in},\mathbf{1}_{A\backslash K}\sigma \right)
\right] \left\Vert k_{K_{in}}^{\ast }\right\Vert _{L^{2}\left(
\omega \right) }\  \notag
\end{eqnarray}%
where the constant $C_{\eta_0 }$ depends on the constant $C_{\gamma }$ in the
statement of the Monotonicity Lemma with $\gamma =\frac{1}{1-\eta_0 }$ since $
\frac{1}{1-\eta_0 }K_{in}\cap \left( A\backslash K\right) =\emptyset $
, and where we have written $\left\{ K^\ell_{in}\right\}
_{\ell =1}^{2^n} $ with $K^{\ell}_{in}$ denoting the innner
grandchildren of $K$.

Thus we see that $\mathsf{P}_{\mathcal{H}}^{\omega ,\mathbf{b}^{\ast }}$ and 
$\mathsf{Q}_{\mathcal{H}}^{\omega ,\mathbf{b}^{\ast }}$ in the Energy Lemma
can be taken to be pseudoprojection onto $K_{in}$, i.e. $\displaystyle\mathsf{P}%
_{K_{in}}^{\omega ,\mathbf{b}^{\ast }}=\!\!\!\!\sum_{J\in \mathcal{G}:\
J\subset K_{in}}\!\!\!\!\square _{J}^{\omega ,\mathbf{b}^{\ast }}$ and $%
\displaystyle\mathsf{Q}_{K_{in}}^{\omega ,\mathbf{b}^{\ast }}=\!\!\!\!\sum_{J\in 
\mathcal{G}:\ J\subset K_{in}}\!\!\!\!\bigtriangleup _{J}^{\omega ,\mathbf{%
b}^{\ast }}$, and we will see below that the cubes $K_{in}$
that arise in subsequent arguments will be pairwise disjoint. Furthermore,
the energy condition will be used to control these full pseudoprojections $\mathsf{P}_{K_{in}}^{\omega,\mathbf{b}^{\ast }}$ when taken over
pairwise disjoint decompositions of cubes by subcubes of the form $K_{in}$.

However, the second line of (\ref{reach'''}) remains problematic because we cannot use any type of testing in $K^\ell_{in}$ with $b^\ast_B$ since $K^\ell_{in}$ does not necessarily belong to $\mathcal{C}_B$, and this
is our point in which we exploit the original testing functions $b_{K^\ell_{in}}^{\ast ,orig}$.

\subsubsection{Return to the original testing functions}

From the discussion above, we recall the identity (\ref{reach'''}) and the
estimate (\ref{after mon'''}). We also have the analogous identity and
estimate with $b_{K^\ell_{in}}^{\ast ,orig}$ in
place of $\mathbf{1}_{K_{in}}b_{B}^{\ast }$:\begin{eqnarray}
&&\!\!\!\!\!\!\!\!\!\!\!\!\!\!\!\!\!\!\!\!\!\!\!\!\!\!\!\!\!\!\!\!\!
\left\langle T_{\sigma }^{\alpha }\left( b_{A}\mathbf{1}_{A\backslash
K}\right) ,b_{K^\ell_{in}}^{\ast ,orig}\right\rangle _{\omega } \notag \\ 
&=&
\left\langle T_{\sigma }^{\alpha }\left( b_{A}\mathbf{1}_{A\backslash K}\right) ,\mathbf{1}_{K^{\ell }_{in
}}\left( b_{K^\ell_{in}}^{\ast ,orig}-\frac{1}{\left\vert K^\ell_{in}\right\vert _{\omega }}\int_{K^{\ell
}_{in }}b_{K^\ell_{in}}^{\ast ,orig}d\omega \right) \right\rangle _{\omega }  \label{reach''''}  \\
&&
\hspace{1cm}+
\left\langle T_{\sigma }^{\alpha }\left( b_{A}\mathbf{1}_{A\backslash
K}\right) ,\mathbf{1}_{K^\ell_{in}}\left( \frac{1}{\left\vert
K^\ell_{in}\right\vert _{\omega }}\int_{K^{\ell }_{in }}b_{K^\ell_{in}}^{\ast ,orig}d\omega
\right) \right\rangle _{\omega }\   \notag
\end{eqnarray}
and
\begin{eqnarray}
&&
\left\vert \left\langle T_{\sigma }^{\alpha }\left( b_{A}\mathbf{1}%
_{A\backslash K}\right) ,\mathbf{1}_{K^\ell_{in}}\left(
b_{K^\ell_{in}}^{\ast ,orig}-\frac{1}{\left\vert
K^\ell_{in}\right\vert _{\omega }}\int_{K^{\ell }_{in }}b_{K^\ell_{in}}^{\ast ,orig}d\omega
\right) \right\rangle _{\omega }\right\vert  \label{after mon''''} \\
&\lesssim &
\mathrm{P}_{\delta }^{\alpha }\mathsf{Q}^{\omega }\left( K^\ell
_{in},\mathbf{1}_{A\backslash K}\sigma \right) \left\Vert 
\mathbf{1}_{K^\ell_{in}}\left( b_{K^{\ell }_{in
}}^{\ast ,orig}-\frac{1}{\left\vert K^{\ell }_{in}\right\vert _{\omega }}\int_{K^\ell_{in}}b_{K^{\ell}_{in}}^{\ast ,orig}d\omega \right) \right\Vert
_{L^{2}\left( \omega \right) } \notag
\end{eqnarray}%
for $1\leq \ell \leq 2^n$, where the implied constants depend on $L^{\infty }$
norms of testing functions and the constant in the Energy Lemma. Using the notation
\begin{equation*}
\left\{ K_{out},K_{in}^{\ell }\right\} ^{orig}\equiv \left\langle T_{\sigma }^{\alpha }b_{A}\mathbf{1}_{K_{out}},b_{K^\ell_{in}}^{\ast ,orig}\right\rangle
_{\omega }\ \text{for } 1\leq \ell \leq 2^n .
\end{equation*}
note that 
\begin{eqnarray*}
&&\left\{ A\backslash K,K_{in}\right\}+ \sum_{\ell=1}^{2^n}\Bigg(\Dl\Bigg) \left\{ K_{out},K_{in}^{\ell }\right\} ^{orig}\\
&=&
\left\{ A\backslash K,K_{in}\right\}- \sum_{\ell=1}^{2^n}\Bigg(\Dl\Bigg) \left\langle T_{\sigma }^{\alpha }\left( b_{A}\mathbf{1}_{A\backslash K}\right) ,b_{K^\ell_{in}}^{\ast ,orig%
}\right\rangle _{\omega } \\
& & 
+ \sum_{\ell=1}^{2^n}\Bigg(\Dl\Bigg)\!\bigg[\!\!\left\langle T_{\sigma }^{\alpha }\left( b_{A}\mathbf{1}_{A}\right) ,b_{K^\ell_{in}}^{\ast ,orig}\right\rangle _{\omega } 
\!\!\!-\!\!
\Big\langle b_A\mathbf{1}_{K_{in}}, T^{\alpha,\ast}_\omega b^{\ast,orig}_{K^\ell_{in}}\Big\rangle_{\!\!\sigma}
\!\bigg]\\
&\equiv&
\mathbf{B}+\mathbf{C}
\end{eqnarray*}
Now for $\mathbf{B}$, using Energy Lemma to the function
$$
\Psi_J^\ell=\Bigg(\Dl\Bigg)b^{\ast,orig}_{K^\ell_{in}}- \left( 
\frac{1}{\left\vert K^\ell_{in}\right\vert _{\omega }}%
\int_{K^\ell_{in}}b_{B}^{\ast }d\omega \right) \mathbf{1}_{K^\ell_{in}}
$$ 
for $1\leq \ell\leq 2^n$ we have 
\begin{eqnarray*}
|\mathbf{B}|
\!\!\!\!\!&=&\!\!\!\!\!
\bigg|\left\langle T_{\sigma }^{\alpha }\left( b_{A}\mathbf{1}_{A\backslash
K}\right) ,\mathbf{1}_{K_{in}}b_{B}^{\ast }\right\rangle _{\omega} \!\! - \!\!
\sum_{\ell=1}^{2^n}  \left( 
\frac{1}{\left\vert K^\ell_{in}\right\vert _{\!\!\omega }}
\int_{K^\ell_{in}}b_{B}^{\ast }d\omega \right)\!\!\Big\langle T_{\sigma }^{\alpha }\left( b_{A}\mathbf{1}_{A\backslash K}\right) ,\mathbf{1}_{K^\ell_{in}}\Big\rangle_\omega \bigg|\\
& &
+O\bigg[ \sum_{\ell=1}^{2^n} \bigg(\frac{\mathrm{P}^\alpha(K^\ell_{in}\mathbf{1}_{A\backslash K}\sigma)}{|K^\ell_{in}|}\Big|\Big|\mathsf{Q}%
_{K^\ell_{in}}^{\omega ,\mathbf{b}^{\ast }}x\Big|\Big|^\spadesuit_{L^2(\omega)}\bigg)\bigg]\sqrt{|K_{in}|_\omega}\\
& &
+O\bigg[ \sum_{\ell=1}^{2^n} \bigg(\frac{\mathrm{P}^\alpha_{1+\delta}(K^\ell_{in}\mathbf{1}_{A\backslash K}\sigma)}{|K^\ell_{in}|}\Big|\Big|x-m_{K^\ell_{in}}\Big|\Big|_{L^2(\omega)}\bigg)\bigg]\sqrt{|K_{in}|_\omega}\\
&\lesssim &
\left[ \sum_{\ell=1}^{2^n} \mathrm{P}_{\delta }^{\alpha }\mathsf{Q}^{\omega }\left(
K^\ell_{in},\mathbf{1}_{A\backslash K}\sigma \right) \right] \sqrt{%
\left\vert K_{in}\right\vert _{\omega }}\ 
\end{eqnarray*}
having used the triangle inequality to get
$$
\big|\big|\Psi_J^\ell\big|\big|_{L^2(\omega)}\lesssim \Bigg|
\frac{\int_{K^\ell_{in}}b^\ast_Bd\omega}{\int_{K^\ell_{in}}b^{\ast,orig}_Bd\omega}\Bigg|\sqrt{|K^\ell_{in}|_\omega}+\sqrt{|K^\ell_{in}|_\omega}\lesssim \sqrt{|K_{in}| _{\omega }}, \quad 1\leq\ell\leq2^n
$$
and
\begin{eqnarray*}
&&\!\!\!\!\!\!\!\!\!\!
\bigg|\left\langle T_{\sigma }^{\alpha }\left( b_{A}\mathbf{1}_{A\backslash
K}\right) ,\mathbf{1}_{K_{in}}b_{B}^{\ast }\right\rangle _{\omega}-
\sum_{\ell=1}^{2^n}  \left( 
\frac{1}{\left\vert K^\ell_{in}\right\vert _{\omega }}%
\int_{K^\ell_{in}}b_{B}^{\ast }d\omega \right)\Big\langle T_{\sigma }^{\alpha }\left( b_{A}\mathbf{1}_{A\backslash K}\right) ,\mathbf{1}_{K^\ell_{in}}\Big\rangle_\omega \bigg|\\
&\lesssim&
\left[\sum_{\ell=1}^{2^n} \mathrm{P}_{\delta }^{\alpha }\mathsf{Q}^{\omega }\left(
K^\ell_{in},\mathbf{1}_{A\backslash K}\sigma \right) \right] 
\left|\left|1_{K_\ell''}\sum_{\ell=1}^{2^n}\left(b^*_B-\frac{1}{|K_l''|_\omega}\int_{k_\ell''}b^*_Bd\omega\right)\right|\right|_{L^2(\omega)}
\\
&\lesssim &
\left[\sum_{\ell=1}^{2^n} \mathrm{P}_{\delta }^{\alpha }\mathsf{Q}^{\omega }\left(
K^\ell_{in},\mathbf{1}_{A\backslash K}\sigma \right)  \right] \sqrt{%
\left\vert K_{in}\right\vert _{\omega }}\
\end{eqnarray*}
where in the last inequality we used accretivity and triangle inequality. 
We turn our attention in term $\mathbf{C}$. We have that
\begin{eqnarray*}
&&
\Bigg|\sum_{\ell=1}^{2^n}\Bigg(\Dl\Bigg)\left\langle T_{\sigma }^{\alpha }\left( b_{A}\mathbf{1}_{A}\right) ,b_{K^\ell_{in}}^{\ast ,orig}\right\rangle _{\omega } \Bigg|\\
&\lesssim&
\sum_{\ell=1}^{2^n}\sqrt{\int_{K_{\ell }^{\prime \prime
}}\left\vert T_{\sigma }^{\alpha }b_{A}\right\vert ^{2}d\omega }\sqrt{%
\int_{K^\ell_{in}}\left\vert b_{K_{\ell }^{\prime \prime
}}^{\ast ,orig}\right\vert ^{2}d\omega }\\
&\lesssim& 
\sqrt{\int_{K_{in}}\left\vert T_{\sigma }^{\alpha }b_{A}\right\vert ^{2}d\omega }%
\sqrt{\left\vert K_{in}\right\vert _{\omega }}
\end{eqnarray*}
Also,
$$
\Bigg|\sum_{\ell=1}^{2^n}\Bigg(\Dl\Bigg)\Big\langle b_A\mathbf{1}_{K_{in}}, T^{\alpha,\ast}_\omega b^{\ast,orig}_{K^\ell_{in}}\Big\rangle_\sigma\Bigg|\equiv \mathrm{I+II+III}
$$
where 
$$
\mathrm{I}=\sum_{\ell=1}^{2^n}\Bigg(\Dl\Bigg)\Big\langle b_A\mathbf{1}_{K_{in}},\mathbf{1}_{K^\ell_{in}} T^{\alpha,\ast}_\omega b^{\ast,orig}_{K^\ell_{in}}\Big\rangle_\sigma 
$$
$$
\mathrm{II}= \sum_{\ell=1}^{2^n}\Bigg(\Dl\Bigg)\Big\langle b_A\mathbf{1}_{K_{in}\backslash(1+\delta)K^\ell_{in}}, T^{\alpha,\ast}_\omega b^{\ast,orig}_{K^\ell_{in}}\Big\rangle_\sigma 
$$
$$
\mathrm{III}= \sum_{\ell=1}^{2^n}\Bigg(\Dl\Bigg)\Big\langle b_A\mathbf{1}_{(K_{in}\backslash K^\ell_{in})\cap(1+\delta)K^\ell_{in}}, T^{\alpha,\ast}_\omega b^{\ast,orig}_{K^\ell_{in}}\Big\rangle_\sigma 
$$
The first term $\mathrm{I}$ is bounded using the dual testing condition. Indeed,
$$
\mathrm{I}\leq \big|\big|b_A\mathbf{1}_{K_{in}}\big|\big|_{L^2(\sigma)}\sum_{\ell=1}^{2^n}\mathfrak{T}^*C_{\mathbf{b}^*}\sqrt{|K^\ell_{in}|_\omega}
\leq
2^n\mathfrak{T}^*C_{\mathbf{b}^*}\big|\big|b_A\mathbf{1}_{K_{in}}\big|\big|_{L^2(\sigma)}\sqrt{|K_{in}|_\omega}
$$
The second term $\mathrm{II}$ is bounded using Lemma \ref{lemma1}. Indeed,
\begin{eqnarray*}
\mathrm{II} &\leq& \sum_{\ell=1}^{2^n}\delta^{\alpha-n}\sqrt{\mathfrak{A}_2^\alpha}\sqrt{|K_{in}\backslash(1+\delta)K^\ell_{in}|_\sigma} \sqrt{|K^\ell_{in}|_\omega}\\
&\leq&
2^n \delta^{\alpha-n}\sqrt{\mathfrak{A}_2^\alpha}\sqrt{|K_{in}|_\omega}\sqrt{|K_{in}|_\sigma}
\end{eqnarray*}
Finally, 
\begin{eqnarray*}
\mathrm{III}   &\leq&  \sum_{\ell=1}^{2^n} \Big|\Big|T^\alpha_\sigma\big(b_A\mathbf{1}_{(K_{in}\backslash K^\ell_{in})\cap(1+\delta)K^\ell_{in}}\big)\Big|\Big|_{L^2(\omega)} \Big|\Big|b^{\ast,orig}_{K^\ell_{in}}\Big|\Big|_{L^2(\omega)}\\
&\leq&
\mathfrak{N}_{T^\alpha}\sqrt{C_{\mathbf{b}}C_{\mathbf{b}^*}}\bigg(\sum_{\ell=1}^{2^n}\Big|(K_{in}\backslash K^\ell_{in})\cap(1+\delta)K^\ell_{in}\Big|_\sigma   \bigg)^\frac{1}{2}  \sqrt{|K_{in}|_\omega}\\
&\equiv&
\sqrt{C_{\mathbf{b}}C_{\mathbf{b}^*}}\cdot\mathbf{\Delta}(K)
\end{eqnarray*}
where we have defined 
$$
\mathbf{\Delta}\left(K\right)=\mathfrak{N}_{T^\alpha}\bigg(\sum_{\ell=1}^{2^n}\Big|(K_{in}\backslash K^\ell_{in})\cap(1+\delta)K^\ell_{in}\Big|_\sigma   \bigg)^\frac{1}{2}  \sqrt{|K_{in}|_\omega}
$$
This last term will be iterated and a final random surgery will give us the desired bound.

\subsubsection{A finite iteration and a final random surgery.}
Letting 
\begin{eqnarray}\label{def PHI}
\Phi^{A,B}(K_{in})  &=&  
\Big|\Big|\mathbf{1}_{K_{in}}T^\alpha_\sigma\big(b_A\big)\Big|\Big|_{L^2(\omega)} \sqrt{\left\vert K_{in}\right\vert _{\omega }}\\
& & +
\sum_{\ell=1}^{2^n}\mathrm{P}_{\delta }^{\alpha }\mathsf{Q}^{\omega }\left(
K^{\ell}_{in},\mathbf{1}_{A\backslash K}\sigma \right)\sqrt{\left\vert K_{in}\right\vert _{\omega }}\notag\\
& & +
\Big(\mathfrak{T}^\alpha+\mathfrak{T}^{\alpha,\ast}+\delta^{\alpha-n}\sqrt{\mathfrak{A}_2^\alpha}\Big)\sqrt{\left\vert K_{in}\right\vert _{\sigma }}\sqrt{\left\vert K_{in}\right\vert _{\omega }}\notag
\end{eqnarray}
and simplifying more our notation 
$$
\left\{ K_{out},K_{in}\right\} ^{orig}\equiv \sum_{\ell=1}^{2^n}\Bigg(\Dl\Bigg) \left\{ K_{out},K_{in}^{\ell }\right\} ^{orig}
$$we have so far that (\ref{K,K}) is written as 
\begin{equation*}
\left\{ K,K\right\} =
\left\{ K_{out},K_{in}\right\} ^{orig}+\left\{ K_{out},K_{out}\right\} +\left\{ K_{in},K_{out}\right\}+
O\big(\Phi^{A,B}(K_{in})+\mathbf{\Delta}(K)\big) \\
\end{equation*}
Now 
$$
\left\{ K_{out},K_{out}\right\}=\sum_{\ell} \{K^{\ell}_{out},K^{\ell}_{out}\}
+\!\!\!\!
\sum_{\substack{m\neq\ell\\K^{m}_{out}\cap K^{\ell}_{out}=\emptyset}}\!\!\!\!\!\!\!\! \{K^{\ell}_{out},K^{m}_{out}\} 
+\!\!\!\!\!\!\!\!  
\sum_{\substack{m\neq\ell\\K^{m}_{out}\cap K^{\ell}_{out}\neq\emptyset}}\!\!\!\!\!\!\!\! \{K^{\ell}_{out},K^{m}_{out}\}
$$
where $K_{out}^\ell, 1\leq\ell\leq 4^n-2^n$, are the outer grandchildren of $K$. For the second sum above, we get 
\begin{eqnarray*}
\bigg|\!\!\!\!\!\!\!\! \sum_{\substack{m\neq\ell\\K^{m}_{out}\cap K^{\ell}_{out}=\emptyset}}\!\!\!\!\!\!\!\! \{K^{\ell}_{out},K^{m}_{out}\}\bigg| 
&\lesssim&
\sqrt{\mathfrak{A}_2^\alpha} \sum_\ell \sqrt{|K^{\ell}_{out}|_\sigma}\!\!\!\!\!\!\!\!\!\!  \sum_{\substack{m\neq\ell\\K^{m}_{out}\cap K^{\ell}_{out}=\emptyset}} \!\!\!\!\!\!\!\! \sqrt{|K^{m}_{out}|_\omega}\\
&\lesssim&
\sqrt{\mathfrak{A}_2^\alpha} \sqrt{|K_{out}|_\sigma} \sqrt{|K_{out}|_\omega}
\end{eqnarray*}
where the implied constant depends on dimension and the accretivity of functions involved and since $\dist(K^{\ell}_{out},K^{m}_{out})\geq \ell(K^\ell_{out})$
there is no $\delta$. For the third sum, we need to use random surgery again. Using Lemma \ref{lemma1},
\begin{eqnarray*}
&&
\!\!\!\!|\{K^{\ell}_{out},K^{m}_{out}\}| =
\bigg|\left\langle T_{\sigma }^{\alpha }\left( b_{A}\mathbf{1}_{K^{\ell}_{out}}\right) ,\mathbf{1}_{K^{m}_{out}}b_{B}^{\ast }\right\rangle _{\omega}\bigg|\\
&\leq&
\!\!\!\!\bigg|\!\!\left\langle T_{\sigma }^{\alpha }\left( b_{A}\mathbf{1}_{K^{\ell}_{out}\backslash(1+\delta)K^{m}_{out}}\right) ,\mathbf{1}_{K^{m}_{out}}b_{B}^{\ast }\right\rangle _{\!\!\omega}\!\!\bigg| +
\bigg|\!\!\left\langle T_{\sigma }^{\alpha }\left( b_{A}\mathbf{1}_{K^{\ell}_{out}\cap(1+\delta)K^{m}_{out}}\right),\mathbf{1}_{K^{m}_{out}}b_{B}^{\ast }\right\rangle _{\!\!\omega}\!\!\bigg|\\
&\leq&
\!\!\!\!\delta^{\alpha-n}\sqrt{\mathfrak{A}_2^\alpha}\sqrt{|K^{\ell}_{out}|_\sigma}\sqrt{|K^{m}_{out}|_\omega} +
\mathfrak{N}_{T^\alpha}\sqrt{|K^{m}_{out}|_\omega}\sqrt{|K^{\ell}_{out}\cap(1+\delta)K^{m}_{out}|_\sigma}
\end{eqnarray*}
Thus, summing 
\begin{eqnarray}
&&\qquad 
\sum_\ell\!\! \!\!\!\sum_{\substack{m\neq\ell\\K^{m}_{out}\cap K^{\ell}_{out}\neq\emptyset}}\!\!\!\!\!\!\!\! |\{K^{\ell}_{out},K^{m}_{out}\}|\\
&\lesssim&\notag
\!\!\!\!\delta^{\alpha-n}\sqrt{\mathfrak{A}_2^\alpha}\sqrt{|K_{out}|_\sigma}\sqrt{|K_{out}|_\omega} +
\mathfrak{N}_{T^\alpha} \!\!\!\sum_\ell\!\! \!\!\!\!\!\!\!\!\!\sum_{\substack{m\neq\ell\\K^{m}_{out}\cap K^{\ell}_{out}\neq\emptyset}}\!\!\!\!\!\!\!\! \!\!\!\sqrt{|K^{m}_{out}|_\omega}\sqrt{|K^{\ell}_{out}\cap(1+\delta)K^{m}_{out}|_\sigma}\nonumber \\
&\leq&\notag
\!\!\!\!\!\!\delta^{\alpha-n}\sqrt{\mathfrak{A}_2^\alpha}\sqrt{|K_{out}|_\sigma}\sqrt{|K_{out}|_\omega} +
\mathfrak{N}_{T^\alpha} \!\!\!\sum_\ell\bigg(\sum_{m\neq\ell}|K^{\ell}_{out}\cap(1+\delta)K^{m}_{out}|_\sigma\bigg)^{\!\!\frac{1}{2}}\!\!\!\sqrt{|K_{out}|_\omega} \nonumber
\end{eqnarray}
Let 
$$
\mathbf{E}(K)= \mathfrak{N}_{T^\alpha}\sum_\ell\bigg(\sum_{m\neq\ell}|K^{\ell}_{out}\cap(1+\delta)K^{m}_{out}|_\sigma\bigg)^\frac{1}{2}\sqrt{|K_{out}|_\omega}
$$
We will iterate this term below and we will the necessary bound. We now turn to $\left\{K_{in},K_{out}\right\}$ and we have
\begin{eqnarray*}
&&
|\left\{K_{in},K_{out}\right\}| \\
&\leq &
\bigg|\left\langle T_{\sigma }^{\alpha }\left( b_{A}\mathbf{1}_{K_{out}\backslash(1+\delta)K_{in}}\right) ,\mathbf{1}_{K_{in}}b_{B}^{\ast }\right\rangle _{\!\omega}\!\bigg| +
\bigg|\left\langle T_{\sigma }^{\alpha }\left( b_{A}\mathbf{1}_{K_{out}\cap(1+\delta)K_{in}}\right) ,\mathbf{1}_{K_{in}}b_{B}^{\ast }\right\rangle _{\!\omega}\!\bigg|\\
&\lesssim&
\delta^{\alpha-n}\sqrt{\mathfrak{A}_2^\alpha}\sqrt{|K_{out}|_\sigma}\sqrt{|K_{in}|_\omega} +
\mathfrak{N}_{T^\alpha}\sqrt{|K_{in}|_\omega} \sqrt{|K_{out}\cap(1+\delta)K_{in}|_\sigma}
\end{eqnarray*}
and similarly $|\left\{ K_{out},K_{in}\right\} ^{orig}|$ is bounded by 
\begin{eqnarray*}
&\lesssim& 
\delta^{\alpha-n}\sqrt{\mathfrak{A}_2^\alpha}\sqrt{|K_{out}|_\sigma}\sqrt{|K_{in}|_\omega} +
\mathfrak{N}_{T^\alpha}\sqrt{|K_{in}|_\omega} \sqrt{|K_{out}\cap(1+\delta)K_{in}|_\sigma}
\end{eqnarray*}
Let 
$$
\mathbf{F}(K)=\mathfrak{N}_{T^\alpha}\sqrt{ |K_{out}\cap(1+\delta)K_{in}|_\sigma}\sqrt{|K_{in}|_\omega}
$$
Using the bounds we found above we have from (\ref{K,K}),
\begin{eqnarray*}
|\left\{ K,K\right\}|&\lesssim& \sum_{\ell=1}^{4^n-2^n} |\{K^{\ell}_{out},K^{\ell}_{out}\}|+O\big(\Phi^{A,B}(K_{in})\big)\\
& &+ \mathbf{\Delta}(K)+\mathbf{E}(K)+\mathbf{F}(K)+C_{\delta,\eta_0,\mathbf{b},\mathbf{b}^*}\sqrt{\mathfrak{A}_2^\alpha} \sqrt{|K|_\sigma}\sqrt{|K|_\omega}
\end{eqnarray*}
Iterating the first term above a \textit{finite} number of times, using again the norm inequality and a final random surgery we get the bound we need. Indeed, for $\nu\in\mathbb{N}$
\begin{eqnarray}
|\left\{ K,K\right\}|
&\leq&\!\!\!\!\!\!\!
\sum_{M\in \mathcal{M}_{\nu}}\!\!\!\!|\left\{ M,M\right\}|
+O\left( \sum_{M\in \mathcal{M}_{\nu}^{\ast }}\!\!\!\left[ \Phi ^{A,B}\left( M_{{in}}\right) \right] +\mathbf{\Delta}(M)+\mathbf{E}(M)+\mathbf{F}(M)\! \right) \notag \\
&&+
C_{\delta,\eta_0,\mathbf{b},\mathbf{b}^*}\sqrt{\mathfrak{A}_{2}^{\alpha }}\sum_{M\in \mathcal{M}_{\nu}^{\ast }}\sqrt{\left\vert M\right\vert _{\sigma }}%
\sqrt{\left\vert M\right\vert _{\omega }}  \notag \\
&\equiv &
A\left( K\right) +B\left( K\right) +C\left( K\right) =A_{\left(I^{\prime },J^{\prime }\right) }\left( K\right) +B_{\left( I^{\prime
},J^{\prime }\right) }\left( K\right) +C_{\left( I^{\prime },J^{\prime
}\right) }\left( K\right) ,  \label{K iterated}
\end{eqnarray}
where the collections of cubes $\mathcal{M}_{\nu}=\mathcal{M}_{\nu}\left(
K\right) $ and $\mathcal{M}_{\nu}^{\ast }=\mathcal{M}_{\nu}^{\ast }\left(
K\right) $ are defined recursively by 
\begin{eqnarray*}
\mathcal{M}_{0} 
&\equiv &
\left\{ K\right\} , \\
\mathcal{M}_{k+1} 
&\equiv &
\bigcup_{M\in \mathcal{M}_{k}} \left\{ M^\ell_{out}\right\} ,\ \ \ \ \ k\geq 0, \\
\mathcal{M}_{\nu}^{\ast } &\equiv &\bigcup_{k=0}^{\nu}\mathcal{M}_{k}\ .
\end{eqnarray*}%
We will include the subscript $\left( I^{\prime },J^{\prime }\right) $ in
the notation when we want to indicate the pair $\left( I^{\prime },J^{\prime
}\right) $ that are defined after \eqref{defIJ}. Now the term $C\left( K\right) $
can be estimated by  
\begin{equation}
C\left( K\right) =C_{\delta,\eta_0,\mathbf{b},\mathbf{b}^*}\sqrt{\mathfrak{A}_{2}^{\alpha }}\sum_{M\in \mathcal{M}_{\nu}^{\ast }}\sqrt{\left\vert M\right\vert _{\sigma }}\sqrt{\left\vert
M\right\vert _{\omega }}\leq \nu C_{\delta,\eta_0,\mathbf{b},\mathbf{b}^*}\sqrt{\mathfrak{A}_{2}^{\alpha }}\sqrt{%
\left\vert K\right\vert _{\sigma }}\sqrt{\left\vert K\right\vert _{\omega }}%
  \label{C est}
\end{equation}%
where $\nu$ is chosen below depending on $\eta_{0}$. For the first term $%
A\left( K\right) $, we will apply the norm inequality and use probability,
namely%
\begin{eqnarray*}
\left\vert A\left( K\right) \right\vert &\leq &\sqrt{C_{\mathbf{b}}C_{%
\mathbf{b}^{\ast }}}\mathfrak{N}_{T^{\alpha }}\sum_{M\in \mathcal{M}_{\nu}}%
\sqrt{\left\vert M\right\vert _{\sigma }}\sqrt{\left\vert M\right\vert
_{\omega }} \\
&\leq &\sqrt{C_{\mathbf{b}}C_{\mathbf{b}^{\ast }}}\mathfrak{N}_{T^{\alpha }}%
\sqrt{\sum_{M\in \mathcal{M}_{\nu}}\left\vert M\right\vert _{\sigma }}\sqrt{%
\sum_{M\in \mathcal{M}_{\nu}}\left\vert M\right\vert _{\omega }} \\ &\leq& \sqrt{C_{%
\mathbf{b}}C_{\mathbf{b}^{\ast }}}\mathfrak{N}_{T^{\alpha }}\sqrt{\sum_{M\in 
\mathcal{M}_{\nu}}\left\vert M\right\vert _{\sigma }}\sqrt{\left\vert
K\right\vert _{\omega }},
\end{eqnarray*}%
where $\sqrt{C_{\mathbf{b}}C_{\mathbf{b}^{\ast }}}$ is an upper bound for
the testing functions involved, followed by 
\begin{equation*}
\boldsymbol{E}_{\Omega }^{\mathcal{G}}\left( \sum_{M\in \mathcal{M}_{\nu}}\left\vert M\right\vert _{\sigma }\right) \leq \varepsilon \left\vert
I^{\prime }\right\vert _{\sigma }\ ,
\end{equation*}%
for a sufficiently small $\varepsilon >0$, where \emph{roughly speaking}, we
use the fact that the cubes $M\in \mathcal{M}_{\nu}$ depend on the grid $%
\mathcal{G}$ and form a relatively small proportion of $I^{\prime }$, which
captures only a small amount of the total mass $\left\vert I^{\prime
}\right\vert _{\sigma }$ as the grid is translated relative to the grid $%
\mathcal{D}$ that contains $I^{\prime }$. 

Here are the details. Recall that the cubes $K$ are taken from
the set of consecutive cubes $\left\{ K_{i}\right\} _{i=1}^{B}$ that lie
in $I^{\prime }\cap J^{\prime }$, that the cubes $M\in \mathcal{M}_{\nu}\left( K_{i}\right) $ have length $\frac{1}{4^{n\nu}}\ell \left(
K_{i}\right) $, and that there are $(4^n-2^{n})^\nu$ such cubes in $\mathcal{M}
_{\nu}\left( K_{i}\right) $ for each $i$. Thus we have 
\begin{eqnarray*}
\sum_{M\in \mathcal{M}_{\nu}\left( K\right) }\left\vert M\right\vert
=\sum_{M\in \mathcal{M}_{\nu}\left( K\right) }\frac{1}{4^{n\nu}}\left\vert
K\right\vert =\left(4^n-2^{n}\right)^\nu\frac{1}{4^{n\nu}}\left\vert K\right\vert   \label{future prob} 
\end{eqnarray*}%
and $\displaystyle\left(\frac{4^n-2^n}{4^n}\right)^\nu \rightarrow 0$ as $\nu \rightarrow \infty$, which implies
$$
 \boldsymbol{E}_{\Omega }^{\mathcal{G}}\left(
\sum_{i=1}^{B}\sum_{M\in \mathcal{M}_{\nu}\left( K_{i}\right) }\left\vert
M\right\vert _{\sigma }\right) 
\leq 
B\left(\frac{4^n-2^n}{4^n}\right)^\nu%
\left\vert I^{\prime }\right\vert _{\sigma }\leq \varepsilon\left\vert
I^{\prime }\right\vert _{\sigma }
$$
where we have used that the variable $B$ is at most $2^{nm}$
 and where the final inequality holds if $\nu$ is chosen large enough such that $%
B\left(\frac{4^n-2^n}{4^n}\right)^\nu\leq \varepsilon$. Then we have by Cauchy-Schwarz applied
first to $\displaystyle \sum_{i=1}^{B}\sum_{M\in \mathcal{M}_{\nu}\left( K_{i}\right) }$ and
then to $\boldsymbol{E}_{\Omega }^{\mathcal{G}}$, 
\begin{eqnarray}
&&
\boldsymbol{E}_{\Omega }^{\mathcal{G}}\left( \sum_{i=1}^{B}\left\vert
A\left( K_{i}\right) \right\vert \right) \leq \boldsymbol{E}_{\Omega }^{%
\mathcal{G}}\sqrt{C_{\mathbf{b}}C_{\mathbf{b}^{\ast }}}\mathfrak{N}%
_{T^{\alpha }}\sqrt{\sum_{i=1}^{B}\sum_{M\in \mathcal{M}_{\nu}\left(
K_{i}\right) }\left\vert M\right\vert _{\sigma }}\sqrt{\left\vert J^{\prime
}\right\vert _{\omega }}  \label{A est} \\
&\leq &
\sqrt{C_{\mathbf{b}}C_{\mathbf{b}^{\ast }}}\mathfrak{N}_{T^{\alpha }}%
\sqrt{\boldsymbol{E}_{\Omega }^{\mathcal{G}}\sum_{i=1}^{B}\sum_{M\in 
\mathcal{M}_{\nu}\left( K_{i}\right) }\left\vert M\right\vert _{\sigma }}\sqrt{%
\left\vert J^{\prime }\right\vert _{\omega }}  \notag \\
&\leq &
\sqrt{C_{\mathbf{b}}C_{\mathbf{b}^{\ast }}}\mathfrak{N}_{T^{\alpha }}%
\sqrt{\varepsilon\left\vert I^{\prime }\right\vert _{\sigma }}\sqrt{\left\vert
J^{\prime }\right\vert _{\omega }}
=
\sqrt{C_{\mathbf{b}}C_{\mathbf{b}^{\ast }}%
}\sqrt{\varepsilon}\mathfrak{N}_{T^{\alpha }}\sqrt{\left\vert I^{\prime
}\right\vert _{\sigma }}\sqrt{\left\vert J^{\prime }\right\vert _{\omega }},
\notag
\end{eqnarray}%
as required.

Now we turn to summing up the remaining terms 
\\$\displaystyle B\left( K\right) =C\sum_{M\in 
\mathcal{M}_{\nu}^{\ast }}\Phi ^{A,B}\left( M_{{in}}\right)+\mathbf{\Delta}(M)+\mathbf{E}(M)+\mathbf{F}(M) $ above.
In the case when the cube $I^{\prime }$ is a 
\emph{natural} child of $I$, i.e. $I^{\prime }\in \mathfrak{C}_{{%
nat}}\left( I\right) $ so that $I^{\prime }\in \mathcal{C}_{A}^{\mathcal{%
A}}$, we have 
\begin{eqnarray*}
\sum_{M\in \mathcal{M}_{\nu}^{\ast }\left( K\right) }\left\Vert \mathbf{1}_{M_{{in}}}T_{\sigma }^{\alpha }b_{A}\right\Vert _{L^{2}\left( \omega\right) }^{2}
=
\sum_{M\in \mathcal{M}_{\nu}^{\ast }\left( K\right) }\int_{M_{
{in}}}\left\vert T_{\sigma }^{\alpha }b_{A}\right\vert ^{2}d\omega
\leq
 \int_{I^{\prime }}\left\vert T_{\sigma }^{\alpha}b_{A}\right\vert^{2}d\omega \lesssim \left( \mathfrak{T}_{T^{\alpha }}^{\mathbf{b}}\right)
^{2}\left\vert I^{\prime }\right\vert _{\sigma }
\end{eqnarray*}
by the weak testing condition for $I^{\prime }$ in the corona $\mathcal{C}_{A}$. Also,
\begin{equation*}
\sum_{M\in \mathcal{M}_{\nu}^{\ast }\left( K\right) }\left\vert M_{{in}%
}\right\vert _{\omega }\leq \left\vert K\right\vert _{\omega }\leq
\left\vert J^{\prime }\right\vert _{\omega }
\end{equation*}
because of the crucial fact that the cubes $%
\left\{ M_{{in}}\right\} _{M\in \mathcal{M}_{\nu}^{\ast }\left(
K\right) }$ form a pairwise disjoint subdecomposition of $K\subset I^{\prime
}\cap J^{\prime }$ (for any $\nu\geq 1$). Of course, this implies
\begin{eqnarray*}
\left( \sum_{M\in \mathcal{M}_{\nu}^{\ast }\left( K\right) }\!\!\!\!\!\left( \mathfrak{T}_{T^{\alpha ,\ast }}+\mathfrak{A}_{2}^{\alpha }\right)
^{2}\left\vert M_{{in}}\right\vert _{\sigma }\right) ^{\!\!\!\frac{1}{2}}
\!\!\!\left( \sum_{M\in \mathcal{M}_{\nu}^{\ast }\left( K\right) }\left\vert M_{
{in}}\right\vert _{\omega }\right) ^{\!\!\!\frac{1}{2}}  \label{next claim}
\!\!\!\!\!\!&\lesssim &\!\!\!\!\!
\Big( \mathfrak{T}_{T^{\alpha ,\ast }}+
\mathfrak{A}_{2}^{\alpha }\Big)\sqrt{\left\vert I^{\prime }\right\vert _{\sigma }\left\vert J^{\prime
}\right\vert _{\omega }}
\end{eqnarray*}
and using the definition of $\mathrm{P}_{\delta }^{\alpha }\mathsf{Q}%
^{\omega }\left( J,\upsilon \right) $ in (\ref{def compact}), 
\begin{eqnarray*}
&&
\sum_{M\in \mathcal{M}_{\nu}^{\ast }\left( K\right) }\sum_{\ell=1}^{2^n}\mathrm{P}_{\delta }^{\alpha }\mathsf{Q}^{\omega }\left(M^{\ell}_{in},\mathbf{1}_{A\backslash K}\sigma \right)^2 \\
&\lesssim &
\sum_{M\in \mathcal{M}_{\nu}^{\ast }\left( K\right) } \sum_{\ell=1}^{2^n}\left( 
\frac{\mathrm{P}^{\alpha }\left( M_{{in}}^\ell,\boldsymbol{1}_{A}\sigma \right) }{\left\vert M_{{in}}^{\ell}\right\vert }\right) ^{2}\left\Vert x-m_{M_{{in}}^{\ell}}\right\Vert _{L^{2}\left( \mathbf{1}_{M_{{in}}^{\ell
}}\omega \right) }^{2}   \\ 
&\lesssim &
(\mathcal{E}_2^{\alpha }+\mathfrak{A}_2^\alpha)\left\vert I^{\prime }\right\vert _{\sigma}
\end{eqnarray*}
upon using the stopping energy condition for $I^{\prime }$ in the corona $\mathcal{C}_{A}$, i.e. the failure of (\ref{def stop 3}), in the corona $\mathcal{C}_{A}$ with the subdecomposition 
$$
I^{\prime }\supset \overset{\cdot }{\bigcup_{M\in \mathcal{M}_{\nu}^{\ast }\left( K\right) }}
\bigcup_{\ell=1}^{2^n} M_{{in}}^{\ell}
$$
Combining these four bounds together with the definition of $\Phi ^{A,B}$ in (\ref{def PHI}), after applying Cauchy-Schwarz, gives 
\begin{equation*}
\sum_{M\in \mathcal{M}_{\nu}^{\ast }\left( K\right) }\Phi ^{A,B}\left( M_{{in}}\right)  
\lesssim
\delta^{\alpha-n}\cdot\mathcal{NTV}_{\alpha }\sqrt{\left\vert
I^{\prime }\right\vert _{\sigma }\left\vert J^{\prime }\right\vert _{\omega }%
} \notag
\end{equation*}

 In particular then, if we
now sum over \emph{natural} children $I^{\prime }$ of $I$ $\in \mathcal{C}%
_{A}$ and the associated children $J^{\prime }$ of $J\in \mathcal{N}(I) $, where
\begin{equation*}
\mathcal{N}\left( I\right) \equiv \left\{
J\in \mathcal{G}:2^{-\mathbf{r}}\ell \left( I\right) <\ell \left( J\right)
\leq \ell \left( I\right) \text{ and }d\left( J,I\right) \leq 2\ell \left(
J\right) ^{\varepsilon }\ell \left( I\right) ^{1-\varepsilon }\right\}.
\end{equation*}%
we obtain the following corona estimate, using the collection of $K$ that is defined after \eqref{defIJ}, 
\begin{eqnarray}
&& \label{corona natural}\\
&&
\sum_{\substack{I\in \mathcal{C}_{A}\\ J\in \mathcal{N}(I) }}\sum_{\substack{ I^{\prime }\in \mathfrak{C
}_{{nat}}\left( I\right) \& J^{\prime }\in \mathfrak{C}
\left( J\right)  \\ K\in \mathcal{K}\left( I^{\prime },J^{\prime }\right) }}  
\!\!\!\!\!\!\left\vert E_{I^{\prime }}^{\sigma }\left( \widehat{\square }_{I}^{\sigma
,\flat ,\mathbf{b}}f\right) \right\vert  \left\vert B_{\left( I^{\prime
},J^{\prime }\right) }\left( K\right) \right\vert\!\! \ \left\vert E_{J^{\prime
}}^{\omega }\left( \widehat{\square }_{J}^{\omega ,\flat ,\mathbf{b}^{\ast
}}g\right) \right\vert  \notag \\
&\lesssim &
\notag 
\delta^{\alpha-n}\cdot B\cdot\mathcal{NTV}_{\alpha }\!\!\sum_{\substack{I\in \mathcal{C}_{A}\\ J\in \mathcal{N}(I) }}\sum_{\substack{I'\in \mathfrak{C}_{{nat}}\left(
I\right) \\ J'\in \mathfrak{C}(J) }} \!\!\!\!\!\!\!\!\left\vert
E_{I^{\prime }}^{\sigma }\left( \widehat{\square }_{I}^{\sigma ,\flat ,%
\mathbf{b}}f\right) \right\vert \ \sqrt{\left\vert I^{\prime }\right\vert
_{\sigma }\left\vert J^{\prime }\right\vert _{\omega }}\ \left\vert
E_{J^{\prime }}^{\omega }\left( \widehat{\square }_{J}^{\omega ,\flat ,%
\mathbf{b}^{\ast }}g\right) \right\vert  \notag \\
&\lesssim &\notag 
\delta^{\alpha-n}\cdot B\cdot\mathcal{NTV}_{\alpha }\left( \sum_{I\in 
\mathcal{C}_{A}}\sum_{I^{\prime }\in \mathfrak{C}_{{nat}}\left(
I\right) }\!\!\!\left\vert I^{\prime }\right\vert _{\sigma }\left\vert
E_{I^{\prime }}^{\sigma }\left( \widehat{\square }_{I}^{\sigma ,\flat ,
\mathbf{b}}f\right) \right\vert ^{2}\right) ^{\!\!\frac{1}{2}} \cdot\\ 
&&
\hspace{3cm}\cdot
\left( \sum_{I\in \mathcal{C}_{A}}
\sum_{J\in \mathcal{N}(I) }
\sum_{J'\in \mathfrak{C}\left( J\right)}\left\vert J^{\prime }\right\vert _{\omega }\ \left\vert E_{J^{\prime
}}^{\omega }\left( \widehat{\square }_{J}^{\omega ,\flat ,\mathbf{b}^{\ast
}}g\right) \right\vert ^{2}\right) ^{\frac{1}{2}}  \notag \\
&\lesssim &
\delta^{\alpha-n}\cdot          B\cdot \mathcal{NTV}_{\alpha }\left\Vert \mathsf{P}_{\mathcal{C}_{A}}^{\sigma }f\right\Vert _{L^{2}\left( \sigma \right)
}^{\bigstar }\left\Vert \mathsf{P}_{\mathcal{C}_{A}^{\mathcal{G},{%
nearby}}}^{\omega }g\right\Vert _{L^{2}\left( \sigma \right) }^{\bigstar }
\notag
\end{eqnarray}
where $\mathcal{C}_{A}^{\mathcal{G},{nearby}}=\bigcup\limits_{I\in 
\mathcal{C}_{A}}\mathcal{N}(I) $, and the final line uses (\ref{PLBP removed}) to obtain 
\begin{eqnarray*}
\sum_{I\in \mathcal{C}_{A}}\sum_{I^{\prime }\in \mathfrak{C}_{{nat}}\left( I\right) }\left\vert I^{\prime }\right\vert _{\sigma
}\left\vert E_{I^{\prime }}^{\sigma }\left( \widehat{\square }_{I}^{\sigma
,\flat ,\mathbf{b}}f\right) \right\vert ^{2}
=
\sum_{I\in \mathcal{C}_{A}}\left\Vert \widehat{\square }_{I}^{\sigma ,\flat ,\mathbf{b}%
}f\right\Vert _{L^{2}\left( \sigma \right) }^{2}
\lesssim
\sum_{I\in \mathcal{C}_{A}}\left\Vert \square _{I}^{\sigma ,\mathbf{b}}f\right\Vert
_{L^{2}\left( \sigma \right) }^{2}\leq \left\Vert \mathsf{P}_{\mathcal{C}_{A}}^{\sigma }f\right\Vert _{L^{2}\left( \sigma \right) }^{\bigstar 2}
\end{eqnarray*}
and similarly for the sum in $J$ and $J^{\prime }$, once we note that given $J\in \mathcal{C}_{A}^{\mathcal{G},{nearby}}$, there are only
boundedly many $I\in \mathcal{C}_{A}$ for which $J\in \mathcal{N}(I)$.

In order to deal with this sum in the case when the child $I^{\prime }$ is
broken, we must take the estimate one step further and sum over those broken
cubes $I^{\prime }$ whose parents belong to the corona $\mathcal{C}_{A}$, i.e. $\left\{ I^{\prime }\in \mathcal{D}:I^{\prime }\in \mathfrak{C}_{{brok}}\left( I\right) \text{ for some }I\in \mathcal{C}%
_{A}\right\} $. Of course this collection is precisely the set of $\mathcal{A}$ -children of $A$, i.e.
\begin{equation}
\left\{ I^{\prime }\in \mathcal{D}:I^{\prime }\in \mathfrak{C}_{{brok}}\left( I\right) \text{ for some }I\in \mathcal{C}_{A}\right\} =%
\mathfrak{C}_{\mathcal{A}}\left( A\right) .  \label{precisely}
\end{equation}

To obtain the same corona estimate when summing over broken $I^{\prime }$,
we will exploit the fact that the cubes $A^{\prime }\in \mathfrak{C}_{%
\mathcal{A}}\left( A\right) $ are pairwise disjoint. But first we note that
when $I^{\prime }$ is a broken child, neither weak testing nor stopping
energy is available. But if we sum over such broken $I^{\prime }$, and use (\ref{precisely}) to see that the broken children are pairwise disjoint, we
obtain the following estimate where for convenience we use the notation $\overset{\sim}{\mathcal{M}}_{\nu} \equiv \bigcup\limits_{K\in \mathcal{K}\left( I^{\prime
},J^{\prime }\right) }\mathcal{M}_{\nu}^{\ast }\left( K\right) $: 
\begin{eqnarray*}
&&
\sum_{\substack{I\in \mathcal{C}_{A}\\ J\in \mathcal{N}(I) }}\sum_{\substack{ I^{\prime }\in \mathfrak{C
}_{{brok}}\left( I\right) \& J^{\prime }\in \mathfrak{C}
\left( J\right)  \\ K\in \mathcal{K}\left( I^{\prime },J^{\prime }\right) }}  
\!\!\!\!\!\!\left\vert E_{I^{\prime }}^{\sigma }\left( \widehat{\square }_{I}^{\sigma
,\flat ,\mathbf{b}}f\right) \right\vert  \left\vert B_{\left( I^{\prime
},J^{\prime }\right) }\left( K\right) \right\vert\!\! \ \left\vert E_{J^{\prime
}}^{\omega }\left( \widehat{\square }_{J}^{\omega ,\flat ,\mathbf{b}^{\ast
}}g\right) \right\vert   \\
&\lesssim &\!\!\!\!\!
\delta^{\alpha-n}\cdot B\cdot \mathcal{NTV}_{\alpha }\!\!\!\sum_{\substack{I\in \mathcal{C}_{A}\\J\in\mathcal{N}(I)}} 
\sum_{\substack{I^{\prime }\in \mathfrak{C}_{{brok}}\left( I\right) \\ J'\in \mathfrak{C}\left( J\right) }}\!\!\!\!
\left\vert E_{I^{\prime}}^{\sigma }\left( \widehat{\square }_{I}^{\sigma ,\flat ,\mathbf{b}}f\right) \right\vert
\sqrt{\left\vert J^{\prime }\right\vert _{\omega }}\ \left\vert E_{J^{\prime
}}^{\omega }\left( \widehat{\square }_{J}^{\omega ,\flat ,\mathbf{b}^{\ast
}}g\right) \right\vert \cdot\\
&&\cdot
\left( \sum_{M\in \overset{\sim}{\mathcal{M}}_{\nu}}\!\!\!\!\!\left\Vert \mathbf{1}_{M_{{in}}}T_{\sigma
}^{\alpha }b_{A}\right\Vert _{L^{2}\left( \omega \right) }^{2}
+\!\!\!\sum_{M\in 
\overset{\sim}{\mathcal{M}}_{\nu}}\sum_{\ell=1}^{2^n} \mathrm{P}_{\delta }^{\alpha }\mathsf{Q}^{\omega }\left( M_{{in}}^\ell,\mathbf{1}_{A}\sigma \right) ^{2}
+\!\!\!
\sum_{M\in \overset{\sim}{\mathcal{M}}_{\nu}} \!\!\!\left\vert M_{{in}}\right\vert _{\sigma } \right)^{\!\!\!1/2}\\
&\lesssim &  \notag
\!\!\!\!\!B\delta^{\alpha-n}\mathcal{NTV}_{\alpha }\!\!\left(\! \sum_{\substack{ I\in \mathcal{C}_{A}  \\ I^{\prime }\in \mathfrak{C}_{{brok}}\left( I\right) }}
\!\sum _{\substack{ J\in \mathcal{N}(I) \\ J^{\prime }\in \mathfrak{C}\left( J\right) }}\!\sum_{M\in 
\overset{\sim}{\mathcal{M}}_{\nu}}  \!\!\!\!\!\left\{
\!\left\Vert \mathbf{1}_{M_{{in}}}T_{\sigma }^{\alpha
}b_{A}\right\Vert _{L^{2}\left( \omega \right) }^{2}
\!\!+\!\!\sum_{\ell=1}^{} \mathrm{P}_{\delta }^{\alpha }\mathsf{Q}^{\omega }\!\!\left( M_{{in}}^\ell,\mathbf{1}_{A}\sigma \right) ^{2}
\left\vert M_{{in}
}\right\vert _{\sigma }\right\} \right)^{\!\!\!\frac{1}{2}} \cdot
\\
&&\hspace{3.5cm}\cdot\left( \frac{1}{\left\vert A\right\vert _{\sigma }}\int_{A}\left\vert
f\right\vert d\sigma \right)\!\!\! \left(\! \sum_{\substack{ J\in \mathcal{C}_{A}^{\mathcal{G},{nearby}}  \\ J^{\prime }\in \mathfrak{C}\left( J\right) 
}}\sum_{\substack{ I\in \mathcal{C}_{A}:\ J\in \mathcal{N}(I)  \\ I^{\prime }\in \mathfrak{C}_{{
brok}}\left( I\right) }}\!\!\!\!\!\!\!\!\left\vert J^{\prime }\right\vert _{\omega
}\left\vert E_{J^{\prime }}^{\omega }\left( \widehat{\square }_{J}^{\omega
,\flat ,\mathbf{b}^{\ast }}g\right) \right\vert ^{2}\right) ^{\!\!\!\frac{1}{2}} 
\notag \\
\end{eqnarray*}
which gives that
\begin{eqnarray}
&&  \label{prelim corona broken} \\
&&\sum_{\substack{I\in \mathcal{C}_{A}\\J\in \mathcal{N}(I)}} \sum_{\substack{ I^{\prime }\in \mathfrak{C}_{{brok}}\left( I\right) \& J^{\prime }\in \mathfrak{C}\left( J\right)  \\ K\in \mathcal{K}\left( I^{\prime },J^{\prime }\right) }}
\left\vert E_{I^{\prime }}^{\sigma }\left( \widehat{\square }_{I}^{\sigma
,\flat ,\mathbf{b}}f\right) \right\vert \ \left\vert B_{\left( I^{\prime
},J^{\prime }\right) }\left( K\right) \right\vert \ \left\vert E_{J^{\prime
}}^{\omega }\left( \widehat{\square }_{J}^{\omega ,\flat ,\mathbf{b}^{\ast}}g\right) \right\vert  \notag \\
&\lesssim &
\mathcal{NTV}_{\alpha }\sqrt{\left\vert A\right\vert _{\sigma
}\left( \frac{1}{\left\vert A\right\vert _{\sigma }}\int_{A}\left\vert
f\right\vert d\sigma \right) ^{2}}\left\Vert \mathsf{P}_{\mathcal{C}_{A}^{%
\mathcal{G},{nearby}}}^{\omega }g\right\Vert _{L^{2}\left( \sigma
\right) }^{\bigstar }  \notag
\end{eqnarray}
because
\begin{equation*}
\left\vert E_{I^{\prime }}^{\sigma }\left( \widehat{\square }_{I}^{\sigma
,\flat ,\mathbf{b}}f\right) \right\vert
=
\left\vert \frac{1}{\int_{I}b_{I}d\sigma }\int_{I}fd\sigma \right\vert 
\lesssim 
\frac{1}{\left\vert I\right\vert _{\sigma }}\int_{I}\left\vert f\right\vert d\sigma
\lesssim 
\frac{1}{\left\vert A\right\vert _{\sigma }}\int_{A}\left\vert
f\right\vert d\sigma
\end{equation*}
if $I^{\prime }\in \mathfrak{C}_{{brok}}\left( I\right) $ and $%
I\in \mathcal{C}_{A}$, and because 
\begin{eqnarray}
&&  \label{last inequ} \\
&&\notag\!\!\!\!\!\!\!\!\!
\sum_{\substack{ I\in \mathcal{C}_{A}  \\ I^{\prime }\in \mathfrak{C}_{{brok}}\left( I\right) }}
\sum_{\substack{ J\in \mathcal{N}(I)  \\ J^{\prime }\in \mathfrak{C}\left( J\right) }}\sum_{M\in \overset{\sim}{\mathcal{M}}_{\nu}} 
\left\{ \left\Vert \mathbf{1}_{M_{{in}
}}T_{\sigma }^{\alpha }b_{A}\right\Vert _{L^{2}\left( \omega \right) }^{2}+
\sum_{\ell=1}^{2^n} \mathrm{P}_{\delta }^{\alpha }\mathsf{Q}^{\omega }\left( M_{{in}}^{\ell},\mathbf{1}_{A}\sigma \right) ^{2}+\left\vert M_{{in}}\right\vert
_{\sigma }\right\}  \notag \\
&\lesssim &
\left( \mathfrak{T}_{T^{\alpha }}^{\mathbf{b}}+%
\mathfrak{E}_{2}^{\alpha }+1\right) ^{2}\left\vert A\right\vert _{\sigma }
\notag
\end{eqnarray}

Indeed, in this last inequality (\ref{last inequ}), we have used first the
testing condition,
\begin{eqnarray*}
\sum_{\substack{ I\in \mathcal{C}_{A}  \\ I^{\prime }\in \mathfrak{C}_{{brok}}\left( I\right) }}
\sum_{\substack{ J\in \mathcal{N}(I)  \\ J^{\prime }\in \mathfrak{C}\left( J\right) }}\sum_{M\in \overset{\sim}{\mathcal{M}}_{\nu}}\!\!\!\left\Vert \mathbf{1}_{M_{{in}
}}T_{\sigma }^{\alpha }b_{A}\right\Vert _{L^{2}\left( \omega \right) }^{2}
&\leq& 
\mathfrak{T}_{T^\alpha}^{\mathbf{b}}\!\!\!\sum_{\substack{ I\in \mathcal{C}_{A}  \\ I^{\prime }\in \mathfrak{C}_{{brok}}\left( I\right) }}
\sum_{\substack{ J\in \mathcal{N}(I)  \\ J^{\prime }\in \mathfrak{C}\left( J\right) }}|I|_\sigma\\
&\lesssim&
\mathfrak{T}_{T^\alpha}^{\mathbf{b}} \sum_{\substack{ I\in \mathcal{C}_{A}  \\ I^{\prime }\in \mathfrak{C}_{{brok}}\left( I\right) }}|I|_\sigma
\leq 
\mathfrak{T}_{T^\alpha}^{\mathbf{b}}|A|_\sigma
\end{eqnarray*}
where in the first inequality we used the fact that the $M_{in}$ that appear are all disjoint and form a subdecomposition of $I'\subset I$ and then used testing. On the second inequality we used the bounded overlap of $J$ for any given $I$, since we are in the case of nearby cubes, and we get the last inequality because the $I\in \mathcal{C}_A$, which have a broken child $I'$, are disjoint and form a subdecomposition of $A$. The same argument can be applied for the second sum of \eqref{last inequ} upon using the energy condition for all $I\in C_A$ which have a broken child $I'$ and using the finite repetition again since we are in the nearby form.

The inequality (\ref{prelim corona broken}) is a suitable estimate since
\begin{equation*}
\sum_{A\in \mathcal{A}}\sqrt{\left\vert A\right\vert _{\sigma }\left( \frac{1%
}{\left\vert A\right\vert _{\sigma }}\int_{A}\left\vert f\right\vert d\sigma
\right) ^{2}}\left\Vert \mathsf{P}_{\mathcal{C}_{A}^{\mathcal{G},{%
nearby}}}^{\omega }g\right\Vert _{L^{2}\left( \sigma \right) }^{\bigstar
}\lesssim \left\Vert f\right\Vert _{L^{2}\left( \sigma \right) }\left\Vert
g\right\Vert _{L^{2}\left( \sigma \right) }
\end{equation*}
by quasiorthogonality and the frame inequalities (\ref{Car
embed}) and (\ref{FRAME}), together with the bounded overlap of the
`nearby' coronas $\left\{ \mathcal{C}_{A}^{\mathcal{G},{nearby}%
}\right\} _{A\in \mathcal{A}}$. We are left with estimating $\mathbf{\Delta},\mathbf{E},\mathbf{F}$ that we get after the  iteration.

Let us first deal with $\mathbf{\Delta}$. By $K^j_{i,\ell}$ we mean a grandchild of a cube $K_i^j$ and $K_i^j$ comes from $K_i$ after having iterated $j$ times, so $K^j_{i,\ell}$ is a ($2j+2$)-child of $K_i$.  We have 
\begin{eqnarray*}
&& \sum_{i=1}^B\sum_{j=1}^\nu \sum_{\ell=1}^{4^n-2^n}\mathbf{\Delta}(K^j_{i,\ell})\\
&\leq&
\mathfrak{N}_{T^\alpha}C_{\mathbf{b},\mathbf{b}^*,\nu}\sum_{i=1}^B\sum_{j=1}^\nu \sum_{\ell=1}^{4^n-2^n}
\bigg(\sum_{q=1}^{2^n}\Big|\big(K^j_{i,\ell,in}\backslash K^{j,q}_{i,\ell,in}\big)\cap(1+\delta)K^{j,q}_{i,\ell,in}\Big|_\sigma \bigg)^{\!\!\frac{1}{2}}\!\!\sqrt{|K^j_{i,\ell}|_\omega}\\
&\leq&
\mathfrak{N}_{T^\alpha}C_{\mathbf{b},\mathbf{b}^*,\nu}\bigg(\sum_{i=1}^B\sum_{j=1}^\nu\sum_{\ell=1}^{4^n-2^n}\sum_{q=1}^{2^n}\Big|\big(K^j_{i,\ell,in}\backslash K^{j,q}_{i,\ell,in}\big)\cap(1+\delta)K^{j,q}_{i,\ell,in}\Big|_\sigma\bigg)^\frac{1}{2}\sqrt{|J'|_\omega}
\end{eqnarray*}
where $K_{i,\ell,in}^{j,q}$ is one of the inner grandchildren of $K_{i,\ell,in}^{j}$.
Now fixing $q=q_0$ and taking averages over the grid $\mathcal{G}$ we get
$$
\boldsymbol{E}^\mathcal{G}_\Omega
\sum_{i=1}^B\sum_{j=1}^\nu\sum_{\ell=1}^{4^n-2^n}\Big|\big(K^j_{i,\ell,in}\backslash K^{j,q}_{i,\ell,in}\big)\cap(1+\delta)K^{j,q}_{i,\ell,in}\Big|_\sigma  \leq C_n \delta |I|_\sigma
$$
the constant depends on dimension since for the same $i,j$ we can have intersection as $\ell$ moves. Adding the different $q$ we get finally 

\begin{eqnarray}\label{Delta}
 \boldsymbol{E}^\mathcal{G}_\Omega\sum_{i=1}^B\sum_{j=1}^\nu\sum_{\ell=1}^{4^n-2^n}\mathbf{\Delta}(K^j_{i,\ell})\leq
\mathfrak{N}_{T^\alpha}C_{\mathbf{b},\mathbf{b}^*,\nu,n}\sqrt{\delta}\sqrt{|I'|_\sigma}\sqrt{|J'|_\omega}.
\end{eqnarray}
For $\mathbf{F}$ we get, 
\begin{eqnarray*}
 \sum_{i=1}^B\sum_{j=1}^\nu\sum_{\ell=1}^{4^n-2^n}\mathbf{F}(K^j_{i,\ell})
\leq
\mathfrak{N}_{T^\alpha} C_{\mathbf{b},\mathbf{b}^*}\bigg(\sum_{i=1}^B\sum_{j=1}^\nu\sum_{\ell=1}^{4^n-2^n}\Big|K^j_{i,\ell,out}\cap(1+\delta)K^j_{i,\ell,in}\Big|_\sigma\bigg)^\frac{1}{2}\sqrt{|J'|_\omega}
\end{eqnarray*}
and again averaging over grids $\mathcal{G}$, we get the bound 
\begin{equation}\label{F}
\boldsymbol{E}^\mathcal{G}_\Omega \sum_{i=1}^B\sum_{j=1}^\nu\sum_{\ell=1}^{4^n-2^n}\mathbf{F}(K^j_{i,\ell})\leq \mathfrak{N}_{T^\alpha} C_{\mathbf{b},\mathbf{b}^*}\sqrt{\delta}\sqrt{|I'|_\sigma}\sqrt{|J'|_\omega}
\end{equation}
Note here that upon choosing $\delta$ small enough there is no repetition in the different terms that arise. Finally, for $\mathbf{E}$, we have

\begin{eqnarray}\label{E}
&& 
\sum_{i=1}^B\sum_{j=1}^\nu\sum_{\ell=1}^{4^n-2^n}\mathbf{E}(K^j_{i,\ell})\\
&\leq&\notag
\mathfrak{N}_{T^\alpha}\sum_{i=1}^B\sum_{j=1}^\nu\sum_{\ell=1}^{4^n-2^n}\sum_{q=1}^{4^n-2^n}\bigg(\sum_{r>q}\Big|K^{j,q}_{i,\ell,out}\cap(1+\delta)K^{j,r}_{i,\ell,out}\Big|_\sigma\bigg)^\frac{1}{2}\sqrt{\Big|K^j_{i,\ell,out}\Big|_\omega}\\
&\leq&\notag
\mathfrak{N}_{T^\alpha}\bigg(\sum_{i=1}^B\sum_{j=1}^\nu\sum_{\ell=1}^{4^n-2^n}\sum_{q=1}^{4^n-2^n}\sum_{r>q}\Big|K^{j,q}_{i,\ell,out}\cap(1+\delta)K^{j,r}_{i,\ell,out}\Big|_\sigma\bigg)^\frac{1}{2}\cdot\\
&& \hspace{5cm}\notag\cdot \bigg(\sum_{i=1}^B\sum_{j=1}^\nu\sum_{\ell=1}^{4^n-2^n}\sum_{q=1}^{4^n-2^n}\sum_{r>q}\Big|K^j_{i,\ell,out}\Big|_\omega\bigg)^\frac{1}{2}\\
&\leq&\notag
\mathfrak{N}_{T^\alpha}\cdot C_{n,\nu}\bigg(\sum_{i=1}^B\sum_{j=1}^\nu\sum_{\ell=1}^{4^n-2^n}\sum_{q=1}^{4^n-2^n}\sum_{r>q}\Big|K^{j,q}_{i,\ell,out}\cap(1+\delta)K^{j,r}_{i,\ell,out}\Big|_\sigma\bigg)^\frac{1}{2}\sqrt{|J'|_\omega}
\end{eqnarray}
Taking averages,
$$
\boldsymbol{E}^\mathcal{G}_\Omega
\sum_{i=1}^B\sum_{j=1}^\nu\sum_{\ell=1}^{4^n-2^n}\mathbf{E}(K^j_{i,\ell})
\leq 
\mathfrak{N}_{T^\alpha}\cdot C_{n,\nu}\sqrt{\delta}\sqrt{|I'|_\sigma}\sqrt{|J'|_\omega}
$$
The constant $C_{n,\nu}$ comes from the intersection of the sets $K^j_{i,\ell,out}$.

Recall that after splitting in the cases of $\delta$-seperated and $\delta$-close cubes, we got the bound \eqref{first bound} in the separated case and after an initial application of random surgery, we reduced the proof of Proposition \ref{The prop} to establishing inequality (\ref{must show final}). Then using the bounds in \eqref{reduction to intersection}, \eqref{ML inequalities}, \eqref{LM surgery}, \eqref{LL surgery}, \eqref{KK disjoint}, \eqref{KK adjacent} we reduced $P\left( I,J\right)$ to getting a bound for $\left\{ K,K\right\} $
in the notation used in (\ref{def E,F}). Then using the estimates in (\ref{C est}), (\ref{A est}), (\ref{corona natural}) and (\ref{prelim corona broken}) together with (\ref{K iterated}), (\ref{Delta}), (\ref{F}) and (\ref{E}) establishes
probabilistic control of the sum of all the inner products $\left\{
K,K\right\} $ taken over appropriate cubes $K$, yielding (\ref{must show final}) as required if we choose $\varepsilon, \lambda $, $\eta_{0}$ and $\delta$ sufficiently
small. And combining all the above bounds we proved proposition \ref{The prop}, namely we got the bound
\begin{eqnarray*}
\ \boldsymbol{E}_{\Omega }^{\mathcal{D}}\boldsymbol{E}_{\Omega }^{\mathcal{G}%
}\sum_{I\in \mathcal{D}}\sum_{\substack{ J\in \mathcal{G}:\ 2^{-\mathbf{r} n}|I| <|J| \leq |I|  \\ %
d\left( J,I\right) \leq 2\ell(J) ^{\varepsilon }\ell(I) ^{1-\varepsilon }}}\!\!\!\left\vert \left\langle T_{\sigma }^{\alpha
}\left( \square _{I}^{\sigma ,\mathbf{b}}f\right) ,\square _{J}^{\omega ,%
\mathbf{b}^{\ast }}g\right\rangle _{\omega }\right\vert  
\lesssim 
\left( C_{\theta }\mathcal{NTV}_{\alpha }+\sqrt{\theta }\mathfrak{%
N}_{T^{\alpha }}\right) \left\Vert f\right\Vert _{L^{2}\left( \sigma \right)
}\left\Vert g\right\Vert _{L^{2}\left( \omega \right) }
\end{eqnarray*}%

\section{Main below form}\label{Sec Main below}

Now we turn to controlling the main below form (\ref{def Theta 2 good}),%
\begin{equation*}
\Theta _{2}^{{good}}\left( f,g\right) =\sum_{I\in \mathcal{D}%
}\sum_{J^{\maltese }\subsetneqq I:\ \ell \left( J\right) \leq 2^{-\mathbf{\rho }%
}\ell \left( I\right) }\int \left( T_{\sigma }\square _{I}^{\sigma ,\mathbf{b%
}}f\right) \square _{J}^{\omega ,\mathbf{b}^{\ast }}gd\omega .
\end{equation*}%

To control $\Theta _{2}^{{good}}\left( f,g\right) \equiv \mathsf{B}%
_{\Subset _{\mathbf{\rho }}}\left( f,g\right) $ we first perform the \emph{%
canonical corona splitting} of $\mathsf{B}_{\Subset _{\mathbf{\rho }}}\left(
f,g\right) $ into a diagonal form and a far below form, namely $\mathsf{T}_{{diagonal}}\left(
f,g\right) $ and $\mathsf{T}_{{far}{below}%
}\left( f,g\right) $ as in \cite{SaShUr6}. This \emph{canonical splitting}
of the form $\mathsf{B}_{\Subset _{\mathbf{\rho }}}\left( f,g\right) $
involves the corona pseudoprojections $\mathsf{P}_{\mathcal{C}_{A}^{\mathcal{%
D}}}^{\sigma ,\mathbf{b}}$ acting on $f$ and the \emph{shifted} corona
pseudoprojections $\mathsf{P}_{\mathcal{C}_{B}^{\mathcal{G},{shift}%
}}^{\omega ,\mathbf{b}^{\ast }}$ acting on $g$, where $B$ is a stopping cube
in $\mathcal{A}$. The stopping cubes $\mathcal{B}$ constructed relative to $%
g\in L^{2}\left( \omega \right) $ play no role in the analysis here, except
to guarantee that the frame and weak Riesz inequalities hold for $g$ and $%
\left\{ \square _{J}^{\omega ,\mathbf{b}^{\ast }}g\right\} _{J\in \mathcal{G}%
}$. Here the shifted corona $\mathcal{C}_{B}^{\mathcal{G},{shift}}$
is defined to include those cubes $J\in \mathcal{G}$ such $J^{\maltese }\in 
\mathcal{C}_{B}^{\mathcal{D}}$. Recall that the parameters $\mathbf{\tau }$
and $\mathbf{\rho }$ are fixed to satisfy 
\begin{equation*}
\mathbf{\tau }>\mathbf{r}\text{ and }\mathbf{\rho }>\mathbf{r}+\mathbf{\tau }%
,
\end{equation*}%
where $\mathbf{r}$ is the goodness parameter already fixed in (\ref{choice
of r}).

\begin{dfn}
\label{shifted corona}For $B\in \mathcal{A}$ we define the shifted $\mathcal{%
G}$-corona by%
\begin{equation*}
\mathcal{C}_{B}^{\mathcal{G},{shift}}=\left\{ J\in \mathcal{G}%
:J^{\maltese }\in \mathcal{C}_{B}^{\mathcal{D}}\right\} .
\end{equation*}
\end{dfn}

We will use repeatedly the fact that the shifted coronas $\mathcal{C}_{B}^{%
\mathcal{G},{shift}}$ are pairwise disjoint in $B$:%
\begin{equation}
\sum_{B\in \mathcal{A}}\mathbf{1}_{\mathcal{C}_{B}^{\mathcal{G},{%
shift}}}\left( J\right) \leq \mathbf{1},\ \ \ \ \ J\in \mathcal{D}.
\label{tau overlap}
\end{equation}%
The forms $\mathsf{B}_{\Subset _{\mathbf{\rho },\varepsilon }}\left(
f,g\right) $ are no longer linear in $f$ and $g$ as the `cut' is determined
by the coronas $\mathcal{C}_{A}^{\mathcal{D}}$ and $\mathcal{C}_{B}^{%
\mathcal{G},{shift}}$, which depend on $f$ as well as the measures $%
\sigma $ and $\omega $. However, if the coronas are held fixed, then the
forms can be considered bilinear in $f$ and $g$. It is convenient at this
point to introduce the following shorthand notation:%
\begin{equation}\label{def shorthand}
\left\langle T_{\sigma }^{\alpha }\left( \mathsf{P}_{\mathcal{C}_{A}^{%
\mathcal{D}}}^{\sigma ,\mathbf{b}}f\right) ,\mathsf{P}_{\mathcal{C}_{B}^{%
\mathcal{G},{shift}}}^{\omega ,\mathbf{b}^{\ast }}g\right\rangle
_{\omega }^{\Subset _{\mathbf{\rho },\varepsilon }}
\equiv\!\!\!\!
\sum_{\substack{ %
I\in \mathcal{C}_{A}^{\mathcal{D}}\text{ and }J\in \mathcal{C}_{B}^{\mathcal{%
G},{shift}}:\ J^{\maltese }\subsetneqq I  \\ \ell \left( J\right) \leq
2^{-\mathbf{\rho }}\ell \left( I\right) }} 
\!\!\!\!\!\left\langle T_{\sigma }^{\alpha
}\left( \square _{I}^{\sigma ,\mathbf{b}}f\right) ,\square _{J}^{\omega ,%
\mathbf{b}^{\ast }}g\right\rangle _{\omega }\ .
\end{equation}

\begin{description}
\item[Caution] One must not assume, from the notation on the left hand side
above, that the function $T_{\sigma }^{\alpha }\left( \mathsf{P}_{\mathcal{C}%
_{A}}^{\sigma }f\right) $ is simply integrated against the function $\mathsf{%
P}_{\mathcal{C}_{B}^{\mathcal{G},{shift}}}^{\omega }g$. Indeed, the
sum on the right hand side is taken over pairs $\left( I,J\right) $ such
that $J^{\maltese }\in \mathcal{C}_{B}^{\mathcal{D}}$ and $J^{\maltese
}\subsetneqq I$ and $\ell \left( J\right) \leq 2^{-\mathbf{\rho }}\ell \left(
I\right) $.
\end{description}

\subsection{The canonical splitting and local below forms}

We then have the canonical splitting determined by the coronas $\mathcal{C}%
_{A}^{\mathcal{D}}$ for $A\in \mathcal{A}$ (the stopping times $\mathcal{B}$ play no explicit role in the canonical splitting of the below
form, other than to guarantee the weak Riesz inequalities for the dual
martingale pseudoprojections $\square _{J}^{\omega ,\mathbf{b}^{\ast }}$)

\begin{eqnarray}
&&\mathsf{B}_{\Subset _{\mathbf{\rho },\varepsilon }}\left( f,g\right)
\label{parallel corona decomp'} \\
&=&\sum_{A,B\in \mathcal{A}}\left\langle T_{\sigma }^{\alpha }\left( \mathsf{%
P}_{\mathcal{C}_{A}}^{\sigma ,\mathbf{b}}f\right) ,\mathsf{P}_{\mathcal{C}%
_{B}^{\mathcal{G},{shift}}}^{\omega ,\mathbf{b}^{\ast
}}g\right\rangle _{\omega }^{\Subset _{\mathbf{\rho },\varepsilon }}  \notag
\\
&=&
\sum_{A\in \mathcal{A}}\left\langle T_{\sigma }^{\alpha }\left( \mathsf{P}%
_{\mathcal{C}_{A}}^{\sigma ,\mathbf{b}}f\right) ,\mathsf{P}_{\mathcal{C}%
_{A}^{\mathcal{G},{shift}}}^{\omega ,\mathbf{b}^{\ast
}}g\right\rangle _{\omega }^{\Subset _{\mathbf{\rho },\varepsilon }}+\sum 
_{\substack{ A,B\in \mathcal{A}  \\ B\subsetneqq A}}\left\langle T_{\sigma
}^{\alpha }\left( \mathsf{P}_{\mathcal{C}_{A}}^{\sigma ,\mathbf{b}}f\right) ,%
\mathsf{P}_{\mathcal{C}_{B}^{\mathcal{G},{shift}}}^{\omega ,\mathbf{b%
}^{\ast }}g\right\rangle _{\omega }^{\Subset _{\mathbf{\rho },\varepsilon }}
\notag \\
&&
+\sum_{\substack{ A,B\in \mathcal{A}  \\ B\supsetneqq A}}\left\langle
T_{\sigma }^{\alpha }\left( \mathsf{P}_{\mathcal{C}_{A}}^{\sigma ,\mathbf{b}%
}f\right) ,\mathsf{P}_{\mathcal{C}_{B}^{\mathcal{G},{shift}%
}}^{\omega ,\mathbf{b}^{\ast }}g\right\rangle _{\omega }^{\Subset _{\mathbf{%
\rho },\varepsilon }}+\!\!\!\sum_{\substack{ A,B\in \mathcal{A}  \\ A\cap
B=\emptyset }}\left\langle T_{\sigma }^{\alpha }\left( \mathsf{P}_{\mathcal{C%
}_{A}}^{\sigma ,\mathbf{b}}f\right) ,\mathsf{P}_{\mathcal{C}_{B}^{\mathcal{G}%
,{shift}}}^{\omega ,\mathbf{b}^{\ast }}g\right\rangle _{\omega
}^{\Subset _{\mathbf{\rho },\varepsilon }}  \notag \\
&\equiv &
\mathsf{T}_{{diagonal}}\left( f,g\right) +\mathsf{T}_{%
{far}{below}}\left( f,g\right) +\mathsf{T}_{{far}%
{above}}\left( f,g\right) +\mathsf{T}_{{disjoint}}\left(
f,g\right) .  \notag
\end{eqnarray}%
Now the final two terms $\mathsf{T}_{{far}{above}}\left(
f,g\right) $ and $\mathsf{T}_{{disjoint}}\left( f,g\right) $ each
vanish since there are no pairs $\left( I,J\right) \in \mathcal{C}_{A}^{%
\mathcal{D}}\times \mathcal{C}_{B}^{\mathcal{G},{shift}}$ with both (%
\textbf{i}) $J^\maltese \subsetneqq I$ and (\textbf{ii}) either $B\subsetneqq A$ or $B\cap A=\emptyset $%
. The far below form $\mathsf{T}_{{far}{below}}\left(
f,g\right) $ requires functional energy, which we discuss in a moment.

Next we follow this splitting by a further decomposition of the diagonal
form into local below forms $\mathsf{B}_{\Subset _{\mathbf{\rho }%
}}^{A}\left( f,g\right) $ given by the individual corona pieces 
\begin{equation}\label{def local}
\mathsf{B}_{\Subset _{\mathbf{\rho },\varepsilon }}^{A}\left( f,g\right)
\equiv \left\langle T_{\sigma }^{\alpha }\left( \mathsf{P}_{\mathcal{C}%
_{A}}^{\sigma ,\mathbf{b}}f\right) ,\mathsf{P}_{\mathcal{C}_{A}^{\mathcal{G},%
{shift}}}^{\omega ,\mathbf{b}^{\ast }}g\right\rangle _{\omega
}^{\Subset _{\mathbf{\rho },\varepsilon }}
\end{equation}%
and prove the following estimate:
\begin{equation*}
\left\vert \mathsf{B}_{\Subset _{\mathbf{\rho },\varepsilon }}^{A}\left(
f,g\right) \right\vert 
\lesssim 
\mathcal{NTV}_{\alpha } \ \left( \alpha _{\mathcal{A}}\left( A\right) \sqrt{\left\vert A\right\vert _{\sigma }}+\left\Vert \mathsf{P}_{\mathcal{C}%
_{A}}^{\sigma ,\mathbf{b}}f\right\Vert _{L^{2}\left( \sigma \right)
}^{\bigstar }\right) \ \left\Vert \mathsf{P}_{\mathcal{C}_{A}^{\mathcal{G},%
{shift}}}^{\omega ,\mathbf{b}^{\ast }}g\right\Vert _{L^{2}\left(
\omega \right) }^{\bigstar }\ .
\end{equation*}%
This reduces matters to the local forms since we then have from
Cauchy-Schwarz that%
\begin{eqnarray*}
\sum_{A\in \mathcal{A}}\left\vert \mathsf{B}_{\Subset _{\mathbf{\rho }%
,\varepsilon }}^{A}\left( f,g\right) \right\vert
&\lesssim&
\mathcal{NTV}%
_{\alpha }\ \left( \sum_{A\in \mathcal{A}}\alpha _{\mathcal{A}}\left(
A\right) ^{2}\left\vert A\right\vert _{\sigma }+\left\Vert \mathsf{P}_{%
\mathcal{C}_{A}^{\mathcal{D}}}^{\sigma ,\mathbf{b}}f\right\Vert
_{L^{2}\left( \sigma \right) }^{\bigstar 2}\right) ^{\frac{1}{2}}\!\! \left(\sum_{A\in \mathcal{A}}\left\Vert \mathsf{P}_{\mathcal{C}_{A}^{\mathcal{G},%
{shift}}}^{\omega ,\mathbf{b}^{\ast }}g\right\Vert _{L^{2}\left(
\omega \right) }^{\bigstar 2}\right) ^{\frac{1}{2}} \\
&\lesssim &
\mathcal{NTV}_{\alpha }\ \left\Vert f\right\Vert _{L^{2}\left( \sigma
\right) }\left\Vert g\right\Vert _{L^{2}\left( \omega \right) }\ .
\end{eqnarray*}
by the lower frame inequalities
$$
\sum\limits_{A\in \mathcal{A}}\left\Vert \mathsf{P}_{\mathcal{C}_{A}}^{\sigma ,\mathbf{b}}f\right\Vert _{L^{2}\left( \sigma
\right) }^{\bigstar 2}\lesssim \left\Vert f\right\Vert _{L^{2}\left( \sigma\right) }^{2}
\text{ and } 
\sum\limits_{A\in \mathcal{A}}\left\Vert \mathsf{P}_{\mathcal{C}%
_{A}^{\mathcal{G},{shift}}}^{\omega ,\mathbf{b}^{\ast }}g\right\Vert
_{L^{2}\left( \omega \right) }^{\bigstar 2}\lesssim \left\Vert g\right\Vert
_{L^{2}\left( \omega \right) }^{2}
$$
using also quasi-orthogonality $\sum\limits_{A\in 
\mathcal{A}}\alpha _{\mathcal{A}}\left( f\right) ^{2}\left\vert A\right\vert
_{\sigma }\lesssim \left\Vert f\right\Vert _{L^{2}\left( \sigma \right)
}^{2} $ in the stopping cubes $\mathcal{A}$, and the pairwise
disjointedness of the shifted coronas $\mathcal{C}_{A}^{\mathcal{G},{%
shift}}$: 
\begin{equation*}
\sum_{A\in \mathcal{A}}\mathbf{1}_{\mathcal{C}_{A}^{\mathcal{G},{%
shift}}}\leq \mathbf{1}_{\mathcal{D}}.
\end{equation*}%
From now on we will often write $\mathcal{C}_{A}$ in place of $\mathcal{C}_{A}^{\mathcal{D}}$ when no confusion is possible.

Finally, the local forms $\mathsf{B}_{\Subset _{\mathbf{\rho },\varepsilon
}}^{A}\left( f,g\right) $ are decomposed into stopping $\mathsf{B}_{{%
stop}}^{A}\left( f,g\right) $, paraproduct $\mathsf{B}_{{paraproduct}%
}^{A}\left( f,g\right) $ and neighbour $\mathsf{B}_{{neighbour}%
}^{A}\left( f,g\right) $ forms. The paraproduct and neighbour terms are
handled as in \cite{SaShUr6}, which in turn follows the treatment
originating in \cite{NTV3}, and this leaves only the stopping form $\mathsf{B%
}_{{stop}}^{A}\left( f,g\right) $ to be bounded, which we treat last by adapting the bottom/up stopping time and recursion of M.
Lacey in \cite{Lac}.

However, in order to obtain the required bounds of the above forms into
which the below form $\mathsf{B}_{\Subset _{\mathbf{\rho }}}\left(
f,g\right) $ was decomposed, we need functional energy. Recall that the
vector-valued function $\mathbf{b}$ in the accretive coronas `breaks' only
at a collection of cubes satisfying a Carleson condition. We define $%
\mathcal{M}_{\left( \mathbf{r},\varepsilon \right) -{deep}}\left(
F\right) $ to consist of the \emph{maximal} $\mathbf{r}$-deeply embedded
dyadic $\mathcal{G}$-subcubes of a $\mathcal{D}$-cube $F$ - see (\ref{def
M_r-deep}) in Appendix  for more detail.

\begin{dfn}
\label{functional energy n}Let $\mathfrak{F}_{\alpha }=%
\mathfrak{F}_{\alpha }\left( \mathcal{D},\mathcal{G}\right) $ be the
smallest constant in the `functional energy' inequality below,
holding for all $h\in L^{2}\left( \sigma \right) $ and all $\sigma $%
-Carleson collections $\mathcal{F}\subset \mathcal{D}$ with Carleson norm $%
C_{\mathcal{F}}$ bounded by a fixed constant $C$: 
\begin{equation}
\sum_{F\in \mathcal{F}}\sum_{M\in \mathcal{M}_{\left( \mathbf{r},1\right) -%
{deep},\mathcal{D}}\left( F\right) }\left( \frac{\mathrm{P}^{\alpha
}\left( M,h\sigma \right) }{\left\vert M\right\vert ^{\frac{1}{n}}}\right)
^{2}\left\Vert \mathsf{Q}_{\mathcal{C}_{F}^{\mathcal{G},{shift}%
}; M}^{\omega ,\mathbf{b}}x\right\Vert _{L^{2}\left( \omega \right)
}^{\spadesuit 2}\leq \mathfrak{F}_{\alpha }\lVert h\rVert _{L^{2}\left( \sigma
\right) }\,,  \label{e.funcEnergy n}
\end{equation}
\end{dfn}

The main ingredient used in reducing control of the below form $\mathsf{B}%
_{\Subset _{\mathbf{\rho }}}\left( f,g\right) $ to control of the functional
energy $\mathfrak{F}_{\alpha }$ constant and the stopping form $\mathsf{B}_{%
{stop}}^{A}\left( f,g\right) $, is the Intertwining Proposition from 
\cite{SaShUr6}. The control of the functional energy condition by the energy
and Muckenhoupt conditions must also be adapted in light of the $p$-weakly
accretive function $\mathbf{b}$ that only `breaks' at a collection of cubes
satisfying a Carleson condition, but this poses no real difficulties. The
fact that the usual Haar bases are orthonormal is here replaced by the
weaker condition that the corresponding broken Haar `bases' are merely
frames satisfying certain lower and weak upper Riesz inequalities, but again
this poses no real difference in the arguments. Finally, the fact that
goodness for $J$ has been replaced with weak goodness, namely $J^{\maltese
}\subsetneqq I$, again forces no real change in the arguments.

We then use the paraproduct / neighbour / stopping splitting mentioned above
to reduce boundedness of $\mathsf{B}_{\Subset _{\mathbf{\rho },\varepsilon
}}^{A}\left( f,g\right) $ to boundedness of the associated stopping form 
\begin{equation}
\mathsf{B}_{{stop}}^{A}\left( f,g\right) \equiv \sum_{I\in \mathcal{C%
}_{A}}\sum_{\substack{ J\in \mathcal{C}_{A}^{\mathcal{G},{shift}}:\
J^{\maltese }\subsetneqq I  \\ \ell \left( J\right) \leq 2^{-\mathbf{\rho }}\ell
\left( I\right) }}\left( E_{I_{J}}^{\sigma }\square _{I}^{\sigma ,%
\mathbf{b}}f\right) \ \left\langle T_{\sigma }^{\alpha }\mathbf{1}%
_{A\backslash I_{J}}b_{A},\square _{J}^{\omega ,\mathbf{b}^{\ast
}}g\right\rangle _{\omega }  \label{bounded stopping form}
\end{equation}%
where $f$ is supported in the cube $A$ and its expectations $E_{I}^{\sigma }\left\vert f\right\vert $ are bounded by $\alpha _{\mathcal{A}%
}\left( A\right) $ for $I\in \mathcal{C}_{A}^{\sigma }$, the dual martingale
support of $f$ is contained in the corona $\mathcal{C}_{A}^{\sigma }$, and
the dual martingale support of $g$\ is contained in $\mathcal{C}_{A}^{%
\mathcal{G},{shift}}$, and where $I_{J}$ is the $\mathcal{D}$-child
of $I$ that contains $J$.

\subsection{Diagonal and far below forms}

Now we turn to the \emph{diagonal} and the \emph{far below} terms  $\mathsf{T}_{{diagonal}}\left( f,g\right)$ and $\mathsf{T}_{{far}{below}}\left( f,g\right)$, where in \cite{SaShUr6} the far below
terms were bounded using the Intertwining Proposition and the control of
functional energy condition by the energy conditions, but of course under
the restriction there that the cubes $J$ were good. Here we write
\begin{eqnarray}
&&\label{Tfarbelow}\\ 
\notag \mathsf{T}_{{far}{below}}\left( f,g\right) &=&
\sum_{\substack{ A,B\in \mathcal{A}  \\ B\subsetneqq A}}\sum_{\substack{ I\in 
\mathcal{C}_{A}\text{ and }J\in \mathcal{C}_{B}^{\mathcal{G},{shift}}  \\ J^{\maltese }\subsetneqq I\text{ and }\ell \left( J\right) \leq 2^{-\mathbf{\mathbf{r}}}\ell \left( I\right) }}\left\langle T_{\sigma }^{\alpha
}\left( \square _{I}^{\sigma ,\mathbf{b}}f\right) ,\left( \square
_{J}^{\omega ,\mathbf{b}^{\ast }}g\right) \right\rangle _{\omega } \\
&=&\notag
\sum_{B\in \mathcal{A}}\sum_{I\in \mathcal{D}:\ B\subsetneqq
I}\left\langle T_{\sigma }^{\alpha }\left( \square _{I}^{\sigma ,\mathbf{b}}f\right) ,\sum_{J\in \mathcal{C}_{B}^{\mathcal{G},{shift}}}\square
_{J}^{\omega ,\mathbf{b}^{\ast }}g\right\rangle _{\omega } \notag \\ 
&&\notag
- \sum_{B\in 
\mathcal{A}}\sum_{I\in \mathcal{D}:\ B\subsetneqq I}\left\langle T_{\sigma
}^{\alpha }\left( \square _{I}^{\sigma ,\mathbf{b}}f\right) ,\sum_{\substack{
J\in \mathcal{C}_{B}^{\mathcal{G},{shift}}  \\ \ell \left( J\right)
>2^{-\mathbf{\mathbf{r}}}\ell \left( I\right) }}\square _{J}^{\omega ,\mathbf{b}%
^{\ast }}g\right\rangle_{\omega } \\
&=&\notag
\mathsf{T}_{{far}{below}}^{1}\left( f,g\right) -\mathsf{T}%
_{{far}{below}}^{2}\left( f,g\right) .
\end{eqnarray}
since if $I\in \mathcal{C}_{A}$ and $J\in \mathcal{C}_{B}^{\mathcal{G},%
{shift}}$, with $J^{\maltese }\subsetneqq I$ and $B\subsetneqq A$,
then we must have $B\subsetneqq I$. First, we note that expectation of the
second sum $\mathsf{T}_{{far}{below}}^{2}\left( f,g\right)$ is controlled by (\ref{delta near}) in Proposition \ref{The prop} , i.e.
\begin{eqnarray*}
&&\boldsymbol{E}_{\Omega }^{\mathcal{D}}\boldsymbol{E}_{\Omega }^{\mathcal{G}%
}\left\vert \sum_{B\in \mathcal{A}}\sum_{I\in \mathcal{D}:\ B\subsetneqq
I}\left\langle T_{\sigma }^{\alpha }\left( \square _{I}^{\sigma ,\mathbf{b}%
}f\right) ,\sum_{\substack{ J\in \mathcal{C}_{B}^{\mathcal{G},{shift}%
}  \\ \ell \left( J\right) >2^{-\mathbf{r}}\ell \left( I\right) }}\square
_{J}^{\omega ,\mathbf{b}^{\ast }}g\right\rangle _{\omega }\right\vert \\
&\lesssim &\boldsymbol{E}_{\Omega }^{\mathcal{D}}\boldsymbol{E}_{\Omega }^{%
\mathcal{G}}\sum_{I\in \mathcal{D}}\sum_{\substack{ J\in \mathcal{G}:\ 2^{-%
\mathbf{r}}\ell \left( I\right) <\ell \left( J\right) \leq \ell \left(
I\right)  \\ d\left( J,I\right) \leq 2\ell \left( J\right) ^{\varepsilon
}\ell \left( I\right) ^{1-\varepsilon }}}\left\vert \left\langle T_{\sigma
}^{\alpha }\left( \square _{I}^{\sigma ,\mathbf{b}}f\right) ,\square
_{J}^{\omega ,\mathbf{b}^{\ast }}g\right\rangle _{\omega }\right\vert \\
&\lesssim &\left( C_{\theta }\mathcal{NTV}_{\alpha }+\sqrt{\theta }\mathfrak{%
N}_{T^{\alpha }}\right) \left\Vert f\right\Vert _{L^{2}\left( \sigma \right)
}\left\Vert g\right\Vert _{L^{2}\left( \omega \right) }\ .
\end{eqnarray*}%

The form $\mathsf{T}_{{far}{below}}^{1}\left( f,g\right) $
can be written as%
\begin{eqnarray*}
\mathsf{T}_{{far}{below}}^{1}\left( f,g\right) &=&\sum_{B\in 
\mathcal{A}}\sum_{I\in \mathcal{D}:\ B\subsetneqq I}\left\langle T_{\sigma
}^{\alpha }\left( \square _{I}^{\sigma,\mathbf{b}}f\right),g_{B}\right\rangle
_{\omega }; \\
\text{where }g_{B} &\equiv &\sum_{J\in \mathcal{C}_{B}^{\mathcal{G},{%
shift}}}\square _{J}^{\omega,\mathbf{b}* }g=\mathsf{P}_{\mathcal{C}_{F}^{%
\mathcal{G},{shift}}}^{\omega ,\mathbf{b}^{\ast }}g
\end{eqnarray*}%
and the Intertwining Proposition \ref{strongly adapted} can now be applied
to this latter form to show that it is bounded by $\mathcal{NTV}_{\alpha }+%
\mathfrak{F}_{\alpha }$. Then Proposition \ref{func ener control}\
can be applied to show that $\mathfrak{F}_{\alpha }\lesssim 
\mathfrak{A}_{2}^{\alpha }+\mathcal{E}_{2}^{\alpha }$, which completes the
proof that%
\begin{equation}
\left\vert \mathsf{T}_{{far}{below}}\left( f,g\right)
\right\vert \lesssim \mathcal{NTV}_{\alpha }\ \left\Vert f\right\Vert
_{L^{2}\left( \sigma \right) }\left\Vert g\right\Vert _{L^{2}\left( \omega
\right) }\ .  \label{far below bound}
\end{equation}

\subsection{Intertwining Proposition}

First we adapt the relevant definitions and theorems from \cite{SaShUr6}.

\begin{dfn}
\label{sigma carleson n}A collection $\mathcal{F}$ of dyadic cubes is $%
\sigma $\emph{-Carleson} if
\begin{equation*}
\sum_{F\in \mathcal{F}:\ F\subset S}\left\vert F\right\vert _{\sigma }\leq
C_{\mathcal{F}}\left\vert S\right\vert _{\sigma },\ \ \ \ \ S\in \mathcal{F}.
\end{equation*}%
The constant $C_{\mathcal{F}}$ is referred to as the Carleson norm of $%
\mathcal{F}$.
\end{dfn}

\begin{dfn}\label{def shift}
Let $\mathcal{F}$ be a collection of dyadic cubes in a grid $\mathcal{D}$.
Then for $F\in \mathcal{F}$, we define the shifted corona $\mathcal{C}_{F}^{%
\mathcal{G},{shift}}$ in analogy with Definition \ref{shifted corona}
by%
\begin{equation*}
\mathcal{C}_{F}^{\mathcal{G},{shift}}=\left\{ J\in \mathcal{G}%
:J^{\maltese }\in \mathcal{C}_{F}\right\} .
\end{equation*}
\end{dfn}

Note that the collections $\mathcal{C}_{F}^{\mathcal{G},{shift}}$
are pairwise disjoint in $F$. Let $\mathfrak{C}_{\mathcal{F}}\left( F\right) 
$ denote the set of $\mathcal{F}$-children of $F$. Given any collection $%
\mathcal{H}\subset \mathcal{G}$ of cubes, a family $\mathbf{b}^{\ast }$
of dual testing functions, and an arbitrary cube $K\in \mathcal{P}$, we
define the corresponding dual pseudoprojection $\mathsf{P}_{\mathcal{H}%
}^{\omega ,\mathbf{b}^{\ast }}$ and its localization $\mathsf{P}_{\mathcal{H}%
;K}^{\omega ,\mathbf{b}^{\ast }}$ to $K$ by%
\begin{equation}
\mathsf{Q}_{\mathcal{H}}^{\omega ,\mathbf{b}^{\ast }}=\sum_{H\in \mathcal{H}%
}\bigtriangleup _{H}^{\omega ,\mathbf{b}^{\ast }}\text{ and }\mathsf{Q}_{%
\mathcal{H};K}^{\omega ,\mathbf{b}^{\ast }}=\sum_{H\in \mathcal{H}:\
H\subset K}\bigtriangleup _{H}^{\omega ,\mathbf{b}^{\ast }}\ .
\label{def localization}
\end{equation}%
Recall from Definition \ref{functional energy n} that $\mathfrak{F}_{\alpha
}=\mathfrak{F}_{\alpha }\left( \mathcal{D},\mathcal{G}\right) =\mathfrak{F}%
_{\alpha }^{\mathbf{b}^{\ast }}\left( \mathcal{D},\mathcal{G}\right) $ is
the best constant in (\ref{e.funcEnergy n}), i.e. 
\begin{equation*}
\sum_{F\in \mathcal{F}}\sum_{M\in \mathcal{M}_{\left( \mathbf{r},1\right) -%
{deep},\mathcal{D}}\left( F\right) }\left( \frac{\mathrm{P}^{\alpha
}\left( M,h\sigma \right) }{\left\vert M\right\vert ^{\frac{1}{n}}}\right)
^{2}\left\Vert \mathsf{Q}_{\mathcal{C}_{F}^{\mathcal{G},{shift}%
}; M}^{\omega ,\mathbf{b}}x\right\Vert _{L^{2}\left( \omega \right)
}^{\spadesuit 2}\leq \mathfrak{F}_{\alpha }\lVert h\rVert _{L^{2}\left( \sigma
\right) }\,. 
\end{equation*}

\begin{rem}
\label{explaining funct ener}If in (\ref{e.funcEnergy n}), we take $%
h=\mathbf{1}_{I}$ and $\mathcal{F}$ to be the trivial Carleson collection $%
\left\{ I_{r}\right\} _{r=1}^{\infty }$ where the cubes $I_{r}$ are pairwise
disjoint in $I$, then we obtain the deep energy condition in Definition \ref%
{energy condition}, but with $\mathsf{P}_{\mathcal{C}_{F}^{\mathcal{G},%
{shift}}; M}^{\omega ,\mathbf{b}^{\ast }}$ in place of $\mathsf{P%
}_{J}^{{weak}{good},\omega }$. However, the pseudoprojection 
$\mathsf{P}_{J}^{{weak}{good},\omega }$ is larger than $%
\mathsf{P}_{\mathcal{C}_{F}^{\mathcal{G},{shift}};J}^{\omega ,%
\mathbf{b}^{\ast }}$, and so we just miss obtaining the deep energy
condition as a consequence of the functional energy condition. Nevertheless,
this near miss with $h=\mathbf{1}_{I}$ explains the terminology `functional'
energy.
\end{rem}

We will need the following `indicator' version of the estimates proved above
for the disjoint form.

\begin{lem}
\label{standard indicator}Suppose $T^{\alpha }$ is a standard fractional
singular integral with $0\leq \alpha <1$, that $\mathbf{\rho }>\mathbf{r}$,
that $f\in L^{2}\left( \sigma \right) $ and $g\in L^{2}\left( \omega \right) 
$, that $\mathcal{F}\subset \mathcal{D}^{\sigma }$ and $\mathcal{G}\subset 
\mathcal{D}^{\omega }$ are $\sigma $-Carleson and $\omega $-Carleson
collections respectively, i.e.,%
\begin{equation*}
\sum_{F^{\prime }\in \mathcal{F}:\ F^{\prime }\subset F}\left\vert F^{\prime
}\right\vert _{\sigma }\lesssim \left\vert F\right\vert _{\sigma },\ \ \ \ \
F\in \mathcal{F},\text{ and }\sum_{G^{\prime }\in \mathcal{G}:\ G^{\prime
}\subset G}\left\vert G^{\prime }\right\vert _{\omega }\lesssim \left\vert
G\right\vert _{\omega },\ \ \ \ \ G\in \mathcal{G},
\end{equation*}%
that there are numerical sequences $\left\{ \alpha _{\mathcal{F}}\left(
F\right) \right\} _{F\in \mathcal{F}}$ and $\left\{ \beta _{\mathcal{G}%
}\left( G\right) \right\} _{G\in \mathcal{G}}$ such that%
\begin{equation}
\sum_{F\in \mathcal{F}}\alpha _{\mathcal{F}}\left( F\right) ^{2}\left\vert
F\right\vert _{\sigma }\leq \left\Vert f\right\Vert _{L^{2}\left( \sigma
\right) }^{2}\text{ and }\sum_{G\in \mathcal{G}}\beta _{\mathcal{G}}\left(
G\right) ^{2}\left\vert G\right\vert _{\sigma }\leq \left\Vert g\right\Vert
_{L^{2}\left( \sigma \right) }^{2}\ ,  \label{qo}
\end{equation}%
Then
\begin{eqnarray}
&&
\sum_{F\in\mathcal{F}} \sum_{\substack{J \in \mathcal{G}:\ \ell(J)\leq \ell(F)  \\ d(J,F)>2\ell(J)^\varepsilon\ell(F)^{1-\varepsilon} }}\left\vert
\left\langle T_{\sigma }^{\alpha }\left(\mathbf{1}_{F}\alpha _{%
\mathcal{F}}\left( F\right) \right) ,\square _{J}^{\omega ,\mathbf{b}^{\ast
}}g\right\rangle _{\omega }\right\vert  \label{indicator far} \\
&&
+
\sum_{G\in\mathcal{G}}\sum_{\substack{ I\in \mathcal{D}:\ell(I)\leq\ell(G)  \\ d(I,G)>2\ell(I)^\varepsilon\ell(G)^{1-\varepsilon}}}\left\vert
\left\langle T_{\sigma }^{\alpha }\left( \square _{I}^{\sigma ,\mathbf{b}%
}f\right),\mathbf{1}_{G}\beta _{\mathcal{G}}\left( G\right)
\right\rangle _{\omega }\right\vert  \notag \\
&\lesssim &
\sqrt{\mathfrak{A}_{2}^{\alpha }}\left\Vert f\right\Vert _{L^{2}\left(
\sigma \right) }\left\Vert g\right\Vert _{L^{2}\left( \omega \right) }. 
\notag
\end{eqnarray}
\end{lem}

The proof of this lemma is similar to those of Lemmas \ref{delta long} and %
\ref{delta short} in Section \ref{Sec disj form}\ above, using the square
function inequalities for $\square _{I}^{\sigma ,\mathbf{b}}$, $\nabla
_{I,\mathcal{F}}^\sigma$ and $\square _{J}^{\omega ,\mathbf{b}^{\ast }}$, $%
\nabla_{J,\mathcal{G}}^{\omega}$.

\begin{prop}[The Intertwining Proposition]
\label{strongly adapted}Let $\mathcal{D}$ and $\mathcal{G}$ be grids, and
suppose that $\mathbf{b}$ and $\mathbf{b}^{\ast }$ are $\infty $-weakly $%
\sigma $-accretive families of cubes in $\mathcal{D}$ and $\mathcal{G}$
respectively. Suppose that $\mathcal{F}\subset \mathcal{D}$ is $\sigma $%
-Carleson and that the $\mathcal{F}$-coronas 
\begin{equation*}
\mathcal{C}_{F}\equiv \left\{ I\in \mathcal{D}:I\subset F\text{ but }%
I\not\subset F^{\prime }\text{ for }F^{\prime }\in \mathfrak{C}_{\mathcal{F}%
}\left( F\right) \right\}
\end{equation*}%
satisfy%
\begin{equation*}
E_{I}^{\sigma }\left\vert f\right\vert \lesssim E_{F}^{\sigma }\left\vert
f\right\vert \text{ and }b_{I}=\mathbf{1}_{I}b_{F},\ \ \ \ \ \text{for all }%
I\in \mathcal{C}_{F}\mathfrak{\ ,\ }F\in \mathcal{F}.
\end{equation*}%
Then 
\begin{equation*}
\boldsymbol{E}_\Omega^\mathcal{D}\left\vert \sum_{F\in \mathcal{F}}\ \sum_{I:\ I\supsetneqq F}\ \left\langle
T_{\sigma }^{\alpha }\square _{I}^{\sigma ,\mathbf{b}}f,\mathsf{P}_{\mathcal{%
C}_{F}^{\mathcal{G},{shift}}}^{\omega ,\mathbf{b}^{\ast
}}g\right\rangle _{\omega }\right\vert \lesssim \left( \mathfrak{F}_{\alpha} \!+\!\mathfrak{T}^\textbf{b}_{T^\alpha}\!+\!\sqrt{\mathfrak{A}_2^\alpha}\delta^{\alpha-n}\! +\!\delta \mathfrak{N}_{T^\alpha}\right) \left\Vert f\right\Vert _{L^{2}\left( \sigma \right)
}\left\Vert g\right\Vert _{L^{2}\left( \omega \right) },
\end{equation*}%
where the implied constant depends on the $\sigma $-Carleson norm $C_{%
\mathcal{F}}$ of the family $\mathcal{F}$.
\end{prop}

\begin{proof}
We write the sum on the left hand side of the display above as%
\begin{eqnarray*}
\sum_{F\in \mathcal{F}}\ \sum_{I:\ I\supsetneqq F}\ \left\langle T_{\sigma}^{\alpha }\square _{I}^{\sigma ,\mathbf{b}}f,\mathsf{P}_{\mathcal{C}_{F}^{\mathcal{G},{shift}}}^{\omega }g\right\rangle _{\omega }
&\!\!\!\!\!=&\!\!\!
\sum_{F\in 
\mathcal{F}}\ \left\langle T_{\sigma }^{\alpha }\left( \sum_{I:\
I\supsetneqq F}\square _{I}^{\sigma ,\mathbf{b}}f\right) ,\mathsf{P}_{%
\mathcal{C}_{F}^{\mathcal{G},{shift}}}^{\omega }g\right\rangle
_{\omega }\\
&\!\!\!\!\!=&\!\!\!
\sum_{F\in \mathcal{F}}\ \left\langle T_{\sigma }^{\alpha }\left(
f^*_{F}\right) ,g_{F}\right\rangle _{\omega };
\end{eqnarray*}
where $\displaystyle f^*_{F}\equiv \sum_{I:\ I\supsetneqq F}\square _{I}^{\sigma ,%
\mathbf{b}}f\text{ and }g_{F}\equiv \mathsf{P}_{\mathcal{C}_{F}^{\mathcal{G},{shift}}}^{\omega }g$.

Note that $g_{F}$ is supported
in $F$. By the telescoping identity for $\square _{I}^{\sigma ,\mathbf{b}}$,
the function $f_{F}^{\ast }$ satisfies%
\begin{equation*}
\mathbf{1}_{F}f_{F}^{\ast }=\sum_{I:\ I_{\infty }\supset I\supsetneqq
F}\square _{I}^{\sigma ,\mathbf{b}}f=\mathbb{F}_{F}^{\sigma ,\mathbf{b}}f-%
\mathbf{1}_{F}\mathbb{F}_{I_{\infty }}^{\sigma ,\mathbf{b}}f=b_{F}\frac{%
E_{F}^{\sigma }f}{E_{F}^{\sigma }b_{F}}-\mathbf{1}_{F}b_{I_{\infty }}\frac{%
E_{I_{\infty }}^{\sigma }f}{E_{I_{\infty }}^{\sigma }b_{I_{\infty }}}\ .
\end{equation*}%
where $I_\infty$ is the starting cube for corona constructions in $\mathcal{D}$. However, we cannot apply the testing condition to the function $\mathbf{1}%
_{F}b_{I_{\infty }}$, and since $E_{I_{\infty }}^{\sigma }f$ does not vanish
in general, we will instead add and subtract the term $\mathbb{F}_{I_{\infty
}}^{\sigma ,\mathbf{b}}f$ to get 
\begin{eqnarray}\label{f*}
\\ \notag
\sum_{F\in \mathcal{F}}\ \left\langle T_{\sigma }^{\alpha }\left(
f_{F}^{\ast }\right) ,g_{F}\right\rangle _{\omega } 
&=&
\sum_{F\in \mathcal{F}%
}\ \left\langle T_{\sigma }^{\alpha }\left( \sum_{I:\ I_{\infty }\supset
I\supsetneqq F}\square _{I}^{\sigma ,\mathbf{b}}f\right) ,\mathsf{P}_{%
\mathcal{C}_{F}^{\mathcal{G},{shift}}}^{\omega }g\right\rangle
_{\omega } \\
&=&\notag
\sum_{F\in \mathcal{F}}\ \left\langle T_{\sigma }^{\alpha }\left( \mathbb{%
F}_{I_{\infty }}^{\sigma ,\mathbf{b}}f+\sum_{I:\ I_{\infty }\supset
I\supsetneqq F}\square _{I}^{\sigma ,\mathbf{b}}f\right) ,\mathsf{P}_{%
\mathcal{C}_{F}^{\mathcal{G},{shift}}}^{\omega }g\right\rangle
_{\omega } \\ 
&&\notag
-\sum_{F\in \mathcal{F}}\ \left\langle T_{\sigma }^{\alpha }\left( 
\mathbb{F}_{I_{\infty }}^{\sigma ,\mathbf{b}}f\right) ,\mathsf{P}_{\mathcal{C%
}_{F}^{\mathcal{G},{shift}}}^{\omega }g\right\rangle _{\omega }\ ,
\end{eqnarray}%
where the second sum on the right hand side of the identity satisfies
$$
\boldsymbol{E}_\Omega^\mathcal{D}\left\vert \sum_{F\in \mathcal{F}}\!\!\! \left\langle T_{\sigma }^{\alpha }\left( 
\mathbb{F}_{I_{\infty }}^{\sigma ,\mathbf{b}}f\right) ,\mathsf{P}_{\mathcal{C%
}_{F}^{\mathcal{G},{shift}}}^{\omega }g\right\rangle _{\omega
}\right\vert \!\!  \lesssim \!\! \left( \mathfrak{T}^\textbf{b}_{T^\alpha}\!+\!\sqrt{\mathfrak{A}_2^\alpha}\delta^{\alpha-n}\!+\!\delta\mathfrak{N}_{T^\alpha}\!\!\right)\left\Vert f\right\Vert_{L^2(\sigma)}\!\left\Vert g\right\Vert_{L^2(\omega)}
$$
Indeed, as 
\begin{eqnarray*}
 &&\sum_{F\in \mathcal{F}}\!\!\! \left\langle T_{\sigma }^{\alpha }\left( 
\mathbb{F}_{I_{\infty }}^{\sigma ,\mathbf{b}}f\right) ,\mathsf{P}_{\mathcal{C%
}_{F}^{\mathcal{G},{shift}}}^{\omega }g\right\rangle _{\omega
}\\
&= &
\left[\int_{I_\infty \cap J_\infty}+\int_{J_\infty \cap\left( (1+\delta)I_\infty\backslash I_\infty\right)} + \int_{J_\infty \backslash (1+\delta)I_\infty}   
\right]
\left( \sum_{F\in \mathcal{F}}\mathsf{P}_{\mathcal{C%
}_{F}^{\mathcal{G},{shift}}}^{\omega }g\right) T_{\sigma }^{\alpha }\left( 
\mathbb{F}_{I_{\infty }}^{\sigma ,\mathbf{b}}f\right) d\omega \\
&\equiv& A_1+A_2+A_3
\end{eqnarray*}
by Cauchy-Schwarz and Riesz inequalities, the term $A_1$ is controlled by testing, the term $A_3$ by Muckenhoupt's condition using lemma \ref{lemma1} and finally
$$
\boldsymbol{E}_\Omega^\mathcal{D}A_2 \leq \left( C\delta \int_{I_{\infty }}\left\vert\sum_{F\in \mathcal{F}}\mathsf{P}_{\mathcal{C%
}_{F}^{\mathcal{G},{shift}}}^{\omega }g\right\vert
^{2}d\omega \right) ^{\frac{1}{2}}\left( \mathfrak{N}_{T^{\alpha }}\int
\left\vert f\right\vert ^{2}d\sigma \right) ^{\frac{1}{2}}
\leq
\sqrt{C\delta 
\mathfrak{N}_{T^{\alpha }}}\left\Vert f\right\Vert _{L^{2}\left( \sigma \right)}\left\Vert g\right\Vert _{L^{2}\left( \omega \right)}\ .
$$
The advantage now is that with
\begin{equation*}
f_{F}\equiv \mathbb{F}_{I_{\infty }}^{\sigma ,\mathbf{b}}f+f_{F}^{\ast }=%
\mathbb{F}_{I_{\infty }}^{\sigma ,\mathbf{b}}f+\sum_{I:\ I_{\infty }\supset
I\supsetneqq F}\square _{I}^{\sigma ,\mathbf{b}}f
\end{equation*}%
then in the first term on the right hand side of (\ref{f*}), the
telescoping identity gives%
\begin{equation*}
\mathbf{1}_{F}f_{F}=\mathbf{1}_{F}\left( \mathbb{F}_{I_{\infty }}^{\sigma ,%
\mathbf{b}}f+\sum_{I:\ I_{\infty }\supset I\supsetneqq F}\square
_{I}^{\sigma ,\mathbf{b}}f\right) =\mathbb{F}_{F}^{\sigma ,\mathbf{b}}f=b_{F}%
\frac{E_{F}^{\sigma }f}{E_{F}^{\sigma }b_{F}},
\end{equation*}%
which shows that $f_{F}$ is a controlled constant times $b_{F}$ on $F$.

The cubes $I$ occurring in this sum are linearly and consecutively
ordered by inclusion, along with the cubes $F^{\prime }\in \mathcal{F}$
that contain $F$. More precisely we can write%
\begin{equation*}
F\equiv F_{0}\subsetneqq F_{1}\subsetneqq F_{2}\subsetneqq ...\subsetneqq
F_{n}\subsetneqq F_{n+1}\subsetneqq ...F_{N}=I_{\infty }
\end{equation*}%
where $F_{m}=\pi _{\mathcal{F}}^{m}F$ for all $m\geq 1$. We can also write%
\begin{equation*}
F=F_{0}\equiv I_{0}\subsetneqq I_{1}\subsetneqq I_{2}\subsetneqq
...\subsetneqq I_{k}\subsetneqq I_{k+1}\subsetneqq ...\subsetneqq
I_{K}=F_{N}=I_{\infty }
\end{equation*}%
where $I_{k}=\pi _{\mathcal{D}}^{k}F$ for all $k\geq 1$. There is a (unique)
subsequence $\left\{ k_{m}\right\} _{m=1}^{N}$ such that%
\begin{equation*}
F_{m}=I_{k_{m}},\ \ \ \ \ 1\leq m\leq N.
\end{equation*}

Then we have%
\begin{eqnarray*}
f_{F}\left( x\right) \!\equiv \mathbb{F}_{I_{\infty }}^{\sigma ,\mathbf{b}%
}f\left( x\right) +\sum_{\ell =1}^{K}\square _{I_{\ell }}^{\sigma ,\mathbf{b}%
}f\left( x\right)\ \ \text{ and } \ \ 
g_{F} \equiv \!\!\!\!\!\!\!\! \sum_{J\in \mathcal{C}_{F}^{\mathcal{G},{shift}%
}}\square _{J}^{\omega ,\mathbf{b}^{\ast }}g.
\end{eqnarray*}%
Assume now that $k_{m}\leq k<k_{m+1}$. We denote by $\theta \left( I\right) $ the $2^n-1$ siblings of $I$, i.e. $ \tilde{I}\in \theta \left( I\right)$ implies  $\tilde{I} \in %
\mathfrak{C}_{\mathcal{D}}\left( \pi _{\mathcal{D}}I\right) \backslash
\left\{ I\right\} $. There are two cases to consider here:%
\begin{equation*}
 \tilde{I}_k \notin \mathcal{F}\text{ and }\tilde{I}_k \in \mathcal{F}.
\end{equation*}%
We first note that in either case, using a telescoping sum, we compute that
for 
\begin{equation*}
x\in \tilde{I}_k \subset
F_{m+1}\backslash F_{m},
\end{equation*}%
we have the formula 
\begin{eqnarray*}
f_{F}\left( x\right) &=&\mathbb{F}_{I_{\infty }}^{\sigma ,\mathbf{b}}f\left(
x\right) +\sum_{\ell = k+1}^K \square _{I_{\ell }}^{\sigma ,\mathbf{b}%
}f\left( x\right) \\
&=&
\mathbb{F}_{\tilde{I}_k }^{\sigma ,\mathbf{b}}f\left(
x\right) -\mathbb{F}_{I_{k+1}}^{\sigma ,\mathbf{b}}f\left( x\right)
+
\sum_{\ell =k+1}^{K-1}\left( \mathbb{F}_{I_{\ell }}^{\sigma ,\mathbf{b}%
}f\left( x\right) -\mathbb{F}_{I_{\ell +1}}^{\sigma ,\mathbf{b}}f\left(
x\right) \right) +\mathbb{F}_{I_{\infty }}^{\sigma ,\mathbf{b}}f\left(
x\right) \\
&=&
\mathbb{F}_{\tilde{I}_k }^{\sigma ,\mathbf{b}%
}f\left( x\right) \ .
\end{eqnarray*}%
Now fix $x\in \tilde{I}_k $. If $\tilde{I}_k
\notin \mathcal{F}$, then $\tilde{I}_k \in \mathcal{C}%
_{F_{m+1}}$, and we have 
\begin{equation}
\left\vert f_{F}\left( x\right) \right\vert =\left\vert \mathbb{F}_{\tilde{I}_k }^{\sigma ,\mathbf{b}}f\left( x\right) \right\vert
\lesssim \left\vert b_{\tilde{I}_k }\left( x\right)
\right\vert \ \frac{E_{\tilde{I}_k }^{\sigma }\left\vert
f\right\vert }{\left\vert E_{\tilde{I}_k }^{\sigma }b_{\theta
\left( I_{k}\right) }\right\vert }\lesssim E_{F_{m+1}}^{\sigma }\left\vert
f\right\vert \ ,  \label{bound for f_F}
\end{equation}%
since the testing functions $b_{\tilde{I}_k }$ are bounded
and accretive, and $E_{\tilde{I}_k }^{\sigma }\left\vert
f\right\vert \lesssim E_{F_{m+1}}^{\sigma }\left\vert f\right\vert $ by
hypothesis. On the other hand, if $\tilde{I}_k \in \mathcal{F}
$, then $I_{k+1}\in \mathcal{C}_{F_{m+1}}$ and we have%
\begin{equation*}
\left\vert f_{F}\left( x\right) \right\vert =\left\vert \mathbb{F}_{\tilde{I}_k }^{\sigma ,\mathbf{b}}f\left( x\right) \right\vert
\lesssim E_{\tilde{I}_k }^{\sigma }\left\vert f\right\vert \ .
\end{equation*}%
Note that $\displaystyle F^{c}=\overset{\cdot }{\bigcup_{k \geq 0}}\theta\left(I_k\right) $. Now we write%
\begin{eqnarray*}
f_{F}
\!\!\!\!\!&=&\!\!\!\!\!
\varphi _{F}+\psi _{F}, \\
\varphi _{F} \equiv \sum_{k\geq 0}\sum_{\substack{\tilde{I}_k \in \theta (I_k)\\\tilde{I}_k\in \mathcal{F}}}
\mathbb{F}_{\tilde{I}_k}^{\sigma ,\mathbf{b}}f
&\text{ and }&
\psi _{F}=f_{F}-\varphi _{F}\ ; \\
\sum_{F\in \mathcal{F}}\ \left\langle T_{\sigma }^{\alpha
}f_{F},g_{F}\right\rangle _{\omega } 
\!\!\!&=&\!\!\!
\sum_{F\in \mathcal{F}}\
\left\langle T_{\sigma }^{\alpha }\varphi _{F},g_{F}\right\rangle _{\omega
}+\sum_{F\in \mathcal{F}}\ \left\langle T_{\sigma }^{\alpha }\psi
_{F},g_{F}\right\rangle _{\omega }\ ,
\end{eqnarray*}%
and note that $\varphi _{F}=0$ on $F$, and $\psi _{F}=b_{F}\frac{%
E_{F}^{\sigma }f}{E_{F}^{\sigma }b_{F}}$ on $F$. We can apply the first line
in (\ref{indicator far}) using $\tilde{I}_k \in \mathcal{F}$
to the first sum above since $J\in \mathcal{C}_{F}^{\mathcal{G},{%
shift}}$ implies $J\subset J^{\maltese }\subset F\subset I_{k}$, which implies that $d( J,\tilde{I}_k
) >2\ell \left( J\right) ^{\varepsilon }\ell ( \tilde{I}_k ) ^{1-\varepsilon }$. Thus we obtain after substituting $%
F^{\prime }$ for $\tilde{I}_k$ below, 
\begin{eqnarray*}
\left\vert \sum_{F\in \mathcal{F}}\ \left\langle T_{\sigma }^{\alpha
}\varphi _{F},g_{F}\right\rangle _{\omega }\right\vert &=&
\left\vert
\sum_{F\in \mathcal{F}}\sum_{J\in \mathcal{C}_{F}^{\mathcal{G},{shift%
}}}\left\langle T_{\sigma }^{\alpha }\left( \sum_{k\geq 0}\sum_{\substack{\tilde{I}_k \in \theta (I_k)\\\tilde{I}_k\in \mathcal{F}}}
\mathbb{F}_{\tilde{I}_k}^{\sigma ,\mathbf{b}}f\right) ,\square _{J}^{\omega ,\mathbf{b}^{\ast
}}g\right\rangle _{\omega }\right\vert \\
&\leq &
\sum_{F\in \mathcal{F}}\sum_{J\in \mathcal{C}_{F}^{\mathcal{G},%
{shift}}} \sum_{k\geq 0}\sum_{\substack{\tilde{I}_k \in \theta (I_k)\\\tilde{I}_k\in \mathcal{F}}}\left\vert \left\langle T_{\sigma }^{\alpha }\left( \mathbb{F}_{\tilde{I}_k }^{\sigma ,\mathbf{b}}f\right) ,\square _{J}^{\omega ,
\mathbf{b}^{\ast }}g\right\rangle _{\omega }\right\vert \\
&\leq &
\sum_{F^{\prime }\in \mathcal{F}}\sum_{\substack{ J\in \mathcal{G}:\
\ell \left( J\right) \leq \ell \left( F^{\prime }\right)  \\ d\left(
J,F^{\prime }\right) >2\ell \left( J\right) ^{\varepsilon }\ell \left(
F^{\prime }\right) ^{1-\varepsilon }}}\left\vert \left\langle T_{\sigma
}^{\alpha }\left( \mathbb{F}_{F^{\prime }}^{\sigma ,\mathbf{b}}f\right)
,\square _{J}^{\omega ,\mathbf{b}^{\ast }}g\right\rangle _{\omega
}\right\vert \\
&\lesssim &
\sqrt{\mathfrak{A}_{2}^{\alpha }}\left\Vert f\right\Vert _{L^{2}\left(
\sigma \right) }\left\Vert g\right\Vert _{L^{2}\left( \omega \right) }\ .
\end{eqnarray*}

Turning to the second sum, we note that for $k_{m}\leq k<k_{m+1}$ and $x\in
\tilde{I}_k $ with $\tilde{I}_k \notin 
\mathcal{F}$, we have 
\begin{equation*}
\left\vert \psi _{F}\left( x\right) \right\vert \lesssim \left\vert
b_{\tilde{I}_k }\right\vert\  E_{\tilde{I}_k
}^{\sigma }\left\vert f\right\vert \ \mathbf{1}_{\tilde{I}_k} \left( x\right) \lesssim \alpha _{\mathcal{F}}\left( F_{m+1}\right)
\mathbf{1}_{\tilde{I}_k}\left( x\right)
\end{equation*}%
Note that for $\sigma$-almost all $x \in I_\infty$ there exists a unique $F\in \mathcal{F}$ such that $\displaystyle x\in F\backslash \bigcup_{F'\in \mathfrak{C}_\mathcal{F}(F)}F'$ since the family $\mathcal{F}$ is a Carleson family. Also from the stopping criteria we have $\alpha_\mathcal{F}(F)\leq\alpha_\mathcal{F}(F')$ for $F' \subset F$. Hence we get the following inequality for $x\notin F$, 
\begin{equation}
\left\vert \psi _{F}\left( x\right) \right\vert 
\lesssim
\Phi \left( x\right) \ \mathbf{1}%
_{F^{c}}\left( x\right) \ ,  \label{Psi_F bound}
\end{equation}%
where we have defined 
\begin{equation*}
\Phi \equiv \sum_{F\in 
\mathcal{F}}\alpha _{\mathcal{F}}\left( F\right) \mathbf{1}_{F\backslash \cup 
\mathfrak{C}_{\mathcal{F}}\left( F\right) }\ .
\end{equation*}

Now we write%
\begin{equation*}
\sum_{F\in \mathcal{F}}\ \left\langle T_{\sigma }^{\alpha }\psi
_{F},g_{F}\right\rangle _{\omega }=\sum_{F\in \mathcal{F}}\ \left\langle
T_{\sigma }^{\alpha }\left( \mathbf{1}_{F}\psi _{F}\right)
,g_{F}\right\rangle _{\omega }+\sum_{F\in \mathcal{F}}\ \left\langle
T_{\sigma }^{\alpha }\left( \mathbf{1}_{F^{c}}\psi _{F}\right)
,g_{F}\right\rangle _{\omega }\equiv\mathrm{I+II}.
\end{equation*}%
Then by cube
testing, 
\begin{equation*}
\left\vert \left\langle T_{\sigma }^{\alpha }\left( b_{F}\mathbf{1}%
_{F}\right) ,g_{F}\right\rangle _{\omega }\right\vert =\left\vert
\left\langle \mathbf{1}_{F}T_{\sigma }^{\alpha }\left( b_{F}\mathbf{1}%
_{F}\right) ,g_{F}\right\rangle _{\omega }\right\vert \lesssim \mathfrak{T}%
_{T^{\alpha }}\sqrt{\left\vert F\right\vert _{\sigma }}\left\Vert
g_{F}\right\Vert _{L^{2}\left( \omega \right) }^{\bigstar }\ ,
\end{equation*}%
and so quasi-orthogonality, together with the fact that on $F$, $\psi
_{F}=b_{F}\frac{E_{F}^{\sigma }f}{E_{F}^{\sigma }b_{F}}$ is a constant $c=%
\frac{E_{F}^{\sigma }f}{E_{F}^{\sigma }b_{F}}$ times $b_{F}$, where $%
\left\vert c\right\vert $ is bounded by $\alpha _{\mathcal{F}}\left(
F\right) $, give
\begin{eqnarray*}
\left\vert\mathrm{I}\right\vert =\left\vert \sum_{F\in \mathcal{F}}\ \left\langle
T_{\sigma }^{\alpha }\left( \mathbf{1}_{F}cb_{F}\right) ,g_{F}\right\rangle
_{\omega }\right\vert
\lesssim
\sum_{F\in \mathcal{F}}\ \alpha _{\mathcal{F}%
}\left( F\right) \ \Big\vert \left\langle T_{\sigma }^{\alpha
}b_{F},g_{F}\right\rangle _{\omega }\Big\vert
\!\!\!\!&\lesssim &\!\!\!
\sum_{F\in \mathcal{F}}\ \alpha _{\mathcal{F}}\left( F\right) 
\mathfrak{T}_{T^\alpha }\sqrt{\left\vert F\right\vert _{\sigma }}\left\Vert
g_{F}\right\Vert _{L^{2}\left( \omega \right) }^{\bigstar }\\
&\lesssim&\!\!\!
\mathfrak{T}_{T^\alpha }\left\Vert f\right\Vert _{L^{2}\left( \sigma \right) }\left[
\sum_{F\in \mathcal{F}}\left\Vert g_{F}\right\Vert _{L^{2}\left( \omega
\right) }^{\bigstar 2}\right] ^{\frac{1}{2}}
\end{eqnarray*}

Now $\mathbf{1}_{F^{c}}\psi _{F}$ is supported outside $F$, and each $J$ in
the dual martingale support $\mathcal{C}_{F}^{\mathcal{G},{shift}}$
of $g_{F}=\mathsf{P}_{\mathcal{C}_{F}^{\mathcal{G},{shift}}}^{\omega
}g$ is in particular ${good}$ in the cube $F$, and as a
consequence, each such cube $J$ as above is contained in some cube $%
M $ for $M\in \mathcal{W}\left( F\right) $. This containment will be used in
the analysis of the term $\mathrm{II_G}$ below.

In addition, each $J$ in the dual martingale support $\mathcal{C}_{F}^{%
\mathcal{G},{shift}}$ of $g_{F}=\mathsf{P}_{\mathcal{C}_{F}^{%
\mathcal{G},{shift}}}^{\omega }g$ is $\left( \left[ \frac{3}{%
\varepsilon }\right] ,\varepsilon \right) $-deeply embedded in $F$, i.e. $%
J\Subset _{\left[ \frac{3}{\varepsilon }\right] ,\varepsilon }F$ the definition of $\mathcal{C}_{F}^{\mathcal{G},{shift}}$. As a consequence, each such
cube $J$ as above is contained in some cube $M$ for $M\in \mathcal{M}%
_{\left( \left[ \frac{3}{\varepsilon }\right] ,\varepsilon \right) -{%
deep},\mathcal{D}}\left( F\right) $. This containment will be used in the
analysis of the term $\mathrm{II_B}$ below.

\begin{notation}
\label{Notation rho}
Define $\mathbf{\rho }\equiv \left[ \frac{3}{\varepsilon }\right] $, so that
for every $J\in \mathcal{C}_{F}^{\mathcal{G},{shift}}$, there is $%
M\in \mathcal{M}_{\left( \mathbf{\rho },\varepsilon \right) -{deep},%
\mathcal{G}}\left( F\right) $ such that $J\subset M$.
\end{notation}

The collections $\mathcal{W}\left( F\right) $ and $\mathcal{M}_{\left( 
\mathbf{\rho },\varepsilon \right) -{deep},\mathcal{G}}\left(
F\right) $ used here, and in the display below, are defined in (\ref{def
M_r-deep}) in Appendix. Finally, since the cubes $M\in \mathcal{W}%
\left( F\right) $, as well as the cubes $M\in \mathcal{M}_{\left( \left[ 
\frac{3}{\varepsilon }\right] ,\varepsilon \right) -{deep},\mathcal{G%
}}\left( F\right) $, satisfy $3M\subset F$, we can apply (\ref{estimate}) in
the Monotonicity Lemma \ref{mono} using (\ref{Psi_F bound}) with $\mu =%
\mathbf{1}_{F^{c}}\psi _{F}$ and $J^{\prime }$ in place of $J$ there, to
obtain%
\begin{eqnarray*}
\left\vert \mathrm{II}\right\vert
&=&
\left\vert \sum_{F\in \mathcal{F}}\left\langle
T_{\sigma }^{\alpha }\left( \mathbf{1}_{F^{c}}\psi _{F}\right)
,g_{F}\right\rangle _{\omega }\right\vert =\left\vert \sum_{F\in \mathcal{F}%
}\sum_{J^{\prime }\in \mathcal{C}_{F}^{\mathcal{G},{shift}%
}}\left\langle T_{\sigma }^{\alpha }\left( \mathbf{1}_{F^{c}}\psi
_{F}\right) ,\square _{J^{\prime }}^{\omega ,\mathbf{b}^{\ast
}}g\right\rangle _{\omega }\right\vert \\
&\lesssim &
\sum_{F\in \mathcal{F}}\sum_{J^{\prime }\in \mathcal{C}_{F}^{%
\mathcal{G},{shift}}}\frac{\mathrm{P}^{\alpha }\left( J^{\prime },%
\mathbf{1}_{F^{c}}|\psi _{F}| \sigma \right) }{\left\vert J^{\prime }\right\vert^\frac{1}{n} }%
\left\Vert \bigtriangleup _{J^{\prime }}^{\omega ,\mathbf{b}^{\ast
}}x\right\Vert _{L^{2}\left( \omega \right) }^{\spadesuit }\left\Vert
\square _{J^{\prime }}^{\omega ,\mathbf{b}^{\ast }}g\right\Vert
_{L^{2}\left( \omega \right) }^{\bigstar } \\
&&
+\sum_{F\in \mathcal{F}}\sum_{J^{\prime }\in \mathcal{C}_{F}^{\mathcal{G},%
{shift}}}\frac{\mathrm{P}_{1+\delta }^{\alpha }\left( J^{\prime },%
\mathbf{1}_{F^{c}}|\psi _{F}| \sigma \right) }{\left\vert J^{\prime }\right\vert^\frac{1}{n} }%
\left\Vert x-m_{J^{\prime }}\right\Vert _{L^{2}\left( \omega \right)
}\left\Vert \square _{J^{\prime }}^{\omega ,\mathbf{b}^{\ast }}g\right\Vert
_{L^{2}\left( \omega \right) }^{\bigstar } \\
&\lesssim &
\sum_{F\in \mathcal{F}}\sum_{M\in \mathcal{W}\left( F\right) }%
\frac{\mathrm{P}^{\alpha }\left( M,\mathbf{1}_{F^{c}}\Phi \sigma \right) }{%
\left\vert M\right\vert^\frac{1}{n} }\left\Vert \mathsf{Q}_{\mathcal{C}_{F;M}^{\mathcal{G%
},{shift}}}^{\omega ,\mathbf{b}^{\ast }}x\right\Vert _{L^{2}\left(
\omega \right) }^{\spadesuit }\left\Vert g_{F;M}\right\Vert _{L^{2}\left(
\omega \right) }^{\bigstar } \\
&&\!\!\!\!\!\!\!+
\sum_{F\in \mathcal{F}}\sum_{J\in \mathcal{M}^{deep}_{\left( \mathbf{\rho }%
,\varepsilon \right),\mathcal{G}}\left( F\right)
}\sum_{J^{\prime }\in \mathcal{C}_{F;J}^{\mathcal{G},{shift}}}
\!\!\!\!\!\!\!\!\frac{
\mathrm{P}_{1+\delta }^{\alpha }\left( J^{\prime },\mathbf{1}_{F^{c}}|\psi _{F}|
\sigma \right) }{\left\vert J^{\prime }\right\vert^\frac{1}{n} }\!\left\Vert
x-m_{J^{\prime }}\right\Vert\! _{L^{2}\left( \mathbf{1}_{J^{\prime }}\omega
\right) }\!\!\left\Vert \square _{J^{\prime }}^{\omega ,\mathbf{b}^{\ast
}}g\right\Vert\! _{L^{2}\left(\omega\! \right) }^{\bigstar } \\
&\equiv &\mathrm{II_G}+\mathrm{II_B}\ .
\end{eqnarray*}%
where $g_{F;M}$ denotes the pseudoprojection $g_{F;M}=\!\!\!\!\!\displaystyle \sum_{J^{\prime }\in 
\mathcal{C}_{F}^{\mathcal{G},{shift}}:\ J^{\prime }\subset M}\square
_{J^{\prime }}^{\omega ,\mathbf{b}^{\ast }}g$.

\textbf{Note}: We could also bound $\mathrm{II_G}$ by using the decomposition $%
\mathcal{M}_{\left( \mathbf{\rho },\varepsilon \right) -{deep},%
\mathcal{G}}\left( F\right) $ of $F$ into certain maximal $\mathcal{G}$%
-cubes, but the `smaller' choice $\mathcal{W}\left( F\right) $ of $%
\mathcal{D}$-cubes is needed for $\mathrm{II_G}$ in order to bound it by the
corresponding functional energy constant $\mathfrak{F}_{\alpha }$, which can then be controlled by the energy and Muckenhoupt
constants in Appendix .

Then from Cauchy-Schwarz, the functional energy condition, and 
\begin{equation*}
\left\Vert \Phi \right\Vert _{L^{2}\left( \sigma \right) }^{2}\leq
\sum_{F\in \mathcal{F}}\alpha _{\mathcal{F}}\left( F\right) ^{2}\left\vert
F\right\vert _{\sigma }\lesssim \left\Vert f\right\Vert _{L^{2}\left( \sigma
\right) }^{2}\ ,
\end{equation*}%
we obtain%
\begin{eqnarray*}
\left\vert \mathrm{II_G}\right\vert 
\!\!\!\!\!&\leq &\!\!\!\!\!
\left(\!\sum_{F\in \mathcal{F}}\!\sum_{M\in 
\mathcal{W}\left( F\right) }\!\!\!\!\left( \frac{\mathrm{P}^{\alpha }\left( M,
\mathbf{1}_{F^{c}}\Phi \sigma \right) }{\left\vert M\right\vert }\right)
^{\!\!2}\!\left\Vert \mathsf{Q}_{\mathcal{C}_{F;M}^{\mathcal{G},{shift}
}}^{\omega ,\mathbf{b}^{\ast }}x\right\Vert _{L^{2}\left( \omega \right)
}^{\spadesuit 2}\!\!\right) ^{\!\!\!\!\frac{1}{2}} 
\!\!\!\left(\!\sum_{F\in \mathcal{F}}\!\sum_{M\in \mathcal{W}\left( F\right)
}\!\!\!\!\!\!\left\Vert g_{F;M}\right\Vert _{L^{2}\left( \omega \right) }^{\bigstar
2}\!\!\!\right) ^{\!\!\!\!\frac{1}{2}} \\
&\lesssim &
\mathfrak{F}_{\alpha }\left\Vert \Phi
\right\Vert _{L^{2}\left( \sigma \right) }\left[ \sum_{F\in \mathcal{F}
}\left\Vert g_{F}\right\Vert _{L^{2}\left( \omega \right) }^{\bigstar 2}
\right] ^{\frac{1}{2}}
\lesssim 
\mathfrak{F}_{\alpha }\left\Vert f\right\Vert _{L^{2}\left( \sigma \right) }\left\Vert
g\right\Vert _{L^{2}\left( \omega \right) },
\end{eqnarray*}
by the pairwise disjointedness of the coronas $\mathcal{C}_{F;M}^{\mathcal{G}%
,{shift}}$ jointly in $F$ and $M$, which in turn follows from the
pairwise disjointedness (\ref{tau overlap}) of the shifted coronas $\mathcal{%
C}_{F}^{\mathcal{G},{shift}}$ in $F$, together with the pairwise
disjointedness of the cubes $M$. Thus we obtain the pairwise disjointedness
of both of the pseudoprojections $\mathsf{P}_{\mathcal{C}_{F;M}^{\mathcal{G},%
{shift}}}^{\omega ,\mathbf{b}^{\ast }}$ and $\mathsf{Q}_{\mathcal{C}%
_{F;M}^{\mathcal{G},{shift}}}^{\omega ,\mathbf{b}^{\ast }}$ jointly
in $F$ and $M$.

In term $\mathrm{II_B}$ the quantities $\left\Vert x-m_{J^{\prime }}\right\Vert
_{L^{2}\left( \mathbf{1}_{J^{\prime }}\omega \right) }^{2}$ are no longer
additive except when the cubes $J^{\prime }$ are pairwise disjoint. As a
result we will use (\ref{Haar trick}) in the form,%
\begin{eqnarray}
\\
\notag \sum_{J^{\prime }\subset J}\left( \frac{\mathrm{P}_{1+\delta }^{\alpha
}\left( J^{\prime },\nu \right) }{\left\vert J^{\prime }\right\vert^\frac{1}{n} }\right)
^{2}\left\Vert x-m_{J^{\prime }}\right\Vert _{L^{2}\left( \mathbf{1}%
_{J^{\prime }} \right) }^{2} 
&\lesssim &
\frac{1}{\gamma ^{2\delta
^{\prime }}}\left( \frac{\mathrm{P}_{1+\delta ^{\prime }}^{\alpha }\left(
J,\nu \right) }{\left\vert J\right\vert^\frac{1}{n} }\right) ^{2}\sum_{J^{\prime \prime
}\subset J}\left\Vert \bigtriangleup _{J^{\prime \prime }}^{\omega
}x\right\Vert _{L^{2} }^{2}  \label{Haar trick'} \\
&\lesssim &
\left( \frac{\mathrm{P}_{1+\delta ^{\prime }}^{\alpha }\left(
J,\nu \right) }{\left\vert J\right\vert ^\frac{1}{n}}\right) ^{2}\left\Vert
x-m_{J}\right\Vert _{L^{2}\left( \mathbf{1}_{J} \right) }^{2}\ , 
\notag
\end{eqnarray}%
and exploit the decay in the Poisson integral $\mathrm{P}_{1+\delta ^{\prime
}}^{\alpha }$ along with weak goodness of the cubes $J$. As a
consequence we will be able to bound $\mathrm{II_B}$ \emph{directly} by the strong
energy condition (\ref{strong b* energy}), without having to invoke the more
difficult functional energy condition. For the decay we compute that for $%
J\in \mathcal{M}_{\left( \mathbf{\rho },\varepsilon \right) -{deep},%
\mathcal{G}}\left( F\right) $%
\begin{eqnarray*}
\frac{\mathrm{P}_{1+\delta ^{\prime }}^{\alpha }\left( J,\mathbf{1}%
_{F^{c}}|\psi _{F}| \sigma \right) }{\left\vert J\right\vert^\frac{1}{n} } 
&\approx&
\int_{F^{c}}\frac{\left\vert J\right\vert ^\frac{\delta'}{n}}{\left\vert
y-c_{J}\right\vert ^{n+1+\delta ^{\prime }-\alpha }}|\psi _{F}| \left( y\right)
d\sigma  \\
&\leq &
\sum_{t=0}^{\infty }\int_{\pi _{\mathcal{F}}^{t+1}F\backslash \pi _{%
\mathcal{F}}^{t}F}\left( \frac{\left\vert J\right\vert^\frac{1}{n} }{\dist
\left( c_{J},\left( \pi _{\mathcal{F}}^{t}F\right) ^{c}\right) }\right)
^{\delta ^{\prime }}\!\!\!\frac{|\psi _{F}|\left( y\right)}{\left\vert y-c_{J}\right\vert ^{n+1-\alpha }} d\sigma  \\
&\lesssim &
\sum_{t=0}^{\infty }\left( \frac{\left\vert J\right\vert^\frac{1}{n} }{%
\dist\left( c_{J},\left( \pi _{\mathcal{F}}^{t}F\right) ^{c}\right) 
}\right) ^{\delta ^{\prime }}\frac{\mathrm{P}^{\alpha }\left( J,\mathbf{1}
_{\pi _{\mathcal{F}}^{t+1}F\backslash \pi _{\mathcal{F}}^{t}F}|\psi _{F}| \sigma
\right) }{\left\vert J\right\vert^\frac{1}{n} },
\end{eqnarray*}%
and then use the weak goodness inequality and the fact that $J\subset F$
\begin{equation*}
\dist\left( c_{J},\left( \pi _{\mathcal{F}}^{t}F\right) ^{c}\right)
\geq 2\ell \left( \pi _{\mathcal{F}}^{t}F\right) ^{1-\varepsilon }\ell
\left( J\right) ^{\varepsilon }\geq 2\cdot 2^{t\left( 1-\varepsilon \right)
}\ell \left( F\right) ^{1-\varepsilon }\ell \left( J\right) ^{\varepsilon
}\geq 2^{t\left( 1-\varepsilon \right) +1}\ell \left( J\right) ,
\end{equation*}%
to conclude that%
\begin{eqnarray}
\ \ \ \ \ \left( \frac{\mathrm{P}_{1+\delta ^{\prime }}^{\alpha }\left( J,\mathbf{1}%
_{F^{c}}|\psi _{F}| \sigma \right) }{\left\vert J\right\vert^\frac{1}{n} }\right) ^{2}
\!\!\!\!\!\!&\lesssim &\!\!\!\!\!\!
\left( \sum_{t=0}^{\infty }2^{-t\delta ^{\prime }\left(
1-\varepsilon \right) }\frac{\mathrm{P}^{\alpha }\left( J,\mathbf{1}_{\pi _{%
\mathcal{F}}^{t+1}F\backslash \pi _{\mathcal{F}}^{t}F}|\psi_{F}| \sigma \right) }{%
\left\vert J\right\vert^\frac{1}{n} }\right) ^{2}  \label{decay in t} \\
&\lesssim &
\sum_{t=0}^{\infty }2^{-t\delta ^{\prime }\left( 1-\varepsilon
\right) }\left( \frac{\mathrm{P}^{\alpha }\left( J,\mathbf{1}_{\pi _{%
\mathcal{F}}^{t+1}F\backslash \pi _{\mathcal{F}}^{t}F}|\psi _{F}| \sigma \right) }{\left\vert J\right\vert^\frac{1}{n} }\right) ^{2}.  \notag
\end{eqnarray}%
where in the last inequality we used the Cauchy-Schwarz inequality. Now we again apply Cauchy-Schwarz and (\ref{Haar trick'}) to obtain%

\begin{eqnarray*}
\mathrm{II_B} 
\!\!\!\!&=&\!\!\!\!
\sum_{F\in \mathcal{F}}\sum_{J\in \mathcal{M}^{deep}_{\left( \mathbf{\rho },\varepsilon \right),\mathcal{G}}\left( F\right)
}\sum_{J^{\prime }\in \mathcal{C}_{F;J}^{\mathcal{G},{shift}}}
\!\!\!\!\!\!\!\frac{\mathrm{P}_{1+\delta }^{\alpha }\left(J^{\prime},\mathbf{1}_{F^{c}}|\psi _{F}|
\sigma \right) }{\left\vert J^{\prime }\right\vert^\frac{1}{n} }\!\left\Vert
x-m_{J^{\prime }}\right\Vert _{L^{2}\left( \mathbf{1}_{J^{\prime }}\omega
\right) }\!\left\Vert \square _{J^{\prime }}^{\omega ,\mathbf{b}^{\ast
}}g\right\Vert _{L^{2}\left( \omega \right) }^{\bigstar } \\
&\leq & 
\!\!\!\left( \sum_{F\in \mathcal{F}}\sum_{J\in \mathcal{M}^{deep}_{\left( \mathbf{\rho },\varepsilon \right),\mathcal{G}}\left( F\right)
}\sum_{J^{\prime }\in \mathcal{C}_{F;J}^{\mathcal{G},{shift}}} 
\!\!\left( \frac{\mathrm{P}_{1+\delta }^{\alpha }\left( J^{\prime },\mathbf{1}
_{F^{c}}|\psi _{F}| \sigma \right) }{\left\vert J^{\prime }\right\vert^\frac{1}{n} }\right)
^{2}\left\Vert x-m_{J^{\prime }}\right\Vert _{L^{2}\left( \mathbf{1}%
_{J^{\prime }}\omega \right) }^{2} \right)^{\frac{1}{2}}\!\!
\left[
\sum_{F\in\mathcal{F}}\left\Vert g_{F}\right\Vert _{L^{2}\left( \omega \right) }^{\bigstar
2}\right] ^{\!\frac{1}{2}} \\
&\leq &
\!\!\!\left( \sum_{F\in \mathcal{F}}\sum_{J\in \mathcal{M}_{\left( \mathbf{%
\rho },\varepsilon \right) -{deep},\mathcal{G}}\left( F\right)
}\left( \frac{\mathrm{P}_{1+\delta ^{\prime }}^{\alpha }\left( J,\mathbf{1}
_{F^{c}}|\psi _{F}| \sigma \right) }{\left\vert J\right\vert ^\frac{1}{n}}\right)
^{2}\left\Vert x-m_{J}\right\Vert _{L^{2}\left( \mathbf{1}_{J}\omega \right)
}^{2}\right) ^{\!\!\frac{1}{2}}\!\!\left\Vert g\right\Vert _{L^{2}(\omega)} \\
&\equiv &
\sqrt{\mathrm{II_{energy}}}\left\Vert g\right\Vert _{L^{2}\left(
\omega \right) },
\end{eqnarray*}%
and it remains to estimate $\mathrm{II_{energy}}$. From (\ref{decay in t})
and the strong energy condition (\ref{strong b* energy}), we have%
\begin{eqnarray*}
\!\!\!&&\!\!\!
\mathrm{II_{energy}}=\sum_{F\in \mathcal{F}}\sum_{J\in \mathcal{M}%
_{\left( \mathbf{\rho },\varepsilon \right) -{deep},\mathcal{G}%
}\left( F\right) }\left( \frac{\mathrm{P}_{1+\delta ^{\prime }}^{\alpha
}\left( J,\mathbf{1}_{F^{c}}|\psi _{F}| \sigma \right) }{\left\vert J\right\vert^\frac{1}{n} }%
\right) ^{2}\left\Vert x-m_{J}\right\Vert _{L^{2}\left( \mathbf{1}_{J}\omega
\right) }^{2} \\
&\leq &
\!\!\!\!\!\sum_{F\in \mathcal{F}}\sum_{J\in \mathcal{M}^{deep}_{\left( \mathbf{\rho }%
,\varepsilon \right),\mathcal{G}}\left( F\right)
}\sum_{t=0}^{\infty }2^{-t\delta ^{\prime }\left( 1-\varepsilon \right)
}
\left( \frac{\mathrm{P}^{\alpha }\left( J,\mathbf{1}_{\pi _{\mathcal{F}%
}^{t+1}F\backslash \pi _{\mathcal{F}}^{t}F}|\psi _{F}| \sigma \right) }{\left\vert
J\right\vert^\frac{1}{n} }\right) ^{2}\left\Vert x-m_{J}\right\Vert _{L^{2}\left( 
\mathbf{1}_{J}\omega \right) }^{2} \\
&=&
\!\!\!\!\!\sum_{t=0}^{\infty }\!2^{-t\delta ^{\prime }\left( 1-\varepsilon \right)
}\!\!\sum_{G\in \mathcal{F}}\!\sum_{F\in \mathfrak{C}_{\mathcal{F}}^{\left(
t+1\right) }\left( G\right) }\!\sum_{J\in \mathcal{M}^{deep}_{\left( \mathbf{\rho },\varepsilon \right),\mathcal{G}}\left( F\right) }
\!\!\!\!\left(\!\! 
\frac{\mathrm{P}^{\alpha }\left( J,\mathbf{1}_{G\backslash \pi _{\mathcal{F}%
}^{t}F}|\psi _{F}| \sigma \right) }{\left\vert J\right\vert ^\frac{1}{n}}\!\!\right) ^{\!\!\!2}\!\!\left\Vert
x-m_{J}\right\Vert_{L^{\!2}\left( \mathbf{1}_{J}\omega\! \right) }^{2} \\
&\lesssim &
\!\!\!\!\sum_{t=0}^{\infty }\!2^{-t\delta ^{\prime }\left( 1-\varepsilon
\right) }\!\!\sum_{G\in \mathcal{F}}\!\!\!\alpha _{\mathcal{F}}\left( G\right)
^{2}\!\!\!\!\!\!\!\!\!\!\sum_{F\in \mathfrak{C}_{\mathcal{F}}^{\left( t+1\right) }\left(
G\right) }\sum_{J\in \mathcal{M}^{deep}_{\left( \mathbf{\rho },\varepsilon \right)}\left( F\right) }
\!\!\!\!\left(\!\! \frac{\mathrm{P}^{\alpha }\left( J,%
\mathbf{1}_{G\backslash \pi _{\mathcal{F}}^{t}F}\sigma \right) }{\left\vert
J\right\vert^\frac{1}{n} }\!\!\right) ^{\!\!\!2}\!\!\left\Vert x-m_{J}\right\Vert _{L^{2}\left( 
\mathbf{1}_{J}\omega\! \right) }^{2} \\
&\lesssim &
\!\!\!\!\sum_{t=0}^{\infty }2^{-t\delta ^{\prime }\left( 1-\varepsilon
\right) }\sum_{G\in \mathcal{F}}\alpha _{\mathcal{F}}\left( G\right)
^{2}\left( \mathcal{E}_{2}^{\alpha }\right) ^{2}\left\vert G\right\vert
_{\sigma }\lesssim \left( \mathcal{E}_{2}^{\alpha }\right) ^{2}\left\Vert
f\right\Vert _{L^{2}\left( \sigma \right) }^{2}.
\end{eqnarray*}

This completes the proof of the Intertwining Proposition \ref{strongly
adapted}.
\end{proof}

The task of controlling functional energy is taken up in Appendix  below.

\subsection{Paraproduct, neighbour and broken forms}

In this subsection we reduce boundedness of the local below form $\mathsf{B}%
_{\Subset _{\mathbf{r},\varepsilon }}^{A}\left( f,g\right) $ defined in (\ref{def local}) to boundedness of the associated stopping form%
\begin{equation}
\mathsf{B}_{{stop}}^{A}\left( f,g\right) \equiv \sum_{\substack{ %
I\in \mathcal{C}_{A}^{\mathcal{D}}\text{ and }J\in \mathcal{C}_{A}^{\mathcal{%
G},{shift}}  \\ J^{\maltese }\subsetneqq I\text{ and }\ell \left(
J\right) \leq 2^{-\mathbf{r}}\ell \left( I\right) }}\left( E_{I_{J}}^{\sigma
}\widehat{\square }_{I}^{\sigma ,\flat ,\mathbf{b}}f\right) \left\langle
T_{\sigma }^{\alpha }\left( \mathbf{1}_{A\backslash I_{J}}b_{A}\right)
,\square _{J}^{\omega ,\mathbf{b}^{\ast }}g\right\rangle _{\omega }\ ,
\label{def stop}
\end{equation}%
where the modified difference $\widehat{\square }_{I}^{\sigma ,\flat ,%
\mathbf{b}}$ must be carefully chosen in order to control the corresponding paraproduct form
below. Indeed, below we will decompose 
\begin{equation*}
\mathsf{B}_{\Subset _{\mathbf{r},\varepsilon }}^{A}\left( f,g\right) =%
\mathsf{B}_{{paraproduct}}^{A}\left( f,g\right) -\mathsf{B}_{%
{stop}}^{A}\left( f,g\right) +\mathsf{B}_{{neighbour}%
}^{A}\left( f,g\right) +\mathsf{B}_{{brok}}^{A}\left( f,g\right) ,
\end{equation*}%
and we will show that 
$$
\sum_{A\in \mathcal{A}}\left\vert \mathsf{B}_{\Subset _{\mathbf{r}%
,\varepsilon }}^{A}\left( f,g\right) +\mathsf{B}_{{stop}}^{A}\left(
f,g\right) \right\vert\lesssim \left( \mathfrak{T}_{T^{\alpha }}^{\mathbf{b}}+\sqrt{\mathfrak{A}%
_{2}^{\alpha }}\right) \left\Vert f\right\Vert _{L^{2}\left( \sigma \right)
}\left\Vert g\right\Vert _{L^{2}\left( \omega \right) }
$$
and the bound of $\mathsf{B}_{{stop}}^{A}\left(
f,g\right)$ will be the main subject of the next section.

Note that the modified dual martingale differences $\square _{I}^{\sigma ,\flat ,%
\mathbf{b}}$ and $\widehat{\square }_{I}^{\sigma ,\flat ,\mathbf{b}}$, 
\begin{equation*}
\square _{I}^{\sigma ,\flat ,\mathbf{b}}f\equiv \square _{I}^{\sigma ,%
\mathbf{b}}f-\sum_{I^{\prime }\in \mathfrak{C}_{{brok}}\left(
I\right) }\mathbb{F}_{I^{\prime }}^{\sigma ,\mathbf{b}}f=b_{A}\sum_{I^{%
\prime }\in \mathfrak{C}\left( I\right) }\mathbf{1}_{I^{\prime
}}E_{I^{\prime }}^{\sigma }\left( \widehat{\square }_{I}^{\sigma ,\flat ,%
\mathbf{b}}f\right) =b_{A}\widehat{\square }_{I}^{\sigma ,\flat ,\mathbf{b}%
}f,
\end{equation*}%
satisfy the following telescoping property for all $K\in \big( \mathcal{C}_{A}\backslash \left\{ A\right\} \big) \cup \left( \bigcup\limits_{A^{\prime }\in 
\mathfrak{C}_{\mathcal{A}}\left( A\right) }A^{\prime }\right) $ and $L\in 
\mathcal{C}_{A}$ with $K\subset L$:%
\begin{equation*}
\sum_{I:\ \pi K\subset I\subset L}E_{I}^{\sigma }\left( \widehat{\square 
}_{I}^{\sigma ,\flat ,\mathbf{b}}f\right) =\left\{ 
\begin{array}{ccc}
-E_{L}^{\sigma }\widehat{\mathbb{F}}_{L}^{\sigma ,\mathbf{b}}f & \text{ if }
& K\in \mathfrak{C}_{\mathcal{A}}\left( A\right) \\ 
E_{K}^{\sigma }\widehat{\mathbb{F}}_{K}^{\sigma ,\mathbf{b}}f-E_{L}^{\sigma }%
\widehat{\mathbb{F}}_{L}^{\sigma ,\mathbf{b}}f & \text{ if } & K\in \mathcal{%
C}_{A}%
\end{array}%
\right. .
\end{equation*}%
Fix $I\in \mathcal{C}_{A}$ for the moment. We will use%
\begin{eqnarray*}
\mathbf{1}_{I} &=&\mathbf{1}_{I_{J}}+\sum_{\tilde{I} \in \theta \left( I_{J}\right)}\mathbf{1}_{\tilde{I}
}\ , \\
\mathbf{1}_{I_{J}} &=&\mathbf{1}_{A}-\mathbf{1}_{A\backslash I_{J}}\ ,
\end{eqnarray*}%
where $\theta \left( I_{J}\right)$ denotes the $2^n-1\ \mathcal{D}$-children of $I$
other than the child $I_{J}$ that contains $J$. We begin with the splitting%
\begin{eqnarray*}
&&
\left\langle T_{\sigma }^{\alpha }\square _{I}^{\sigma ,\mathbf{b}%
}f,\square _{J}^{\omega ,\mathbf{b}^{\ast }}g\right\rangle _{\omega
}\\
&=&
\left\langle T_{\sigma }^{\alpha }\left( \mathbf{1}_{I_{J}}\square
_{I}^{\sigma ,\mathbf{b}}f\right) ,\square _{J}^{\omega ,\mathbf{b}^{\ast
}}g\right\rangle _{\omega }+
\sum_{\tilde{I}\in \theta \left( I_{J}\right)}\!\!\!\! \left\langle T_{\sigma }^{\alpha }\left( \mathbf{%
1}_{\tilde{I}}\square _{I}^{\sigma ,\mathbf{b}}f\right)
,\square _{J}^{\omega ,\mathbf{b}^{\ast }}g\right\rangle _{\omega } \\
&=&
\left\langle T_{\sigma }^{\alpha }\left( \mathbf{1}_{I_{J}}\square
_{I}^{\sigma ,\flat ,\mathbf{b}}f\right) ,\square _{J}^{\omega ,\mathbf{b}%
^{\ast }}g\right\rangle _{\omega }+ \left\langle T_{\sigma }^{\alpha }\left( 
\mathbf{1}_{I_{J}}\sum_{I^{\prime }\in \mathfrak{C}_{{brok}}\left(
I\right) }\mathbb{F}_{I^{\prime }}^{\sigma ,\mathbf{b}}f\right) ,\square
_{J}^{\omega ,\mathbf{b}^{\ast }}g\right\rangle _{\omega }\\
&& 
\hspace{4.59cm}+\!\!\!\sum_{\tilde{I}\in \theta \left( I_{J}\right)}\!\!\!\! \left\langle T_{\sigma }^{\alpha }\left( \mathbf{%
1}_{\tilde{I}}\square _{I}^{\sigma ,\mathbf{b}}f\right)
,\square _{J}^{\omega ,\mathbf{b}^{\ast }}g\right\rangle _{\omega } \\
&\equiv &\mathrm{I+II+III}\ .
\end{eqnarray*}%
From (\ref{factor b_A}) we have%
\begin{eqnarray*}
\mathrm{I} &=&
\left\langle T_{\sigma }^{\alpha }\left( \mathbf{1}_{I_{J}}\square
_{I}^{\sigma ,\flat ,\mathbf{b}}f\right) ,\square _{J}^{\omega ,\mathbf{b}%
^{\ast }}g\right\rangle _{\omega }=\left\langle T_{\sigma }^{\alpha }\left[
b_{A}\left( \mathbf{1}_{I_{J}}\widehat{\square }_{I}^{\sigma ,\flat ,\mathbf{%
b}}f\right) \right] ,\square _{J}^{\omega ,\mathbf{b}^{\ast }}g\right\rangle
_{\omega } \\
&=&
E_{I_{J}}^{\sigma }\left( \widehat{\square }_{I}^{\sigma ,\flat ,\mathbf{b%
}}f\right) \left\langle T_{\sigma }^{\alpha }\left( \mathbf{1}%
_{I_{J}}b_{A}\right) ,\square _{J}^{\omega ,\mathbf{b}^{\ast
}}g\right\rangle _{\omega } \\
&=&
E_{I_{J}}^{\sigma }\left( \widehat{\square }_{I}^{\sigma ,\flat ,\mathbf{b%
}}f\right) \left\langle T_{\sigma }^{\alpha }b_{A},\square _{J}^{\omega ,%
\mathbf{b}^{\ast }}g\right\rangle _{\omega }\!\!\!-E_{I_{J}}^{\sigma }\left( 
\widehat{\square }_{I}^{\sigma ,\flat ,\mathbf{b}}f\right) \left\langle
T_{\sigma }^{\alpha }\left( \mathbf{1}_{A\backslash I_{J}}b_{A}\right)
,\square _{J}^{\omega ,\mathbf{b}^{\ast }}g\right\rangle _{\omega }
\end{eqnarray*}%
Since the function $\mathbb{F}_{I_{J}}^{\sigma ,\mathbf{b}}f$ is a constant
multiple of $b_{I_{J}}$ on $I_{J}$, we can define $\widehat{\mathbb{F}}%
_{I_{J}}^{\sigma ,\mathbf{b}}f\equiv \frac{1}{b_{I_{J}}}\mathbb{F}%
_{I_{J}}^{\sigma ,\mathbf{b}}f$
and then%
\begin{eqnarray*}
\ \mathrm{II}
=\!
\left\langle T_{\sigma }^{\alpha }\!\!\left(\!\! \mathbf{1}_{I_{J}}\!\!\!\!\sum_{I^{\prime }\in \mathfrak{C}_{{brok}}\left( I\right) }\!\!\!\!\!\!\mathbb{F}%
_{I^{\prime }}^{\sigma ,\mathbf{b}}f\!\!\right)\!\!,\square _{J}^{\omega ,\mathbf{b}%
^{\ast }}g\right\rangle _{\omega }\!\!\!\!
=
\mathbf{1}_{\mathfrak{C}_{\mathcal{A}%
}\left( A\right) }\!\!\left( I_{J}\right) E_{I_{J}}^{\sigma }\left( \widehat{%
\mathbb{F}}_{I_{J}}^{\sigma ,\mathbf{b}}f\right) \left\langle T_{\sigma
}^{\alpha }b_{I_{J}},\square _{J}^{\omega ,\mathbf{b}^{\ast }}g\right\rangle
_{\!\omega }
\end{eqnarray*}%
where the presence of the indicator function $\mathbf{1}_{\mathfrak{C}_{%
\mathcal{A}}\left( A\right) }\left( I_{J}\right) $ simply means that term $\mathrm{II}$ vanishes unless $I_{J}$ is an $\mathcal{A}$-child of $A$. We now write
these terms as
\begin{eqnarray*}
\left\langle T_{\sigma }^{\alpha }\square _{I}^{\sigma ,\mathbf{b}}f,\square
_{J}^{\omega ,\mathbf{b}^{\ast }}g\right\rangle _{\omega }
&=&
E_{I_{J}}^{\sigma }\left( \widehat{\square }_{I}^{\sigma ,\flat ,\mathbf{b%
}}f\right) \left\langle T_{\sigma }^{\alpha }b_{A},\square _{J}^{\omega ,%
\mathbf{b}^{\ast }}g\right\rangle _{\omega } \\
&&-E_{I_{J}}^{\sigma }\left( \widehat{\square }_{I}^{\sigma ,\flat ,\mathbf{b%
}}f\right) \left\langle T_{\sigma }^{\alpha }\left( \mathbf{1}_{A\backslash
I_{J}}b_{A}\right) ,\square _{J}^{\omega ,\mathbf{b}^{\ast }}g\right\rangle
_{\omega } \\
&&
+\sum_{\tilde{I} \in \theta \left( I_{J}\right)}\!\!\!\!\left\langle T_{\sigma }^{\alpha }\left( \mathbf{1}_{\tilde{I}}\square _{I}^{\sigma ,\mathbf{b}}f\right) ,\square
_{J}^{\omega ,\mathbf{b}^{\ast }}g\right\rangle _{\omega } \\
&&
+\mathbf{1}_{\left\{ I_{J}\in \mathfrak{C}_{\mathcal{A}}\left( A\right)
\right\} }\ E_{I_{J}}^{\sigma }\left( \widehat{\mathbb{F}}_{I_{J}}^{\sigma ,%
\mathbf{b}}f\right) \ \left\langle T_{\sigma }^{\alpha }b_{I_{J}},\square
_{J}^{\omega ,\mathbf{b}^{\ast }}g\right\rangle _{\omega }\ ,
\end{eqnarray*}%
where the four lines are respectively a paraproduct, stopping, neighbour and
broken term.

The corresponding NTV splitting of $\mathsf{B}_{\Subset _{\mathbf{r}%
,\varepsilon }}^{A}\left( f,g\right) $ using (\ref{def local}) and (\ref{def shorthand}) becomes%
\begin{eqnarray*}
\mathsf{B}_{\Subset _{\mathbf{r},\varepsilon }}^{A}\left( f,g\right)
&=&
\left\langle T_{\sigma }^{\alpha }\left( \mathsf{P}_{\mathcal{C}%
_{A}}^{\sigma }f\right) ,\mathsf{P}_{\mathcal{C}_{A}^{\mathcal{G},{%
shift}}}^{\omega }g\right\rangle _{\omega }^{\Subset _{\mathbf{r}%
,\varepsilon }} \\
&=&
\sum_{\substack{ I\in \mathcal{C}_{A}\text{ and }J\in 
\mathcal{C}_{A}^{\mathcal{G},{shift}}  \\ J^{\maltese }\subsetneqq I
\text{ and }\ell \left( J\right) \leq 2^{-\mathbf{r}}\ell \left( I\right) }}\!\!\!\!
\left\langle T_{\sigma }^{\alpha }\left( \square _{I}^{\sigma ,\mathbf{b}
}f\right) ,\square _{J}^{\omega ,\mathbf{b}^{\ast }}g\right\rangle _{\omega }
\\
&=&\mathsf{B}_{{paraproduct}}^{A}\left( f,g\right) -\mathsf{B}_{%
{stop}}^{A}\left( f,g\right) +\mathsf{B}_{{neighbour}%
}^{A}\left( f,g\right) +\mathsf{B}_{{brok}}^{A}\left( f,g\right) ,
\end{eqnarray*}%
where

\begin{eqnarray*}
\mathsf{B}_{{paraproduct}}^{A}\left( f,g\right) 
&\equiv &
\sum 
_{\substack{ I\in \mathcal{C}_{A}\text{ and }J\in \mathcal{C}_{A}^{\mathcal{G%
},{shift}}  \\ J^{\maltese }\subsetneqq I\text{ and }\ell \left(
J\right) \leq 2^{-\mathbf{r}}\ell \left( I\right) }}E_{I_{J}}^{\sigma
}\left( \widehat{\square }_{I}^{\sigma ,\flat ,\mathbf{b}}f\right)
\left\langle T_{\sigma }^{\alpha }b_{A},\square _{J}^{\omega ,\mathbf{b}%
^{\ast }}g\right\rangle _{\omega } \\
\mathsf{B}_{{stop}}^{A}\left( f,g\right) 
&\equiv &
\sum_{\substack{ %
I\in \mathcal{C}_{A}\text{ and }J\in \mathcal{C}_{A}^{\mathcal{G},{%
shift}}  \\ J^{\maltese }\subsetneqq I\text{ and }\ell \left( J\right) \leq
2^{-\mathbf{r}}\ell \left( I\right) }}\!\!\!\!\!\!\!\! E_{I_{J}}^{\sigma }\left( \widehat{%
\square }_{I}^{\sigma ,\flat ,\mathbf{b}}f\right) \left\langle T_{\sigma
}^{\alpha }\left( \mathbf{1}_{A\backslash I_{J}}b_{A}\right) ,\square
_{J}^{\omega ,\mathbf{b}^{\ast }}g\right\rangle _{\omega } \\
\mathsf{B}_{{neighbour}}^{A}\left( f,g\right) 
&\equiv &
\!\!\!\!\!\!\!\!\sum 
_{\substack{ I\in \mathcal{C}_{A}\text{ and }J\in \mathcal{C}_{A}^{\mathcal{G},{shift}}  \\ J^{\maltese }\subsetneqq I\text{ and }\ell \left(
J\right) \leq 2^{-\mathbf{r}}\ell \left( I\right) }}\!\!
\sum_{\tilde{I}\in \theta \left( I_{J}\right)}\!\!\!\!
\left\langle T_{\sigma
}^{\alpha }\left( \mathbf{1}_{\tilde{I} }\square
_{I}^{\sigma ,\mathbf{b}}f\right) ,\square _{J}^{\omega ,\mathbf{b}^{\ast
}}g\right\rangle _{\!\omega }
\end{eqnarray*}%
correspond to the three original NTV forms associated with $1$-testing, and
where 
\begin{equation}
\mathsf{B}_{{brok}}^{A}\left( f,g\right) \equiv \sum_{\substack{ %
I\in \mathcal{C}_{A}\text{ and }J\in \mathcal{C}_{A}^{\mathcal{G},{%
shift}}  \\ J^{\maltese }\subsetneqq I\text{ and }\ell \left( J\right) \leq
2^{-\mathbf{r}}\ell \left( I\right) }}\!\!\!\!\!\!\mathbf{1}_{\left\{ I_{J}\in \mathfrak{%
C}_{\mathcal{A}}\left( A\right) \right\} }\ E_{I_{J}}^{\sigma }\left( 
\widehat{\mathbb{F}}_{I_{J}}^{\sigma ,\mathbf{b}}f\right) \ \left\langle
T_{\sigma }^{\alpha }b_{I_{J}},\square _{J}^{\omega ,\mathbf{b}^{\ast
}}g\right\rangle _{\omega }  \label{broken vanish}
\end{equation}%
"vanishes" since $J^{\maltese }\subsetneqq I$ and $I_{J}\in \mathfrak{C}_{\mathcal{A}%
}\left( A\right) $ imply that $J^{\maltese }\notin \mathcal{C}_{A}^{\mathcal{G}}$, contradicting $J\in \mathcal{C}_{A}^{\mathcal{G},{shift}}$.

\begin{rem}
The inquisitive reader will note that the pairs $\left( I,J\right) $ arising
in the above sum with $J^{\maltese }\subsetneqq I$ replaced by $J^{\maltese
}=I$ are handled in the probabilistic estimate (\ref{HM bad}) for the bad
form $\Theta _{2}^{{bad}\natural }$ defined in (\ref{Theta_2^bad
sharp}).
\end{rem}

\subsubsection{The paraproduct form}

The paraproduct form $\mathsf{B}_{{paraproduct}}^{A}\left(
f,g\right) $ is easily controlled by the testing condition for $T^{\alpha }$
together with weak Riesz inequalities for dual martingale differences.
Indeed, recalling the telescoping identity (\ref{telescoping}), and that the
collection $\left\{ I\in \mathcal{C}_{A}\text{:\ }\ell \left( J\right) \leq
2^{-\mathbf{r}}\ell \left( I\right) \right\} $ is tree connected for all $%
J\in \mathcal{C}_{A}^{\mathcal{G},{shift}}$, we have%
\begin{eqnarray*}
\mathsf{B}_{{paraproduct}}^{A}\left( f,g\right) &=&
\sum_{\substack{ 
I\in \mathcal{C}_{A}\text{ and }J\in \mathcal{C}_{A}^{\mathcal{G},{%
shift}}  \\ J^{\maltese }\subsetneqq I\text{ and }\ell \left( J\right) \leq
2^{-\mathbf{r}}\ell \left( I\right) }}E_{I_{J}}^{\sigma }\left( \widehat{%
\square }_{I}^{\sigma ,\flat ,\mathbf{b}}f\right) \left\langle T_{\sigma
}^{\alpha }b_{A},\square _{J}^{\omega ,\mathbf{b}^{\ast }}g\right\rangle
_{\omega } \\
&=&
\sum_{J\in \mathcal{C}_{A}^{\mathcal{G},{shift}}}\left\langle
T_{\sigma }^{\alpha }b_{A},\square _{J}^{\omega ,\mathbf{b}^{\ast
}}g\right\rangle _{\omega }\left\{ \sum_{I\in \mathcal{C}_{A}\text{:\ }%
J^{\maltese }\subsetneqq I\text{ and }\ell \left( J\right) \leq 2^{-\mathbf{r%
}}\ell \left( I\right) }E_{I_{J}}^{\sigma }\left( \widehat{\square }%
_{I}^{\sigma ,\flat ,\mathbf{b}}f\right) \right\}
\end{eqnarray*}
\begin{eqnarray*}
&=&
\sum_{J\in \mathcal{C}_{A}^{\mathcal{G},{shift}}}\left\langle
T_{\sigma }^{\alpha }b_{A},\square _{J}^{\omega ,\mathbf{b}^{\ast
}}g\right\rangle _{\omega } \left\{ \mathbf{1}_{\left\{ J:I^{\natural }\left(
J\right) _{J}\in \mathcal{C}_{A}\right\} }E_{I^{\natural }\left( J\right)
_{J}}^{\sigma }\widehat{\mathbb{F}}_{I^{\natural }\left( J\right)
_{J}}^{\sigma ,\mathbf{b}}f-E_{A}^{\sigma }\widehat{\mathbb{F}}_{A}^{\sigma ,%
\mathbf{b}}f\right\} \\
&=&
\left\langle T_{\sigma }^{\alpha }b_{A},\sum_{J\in \mathcal{C}_{A}^{%
\mathcal{G},{shift}}}\left\{ \mathbf{1}_{\left\{ J:I^{\natural
}\left( J\right) _{J}\in \mathcal{C}_{A}\right\} }E_{I^{\natural }\left(
J\right) _{J}}^{\sigma }\widehat{\mathbb{F}}_{I^{\natural }\left( J\right)
_{J}}^{\sigma ,\mathbf{b}}f-E_{A}^{\sigma }\widehat{\mathbb{F}}_{A}^{\sigma ,%
\mathbf{b}}f\right\} \square _{J}^{\omega ,\mathbf{b}^{\ast }}g\right\rangle
_{\omega }
\end{eqnarray*}
where $I^{\natural }\left( J\right) $ denotes the smallest cube $I\in 
\mathcal{C}_{A}$ such that $J^{\maltese }\subsetneqq I$ and $\ell \left(
J\right) \leq 2^{-\mathbf{r}}\ell \left( I\right) $, and of course $%
I^{\natural }\left( J\right) _{J}$ denotes its child containing $J$. Note
that by construction of the modified difference operator $\square
_{I}^{\sigma ,\flat ,\mathbf{b}}$, the only time the average $\widehat{%
\mathbb{F}}_{I^{\natural }\left( J\right) _{J}}^{\sigma }f$ appears in the
above sum is when $I^{\natural }\left( J\right) _{J}\in \mathcal{C}_{A}$,
since the case $I^{\natural }\left( J\right) _{J}\in \mathcal{A}$ has been
removed to the broken term. This is reflected above with the inclusion of
the indicator $\mathbf{1}_{\left\{ J:I^{\natural }\left( J\right) _{J}\in 
\mathcal{C}_{A}\right\} }$. It follows that we have the bound 
$$
\left\vert 
\mathbf{1}_{\left\{ J:I^{\natural }\left( J\right) _{J}\in \mathcal{C}%
_{A}\right\} }E_{I^{\natural }\left( J\right) _{J}}^{\sigma }\widehat{%
\mathbb{F}}_{I^{\natural }\left( J\right) _{J}}^{\sigma ,\mathbf{b}%
}f\right\vert +\left\vert E_{A}^{\sigma }\widehat{\mathbb{F}}_{A}^{\sigma ,%
\mathbf{b}}f\right\vert \lesssim E_{A}^{\sigma }\left\vert f\right\vert \leq
\alpha _{\mathcal{A}}\left( A\right)
$$

Thus from Cauchy-Schwarz, the upper weak Riesz inequalities for the pseudoprojections $\square _{J}^{\omega ,\mathbf{b}%
^{\ast }}g$ and the bound on the coefficients 
$$
\lambda _{J}\equiv \left( 
\mathbf{1}_{\left\{ J:I^{\natural }\left( J\right) _{J}\in \mathcal{C}%
_{A}\right\} }E_{I^{\natural }\left( J\right) _{J}}^{\sigma }\widehat{%
\mathbb{F}}_{I^{\natural }\left( J\right) _{J}}^{\sigma ,\mathbf{b}%
}f-E_{A}^{\sigma }\widehat{\mathbb{F}}_{A}^{\sigma ,\mathbf{b}}f\right) 
$$
given by $\left\vert \lambda _{J}\right\vert \lesssim \alpha _{\mathcal{A}%
}\left( A\right) $, we have
\begin{eqnarray}
&&  \label{est para} \\
&&\hspace{-1cm}
\left\vert \mathsf{B}_{{paraproduct}}^{A}\left( f,g\right)
\right\vert
=
\left\vert \left\langle T_{\sigma }^{\alpha }b_{A},\!\!\!\!\!\!\sum_{J\in 
\mathcal{C}_{A}^{\mathcal{G},{shift}}}\left\{ \left( \mathbf{1}%
_{\left\{ J:I^{\natural }\left( J\right) _{J}\in \mathcal{C}_{A}\right\}
}E_{I^{\natural }\left( J\right) _{J}}^{\sigma }\widehat{\mathbb{F}}
_{I^{\natural }\left( J\right) _{J}}^{\sigma ,\mathbf{b}}f-E_{A}^{\sigma }%
\widehat{\mathbb{F}}_{A}^{\sigma ,\mathbf{b}}f\right) \right\} \square
_{J}^{\omega ,\mathbf{b}^{\ast }}g\right\rangle _{\omega }\right\vert  \notag
\\
&\leq &
\left\Vert \mathbf{1}_{A}T_{\sigma }^{\alpha }b_{A}\right\Vert
_{L^{2}\left( \omega \right) }\left\Vert \sum_{J\in \mathcal{C}_{A}^{%
\mathcal{G},{shift}}}\lambda _{J}\square _{J}^{\omega ,\mathbf{b}%
^{\ast }}g\right\Vert _{L^{2}\left( \omega \right) }  \notag \\
&\lesssim &
\alpha _{\mathcal{A}}\left( A\right) \ \left\Vert \mathbf{1}%
_{A}T_{\sigma }^{\alpha }b_{A}\right\Vert _{L^{2}\left( \omega \right) }\
\sum_{J\in \mathcal{C}_{A}^{\mathcal{G},{shift}}}\left\Vert \square
_{J}^{\omega ,\mathbf{b}^{\ast }}g\right\Vert _{L^{2}\left( \omega \right)
}  \notag \\
&\leq &
\mathfrak{T}_{T^{\alpha }}^{\mathbf{b}}\ \alpha _{\mathcal{A}}\left(
A\right) \ \sqrt{\left\vert A\right\vert _{\sigma }}\ \left\Vert \mathsf{P}_{%
\mathcal{C}_{A}^{\mathcal{G},{shift}}}^{\omega ,\mathbf{b}^{\ast
}}g\right\Vert _{L^{2}\left( \omega \right) }^{\bigstar }.  \notag
\end{eqnarray}

\subsubsection{The neighbour form}

Next, the neighbour form $\mathsf{B}_{{neighbour}}^{A}\left(
f,g\right) $ is easily controlled by the $\mathfrak{A}_{2}^{\alpha }$
condition using the pivotal estimate in Energy Lemma \ref{ener} and the fact
that the cubes $J\in \mathcal{C}_{A}^{\mathcal{G},{shift}}$ are
good in $I$ and beyond when the pair $\left( I,J\right) $ occurs in the sum.
In particular, the information encoded in the stopping tree $\mathcal{A}$
plays no role here, apart from appearing in the corona projections on the
right hand side of (\ref{est neigh}) below. We have%
\begin{equation}
\mathsf{B}_{{neighbour}}^{A}\left( f,g\right) 
= 
\!\!\!\!\!\!\!\!\sum 
_{\substack{ I\in \mathcal{C}_{A}\text{ and }J\in \mathcal{C}_{A}^{\mathcal{G},{shift}}  \\ J^{\maltese }\subsetneqq I\text{ and }\ell \left(
J\right) \leq 2^{-\mathbf{r}}\ell \left( I\right) }}\!\!
\sum_{\tilde{I}\in\theta \left( I_{J}\right)} \!\!\!
\left\langle T_{\sigma
}^{\alpha }\left( \mathbf{1}_{\tilde{I} }\square
_{I}^{\sigma ,\mathbf{b}}f\right) ,\square _{J}^{\omega ,\mathbf{b}^{\ast
}}g\right\rangle _{\!\omega }
 \label{def neighbour}
\end{equation}%
where we keep in mind that the pairs $\left( I,J\right) \in \mathcal{D}%
\times \mathcal{G}$ that arise in the sum for $\mathsf{B}_{{neighbour%
}}^{A}\left( f,g\right) $ satisfy the property that $J^{\maltese
}\subsetneqq I$, so that $J$ is good with respect to all cubes $K$ of
size at least that of $J^{\maltese }$, which includes $I$. Recall that $%
I_{J} $ is the child of $I$ that contains $J$, and that $\theta \left(
I_{J}\right) $ denotes its $2^n-1$ siblings in $I$, i.e. $\theta \left( I_{J}\right)
= \mathfrak{C}_{\mathcal{D}}\left( I\right) \backslash \left\{
I_{J}\right\} $. Fix $\left( I,J\right) $ momentarily, and an integer $s\geq 
\mathbf{r}$. Using $\square _{I}^{\sigma ,\mathbf{b}}=\square _{I}^{\sigma
,\flat ,\mathbf{b}}+\square _{I,{brok}}^{\sigma ,\flat ,\mathbf{b}%
} $ and the fact that $\square _{I}^{\sigma ,\flat ,\mathbf{b}}f$ is a
constant multiple of $b_{\tilde{I}}$ on the cube $\tilde{I} $, we have the estimates 
\begin{eqnarray*}
\left\vert \mathbf{1}_{\tilde{I}}\square _{I}^{\sigma
,\flat ,\mathbf{b}}f\right\vert &=&\left\vert \left( E_{\tilde{I}}^{\sigma }\widehat{\square }_{I}^{\sigma ,\flat ,\mathbf{b}%
}f\right) b_{\tilde{I} }\right\vert \leq C_{\mathbf{b}
}\left\vert E_{\tilde{I} }^{\sigma }\widehat{\square }
_{I}^{\sigma ,\flat ,\mathbf{b}}f\right\vert , \\
\left\vert \mathbf{1}_{\tilde{I}}\square _{I,{
brok}}^{\sigma ,\flat ,\mathbf{b}}f\right\vert 
&\leq &
\mathbf{1}_{\mathfrak{C}_{A}\left( A\right) }( \tilde{I} ) \
E_{\tilde{I} }^{\sigma }\left\vert f\right\vert ,
\end{eqnarray*}%
and hence%
\begin{equation}
\mathbf{1}_{\tilde{I} }\left\vert \square _{I}^{\sigma ,
\mathbf{b}}f\right\vert \leq C\mathbf{1}_{\tilde{I}}\left(
\left\vert E_{\tilde{I} }^{\sigma }\widehat{\square }
_{I}^{\sigma ,\flat ,\mathbf{b}}f\right\vert +\mathbf{1}_{\mathfrak{C}
_{A}\left( A\right) }( \tilde{I} ) \ E_{\tilde{I} }^{\sigma }\left\vert f\right\vert \right) ,
\label{box bound}
\end{equation}%
which will be used below after an application of the Energy Lemma. We can write $\mathsf{B}_{{neighbour}}^{A}\left( f,g\right)$ as
\begin{equation*}
\sum_{\substack{ I\in 
\mathcal{C}_{A}\& J\in \mathcal{G}_{\left( \kappa \left(
I_{J},J\right) ,\varepsilon \right) -{good}}^{\mathcal{D}}\cap 
\mathcal{C}_{A}^{\mathcal{G},{shift}}\& J^{\maltese
}\subsetneqq I  \\ d\left( J,\tilde{I} \right) >2\ell
\left( J\right) ^{\varepsilon }\ell \left( \tilde{I}
\right) ^{1-\varepsilon }\text{ and }\ell \left( J\right) \leq 2^{-\mathbf{r}%
}\ell \left( I\right) }}
\sum_{\tilde{I} \in \theta \left( I_{J}\right)
 }\!\!\!\!\left\langle T_{\sigma }^{\alpha }\left( \mathbf{1}%
_{\tilde{I} }\square _{I}^{\sigma ,\mathbf{b}}f\right)
,\square _{J}^{\omega ,\mathbf{b}^{\ast }}g\right\rangle _{\omega }
\end{equation*}%
where we have included the conditions 
$$
J\in \mathcal{G}_{\left( \kappa
\left( I_{J},J\right) ,\varepsilon \right) -{good}}^{\mathcal{D}}
\text{ and }
d( J,\tilde{I}) >2\ell \left( J\right)
^{\varepsilon }\ell ( \tilde{I} )
^{1-\varepsilon }
$$
in the summation since they are already implied the
remaining four conditions, and will be used in estimates below.

We will also use the following fractional analogue of the Poisson inequality
in \cite{Vol}.

\begin{lem}
\label{Poisson inequality}Suppose $0\leq \alpha <1$ and $J\!\subset\! I\!\subset\! K$
and that $d\left( J,\partial I\right)\!>\! 2\ell \left( J\right) ^{\varepsilon
}\!\ell \left( I\right) ^{1-\varepsilon }$ for some $0<\varepsilon <\frac{1}{%
n+1-\alpha }$. Then for a positive Borel measure $\mu $ we have%
\begin{equation}
\mathrm{P}^{\alpha }(J,\mu \mathbf{1}_{K\backslash I})\lesssim \left( \frac{%
\ell \left( J\right) }{\ell \left( I\right) }\right) ^{1-\varepsilon \left(
n+1-\alpha \right) }\mathrm{P}^{\alpha }(I,\mu \mathbf{1}_{K\backslash I}).
\label{e.Jsimeq}
\end{equation}
\end{lem}

\begin{proof}
We have%
\begin{equation*}
\mathrm{P}^{\alpha }\left( J,\mu \mathbf{1}_{K\backslash I}\right) \approx
\sum_{k=0}^{\infty }2^{-k}\frac{1}{\left\vert 2^{k}J\right\vert ^{1-\frac{\alpha}{n} }}%
\int_{\left( 2^{k}J\right) \cap \left( K\backslash I\right) }d\mu ,
\end{equation*}%
and $\left( 2^{k}J\right) \cap \left( K\backslash I\right) \neq \emptyset $
requires%
\begin{equation*}
d\left( J,K\backslash I\right) \leq c2^{k}\ell \left( J\right) ,
\end{equation*}%
for some dimensional constant $c>0$. Let $k_{0}$ be the smallest such $k$.
By our distance assumption we must then have%
\begin{equation*}
2\ell \left( J\right) ^{\varepsilon }\ell \left( I\right) ^{1-\varepsilon
}\leq d\left( J,\partial I\right) \leq c2^{k_{0}}\ell \left( J\right) ,
\end{equation*}%
or%
\begin{equation*}
2^{-k_{0}+1}\leq c\left( \frac{\ell \left( J\right) }{\ell \left( I\right) }%
\right) ^{1-\varepsilon }.
\end{equation*}%
Now let $k_{1}$ be defined by $2^{k_{1}}\equiv \frac{\ell \left( I\right) }{%
\ell \left( J\right) }$. Then assuming $k_{1}>k_{0}$ (the case $k_{1}\leq
k_{0}$ is similar) we have%

\begin{eqnarray*}
&&\mathrm{P}^{\alpha }\left( J,\mu \mathbf{1}_{K\backslash I}\right) 
\approx
\left\{ \sum_{k=k_{0}}^{k_{1}}+\sum_{k=k_{1}}^{\infty }\right\} 2^{-k}\frac{%
1}{\left\vert 2^{k}J\right\vert ^{1-\frac{\alpha}{n} }}\int_{\left( 2^{k}J\right) \cap
\left( K\backslash I\right) }d\mu \\
&\lesssim &
2^{-k_{0}}\frac{\left\vert I\right\vert ^{1-\frac{\alpha}{n} }}{\left\vert
2^{k_{0}}J\right\vert ^{1-\frac{\alpha}{n} }}\left( \frac{1}{\left\vert I\right\vert
^{1-\frac{\alpha}{n} }}\int_{\left( 2^{k_{1}}J\right) \cap \left( K\backslash I\right)
}d\mu \right) +2^{-k_{1}}\mathrm{P}^{\alpha }\left( I,\mu \mathbf{1}
_{K\backslash I}\right) \\
&\lesssim &
\left( \frac{\ell \left( J\right) }{\ell \left( I\right) }\right)
^{\left( 1-\varepsilon \right) \left(n+1-\alpha \right) }\left( \frac{\ell
\left( I\right) }{\ell \left( J\right) }\right) ^{n-\alpha }\mathrm{P}%
^{\alpha }\left( I,\mu \mathbf{1}_{K\backslash I}\right)
+
\frac{\ell \left(
J\right) }{\ell \left( I\right) }\mathrm{P}^{\alpha }\left( I,\mu \mathbf{1}%
_{K\backslash I}\right) ,
\end{eqnarray*}%
which is the inequality (\ref{e.Jsimeq}).
\end{proof}

Now fix $I_{0}=I_{J},I_{\theta } \in \theta \left( I_{J}\right)$ and assume that $J\Subset _{\mathbf{r}%
,\varepsilon }I_{0}$. Let $\frac{\ell \left( J\right) }{\ell \left(
I_{0}\right) }=2^{-s}$ in the pivotal estimate from Energy Lemma \ref{ener}
with $J\subset I_{0}\subset I$ to obtain 
\begin{eqnarray*}
 \left\vert \langle T_{\sigma }^{\alpha }\left( \mathbf{1}_{I_\theta}\square _{I}^{\sigma ,\mathbf{b}}f\right) ,\square
_{J}^{\omega ,\mathbf{b}^{\ast }}g\rangle _{\omega }\right\vert 
\!\!\!\!& \lesssim&
 \left\Vert \square _{J}^{\omega ,\mathbf{b}^{\ast }}g\right\Vert
_{L^{2}\left( \omega \right) }\sqrt{\left\vert J\right\vert _{\omega }}%
\mathrm{P}^{\alpha }\left( J,\mathbf{1}_{I_\theta}\left\vert \square _{I}^{\sigma ,\mathbf{b}}f\right\vert \sigma \right) \\
& \lesssim&
 \left\Vert \square _{J}^{\omega ,\mathbf{b}^{\ast }}g\right\Vert
_{L^{2}\left( \omega \right) }\sqrt{\left\vert J\right\vert _{\omega }}\cdot
2^{-\left( 1-\varepsilon \left( n+1-\alpha \right) \right) s}\mathrm{P}%
^{\alpha }\left( I_{0},\mathbf{1}_{I_\theta }\left\vert
\square _{I}^{\sigma ,\mathbf{b}}f\right\vert \sigma \right) \\
& \lesssim&
 \left\Vert \square _{J}^{\omega ,\mathbf{b}^{\ast }}g\right\Vert
_{L^{2}\left( \omega \right) }\sqrt{\left\vert J\right\vert _{\omega }}\cdot
2^{-\left( 1-\varepsilon \left( n+1-\alpha \right) \right) s} 
 \mathrm{P}%
^{\alpha }\left( I_{0},\mathbf{1}_{I_\theta}\mathbf{E}_{I_\theta }^{\sigma }f\cdot \sigma
\right) 
\end{eqnarray*}%
Here we are using (\ref{e.Jsimeq}) in the third line, which applies since $%
J\subset I_{0}$, and we have used (\ref{box bound}) in the fourth line and the shorthand notation%
\begin{equation*}
\mathbf{E}_{I_\theta}^{\sigma }f\equiv \left\vert
E_{I_\theta }^{\sigma }\widehat{\square }_{I}^{\sigma
,\flat ,\mathbf{b}}f\right\vert +\mathbf{1}_{\mathfrak{C}_{A}\left( A\right)
}\left( I_\theta \right) \ E_{I_\theta 
}^{\sigma }\left\vert f\right\vert
\end{equation*}%
where the cube $I$ on the right hand side is determined
uniquely by the cube $I_\theta \in \theta \left( I_{J}\right) $.

In the sum below, we keep the side lengths of the cubes $J$ fixed at $%
2^{-s}$ times that of $I_{0}$, and of course take $J\subset I_{0}$. We also
keep the underlying assumptions that $J\in \mathcal{C}_{A}^{\mathcal{G},%
{shift}}$ and that $J\in \mathcal{G}_{\left( \kappa \left(
I_{J},J\right) ,\varepsilon \right) -{good}}^{\mathcal{D}}$ in mind
without necessarily pointing to them in the notation. Matters will shortly
be reduced to estimating the following term: 
\begin{align*}
A(I,I_{0},I_{\theta },s)
& \equiv
 \sum_{J\;:\;2^{s+1}\ell \left( J\right)
=\ell \left( I\right) :J\subset I_{0}}\left\vert \langle T_{\sigma }^{\alpha
}\left( \mathbf{1}_{I_{\theta }}\square _{I}^{\sigma ,\mathbf{b}}f\right)
,\square _{J}^{\omega ,\mathbf{b}^{\ast }}g\rangle _{\omega }\right\vert \\
& \leq
 2^{-\left( 1-\varepsilon \left(n+1-\alpha \right) \right) s}\!\left( 
\mathbf{E}_{I_\theta}^{\sigma }f\right) \mathrm{P}^{\alpha }(I_{0},\mathbf{1}_{I_\theta}\sigma
)\!\!\!\!\!\!\!\!\!\!\sum_{\substack{J:  J\subset
I_{0}\\ 2^{s+1}\ell \left( J\right) =\ell \left( I\right)}}\!\!\!\!\!\!\!\!\! \left\Vert \square _{J}^{\omega ,\mathbf{b}^{\ast }}g\right\Vert
_{L^{2}\left( \omega \right) }\!\!\sqrt{\left\vert J\right\vert _{\omega }} \\
& \leq
 2^{-\left( 1-\varepsilon \left(n+1-\alpha \right) \right) s}\!\left( 
\mathbf{E}_{I_\theta  }^{\sigma }f\right)\mathrm{P}%
^{\alpha }(I_{0},\mathbf{1}_{I_\theta}\sigma )\sqrt{%
\left\vert I_{0}\right\vert _{\omega }}\Lambda (I,I_{0},I_{\theta },s)
\end{align*}
$\displaystyle \text{where }\Lambda (I,I_{0},I_{\theta },s)^{2}\equiv \sum_{J\in \mathcal{%
C}_{A}^{\mathcal{G},{shift}}:\;2^{s+1}\ell \left( J\right) =\ell
\left( I\right) :\ J\subset I_{0}}\left\Vert \square _{J}^{\omega ,\mathbf{b}%
^{\ast }}g\right\Vert _{L^{2}\left( \omega \right) }^{2}$\,.

The last line follows upon using the Cauchy-Schwarz inequality and the fact
that $J\in \mathcal{C}_{A}^{\mathcal{G},{shift}}$. We also note that
since $2^{s+1}\ell \left( J\right) =\ell \left( I\right) $, 
\begin{eqnarray}
\sum_{I_{0}\in \mathfrak{C}_{\mathcal{D}}\left( I\right) }\Lambda
(I,I_{0},I_{\theta },s)^{2} &\equiv &\sum_{J\in \mathcal{C}_{A}^{\mathcal{G},%
{shift}}:\;2^{s+1}\ell \left( J\right) =\ell \left( I\right) :\
J\subset I}\left\Vert \square _{J}^{\omega ,\mathbf{b}^{\ast }}g\right\Vert
_{L^{2}\left( \omega \right) }^{2}\ ;  \label{g} \\
\sum_{I\in \mathcal{C}_{A}}\sum_{I_{0}\in \mathfrak{C}_{\mathcal{D}}\left(
I\right) }\Lambda (I,I_{0},I_{\theta },s)^{2} &\leq &\left\Vert \mathsf{P}_{%
\mathcal{C}_{A}^{\mathcal{G},{shift}}}^{\omega ,\mathbf{b}^{\ast
}}g\right\Vert _{L^{2}(\omega )}^{\bigstar 2}  \notag
\end{eqnarray}
Using (\ref{PLBP removed}) we obtain 
\begin{equation}
\left\vert E_{I_{\theta }}^{\sigma }\left( \widehat{\square }_{I}^{\sigma
,\flat ,\mathbf{b}}f\right) \right\vert \leq \sqrt{E_{I_{\theta }}^{\sigma
}\left\vert \widehat{\square }_{I}^{\sigma ,\flat ,\mathbf{b}}f\right\vert
^{2}}\lesssim \left\Vert \square _{I}^{\sigma ,\mathbf{b}}f\right\Vert
_{L^{2}\left( \sigma \right) }^{\bigstar }\ \left\vert I_{\theta
}\right\vert _{\sigma }^{-\frac{1}{2}}  \label{e.haarAvg}
\end{equation}%
and hence%
\begin{eqnarray*}
\mathbf{E}_{I_\theta }^{\sigma }f
&\equiv&
\left\vert
E_{I_\theta \left( I_{J}\right) }^{\sigma }\widehat{\square }_{I}^{\sigma
,\flat ,\mathbf{b}}f\right\vert +\mathbf{1}_{\mathfrak{C}_{A}\left( A\right)}\left( I_\theta\right) \ E_{I_\theta
}^{\sigma }\left\vert f\right\vert \\
&\lesssim&
 \left( \left\Vert \square
_{I}^{\sigma ,\mathbf{b}}f\right\Vert _{L^{2}\left( \sigma \right)
}^{\bigstar }\ +\mathbf{1}_{\mathfrak{C}_{A}\left( A\right) }\left(I_\theta \right) \ \left\vert I_{\theta }\right\vert _{\sigma }^{%
\frac{1}{2}}E_{I_\theta}^{\sigma }\left\vert f\right\vert
\right) \left\vert I_{\theta }\right\vert _{\sigma }^{-\frac{1}{2}}
\end{eqnarray*}%
and thus $A(I,I_{0},I_{\theta },s)$ is bounded by
\begin{eqnarray*}
&&\hspace{-1cm}
2^{-\left( 1-\varepsilon \left( n+1-\alpha \right) \right) s}\left(
\left\Vert \square _{I}^{\sigma ,\mathbf{b}}f\right\Vert _{L^{2}\left(
\sigma \right) }^{\bigstar }+\mathbf{1}_{\mathfrak{C}_{A}\left( A\right)
}\left( I_{\theta }\right) \ \left\vert I_{\theta }\right\vert _{\sigma }^{\frac{1}{2}}E_{I_{\theta }}^{\sigma }\left\vert f\right\vert \right) \Lambda(I,I_{0},I_{\theta },s) \left\vert I_{\theta }\right\vert _{\sigma }^{-\frac{1}{2}}\mathrm{P}^{\alpha }(I_{0},\mathbf{1}_{I_\theta}\sigma )\sqrt{\left\vert I_{0}\right\vert _{\omega }} \\
&\lesssim &
\!\!\!\!\!\!\sqrt{\mathfrak{A}_{2}^{\alpha }}2^{-\left( 1-\varepsilon \left(
n+1-\alpha \right) \right) s}\left( \left\Vert \square _{I}^{\sigma ,\mathbf{b}%
}f\right\Vert _{L^{2}\left( \sigma \right) }^{\bigstar }+\mathbf{1}_{\mathfrak{C}_{A}\left( A\right) }\left( I_{\theta }\right) \ \left\vert
I_{\theta }\right\vert _{\sigma }^{\frac{1}{2}}E_{I_{\theta }}^{\sigma
}\left\vert f\right\vert \right) \Lambda (I,I_{0},I_{\theta },s)
\end{eqnarray*}%
since $\mathrm{P}^{\alpha }(I_{0},\mathbf{1}_{I_\theta
}\sigma )\lesssim \frac{\displaystyle\left\vert I_{\theta }\right\vert _{\sigma }}{\displaystyle
\left\vert I_{\theta }\right\vert ^{1-\frac{\alpha}{n} }}$ shows that 
\begin{equation*}
\left\vert I_{\theta }\right\vert _{\sigma }^{-\frac{1}{2}}\mathrm{P}%
^{\alpha }(I_{0},\mathbf{1}_{I_\theta}\sigma )\ \sqrt{%
\left\vert I_{0}\right\vert _{\omega }}\lesssim \frac{\sqrt{\left\vert
I_{\theta }\right\vert _{\sigma }}\sqrt{\left\vert I_{0}\right\vert _{\omega
}}}{\left\vert I_\theta\right\vert ^{1-\frac{\alpha}{n} }}\lesssim \sqrt{\mathfrak{A}%
_{2}^{\alpha }}
\end{equation*}
where the implied constant depends on $\alpha$ and the dimension. An application of Cauchy-Schwarz to the sum over $I$ using (\ref{g}) then shows that 

\begin{eqnarray*}
&&\sum_{I\in \mathcal{C}_{A}}\sum_{\substack{I_{0},I_{\theta }\in \mathfrak{%
C}_{\mathcal{D}}\left( I\right)  \\ I_{0}\neq I_{\theta }}}%
A(I,I_{0},I_{\theta },s) \\
&\lesssim &\!\!\!\!\!
\sqrt{\mathfrak{A}_{2}^{\alpha }}2^{-\left( 1-\varepsilon \left(
n+1-\alpha \right) \right) s}\sqrt{\sum_{I\in \mathcal{C}_{A}}\left\Vert
\square _{I}^{\sigma ,\mathbf{b}}f\right\Vert _{L^{2}\left( \sigma \right)
}^{\bigstar 2}\!\!+\!\!\sum_{I_{\theta }\in \mathfrak{C}_{A}\left( A\right)
}\!\!\!\left\vert I_{\theta }\right\vert _{\sigma }\left( E_{I_{\theta }}^{\sigma
}\left\vert f\right\vert \right) ^{2}}
 \sqrt{\!\sum_{I\in \mathcal{C}%
_{A}}\!\!\left( \!\sum_{\substack{ I_{0},I_{\theta }\in \mathfrak{C}_{\mathcal{D}%
}\left( I\right)  \\ I_{0}\neq I_{\theta }}}\!\!\!\!\Lambda (I,I_{0},I_{\theta
},s)\!\!\right)^{\!\!2}} \\
&\lesssim &\!\!\!\!\!
\sqrt{\mathfrak{A}_{2}^{\alpha }}2^{-\left( 1-\varepsilon \left(
n+1-\alpha \right) \right) s}\sqrt{\left\Vert \mathsf{P}_{\mathcal{C}%
_{A}}^{\sigma }f\right\Vert _{L^{2}(\sigma )}^{\bigstar 2}+\sum_{A^{\prime
}\in \mathfrak{C}_{A}\left( A\right) }\left\vert A^{\prime }\right\vert
_{\sigma }\left( E_{A^{\prime }}^{\sigma }\left\vert f\right\vert \right)
^{\!2}}
\sqrt{\sum_{I\in \mathcal{C}_{A}}\left( \sum_{\substack{ I_{0}\in 
\mathfrak{C}_{\mathcal{D}}\left( I\right)  \\ I_{0}\neq I_{\theta }}}\Lambda
(I,I_{0},I_{\theta },s)\right) ^{2}} \\
&\lesssim &\!\!\!\!\!
\sqrt{\mathfrak{A}_{2}^{\alpha }}2_{{}}^{-\left( 1-\varepsilon
\left( n+1-\alpha \right) \right) s} 
\left( \lVert \mathsf{P}_{\mathcal{C}%
_{A}}^{\sigma }f\rVert _{L^{2}(\sigma )}^{\bigstar }+\sqrt{\sum_{A^{\prime
}\in \mathfrak{C}_{A}\left( A\right) }\left\vert A^{\prime }\right\vert
_{\sigma }\left( E_{A^{\prime }}^{\sigma }\left\vert f\right\vert \right)
^{2}}\right) \left\Vert \mathsf{P}_{\mathcal{C}_{A}^{\mathcal{G},{%
shift}}}^{\omega ,\mathbf{b}^{\ast }}g\right\Vert _{L^{2}(\omega
)}^{\bigstar }
\end{eqnarray*}%
This estimate is summable in $s\geq \mathbf{r}$ since $\varepsilon <\frac{1}{%
n+1-\alpha }$, and so the proof of 
\begin{eqnarray}
&& \label{est neigh}
 \left\vert \mathsf{B}_{{neighbour}}^{A}\left( f,g\right) \right\vert
\leq 
\sum_{I\in \mathcal{C}_{A}}\sum_{\substack{ I_{0}\text{ and }I_{\theta }\in 
\mathfrak{C}_{\mathcal{D}}\left( I\right)  \\ I_{0}\neq I_{\theta }}}\sum_{s=\mathbf{r}}^{\infty }A(I,I_{0},I_{\theta },s)   \\
&\lesssim &
\sqrt{\mathfrak{A}_{2}^{\alpha }}\left( \left\Vert \mathsf{P}_{%
\mathcal{C}_{A}}^{\sigma }f\right\Vert _{L^{2}(\sigma )}^{\bigstar }+
\sqrt{\sum_{A^{\prime }\in \mathfrak{C}_{A}\left( A\right) }\left\vert A^{\prime
}\right\vert _{\sigma }\alpha _{\mathcal{A}}\left( A^{\prime }\right) ^{2}}%
\right) \left\Vert \mathsf{P}_{\mathcal{C}_{A}^{\mathcal{G},{shift}%
}}^{\omega ,\mathbf{b}^{\ast }}g\right\Vert _{L^{2}(\omega )}^{\bigstar } 
\notag
\end{eqnarray}%
is complete since $E_{A^{\prime }}^{\sigma }\left\vert f\right\vert \lesssim
\alpha _{\mathcal{A}}\left( A^{\prime }\right) $.

Now if we sum in $A\in \mathcal{A}$ the inequalities (\ref{est para}), (\ref{est neigh}) and (\ref{broken vanish}) we get
\begin{eqnarray*}
&&\sum_{A\in \mathcal{A}}\left\vert \mathsf{B}_{\Subset _{\mathbf{r}
,\varepsilon }}^{A}\left( f,g\right) +\mathsf{B}_{{stop}}^{A}\left(
f,g\right) \right\vert \\
&\lesssim &
\left( \mathfrak{T}_{T^{\alpha }}^{\mathbf{b}}+\sqrt{\mathfrak{A}%
_{2}^{\alpha }}\right)
\sqrt{\sum_{A\in 
\mathcal{A}}\left\Vert \mathsf{P}_{\mathcal{C}_{A}^{\mathcal{G},{%
shift}}}^{\omega ,\mathbf{b}^{\ast }}g\right\Vert _{L^{2}\left( \omega
\right) }^{\bigstar 2}}
\cdot \\ 
&&
\cdot\sqrt{\sum_{A\in \mathcal{A}}\left\{ \alpha _{%
\mathcal{A}}\left( A\right) ^{2}\left\vert A\right\vert _{\sigma
}+\left\Vert \mathsf{P}_{\mathcal{C}_{A}}^{\sigma }f\right\Vert
_{L^{2}(\sigma )}^{\bigstar 2}+\sum_{A^{\prime }\in \mathfrak{C}_{\mathcal{A}%
}\left( A\right) }\alpha _{\mathcal{A}}\left( A^{\prime }\right)
^{2}\left\vert A^{\prime }\right\vert _{\sigma }\right\} }  \\ 
&\lesssim &
\left( \mathfrak{T}_{T^{\alpha }}^{\mathbf{b}}+\sqrt{\mathfrak{A}%
_{2}^{\alpha }}\right) \left\Vert f\right\Vert _{L^{2}\left( \sigma \right)
}\left\Vert g\right\Vert _{L^{2}\left( \omega \right) }
\end{eqnarray*}%
The stopping form is the subject of the following section.

\section{The stopping form}\label{Sec stop}

Here we deal with the stopping form. We modify the adaptation of the argument of M. Lacey in to apply in
the setting of a $Tb$ theorem for an $\alpha $-fractional Calder\'{o}n-Zygmund operator $T^{\alpha }$ in $\mathbb{R}^n$ using the Monotonicity
Lemma \ref{mono}, the energy condition, and the weak goodness of Hyt\"{o}nen and Martikainen \cite
{HyMa}. We directly control the pairs $\left( I,J\right) $ in the stopping
form according to the $\mathcal{L}\,$-coronas (constructed from the `bottom
up' with stopping times involving the energies $\left\Vert \square
_{J}^{\omega ,\mathbf{b}^{\ast }}\right\Vert _{L^{2}\left( \omega \right)
}^{2}$) to which $I$ and $J^{\maltese }$ are associated. However, due to the
fact that the cubes $I$ need no longer be good in any sense, we must
introduce an additional top/down `indented' corona construction on top of
the bottom/up construction of M. Lacey, and in connection with this we
introduce a Substraddling Lemma. We then control the stopping form by
absorbing the case when both $I$ and $J^{\maltese }$ belong to the same $
\mathcal{L}\,$-corona, and by using the Straddling and Substraddling Lemmas,
together with the Orthogonality Lemma, to control the case when $I$ and $
J^{\maltese }$ lie in different coronas, with a geometric gain coming from
the separation of the coronas. This geometric gain is where the new
`indented' corona is required.

Apart from this change, the remaining modifications are more cosmetic, such
as

\begin{itemize}
\item the use of the weak goodness of Hyt\"{o}nen and Martikainen \cite{HyMa}
for pairs $\left( I,J\right) $ arising in the stopping form, rather than
goodness for all cubes $J$ that was available in \cite{Lac}, \cite{SaShUr7}, 
\cite{SaShUr9} and \cite{SaShUr10}. For the most part definitions such as
admissible collections are modified to require $J^{\maltese }\subset I$;

\item the pseudoprojections $\square _{I}^{\sigma ,\mathbf{b}}, \square _{J}^{\omega ,\mathbf{b}^{\ast }}$ are used in place of the orthogonal Haar
projections, and the frame and weak Riesz inequalities compensate for the
lack of orthogonality.
\end{itemize}
Fix grids $\mathcal{D}$ and $\mathcal{G}$. We will prove the bound
\begin{equation}
\left\vert \mathsf{B}_{ {stop}}^{A}\left( f,g\right) \right\vert
\lesssim \mathcal{NTV}_{\alpha }\left\Vert \mathsf{P}_{\mathcal{C}_{A}^{
\mathcal{D}}}^{\sigma ,\mathbf{b}}f\right\Vert _{L^{2}\left( \sigma \right)
}^{\bigstar }\left\Vert \mathsf{P}_{\mathcal{C}_{A}^{\mathcal{G}, {
shift}}}^{\sigma ,\mathbf{b}}g\right\Vert _{L^{2}\left( \omega \right)
}^{\bigstar }\ ,  \label{B stop form 3}
\end{equation}
where we recall that the nonstandard `norms' are given  by,
\begin{eqnarray*}
\left\Vert \mathsf{P}_{\mathcal{C}_{A}^{\mathcal{D}}}^{\sigma ,\mathbf{b}
}f\right\Vert _{L^{2}\left( \sigma \right) }^{\bigstar 2} &\equiv
&\sum_{I\in \mathcal{C}_{A}^{\mathcal{D}}}\left\Vert \square _{I}^{\sigma ,
\mathbf{b}}f\right\Vert _{L^{2}\left( \sigma \right) }^{2}, \\
\left\Vert \mathsf{P}_{\mathcal{C}_{A}^{\mathcal{G}, {shift}
}}^{\sigma ,\mathbf{b}}g\right\Vert _{L^{2}\left( \omega \right) }^{\bigstar
2} &\equiv &\sum_{J\in \mathcal{C}_{A}^{\mathcal{G}, {shift}
}}\left\Vert \square _{J}^{\omega ,\mathbf{b}^{\ast }}g\right\Vert
_{L^{2}\left( \omega \right) }^{2},
\end{eqnarray*}
and that the stopping form is given by
\begin{eqnarray*}
\mathsf{B}_{ {stop}}^{A}\left( f,g\right) &\equiv &\sum_{\substack{ 
I\in \mathcal{C}_{A}^{\mathcal{D}}\text{ and }J\in \mathcal{C}_{A}^{\mathcal{
G}, {shift}}  \\ J^{\maltese }\subsetneqq I\text{ and }\ell \left(
J\right) \leq 2^{-\mathbf{\rho }}\ell \left( I\right) }}\left( E
_{I_{J}}^{\sigma }\widehat{\square }_{I}^{\sigma,\flat ,\mathbf{b}}f\right)
\left\langle T_{\sigma }^{\alpha }\left( b_{A}\mathbf{1}_{A\backslash
I_{J}}\right) ,\square _{J}^{\omega ,\mathbf{b}^{\ast }}g\right\rangle
_{\omega }  \label{dummy} \\
&=&
\sum_{\substack{ I:\ \pi I\in \mathcal{C}_{A}^{\mathcal{D}}\text{ and }
J\in \mathcal{C}_{A}^{\mathcal{G}, {shift}}  \\ J^{\maltese }\subsetneqq I
\text{ and }\ell \left( J\right) \leq 2^{-\left( \mathbf{\rho }-1\right)
}\ell \left( I\right) }}\left( E_{I}^{\sigma }\widehat{\square }
_{\pi I}^{\sigma,\flat ,\mathbf{b}}f\right) \left\langle T_{\sigma }^{\alpha
}\left( b_{A}\mathbf{1}_{A\backslash I}\right) ,\square _{J}^{\omega ,\mathbf{
b}^{\ast }}g\right\rangle _{\omega }  \notag
\end{eqnarray*}
where we have made the `change of dummy variable' $I_{J}\rightarrow I$ for
convenience in notation (recall that the child of $I$ that contains $J$ is
denoted $I_{J}$). Changing $\rho-1$ to $\rho$ we have:
\begin{eqnarray*}
\mathsf{B}_{ {stop}}^{A}\left( f,g\right)=\sum_{\substack{ I:\ \pi I\in \mathcal{C}_{A}^{\mathcal{D}}\text{ and }
J\in \mathcal{C}_{A}^{\mathcal{G}, {shift}}  \\ J^{\maltese }\subsetneqq I
\text{ and }\ell \left( J\right) \leq 2^{-\mathbf{\rho}
}\ell \left( I\right) }}\left( E_{I}^{\sigma }\widehat{\square }
_{\pi I}^{\sigma,\flat,\mathbf{b}}f\right) \left\langle T_{\sigma }^{\alpha
}\left( b_{A}\mathbf{1}_{A\backslash I}\right) ,\square _{J}^{\omega ,\mathbf{
b}^{\ast }}g\right\rangle _{\omega },  \notag
\end{eqnarray*}

For $A\in \mathcal{A}$ recall that we have defined the \emph{shifted} $
\mathcal{G}$-corona by 
\begin{equation*}
\mathcal{C}_{A}^{\mathcal{G}, {shift}}\equiv \left\{ J\in \mathcal{G}
:J^{\maltese }\in \mathcal{C}_{A}^{\mathcal{D}}\right\} ,
\end{equation*}
and also defined the \emph{restricted} $\mathcal{D}$-corona by 
\begin{equation*}
\mathcal{C}_{A}^{\mathcal{D}, {restrict}}\equiv \mathcal{C}
_{A}\backslash \left\{ A\right\} \equiv \mathcal{C}_{A}'.
\end{equation*}

\begin{dfn}
Suppose that $A\in \mathcal{A}$ and that $\mathcal{P}\subset \mathcal{C}
_{A}^{\mathcal{D}, {restrict}}\times \mathcal{C}_{A}^{\mathcal{G},
 {shift}}$. We say that the collection of pairs $\mathcal{P}$ is $A$
\emph{-admissible} if
\end{dfn}

\begin{itemize}
\item (good and $\left( \mathbf{\rho},\varepsilon \right) $-deeply
embedded) For every $\left( I,J\right) \in \mathcal{P},$  and $J^{\maltese }\subset I\varsubsetneqq A$. 

\item (tree-connected in the first component) if $I_{1}\subset I_{2}$ and
both $\left( I_{1},J\right) \in \mathcal{P}$ and $\left( I_{2},J\right) \in 
\mathcal{P}$, then $\left( I,J\right) \in \mathcal{P}$ for every $I$ in the
geodesic $\left[ I_{1},I_{2}\right] =\left\{ I\in \mathcal{D}:I_{1}\subset
I\subset I_{2}\right\} $.
\end{itemize}

From now on we often write $\mathcal{C}_{A}$ and $\mathcal{C}_{A}'$ in place of $\mathcal{C}_{A}^{\mathcal{D}}$ and $\mathcal{C}
_{A}^{\mathcal{D}, {restrict}}$ respectively when there is no
confusion. The basic example of an admissible collection of pairs is
obtained from the pairs of cubes summed in the stopping form $\mathsf{B}
_{stop}^{A}\left( f,g\right)$,
\begin{equation}
\mathcal{P}^{A}\equiv \left\{ \left( I,J\right) :I\in \mathcal{C}
_{A}^{\prime }\text{ and }J\in \mathcal{G}_{\left(\rho,\varepsilon \right) - {good}}^{\mathcal{D}}\cap \mathcal{C}_{A}^{
\mathcal{G}, {shift}}\text{ where\ }J\Subset _{\mathbf{\rho}
,\varepsilon }I\right\} .  \label{initial P}
\end{equation}

\begin{dfn}\label{def stop P}
Suppose that $A\in \mathcal{A}$ and that $\mathcal{P}$ is an $A$\emph{
-admissible} collection of pairs. Define the associated \emph{stopping} form 
$\mathsf{B}_{ {stop}}^{A,\mathcal{P}}$ by
\begin{equation*}
\mathsf{B}_{ {stop}}^{A,\mathcal{P}}\left( f,g\right) \equiv
\sum_{\left( I,J\right) \in \mathcal{P}}\left( E_{I}^{\sigma }
\widehat{\square }_{\pi I}^{\sigma,\flat ,\mathbf{b}}f\right) \ \left\langle
T_{\sigma }^{\alpha }\left( b_{A}\mathbf{1}_{A\backslash I}\right) ,\square
_{J}^{\omega ,\mathbf{b}^{\ast }}g\right\rangle _{\omega }\ .
\end{equation*}
\end{dfn}


\begin{prop}
\label{stopping bound}Suppose that $A\in \mathcal{A}$ and that $\mathcal{P}$
is an $A$\emph{-admissible} collection of pairs. Then the stopping form $
\mathsf{B}_{ {stop}}^{A,\mathcal{P}}$ satisfies the bound
\begin{equation}
\left\vert \mathsf{B}_{ {stop}}^{A,\mathcal{P}}\left( f,g\right)
\right\vert \lesssim \left( \mathcal{E}^{\alpha }_2+\sqrt{
\mathfrak{A}_{2}^{\alpha }}\right)  \left\Vert \mathsf{P}_{\mathcal{C}
_{A}}^{\sigma ,\mathbf{b}}f\right\Vert _{L^{2}\left( \sigma \right)
}^{\bigstar } \left\Vert \mathsf{P}_{\mathcal{C}_{A}^{
\mathcal{G}, {shift}}}^{\omega ,\mathbf{b}^{\ast }}g\right\Vert
_{L^{2}\left( \omega \right) }^{\bigstar }\  \label{B stop form}
\end{equation}
\end{prop}

With the above proposition in hand, we can complete the proof of (\ref{B
stop form 3}) by summing over the stopping cubes $A\in \mathcal{A}$ with the
choice $\mathcal{P}^{A}$ of $A$-admissible pairs for each $A$: 
\begin{eqnarray*}
&&\sum_{A\in \mathcal{A}}\left\vert \mathsf{B}_{ {stop}}^{A,\mathcal{P
}^{A}}\left( f,g\right) \right\vert \\
&\lesssim &
\sum_{A\in \mathcal{A}}\left( \mathcal{E}^{\alpha }_2+\sqrt{
\mathfrak{A}_{2}^{\alpha }}\right) 
\left\Vert \mathsf{P}_{\mathcal{C}
_{A}}^{\sigma ,\mathbf{b}}f\right\Vert _{L^{2}\left( \sigma \right)
}^{\bigstar } \left\Vert \mathsf{P}_{\mathcal{C}_{A}^{
\mathcal{G}, {shift}}}^{\omega ,\mathbf{b}^{\ast }}g\right\Vert
_{L^{2}\left( \omega \right) }^{\bigstar } \\
&\lesssim &
\left( \mathcal{E}^{\alpha }_2+\sqrt{\mathfrak{A}_{2}^{\alpha }}\right) \left( \sum_{A\in \mathcal{A}}
\left\Vert \mathsf{P}_{
\mathcal{C}_{A}}^{\sigma ,\mathbf{b}}f\right\Vert _{L^{2}\left( \sigma
\right) }^{\bigstar 2} \right) ^{\frac{1}{2}}\left( \sum_{A\in 
\mathcal{A}}\left\Vert \mathsf{P}_{\mathcal{C}_{A}^{\mathcal{G}, {
shift}}}^{\omega ,\mathbf{b}^{\ast }}g\right\Vert _{L^{2}\left( \omega
\right) }^{\bigstar 2}\right) ^{\frac{1}{2}} \\
&\lesssim &
\left( \mathcal{E}^{\alpha }_2+\sqrt{\mathfrak{A}_{2}^{\alpha }}\right) \left\Vert f\right\Vert _{L^{2}\left( \sigma \right) }\left\Vert
g\right\Vert _{L^{2}\left( \omega \right) }
\end{eqnarray*}
by the lower Riesz inequality $\displaystyle\sum_{A\in \mathcal{A}}\left\Vert \mathsf{P}_{
\mathcal{C}_{A}}^{\sigma ,\mathbf{b}}f\right\Vert _{L^{2}\left( \sigma
\right) }^{\bigstar 2}\lesssim \left\Vert f\right\Vert _{L^{2}\left( \sigma
\right) }^{2}$, quasi-orthogonality $\displaystyle\sum_{A\in \mathcal{A}}\alpha _{
\mathcal{A}}\left( f\right) ^{2}\left\vert A\right\vert _{\sigma }\lesssim
\left\Vert f\right\Vert _{L^{2}\left( \sigma \right) }^{2}$ in the stopping
cubes $\mathcal{A}$, and by the pairwise disjointedness of the shifted
coronas $\mathcal{C}_{A}^{\mathcal{G}, {shift}}$: 
$\displaystyle
\sum_{A\in \mathcal{A}}\mathbf{1}_{\mathcal{C}_{A}^{\mathcal{G}, {
shift}}}\leq \mathbf{1}_{\mathcal{D}}.
$

To prove Proposition \ref{stopping bound}, we begin by letting 
\begin{eqnarray*}
\Pi _{1}\mathcal{P} &\equiv &\left\{ I\in \mathcal{C}_{A}^{\mathcal{D},restrict
}:\left( I,J\right) \in \mathcal{P}\text{ for some }J\in 
\mathcal{C}_{A}^{\mathcal{G},{shift}}\right\} , \\
\Pi _{2}\mathcal{P} &\equiv &\left\{ J\in \mathcal{C}_{A}^{\mathcal{G},{
shift}}:\left( I,J\right) \in \mathcal{P}\text{ for some }I\in \mathcal{C}
_{A}'\right\} ,
\end{eqnarray*}
consist of the first and second components respectively of the pairs in $
\mathcal{P}$, and writing
\begin{eqnarray*}
\mathsf{B}_{ {stop}}^{A,\mathcal{P}}\left( f,g\right) &=&\sum_{J\in
\Pi _{2}\mathcal{P}}\left\langle T_{\sigma }^{\alpha }\varphi _{J}^{\mathcal{
P}},\square _{J}^{\omega ,\mathbf{b}^{\ast }}g\right\rangle _{\omega }; \\
\text{where }\varphi _{J}^{\mathcal{P}} &\equiv &\sum_{I\in \mathcal{C}
_{A}^{\prime }:\ \left( I,J\right) \in \mathcal{P}}b_{A}E
_{I}^{\sigma }\left( \widehat{\square }_{\pi I}^{\sigma,\flat ,\mathbf{b}}f\right)
\ \mathbf{1}_{A\backslash I}\ \text{(since }b_{I}=b_{A}\text{ for }I\in 
\mathcal{C}_{A}\text{)}.
\end{eqnarray*}
By the tree-connected property of $\mathcal{P}$, and the telescoping
property of dual martingale differences, together with the bound $\alpha _{
\mathcal{A}}\left( A\right) $ on the averages of $f$ in the corona $\mathcal{
C}_{A}$, we have
\begin{equation}
\left\vert \varphi _{J}^{\mathcal{P}}\right\vert \lesssim \alpha _{\mathcal{A
}}\left( A\right) 1_{A\backslash I_{\mathcal{P}}\left( J\right) },
\label{phi bound}
\end{equation}
where $I_{\mathcal{P}}\left( J\right) \equiv \bigcap \left\{ I:\left(
I,J\right) \in \mathcal{P}\right\} $ is the smallest cube $I$ for which $
\left( I,J\right) \in \mathcal{P}$. It is important to note that $J$ is good
with respect to $I_{\mathcal{P}}\left( J\right) $ by our infusion of weak
goodness above. Another important property of these functions is the
sublinearity:
\begin{equation}
\left\vert \varphi _{J}^{\mathcal{P}}\right\vert \leq \left\vert \varphi
_{J}^{\mathcal{P}_{1}}\right\vert +\left\vert \varphi _{J}^{\mathcal{P}
_{2}}\right\vert ,\ \ \ \ \ \mathcal{P}=\mathcal{P}_{1}\dot{\cup}\mathcal{P}
_{2}\ .  \label{phi sublinear}
\end{equation}
Now apply the Monotonicity Lemma \ref{mono} to the inner product $
\left\langle T_{\sigma }^{\alpha }\varphi _{J},\square _{J}^{\omega
}g\right\rangle _{\omega }$ to obtain
\begin{eqnarray*}
\left\vert \left\langle T_{\sigma }^{\alpha }\varphi _{J},\square
_{J}^{\omega ,\mathbf{b}^{\ast }}g\right\rangle _{\omega }\right\vert
&\lesssim &
\frac{\mathrm{P}^{\alpha }\left( J,\left\vert \varphi
_{J}\right\vert \mathbf{1}_{A\backslash I_{\mathcal{P}}\left( J\right)
}\sigma \right) }{\left\vert J\right\vert ^{\frac{1}{n}}}\left\Vert
\bigtriangleup _{J}^{\omega ,\mathbf{b}^{\ast }}x\right\Vert _{L^{2}\left(
\omega \right) }^{\spadesuit }\left\Vert \square _{J}^{\omega ,\mathbf{b}
^{\ast }}g\right\Vert _{L^{2}\left( \omega \right) }^{\bigstar }  \\
&+&
\frac{\mathrm{P}_{1+\delta }^{\alpha }\left( J,\left\vert \varphi
_{J}\right\vert \mathbf{1}_{A\backslash I_{\mathcal{P}}\left( J\right)
}\sigma \right) }{\left\vert J\right\vert ^{\frac{1}{n}}}\left\Vert
x-m_{J}^{\omega }\right\Vert _{L^{2}\left( \mathbf{1}_{J}\omega \right)
}\left\Vert \square _{J}^{\omega ,\mathbf{b}^{\ast }}g\right\Vert
_{L^{2} }^{\bigstar }
\end{eqnarray*}
Thus we have
\begin{eqnarray}
&& \label{def split} \\
\left\vert \mathsf{B}_{ {stop}}^{A,\mathcal{P}}\left( f,g\right)
\right\vert 
&\leq &
\sum_{J\in \Pi _{2}\mathcal{P}}\frac{\mathrm{P}^{\alpha }\left( J,\left\vert \varphi
_{J}\right\vert \mathbf{1}_{A\backslash I_{\mathcal{P}}\left( J\right)
}\sigma \right) }{\left\vert J\right\vert ^{\frac{1}{n}}}\left\Vert
\bigtriangleup _{J}^{\omega ,\mathbf{b}^{\ast }}x\right\Vert _{L^{2}\left(
\omega \right) }^{\spadesuit }\left\Vert \square _{J}^{\omega ,\mathbf{b}
^{\ast }}g\right\Vert _{L^{2}\left( \omega \right) }^{\bigstar }
\notag \\
&+&
\sum_{J\in \Pi _{2}\mathcal{P}}\frac{\mathrm{P}_{1+\delta }^{\alpha }\left( J,\left\vert \varphi
_{J}\right\vert \mathbf{1}_{A\backslash I_{\mathcal{P}}\left( J\right)
}\sigma \right) }{\left\vert J\right\vert ^{\frac{1}{n}}}\left\Vert
x-m_{J}^{\omega }\right\Vert _{L^{2}\left( \mathbf{1}_{J}\omega \right)
}\left\Vert \square _{J}^{\omega ,\mathbf{b}^{\ast }}g\right\Vert
_{L^{2} }^{\bigstar } \notag \\
&\equiv &
\left\vert \mathsf{B}\right\vert _{ {stop},1,\bigtriangleup
^{\omega }}^{A,\mathcal{P}}\left( f,g\right) +\left\vert \mathsf{B}
\right\vert _{ {stop},1+\delta ,\mathsf{P}^{\omega }}^{A,\mathcal{P}
}\left( f,g\right) ,  \notag
\end{eqnarray}
where we have dominated the stopping form by two sublinear stopping forms
that involve the Poisson integrals of order $1$ and $1+\delta $
respectively, and where the smaller Poisson integral $\mathrm{P}_{1+\delta
}^{\alpha }$ is multiplied by the larger quantity $\left\Vert
x-m_{J}^{\omega }\right\Vert _{L^{2}}\left( \mathbf{1}_{J}\omega \right)$. This splitting turns out to be successful in separating the two
energy terms from the right hand side of the Energy Lemma, because of the
two properties (\ref{phi bound}) and (\ref{phi sublinear}) above. It remains
to show the two inequalities:
\begin{equation}
\left\vert \mathsf{B}\right\vert _{{stop},\bigtriangleup ^{\omega
}}^{A,\mathcal{P}}\left( f,g\right) \lesssim \left( \mathcal{E}^{\alpha }_2+\sqrt{\mathfrak{A}_{2}^{\alpha }}\right) \left\Vert \mathsf{P}_{\pi \left(
\Pi _{1}\mathcal{P}\right) }^{\sigma ,\mathbf{b}}f\right\Vert _{L^{2}\left(
\sigma \right) }^{\bigstar }\left\Vert \mathsf{P}_{\Pi _{2}\mathcal{P}
}^{\omega ,\mathbf{b}^{\ast }}g\right\Vert _{L^{2}\left( \omega \right)
}^{\bigstar },  \label{First inequality}
\end{equation}
for $f\in L^{2}\left( \sigma \right) $ satisfying where $E
_{I}^{\sigma }\left\vert f\right\vert \leq \alpha _{\mathcal{A}}\left(
A\right) $ for all $I\in \mathcal{C}_{A}$; and where $\pi \left( \Pi _{1}\mathcal{P}\right) \equiv
\left\{ \pi _{\mathcal{D}}I:I\in \Pi _{1}\mathcal{P}\right\} $; and
\begin{equation}
\left\vert \mathsf{B}\right\vert _{ {stop},1+\delta ,\mathsf{P}
^{\omega }}^{A,\mathcal{P}}\left( f,g\right) \lesssim \left( \mathcal{E}^{\alpha }_2+\sqrt{\mathfrak{A}_{2}^{\alpha }}\right) \left\Vert \mathsf{P}_{\mathcal{C}_{A}^{\mathcal{D}
}}^{\sigma ,\mathbf{b}}f\right\Vert _{L^{2}\left( \sigma \right) }\left\Vert 
\mathsf{P}_{\mathcal{C}_{A}^{\mathcal{G},{shift}}}^{\omega ,\mathbf{b}
^{\ast }}g\right\Vert _{L^{2}\left( \omega \right) } \label{Second inequality}
\end{equation}
where we only need the case $\mathcal{P}=\mathcal{P}^{A}$ in this latter
inequality as there is no recursion involved in treating this second
sublinear form. We consider first the easier inequality (\ref{Second
inequality}) that does not require recursion.

\subsection{The bound for the second sublinear inequality}

Now we turn to proving (\ref{Second inequality}), i.e.
\begin{equation*}
\left\vert \mathsf{B}\right\vert _{ {stop},1+\delta ,\mathsf{P}
^{\omega }}^{A,\mathcal{P}}\left( f,g\right)
\lesssim
\left( \mathcal{E}_{2}^{\alpha }+\sqrt{\mathfrak{A}_{2}^{\alpha }}\right)
\left\Vert \mathsf{P}_{\mathcal{C}_{A}^{\mathcal{D}}}^{\sigma ,\mathbf{b}
}f\right\Vert _{L^{2}\left( \sigma \right) }\left\Vert \mathsf{P}_{\mathcal{C
}_{A}^{\mathcal{G},{shift}}}^{\omega ,\mathbf{b}^{\ast }}g\right\Vert
_{L^{2}\left( \omega \right) }
\end{equation*}
where since 
\begin{equation*}
\left\vert \varphi _{J}\right\vert =\left\vert \sum_{I\in \mathcal{C}
_{A}^{\prime }:\ \left( I,J\right) \in \mathcal{P}}E_{I}^{\sigma
}\left( \widehat{\square }_{\pi I}^{\sigma,\flat ,\mathbf{b}}f\right) \ b_{A}
\mathbf{1}_{A\backslash I}\right\vert \leq \sum_{I\in \mathcal{C}_{A}^{\prime
}:\ \left( I,J\right) \in \mathcal{P}}\left\vert E_{I}^{\sigma
}\left( \widehat{\square }_{\pi I}^{\sigma,\flat ,\mathbf{b}}f\right) b_{A}\ 
\mathbf{1}_{A\backslash I}\right\vert ,
\end{equation*}
the sublinear form $\left\vert \mathsf{B}\right\vert _{ {stop}
,1+\delta ,\mathsf{P}^{\omega }}^{A,\mathcal{P}}$ can be dominated and then
decomposed by pigeonholing the ratio of side lengths of $J$ and $I$:
\begin{eqnarray*}
&&\left\vert \mathsf{B}\right\vert _{{stop},1+\delta }^{A,\mathcal{P}
}\left( f,g\right)\\
&=&
\sum_{J\in \Pi _{2}\mathcal{P}}\frac{\mathrm{P}
_{1+\delta }^{\alpha }\left( J,\left\vert \varphi _{J}\right\vert \mathbf{1}
_{A\backslash I_{\mathcal{P}}\left( J\right) }\sigma \right) }{\left\vert
J\right\vert ^{\frac{1}{n}}}\left\Vert x-m_{J}\right\Vert _{L^{2}\left( 
\mathbf{1}_{J}\omega \right) }\left\Vert \square _{J}^{\omega ,\mathbf{b}
^{\ast }}g\right\Vert _{L^{2}\left( \omega \right) }^{\bigstar } \\
&\leq&
\sum_{\left( I,J\right) \in \mathcal{P}}\frac{\mathrm{P}_{1+\delta
}^{\alpha }\left( J,\left\vert E_{I}^{\sigma }\left( \square _{\pi
I}^{\sigma ,\flat ,\mathbf{b}}f\right) \right\vert \mathbf{1}_{A\backslash
I}\sigma \right) }{\left\vert J\right\vert ^{\frac{1}{n}}}\left\Vert
x-m_{J}\right\Vert _{L^{2}\left( \mathbf{1}_{J}\omega \right) }\left\Vert
\square _{J}^{\omega ,\mathbf{b}^{\ast }}g\right\Vert _{L^{2}\left( \omega
\right) }^{\bigstar } \\
&\equiv&
\sum_{s=0}^{\infty }\left\vert \mathsf{B}\right\vert _{{stop
},1+\delta }^{A,\mathcal{P};s}\left( f,g\right) ; \\
\end{eqnarray*}
We will now adapt the argument for the stopping term starting on page 42 of \cite
{LaSaUr2}, where the geometric gain from the assumed `Energy Hypothesis'
there will be replaced by a geometric gain from the smaller Poisson integral 
$\mathrm{P}_{1+\delta }^{\alpha }$ used here.

First, we exploit the additional decay in the Poisson integral $\mathrm{P}
_{1+\delta }^{\alpha }$ as follows. Suppose that $\left( I,J\right) \in 
\mathcal{P}$ with $\ell \left( J\right) =2^{-s}\ell \left( I\right) $. We
then compute
\begin{eqnarray*}
\frac{\mathrm{P}_{1+\delta }^{\alpha }\left( J,\left\vert b_{A}\right\vert 
\mathbf{1}_{A\backslash I}\sigma \right) }{\left\vert J\right\vert ^{\frac{1}{
n}}}
&\approx &
\int_{A\backslash I}\frac{\left\vert J\right\vert ^{\frac{
\delta }{n}}}{\left\vert y-c_{J}\right\vert ^{n+1+\delta -\alpha }}
\left\vert b_{A}\left( y\right) \right\vert d\sigma \left( y\right) \\
&\leq &
\int_{A\backslash I}\left( \frac{\left\vert J\right\vert ^{\frac{1}{n}}
}{\dist\left( c_{J},I^{c}\right) }\right) ^{\delta }\frac{1}{\left\vert
y-c_{J}\right\vert ^{n+1-\alpha }}\left\vert b_{A}\left( y\right)
\right\vert d\sigma \left( y\right) \\
&\lesssim &
\left( \frac{\left\vert J\right\vert ^{\frac{1}{n}}}{\dist\left(
c_{J},I^{c}\right) }\right) ^{\delta }\frac{\mathrm{P}^{\alpha }\left(
J,\left\vert b_{A}\right\vert \mathbf{1}_{A\backslash I}\sigma \right) }{
\left\vert J\right\vert ^{\frac{1}{n}}},
\end{eqnarray*}
and using the goodness of $J$ in $I$,
\begin{equation*}
d\left( c_{J},I^{c}\right) \geq 2\ell \left( I\right)
^{1-\varepsilon }\ell \left( J\right) ^{\varepsilon }\geq 2\cdot 2^{s\left( 1-\varepsilon \right) }\ell \left( J\right) ,
\end{equation*}
to conclude, using accretivity, that
\begin{equation}
\left( \frac{\mathrm{P}_{1+\delta }^{\alpha }\left( J, |b_A| \mathbf{1}_{A\backslash I}\sigma \right) }{\left\vert
J\right\vert ^{\frac{1}{n}}}\right) \lesssim 2^{-s\delta \left(
1-\varepsilon \right) }\frac{\mathrm{P}^{\alpha }\left( J, \mathbf{1}_{A\backslash I}\sigma \right) }{\left\vert
J\right\vert ^{\frac{1}{n}}}.  \label{Poisson decay}
\end{equation}

We next claim that for $s\geq 0$ an integer,
\begin{eqnarray*}
\left\vert \mathsf{B}\right\vert _{ {stop},1+\delta ,\mathsf{P}
^{\omega }}^{A,\mathcal{P};s}\left( f,g\right)
&\lesssim &
2^{-s\delta \left( 1-\varepsilon \right) }\ \left( \mathcal{E}_{2}^{\alpha }+\sqrt{\mathfrak{A}_{2}^{\alpha }}\right)
\left\Vert \mathsf{P}_{\mathcal{C}_{A}^{\mathcal{D}}}^{\sigma ,\mathbf{b}
}f\right\Vert _{L^{2}\left( \sigma \right) }\left\Vert \mathsf{P}_{\mathcal{C
}_{A}^{\mathcal{G},{shift}}}^{\omega ,\mathbf{b}^{\ast }}g\right\Vert
_{L^{2}\left( \omega \right) }
\end{eqnarray*}
from which (\ref{Second inequality}) follows upon summing in $s\geq 0$. Now
using both
$$
\left\vert E_{I}^{\sigma }\left( \widehat{\square }_{\pi I}^{\sigma,\flat
,\mathbf{b}}f\right) \right\vert\frac{1}{\left\vert I\right\vert _{\sigma }}
\int_{I}\left\vert \square _{\pi I}^{\sigma,\flat ,\mathbf{b}}f\right\vert d\sigma
\leq \left\Vert \square _{\pi I}^{\sigma,\flat ,\mathbf{b}}f\right\Vert
_{L^{2}\left( \sigma \right) }\frac{1}{\sqrt{\left\vert I\right\vert
_{\sigma }}},
$$
$$
\sum_{I\in \mathcal{D}}\left\Vert \square _{\pi I}^{\sigma ,\flat ,\mathbf{b}
}f\right\Vert _{L^{2}\left( \sigma \right) }^{2} 
\lesssim 
\sum_{I\in\mathcal{D}}\left( \left\Vert \square _{\pi I}^{\sigma ,\mathbf{b}
}f\right\Vert _{L^{2}\left( \sigma \right) }^{2}+\left\Vert \nabla _{\pi
I}^{\sigma }f\right\Vert _{L^{2}\left( \sigma \right) }^{2}\right) \approx
\left\Vert f\right\Vert _{L^{2}(\sigma )}^{2}\ ,
$$
we apply Cauchy-Schwarz in the $I$ variable above to see that 
\begin{eqnarray*}
&&\!\!\!\!\left[ \left\vert \mathsf{B}\right\vert _{ {stop},1+\delta ,\mathsf{
P}^{\omega }}^{A,\mathcal{P};s}\left( f,g\right) \right] ^{2} \\
&&\hspace{-0.75cm}\lesssim\!\!
\left\Vert \mathsf{P}_{\mathcal{C}_{A}^{\mathcal{D}}}^{\sigma ,\mathbf{b}
}f\right\Vert _{L^{2}\left( \sigma \right) }\!\!\! \left[ \sum_{I\in \mathcal{C}_{A}^{\prime }}\!\!\!\left(\!\!\!  \frac{1}{\sqrt{
\left\vert I\right\vert _{\sigma }}}\!\!\!\!\!\!\!\sum_{\substack{\ \ \  J:\ \left( I,J\right)
\in \mathcal{P}  \\ \ell \left( J\right) =2^{-s}\ell \left( I\right) }}\!\!\!\!\!\!\!\!\!\!\frac{
\mathrm{P}_{1+\delta }^{\alpha }\left( J, \mathbf{
1}_{A\backslash I}\sigma \right) }{\left\vert J\right\vert ^{\frac{1}{n}}}
\!\!\left\Vert
x-m_{J}\right\Vert _{L^{2}\left( \mathbf{1}_{J}\omega \right) }\!\!\left\Vert
\square _{J}^{\omega ,\mathbf{b}^{\ast }}g\right\Vert _{L^{2}\left( \omega
\right) }^{\bigstar }\!\!\! \right) ^{\!\!\! 2}\right]^{\!\!\frac{1}{2}}
\end{eqnarray*}
Using the frame inequality for $\square _{J}^{\omega ,\mathbf{b}^{\ast }}$
we can then estimate the sum inside the square brackets by
\begin{eqnarray*}
\sum_{I\in \mathcal{C}_{A}^{\prime }}
\!\!\!\left\{ \sum_{\substack{ J:\ \left(
I,J\right) \in \mathcal{P}  \\ \ell \left( J\right) =2^{-s}\ell \left(
I\right) }}\!\!\!\!\!\!\left\Vert \square _{J}^{\omega ,\mathbf{b}^{\ast }}g\right\Vert
_{L^{2}\left( \omega \right) }^{\bigstar2}\!\!\!\right\} \!\!\sum_{\substack{ J:\ \left(
I,J\right) \in \mathcal{P}  \\ \ell \left( J\right) =2^{-s}\ell \left(
I\right) }}
\!\!\!\!&\displaystyle\!\!\!\!\!\frac{1}{\left\vert I\right\vert _{\sigma }}&\!\!\!\!\!\!
\left( \frac{\mathrm{P}_{1+\delta }^{\alpha }\left( J,\mathbf{1}_{A\backslash I}\sigma \right) }{\left\vert J\right\vert ^{\frac{1}{n}}}\right) ^{\!\!2}
\!\!\!\left\Vert x-m_J\right\Vert _{L^{2}\left( \mathbf{1}_J\omega \right) }^{2}\\
&\lesssim& 
\left\Vert \mathsf{P}_{\Pi_2\mathcal{P}}^{\omega ,\mathbf{b}^{\ast }}g\right\Vert
_{L^{2}\left( \omega \right) }^{\bigstar2}A\left( s\right) ^{2},
\end{eqnarray*}
where
\begin{equation*}
A\left( s\right) ^{2}\equiv \sup_{I\in \mathcal{C}_{A}^{\prime }}\sum_{\substack{ J:\ \left(
I,J\right) \in \mathcal{P}  \\ \ell \left( J\right) =2^{-s}\ell \left(
I\right) }}\frac{1}{\left\vert I\right\vert _{\sigma }}
\left( \frac{\mathrm{P}_{1+\delta }^{\alpha }\left( J, \mathbf{1}
_{A\backslash I}\sigma \right) }{\left\vert J\right\vert ^{\frac{1}{n}}}
\right) ^{2}\left\Vert x-m_J\right\Vert _{L^{2}\left( \mathbf{1}_J\omega \right) }^{2}
\end{equation*}
Finally then we turn to the analysis of the supremum in last display. From
the Poisson decay (\ref{Poisson decay}) we have 
\begin{eqnarray*}
A\left( s\right) ^{2} 
&\lesssim &
\sup_{I\in \mathcal{C}_{A}^{\prime }}\frac{1
}{\left\vert I\right\vert _{\sigma }}2^{-2s\delta \left( 1-\varepsilon
\right) }\!\!\!\!
\sum_{\substack{ J:\ \left( I,J\right) \in \mathcal{P}  \\ \ell
\left( J\right) =2^{-s}\ell \left( I\right) }}\left( \frac{\mathrm{P}
^{\alpha }\left( J, \mathbf{1}_{A\backslash
I}\sigma \right) }{\left\vert J\right\vert ^{\frac{1}{n}}}\right)
^{2}\left\Vert x-m_J\right\Vert _{L^{2}\left( \mathbf{1}_J\omega \right) }^{2} \\
&\lesssim &
2^{-2s\delta \left( 1-\varepsilon \right) }\left[ \left( \mathcal{
E}^{\alpha }_2\right) ^{2}+\mathfrak{A}_{2}^{\alpha }\right]
\,,
\end{eqnarray*}
Indeed, from Definition~\ref{def energy corona 3}, as $(I,J)\in \mathcal{P}$
, we have that $I$ is \emph{not} a stopping cube in $\mathcal{A}$, and hence
that (\ref{def stop 3}) \emph{fails} to hold, delivering the estimate above
since $J\Subset _{\mathbf{\rho},\varepsilon }I$ good must be contained in
some $K\in \mathcal{M}_{\left( \mathbf{r},\varepsilon \right) - {deep}
}\left( I\right) $, and since $\frac{\mathrm{P}^{\alpha }\left( J,\left\vert
b_{I}\right\vert \mathbf{1}_{A\backslash I}\sigma \right) }{\left\vert
J\right\vert ^{\frac{1}{n}}}\approx \frac{\mathrm{P}^{\alpha }\left(
K,\left\vert b_{I}\right\vert \mathbf{1}_{A\backslash I}\sigma \right) }{
\left\vert K\right\vert ^{\frac{1}{n}}}$. The terms $\left\Vert \mathsf{P}
_{J}^{\omega }x\right\Vert _{L^{2}\left( \omega \right) }^{2}$ are additive
since the $J^{\prime }s$ are pigeonholed by $\ell \left( J\right)
=2^{-s}\ell \left( I\right) $.

\subsection{The bound for the first sublinear inequality}

Now we turn to proving the more difficult inequality \eqref{First inequality}. Denote by $\mathfrak{N}_{ {stop},\bigtriangleup ^{\omega }}^{A,\mathcal{P}
}$ the best constant in 
\begin{equation}\label{best hat}
\left\vert \mathsf{B}\right\vert _{ {stop},\bigtriangleup ^{\omega }}^{A,
\mathcal{P}}\left( f,g\right) 
\leq 
\mathfrak{N}_{ {stop},\bigtriangleup
^{\omega }}^{A,\mathcal{P}}\left\Vert \mathsf{P}_{\pi(\Pi _{1}\mathcal{P})}^{\sigma ,\mathbf{b}}f\right\Vert _{L^{2}\left( \sigma \right) }^{\bigstar
} \left\Vert \mathsf{P}_{\Pi _{2}\mathcal{P}}^{\omega ,\mathbf{b}
^{\ast }}g\right\Vert _{L^{2}\left( \omega \right) }^{\bigstar },
\end{equation}
where $f\in L^{2}\left( \sigma \right) $ satisfies $E_{I}^{\sigma
}\left\vert f\right\vert \leq \alpha _{\mathcal{A}}\left( A\right) $ for all 
$I\in \mathcal{C}_{A}$, and $g\in L^{2}\left( \omega \right) $ and 
$
\pi(\Pi_1\mathcal{P})=\{\pi I: I\in\Pi_1P\}
$. 
We refer to $\mathfrak{N}_{ {stop},\bigtriangleup
^{\omega }}^{A,\mathcal{P}}$ as the \textit{restricted} norm relative to the collection $\mathcal{P}$. Inequality (\ref{First inequality}) follows once we have shown that $\mathfrak{N}_{ {stop},\bigtriangleup
^{\omega }}^{A,\mathcal{P}}\lesssim \mathcal{E}_{2}^{\alpha }+
\sqrt{\mathfrak{A}_{2}^{\alpha }}$.

The following general result on mutually orthogonal admissible collections
will prove very useful in establishing (\ref{First inequality}). Given a set 
$\left\{ \mathcal{Q}_{m}\right\} _{m=0}^{\infty }$ of admissible collections
for $A$, we say that the collections $\mathcal{Q}_{m}$ are \emph{mutually
orthogonal}, if each collection $\mathcal{Q}_{m}$ satisfies
\begin{equation*}
\mathcal{Q}_{m}\subset \bigcup\limits_{j=0}^{\infty }\left\{ \mathcal{A}
_{m,j}\times \mathcal{B}_{m,j}\right\} \
\end{equation*}
where the sets $\left\{ \mathcal{A}_{m,j}\right\} _{m,j}$ and $\left\{ 
\mathcal{B}_{m,j}\right\} _{m,j}$ are each pairwise disjoint in their respective dyadic
grids $\mathcal{D}$ and $\mathcal{G}$: 
\begin{equation*}
\sum_{m,j=0}^{\infty }\mathbf{1}_{\mathcal{A}_{m,j}}\leq \mathbf{1}_{
\mathcal{D}}\text{ and }\sum_{m,j=0}^{\infty }\mathbf{1}_{\mathcal{B}
_{m,j}}\leq \mathbf{1}_{\mathcal{G}}.
\end{equation*}

\begin{lem}
\label{mut orth}Suppose that $\left\{ \mathcal{Q}_{m}\right\} _{m=0}^{\infty
}$ is a set of admissible collections for $A$ that are \emph{mutually
orthogonal}. Then $\mathcal{Q}\equiv \bigcup\limits_{m=0}^{\infty }\mathcal{
Q}_{m}$ is admissible, and the sublinear stopping form $\left\vert \mathsf{B}
\right\vert _{ {stop},\bigtriangleup ^{\omega }}^{A,\mathcal{Q}}\left(
f,g\right) $ has its restricted norm $\mathfrak{N}_{ {stop},\bigtriangleup
^{\omega }}^{A,\mathcal{Q}}$ controlled by the \emph{supremum} of the
restricted norms $\mathfrak{N}_{ {stop},\bigtriangleup ^{\omega }}^{A,
\mathcal{Q}_{m}}$: 
\begin{equation*}
\mathfrak{N}_{ {stop},\bigtriangleup ^{\omega }}^{A,\mathcal{Q}}\leq \sup_{m\geq 0}\mathfrak{N}_{ {stop},\bigtriangleup ^{\omega }}^{A,\mathcal{Q
}_{m}}.
\end{equation*}
\end{lem}

\begin{proof}
If $J\in \Pi _{2}\mathcal{Q}_{m}$, then $\varphi _{J}^{\mathcal{Q}}=\varphi
_{J}^{\mathcal{Q}_{m}}$ and $I_{\mathcal{Q}}\left( J\right) =I_{\mathcal{Q}
_{m}}\left( J\right) $, since the collection $\left\{ \mathcal{Q}
_{m}\right\} _{m=0}^{\infty }$ is mutually orthogonal. Thus we have 
\begin{eqnarray*}
&&
\left\vert \mathsf{B}\right\vert _{{stop},\bigtriangleup ^{\omega
}}^{A,\mathcal{Q}}\left( f,g\right) 
=
\sum_{J\in \Pi _{2}\mathcal{Q}}\frac{
\mathrm{P}^{\alpha }\left( J,\left\vert \varphi _{J}^{\mathcal{Q}
}\right\vert \mathbf{1}_{A\backslash I_{\mathcal{Q}}\left( J\right) }\sigma
\right) }{\left\vert J\right\vert^\frac{1}{n} }\left\Vert \bigtriangleup _{J}^{\omega ,
\mathbf{b}^{\ast }}x\right\Vert _{L^{2}\left( \omega \right) }^{\spadesuit
}\left\Vert \square _{J}^{\omega ,\mathbf{b}^{\ast }}g\right\Vert
_{L^{2}\left( \omega \right) }^{\bigstar } \\
&=&
\sum_{m\geq 0}\sum_{J\in \Pi _{2}\mathcal{Q}_{m}}\frac{\mathrm{P}^{\alpha
}\left( J,\left\vert \varphi _{J}^{\mathcal{Q}_{m}}\right\vert \mathbf{1}
_{A\backslash I_{\mathcal{Q}_{m}}\left( J\right) }\sigma \right) }{\left\vert
J\right\vert^\frac{1}{n} }\left\Vert \bigtriangleup _{J}^{\omega ,\mathbf{b}^{\ast
}}x\right\Vert _{L^{2}\left( \omega \right) }^{\spadesuit }\left\Vert
\square _{J}^{\omega ,\mathbf{b}^{\ast }}g\right\Vert _{L^{2}\left( \omega
\right) }^{\bigstar }\\
&=&
\sum_{m\geq 0}\left\vert \mathsf{B}\right\vert _{
{stop},\bigtriangleup ^{\omega }}^{A,\mathcal{Q}_{m}}\left(
f,g\right) ,
\end{eqnarray*}
and we can continue with the definition of $\widehat{\mathfrak{N}}_{{
stop},\bigtriangleup ^{\omega }}^{A,\mathcal{Q}_{m}}$ and Cauchy-Schwarz to
obtain
\begin{eqnarray*}
\left\vert \mathsf{B}\right\vert _{{stop},\bigtriangleup ^{\omega
}}^{A,\mathcal{Q}}\left( f,g\right) 
&\!\!\!\!\!\leq &
\sum_{m\geq 0}\widehat{\mathfrak{N}}_{{stop},\bigtriangleup ^{\omega }}^{A,\mathcal{Q}_{m}}\left\Vert 
\mathsf{P}_{\pi \left( \Pi _{1}\mathcal{Q}_{m}\right) }^{\sigma ,\mathbf{b}
}f\right\Vert _{L^{2}\left( \sigma \right) }^{\bigstar }\left\Vert \mathsf{P}
_{\Pi _{2}\mathcal{Q}_{m}}^{\omega ,\mathbf{b}^{\ast }}g\right\Vert
_{L^{2}\left( \omega \right) }^{\bigstar } \\
&\!\!\!\!\!\leq &
\left( \sup_{m\geq 0}\widehat{\mathfrak{N}}_{{stop}
,\bigtriangleup ^{\omega }}^{A,\mathcal{Q}_{m}}\right) \sqrt{\sum_{m\geq
0}\left\Vert \mathsf{P}_{\pi \left( \Pi _{1}\mathcal{Q}_{m}\right) }^{\sigma
,\mathbf{b}}f\right\Vert _{L^{2}\left( \sigma \right) }^{\bigstar 2}}\sqrt{
\sum_{m\geq 0}\left\Vert \mathsf{P}_{\Pi _{2}\mathcal{Q}_{m}}^{\omega ,
\mathbf{b}^{\ast }}g\right\Vert _{L^{2}\left( \omega \right) }^{\bigstar 2}}
\\
&\!\!\!\!\!\leq &
\left( \sup_{m\geq 0}\widehat{\mathfrak{N}}_{{stop}
,\bigtriangleup ^{\omega }}^{A,\mathcal{Q}_{m}}\right) \sqrt{\left\Vert 
\mathsf{P}_{\pi \left( \Pi _{1}\mathcal{Q}\right) }^{\sigma ,\mathbf{b}
}f\right\Vert _{L^{2}\left( \sigma \right) }^{\bigstar 2}}\sqrt{\left\Vert 
\mathsf{P}_{\Pi _{2}\mathcal{Q}}^{\omega ,\mathbf{b}^{\ast }}g\right\Vert
_{L^{2}\left( \omega \right) }^{\bigstar 2}}.
\end{eqnarray*}
\end{proof}

Now we turn to proving inequality (\ref{First inequality}) for the sublinear
form $\left\vert \mathsf{B}\right\vert _{{stop},\bigtriangleup
^{\omega }}^{A,\mathcal{P}}\left( f,g\right) $, i.e.
\begin{eqnarray*}
\left\vert \mathsf{B}\right\vert _{{stop},\bigtriangleup ^{\omega
}}^{A,\mathcal{P}}\left( f,g\right) &\equiv &\sum_{J\in \Pi _{2}\mathcal{P}}
\frac{\mathrm{P}^{\alpha }\left( J,\left\vert \varphi _{J}\right\vert 
\mathbf{1}_{A\backslash I_{\mathcal{P}}\left( J\right) }\sigma \right) }{
\left\vert J\right\vert }\left\Vert \bigtriangleup _{J}^{\omega ,\mathbf{b}
^{\ast }}x\right\Vert _{L^{2}\left( \omega \right) }^{\spadesuit }\left\Vert
\square _{J}^{\omega ,\mathbf{b}^{\ast }}g\right\Vert _{L^{2}\left( \omega
\right) }^{\bigstar } \\
&\lesssim &\left( \mathcal{E}_{2}^{\alpha }+\sqrt{\mathfrak{A}_{2}^{\alpha }}
\right) \left\Vert \mathsf{P}_{\pi \left( \Pi _{1}\mathcal{P}\right)
}^{\sigma ,\mathbf{b}}f\right\Vert _{L^{2}\left( \sigma \right) }^{\bigstar
}\left\Vert \mathsf{P}_{\Pi _{2}\mathcal{P}}^{\omega ,\mathbf{b}^{\ast
}}g\right\Vert _{L^{2}\left( \omega \right) }^{\bigstar };
\end{eqnarray*}
\begin{eqnarray*}
\text{where }\varphi _{J} 
&\equiv &
\sum_{I\in \mathcal{C}_{A}^{{'}}:\ \left( I,J\right) \in \mathcal{P}}\left( E_{I}^{\sigma }
\widehat{\square }_{\pi I}^{\sigma ,\flat ,\mathbf{b}}f\right) b_{A}\ 
\mathbf{1}_{A\backslash I}\ \text{\ is supported in }A\backslash I_{\mathcal{P}
}\left( J\right) 
\end{eqnarray*}
and $I_{\mathcal{P}}\left( J\right) $ denotes the smallest cube $I\in 
\mathcal{D}$ for which $\left( I,J\right) \in \mathcal{P}$. We recall the
stopping energy from (\ref{def stopping energy 3}),
\begin{equation*}
\mathbf{X}_{\alpha }\left( \mathcal{C}_{A}\right) ^{2}\equiv \sup_{I\in 
\mathcal{C}_{A}}\frac{1}{\left\vert I\right\vert _{\sigma }}\sup_{I\supset 
\overset{\cdot }{\cup }J_{r}}\sum_{r=1}^{\infty }\left( \frac{\mathrm{P}
^{\alpha }\left( J_{r},\mathbf{1}_{A}\sigma \right) }{\left\vert
J_{r}\right\vert }\right) ^{2}\left\Vert x-m_{J_{r}}\right\Vert
_{L^{2}\left( \mathbf{1}_{J_{r}}\omega \right) }^{2}\ ,
\end{equation*}
where the cubes $J_{r}\in \mathcal{G}$ are pairwise disjoint in $I$.

What now follows is an adaptation to our sublinear form $\left\vert \mathsf{B
}\right\vert _{{stop},\bigtriangleup ^{\omega }}^{A,\mathcal{P}}$ of the
arguments of M. Lacey in \cite{Lac}, together with an additional `indented'
corona construction. We have the following Poisson inequality for cubes $
B\subset A\subset I$:
\begin{eqnarray}
\frac{\mathrm{P}^{\alpha }\left( A,\mathbf{1}_{I\backslash A}\sigma \right) }{
\left\vert A\right\vert^\frac{1}{n} } 
&\approx &
\int_{I\backslash A}\frac{1}{\left(
\left\vert y-c_{A}\right\vert \right) ^{n+1-\alpha }}d\sigma \left( y\right)
\label{BAI} \\
&\lesssim &
\int_{I\backslash A}\frac{1}{\left( \left\vert y-c_{B}\right\vert
\right) ^{n+1-\alpha }}d\sigma \left( y\right) \approx \frac{\mathrm{P}
^{\alpha }\left( B,\mathbf{1}_{I\backslash A}\sigma \right) }{\left\vert
B\right\vert^\frac{1}{n} }  \notag
\end{eqnarray}
where the implied constants depend on $n,\alpha$.

Fix $A\in \mathcal{A}$. Following \cite{Lac} we will use a `decoupled'
modification of the stopping energy $\mathbf{X}_{\alpha }\left( \mathcal{C}
_{A}\right) $ to define a `size functional' of an $A$-admissible collection $
\mathcal{P}$. So suppose that $\mathcal{P}$ is an $A$-admissible collection
of pairs of cubes, and recall that $\Pi _{1}\mathcal{P}$ and $\Pi _{2}
\mathcal{P}$ denote the cubes in the first and second components of the
pairs in $\mathcal{P}$ respectively.

\begin{dfn}
\label{rest K}For an $A$-admissible collection of pairs of cubes $
\mathcal{P}$, and a cube $K\in \Pi _{1}\mathcal{P}$, define the
projection of $\mathcal{P}$ `relative to $K$' by 
\begin{equation*}
\Pi _{2}^{K}\mathcal{P}\equiv \left\{ J\in \Pi _{2}\mathcal{P}:\ J^{\maltese
}\subset K\right\} ,
\end{equation*}
where we have suppressed dependence on $A$.
\end{dfn}

\begin{dfn}
\label{Pi below}We will use as the `size testing collection' of cubes
for $\mathcal{P}$ the collection 
\begin{equation*}
\Pi _{1}^{{below}}\mathcal{P}\equiv \left\{ K\in \mathcal{D}
:K\subset I\text{ for some }I\in \Pi _{1}\mathcal{P}\right\} ,
\end{equation*}
which consists of all cubes contained in a cube from $\Pi _{1}
\mathcal{P}$.
\end{dfn}

Continuing to follow Lacey \cite{Lac}, we define two `size functionals' of $
\mathcal{P}$ as follows. Recall 
that for a pseudoprojection $\mathsf{Q}_{\mathcal{H}
}^{\omega }$ on $x$ we have 
\begin{eqnarray*}
\left\Vert \mathsf{Q}_{\mathcal{H}}^{\omega ,\mathbf{b}^{\ast }}x\right\Vert
_{L^{2}\left( \omega \right) }^{\spadesuit 2}
&=&
\sum_{J\in \mathcal{H}
}\left\Vert \bigtriangleup _{J}^{\omega ,\mathbf{b}^{\ast }}x\right\Vert
_{L^{2}\left( \omega \right) }^{\spadesuit 2}\\
&=&\sum_{J\in \mathcal{H}}\left(
\left\Vert \bigtriangleup _{J}^{\omega ,\mathbf{b}^{\ast }}x\right\Vert
_{L^{2}\left( \omega \right) }^{2}+\inf_{z\in \mathbb{R}^n}\sum_{J^{\prime
}\in \mathfrak{C}_{{brok}}\left( J\right) }\left\vert J^{\prime
}\right\vert _{\omega }\left( E_{J^{\prime }}^{\omega }\left\vert
x-z\right\vert \right) ^{2}\right)  \label{h}
\end{eqnarray*}

\begin{dfn}
\label{def ext size}If $\mathcal{P}$ is $A$-admissible, define an \emph{
initial} size condition $\mathcal{S}_{{init}{size}}^{\alpha
,A}\left( \mathcal{P}\right) $ by
\begin{equation}
\mathcal{S}_{{init}{size}}^{\alpha ,A}\left( \mathcal{P}
\right) ^{2}\equiv \sup_{K\in \Pi _{1}^{{below}}\mathcal{P}}\frac{1}{
\left\vert K\right\vert _{\sigma }}\left( \frac{\mathrm{P}^{\alpha }\left( K,
\mathbf{1}_{A\backslash K}\sigma \right) }{\left\vert K\right\vert^\frac{1}{n} }\right)
^{2}\left\Vert \mathsf{Q}_{\Pi _{2}^{K}\mathcal{P}}^{\omega ,\mathbf{b}
^{\ast }}x\right\Vert _{L^{2}\left( \omega \right) }^{\spadesuit 2}.
\label{def P stop energy 3}
\end{equation}
\end{dfn}
The following key fact is essential.

\textbf{Key Fact \#1:}
\begin{equation}
\text{If }K\subset A\text{ and }K\notin \mathcal{C}_{A}, \text{ then } \Pi _{2}^{K}
\mathcal{P}=\emptyset \ .  \label{later use}
\end{equation}
To see this, suppose that $K\subset A$ and $K\notin \mathcal{C}_{A}$. Then $
K\subset A^{\prime }$ for some $A^{\prime }\in \mathfrak{C}_{\mathcal{A}
}\left( A\right) $, and so if there is $J^{\prime }\in \Pi _{2}^{K}\mathcal{P
}$, then $\left( J^{\prime }\right) ^{\maltese }\subset K\subset A^{\prime }$
, which implies that $J^{\prime }\notin \mathcal{C
}_{A}^{\mathcal{G},{shift}}$, which contradicts $\Pi _{2}^{K}
\mathcal{P}\subset \mathcal{C}_{A}^{\mathcal{G},{shift}}$. We now
observe from (\ref{later use}) that we may also write the initial size
functional as
\begin{equation}
\mathcal{S}_{{init}{size}}^{\alpha ,A}\left( \mathcal{P}
\right) ^{2}\equiv \sup_{K\in \Pi _{1}^{{below}}\mathcal{P}\cap
\mathcal{C}_{A}'}\frac{1}{\left\vert K\right\vert _{\sigma }}
\left( \frac{\mathrm{P}^{\alpha }\left( K,\mathbf{1}_{A\backslash K}\sigma
\right) }{\left\vert K\right\vert^\frac{1}{n} }\right) ^{2}\left\Vert \mathsf{Q}_{\Pi
_{2}^{K}\mathcal{P}}^{\omega ,\mathbf{b}^{\ast }}x\right\Vert _{L^{2}\left(
\omega \right) }^{\spadesuit 2}.  \label{rewrite size}
\end{equation}

However, we will also need to control certain pairs $\left( I,J\right) \in 
\mathcal{P}$ using testing cubes $K$ which are strictly smaller than $
J^{\maltese }$, namely those $K\in \mathcal{C}_{A}$ such that $K\subset
J^{\maltese }\subset \pi _{\mathcal{D}}^{\left( 2\right) }K$. For this, we
need a second key fact regarding the cubes $J^{\maltese }$, that will
also play a crucial role in controlling pairs in the indented corona below,
and which is that $J$ is always contained in one of the \emph{inner} $2^n$
grandchildren of $J^{\maltese }$. For $M\in \mathcal{D}$, denote by $M_\searrow$ and $M_\nearrow$ any of the inner and outer respectively grandchildren of $M$
and by $M_{J}$ and $M^{\flat}$ the child and grandchild respectively of $M$
that contains $J$, provided they exist.

\textbf{Key Fact \#2:}
\begin{eqnarray}
&&3J\subset J^{\flat } \text{ and $J^\flat$ is an inner grandchild of $J^\maltese$}\label{indentation} 
\end{eqnarray}

To see this, suppose that the child $
J_{J}^{\maltese }$ of $J^{\maltese }$ contains $J$ ($J_{J}^{\maltese }$ exists because $J$ is good in $J^{\maltese }$).
Then observe that $J$ is by definition $\varepsilon -{bad}$ in $J_{J}^{\maltese }$, i.e. 
$$
{\dist}\left( J,{\body}J_{J}^{\maltese }\right) \leq 2\left\vert J\right\vert ^{\frac{\varepsilon}{n}
}\left\vert J_{J}^{\maltese }\right\vert ^{\frac{1-\varepsilon}{n} }
$$
and so cannot lie in any of the $4^n-2^n$ outermost grandchildren $J_{\nearrow}^{\maltese }$. Indeed, if $J\subset
J_{\nearrow}^{\maltese }$, then
\begin{eqnarray*}
\dist\left( J,{\body}J^{\maltese }\right) 
&=&
\dist\left( J,{\body}J_{J}^{\maltese }\right) \leq 
2\left\vert J\right\vert ^{\frac{\varepsilon}{n}
}\left\vert J_{J}^{\maltese }\right\vert ^{\frac{1-\varepsilon}{n} }
\\
&=&
2^{\varepsilon }\left\vert J\right\vert ^{\frac{\varepsilon}{n}
}\left\vert J^{\maltese }\right\vert ^{\frac{1-\varepsilon}{n} }
<
2\left\vert J\right\vert ^{\frac{\varepsilon}{n}
}\left\vert J^{\maltese }\right\vert ^{\frac{1-\varepsilon}{n} }
\end{eqnarray*}
contradicting the fact that $J$ is $\varepsilon -{good}$ in $
J^{\maltese }$. Thus we must have $J\subset J^\flat$, and of course we get that $J^\flat$ is an inner grandchild of $J^\maltese$, (where the
body of $J^{\maltese }$ does not intersect the interior of $
J^\flat$, thus permitting $J$ to be $\varepsilon -{good}$
in $J^{\maltese }$). Finally, the fact that $J$ is $\varepsilon -{
good}$ in $J^{\maltese }$ implies that $3J\subset J^\flat$.

This second key fact is what underlies the construction of the indented
corona below, and motivates the next definition of augmented projection, in
which we allow cubes $K$ satisfying $J\subset K\subsetneqq J^{\maltese
}\subset \pi _{\mathcal{D}}^{\left( 2\right) }K$, as well as $K\in C_{A}$,
to be tested over in the augmented size condition below.

\begin{dfn}
\label{augs}Suppose $\mathcal{P}$ is an $A$-admissible collection.

\begin{enumerate}[(1).]
\item For $K\in \Pi _{1}\mathcal{P}$, define the \emph{augmented} projection
of $\mathcal{P}$ relative to $K$ by
\begin{equation*}
\Pi _{2}^{K,{aug}}\mathcal{P}\equiv \left\{ J\in \Pi _{2}\mathcal{P}
:J\subset K\text{ and }J^{\maltese }\subset \pi _{\mathcal{D}}^{\left(
2\right) }K\right\}.
\end{equation*}

\item Define the corresponding \emph{augmented} size functional $\mathcal{S}
_{{aug}{size}}^{\alpha ,A}\left( \mathcal{P}\right) $ by 
\begin{equation*}
\mathcal{S}_{{aug}{size}}^{\alpha ,A}\left( \mathcal{P}
\right) ^{2}
\equiv\!\!\!\!\!\!\!
\sup_{K\in \Pi _{1}^{{below}}\mathcal{P}\cap
\mathcal{C}_{A}'}\frac{1}{\left\vert K\right\vert _{\sigma }}
\left( \frac{\mathrm{P}^{\alpha }\left( K,\mathbf{1}_{A\backslash K}\sigma
\right) }{\left\vert K\right\vert }\right) ^{\!\!\!2}
\!\left\Vert \mathsf{Q}_{\Pi
_{2}^{K,{aug}}\mathcal{P}}^{\omega ,\mathbf{b}^{\ast }}x\right\Vert
_{L^{2}\left( \omega \right) }^{\spadesuit 2}\
\end{equation*}
\end{enumerate}
\end{dfn}

We note that the augmented projection $\Pi _{2}^{K,{aug}}\mathcal{P}$
includes cubes $J$ for which $J\subset K\subsetneqq J^{\maltese }\subset
\pi _{\mathcal{D}}^{\left( 2\right) }K$, and hence $J$ need not be $
\varepsilon -{good}$ inside $K$.  Then by the second key fact (\ref
{indentation}), and using that the boundaries of  $J_{\searrow}^{\maltese }$
 lie in the ${body}$ of $J^{\maltese }$, we
have two consequences,
\begin{equation*}
K\in \left\{ J_{J}^{\maltese },J^\flat\right\} \text{ and 
}3J\subset J^\flat\subset 3J^\flat\subset J^{\maltese }\text{ for }J\in \Pi _{2}^{K,{aug}}\mathcal{P},
\end{equation*}
which will play an important role below.

The augmented size functional $\mathcal{S}_{{aug}{size}
}^{\alpha ,A}\left( \mathcal{P}\right) $ is a `decoupled' form of the
stopping energy $\mathbf{X}_{\alpha }\left( \mathcal{C}_{A}\right) $
restricted to $\mathcal{P}$, in which the cubes $J$ appearing in $
\mathbf{X}_{\alpha }\left( \mathcal{C}_{A}\right) $ no longer appear in the
Poisson integral in $\mathcal{S}_{{aug}{size}}^{\alpha
,A}\left( \mathcal{P}\right) $, and it plays a crucial role in Lacey's
argument in \cite{Lac}. We note two essential properties of this definition
of size functional:

\begin{enumerate}
\item \textbf{Monotonicity} of size: $\mathcal{S}_{{aug}{size
}}^{\alpha ,A}\left( \mathcal{P}\right) \leq \mathcal{S}_{{aug}
{size}}^{\alpha ,A}\left( \mathcal{Q}\right) $ if $\mathcal{P}
\subset \mathcal{Q}$,

\item \textbf{Control} by energy and Muckenhoupt conditions: $\mathcal{S}_{
{aug}{size}}^{\alpha ,A}\left( \mathcal{P}\right) \lesssim 
\mathcal{E}_{2}^{\alpha }+\sqrt{\mathfrak{A}_{2}^{\alpha }}$.
\end{enumerate}

The monotonicity property follows from $\Pi _{1}^{{below}}\mathcal{P}
\subset \Pi _{1}^{{below}}\mathcal{Q}$ and $\Pi _{2}^{K}\mathcal{P}
\subset \Pi _{2}^{K}\mathcal{Q}$. The control property is contained in the
next lemma, which uses the stopping energy control for the form $\mathsf{B}_{
{stop}}^{A}\left( f,g\right) $ associated with $A$.

\begin{lem}
\label{energy control}If $\mathcal{P}^{A}$ is as in (\ref{initial P}) and $
\mathcal{P}\subset \mathcal{P}^{A}$, then 
\begin{equation*}
\mathcal{S}_{{aug}{size}}^{\alpha ,A}\left( \mathcal{P}
\right) \leq \mathbf{X}_{\alpha }\left( \mathcal{C}_{A}\right) \lesssim 
\mathcal{E}_{2}^{\alpha }+\sqrt{\mathfrak{A}_{2}^{\alpha }}.
\end{equation*}
\end{lem}

\begin{proof}
We have
\begin{eqnarray*}
\mathcal{S}_{{aug}{size}}^{\alpha ,A}\left( \mathcal{P}
\right) ^{2} 
&=&\!\!\!\!\!\!\!\!\!\!
\sup_{K\in \Pi _{1}^{{below}}\mathcal{P}\cap 
\mathcal{C}_{A}'}\frac{1}{\left\vert K\right\vert
_{\sigma }}\left( \frac{\mathrm{P}^{\alpha }\left( K,\mathbf{1}_{A\backslash
K}\sigma \right) }{\left\vert K\right\vert^\frac{1}{n} }\right) ^{2}\left\Vert \mathsf{Q}
_{\Pi _{2}^{K}\mathcal{P}\cup \Pi _{2}^{K,{aug}}\mathcal{P}}^{\omega
,\mathbf{b}^{\ast }}x\right\Vert _{L^{2}\left( \omega \right) }^{\spadesuit
2} \\
&\lesssim &
\sup_{K\in \mathcal{C}_{A}'}\frac{1}{
\left\vert K\right\vert _{\sigma }}\left( \frac{\mathrm{P}^{\alpha }\left( K,
\mathbf{1}_{A}\sigma \right) }{\left\vert K\right\vert^\frac{1}{n} }\right)
^{2}\left\Vert x-m_{K}\right\Vert _{L^{2}\left( \mathbf{1}_{K}\omega \right)
}^{2}\leq \mathbf{X}_{\alpha }\left( \mathcal{C}_{A}\right) ^{2},
\end{eqnarray*}
which is the first inequality in the statement of the lemma. The second
inequality follows from (\ref{def stopping bounds 3}).
\end{proof}

There is an important special circumstance, introduced by M. Lacey in \cite
{Lac}, in which we can bound our forms by the size functional, namely when
the pairs all straddle a subpartition of $A$, and we present this in the
next subsection. In order to handle the fact that the cubes in $\Pi
_{1}^{{below}}\mathcal{P}\cap \mathcal{C}_{A}'$ need no
longer enjoy any goodness, we will need to formulate a Substraddling Lemma
to deal with this situation as well. See \textbf{Remark on lack of usual
goodness} after (\ref{N_L}), where it is explained how this applies to the
proof of (\ref{rem}). Then in the following subsection, we use the bottom/up
stopping time construction of M. Lacey, together with an additional
`indented' top/down corona construction, to reduce control of the sublinear
stopping form $\left\vert \mathsf{B}\right\vert _{{stop}
,\bigtriangleup ^{\omega }}^{A,\mathcal{P}}\left( f,g\right) $ in inequality
(\ref{First inequality}) to the three special cases addressed by the
Orthogonality Lemma, the Straddling Lemma and the Substraddling Lemma.

\subsection{$\flat$Straddling, Substraddling, Corona-Straddling Lemmas}

We begin with the Corona-straddling Lemma in which the straddling collection
is the set of $\mathcal{A}$-children of $A$, and applies to the `corona
straddling' subcollection of the initial admissible collection $\mathcal{P}
^{A}$ - see (\ref{initial P}). Define the `corona straddling' collection $
\mathcal{P}_{{cor}}^{A}$ by
\begin{equation}
\mathcal{P}_{{cor}}^{A}\equiv \bigcup\limits_{A^{\prime }\in \mathfrak{
C}_{\mathcal{A}}\left( A\right) }\left\{ \left( I,J\right) \in \mathcal{P}
^{A}:J\subset A^{\prime }\varsubsetneqq J^{\maltese }\subset \pi _{\mathcal{D
}}^{\left( 2\right) }A^{\prime }\right\} .  \label{def cor}
\end{equation}
Note that $\mathcal{P}_{{cor}}^{A}$ is an $A$-admissible collection
that consists of just those pairs $\left( I,J\right) $ for which $
J^{\maltese }$ is either the $\mathcal{D}$-parent or the $\mathcal{D}$
-grandparent of a stopping cube $A^{\prime }\in \mathfrak{C}_{\mathcal{A}
}\left( A\right) $. The bound for the norm of the corresponding form is
controlled by the energy condition.

\begin{lem}
\label{cor strad 1}We have the sublinear form bound
\begin{equation*}
\mathfrak{N}_{{stop},\bigtriangleup ^{\omega }}^{A,\mathcal{P}_{
{cor}}^{A}}\leq C\mathcal{E}_{2}^{\alpha }.
\end{equation*}
\end{lem}

\begin{proof}
The key point here is our assumption that $J\subset A^{\prime
}\varsubsetneqq J^{\maltese }\subset \pi _{\mathcal{D}}^{\left( 2\right)
}A^{\prime }$ for $\left( I,J\right) \in \mathcal{P}_{{cor}}^{A}$,
which implies that in fact $3J\subset A^{\prime }$ since $J\cap \body\left( \pi _{\mathcal{D}}^{\left( 2\right) }A^{\prime }\right) =\emptyset $
because $J$ is $\varepsilon -{good}$ in $\pi _{\mathcal{D}}^{\left(
2\right) }A^{\prime }$. We start with
\begin{eqnarray*}
\left\vert \mathsf{B}\right\vert _{{stop},\bigtriangleup ^{\omega
}}^{A,\mathcal{P}_{{cor}}^{A}}\left( f,g\right) 
\!\!\!&=&\!\!\!\!\!\!\!\!\!
\sum_{J\in \Pi _{2}
\mathcal{P}_{{cor}}^{A}}\!\!\!\!\!\!
\frac{\mathrm{P}^{\alpha }\left( J,\left\vert
\varphi _{J}^{\mathcal{P}_{{cor}}^{A}}\right\vert \mathbf{1}
_{A\backslash I_{\mathcal{P}_{{cor}}^{A}}\left( J\right) }\sigma \right) 
}{\left\vert J\right\vert }\left\Vert \bigtriangleup _{J}^{\omega ,\mathbf{b}
^{\ast }}x\right\Vert _{L^{2}\left( \omega \right) }^{\spadesuit }\left\Vert
\square _{J}^{\omega ,\mathbf{b}^{\ast }}g\right\Vert _{L^{2}\left( \omega
\right) }^{\bigstar }
\end{eqnarray*}
\begin{equation*}
\quad=
\sum_{A^{\prime }\in \mathfrak{C}_{\mathcal{A}}\left( A\right)
}\!\!\!\!\sum_{\,\,\,\,\,\substack{J\in \Pi _{2}\mathcal{P}_{{cor}}^{A}\\ \ 3J\subset A^{\prime }}}\!\!\!\!\!\!\!
\frac{\mathrm{P}^{\alpha }\left( J,\left\vert \varphi _{J}^{\mathcal{P}_{
{cor}}^{A}}\right\vert \mathbf{1}_{A\backslash I_{\mathcal{P}_{{cor}
}^{A}}\left( J\right) }\sigma \right) }{\left\vert J\right\vert }\left\Vert
\bigtriangleup _{J}^{\omega ,\mathbf{b}^{\ast }}x\right\Vert _{L^{2}\left(
\omega \right) }^{\spadesuit }\left\Vert \square _{J}^{\omega ,\mathbf{b}
^{\ast }}g\right\Vert _{L^{2}\left( \omega \right) }^{\bigstar } 
\end{equation*}
where 
$$
\varphi _{J}^{\mathcal{P}_{{cor}}^{A}} \equiv \sum_{I\in
\Pi _{1}\mathcal{P}_{{cor}}^{A}:\mathcal{\ }\left( I,J\right) \in 
\mathcal{P}_{{cor}}^{A}}b_{A}E_{I}^{\sigma }\left( \widehat{\square }
_{\pi I}^{\sigma ,\flat ,\mathbf{b}}f\right) \ \mathbf{1}_{A\backslash I}\ .
$$
If $J\in \Pi _{2}\mathcal{P}_{{cor}}^{A}$ and $J\subset A^{\prime }\in 
\mathfrak{C}_{\mathcal{A}}\left( A\right) $, then either $A^{\prime
}=J^\flat$ or $A^{\prime
}=J_J^{\maltese }$ and we have 
\begin{equation*}
\frac{\mathrm{P}^{\alpha }\left( J,\mathbf{1}_{A\backslash I_{\mathcal{P}_{
{cor}}^{A}}\left( J\right) }\sigma \right) }{\left\vert J\right\vert^\frac{1}{n} }
\approx \left\{ 
\begin{array}{ccc}
\frac{\mathrm{P}^{\alpha }\left( A^{\prime },\mathbf{1}_{A\backslash I_{
\mathcal{P}_{{cor}}^{A}}}\sigma \right) }{\left\vert A^{\prime
}\right\vert^\frac{1}{n} }
\leq 
\frac{\mathrm{P}^{\alpha }\left( A^{\prime },\mathbf{1}
_{A}\sigma \right) }{\left\vert A^{\prime }\right\vert^\frac{1}{n} } & \text{if} & 
A^{\prime }=J^{\flat } \\ 
\frac{\mathrm{P}^{\alpha }\left( A_{J}^{\prime },\mathbf{1}_{A\backslash I_{
\mathcal{P}_{{cor}}^{A}}}\sigma \right) }{\left\vert A_{J}^{\prime
}\right\vert^\frac{1}{n} }
\lesssim 
\frac{\mathrm{P}^{\alpha }\left( A^{\prime },\mathbf{1
}_{A}\sigma \right) }{\left\vert A^{\prime }\right\vert^\frac{1}{n} } & \text{if} & 
A^{\prime }=J_J^{\maltese}
\end{array}
\right.
\end{equation*}
Since $\left\vert \varphi _{J}^{\mathcal{P}_{{cor}}^{A}}\right\vert
\lesssim \alpha _{\mathcal{A}}\left( A\right) \mathbf{1}_{A}$ by (\ref{phi
bound}), we can then bound $\left\vert \mathsf{B}\right\vert _{{stop}
,\bigtriangleup ^{\omega }}^{A,\mathcal{P}_{{cor}}^{A}}\left(
f,g\right) $ by
\begin{eqnarray*}
&&\alpha _{\mathcal{A}}\left( A\right) \sum_{A^{\prime }\in \mathfrak{C}_{
\mathcal{A}}\left( A\right) }\left( \frac{\mathrm{P}^{\alpha }\left(
A^{\prime },\mathbf{1}_{A}\sigma \right) }{\left\vert A^{\prime }\right\vert^\frac{1}{n} 
}\right) \left\Vert \mathsf{Q}_{\Pi _{2}\mathcal{P}_{{cor}
}^{A};A^{\prime }}^{\omega ,\mathbf{b}^{\ast }}x\right\Vert _{L^{2}\left(
\omega \right) }^{\spadesuit }\left\Vert \mathsf{P}_{\Pi _{2}\mathcal{P}_{
{cor}}^{A};A^{\prime }}^{\omega ,\mathbf{b}^{\ast }}g\right\Vert
_{L^{2}\left( \omega \right) }^{\bigstar } \\
&\leq &
\alpha _{\mathcal{A}}\left( A\right) \left( \sum_{A^{\prime }\in 
\mathfrak{C}_{\mathcal{A}}\left( A\right) }\left( \frac{\mathrm{P}^{\alpha
}\left( A^{\prime },\mathbf{1}_{A}\sigma \right) }{\left\vert A^{\prime
}\right\vert^\frac{1}{n} }\right) ^{2}\left\Vert x-m_{A^{\prime }}^{\sigma }\right\Vert
_{L^{2}\left( \mathbf{1}_{A^{\prime }}\sigma \right) }^{ 2}\right)
^{\frac{1}{2}}\cdot \\
&&\hspace{5cm}\cdot\left( \sum_{A^{\prime }\in \mathfrak{C}_{\mathcal{A}
}\left( A\right) }\left\Vert \mathsf{P}_{\Pi _{2}\mathcal{P}_{{cor}
}^{A};A^{\prime }}^{\omega ,\mathbf{b}^{\ast }}g\right\Vert _{L^{2}\left(
\omega \right) }^{\bigstar 2}\right) ^{\frac{1}{2}} \\
&\leq &
\mathcal{E}_{2}^{\alpha }\alpha _{\mathcal{A}}\left( A\right) \sqrt{
\left\vert A\right\vert _{\sigma }}\left\Vert \mathsf{P}_{\Pi _{2}\mathcal{P}
_{{cor}}^{A}}^{\omega ,\mathbf{b}^{\ast }}g\right\Vert _{L^{2}\left(
\omega \right) }^{\bigstar }\\
&\leq &
\mathcal{E}_{2}^{\alpha }\alpha _{\mathcal{A
}}\left( A\right) \sqrt{\left\vert A\right\vert _{\sigma }}\left\Vert 
\mathsf{P}_{\mathcal{C}_{A}^{{shift}}}^{\omega ,\mathbf{b}^{\ast
}}g\right\Vert _{L^{2}\left( \omega \right) }^{\bigstar }
\end{eqnarray*}
where in the last line we have used the strong energy constant $\mathcal{E}
_{2}^{\alpha }$ in (\ref{strong b* energy}).
\end{proof}

\begin{dfn}
We say that an admissible collection of pairs $\mathcal{P}$ is \emph{reduced}
if it contains no pairs from $\mathcal{P}_{{cor}}^{A}$, i.e.
\begin{equation*}
\mathcal{P}\cap \mathcal{P}_{{cor}}^{A}=\emptyset .
\end{equation*}
\end{dfn}

Recall that in terms of $J^{\flat }$ we rewrite
\begin{eqnarray*}
\Pi _{2}^{K,{aug}}\mathcal{P}&=&\left\{ J\in \Pi _{2}\mathcal{P}
:J\subset K\text{ and }J^{\maltese }\subset \pi _{\mathcal{D}}^{\left(
2\right) }K\right\} \\
&=&\left\{ J\in \Pi _{2}\mathcal{P}:J\subset K\text{ and }
J^{\flat }\subset K\right\} 
\end{eqnarray*}

\begin{dfn}
\label{flat straddles}Given a \emph{reduced admissible} collection of pairs $
\mathcal{Q}$ for $A$, and a subpartition $\mathcal{S}\subset \Pi _{1}^{
{below}}\mathcal{Q}\cap \mathcal{C}_{A}'$ of
pairwise disjoint cubes in $A$, we say that $\mathcal{Q}$ $\flat$
\textbf{straddles} $\mathcal{S}$ if for every pair $\left( I,J\right) \in 
\mathcal{Q}$ there is $S\in \mathcal{S}\cap \left[ J,I\right] $ with $
J^{\flat }\subset S$. To avoid trivialities, we further assume that for
every $S\in \mathcal{S}$, there is at least one pair $\left( I,J\right) \in 
\mathcal{Q}$ with $J^{\flat }\subset S\subset I$. Here $\left[ J,I\right] $
denotes the geodesic in the dyadic tree $\mathcal{D}$ that connects $J^{
\mathcal{D}}$ to $I$, where $J^{\mathcal{D}}$ is the minimal cube in $
\mathcal{D}$ that contains $J$.
\end{dfn}

\begin{dfn}
\label{def Whit}For any dyadic cube $S\in \mathcal{D}$, define the
Whitney collection $\mathcal{W}\left( S\right) $ to consist of the maximal
subcubes $K$ of $S$ whose triples $3K$ are contained in $S$. Then set $
\mathcal{W}^{\ast }\left( S\right) \equiv \mathcal{W}\left( S\right) \cup
\left\{ S\right\} $.
\end{dfn}

The following geometric proposition will prove useful in proving the $\flat $
Straddling Lemma \ref{straddle 3 ref} below. For $S\in \mathcal{S}$, let $\mathcal{Q}^{S}\equiv \left\{ \left( I,J\right)
\in \mathcal{Q}:J^{\flat }\subset S\subset I\right\} $.

\begin{prop}
\label{flatness}Suppose $\mathcal{Q}$ is reduced admissible and $\flat $
straddles a subpartition $\mathcal{S}$ of $A$. Fix $S\in \mathcal{S}$.
Define 
\begin{equation*}
\varphi _{J}^{\mathcal{Q}^{S}}\left[ h\right] \equiv \sum_{I\in \Pi _{1}
\mathcal{Q}^{S}:\mathcal{\ }\left( I,J\right) \in \mathcal{Q}
^{S}}b_{A}E_{I}^{\sigma }\left( \widehat{\square }_{\pi I}^{\sigma ,\flat ,
\mathbf{b}}h\right) \ \mathbf{1}_{A\backslash I}\ ,
\end{equation*}
assume that $h\in L^{2}\left( \sigma \right) $ is supported in the cube $
A$, and that there is a cube $H\in \mathcal{C}_{A}$ with $H\supset S$
such that 
\begin{equation*}
E_{I}^{\sigma }\left\vert h\right\vert \leq CE_{H}^{\sigma }\left\vert
h\right\vert ,\ \ \ \ \ \text{for all }I\in \Pi _{1}^{{below}}
\mathcal{Q}\cap \mathcal{C}_{A}'\text{ with }I\supset S.
\end{equation*}
Then
\begin{eqnarray*}
&&\sum_{J\in \Pi _{2}\mathcal{Q}:\ J^{\flat }\subset S}\frac{\mathrm{P}
^{\alpha }\left( J,\left\vert \varphi _{J}^{\mathcal{Q}}\left[ h\right]
\right\vert \mathbf{1}_{A\backslash I_{\mathcal{Q}}\left( J\right) }\sigma
\right) }{\left\vert J\right\vert }\left\Vert \bigtriangleup _{J}^{\omega ,
\mathbf{b}^{\ast }}x\right\Vert _{L^{2}\left( \omega \right) }^{\spadesuit
}\left\Vert \square _{J}^{\omega ,\mathbf{b}^{\ast }}g\right\Vert
_{L^{2}\left( \omega \right) }^{\bigstar } \\
&\lesssim &\alpha _{\mathcal{H}}\left( H\right) \frac{\mathrm{P}^{\alpha
}\left( S,\mathbf{1}_{A\backslash S}\sigma \right) }{\left\vert S\right\vert }
\left\Vert \mathsf{Q}_{\Pi _{2}^{S,{aug}}\mathcal{Q}}^{\omega ,
\mathbf{b}^{\ast }}x\right\Vert _{L^{2}\left( \omega \right) }^{\spadesuit
}\left\Vert \mathsf{P}_{\Pi _{2}^{S,{aug}}\mathcal{Q}}^{\omega ,
\mathbf{b}^{\ast }}g\right\Vert _{L^{2}\left( \omega \right) }^{\bigstar } \\
&&+\alpha _{\mathcal{H}}\left( H\right) \sum\limits_{K\in \mathcal{W}\left(
S\right) }\frac{\mathrm{P}^{\alpha }\left( K,\mathbf{1}_{A\backslash K}\sigma
\right) }{\left\vert K\right\vert }\left\Vert \mathsf{Q}_{\Pi _{2}^{K,
{aug}}\mathcal{Q}}^{\omega ,\mathbf{b}^{\ast }}x\right\Vert
_{L^{2}\left( \omega \right) }^{\spadesuit }\left\Vert \mathsf{P}_{\Pi
_{2}^{K,{aug}}\mathcal{Q}}^{\omega ,\mathbf{b}^{\ast }}g\right\Vert
_{L^{2}\left( \omega \right) }^{\bigstar }\ .
\end{eqnarray*}
The sum over Whitney cubes $K\in \mathcal{W}\left( S\right) $ is only
required to bound the sum of those terms on the left for which $J^{\flat
}\subset S^{\prime \prime }$ for some $S^{\prime \prime }\in \mathfrak{C}_{
\mathcal{D}}^{\left( 2\right) }\left( S\right) $.
\end{prop}

\begin{proof}
Suppose first that $J^{\flat }=S\in \mathcal{C}_{A}'$.
Then $3S=3J^{\flat }\subset J^{\maltese }\subset I_{\mathcal{Q}}\left(
J\right) $ and using (\ref{phi bound}) with $\alpha _{\mathcal{H}}\left(
H\right) $ in place of $\alpha _{\mathcal{A}}\left( A\right) $, we have
\begin{eqnarray*}
\frac{\mathrm{P}^{\alpha }\left( J,\left\vert \varphi _{J}^{\mathcal{Q}
}\right\vert \mathbf{1}_{A\backslash I_{\mathcal{Q}}\left( J\right) }\sigma
\right) }{\left\vert J\right\vert^\frac{1}{n} } 
&\lesssim &
\alpha _{\mathcal{H}}\left(
H\right) \frac{\mathrm{P}^{\alpha }\left( J,\mathbf{1}_{A\backslash
J^{\maltese }}\sigma \right) }{\left\vert J\right\vert^\frac{1}{n} } \\
&\lesssim &
\alpha _{\mathcal{H}}\left( H\right) \frac{\mathrm{P}^{\alpha
}\left( S,\mathbf{1}_{A\backslash J^{\maltese }}\sigma \right) }{\left\vert
S\right\vert^\frac{1}{n} }
\leq
\alpha _{\mathcal{H}}\left( H\right) \frac{\mathrm{P}
^{\alpha }\left( S,\mathbf{1}_{A\backslash S}\sigma \right) }{\left\vert
S\right\vert^\frac{1}{n} }.
\end{eqnarray*}
Suppose next that $J^{\flat }=S^{\prime }\in \mathfrak{C}_{\mathcal{D}
}\left( S\right) $. Then $3S^{\prime }=3J^{\flat }\subset J^{\maltese
}\subset I_{\mathcal{Q}}\left( J\right) $ and (\ref{phi bound}) give
\begin{eqnarray*}
\frac{\mathrm{P}^{\alpha }\left( J,\left\vert \varphi _{J}^{\mathcal{Q}
}\right\vert \mathbf{1}_{A\backslash I_{\mathcal{Q}}\left( J\right) }\sigma
\right) }{\left\vert J\right\vert ^\frac{1}{n}} 
&\lesssim &
\alpha _{\mathcal{H}}\left(
H\right) \frac{\mathrm{P}^{\alpha }\left( J,\mathbf{1}_{A\backslash
J^{\maltese }}\sigma \right) }{\left\vert J\right\vert^\frac{1}{n} } \\
&\lesssim &
\alpha _{\mathcal{H}}\left( H\right) \frac{\mathrm{P}^{\alpha
}\left( S^{\prime },\mathbf{1}_{A\backslash J^{\maltese }}\sigma \right) }{
\left\vert S^{\prime }\right\vert ^\frac{1}{n}} \\
&\leq &
\alpha _{\mathcal{H}}\left( H\right) \frac{\mathrm{P}^{\alpha }\left(
S^{\prime },\mathbf{1}_{A\backslash S}\sigma \right) }{\left\vert S^{\prime
}\right\vert^\frac{1}{n} }
\approx
\alpha _{\mathcal{H}}\left( H\right) \frac{\mathrm{P}^{\alpha }\left( S,\mathbf{1}_{A\backslash S}\sigma \right) }{\left\vert
S\right\vert^\frac{1}{n} }.
\end{eqnarray*}
Thus in these two cases, by Cauchy-Schwarz, the left hand side of our
conclusion is bounded by a multiple of
\begin{equation*}
\alpha _{\mathcal{H}}\left( H\right) \frac{\mathrm{P}^{\alpha }\left( S,
\mathbf{1}_{A\backslash S}\sigma \right) }{\left\vert S\right\vert^\frac{1}{n} }\left(
\sum_{\substack{J\in \Pi _{2}\mathcal{Q}\\ \ J^{\flat }\subset S}}\left\Vert
\bigtriangleup _{J}^{\omega ,\mathbf{b}^{\ast }}x\right\Vert _{L^{2}\left(
\omega \right) }^{\spadesuit 2}\right) ^{\!\!\!\frac{1}{2}} 
\left( \sum_{{\substack{J\in \Pi _{2}\mathcal{Q}\\ \ J^{\flat }\subset S}}}\left\Vert \square _{J}^{\omega ,\mathbf{b}^{\ast }}g\right\Vert _{L^{2}\left( \omega \right) }^{\bigstar
2}\right) ^{\!\!\!\frac{1}{2}} 
\end{equation*}
\begin{equation*}
=
\alpha _{\mathcal{H}}\left( H\right) \frac{\mathrm{P}^{\alpha }\left( S,
\mathbf{1}_{A\backslash S}\sigma \right) }{\left\vert S\right\vert^\frac{1}{n} }
\left\Vert \mathsf{Q}_{\Pi _{2}^{S,{aug}}\mathcal{Q}}^{\omega ,
\mathbf{b}^{\ast }}x\right\Vert _{L^{2}\left( \omega \right) }^{\spadesuit
}\left\Vert \mathsf{P}_{\Pi _{2}^{S,{aug}}\mathcal{Q}}^{\omega ,
\mathbf{b}^{\ast }}g\right\Vert _{L^{2}\left( \omega \right) }^{\bigstar }\ 
\end{equation*}
Finally, suppose that $J^{\flat }\subset S^{\prime \prime }$ for some $
S^{\prime \prime }\in \mathfrak{C}_{\mathcal{D}}^{\left( 2\right) }\left(
S\right) $. Then $J^{\maltese }\subset S$, and Key Fact \#2 in (\ref
{indentation}) shows that $3J^{\flat }\subset J^{\maltese }$, so that $
3J^{\flat }\subset J^{\maltese }\subset S\subset I_{\mathcal{Q}}\left(
J\right) $. Thus we have $J^{\flat }\subset K=K\left[ J\right] $ for some $
K\in \mathcal{W}\left( S\right) $ and so by (\ref{phi bound}) again,
\begin{eqnarray*}
\frac{\mathrm{P}^{\alpha }\left( J,\left\vert \varphi _{J}^{\mathcal{Q}
}\right\vert \mathbf{1}_{A\backslash I_{\mathcal{Q}}\left( J\right) }\sigma
\right) }{\left\vert J\right\vert^\frac{1}{n} } 
&\lesssim &
\alpha _{\mathcal{H}}\left(
H\right) \frac{\mathrm{P}^{\alpha }\left( J,\mathbf{1}_{A\backslash S}\sigma
\right) }{\left\vert J\right\vert^\frac{1}{n} } \\
&\lesssim &
\alpha _{\mathcal{H}}\left( H\right) \frac{\mathrm{P}^{\alpha
}\left( K,\mathbf{1}_{A\backslash S}\sigma \right) }{\left\vert K\right\vert^\frac{1}{n} }
\leq
\alpha _{\mathcal{H}}\left( H\right) \frac{\mathrm{P}^{\alpha }\left( K,
\mathbf{1}_{A\backslash K}\sigma \right) }{\left\vert K\right\vert ^\frac{1}{n}}.
\end{eqnarray*}
Now we apply Cauchy-Schwarz again, but noting that $J^{\flat }\subset K$
this time, to obtain that the left hand side of our conclusion is bounded by
a multiple of
\begin{equation*}
\alpha _{\mathcal{H}}\left( H\right) \sum\limits_{K\in \mathcal{W}\left(
S\right) }\!\!\!\!\frac{\mathrm{P}^{\alpha }\left( K,\mathbf{1}_{A\backslash K}\sigma
\right) }{\left\vert K\right\vert^\frac{1}{n} }\left( \sum_{\substack{J\in \Pi _{2}\mathcal{Q}\\\
J^{\flat }\subset K}}\left\Vert \bigtriangleup _{J}^{\omega ,\mathbf{b}^{\ast
}}x\right\Vert _{L^{2}\left( \omega \right) }^{\spadesuit 2}\right) ^{\!\!\!\frac{1}{2}} 
\left( \sum_{\substack{J\in \Pi _{2}\mathcal{Q}\\\ J^{\flat }\subset K}}\left\Vert
\square _{J}^{\omega ,\mathbf{b}^{\ast }}g\right\Vert _{L^{2}\left( \omega
\right) }^{\bigstar 2}\right) ^{\!\!\!\frac{1}{2}} \end{equation*}
\begin{equation*}
=\alpha _{\mathcal{H}}\left( H\right) \sum\limits_{K\in \mathcal{W}\left(
S\right) }\frac{\mathrm{P}^{\alpha }\left( K,\mathbf{1}_{A\backslash K}\sigma
\right) }{\left\vert K\right\vert^\frac{1}{n} }\left\Vert \mathsf{Q}_{\Pi _{2}^{K,
{aug}}\mathcal{Q}}^{\omega ,\mathbf{b}^{\ast }}x\right\Vert
_{L^{2}\left( \omega \right) }^{\spadesuit }\left\Vert \mathsf{P}_{\Pi
_{2}^{K,{aug}}\mathcal{Q}}^{\omega ,\mathbf{b}^{\ast }}g\right\Vert
_{L^{2}\left( \omega \right) }^{\bigstar }\ .
\end{equation*}
This completes the proof of Proposition \ref{flatness}.
\end{proof}

Recall the family of operators $\left\{ \square _{I}^{\sigma ,\pi ,\mathbf{b}
}\right\} _{I\in \mathcal{C}_{A}^{\mathcal{A}}}$, where for $I\in \mathcal{C}
_{A}^{\mathcal{A}}$, the dual martingale difference $\square _{I}^{\sigma
,\pi ,\mathbf{b}}$ is defined in (\ref{def pi box}), and
satisfies
\begin{equation*}
\square _{I}^{\sigma ,\pi ,\mathbf{b}}f=\left[ \sum_{I^{\prime }\in 
\mathfrak{C}\left( I\right) }\mathbb{F}_{I^{\prime }}^{\sigma ,\pi ,\mathbf{b
}}f\right] -\mathbb{F}_{I}^{\sigma ,\mathbf{b}}f=\sum_{I^{\prime }\in 
\mathfrak{C}\left( I\right) }\mathbb{F}_{I^{\prime }}^{\sigma ,b_{A}}f-
\mathbb{F}_{I}^{\sigma ,b_{A}}f\ .
\end{equation*}
Since $\square _{I}^{\sigma ,\pi ,\mathbf{b}}$ is the transpose of $
\triangle _{I}^{\sigma ,\pi ,\mathbf{b}}$ for $I\in \mathcal{C}_{A}^{
\mathcal{A}}$, the first line of Lemma \ref{b proj} (where the superscript $
\pi $ is suppressed for convenience) shows that $\left\{ \square
_{I}^{\sigma ,\pi ,\mathbf{b}}\right\} _{I\in \mathcal{C}_{A}^{\mathcal{A}}}$
is a family of projections, and the second line of Lemma \ref{b proj} shows
it is an orthogonal family, i.e. 
\begin{equation*}
\square _{I}^{\sigma ,\pi ,\mathbf{b}}\square _{K}^{\sigma ,\pi ,\mathbf{b}
}=\left\{ 
\begin{array}{ccc}
\square _{I}^{\sigma ,\pi ,\mathbf{b}} & \text{ if } & I=K \\ 
0 & \text{ if } & I\not=K
\end{array}
\right. ,\ \ \ \ \ I,K\in \mathcal{C}_{A}^{\mathcal{A}}.
\end{equation*}
The orthogonal projections 
\begin{eqnarray*}
\mathsf{P}_{\pi \left( \Pi _{1}\mathcal{Q}\right) }^{\sigma ,\pi ,\mathbf{b}
} &\equiv &\sum_{I\in \pi \left( \Pi _{1}\mathcal{Q}\right) }\square
_{I}^{\sigma ,\pi ,\mathbf{b}}=\sum_{I\in \Pi _{1}\mathcal{Q}}\square _{\pi
I}^{\sigma ,\pi ,\mathbf{b}}, \\
\text{where }\pi \left( \Pi _{1}\mathcal{Q}\right) &\equiv &\left\{ \pi _{
\mathcal{D}}I:I\in \Pi _{1}\mathcal{Q}\right\} \text{ and }\Pi _{1}\mathcal{Q
}\subset \mathcal{C}_{A}^{\mathcal{A},'}\ ,
\end{eqnarray*}
thus satisfy the equalities
\begin{equation}
\square _{\pi I}^{\sigma ,\pi ,\mathbf{b}}f=\square _{\pi I}^{\sigma ,\pi ,
\mathbf{b}}\mathsf{P}_{\pi \left( \Pi _{1}\mathcal{Q}\right) }^{\sigma ,\pi ,
\mathbf{b}}f\text{ and }\widehat{\square }_{\pi I}^{\sigma ,\pi ,\mathbf{b}
}f=\widehat{\square }_{\pi I}^{\sigma ,\pi ,\mathbf{b}}\mathsf{P}_{\pi
\left( \Pi _{1}\mathcal{Q}\right) }^{\sigma ,\pi ,\mathbf{b}}f  \label{sat}
\end{equation}
 for $I\in \Pi _{1}\mathcal{Q}\subset \mathcal{C}_{A}^{\mathcal{A}
{restrict}}$, which will permit us to apply certain projection tricks used for Haar
projections in the proof of $T1$ theorems.

However, in our sublinear stopping form $\left\vert \mathsf{B}\right\vert _{
{stop},\bigtriangleup ^{\omega }}^{A,\mathcal{Q}}$, the dual
martingale projections in use in the function
\begin{equation}
\varphi _{J}^{\mathcal{Q}^{S}}\equiv \sum_{I\in \Pi _{1}\mathcal{Q}^{S}:
\mathcal{\ }\left( I,J\right) \in \mathcal{Q}^{S}}b_{A}E_{I}^{\sigma }\left( 
\widehat{\square }_{\pi I}^{\sigma ,\flat ,\mathbf{b}}f\right) \ \mathbf{1}
_{A\backslash I}\ ,  \label{in use}
\end{equation}
given in Proposition \ref{flatness} above, are the modified pseudoprojections $
\left\{ \widehat{\square }_{\pi I}^{\sigma ,\flat ,\mathbf{b}}\right\}
_{I\in \Pi _{1}\mathcal{Q}}$, where $\square _{\pi I}^{\sigma ,\flat ,
\mathbf{b}}$ differs from the orthogonal projection $\square _{\pi
I}^{\sigma ,\pi ,\mathbf{b}}$ for $I\in \Pi _{1}\mathcal{Q}$ by
\begin{eqnarray*}
&& \square _{\pi I}^{\sigma ,\flat ,\mathbf{b}}f-\square _{\pi I}^{\sigma ,\pi ,
\mathbf{b}}f\\
&=&
\left\{ \left( \sum_{I^{\prime }\in \mathfrak{C}_{{
nat}}\left( \pi I\right) }\mathbb{F}_{I^{\prime }}^{\sigma
,b_{A}}f\right) -\mathbb{F}_{\pi I}^{\sigma ,b_{A}}f\right\} -\left\{ \left(
\sum_{I^{\prime }\in \mathfrak{C}\left( \pi I\right) }\mathbb{F}_{I^{\prime
}}^{\sigma ,b_{A}}f\right) -\mathbb{F}_{\pi I}^{\sigma ,b_{A}}f\right\}\\
&=&
-\sum_{I^{\prime }\in \mathfrak{C}_{{brok}}\left( \pi I\right) }
\mathbb{F}_{I^{\prime }}^{\sigma ,b_{A}}f.
\end{eqnarray*}
But the "box support" ${Supp}_{{box}}$ of this last expression $\displaystyle
\sum_{I^{\prime }\in \mathfrak{C}_{{brok}}\left( \pi I\right) }
\mathbb{F}_{I^{\prime }}^{\sigma ,b_{A}}f$ consists of the broken children
of $\pi I$, $\mathfrak{C}_{{brok}}\left( \pi I\right) $, and is
contained in the set 
$$
\bigcup\limits_{I\in \mathcal{C}_{A}'} \bigcup\limits_{I^{\prime }\in \mathfrak{C}_{\mathcal{A}}\left(
A\right) \cap \mathfrak{C}_{\mathcal{D}}\left( \pi I\right) }  \left\{I^{\prime }\right\} 
$$
i.e. 
\begin{eqnarray*}
&&\!\!\!\!\!\!\!\!{Supp}_{{box}}\!\!\!\left( \sum_{I^{\prime }\in \mathfrak{C}_{
{brok}}\left( \pi I\right) }\mathbb{F}_{I^{\prime }}^{\sigma
,b_{A}}f\right) 
\subset
\left\{ I^{\prime }\in \mathfrak{C}_{\mathcal{A}
}\left( A\right) :I^{\prime }\in \mathfrak{C}_{{brok}}\left( \pi
I\right) \text{ for some }I\in \mathcal{C}_{A}'\right\}
\\
&&
\hspace{4.1cm}=\bigcup\limits_{I\in \mathcal{C}_{A}'}\bigcup\limits_{I^{\prime }\in \mathfrak{C}_{\mathcal{A}}\left( A\right)
\cap \mathfrak{C}_{\mathcal{D}}\left( \pi I\right) }\left\{ I^{\prime
}\right\} .
\end{eqnarray*}
But $I\in \Pi _{1}\mathcal{Q}^{S}\subset \mathcal{C}_{A}' $ is a \emph{natural} child of $\pi I$, and so
\begin{equation*}
I\cap {Supp}_{{box}}\left( \sum_{I^{\prime }\in \mathfrak{C}_{{brok}}\left( \pi I\right) }\mathbb{F}_{I^{\prime }}^{\sigma
,b_{A}}f\right) =\emptyset 
\end{equation*}
It now follows that we have
\begin{equation}
E_{I}^{\sigma }\left( \widehat{\square }_{\pi I}^{\sigma ,\flat ,\mathbf{b}
}f\right) =E_{I}^{\sigma }\left( \widehat{\square }_{\pi I}^{\sigma ,\pi ,
\mathbf{b}}f\right) ,\ \ \ \ \ \text{for }I\in \mathcal{C}_{A}' \label{fol}
\end{equation}

Returning to (\ref{in use}), we have from (\ref{sat}) and (\ref{fol}) the
identity 
\begin{eqnarray}
\varphi _{J}^{\mathcal{Q}^{S}} 
&\equiv &
\sum_{I\in \Pi _{1}\mathcal{Q}^{S}:
\mathcal{\ }\left( I,J\right) \in \mathcal{Q}^{S}}b_{A}E_{I}^{\sigma }\left( 
\widehat{\square }_{\pi I}^{\sigma ,\pi ,\mathbf{b}}f\right) \ \mathbf{1}
_{A\backslash I}  \label{iden} \\
&=&
\sum_{I\in \Pi _{1}\mathcal{Q}^{S}:\mathcal{\ }\left( I,J\right) \in 
\mathcal{Q}^{S}}b_{A}E_{I}^{\sigma }\left( \widehat{\square }_{\pi
I}^{\sigma ,\pi ,\mathbf{b}}\left( \mathsf{P}_{\pi \left( \Pi _{1}\mathcal{Q}
\right) }^{\sigma ,\pi ,\mathbf{b}}f\right) \right) \ \mathbf{1}_{A\backslash
I}\ 
\notag
\end{eqnarray}
which will play a critical role in proving the following $\flat $Straddling
and Substraddling lemmas. The $\flat $Straddling Lemma is an adaptation of
Lemmas 3.19 and 3.16 in \cite{Lac}.

\begin{lem}
\label{straddle 3 ref}Let $\mathcal{Q}$ be a reduced admissible collection
of pairs for $A$, and suppose that $\mathcal{S}\subset \Pi _{1}^{{
below}}\mathcal{Q}\cap \mathcal{C}_{A}'$ is a
subpartition of $A$ such that $\mathcal{Q}$ $\flat $straddles $\mathcal{S}$.
Then we have the restricted sublinear norm bound
\begin{equation}
\widehat{\mathfrak{N}}_{{stop},\bigtriangleup ^{\omega }}^{A,
\mathcal{Q}}\leq C_{\mathbf{r}}\sup_{S\in \mathcal{S}}\mathcal{S}_{{
loc}{size}}^{\alpha ,A;S}\left( \mathcal{Q}\right) \leq C_{\mathbf{r}
}\mathcal{S}_{{aug}{size}}^{\alpha ,A}\left( \mathcal{Q}
\right) ,  \label{sub loc bound}
\end{equation}
where $\mathcal{S}_{{loc}{size}}^{\alpha ,A;S}$ is an $S$-localized size condition with an $S$-hole given by
\begin{equation}
\mathcal{S}_{{loc}{size}}^{\alpha ,A;S}\left( \mathcal{Q}
\right) ^{2}\equiv \sup_{K\in \mathcal{W}^{\ast }\left( S\right) \cap 
\mathcal{C}_{A}'}\frac{1}{\left\vert K\right\vert
_{\sigma }}\left( \frac{\mathrm{P}^{\alpha }\left( K,\mathbf{1}_{A\backslash
S}\sigma \right) }{\left\vert K\right\vert^\frac{1}{n} }\right) ^{2}\sum_{J\in \Pi
_{2}^{K,{aug}}\mathcal{Q}}\left\Vert \bigtriangleup _{J}^{\omega ,
\mathbf{b}^{\ast }}x\right\Vert _{L^{2}\left( \omega \right) }^{\spadesuit
2}  \label{localized size ref}
\end{equation}
\end{lem}

\begin{proof}
 We begin by using
that the reduced collection $\mathcal{Q}$ $\flat$straddles $\mathcal{S}$ 
to write
\begin{eqnarray*}
&&\!\!\!\!\!\!\left\vert \mathsf{B}\right\vert _{{stop},\bigtriangleup ^{\omega
}}^{A,\mathcal{Q}}\left( f,g\right) =\sum_{J\in \Pi _{2}\mathcal{Q}}\frac{
\mathrm{P}^{\alpha }\left( J,\left\vert \varphi _{J}^{\mathcal{Q}
}\right\vert \mathbf{1}_{A\backslash I_{\mathcal{Q}}\left( J\right) }\sigma
\right) }{\left\vert J\right\vert^\frac{1}{n} }\left\Vert \bigtriangleup _{J}^{\omega ,
\mathbf{b}^{\ast }}x\right\Vert _{L^{2}\left( \omega \right) }^{\spadesuit
}\left\Vert \square _{J}^{\omega ,\mathbf{b}^{\ast }}g\right\Vert
_{L^{2}\left( \omega \right) }^{\bigstar } \\
&=&
\sum_{S\in \mathcal{S}}\sum_{J\in \Pi _{2}^{S,{aug}}\mathcal{Q
}}\frac{\mathrm{P}^{\alpha }\left( J,\left\vert \varphi _{J}^{\mathcal{Q}
^{S}}\right\vert \mathbf{1}_{A\backslash I_{\mathcal{Q}}\left( J\right)
}\sigma \right) }{\left\vert J\right\vert^\frac{1}{n} }\left\Vert \bigtriangleup
_{J}^{\omega ,\mathbf{b}^{\ast }}x\right\Vert _{L^{2}\left( \omega \right)
}^{\spadesuit }\left\Vert \square _{J}^{\omega ,\mathbf{b}^{\ast
}}g\right\Vert _{L^{2}\left( \omega \right) }^{\bigstar } \\
\end{eqnarray*}
\begin{eqnarray*}
\text{where }\varphi _{J}^{\mathcal{Q}^{S}} &\equiv &\sum_{I\in \Pi _{1}
\mathcal{Q}^{S}:\mathcal{\ }\left( I,J\right) \in \mathcal{Q}
^{S}}b_{A}E_{I}^{\sigma }\left( \widehat{\square }_{\pi I}^{\sigma ,\flat ,
\mathbf{b}}f\right) \ \mathbf{1}_{A\backslash I}\ .
\end{eqnarray*}

At this point we invoke the identity (\ref{iden}),
\begin{equation*}
\varphi _{J}^{\mathcal{Q}^{S}}=\sum_{I\in \Pi _{1}\mathcal{Q}^{S}:\mathcal{\ 
}\left( I,J\right) \in \mathcal{Q}^{S}}b_{A}E_{I}^{\sigma }\left( \widehat{
\square }_{\pi I}^{\sigma ,\pi ,\mathbf{b}}\left( \mathsf{P}_{\pi \left( \Pi
_{1}\mathcal{Q}\right) }^{\sigma ,\pi ,\mathbf{b}}f\right) \right) \ \mathbf{
1}_{A\backslash I}\ ,
\end{equation*}
so that
\begin{equation*}
\left\vert \mathsf{B}\right\vert _{{stop},\bigtriangleup ^{\omega
}}^{A,\mathcal{Q}}\left( f,g\right) =\left\vert \mathsf{B}\right\vert _{
{stop},\bigtriangleup ^{\omega }}^{A,\mathcal{Q}}\left( h,g\right)
,\ \ \ \ \ \text{where }h\equiv \mathsf{P}_{\pi \left( \Pi _{1}\mathcal{Q}
\right) }^{\sigma ,\pi ,\mathbf{b}}f\ .
\end{equation*}
We will treat the sublinear form $\left\vert \mathsf{B}\right\vert _{
{stop},\bigtriangleup ^{\omega }}^{A,\mathcal{Q}}\left( h,g\right) $
with $h=\mathsf{P}_{\pi \left( \Pi _{1}\mathcal{Q}\right) }^{\sigma ,\pi ,
\mathbf{b}}f$ using a small variation on the corresponding argument in Lacey 
\cite{Lac}
. Namely, we will apply a Calder\'{o}n-Zygmund stopping time decomposition to the
function $h=\mathsf{P}_{\pi \left( \Pi _{1}\mathcal{Q}\right) }^{\sigma ,\pi
,\mathbf{b}}f$ on the cube $A$ \ with `obstacle' $\mathcal{S}\cup 
\mathfrak{C}_{A}$ $\left( A\right) $, to obtain stopping times $\mathcal{H}$ 
$\subset \mathcal{C}_{A}$ with the property that for all $H\in \mathcal{H}
\backslash \left\{ A\right\} $ we have 
\begin{eqnarray*}
&&H\in \mathcal{C}_{A}\text{ is not strictly contained in any cube from }
\mathcal{S}, \\
&&E_{H}^{\sigma }\left\vert h\right\vert >\Gamma E_{\pi _{\mathcal{H}
}H}^{\sigma }\left\vert h\right\vert , \\
&&E_{H^{\prime }}^{\sigma }\left\vert h\right\vert \leq \Gamma E_{\pi _{
\mathcal{H}}H}^{\sigma }\left\vert h\right\vert \text{ for all }H\subsetneqq
H^{\prime }\subset \pi _{\mathcal{H}}H\text{ with }H^{\prime }\in \mathcal{C}
_{A}.
\end{eqnarray*}
More precisely, define generation $0$ of $\mathcal{H}$ to consist of the
single cube $A$. Having defined generation $n$, let generation $n+1$
consist of the union over all cubes $M$ in generation $n$ of the maximal
cubes $M^{\prime }$ in $\mathcal{C}_{A}$ that are contained in $M$ with $
E_{M^{\prime }}^{\sigma }\left\vert h\right\vert >\Gamma E_{M}^{\sigma
}\left\vert h\right\vert $, but are \emph{not} strictly contained in any
cube $S$ from $\mathcal{S}$ or contained in any cube $A^{\prime }$
from $\mathfrak{C}_{A}$ $\left( A\right) $ - thus the construction stops at the obstacle $\mathcal{S}\cup \mathfrak{C}_{A}(A)$. Then $
\mathcal{H}$ is the union of all generations $n\geq 0$.

Denote by 
\begin{equation*}
\mathcal{C}_{H}^{\mathcal{H}}\equiv \left\{ H^{\prime }\in \mathcal{C}
_{A}:H^{\prime }\subset H\text{ but }H^{\prime }\not\subset H^{\prime \prime
}\text{ for any }H^{\prime \prime }\in \mathfrak{C}_{\mathcal{H}}\left(
H\right) \right\}
\end{equation*}
the usual $\mathcal{H}$-corona associated with the stopping cube $H$,
but restricted to $\mathcal{C}_{A}$, and let $\alpha _{\mathcal{H}}\left(
H\right) =E_{H}^{\sigma }\left\vert f\right\vert $ as is customary for a
Calder\'{o}n-Zygmund corona. Since these coronas $\mathcal{C}_{H}^{\mathcal{H
}}$ are all contained in $\mathcal{C}_{A}$, we have the stopping energy from
the $\mathcal{A}$-corona $\mathcal{C}_{A}$ at our disposal, which is crucial for the argument. Furthermore, we denote by
\begin{equation}
\mathcal{Q}_{H}\equiv \left\{ \left( I,J\right) \in \mathcal{Q}:J\in 
\mathcal{C}_{H}^{\mathcal{H},\flat {shift}}\right\},\ \ \ \text{with }\mathcal{C}_{H}^{\mathcal{H},\flat {shift}}\equiv \left\{ J\in
\Pi _{2}\mathcal{Q}:J^{\flat }\in \mathcal{C}_{H}^{\mathcal{H}}\right\}
\label{def Q H}
\end{equation}
the restriction of the pairs $\left( I,J\right) $ in $\mathcal{Q}$ to those
for which $J$ lies in the flat shifted $\mathcal{H}$-corona $\mathcal{C}
_{H}^{\mathcal{H},\flat {shift}}$. Since the $\mathcal{H}$-stopping
cubes satisfy a $\sigma $-Carleson condition for $\Gamma $ chosen large
enough, we have the quasiorthogonal inequality 
\begin{equation}
\sum_{H\in \mathcal{H}}\alpha _{\mathcal{H}}\left( H\right) ^{2}\left\vert
H\right\vert _{\sigma }\lesssim \left\Vert h\right\Vert _{L^{2}\left( \sigma
\right) }^{2},  \label{qor}
\end{equation}
which below we will see reduces matters to proving inequality (\ref{sub loc
bound}) for the family of reduced admissible collections $\left\{ \mathcal{Q}
_{H}\right\} _{H\in \mathcal{H}}$ with constants independent of $H$:
\begin{equation*}
\widehat{\mathfrak{N}}_{{stop},\bigtriangleup ^{\omega }}^{A,
\mathcal{Q}_{H}}\leq C_{\mathbf{r}}\sup_{S\in \mathcal{S}}\mathcal{S}_{
{loc}{size}}^{\alpha ,A;S}\left( \mathcal{Q}_{H}\right) \leq
C_{\mathbf{r}}\mathcal{S}_{{aug}{size}}^{\alpha ,A}\left( 
\mathcal{Q}_{H}\right) ,\ \ \ \ \ H\in \mathcal{H}.
\end{equation*}

Given $S\in \mathcal{S}$, define $H_{S}\in \mathcal{H}$ to be the minimal
cube in $\mathcal{H}$ that contains $S$, and then define 
\begin{equation*}
\mathcal{H}_{\mathcal{S}}\equiv \left\{ H_{S}\in \mathcal{H}:S\in \mathcal{S}
\right\} .
\end{equation*}
Note that a given $H\in \mathcal{H}_{\mathcal{S}}$ may have many cubes $
S\in \mathcal{S}$ such that $H=H_{S}$, and we denote the collection of these
cubes by $\mathcal{S}_{H}\equiv \left\{ S\in \mathcal{S}:H_{S}=H\
\right\} $. We will organize the straddling cubes $\mathcal{S}$ as
\begin{equation*}
\mathcal{S}=\bigcup\limits_{H\in \mathcal{H}_{\mathcal{S}
}}\bigcup\limits_{S\in \mathcal{S}_{H}}S
\end{equation*}
where each $S\in \mathcal{S}$ occurs exactly once in the union on the right
hand side, i.e. the collections $\left\{ \mathcal{S}_{H}\right\} _{H\in 
\mathcal{H}_{\mathcal{S}}}$ are pairwise disjoint.

We now momentarily fix $H\in \mathcal{H}_{\mathcal{S}}$, and consider the
reduced admissible collection $\mathcal{Q}_{H}$, so that its projection onto
the second component $\Pi _{2}\mathcal{Q}_{H}$ of $\mathcal{Q}_{H}$\ is 
\emph{contained} in the corona $\mathcal{C}_{H}^{\mathcal{H},\flat {
shift}}$. Then the collection $\mathcal{Q}_{H}$ $\flat $straddles the set $
\mathcal{S}_{H}=\left\{ S\in \mathcal{S}:H_{S}=H\ \right\} $. Moreover, $
\mathcal{Q}_{H}=\bigcup\limits_{S\in \mathcal{S}:\ S\subset H}\mathcal{Q}
_{H}^{S}$ and $\Pi _{2}\mathcal{Q}_{H}^{S}=\Pi _{2}^{S,{aug}}
\mathcal{Q}_{H}$.

Recall that a Whitney cube $K$ was required in the right hand side of
the conclusion of Proposition \ref{flatness} only in the case that $J^{\flat
}\subset S^{\prime \prime }$ for some $S^{\prime \prime }\in \mathfrak{C}_{
\mathcal{D}}^{\left( 2\right) }\left( S\right) $, which of course implies $
3J^{\flat }\subset J^{\maltese }\subset S$. In this case we claim that $K\in 
\mathcal{C}_{A}$. Indeed, suppose in order to derive a contradiction, that $
K\not\in \mathcal{C}_{A}$. Then $J^{\maltese }\not\subset K$, and hence $
3J^{\maltese }\not\subset S$. Since $J^{\maltese }\subset S$, it follows
that $J^{\maltese }$ shares a common part of the boundary with $S$ (since if not, then $
3J^{\maltese }\subset S$, a contradiction). Now Key Fact \#2 in (\ref
{indentation}) implies that the inner grandchild containing $J$, $J^\flat$, is contained in $K$ where $
K\not\in \mathcal{C}_{A}$. This then implies that the pair $\left(
I,J\right) $ belongs to the corona straddling subcollection $\mathcal{P}_{
{cor}}^{A}$, contradicting the assumption that $\mathcal{Q}$ is reduced.

Thus we have $S\in \Pi _{1}^{{below}}\mathcal{Q}\cap \mathcal{C}_{A}'$ and $K\in \mathcal{W}\left( S\right) \cap 
\mathcal{C}_{A}'$ and we can use Proposition \eqref{flatness} with $H=H_{S}$ to bound $\left\vert \mathsf{B}\right\vert _{
{stop},\bigtriangleup ^{\omega }}^{A,\mathcal{Q}}\left( f,g\right) $
by first summing over $H\in \mathcal{H}_{\mathcal{S}}$ and then over $S\in 
\mathcal{S}_{H}$. Indeed, $\mathcal{Q}_{H}$ $\flat $straddles $\mathcal{S}
_{H}\equiv \left\{ S\in \mathcal{S}:H_{S}=H\ \right\} $, so that $\left\vert
\varphi _{J}^{\mathcal{Q}_{H}}\right\vert \lesssim \alpha _{\mathcal{H}
}\left( H\right) \mathbf{1}_{A\backslash I_{\mathcal{Q}_{H}}\left( J\right) }$
by (\ref{phi bound}), and so the sum over $S\in \mathcal{S}_{H}$ of the
first term on the right side of the conclusion of Proposition \eqref{flatness}
is bounded by
\begin{eqnarray*}
&&\!\!\!\!\!\!\!\!\!\!\!\!\!\!\!\!\!
\alpha _{\mathcal{H}}\left( H\right) \sum_{S\in \mathcal{S}_{H}}\!\!\!\!\sqrt{\left\vert S\right\vert _{\sigma }}\frac{1}{\sqrt{\left\vert S\right\vert
_{\sigma }}}\left( \frac{\mathrm{P}^{\alpha }\left( S,\mathbf{1}_{A\backslash
S}\sigma \right) }{\left\vert S\right\vert^\frac{1}{n} }\right)\left\Vert \mathsf{Q}
_{\Pi _{2}^{S,{aug}}\mathcal{Q}_{H}}^{\omega ,\mathbf{b}^{\ast
}}x\right\Vert _{L^{2}\left( \omega \right) }^{\spadesuit }\!\!\left\Vert 
\mathsf{P}_{\Pi _{2}^{S,{aug}}\mathcal{Q}_{H}}^{\omega ,\mathbf{b}
^{\ast }}g\right\Vert _{L^{2}\left( \omega \right) }^{\bigstar } \\
&\leq &
\alpha _{\mathcal{H}}\left( H\right) \left\{ \sup_{S\in \mathcal{S}
_{H}}\frac{1}{\sqrt{\left\vert S\right\vert _{\sigma }}}\left( \frac{\mathrm{
P}^{\alpha }\left( S,\mathbf{1}_{A\backslash S}\sigma \right) }{\left\vert
S\right\vert^\frac{1}{n} }\right) \left\Vert \mathsf{Q}_{\Pi _{2}^{S,{aug}}
\mathcal{Q}_{H}}^{\omega ,\mathbf{b}^{\ast }}x\right\Vert _{L^{2}\left(
\omega \right) }^{\spadesuit }\right\} \cdot \\
&& 
\hspace{6.65cm}\cdot \sum_{S\in \mathcal{S}_{H}}\sqrt{
\left\vert S\right\vert _{\sigma }}\left\Vert \mathsf{P}_{\Pi _{2}^{S,
{aug}}\mathcal{Q}_{H}}^{\omega ,\mathbf{b}^{\ast }}g\right\Vert
_{L^{2}\left( \omega \right) }^{\bigstar } \\
&\leq &\alpha _{\mathcal{H}}\left( H\right) \left\{ \sup_{S\in \mathcal{S}
_{H}}\mathcal{S}_{{loc}{size}}^{\alpha ,A;S}\left( \mathcal{Q
}_{H}\right) \right\} \sqrt{\left\vert H\right\vert _{\sigma }}\left\Vert 
\mathsf{P}_{\Pi _{2}\mathcal{Q}_{H}}^{\omega ,\mathbf{b}^{\ast
}}g\right\Vert _{L^{2}\left( \omega \right) }^{\bigstar }\
\end{eqnarray*}
where $\Pi _{2}^{K,{aug}}\mathcal{Q}_{H}$ is as in Definition \ref
{augs}, and the corresponding sum over $S\in \mathcal{S}_{H}$ of the
second term is bounded by
\begin{eqnarray*}
&&\!\!\!\!\!\!\!\!\!\!\!\!\!\!\!\!\!
\alpha _{\mathcal{H}}\left( H\right)\!\!\! \sum_{S\in \mathcal{S}_{H}}\!\sum_{K\in 
\mathcal{W}\left( S\right) \cap \mathcal{C}_{A}'}\!\!\!\frac{\sqrt{\left\vert K\right\vert _{\sigma }}}{\sqrt{\left\vert K\right\vert
_{\sigma }}}\frac{\mathrm{P}^{\alpha }\left( K,\mathbf{1}_{A\backslash
S}\sigma \right) }{\left\vert K\right\vert^\frac{1}{n} } \!\!\left\Vert \mathsf{Q}
_{\Pi _{2}^{K,{aug}}\mathcal{Q}_{H}^{S}}^{\omega ,\mathbf{b}^{\ast
}}x\right\Vert _{L^{2}\left(\omega\! \right) }^{\spadesuit }\!\!\left\Vert 
\mathsf{P}_{\Pi _{2}^{K,{aug}}\mathcal{Q}_{H}^{S}}^{\omega ,\mathbf{b
}^{\ast }}g\right\Vert _{L^{2}\left(\omega\! \right) }^{\bigstar } \\
&\lesssim &
\alpha _{\mathcal{H}}\left( H\right) \sup_{S\in \mathcal{S}_{
\mathcal{H}}}\mathcal{S}_{{loc}{size}}^{\alpha ,A;S}\left( 
\mathcal{Q}_{H}\right) \left( \sum_{S\in \mathcal{S}}\sum_{K\in \mathcal{W}
\left( S\right) }\left\vert K\right\vert _{\sigma }\right) ^{\frac{1}{2}
}\left\Vert \mathsf{P}_{\Pi _{2}\mathcal{Q}_{H}}^{\omega ,\mathbf{b}^{\ast
}}g\right\Vert _{L^{2}\left( \omega \right) }^{\bigstar } \\
&\leq &\left\{ \sup_{S\in \mathcal{S}_{\mathcal{H}}}\mathcal{S}_{{loc
}{size}}^{\alpha ,A;S}\left( \mathcal{Q}_{H}\right) \right\} \alpha
_{\mathcal{H}}\left( H\right) \sqrt{\left\vert H\right\vert _{\sigma }}
\left\Vert \mathsf{P}_{\Pi _{2}\mathcal{Q}_{\mathcal{H}}}^{\omega ,\mathbf{b}
^{\ast }}g\right\Vert _{L^{2}\left( \omega \right) }^{\bigstar }
\end{eqnarray*}

Using the definition of $\left\vert \mathsf{B}\right\vert _{{stop}
,\bigtriangleup ^{\omega }}^{A,\mathcal{Q}}\left( f,g\right) $, we now sum the previous inequalities over the cubes $H\in 
\mathcal{H}_{\mathcal{S}}$ to obtain the following string of inequalities
(explained in detail after the display)
\begin{eqnarray}
&&\!\!\!\!\!\!\!\!\!\!\!\!\!\!\!\!\!\!\!\!\!\!\!\!\!\!\!\!\!\!
\left\vert \mathsf{B}\right\vert _{{stop},\bigtriangleup ^{\omega
}}^{A,\mathcal{Q}}\left( f,g\right) \leq \left\{ \sup_{S\in \mathcal{S}}
\mathcal{S}_{{loc}{size}}^{\alpha ,A;S}\left( \mathcal{Q}
\right) \right\} \sum_{H\in \mathcal{H}_{\mathcal{S}}}\alpha _{\mathcal{H}
}\left( H\right) \sqrt{\left\vert H\right\vert _{\sigma }}\left\Vert \mathsf{
P}_{\Pi _{2}\mathcal{Q}_{H}}^{\omega ,\mathbf{b}^{\ast }}g\right\Vert
_{L^{2}\left( \omega \right) }^{\bigstar }  \notag \label{unfix} \\
&\leq &\left\{ \sup_{S\in \mathcal{S}}\mathcal{S}_{{loc}{size
}}^{\alpha ,A;S}\left( \mathcal{Q}\right) \right\} \sqrt{\sum_{H\in \mathcal{
H}_{\mathcal{S}}}\alpha _{\mathcal{H}}\left( H\right) ^{2}\left\vert
H\right\vert _{\sigma }}\sqrt{\sum_{H\in \mathcal{H}_{\mathcal{S}
}}\left\Vert \mathsf{P}_{\Pi _{2}\mathcal{Q}_{H}}^{\omega ,\mathbf{b}^{\ast
}}g\right\Vert _{L^{2}\left( \omega \right) }^{\bigstar 2}}  \notag \\
&\lesssim &\left\{ \sup_{S\in \mathcal{S}}\mathcal{S}_{{loc}{
size}}^{\alpha ,A;S}\left( \mathcal{Q}\right) \right\} \left\Vert
h\right\Vert _{L^{2}\left( \sigma \right) }\sqrt{\sum_{H\in \mathcal{H}_{
\mathcal{S}}}\left\Vert \mathsf{P}_{\Pi _{2}\mathcal{Q}_{H}}^{\omega ,
\mathbf{b}^{\ast }}g\right\Vert _{L^{2}\left( \omega \right) }^{\bigstar 2}}
\notag \\
&\leq &\left\{ \sup_{S\in \mathcal{S}}\mathcal{S}_{{loc}{size
}}^{\alpha ,A;S}\left( \mathcal{Q}\right) \right\} \left\Vert \mathsf{P}
_{\pi \left( \Pi _{1}\mathcal{Q}\right) }^{\sigma ,\pi ,\mathbf{b}
}f\right\Vert _{L^{2}\left( \sigma \right) }\left\Vert \mathsf{P}_{\Pi _{2}
\mathcal{Q}}^{\omega ,\mathbf{b}^{\ast }}g\right\Vert _{L^{2}\left( \omega
\right) }^{\bigstar }  \notag \\
&\lesssim &\left\{ \sup_{S\in \mathcal{S}}\mathcal{S}_{{loc}{
size}}^{\alpha ,A;S}\left( \mathcal{Q}\right) \right\} \left\Vert \mathsf{P}
_{\pi \left( \Pi _{1}\mathcal{Q}\right) }^{\sigma ,\mathbf{b}}f\right\Vert
_{L^{2}\left( \sigma \right) }^{\bigstar }\left\Vert \mathsf{P}_{\Pi _{2}
\mathcal{Q}}^{\omega ,\mathbf{b}^{\ast }}g\right\Vert _{L^{2}\left( \omega
\right) }^{\bigstar }\   \notag
\end{eqnarray}
where in the first line we have used $\mathcal{Q}=\bigcup\limits_{H\in 
\mathcal{H}_{\mathcal{S}}}\mathcal{Q}_{H}$, which follows from the fact that
each $J^{\flat }$ is contained in a unique $S\in \mathcal{S}$; in the third
line we have used the quasiorthogonal inequality (\ref{qor}); in the fourth
line we have used that the sets $\Pi _{2}\mathcal{Q}_{H}\subset \mathcal{C}
_{H}^{\mathcal{H},\flat {shift}}$ are pairwise disjoint in $H$ and have
union $\displaystyle \Pi _{2}\mathcal{Q}=\overset{\cdot }{\bigcup_{H\in \mathcal{H}_{\mathcal{S}}}}\Pi _{2}\mathcal{Q}_{H}$. In the final line, we have used first
the equality (\ref{box pi equals}), second the fact that the functions $
\square _{I,{brok}}^{\sigma ,\pi ,\mathbf{b}}f$ have pairwise
disjoint supports, third the upper weak Riesz inequality and fourth the estimate (\ref{F est}) - which relies on the
reverse H\"{o}lder property for children in Lemma \ref{prelim control of
corona} - to obtain 
\begin{eqnarray}
\left\Vert \mathsf{P}_{\pi \left( \Pi _{1}\mathcal{Q}\right) }^{\sigma ,\pi ,
\mathbf{b}}f\right\Vert _{L^{2}\left( \sigma \right) }^{2} &=&\left\Vert
\sum_{I\in \pi \left( \Pi _{1}\mathcal{Q}\right) }\square _{I}^{\sigma ,
\mathbf{b}}f-\sum_{I\in \pi \left( \Pi _{1}\mathcal{Q}\right) }\square _{I,
{brok}}^{\sigma ,\pi ,\mathbf{b}}f\right\Vert _{L^{2}\left( \sigma
\right) }^{2}  \label{needed for unfix} \notag \\
&\lesssim &\left\Vert \sum_{I\in \pi \left( \Pi _{1}\mathcal{Q}\right)
}\square _{I}^{\sigma ,\mathbf{b}}f\right\Vert _{L^{2}\left( \sigma \right)
}^{2}+\left\Vert \sum_{I\in \pi \left( \Pi _{1}\mathcal{Q}\right) }\square
_{I,{brok}}^{\sigma ,\pi ,\mathbf{b}}f\right\Vert _{L^{2}\left(
\sigma \right) }^{2}  \notag \\
&\lesssim &\left\Vert \mathsf{P}_{\pi \left( \Pi _{1}\mathcal{Q}\right)
}^{\sigma ,\mathbf{b}}f\right\Vert _{L^{2}\left( \sigma \right)
}^{2}+\sum_{I\in \pi \left( \Pi _{1}\mathcal{Q}\right) }\left\Vert \square
_{I,{brok}}^{\sigma ,\pi ,\mathbf{b}}f\right\Vert _{L^{2}\left(
\sigma \right) }^{2}  \notag \\
&\lesssim &\sum_{I\in \pi \left( \Pi _{1}\mathcal{Q}\right) }\left\Vert
\square _{I}^{\sigma ,\mathbf{b}}f\right\Vert _{L^{2}\left( \sigma \right)
}^{2}+\sum_{I\in \pi \left( \Pi _{1}\mathcal{Q}\right) }\left\Vert
\bigtriangledown _{I}^{\sigma }f\right\Vert _{L^{2}\left( \sigma \right)
}^{2} \\ &\lesssim& \left\Vert \mathsf{P}_{\pi \left( \Pi _{1}\mathcal{Q}\right)
}^{\sigma ,\mathbf{b}}f\right\Vert _{L^{2}\left( \sigma \right) }^{\bigstar
2}  \notag
\end{eqnarray}

\medskip

We now use the fact that the supremum in the definition of $\mathcal{S}_{
{loc}{size}}^{\alpha ,A;S}\left( \mathcal{Q}\right) $ is
taken over $K\in \mathcal{W}^{\ast }\left( S\right) \cap \mathcal{C}_{A}'$ to conclude that 
\begin{equation*}
\sup_{S\in \mathcal{S}}\mathcal{S}_{{loc}{size}}^{\alpha
,A;S}\left( \mathcal{Q}\right) \leq \mathcal{S}_{{aug}{size}
}^{\alpha ,A}\left( \mathcal{Q}\right) ,
\end{equation*}
and this completes the proof of Lemma \ref{straddle 3 ref}.
\end{proof}

In a similar fashion we can obtain the following Substraddling Lemma.

\begin{dfn}
\label{def substraddles}Given a \emph{reduced admissible} collection of
pairs $\mathcal{Q}$ for $A$, and a $\mathcal{D}$-cube $L$ contained in $
A $, we say that $\mathcal{Q}$ \textbf{substraddles} $L$ if for every pair $
\left( I,J\right) \in \mathcal{Q}$ there is $K\in \mathcal{W}\left( L\right)
\cap \mathcal{C}_{A}'$ with $J\subset K\subset 3K\subset
I\subset L$.
\end{dfn}

\begin{lem}
\label{substraddle ref}Let $L$ be a $\mathcal{D}$-cube contained in $A$,
and suppose that $\mathcal{Q}$ is an admissible collection of pairs that
substraddles $L$. Then we have the sublinear form bound
\begin{equation*}
\widehat{\mathfrak{N}}_{{stop},\bigtriangleup ^{\omega }}^{A,
\mathcal{Q}}\leq C\mathcal{S}_{{aug}{size}}^{\alpha
,A}\left( \mathcal{Q}\right) .
\end{equation*}
\end{lem}

\begin{proof}
We will show that $\mathcal{Q}$ $\flat $straddles the subset $\mathcal{W}
_{L} $ of Whitney cubes for $L$ given by
\begin{equation*}
\mathcal{W}^{\mathcal{Q}}\left( L\right) \equiv \left\{ K\in \mathcal{W}
\left( L\right) \cap \mathcal{C}_{A}':J\subset K\subset
3K\subset I\subset L\text{ for some }\left( I,J\right) \in \mathcal{Q}
\right\} .
\end{equation*}
It is clear that $\mathcal{W}^{\mathcal{Q}}\left( L\right) \subset \Pi _{1}^{
{below}}\mathcal{Q}\cap \mathcal{C}_{A}'$ is a
subpartition of $A$. It remains to show that for every pair $\left(
I,J\right) \in \mathcal{Q}$ there is $K\in \mathcal{W}^{\mathcal{Q}}\left(
L\right) \cap \left[ J,I\right] $ such that $J^{\flat }\subset K$. But our
hypothesis implies that there is $K\in \mathcal{W}^{\mathcal{Q}}\left(
L\right) $ with $J\subset K\subset 3K\subset I\subset L$. We now consider
two cases.

\textbf{Case 1}: If $\pi _{\mathcal{D}}^{\left( 3\right) }K\subset L$, then since $K$ is maximal Whitney cube, it is contained in an \emph{outer}
grandchild of $\pi _{\mathcal{D}}^{\left( 3\right) }K$ and $\pi _{\mathcal{D}}^{\left( 1\right) }K$ has to share an endpoint with $L$. Then so does $\pi_{\mathcal{D}}^{\left( 3\right) }K$. Recall, from Key Fact \#2 in (\ref{indentation}), $3J\subset J^\flat$, an \emph{inner} grandchild of $J^{\maltese }$. We thus have $J^{\maltese }\subset \pi
_{\mathcal{D}}^{\left( 2\right) }K$ (If not; $\pi
_{\mathcal{D}}^{\left( 2\right) }K\subset J^{\maltese}$ which implies that $J^\flat$ has the same endpoint as $L$, a contradiction). This implies that $J^{\flat }\subset K$
.

\textbf{Case 2}: If $\pi _{\mathcal{D}}^{\left( 3\right) }K\varsupsetneqq L$
, then $K\subset 3K\subset I\subset L$ implies that $I=L=\pi _{\mathcal{D}
}^{\left( 2\right) }K$. Thus we have $J^{\maltese }\subset I=\pi _{\mathcal{D
}}^{\left( 2\right) }K$, which again gives $J^{\flat }\subset K$.

Now that we know $\mathcal{Q}$ $\flat $straddles the subset $\mathcal{W}^{
\mathcal{Q}}\left( L\right) $, we can apply Lemma \ref{straddle 3 ref} to
obtain the required bound $\widehat{\mathfrak{N}}_{{stop}
,\bigtriangleup ^{\omega }}^{A,\mathcal{Q}}\leq C\mathcal{S}_{{aug}
{size}}^{\alpha ,A}\left( \mathcal{Q}\right) $.
\end{proof}

\subsection{The bottom/up stopping time argument of M. Lacey}

Before introducing Lacey's stopping times, we note that the
Corona-straddling Lemma \ref{cor strad 1} allows us to remove the `corona
straddling' collection $\mathcal{P}_{{cor}}^{A}$ of pairs of cubes
in (\ref{def cor}) from the collection $\mathcal{P}^{A}$ in (\ref{initial P}
) used to define the stopping form $\mathsf{B}_{{stop}}^{A}\left(
f,g\right) $. The collection $\mathcal{P}^{A}\backslash \mathcal{P}_{{cor
}}^{A}$ is of course also $A$-admissible.

\label{assume}We assume for the remainder of the proof that all admissible
collections $\mathcal{P}$ are reduced, i.e. 
\begin{equation}
\mathcal{P}^{A}\cap \mathcal{P}_{{cor}}^{A}=\emptyset ,\text{ as well
as }\mathcal{P}\cap \mathcal{P}_{{cor}}^{A}=\emptyset \text{ for all }A
\text{-admissible }\mathcal{P}.  \label{empty assumption}
\end{equation}

For a cube $K\in \mathcal{D}$, we define
\begin{equation*}
\mathcal{G}\left[ K\right] \equiv \left\{ J\in \mathcal{G}:J\subset K\right\}
\end{equation*}
to consist of all cubes $J$ in the other grid $\mathcal{G}$ that are
contained in $K$. For an $A$-admissible collection $\mathcal{P}$ of pairs,
define two atomic measures $\omega _{\mathcal{P}}$ and $\omega _{\flat 
\mathcal{P}}$ in the upper half space $\mathbb{R}_{+}^{n+1}$ by
\begin{equation}
\omega _{\mathcal{P}}\equiv \sum_{J\in \Pi _{2}\mathcal{P}}\left\Vert
\bigtriangleup _{J}^{\omega ,\mathbf{b}^{\ast }}x\right\Vert _{L^{2}\left(
\omega \right) }^{\spadesuit 2}\ \delta _{\left( c_{J^{\maltese }},\ell
\left( J^{\maltese }\right) \right) }
\end{equation}
and
\begin{equation}
\omega _{\flat \mathcal{P}
}\equiv \sum_{J\in \Pi _{2}\mathcal{P}}\left\Vert \bigtriangleup
_{J}^{\omega ,\mathbf{b}^{\ast }}x\right\Vert _{L^{2}\left( \omega \right)
}^{\spadesuit 2}\ \delta _{\left( c_{J^{\flat }},\ell \left( J^{\flat
}\right) \right) },  \label{def atomic}
\end{equation}
where $J^\flat$ is the inner grandchild of $
J^{\maltese }$ that contains $J$

Note that each cube $J\in \Pi _{2}\mathcal{P}$ has its
`energy' $\left\Vert \bigtriangleup _{J}^{\omega ,\mathbf{b}^{\ast
}}x\right\Vert _{L^{2}\left( \omega \right) }^{\spadesuit 2}$ in the measure 
$\omega _{\flat \mathcal{P}}$ assigned to exactly one of the $2^n$ points $
\left( c_{J^{\flat }},\frac{1}{4}\ell \left( J^{\maltese }\right)
\right) $ in the upper half plane $\mathbb{R}_{+}^{n+1}$
since $J$ is contained in one of $J_{\searrow}^{\maltese }$, namely in $J^\flat$, by Key Fact \#2 in (\ref{indentation}). Note also that
the atomic measure $\omega _{\flat \mathcal{P}}$ differs from the measure $
\mu $ in (\ref{def mu n}) in Appendix  below - which is used there to
control the functional energy condition - in that here we bundle together
all the $J^{\prime }s$ having a common $J^{\flat }$. This is in order to
rewrite the \emph{augmented} size functional in terms of the measure $\omega
_{\flat \mathcal{P}}$. We can get away with this here, as opposed to in
Appendix , due to the `smaller and decoupled' nature of the augmented size
functional to which we will relate $\omega _{\flat \mathcal{P}}$.

Define the tent $\mathbf{T}\left( L\right) $ over a cube $L$ to be the
convex hull of the cube $L\times \left\{ 0\right\} $ and the point $
\left( c_{L},\ell \left( L\right) \right) \in \mathbb{R}_{+}^{n+1}$. Then for $
J\in \Pi _{2}\mathcal{P}$ we have $J\in \Pi _{2}^{K,{aug}}\mathcal{P}
$ \emph{iff} $\left\{ J\subset K\text{ and }J^{\maltese }\subset \pi _{
\mathcal{D}}^{\left( 2\right) }K\right\} $ \emph{iff} $J^\flat\subset K$ \emph{iff} $\left( c_{J^{\flat
}},\ell \left( J^{\flat }\right) \right) \in \mathbf{T}\left( K\right) $. We
can now rewrite the augmented size functional of $\mathcal{P}$ in Definition 
\ref{augs}\ as
\begin{equation}
\mathcal{S}_{{aug}{size}}^{\alpha ,A}\left( \mathcal{P}
\right) ^{2}\equiv \sup_{K\in \Pi _{1}^{{below}}\mathcal{P}\cap 
\mathcal{C}_{A}'}\frac{1}{\left\vert K\right\vert
_{\sigma }}\left( \frac{\mathrm{P}^{\alpha }\left( K,\mathbf{1}_{A\backslash
K}\sigma \right) }{\left\vert K\right\vert^\frac{1}{n} }\right) ^{2}\omega _{\flat 
\mathcal{P}}\left( \mathbf{T}\left( K\right) \right) .
\label{def P stop energy' 3}
\end{equation}
It will be convenient to write
\begin{equation*}
\Psi ^{\alpha }\left( K;\mathcal{P}\right) ^{2}\equiv \left( \frac{\mathrm{P}
^{\alpha }\left( K,\mathbf{1}_{A\backslash K}\sigma \right) }{\left\vert
K\right\vert^\frac{1}{n} }\right) ^{2}\omega _{\flat \mathcal{P}}\left( \mathbf{T}\left(
K\right) \right) ,
\end{equation*}
so that we have simply
\begin{equation*}
\mathcal{S}_{{aug}{size}}^{\alpha ,A}\left( \mathcal{P}
\right) ^{2}=\sup_{K\in \Pi _{1}^{{below}}\mathcal{P}\cap \mathcal{C}
_{A}'}\frac{\Psi ^{\alpha }\left( K;\mathcal{P}\right)
^{2}}{\left\vert K\right\vert _{\sigma }}.
\end{equation*}

\begin{rem}
The functional $\omega _{\flat \mathcal{P}}\left( \mathbf{T}\left( K\right)
\right) $ is increasing in $K$, while the functional $\frac{\mathrm{P}
^{\alpha }\left( K,\mathbf{1}_{A\backslash K}\sigma \right) }{\left\vert
K\right\vert^\frac{1}{n} }$ is `almost decreasing' in $K$: if $K_{0}\subset K$ then
\begin{eqnarray*}
\frac{\mathrm{P}^{\alpha }\left( K,\mathbf{1}_{A\backslash K}\sigma \right) }{
\left\vert K\right\vert ^\frac{1}{n}} 
&=&
\int_{A\backslash K}\frac{d\sigma \left(
y\right) }{\left( \left\vert K\right\vert^\frac{1}{n} +\left\vert y-c_{K}\right\vert
\right) ^{n+1-\alpha }} \\
&\lesssim &
\int_{A\backslash K}\frac{(\sqrt{n})^{n+1-\alpha}d\sigma \left( y\right) }{\left(
\left\vert K_{0}\right\vert^\frac{1}{n} +\left\vert y-c_{K_{0}}\right\vert \right)
^{n+1-\alpha }} \\
&\leq &
\int_{A\backslash K_{0}}\frac{C_{\alpha,n }\ \ d\sigma \left( y\right) }{
\left( \left\vert K_{0}\right\vert^\frac{1}{n} +\left\vert y-c_{K_{0}}\right\vert
\right) ^{n+1-\alpha }}
=
C_{\alpha,n}\frac{\mathrm{P}^{\alpha }\left( K_{0},
\mathbf{1}_{A\backslash K_{0}}\sigma \right) }{\left\vert K_{0}\right\vert^\frac{1}{n} }
\end{eqnarray*}
since $\left\vert K_{0}\right\vert +\left\vert y-c_{K_{0}}\right\vert \leq
\left\vert K\right\vert +\left\vert y-c_{K}\right\vert +\frac{1}{2}\diam\left( K\right) $ for $y\in A\backslash K$.
\end{rem}

Recall that if $\mathcal{P}$ is an admissible collection for a dyadic
cube $A$, the corresponding sublinear form in (\ref{First inequality}) is given by
\begin{eqnarray*}
\left\vert \mathsf{B}\right\vert _{{stop},\bigtriangleup ^{\omega
}}^{A,\mathcal{P}}\left( f,g\right) &\equiv &\sum_{J\in \Pi _{2}\mathcal{P}}
\frac{\mathrm{P}^{\alpha }\left( J,\left\vert \varphi _{J}^{\mathcal{P}
}\right\vert \mathbf{1}_{A\backslash I_{\mathcal{P}}\left( J\right) }\sigma
\right) }{\left\vert J\right\vert^\frac{1}{n} }\left\Vert \bigtriangleup _{J}^{\omega ,
\mathbf{b}^{\ast }}x\right\Vert _{L^{2}\left( \omega \right) }^{\spadesuit
}\left\Vert \square _{J}^{\omega ,\mathbf{b}^{\ast }}g\right\Vert
_{L^{2}\left( \omega \right) }^{\bigstar }; \\
\text{where }\varphi _{J}^{\mathcal{P}} &\equiv &\sum_{I\in \mathcal{C}_{A}^{
'}:\ \left( I,J\right) \in \mathcal{P}}b_{A}E_{I}^{\sigma
}\left( \widehat{\square }_{\pi I}^{\sigma ,\flat ,\mathbf{b}}f\right) \ 
\mathbf{1}_{A\backslash I}\ .
\end{eqnarray*}
In the notation for $\left\vert \mathsf{B}\right\vert _{{stop}
,\bigtriangleup ^{\omega }}^{A,\mathcal{P}}$, we are omitting dependence on
the parameter $\alpha $, and to avoid clutter, we will often do so from now
on when the dependence on $\alpha $ is inconsequential.

Recall further that the `size testing collection' of cubes $\Pi _{1}^{
{below}}\mathcal{P}$ for the initial size testing functional $
\mathcal{S}_{{init}{size}}^{\alpha ,A}\left( \mathcal{P}
\right) $ is the collection of all subcubes of cubes in $\Pi _{1}
\mathcal{P}$, and moreover, by Key Fact \#1 in (\ref{later use}), that we
can restrict the collection to $\Pi _{1}^{{below}}\mathcal{P}\cap 
\mathcal{C}_{A}'$. This latter set is used for the
augmented size functional.

\begin{description}
\item[Assumption]
$ $\\
We may assume that the corona $\mathcal{C}_{A}$ is finite,
and that each $A$-admissible collection $\mathcal{P}$ is a finite
collection, and hence so are $\Pi _{1}\mathcal{P}$, $\Pi _{1}^{{below
}}\mathcal{P}\cap \mathcal{C}_{A}'$ and $\Pi _{2}
\mathcal{P}$, provided all of the bounds we obtain are independent of the
cardinality of these latter collections.
\end{description}

Consider $0<\varepsilon <1$, where $\rho =1+\varepsilon $ will be chosen
later in (\ref{choose rho}). Begin by defining the collection $\mathcal{L}
_{0}$ to consist of the \emph{minimal} dyadic cubes $K$ in $\Pi _{1}^{
{below}}\mathcal{P}\cap \mathcal{C}_{A}'$ such
that
\begin{equation*}
\frac{\Psi ^{\alpha }\left( K;\mathcal{P}\right) ^{2}}{\left\vert
K\right\vert _{\sigma }}\geq \varepsilon \mathcal{S}_{{aug}{
size}}^{\alpha ,A}\left( \mathcal{P}\right) ^{2}.
\end{equation*}
where we recall that
\begin{equation*}
\Psi ^{\alpha }\left( K;\mathcal{P}\right) ^{2}\equiv \left( \frac{\mathrm{P}
^{\alpha }\left( K,\mathbf{1}_{A\backslash K}\sigma \right) }{\left\vert
K\right\vert^\frac{1}{n} }\right) ^{2}\omega _{\flat \mathcal{P}}\left( \mathbf{T}\left(
K\right) \right) .
\end{equation*}
Note that such minimal cubes exist when $0<\varepsilon <1$ because $
\mathcal{S}_{{aug}{size}}^{\alpha ,A}\left( \mathcal{P}
\right) ^{2}$ is the supremum over $K\in \Pi _{1}^{{below}}\mathcal{P
}\cap \mathcal{C}_{A}'$ of $\frac{\Psi ^{\alpha }\left(
K;\mathcal{P}\right) ^{2}}{\left\vert K\right\vert _{\sigma }}$. A key
property of the minimality requirement is that
\begin{equation}
\frac{\Psi ^{\alpha }\left( K^{\prime };\mathcal{P}\right) ^{2}}{\left\vert
K^{\prime }\right\vert _{\sigma }}<\varepsilon \mathcal{S}_{{aug}
{size}}^{\alpha ,A}\left( \mathcal{P}\right) ^{2},
\label{key property 3}
\end{equation}
whenever there is $K^{\prime }\in \Pi _{1}^{{below}}\mathcal{P}\cap 
\mathcal{C}_{A}'$ with $K^{\prime }\varsubsetneqq K$ and 
$K\in \mathcal{L}_{0}$.

We now perform a stopping time argument `from the bottom up' with respect to
the atomic measure $\omega _{\mathcal{P}}$ in the upper half space. This
construction of a stopping time `from the bottom up', together with the
subsequent applications of the Orthogonality Lemma and the Straddling Lemma,
comprise the key innovations in Lacey's argument \cite{Lac}. However, in our
situation the cubes $I$ belonging to $\Pi _{1}^{{below}}\mathcal{
P}$ are no longer `good' in any sense, and we must include an additional
top/down stopping criterion in the next subsection to accommodate this lack
of `goodness'. The argument in \cite{Lac} will apply to these special
stopping cubes, called `indented' cubes, and the remaining cubes
form towers with a common endpoint, that are controlled using all three
straddling lemmas.

We refer to $\mathcal{L}_{0}$ as the initial or level $0$ generation of
stopping cubes. Set
\begin{equation}
\rho =1+\varepsilon .  \label{def rho}
\end{equation}
As in \cite{SaShUr7}, \cite{SaShUr9} and \cite{SaShUr10}, we follow Lacey 
\cite{Lac} by recursively defining a finite sequence of generations $\left\{ 
\mathcal{L}_{m}\right\} _{m\geq 0}$ by letting $\mathcal{L}_{m}$ consist of
the \emph{minimal} dyadic cubes $L$ in $\Pi _{1}^{{below}}
\mathcal{P}\cap \mathcal{C}_{A}'$ that contain a cube from some previous level $\mathcal{L}_{\ell }$, $\ell <m$, such that
\begin{equation}
\omega _{\flat \mathcal{P}}\left( \mathbf{T}\left( L\right) \right) \geq
\rho \omega _{\flat \mathcal{P}}\left( \bigcup\limits_{L^{\prime }\in
\bigcup\limits_{\ell =0}^{m-1}\mathcal{L}_{\ell }:\ L^{\prime }\subset L}
\mathbf{T}\left( L^{\prime }\right) \right) .  \label{up stopping condition}
\end{equation}
Since $\mathcal{P}$ is finite this recursion stops at some level $M$. We
then let $\mathcal{L}_{M+1}$ consist of all the maximal cubes in $\Pi
_{1}^{{below}}\mathcal{P}\cap \mathcal{C}_{A}'$
that are not already in some $\mathcal{L}_{m}$ with $m\leq M$. Thus $
\mathcal{L}_{M+1}$ will contain either none, some, or all of the maximal
cubes in $\Pi _{1}^{{below}}\mathcal{P}$. We do not of course
have (\ref{up stopping condition}) for $A^{\prime }\in \mathcal{L}_{M+1}$ in
this case, but we do have that (\ref{up stopping condition}) fails for
subcubes $K$ of $A^{\prime }\in \mathcal{L}_{M+1}$ that are not
contained in any other $L\in \mathcal{L}_{m}$ with $m\leq M$, and this is
sufficient for the arguments below.

We now decompose the collection of pairs $\left( I,J\right) $ in $\mathcal{P}
$ into collections $\mathcal{P}^{\flat small}$ and $\mathcal{P}^{\flat big}$
according to the location of $I$ and $J^{\flat }$, but only after
introducing below the indented corona $\mathcal{H}$. The collection $
\mathcal{P}^{\flat big}$ will then essentially consist of those pairs $
\left( I,J\right) \in \mathcal{P}$ for which there are $L^{\prime },L\in 
\mathcal{H}$ with $L^{\prime }\varsubsetneqq L$ and such that $J^{\flat }\in 
\mathcal{C}_{L^{\prime }}^{\mathcal{H}}$ and $I\in \mathcal{C}_{L}^{\mathcal{
H}}$. The collection $\mathcal{P}^{\flat small}$ will consist of the
remaining pairs $\left( I,J\right) \in \mathcal{P}$ for which there is $L\in 
\mathcal{H}$ such that $J^{\flat },I\in \mathcal{C}_{L}^{\mathcal{H}}$,
along with the pairs $\left( I,J\right) \in \mathcal{P}$ such that $I\subset
I_{0}$ for some $I_{0}\in \mathcal{L}_{0}$. This will cover all pairs $
\left( I,J\right) $ in $\mathcal{P}\subset \mathcal{P}_{A}$, since for such
pairs, $I\in \mathcal{C}_{A}'$\ and $J\in \mathcal{C}
_{A}^{\mathcal{G}{shift}}$, which in turn implies $I\in \mathcal{C}
_{L}^{\mathcal{H}}$ and $J^{\flat }\in \mathcal{C}_{L^{\prime }}^{\mathcal{H}
}$ for some $L,L^{\prime }\in \mathcal{H}$. But a considerable amount of
further analysis is required to prove (\ref{First inequality}).

First recall that $\mathcal{L}\equiv \bigcup\limits_{m=0}^{M+1}\mathcal{L}
_{m}$ is the tree of stopping $\omega _{\mathcal{P}}$-energy cubes
defined above. By the construction above, the maximal elements in $\mathcal{L
}$ are the maximal cubes in $\Pi _{1}^{{below}}\mathcal{P}\cap 
\mathcal{C}_{A}'$. For $L\in \mathcal{L}$, denote by $
\mathcal{C}_{L}^{\mathcal{L}}$ the \emph{corona} associated with $L$ in the
tree $\mathcal{L}$,
\begin{equation*}
\mathcal{C}_{L}^{\mathcal{L}}\equiv \left\{ K\in \mathcal{D}:K\subset L\text{
and there is no }L^{\prime }\in \mathcal{L}\text{ with }K\subset L^{\prime
}\subsetneqq L\right\} ,
\end{equation*}
and define the $\flat $\emph{
shifted} $\mathcal{L}$-corona by
\begin{eqnarray*}
\mathcal{C}_{L}^{\mathcal{L},\flat {shift}} 
&\equiv&
\left\{ J\in 
\mathcal{G}:J^{\flat }\in \mathcal{C}_{L}^{\mathcal{L}}\text{ }\right\} .
\end{eqnarray*}

Now the parameter $m$
in $\mathcal{L}_{m}$ refers to the level at which the stopping construction
was performed, but for\thinspace $L\in \mathcal{L}_{m}$, the corona children 
$L^{\prime }$ of $L$ are \emph{not} all necessarily in $\mathcal{L}_{m-1}$,
but may be in $\mathcal{L}_{m-t}$ for $t$ large.

At this point we introduce the notion of geometric depth $d$ in the tree $
\mathcal{L}$ by defining
\begin{eqnarray}
\mathcal{G}_{0} &\equiv &\left\{ L\in \mathcal{L}:L\text{ is maximal}
\right\} ,  \label{geom depth} \\
\mathcal{G}_{1} &\equiv &\left\{ L\in \mathcal{L}:L\text{ is maximal wrt }
L\subsetneqq L_{0}\text{ for some }L_{0}\in \mathcal{G}_{0}\right\} ,  \notag
\\
&&\vdots  \notag \\
\mathcal{G}_{d+1} &\equiv &\left\{ L\in \mathcal{L}:L\text{ is maximal wrt }
L\subsetneqq L_{d}\text{ for some }L_{d}\in \mathcal{G}_{d}\right\} ,  \notag
\\
&&\vdots  \notag
\end{eqnarray}
We refer to $\mathcal{G}_{d}$ as the $d^{th}$ generation of cubes in the
tree $\mathcal{L}$, and say that the cubes in $\mathcal{G}_{d}$ are at
depth $d$ in the tree $\mathcal{L}$ (the generations $\mathcal{G}_{d}$ here
are \emph{not} related to the grid $\mathcal{G}$), and we write $d_{{
geom}}\left( L\right) $ for the geometric depth of $L$. Thus the cubes
in $\mathcal{G}_{d}$ are the stopping cubes in $\mathcal{L}$ that are $d$
levels in the \emph{geometric} sense below the top level. While the
geometric depth $d_{{geom}}$ is about to be superceded by the
`indented' depth $d_{{indent}}$ defined in the next subsection, we
will return to the geometric depth in order to iterate Lacey's bottom/up
stopping criterion when proving the second line in (\ref{rest bounds}) in
Proposition \ref{bottom up 3} below.

\subsection{The indented corona construction}

Now we address the lack of goodness in $\Pi _{1}^{{below}}\mathcal{P}
\cap \mathcal{C}_{A}'$. For this we introduce an
additional top/down stopping time $\mathcal{H}$ over the collection $
\mathcal{L}$. Given the initial generation 
\begin{equation*}
\mathcal{H}_{0}=\left\{ \text{maximal }L\in \mathcal{
L}\right\} =\left\{ \text{maximal }I\in \Pi _{1}^{{below}}\mathcal{P}
\right\} ,
\end{equation*}
define subsequent generations $\mathcal{H}_{k}$ as follows. For $k\geq 1$
and each $H\in \mathcal{H}_{k-1}$, let 
\begin{equation*}
\mathcal{H}_{k}\left( H\right) \equiv \left\{ \text{maximal }L\in \mathcal{L}
:3L\subset H\right\}
\end{equation*}
consist of the next $\mathcal{H}$-generation of $\mathcal{L}$-cubes
below $H$, and set $\mathcal{H}_{k}\equiv \bigcup\limits_{H\in \mathcal{H}
_{k-1}}\mathcal{H}_{k}\left( H\right) $. Finally set $\mathcal{H}\equiv
\bigcup\limits_{k=0}^{\infty }\mathcal{H}_{k}$. We refer to the stopping
cubes $H\in \mathcal{H}$ as \emph{indented} stopping cubes since $
3H\subset \pi _{\mathcal{H}}H$ for all $H\in \mathcal{H}$ at indented
generation one or more, i.e. each successive such $H$ is `indented' in its $
\mathcal{H}$-parent. This property of indentation is precisely what is
required in order to generate geometric decay in indented generations at the
end of the proof. We refer to $k$ as the \emph{indented depth} of the
stopping cube $H\in \mathcal{H}_{k}$, written $k=d_{{indent}
}\left( H\right) $, which is a refinement of the geometric depth $d_{
{geom}}$ introduced above. We will often revert to writing the dummy
variable for cubes in $\mathcal{H}$ as $L$ instead of $H$. For $L\in 
\mathcal{H}$ define the $\mathcal{H}$-corona $\mathcal{C}_{L}^{\mathcal{H}}$ and $\mathcal{H}$-$\flat $shifted corona $
\mathcal{C}_{L}^{\mathcal{H},\flat {shift}}$ by
\begin{eqnarray*}
\mathcal{C}_{L}^{\mathcal{H}} &\equiv &\left\{ I\in \mathcal{D}:I\subset L
\text{ and }I\not\subset L^{\prime }\text{ for any }L^{\prime }\in \mathfrak{
C}_{\mathcal{H}}\left( L\right) \right\} , \\
\mathcal{C}_{L}^{\mathcal{H},\flat {shift}} &\equiv &\left\{ J\in 
\mathcal{G}:J^{\flat }\in \mathcal{C}_{L}^{\mathcal{H}}\right\} .
\end{eqnarray*}
We will also need recourse to the coronas $\mathcal{C}_{L}^{\mathcal{H}}$
restricted to cubes in $\mathcal{L}$, i.e.
\begin{equation*}
\mathcal{C}_{L}^{\mathcal{H}}\left( \mathcal{L}\right) \equiv \mathcal{C}
_{L}^{\mathcal{H}}\cap \mathcal{L}=\left\{ T\in \mathcal{L}:T\subset L\text{
and }T\not\subset L^{\prime }\text{ for any }L^{\prime }\in \mathcal{H}\text{
with }L^{\prime }\subsetneqq L\right\} .
\end{equation*}
and 
\begin{equation*}
\mathcal{T}(L)\equiv \mathcal{C}_{L}^{\mathcal{H},{restrict}}(\mathcal{L})=\mathcal{C}_{L}^{\mathcal{H}}(\mathcal{L})\backslash\{L\} 
\end{equation*}
We emphasize the
distinction `indented generation' as this refers to the indented depth
rather than either the level of initial stopping construction of $\mathcal{L}
$, or the geometric depth. The point of introducing the tree $\mathcal{H}$
of indented stopping cubes, is that the inclusion $3L\subset \pi _{
\mathcal{H}}L$ for all $L\in \mathcal{H}$ with $d_{{indent}}\left(
L\right) \geq 1$ turns out to be an adequate substitute for the standard
`goodness' lost in the process of infusing the weak goodness of Hyt\"{o}nen
and Martikainen in\ Subsection \ref{Subsec HM} above.

\subsubsection{Flat shifted coronas}

We now define the $\flat $shifted admissible
collections of pairs $\mathcal{P}_{L,t}^{\flat \mathcal{H}}$ using the coronas
\begin{equation*}
\mathcal{C}_{L}^{\mathcal{H},\flat {shift}}\equiv \left\{ J\in \Pi
_{2}\mathcal{P}:J^{\flat }\in \mathcal{C}_{L}^{\mathcal{H}}\right\} \text{
and }\mathcal{C}_{L}^{\mathcal{L},\flat {shift}}\equiv \left\{ J\in
\Pi _{2}\mathcal{P}:J^{\flat }\in \mathcal{C}_{L}^{\mathcal{L}}\right\} .
\end{equation*}
In these flat shifted $\mathcal{H}$ and $\mathcal{L}$ coronas, we have
effectively shift the cubes $J$ two levels `up' by requiring $J^{\flat
}\in \mathcal{C}_{L}^{\mathcal{L}}$, but because $\mathcal{P}$ is admissible, we
always have $J^{\maltese }\in \mathcal{C}_{A}^{\mathcal{A},{restrict}}$. We define 
\begin{eqnarray*}
\mathcal{P}_{L,t}^{\flat \mathcal{H}} &\equiv &\left\{ \left( I,J\right) \in 
\mathcal{P}:I\in \mathcal{C}_{L}^{\mathcal{H}}, J\in \mathcal{C}
_{L^{\prime }}^{\mathcal{H},\flat {shift}}\text{ for some }L^{\prime
}\in \mathcal{H}_{d_{{indent}}\left( L\right) +t}
L^{\prime }\subset L\right\} , \\
\mathcal{P}_{L,0}^{\flat \mathcal{H}} &=&\left\{ \left( I,J\right) \in 
\mathcal{P}:I\in \mathcal{C}_{L}^{\mathcal{H}}\text{ and }J\in \mathcal{C}
_{L}^{\mathcal{H},\flat {shift}}\right\} 
\end{eqnarray*}
and
\begin{eqnarray*}
\mathcal{P}_{L,0}^{\flat \mathcal{H}} &=&\mathcal{P}_{L,0}^{\flat \mathcal{H}
-small}\dot{\cup}\mathcal{P}_{L,0}^{\flat \mathcal{H}-big}; \\
\mathcal{P}_{L,0}^{\flat \mathcal{H}-small}
&\equiv &
\left\{ \left(
I,J\right) \in \mathcal{P}_{L,0}^{\flat \mathcal{H}}:\text{there is no }
L^{\prime }\in \mathcal{T}\left( L\right) \text{ with }J^{\flat }\subset
L^{\prime }\subset I\right\} \\
&=&
\left\{ \left( I,J\right) \in \mathcal{P}_{L,0}^{\flat \mathcal{H}}:I\in 
\mathcal{C}_{L^{\prime }}^{\mathcal{L}}\backslash \left\{ L^{\prime }\right\}, J\in \mathcal{C}_{L^{\prime }}^{\mathcal{L},\flat {shift}
}\text{ for some }L^{\prime }\in \mathcal{T}\left( L\right) \right\} , \\
\mathcal{P}_{L,0}^{\flat \mathcal{H}-big} &\equiv &\left\{ \left( I,J\right)
\in \mathcal{P}_{L,0}^{\flat \mathcal{H}}:\text{there is }L^{\prime }\in 
\mathcal{T}\left( L\right) \text{ with }J^{\flat }\subset L^{\prime }\subset
I\right\} ,
\end{eqnarray*}
with one exception: if $L\in \mathcal{H}_{0}$ we set $
\mathcal{P}_{L,0}^{\flat \mathcal{H}-small}\equiv \mathcal{P}_{L,0}^{\flat 
\mathcal{H}}$ and $\mathcal{P}_{L,0}^{\flat \mathcal{H}-big}\equiv \emptyset 
$ since in this case $L$ fails to satisfy (\ref{up stopping condition}) as
pointed out above. Finally, for $L\in \mathcal{H}$ we further decompose $
\mathcal{P}_{L,0}^{\flat \mathcal{H}-small}$ as
\begin{eqnarray*}
\mathcal{P}_{L,0}^{\flat \mathcal{H}-small} 
&=&
\overset{\cdot }{\bigcup 
_{L^{\prime }\in \mathcal{T}\left( L\right)}}\mathcal{P}_{L^{\prime
},0}^{\flat \mathcal{L}-small} \\
\text{where }\mathcal{P}_{L^{\prime },0}^{\flat \mathcal{L}-small} 
&\equiv &
\left\{
\left( I,J\right) \in \mathcal{P}:I\in \mathcal{C}_{L^{\prime }}^{\mathcal{L}
}\backslash \left\{ L^{\prime }\right\} \text{ and }J\in \mathcal{C}
_{L^{\prime }}^{\mathcal{L},\flat {shift}}\right\} 
\end{eqnarray*}
Then we set
\begin{eqnarray}
\mathcal{P}^{\flat big} &\equiv &\left\{ \bigcup\limits_{L\in \mathcal{H}}
\mathcal{P}_{L,0}^{\flat \mathcal{H}-big}\right\} \bigcup \left\{
\bigcup\limits_{t\geq 1}\bigcup\limits_{L\in \mathcal{H}}\mathcal{P}
_{L,t}^{\flat \mathcal{H}}\right\} ;  \label{def big small flat} \\
\mathcal{P}^{\flat small} &\equiv &\bigcup\limits_{L\in \mathcal{L}}
\mathcal{P}_{L,0}^{\flat \mathcal{L}-small}\text{ } \notag
\end{eqnarray}
We observed above that every pair $\left( I,J\right) \in \mathcal{P}$ is
included in either $\mathcal{P}^{small}$ or $\mathcal{P}^{big}$, and it
follows that every pair $\left( I,J\right) \in \mathcal{P}$ is thus included
in either $\mathcal{P}^{\flat small}$ or $\mathcal{P}^{\flat big}$, simply
because the pairs $\left( I,J\right) $ have been shifted up by two dyadic
levels in the cube $J$. Thus the coronas $\mathcal{P}_{L,0}^{\flat 
\mathcal{L}-small}$ are now even \emph{smaller} than the regular coronas $
\mathcal{P}_{L,0}^{\mathcal{L}-small}$, which permits the estimate (\ref
{small claim' 3}) below to hold for the larger augmented size functional. On
the other hand, the coronas $\mathcal{P}_{L,0}^{\flat \mathcal{H}-big}$ and $
\mathcal{P}_{L,t}^{\flat \mathcal{H}}$ are now bigger than before, requiring
the stronger straddling lemmas above in order to obtain the estimates (\ref
{rest bounds}) below. More specifically, we will see that stopping forms
with pairs in $\mathcal{P}^{\flat big}$ will be estimated using the $\flat $
Straddling and Substraddling Lemmas (Substraddling applies to part of $
\mathcal{P}_{L,0}^{\flat \mathcal{H}-big}$ and $\flat $Straddling applies to
the remaining part of $\mathcal{P}_{L,0}^{\flat \mathcal{H}-big}$ and to $
\mathcal{P}_{L,t}^{\flat \mathcal{H}}$), and it is here that the removal of
the corona-straddling collection $\mathcal{P}_{{cor}}^{A}$ is
essential, while forms with pairs in $\mathcal{P}^{\flat small}$ will be
absorbed.

\subsection{Size estimates}

Now we turn to proving the \emph{size estimates} we need for these
collections. Recall that the \emph{restricted} norm $\widehat{\mathfrak{N}}_{
{stop},\bigtriangleup ^{\omega }}^{A,\mathcal{P}}$ is the best
constant in the inequality
\begin{equation*}
\left\vert \mathsf{B}\right\vert _{{stop},\bigtriangleup ^{\omega
}}^{A,\mathcal{P}}\left( f,g\right) \leq \widehat{\mathfrak{N}}_{{
stop},\bigtriangleup ^{\omega }}^{A,\mathcal{P}}\left\Vert \mathsf{P}_{\Pi
_{1}\mathcal{P}}^{\sigma ,\mathbf{b}}f\right\Vert _{L^{2}\left( \sigma
\right) }^{\bigstar }\left\Vert \mathsf{P}_{\Pi _{2}\mathcal{P}}^{\omega ,
\mathbf{b}^{\ast }}g\right\Vert _{L^{2}\left( \omega \right) }^{\bigstar }
\end{equation*}
where $f\in L^{2}\left( \sigma \right) $ satisfies $E_{I}^{\sigma
}\left\vert f\right\vert \leq \alpha _{\mathcal{A}}\left( A\right) $ for all 
$I\in \mathcal{C}_{A}$, and $g\in L^{2}\left( \omega \right) $.

\begin{prop}
\label{bottom up 3}Suppose $\rho $ in (\ref{def rho}) is greater than $1$,
and $\mathcal{P}$ is a \emph{reduced admissible} collection of pairs for a
dyadic cube $A$. Let $\mathcal{P}=\mathcal{P}^{\flat big}\dot{\cup}
\mathcal{P}^{\flat small}$ be the decomposition satisfying above, i.e.
\begin{equation*}
\mathcal{P}=\left\{ \bigcup\limits_{L\in \mathcal{H}}\mathcal{P}
_{L,0}^{\flat \mathcal{H}-big}\right\} \bigcup \left\{
\bigcup\limits_{t\geq 1}\bigcup\limits_{L\in \mathcal{H}}\mathcal{P}
_{L,t}^{\flat \mathcal{H}}\right\} \ \bigcup \ \left( \bigcup_{L\in \mathcal{L}
}\mathcal{P}_{L,0}^{\flat \mathcal{L}-small}\right)
\end{equation*}
Then all of these collections $\mathcal{P}_{L,0}^{\flat \mathcal{L}-small}$, 
$\mathcal{P}_{L,0}^{\flat \mathcal{H}-big}$ and $\mathcal{P}_{L,t}^{\flat 
\mathcal{H}}$ are reduced admissible, and we have the estimate 
\begin{equation}
\mathcal{S}_{{aug}{size}}^{\alpha ,A}\left( \mathcal{P}
_{L,0}^{\flat \mathcal{L}-small}\right) ^{2}\leq \left( \rho -1\right) 
\mathcal{S}_{{aug}{size}}^{\alpha ,A}\left( \mathcal{P}
\right) ^{2},\ \ \ \ \ L\in \mathcal{L} \label{small claim' 3}
\end{equation}
and the localized norm bounds,
\begin{eqnarray}
\widehat{\mathfrak{N}}_{{stop},\bigtriangleup ^{\omega
}}^{A,\bigcup\limits_{L\in \mathcal{H}}\mathcal{P}_{L,0}^{\flat \mathcal{H}
-big}} &\leq &C\mathcal{S}_{{aug}{size}}^{\alpha ,A}\left( 
\mathcal{P}\right) ,  \label{rest bounds} \\
\widehat{\mathfrak{N}}_{{stop},\bigtriangleup ^{\omega
}}^{A,\bigcup\limits_{L\in \mathcal{H}}\mathcal{P}_{L,t}^{\flat \mathcal{H}
}} &\leq &C\rho ^{-\frac{t}{2}}\mathcal{S}_{{aug}{size}
}^{\alpha ,A}\left( \mathcal{P}\right) ,\ \ \ \ \ t\geq 1.  \notag
\end{eqnarray}
\end{prop}

Using this proposition on size estimates, we can finish the proof of (\ref
{First inequality}), and hence the proof of (\ref{B stop form 3}).

\begin{cor}
The sublinear stopping form inequality (\ref{First inequality}) holds.
\end{cor}

\begin{proof}
Recall that $\widehat{\mathfrak{N}}_{{stop},\bigtriangleup ^{\omega
}}^{A,\mathcal{P}}$ is the best constant in the inequality (\ref{best hat}). 
Since $\left\{ \mathcal{P}_{L,0}^{\flat \mathcal{L}-small}\right\} _{L\in 
\mathcal{L}}$ is a mutually orthogonal family of $A$-admissible pairs, the
Orthogonality Lemma \ref{mut orth} implies that
\begin{equation*}
\widehat{\mathfrak{N}}_{{stop},\bigtriangleup ^{\omega
}}^{A,\bigcup\limits_{L\in \mathcal{L}}\mathcal{P}_{L,0}^{\flat \mathcal{L}
-small}}\leq \sup_{L\in \mathcal{L}}\widehat{\mathfrak{N}}_{{stop}
,\bigtriangleup ^{\omega }}^{A,\mathcal{P}_{L,0}^{\flat \mathcal{L}-small}}
\end{equation*}
Using this, together with the decomposition of $\mathcal{P}$ and (\ref{rest
bounds}) above, we obtain
\begin{eqnarray*}
\widehat{\mathfrak{N}}_{{stop},\bigtriangleup ^{\omega }}^{A,
\mathcal{P}} &\leq &\sup_{L\in \mathcal{H}}\widehat{\mathfrak{N}}_{{
stop},\bigtriangleup ^{\omega }}^{A,\bigcup\limits_{L\in \mathcal{H}}
\mathcal{P}_{L,0}^{\flat \mathcal{H}-big}}+\sum_{t=1}^{M+1}\sup_{L\in 
\mathcal{H}}\widehat{\mathfrak{N}}_{{stop},\bigtriangleup ^{\omega
}}^{A,\bigcup\limits_{L\in \mathcal{H}}\mathcal{P}_{L,t}^{\flat \mathcal{H}
}}+\widehat{\mathfrak{N}}_{{stop},\bigtriangleup ^{\omega
}}^{A,\bigcup\limits_{L\in \mathcal{L}}\mathcal{P}_{L,0}^{\flat \mathcal{L}
-small}} \\
&\lesssim &\mathcal{S}_{{aug}{size}}^{\alpha ,A}\left( 
\mathcal{P}\right) +\left( \sum_{t=1}^{M+1}\rho ^{-\frac{t}{2}}\right) 
\mathcal{S}_{{aug}{size}}^{\alpha ,A}\left( \mathcal{P}
\right) +\sup_{L\in \mathcal{L}}\widehat{\mathfrak{N}}_{{stop}
,\bigtriangleup ^{\omega }}^{A,\mathcal{P}_{L,0}^{\flat \mathcal{L}-small}}\
\end{eqnarray*}
Since the admissible collection $\mathcal{P}^{A}$ in (\ref{initial P}) that
arises in the stopping form is finite, we can define $\mathfrak{L}$ to be
the best constant in the inequality
\begin{equation*}
\widehat{\mathfrak{N}}_{{stop},\bigtriangleup ^{\omega }}^{A,
\mathcal{P}}\leq \mathfrak{L}\mathcal{S}_{{aug}{size}
}^{\alpha ,A}\left( \mathcal{P}\right) \text{ for all }A\text{-admissible
collections }\mathcal{P}.
\end{equation*}
Now choose $\mathcal{P}$ so that 
\begin{equation*}
\frac{\widehat{\mathfrak{N}}_{{stop},\bigtriangleup ^{\omega }}^{A,
\mathcal{P}}}{\mathcal{S}_{{aug}{size}}^{\alpha ,A}\left( 
\mathcal{P}\right) }>\frac{1}{2}\mathfrak{L=}\frac{1}{2}\sup_{\mathcal{Q}
\text{ is }A\text{-admissible}}\frac{\widehat{\mathfrak{N}}_{{stop}
,\bigtriangleup ^{\omega }}^{A,\mathcal{Q}}}{\mathcal{S}_{{aug}
{size}}^{\alpha ,A}\left( \mathcal{Q}\right) }\ .
\end{equation*}
Then using $\displaystyle \sum_{t=1}^{M+1}\rho ^{-\frac{t}{2}}\leq \frac{1}{\sqrt{\rho }-1}
$ we have
\begin{eqnarray*}
\mathfrak{L} &<&2\frac{\widehat{\mathfrak{N}}_{{stop},\bigtriangleup
^{\omega }}^{A,\mathcal{P}}}{\mathcal{S}_{{aug}{size}
}^{\alpha ,A}\left( \mathcal{P}\right) }\leq \frac{C\frac{1}{\sqrt{\rho }-1}
\mathcal{S}_{{aug}{size}}^{\alpha ,A}\left( \mathcal{P}
\right) +  \displaystyle
C\sup_{L\in \mathcal{L}}\widehat{\mathfrak{N}}_{{stop}
,\bigtriangleup ^{\omega }}^{A,\mathcal{P}_{L,0}^{\flat \mathcal{L}-small}}}{
\mathcal{S}_{{aug}{size}}^{\alpha ,A}\left( \mathcal{P}
\right) } \\
&\leq &
C\frac{1}{\sqrt{\rho }-1}+C\sup_{L\in \mathcal{L}}\mathfrak{L}\frac{
\mathcal{S}_{{aug}{size}}^{\alpha ,A}\left( \mathcal{P}
_{L,0}^{\flat \mathcal{L}-small}\right) }{\mathcal{S}_{{aug}{size}}^{\alpha ,A}\left( \mathcal{P}\right) }\leq C\frac{1}{\sqrt{\rho }-1}+C
\mathfrak{L}\sqrt{\rho -1}\ 
\end{eqnarray*}
where we have used (\ref{small claim' 3}) in the last line. If we choose $
\rho >1$ so that 
\begin{equation}
C\sqrt{\rho -1}<\frac{1}{2},  \label{choose rho}
\end{equation}
then we obtain $\mathfrak{L}\leq 2C\frac{1}{\sqrt{\rho }-1}$. Together with
Lemma \ref{energy control}, this yields
\begin{equation*}
\widehat{\mathfrak{N}}_{{stop},\bigtriangleup ^{\omega }}^{A,
\mathcal{P}}\leq \mathfrak{L}\mathcal{S}_{{aug}{size}
}^{\alpha ,A}\left( \mathcal{P}\right) \leq 2C\frac{1}{\sqrt{\rho }-1}\left( 
\mathcal{E}_{2}^{\alpha }+\sqrt{\mathfrak{A}_{2}^{\alpha }}\right)
\end{equation*}
as desired, and completes the proof of inequality (\ref{First inequality}).
\end{proof}

Thus, in view of Conclusion \ref{assume}, it remains only to prove
Proposition \ref{bottom up 3} using the Orthogonality and Straddling and
Substraddling Lemmas above, and we now turn to this task.

\begin{proof}[Proof of Proposition \protect\ref{bottom up 3}]
We split the proof into three parts.

\textbf{Proof of (\ref{small claim' 3})}: To prove the inequality (\ref
{small claim' 3}), suppose first that $L\notin \mathcal{L}_{M+1}$. In the
case that $L\in \mathcal{L}_{0}$ is an initial generation cube, then
from (\ref{key property 3}) and the fact that every $I\in \mathcal{P}
_{L,0}^{\flat \mathcal{L}-small}$ satisfies $I\subsetneqq L$, we obtain that
\begin{eqnarray*}
\mathcal{S}_{{aug}{size}}^{\alpha ,A}\left( \mathcal{P}_{L,0}^{\flat \mathcal{L}-small}\right) ^{2}
&=&
\sup_{K^{\prime }\in \Pi _{1}^{{below}}\mathcal{P}_{L,0}^{\flat \mathcal{L}-small}\cap \mathcal{C}_{A}'}\frac{\Psi ^{\alpha }\left( K^{\prime };\mathcal{P}
_{L,0}^{\flat \mathcal{L}-small}\right) ^{2}}{\left\vert K^{\prime
}\right\vert _{\sigma }} \\
&\leq &
\sup_{K^{\prime }\in \Pi _{1}^{{below}}\mathcal{P}\cap 
\mathcal{C}_{A}':\ K^{\prime }\varsubsetneqq L}\frac{
\Psi ^{\alpha }\left( K^{\prime };\mathcal{P}_{L,0}^{\flat \mathcal{L}
-small}\right) ^{2}}{\left\vert K^{\prime }\right\vert _{\sigma }}\\
&\leq&
\varepsilon \mathcal{S}_{{aug}{size}}^{\alpha ,A}\left( 
\mathcal{P}\right) ^{2}
\end{eqnarray*}
Now suppose that $L\not\in \mathcal{L}_{0}$ in addition to $L\notin \mathcal{
L}_{M+1}$. Pick a pair $\left( I,J\right) \in \mathcal{P}_{L,0}^{\flat 
\mathcal{L}-small}$. Then $I$ is in the restricted corona $\mathcal{C}_{L}^{
\mathcal{L},'}$ and $J$ is in the $\flat $\emph{shifted}
corona $\mathcal{C}_{L}^{\mathcal{L},\flat {shift}}$. Since $\mathcal{P}_{L,0}^{\flat \mathcal{L}-small}$ is a finite collection, the
definition of $\mathcal{S}_{{aug}{size}}^{\alpha ,A}\left( 
\mathcal{P}_{L,0}^{\flat \mathcal{L}-small}\right) $ shows that there is an
cube $K\in \Pi _{1}^{{below}}\mathcal{P}_{L,0}^{\flat \mathcal{L}
-small}\cap \mathcal{C}_{A}'$ so that
\begin{equation*}
\mathcal{S}_{{aug}{size}}^{\alpha ,A}\left( \mathcal{P}
_{L,0}^{\flat \mathcal{L}-small}\right) ^{2}=\frac{1}{\left\vert
K\right\vert _{\sigma }}\left( \frac{\mathrm{P}^{\alpha }\left( K,\mathbf{1}
_{A\backslash K}\sigma \right) }{\left\vert K\right\vert^\frac{1}{n} }\right) ^{2}\omega
_{\flat \mathcal{P}}\left( \mathbf{T}\left( K\right) \right) .
\end{equation*}
Note that $K\subsetneqq L$ by definition of $\mathcal{P}_{L,0}^{\flat 
\mathcal{L}-small}$. Now let $t$ be such that $L\in \mathcal{L}_{t}$, and
define 
\begin{equation*}
t^{\prime }=t^{\prime }\left( K\right) \equiv \max \left\{ s:\text{there is }
L^{\prime }\in \mathcal{L}_{s}\text{ with }L^{\prime }\subset K\right\} ,
\end{equation*}
and note that $0\leq t^{\prime }<t$.  
First, suppose that $t^{\prime }=0$ so that $%
K $ does not contain any $L^{\prime }\in \mathcal{L}$. Then it follows from
the construction at level $\ell =0$ that%
\begin{equation*}
\frac{1}{\left\vert K\right\vert _{\sigma }}\left( \frac{\mathrm{P}^{\alpha
}\left( K,\mathbf{1}_{A\setminus K}\sigma \right) }{\left\vert K\right\vert }%
\right) ^{2}\omega _{\flat \mathcal{P}}\left( \mathbf{T}\left( K\right)
\right) <\varepsilon \mathcal{S}_{{aug}{size}}^{\alpha
,A}\left( \mathcal{P}\right) ^{2},
\end{equation*}%
and hence from $\rho =1+\varepsilon $ we obtain 
\begin{equation*}
\mathcal{S}_{{aug}{size}}^{\alpha ,A}\left( \mathcal{P}%
_{L,0}^{\flat \mathcal{L}-small}\right) ^{2}<\varepsilon \mathcal{S}_{%
{aug}{size}}^{\alpha ,A}\left( \mathcal{P}\right)
^{2}=\left( \rho -1\right) \mathcal{S}_{{aug}{size}}^{\alpha
,A}\left( \mathcal{P}\right) ^{2}.
\end{equation*}%
Now suppose that $t^{\prime }\geq 1$.
Then $K$ fails the stopping condition (
\ref{up stopping condition}) with $m=t^{\prime }+1$, since otherwise it
would contain a cube $L^{\prime \prime }\in \mathcal{L}_{t^{\prime }+1}$
contradicting our definition of $t^{\prime }$, and so
\begin{equation*}
\omega _{\flat \mathcal{P}}\left( \mathbf{T}\left( K\right) \right) <\rho
\omega _{\flat \mathcal{P}}\left( \mathbf{V}\left( K\right) \right) \text{
where }\mathbf{V}\left( K\right) \equiv\!\!\!\!\!\!\! \bigcup\limits_{L^{\prime }\in
\bigcup\limits_{\ell =0}^{t^{\prime }}\mathcal{L}_{\ell }:\ L^{\prime
}\subset K}\mathbf{T}\left( L^{\prime }\right) .
\end{equation*}

Now we use the crucial fact that the positive measure $\omega _{\flat 
\mathcal{P}}$ is \emph{additive} and finite to obtain from this that
\begin{equation}
\omega _{\flat \mathcal{P}}\left( \mathbf{T}\left( K\right) \backslash 
\mathbf{V}\left( K\right) \right) =\omega _{\flat \mathcal{P}}\left( \mathbf{
T}\left( K\right) \right) -\omega _{\flat \mathcal{P}}\left( \mathbf{V}
\left( K\right) \right) \leq \left( \rho -1\right) \omega _{\flat \mathcal{P}
}\left( \mathbf{V}\left( K\right) \right) .  \label{additive}
\end{equation}
Now recall that
\begin{equation*}
\mathcal{S}_{{aug}{size}}^{\alpha ,A}\left( \mathcal{Q}
\right) ^{2}\equiv \sup_{K\in \Pi _{1}^{{below}}\mathcal{Q}\cap
C_{A}^{'}}\frac{1}{\left\vert K\right\vert _{\sigma }}
\left( \frac{\mathrm{P}^{\alpha }\left( K,\mathbf{1}_{A\backslash K}\sigma
\right) }{\left\vert K\right\vert^\frac{1}{n} }\right) ^{2}\left\Vert \mathsf{Q}_{\Pi
_{2}^{K,{aug}}\mathcal{Q}}^{\omega ,\mathbf{b}^{\ast }}x\right\Vert
_{L^{2}\left( \omega \right) }^{\spadesuit 2}.
\end{equation*}
We claim it follows that for each $J\in \Pi _{2}^{K,{aug}}\mathcal{P}
_{L,0}^{\flat \mathcal{L}-small}$ the support $\left( c_{J^{\flat }},\ell
\left( J^{\flat }\right) \right) $ of the atom $\delta _{\left( c_{J^{\flat
}},\ell \left( J^{\flat }\right) \right) }$ is contained in the set $\mathbf{
T}\left( K\right) $, but not in the set 
\begin{equation*}
\mathbf{V}\left( K\right) \equiv \bigcup \left\{ \mathbf{T}\left( L^{\prime
}\right) :L^{\prime }\in \bigcup\limits_{\ell =0}^{t^{\prime }}\mathcal{L}
_{\ell }:\ L^{\prime }\subset K\right\} .
\end{equation*}
Indeed, suppose in order to derive a contradiction, that $\left( c_{J^{\flat
}},\ell \left( J^{\flat }\right) \right) \in \mathbf{T}\left( L^{\prime
}\right) $ for some $L^{\prime }\in \mathcal{L}_{\ell }$ with $0\leq \ell
\leq t^{\prime }$. Recall that $L\in \mathcal{L}_{t}$ with $t^{\prime }<t$
so that $L^{\prime }\subsetneqq L$. Thus $\left( c_{J^{\flat }},\ell \left(
J^{\flat }\right) \right) \in \mathbf{T}\left( L^{\prime }\right) $ implies 
$J^{\flat }\subset L^{\prime }$, which contradicts the fact that 
\begin{equation*}
J\!\in\! \Pi _{2}^{K}\mathcal{P}_{L,0}^{\flat \mathcal{L}-small}\subset \Pi _{2}
\mathcal{P}_{L,0}^{\flat \mathcal{L}-small}=\left\{ \left( I,J\right) \in 
\mathcal{P}:I\in \mathcal{C}_{L}^{\mathcal{L}}\backslash \left\{ L\right\} 
\text{ and }J\in \mathcal{C}_{L}^{\mathcal{L},\flat {shift}}\right\}
\end{equation*}
implies $J^{\flat }\in \mathcal{C}_{L}^{\mathcal{L}}$ - because $L^{\prime
}\notin \mathcal{C}_{L}^{\mathcal{L}}$.

Thus from the definition of $\omega _{\flat \mathcal{P}}$ in (\ref{def
atomic}), the `energy' $\left\Vert \mathsf{Q}_{\Pi _{2}^{K,{aug}}
\mathcal{P}_{L,0}^{\flat \mathcal{L}-small}}^{\omega ,\mathbf{b}^{\ast
}}x\right\Vert _{L^{2}\left( \omega \right) }^{\spadesuit 2}$ is at most the 
$\omega _{\flat \mathcal{P}}$-measure of $\mathbf{T}\left( K\right)
\backslash \mathbf{V}\left( K\right) $. Using now
\begin{equation*}
\omega _{\flat \mathcal{P}_{L,0}^{\flat \mathcal{L}-small}}\left( \mathbf{T}
\left( K\right) \right)
=
\omega _{\flat \mathcal{P}_{L,0}^{\flat \mathcal{L}-small}}\left( \mathbf{T}
\left( K\right)\backslash \mathbf{V}\left( K\right) \right)
\leq
\omega _{\flat \mathcal{P}}\left( \mathbf{T}
\left( K\right) \backslash \mathbf{V}\left( K\right) \right) 
\end{equation*}
and (\ref{additive}), we then have
\begin{eqnarray*}
\mathcal{S}_{{aug}{size}}^{\alpha ,A}\left( \mathcal{P}
_{L,0}^{\flat \mathcal{L}-small}\right) ^{2} 
&\leq &
\sup_{K\in \Pi _{1}^{{below}}\mathcal{P}_{L,0}^{\flat 
\mathcal{L}-small}\cap \mathcal{C}_{A}'}\frac{1}{
\left\vert K\right\vert _{\sigma }}\left( \frac{\mathrm{P}^{\alpha }\left( K,
\mathbf{1}_{A\backslash K}\sigma \right) }{\left\vert K\right\vert^\frac{1}{n} }\right)
^{2}\omega _{\flat \mathcal{P}}\left( \mathbf{T}\left( K\right) \backslash 
\mathbf{V}\left( K\right) \right) \\
&\leq &
\left( \rho -1\right) \sup_{K\in \Pi _{1}^{{below}}\mathcal{P}
_{L,0}^{\flat \mathcal{L}-small}\cap \mathcal{C}_{A}'}
\frac{1}{\left\vert K\right\vert _{\sigma }}\left( \frac{\mathrm{P}^{\alpha
}\left( K,\mathbf{1}_{A\backslash K}\sigma \right) }{\left\vert K\right\vert^\frac{1}{n} }
\right) ^{2}\omega _{\flat \mathcal{P}}\left( \mathbf{V}\left( K\right)
\right) 
\end{eqnarray*}
and we can continue with 
\begin{eqnarray*}
\mathcal{S}_{{aug}{size}}^{\alpha ,A}\left( \mathcal{P}
_{L,0}^{\flat \mathcal{L}-small}\right) ^{2} 
&\leq &
\left( \rho -1\right) \sup_{K\in \Pi _{1}^{{below}}\mathcal{P}
\cap \mathcal{C}_{A}'}\frac{1}{\left\vert K\right\vert
_{\sigma }}\left( \frac{\mathrm{P}^{\alpha }\left( K,\mathbf{1}_{A\backslash
K}\sigma \right) }{\left\vert K\right\vert^\frac{1}{n} }\right) ^{2}\omega _{\flat 
\mathcal{P}}\left( \mathbf{T}\left( K\right) \right) \\
&\leq &
\left( \rho -1\right) \mathcal{S}_{{aug}{size}
}^{\alpha ,A}\left( \mathcal{P}\right) ^{2}.
\end{eqnarray*}

In the remaining case where $L\in \mathcal{L}_{M+1}$ we can include $L$ as a
testing cube $K$ and the same reasoning applies. This completes the
proof of (\ref{small claim' 3}).

\bigskip

To prove the other inequality (\ref{rest bounds}) in Proposition \ref{bottom
up 3}, we will use the $\flat$ Straddling and Substraddling Lemmas to bound the norm
of certain `straddled' stopping forms by the augmented size functional $
\mathcal{S}_{{aug}{size}}^{\alpha ,A}$, and the
Orthogonality Lemma to bound sums of `mutually orthogonal' stopping forms.
Recall that 
\begin{eqnarray*}
\mathcal{P}^{\flat big} &=&\left\{ \bigcup\limits_{L\in \mathcal{H}}
\mathcal{P}_{L,0}^{\flat \mathcal{H}-big}\right\} \bigcup \left\{
\bigcup\limits_{t\geq 1}\bigcup\limits_{L\in \mathcal{H}}\mathcal{P}
_{L,t}^{\flat \mathcal{H}}\right\} \equiv \mathcal{Q}_{0}^{\flat \mathcal{H}
-big}\bigcup \mathcal{Q}_{1}^{\flat \mathcal{H}-big}; \\
\mathcal{Q}_{0}^{\flat \mathcal{H}-big} &\equiv &\bigcup\limits_{L\in 
\mathcal{L}}\mathcal{P}_{L,0}^{\flat \mathcal{H}-big}\ ,\ \ \ \ \ \mathcal{Q}
_{1}^{\flat \mathcal{H}-big}\equiv \bigcup\limits_{t\geq 1}\mathcal{P}
_{t}^{\flat \mathcal{H}-big},\ \ \ \ \ \mathcal{P}_{t}^{\flat \mathcal{H}
-big}\equiv \bigcup\limits_{L\in \mathcal{H}}\mathcal{P}_{L,t}^{\flat 
\mathcal{H}}
\end{eqnarray*}

\bigskip

\textbf{Proof of the second line in (\ref{rest bounds})}: We first turn to
the collection
\begin{eqnarray*}
\mathcal{Q}_{1}^{\flat \mathcal{H}-big} &=&\bigcup\limits_{t\geq
1}\bigcup\limits_{L\in \mathcal{H}}\mathcal{P}_{L,t}^{\flat \mathcal{H}
}=\bigcup\limits_{t\geq 1}\mathcal{P}_{t}^{\flat \mathcal{H}-big}; \\
\mathcal{P}_{t}^{\flat \mathcal{H}-big} &\equiv &\bigcup\limits_{L\in 
\mathcal{L}}\mathcal{P}_{L,t}^{\flat \mathcal{H}}\ ,\ \ \ \ \ t\geq 1,
\end{eqnarray*}
where
\begin{equation*}
\mathcal{P}_{L,t}^{\flat \mathcal{H}}=\left\{ \left( I,J\right) \in \mathcal{
P}:I\in \mathcal{C}_{L}^{\mathcal{H}}, J\in \mathcal{C}_{L^{\prime
}}^{\mathcal{H},\flat {shift}}\text{ for some }L^{\prime }\in 
\mathcal{H}_{d_{{indent}}\left( L\right) +t}, L^{\prime
}\subset L\right\} .
\end{equation*}
We now claim that the second line in (\ref{rest bounds}) holds, i.e.
\begin{equation}
\widehat{\mathfrak{N}}_{{stop},\bigtriangleup ^{\omega }}^{A,
\mathcal{P}_{t}^{\flat \mathcal{H}-big}}\leq C\rho ^{-\frac{t}{2}}\mathcal{S}
_{{aug}{size}}^{\alpha ,A}\left( \mathcal{P}\right) ,\ \ \ \
\ t\geq 1,  \label{S big t 3}
\end{equation}
which recovers the key geometric gain obtained by Lacey in \cite{Lac},
except that here we are only gaining this decay relative to the indented
subtree $\mathcal{H}$ of the tree $\mathcal{L}$.

The case $t=1$ can be handled with relative ease since decay is not relevant
here. Indeed, $\mathcal{P}_{L,1}^{\flat \mathcal{H}}$ straddles the
collection $\mathfrak{C}_{\mathcal{H}}\left( L\right) $ of $\mathcal{H}$
-children of $L$, and so the localized $\flat $Straddling Lemma \ref
{straddle 3 ref} applies to give
\begin{equation*}
\widehat{\mathfrak{N}}_{{stop},\bigtriangleup ^{\omega }}^{A,
\mathcal{P}_{L,1}^{\flat \mathcal{H}}}\leq C\mathcal{S}_{{aug}
{size}}^{\alpha ,A}\left( \mathcal{P}_{L,1}^{\flat \mathcal{H}
}\right) \leq C\mathcal{S}_{{aug}{size}}^{\alpha ,A}\left( 
\mathcal{P}\right) ,
\end{equation*}
and then the Orthogonality Lemma \ref{mut orth} applies to give
\begin{equation*}
\widehat{\mathfrak{N}}_{{stop},\bigtriangleup ^{\omega }}^{A,
\mathcal{P}_{1}^{\flat \mathcal{H}-big}}\leq \sup_{L\in \mathcal{H}}
\mathfrak{N}_{{stop},\bigtriangleup ^{\omega }}^{A,\mathcal{P}
_{L,1}^{\flat \mathcal{H}}}\leq C\mathcal{S}_{{aug}{size}
}^{\alpha ,A}\left( \mathcal{P}\right) ,
\end{equation*}
since $\left\{ \mathcal{P}_{L,1}^{\flat \mathcal{H}}\right\} _{L\in \mathcal{
L}}$ is mutually orthogonal as $\mathcal{P}_{L,1}^{\flat \mathcal{H}}\subset 
\mathcal{C}_{L}^{\mathcal{H}}\times \mathcal{C}_{L^{\prime }}^{\mathcal{H}
,\flat {shift}}$ with $L\in \mathcal{H}_{k}$ and $L^{\prime }\in 
\mathcal{H}_{k+1}$ for indented depth $k=k\left( L\right) $. The case $t=2$
is equally easy.

Now we consider the case $t\geq 2$, where it is essential to obtain
geometric decay in $t$. We remind the reader that all of our admissible
collections $\mathcal{P}_{L,t}^{\flat \mathcal{H}}$ are \emph{reduced} by
Conclusion \ref{assume}. We again apply Lemma \ref{straddle 3 ref} to $
\mathcal{P}_{L,t}^{\flat \mathcal{H}}$ with $\mathcal{S}=\mathfrak{C}_{
\mathcal{H}}\left( L\right) $, so that for any $\left( I,J\right) \in 
\mathcal{P}_{L,t}^{\flat \mathcal{H}}$, there is $H^{\prime }\in \mathfrak{C}
_{\mathcal{H}}\left( L\right) $ with $J^{\flat }\subset H^{\prime
}\subsetneqq I\in \mathcal{C}_{L}^{\mathcal{H}}$. But this time we must use
the stronger localized bounds $\mathcal{S}_{{loc}{size}
}^{\alpha ,A;S}$ with an $S$-hole, that give
\begin{eqnarray}
\widehat{\mathfrak{N}}_{{stop},\bigtriangleup ^{\omega }}^{A,
\mathcal{P}_{L,t}^{\flat \mathcal{H}}} 
&\leq &
C\sup_{H^{\prime }\in 
\mathfrak{C}_{\mathcal{H}}\left( L\right) }\mathcal{S}_{{loc}
{size}}^{\alpha ,A;H^{\prime }}\left( \mathcal{P}_{L,t}^{\flat 
\mathcal{H}}\right) ,\ \ \ \ \ t\geq 0;  \label{t,n 3} \notag \\
\mathcal{S}_{{loc}{size}}^{\alpha ,A;H^{\prime }}\left( 
\mathcal{P}_{L,t}^{\flat \mathcal{H}}\right) ^{2} &=&\!\!\!\!\!\!\!\!\!\!
\sup_{K\in \mathcal{W}
^{\ast }\left( H^{\prime }\right) \cap \mathcal{C}_{A}'}
\frac{1}{\left\vert K\right\vert _{\sigma }}\left( \frac{\mathrm{P}^{\alpha
}\left( K,\mathbf{1}_{A\backslash H^{\prime }}\sigma \right) }{\left\vert K\right\vert^\frac{1}{n} }\right) ^{2}
\!\!\!\!\!\sum_{J\in \Pi _{2}^{K,{aug}}\mathcal{P}
_{L,t}^{\flat \mathcal{H}}}\!\!\!\!\!\!\!\!
\left\Vert \bigtriangleup _{J}^{\omega ,\mathbf{b}
^{\ast }}x\right\Vert _{L^{2}\left( \omega \right) }^{\spadesuit 2}\ 
\notag
\end{eqnarray}

It remains to show that
\begin{eqnarray}
&&\sum_{J\in \Pi _{2}^{K,{aug}}\mathcal{P}_{L,t}^{\flat \mathcal{H}
}}\left\Vert \bigtriangleup _{J}^{\omega ,\mathbf{b}^{\ast }}x\right\Vert
_{L^{2}\left( \omega \right) }^{\spadesuit 2}\leq \rho ^{-\left( t-2\right)
}\omega _{\flat \mathcal{P}}\left( \mathbf{T}\left( K\right) \right) ,
\label{rem} \\
&&\text{for}\ t\geq 2,\ K\in \mathcal{W}^{\ast }\left( H^{\prime }\right)
\cap \mathcal{C}_{A}',\ H^{\prime }\in \mathfrak{C}_{
\mathcal{H}}\left( L\right)   \notag
\end{eqnarray}
so that we then have
\begin{eqnarray*}
&&\frac{1}{\left\vert K\right\vert _{\sigma }}\left( \frac{\mathrm{P}
^{\alpha }\left( K,\mathbf{1}_{A\backslash H^{\prime }}\sigma \right) }{
\left\vert K\right\vert^\frac{1}{n} }\right) ^{2}\sum_{J\in \Pi _{2}^{K,{aug}}
\mathcal{P}_{L,t}^{\flat \mathcal{H}}}\left\Vert \bigtriangleup _{J}^{\omega
,\mathbf{b}^{\ast }}x\right\Vert _{L^{2}\left( \omega \right) }^{\spadesuit
2} \\
&\leq &\rho ^{-\left( t-2\right) }\frac{1}{\left\vert K\right\vert _{\sigma }
}\left( \frac{\mathrm{P}^{\alpha }\left( K,\mathbf{1}_{A\backslash K}\sigma
\right) }{\left\vert K\right\vert^\frac{1}{n} }\right) ^{2}\omega _{\flat \mathcal{P}
}\left( \mathbf{T}\left( K\right) \right) \leq \rho ^{-\left( t-2\right) }
\mathcal{S}_{{aug}{size}}^{\alpha ,A}\left( \mathcal{P}
\right) ^{2}
\end{eqnarray*}
by (\ref{def P stop energy' 3}), and hence conclude the required bound for $
\mathfrak{N}_{{stop},\bigtriangleup ^{\omega }}^{A,\mathcal{P}
_{L,t}^{\flat \mathcal{H}}}$, namely that
\begin{eqnarray}
&&\label{N_L}
\widehat{\mathfrak{N}}_{{stop},\bigtriangleup ^{\omega }}^{A,
\mathcal{P}_{L,t}^{\flat \mathcal{H}}} \\
&\leq&
C\!\!\! \!\!\! 
\sup_{H^{\prime }\in 
\mathfrak{C}_{\mathcal{H}}\left( L\right) }\sup_{K\in \mathcal{W}^{\ast
}\left( H^{\prime }\right) \cap \mathcal{C}_{A}'}
\! \sqrt{\frac{1}{\left\vert K\right\vert _{\sigma }}\!\! \left(\! \frac{\mathrm{P}^{\alpha
}\left( K,\mathbf{1}_{A\backslash H^{\prime }}\sigma \right) }{\left\vert
K\right\vert^\frac{1}{n} }\right) ^{2}\!\!\! \!\!\! \!
\sum_{J\in \Pi _{2}^{K,{aug}}\mathcal{P}
_{L,t}^{\flat \mathcal{H}}}\!\!\! \!\!\!\! \left\Vert \bigtriangleup _{J}^{\omega ,\mathbf{b}
^{\ast }}x\right\Vert _{L^{2}\left( \omega \right) }^{\spadesuit 2}}  \notag
\\
&\leq &C\sqrt{\rho ^{-\left( t-2\right) }}\mathcal{S}_{{aug}{
size}}^{\alpha ,A}\left( \mathcal{P}\right) =C^{\prime }\rho ^{-\frac{t}{2}}
\mathcal{S}_{{aug}{size}}^{\alpha ,A}\left( \mathcal{P}
\right) .  \notag
\end{eqnarray}

\medskip

\textbf{Remark on lack of usual goodness}: To prove (\ref{rem}), it is
essential that the cubes $H^{k+2}\in \mathcal{H}_{k+2}$ at the next
indented level down from $H^{k+1}\in \mathfrak{C}_{\mathcal{H}}\left(
L\right) $ are each contained in one of the Whitney cubes $K\in \mathcal{
W}\left( H^{k+1}\right) \cap \mathcal{C}_{A}'$ for some $
H^{k+1}\in \mathfrak{C}_{\mathcal{H}}\left( L\right) $. And this is the
reason we introduced the indented corona - namely so that $3H^{k+2}\subset
H^{k+1}$ for some $H^{k+1}\in \mathfrak{C}_{\mathcal{H}}\left( L\right) $,
and hence $H^{k+2}\subset K$ for some $K\in \mathcal{W}\left( H^{k+1}\right) 
$. In the argument of Lacey in \cite{Lac}, the corresponding cubes were
good in the usual sense, and so the above triple property was automatic.

\medskip

So we begin by fixing $K\in \mathcal{W}^{\ast }\left( H^{k+1}\right) \cap 
\mathcal{C}_{A}'$ with $H^{k+1}\in \mathfrak{C}_{
\mathcal{H}}\left( L\right) $, and noting from the above that each $J\in \Pi
_{2}^{K,{aug}}\mathcal{P}_{L,t}^{\flat \mathcal{H}}$ satisfies 
\begin{equation*}
J^{\flat }\subset H^{k+t}\subset H^{k+t-1}\subset ...\subset H^{k+2}\subset K
\end{equation*}
for $H^{k+j}\in \mathcal{H}_{k+j}$ uniquely determined by $J^{\flat }$. Thus
for $t\geq 2$ we have 
\begin{eqnarray*}
\sum_{J\in \Pi _{2}^{K,{aug}}\mathcal{P}_{L,t}^{\flat \mathcal{H}
}}\left\Vert \bigtriangleup _{J}^{\omega ,\mathbf{b}^{\ast }}x\right\Vert
_{L^{2}\left( \omega \right) }^{\spadesuit 2} 
&=&
\sum_{\substack{H^{k+t}\in \mathcal{H}
_{k+t}\\\ H^{k+t}\subset K}}\sum_{\substack{J\in \Pi _{2}^{K,{aug}}\mathcal{P}
_{L,t}^{\flat \mathcal{H}}\\ \ J^{\flat }\subset H^{k+t}}}\left\Vert
\bigtriangleup _{J}^{\omega ,\mathbf{b}^{\ast }}x\right\Vert _{L^{2}\left(
\omega \right) }^{\spadesuit 2} \\
&\leq &
\sum_{\substack{H^{k+t}\in \mathcal{H}_{k+t}\\ \ H^{k+t}\subset K}}\omega _{\flat 
\mathcal{P}}\left( \mathbf{T}\left( H^{k+t}\right) \right)
\end{eqnarray*}
In the case $t=2$ we are done since the final sum above is at most $\omega
_{\flat \mathcal{P}}\left( \mathbf{T}\left( K\right) \right) $.

Now suppose $t\geq 3$. In order to obtain geometric gain in $t$, we will
apply the stopping criterion (\ref{up stopping condition}) in the following
form,
\begin{equation}
\sum_{L^{\prime }\in \mathfrak{C}_{\mathcal{L}}\left( L_{0}\right) }\!\!\!\!\omega
_{\flat \mathcal{P}}\left( \mathbf{T}\left( L^{\prime }\right) \right)
=\omega _{\flat \mathcal{P}}\left( \bigcup\limits_{L^{\prime }\in \mathfrak{
C}_{\mathcal{L}}\left( L_{0}\right) }\mathbf{T}\left( L^{\prime }\right)
\right) \leq \frac{1}{\rho }\omega _{\flat \mathcal{P}}\left( \mathbf{T}
\left( L_{0}\right) \right) ,\ \ \ \text{for all }L_{0}\in \mathcal{L}
\label{foll form}
\end{equation}
where we have used the fact that the \emph{maximal} cubes $L^{\prime }$
in the collection
$$
\bigcup\limits_{\ell =0}^{m-1}\left\{ L^{\prime }\in 
\mathcal{L}_{\ell }:\ L^{\prime }\subset L_{0}\right\} 
$$
for $L_{0}\in \mathcal{L}_{m}$ (that appears in (\ref{up stopping condition})) are
precisely the $\mathcal{L}$-children of $L_{0}$ in the tree $\mathcal{L}$
(the cubes $L^{\prime }$ above are strictly contained in $L_{0}$ since $
\rho >1$ in (\ref{up stopping condition})), so that
\begin{equation*}
\bigcup\limits_{L^{\prime }\in \Gamma }L^{\prime
}=\bigcup\limits_{L^{\prime }\in \mathfrak{C}_{\mathcal{L}}\left(
L_{0}\right) }L^{\prime }\text{ where }\Gamma \equiv \bigcup\limits_{\ell
=0}^{m-1}\left\{ L^{\prime }\in \mathcal{L}_{\ell }:\ L^{\prime }\subset
L_{0}\right\} .
\end{equation*}

In order to apply (\ref{foll form}), we collect the pairwise disjoint
cubes $H^{k+t}\in \mathcal{H}_{k+t}$ such that$\ H^{k+t}\subset
H^{k+2}\subset K$, into groups according to which cube $L^{k^{\prime
}+t-2}\in \mathcal{G}_{k^{\prime }+t-2}$ they are contained in, where $
k^{\prime }=d_{{geom}}\left( H^{k+2}\right) $ is the geometric depth
of $H^{k+2}$ in the tree $\mathcal{L}$ introduced in (\ref{geom depth}). It
follows that each cube $H^{k+t}\in \mathcal{H}_{k+t}$ is contained in a
unique cube $L^{d_{{geom}}\left( H^{k+2}\right) +t-2}\in 
\mathcal{G}_{d_{{geom}}\left( H^{k+2}\right) +t-2}$. Thus we\ obtain
from the previous inequality that
\begin{eqnarray*}
\sum_{J\in \Pi _{2}^{K,{aug}}\mathcal{P}_{L,t}^{\flat \mathcal{H}
}}\left\Vert \bigtriangleup _{J}^{\omega ,\mathbf{b}^{\ast }}x\right\Vert
_{L^{2}\left( \omega \right) }^{\spadesuit 2}
&\leq& 
\sum_{\substack{H^{k+t}\in 
\mathcal{H}_{k+t}\\ \ H^{k+t}\subset K}}\omega _{\flat \mathcal{P}}\left( 
\mathbf{T}\left( H^{k+t}\right) \right) \\
&\leq &
\sum_{\substack{ H^{k+2}\in \mathcal{H}_{k+2}  \\ H^{k+2}\subset K}}
\sum_{\substack{ L^{k^{\prime }+t-2}\in \mathcal{G}_{k^{\prime }+t-2}\\ \
L^{k^{\prime }+t-2}\subset H^{k+2}  \\ \text{where }k^{\prime }=d_{{
geom}}\left( H^{k+2}\right) }}\!\!\!\!\!\!\!\!
\omega _{\flat \mathcal{P}}\left( \mathbf{T}
\left( L^{k^{\prime }+t-2}\right) \right) 
\end{eqnarray*}
and this last expression is equal to
\begin{eqnarray*}
&&
\sum_{\substack{ H^{k+2}\in \mathcal{H}_{k+2}  \\ H^{k+2}\subset K}}\sum 
_{\substack{ L^{k^{\prime }+t-3}\in \mathcal{G}_{k^{\prime }+t-3}\\ \L^{k^{\prime }+t-3}\subset H^{k+2}  \\ \text{where }k^{\prime }=d_{{geom}}\left( H^{k+2}\right) }}
\left\{ \sum_{\substack{ L^{k^{\prime
}+t-2}\in \mathcal{G}_{k^{\prime }+t-2}\\ \ L^{k^{\prime }+t-2}\subset
L^{k^{\prime }+t-3}  \\ \text{where }k^{\prime }=d_{{geom}}\left(
H^{k+2}\right) }}\omega _{\flat \mathcal{P}}\left( \mathbf{T}\left(
L^{k^{\prime }+t-2}\right) \right) \right\} \\
&\leq &
\sum_{\substack{ H^{k+2}\in \mathcal{H}_{k+2}  \\ H^{k+2}\subset K}}
\sum_{\substack{ L^{k^{\prime }+t-3}\in \mathcal{G}_{k^{\prime }+t-3} \\ \
L^{k^{\prime }+t-3}\subset H^{k+2}  \\ \text{where }k^{\prime }=d_{{
geom}}\left( H^{k+2}\right) }}\left\{ \frac{1}{\rho }\omega _{\flat \mathcal{
P}}\left( \mathbf{T}\left( L^{k^{\prime }+t-3}\right) \right) \right\} 
\end{eqnarray*}
where in the last line we have used (\ref{foll form}) with $
L_{0}=L^{k^{\prime }+t-3}$ on the sum in braces. We then continue (if necessary) with
\begin{eqnarray*}
\sum_{J\in \Pi _{2}^{K,{aug}}\mathcal{P}_{L,t}^{\flat \mathcal{H}
}}\left\Vert \bigtriangleup _{J}^{\omega ,\mathbf{b}^{\ast }}x\right\Vert
_{L^{2}\left( \omega \right) }^{\spadesuit 2} 
&\leq &
\frac{1}{\rho }\sum_{\substack{ H^{k+2}\in \mathcal{H}_{k+2}  \\ H^{k+2}\subset K}}\sum 
_{\substack{ L^{k^{\prime }+t-3}\in \mathcal{G}_{k^{\prime }+t-3}\\ \
L^{k^{\prime }+t-3}\subset H^{k+2}  \\ \text{where }k^{\prime }=d_{{geom}}\left( H^{k+2}\right) }}\omega _{\flat \mathcal{P}}\left( \mathbf{T}
\left( L^{k^{\prime }+t-3}\right) \right)\\
&\leq &
\frac{1}{\rho ^{2}}\sum_{\substack{ H^{k+2}\in \mathcal{H}_{k+2}  \\ 
H^{k+2}\subset K}}\sum_{\substack{ L^{k^{\prime }+t-4}\in \mathcal{G}_{k^{\prime }+t-4}\\ \ L^{k^{\prime }+t-4}\subset H^{k+2}  \\ \text{where }
k^{\prime }=d_{{geom}}\left( H^{k+2}\right) }}\omega _{\flat 
\mathcal{P}}\left( \mathbf{T}\left( L^{k^{\prime }+t-4}\right) \right)
\end{eqnarray*}
\begin{eqnarray*}
&&\vdots \\
&\leq &
\frac{1}{\rho ^{t-2}}\sum_{\substack{ H^{k+2}\in \mathcal{H}_{k+2} 
\\ H^{k+2}\subset K}}\sum_{\substack{ L^{k^{\prime }}\in \mathcal{G}
_{k^{\prime }}:\ L^{k^{\prime }}\subset H^{k+2}  \\ \text{where }k^{\prime
}=d_{{geom}}\left( H^{k+2}\right) }}\omega _{\flat \mathcal{P}
}\left( \mathbf{T}\left( L^{k^{\prime }}\right) \right) 
\end{eqnarray*}
Since $L^{k^{\prime }}\subset H^{k+2}$ implies $L^{k^{\prime }}=H^{k+2}$, we
now obtain
\begin{eqnarray*}
\sum_{J\in \Pi _{2}^{K,{aug}}\mathcal{P}_{L,t}^{\flat \mathcal{H}
}}\left\Vert \bigtriangleup _{J}^{\omega ,\mathbf{b}^{\ast }}x\right\Vert
_{L^{2}\left( \omega \right) }^{\spadesuit 2}
&\leq&
\frac{1}{\rho ^{t-2}}
\sum_{H^{k+2}\in \mathcal{H}_{k+2}:\ H^{k+2}\subset K}\omega _{\flat 
\mathcal{P}}\left( \mathbf{T}\left( H^{k+2}\right) \right)\\
 &\leq&
 \frac{1}{
\rho ^{t-2}}\omega _{\flat \mathcal{P}}\left( \mathbf{T}\left( K\right)\right) 
\end{eqnarray*}
which completes the proof of (\ref{rem}), and hence that of (\ref{N_L}).
Finally, an application of the Orthogonality Lemma \ref{mut orth}\ proves (\ref{S big t 3}).

\bigskip

\textbf{Proof of the first line in (\ref{rest bounds})}: At last we turn to
proving the first line in (\ref{rest bounds}). Recalling that $\mathcal{T}\left( L\right) =\mathcal{C}_{L}^{\mathcal{H}}(\mathcal{L}) \backslash\{ L\} $,
we consider the collection 
\begin{eqnarray*}
&&
\mathcal{Q}_{0}^{\flat \mathcal{H}-big}=\bigcup\limits_{L\in \mathcal{H}}
\mathcal{P}_{L,0}^{\flat \mathcal{H}-big} \\
\text{where } 
&&
\mathcal{P}_{L,0}^{\flat \mathcal{H}-big}=\left\{ \left(
I,J\right) \in \mathcal{P}_{L,0}^{\flat \mathcal{H}}:\text{there is }
L^{\prime }\in \mathcal{T}\left( L\right), J^{\flat }\subset
L^{\prime }\subset I\right\}, L\in \mathcal{H} \\
\text{and } 
&&
\mathcal{P}_{L,0}^{\flat \mathcal{H}}=\left\{ \left(
I,J\right) \in \mathcal{P}:I\in \mathcal{C}_{L}^{\mathcal{H}}
J\in \mathcal{C}_{L}^{\mathcal{H},\flat {shift}}\text{ for some }
L\in \mathcal{H}\right\}, L\in \mathcal{H}
\end{eqnarray*}
and begin by claiming that
\begin{equation}
\widehat{\mathfrak{N}}_{{stop},\bigtriangleup ^{\omega }}^{A,
\mathcal{P}_{L,0}^{\flat \mathcal{H}-big}}\leq C\mathcal{S}_{{aug}
{size}}^{\alpha ,A}\left( \mathcal{P}_{L,0}^{\flat \mathcal{H}
-big}\right) \leq C\mathcal{S}_{{aug}{size}}^{\alpha
,A}\left( \mathcal{P}\right) ,\ \ \ \ \ L\in \mathcal{H}.  \label{big t 3}
\end{equation}
To see this, we fix $L\in \mathcal{H}$ and order the cubes of  $\mathcal{T}\left( L\right) =\left\{
L^{k,i}\right\} _{k,i}$, where $1 \leq i \leq n_k$ where $L^0=L$ and  $L^{1,i}$ are the maximal cubes in $L^0$ and then $L^{k+1,i}$ are the maximal cubes inside a cube $L^{k,j}$ of some previous generation. Then $\mathcal{P}_{L,0}^{\flat 
\mathcal{H}-big}$ can be decomposed as follows, remembering that $J^{\flat
}\subset I\subset L$ for $\left( I,J\right) \in \mathcal{P}_{L,0}^{\flat 
\mathcal{H}-big}\subset \mathcal{P}_{L,0}^{\flat \mathcal{H}}$:
\begin{eqnarray*}
\mathcal{P}_{L,0}^{\flat \mathcal{H}-big}
&=&
\bigcup_{k,i}^{\cdot}\left\{ \mathcal{R}_{L_{{out,out}}^{k,i}}^{\flat\mathcal{L}}
\dot{\cup}\ \
\mathcal{R}_{L_{{out,in}}^{k,i}}^{\flat\mathcal{L}}
\dot{\cup}\ \
\mathcal{R}_{L_{{in}}^{k,i}}^{\flat\mathcal{L}}\right\} \\
&=&
\left( \bigcup_{k,i}^{\cdot}\mathcal{R}_{L_{{out,out}}^{k,i}}^{\flat \mathcal{L}}\right) \dot{\cup}\left( \bigcup_{k,i}^{\cdot }\mathcal{R}_{L_{{out,in}}^{k,i}}^{\flat 
\mathcal{L}}\right) \dot{\cup}\left( \bigcup_{k,i}^{\cdot }\mathcal{R}_{L_{{in}}^{k,i}}^{\flat \mathcal{L}
}\right) \ ; \\
\mathcal{R}_{L_{{out,in}}^{k,i}}^{\flat \mathcal{L}} &\equiv &\left\{\left( I,J\right) \in \mathcal{P}_{L,0}^{\flat \mathcal{H}-big}:I\in 
\mathcal{C}_{L^{k-1,i}}^{\mathcal{L}}\text{ and }J^{\flat }\subset L_{{out,in}}^{k,i}\right\} , \\
\mathcal{R}_{L_{{out,out}}^{k,i}}^{\flat\mathcal{L}} 
&\equiv &
\left\{
\left( I,J\right) \in \mathcal{P}_{L,0}^{\flat \mathcal{H}-big}:I\in 
\mathcal{C}_{L^{k-1,i}}^{\mathcal{L}}\text{ and }J^{\flat }\subset L_{{out,out}}^{k,i}\right\}, \\
\mathcal{R}_{L_{{in}}^{k,i}}^{\flat \mathcal{L}} 
&\equiv&
\left\{ \left( I,J\right) \in \mathcal{P}_{L,0}^{\flat \mathcal{H}
-big}:I\in \mathcal{C}_{L^{k-1,i}}^{\mathcal{L}}\text{ and }J^{\flat }\in 
\mathcal{C}_{L^{k-1,i}}^{\mathcal{L}}\text{ and }J^{\flat }\cap
L^{k,i}=\emptyset \right\} \\
&=&\left\{ \left( I,J\right) \in \mathcal{P}_{L,0}^{\flat \mathcal{H}
-big}:I=L^{k-1,i}\text{ and }J^{\flat }\in \mathcal{C}_{L^{k-1,i}}^{\mathcal{L}}
\text{ and }J^{\flat }\cap L^{k-1,i}_{out}=\emptyset, \right\} ,
\end{eqnarray*}
where by $L^{k,i}_{in}$ we denote the union of the children of $L^{k,i}$ that do not touch the boundary of $L$, by $L^{k,i}_{out,in}$ the union of the grandchildren of $L^{k,i}$ that do not touch the boundary of $L$ while their father does, and by $L^{k,i}_{out,out}$ the grandchildren of $L^{k,i}$ that touch the boundary of $L$ and where in the last line we have used the fact that if $I,J^{\flat }\in 
\mathcal{C}_{L^{k-1,i}}^{\mathcal{L}}$ and there is $L^{\prime }\in \mathcal{T}
\left( L\right) $ with $J^{\flat }\subset L^{\prime }\subset I$, then we
must have $I=L^{k-1,i}$. All of the pairs $\left( I,J\right) \in \mathcal{P}
_{L,0}^{\flat \mathcal{H}-big}$ are included in either $\mathcal{R}_{L_{
{out,in}}^{k,i}}^{\flat \mathcal{L}}$, $\mathcal{R}_{L_{{out,out}
}^{k,i}}^{\flat \mathcal{L}}$ or $\mathcal{R}_{L_{{in}
}^{k,i}}^{\flat \mathcal{L}}$ for some $k$, since if $J^{\flat }\supset L^{k,i}$
, then $J^{\flat }$ shares boundary with $L$, which contradicts the fact
that $3J^{\flat }\subset J^{\maltese }\subset I\subset L$.

We can easily deal with the `in' collection $\mathcal{Q}^{{in}}\equiv \overset{\cdot }{\bigcup }_{k=1}^{\infty }\mathcal{R}_{L_{
{in}}^{k,i}}^{\flat \mathcal{L}}$ by applying a $\emph{trivial}$
case of the $\flat $Straddling Lemma to $\mathcal{R}_{L_{{in}
}^{k,i}}^{\flat \mathcal{L}}$ with a single straddling cube, followed by
an application of the Orthogonality Lemma to $\mathcal{Q}^{{in}
}$. More precisely, every pair $\left( I,J\right) \in \mathcal{R}_{L_{
{in}}^{k,i}}^{\flat \mathcal{L}}$ satisfies $J^{\flat }\subset
L^{k-1,i}=I$, so that the reduced admissible collection $\mathcal{R}_{L_{
{in}}^{k,i}}^{\flat \mathcal{L}}$ $\flat $straddles the trivial
choice $\mathcal{S}=\left\{ L^{k-1,i}\right\} $, the singleton consisting of
just the cube $L^{k-1,i}$. Then the inequality
\begin{equation*}
\widehat{\mathfrak{N}}_{{stop},\bigtriangleup ^{\omega }}^{A,
\mathcal{R}_{L_{{in}}^{k,i}}^{\flat \mathcal{L}}}\leq C\mathcal{S
}_{{aug}{size}}^{\alpha ,A}\left( \mathcal{R}_{L_{{in}}^{k,i}}^{\flat \mathcal{L}}\right) ,
\end{equation*}
follows from $\flat $Straddling Lemma \ref{straddle 3 ref}. The collection $
\left\{ \mathcal{R}_{L_{{in}}^{k,i}}^{\flat \mathcal{L}}\right\}
_{k,i}$ is mutually orthogonal since
\begin{eqnarray*}
\mathcal{R}_{L_{{in}}^{k,i}}^{\flat \mathcal{L}} &\subset &
\mathcal{C}_{L^{k-1,i}}^{\mathcal{L}}\times \mathcal{C}_{L^{k-1,i}}^{\mathcal{L}
,\flat {shift}} \\
\sum\limits_{k=1}^{\infty}\sum\limits_{i=1}^{n_k}\mathbf{1}_{\mathcal{C}_{L^{k-1,i}}^{\mathcal{L}}}\!\!\!\!\!
&\leq &\!\!\!
\mathbf{1}
\text{ and }\sum\limits_{k=1}^{\infty }\sum\limits_{i=1}^{n_k}\mathbf{1}_{
\mathcal{C}_{L^{k-1,i}}^{\mathcal{L},\flat {shift}}}\leq \mathbf{1}.
\end{eqnarray*}
Since $\displaystyle\bigcup_{k,i}^{\cdot }\mathcal{R}_{L_{{in}}^{k,i}}^{\flat \mathcal{L}}$ is reduced and admissible (each $J\in
\Pi _{2}\left(\displaystyle\bigcup_{k,i}^{\cdot }\mathcal{R}_{L_{
{in}}^{k,i}}^{\flat \mathcal{L}}\right) $ is paired with a
single $I$, namely the top of the $\mathcal{L}$-corona to which $J^{\flat }$
belongs), the Orthogonality Lemma \ref{mut orth} applies to obtain the
estimate
\begin{equation}
\widehat{\mathfrak{N}}_{{stop},\bigtriangleup ^{\omega }}^{A,\bigcup_{k,i}\mathcal{R}_{L_{{in}}^{k,i}}^{\flat \mathcal{L}}}
\leq
\sup_{\substack{1 \leq k\\1 \leq i \leq n_k}}\widehat{\mathfrak{N}}_{{stop},\bigtriangleup ^{\omega
}}^{A,\mathcal{R}_{L_{{in}}^{k,i}}^{\flat \mathcal{L}}}
\leq
C\sup_{\substack{1 \leq k\\1 \leq i \leq n_k}}\mathcal{S}_{{aug}{size}}^{\alpha ,A}\left( 
\mathcal{R}_{L_{{in}}^{k,i}}^{\flat \mathcal{L}}\right) \leq C
\mathcal{S}_{{aug}{size}}^{\alpha ,A}\left( \mathcal{P}
_{L,0}^{\flat \mathcal{H}-big}\right)   \label{disjoint bound}
\end{equation}
Now we turn to estimating the norm of the `out-in' collection $\mathcal{Q}^{
{out,in}}\equiv \displaystyle \bigcup_{k,i}\mathcal{R}_{L_{{out,in}
}^{k,i}}^{\flat \mathcal{L}}$. First we note that $L_{{out,in}}^{k,i}\in 
\mathcal{C}_{A}^{\mathcal{A},{restrict}}$ if $\left( I,J\right) \in 
\mathcal{R}_{L_{{out,in}}^{k,i}}^{\flat \mathcal{L}}$ since $\mathcal{R}
_{L_{{out,in}}^{k,i}}^{\flat \mathcal{L}}$ is reduced, i.e. doesn't
contain any pairs $\left( I,J\right) $ with $J^{\flat }\subset A^{\prime }$
for some $A^{\prime }\in \mathfrak{C}_{\mathcal{A}}\left( A\right) $. Next
we note that $\mathcal{Q}^{{out,in}}$ is admissible since if $J\in
\Pi _{2}\mathcal{Q}^{{out,in}}$, then $J\in \Pi _{2}\mathcal{R}_{L_{{out,in}}^{k,i}}^{\flat \mathcal{L}}$ for a unique index $(k,i)$, and of
course $\mathcal{R}_{L_{{out,in}}^{k,i}}^{\flat \mathcal{L}}$ is
admissible, so that the cubes $I$ that are paired with $J$ are
tree-connected. Thus we can apply the Straddling Lemma \ref{straddle 3 ref}
to the reduced admissible collection $\mathcal{Q}^{{out,in}}$ with
the `straddling' set $\mathcal{S}\equiv\left( \bigcup_{k,i}\bigcup_{L' \in L^{k,i}} L'\right)
\cap \mathcal{C}_{A}^{\mathcal{A},{restrict}}$ to obtain the estimate
\begin{equation}
\widehat{\mathfrak{N}}_{{stop},\bigtriangleup ^{\omega
}}^{A,\bigcup_{k=1}^{\infty }\mathcal{R}_{L_{{out,in}}^{k,i}}^{\flat 
\mathcal{L}}}=\widehat{\mathfrak{N}}_{{stop},\bigtriangleup ^{\omega
}}^{A,\mathcal{Q}^{{out,in}}}\leq C\mathcal{S}_{{aug}{
size}}^{\alpha ,A}\left( \mathcal{Q}^{{out,in}}\right) \leq C\mathcal{
S}_{{aug}{size}}^{\alpha ,A}\left( \mathcal{P}_{L,0}^{\flat 
\mathcal{H}-big}\right) \   \label{right bound}
\end{equation}
As for the remaining `out-out' form $\left\vert \mathsf{B}\right\vert _{
{stop},\bigtriangleup ^{\omega }}^{A,\bigcup_{k,i}\mathcal{
R}_{L_{{out,out}}^{k,i}}^{\flat \mathcal{L}}}\left( f,g\right) $, if the
cube pair $\left( I,J\right) \in \mathcal{R}_{L_{{out,out}
}^{k,i}}^{\flat \mathcal{L}}$, then either $J^{\flat }\subset L' \in  L_{{out,out}}^{k,i}\subsetneqq J^{\maltese }$ or $J^{\maltese }\subset L' \in L_{{out,out}
}^{k,i}$. But $J^{\flat }\subset L'\subsetneqq J^{\maltese
}$ implies that either $J^{\flat }=L'\subsetneqq
J^{\maltese }\subset I\subset L$, which is impossible since $J^{\flat }$
cannot share an endpoint with $L$, or that $J^{\flat }=L''\in L'_{in} $ and $J^{\maltese }=L^{k,i}$.
So we conclude that if $\left( I,J\right) \in \mathcal{R}_{L_{{out,out}}^{k,i}}^{\flat \mathcal{L}}$, then \begin{equation}
\text{either }J^{\maltese }\subset L_{{out,out}}^{k,i}
\text{ or \{}
J^{\maltese }=L^{k,i}\text{ and }J\subset L_{{out,out}}^{k,i}\text{\}}.
\label{either or}
\end{equation}
In either case in (\ref{either or}), there is a unique cube $K\left[ J\right] \in \mathcal{W}\left( L\right) $ that contains $J$. It follows that
there are now two remaining cases:

\textbf{Case 1}: $K\left[ J\right] \in \mathcal{C}_{A}'$,

\textbf{Case 2}: $K\left[ J\right] \subset A^{\prime }\subsetneqq I$ for
some $A^{\prime }\in \mathfrak{C}_{\mathcal{A}}\left( A\right) $.

However, since $J^\flat\subset K[J]$, as $K[J]$ is the maximal cube whose triple is contained in $L$, and since  $\mathcal{R}_{L_{{out,out}}^{k,i}}^{\flat \mathcal{L}}$ is reduced, the pairs $\left( I,J\right) $ in \textbf{Case 2} lie in the `corona straddling'
collection $P_{{cor}}^{A}$ that was removed from all $A$-admissible
collections in (\ref{empty assumption}) of Conclusion \ref{assume}\ above,
and thus there are no pairs in \textbf{Case 2} here. Thus we conclude that $K
\left[ J\right] \in \mathcal{C}_{A}'$.

We now claim that $3K\left[ J\right] \subset I$ for all pairs $\left(
I,J\right) \in \bigcup_{k,i}\mathcal{R}_{L_{{out,out}
}^{k}}^{\flat \mathcal{L}}$. To see this, suppose that $\left( I,J\right)
\in \mathcal{R}_{L_{{out,out}}^{k,i}}^{\flat \mathcal{L}}$ for some $
k\geq 1$, $1 \leq i \leq n_k$. Then by (\ref{either or}) we have both that $K\left[ J\right]
\subset L_{{out,out}}^{k,i}$ and $L^{k,i}\subsetneqq I$. But then $K\left[ J
\right] \subset L_{{out,out}}^{k,i}$ implies that $3K\left[ J\right]
\subset L^{k,i}\subset I$ as claimed.

Now the `out-out' collection $\displaystyle\mathcal{Q}^{{out,out}}\equiv
\bigcup_{k,i}\mathcal{R}_{L_{{out,out}}^{k,i}}^{\flat \mathcal{L
}}$ is admissible, since if $J\in \Pi _{2}\mathcal{Q}^{{out,out}}$ and $I_{j}\in \Pi _{1}\mathcal{Q}^{{out,out}}$ with $\left( I_{j},J\right)
\in \mathcal{Q}^{{out,out}}$ for $j=1,2$, then $I_{j}\in \mathcal{C}_{L^{k_{j}-1,i}}^{\mathcal{L}}$ for some $k_{j}$ and $i$ and all of the cubes $I\in \left[ I_{1},I_{2}\right] \ $lie in one of the coronas $\mathcal{C}
_{L^{k-1,i}}^{\mathcal{L}}$ for $k$ between $k_{1}$ and $k_{2}$. And of course for those coronas we have $J \in L^{k,i}_{out,out}$. Thus $\left(
I,J\right) \in \mathcal{R}_{L_{{out,out}}^{k}}^{\flat \mathcal{L}}\subset \mathcal{Q}^{{out,out}}$ and we have proved the required
connectedness. From the containment $3K\left[ J\right] \subset I\subset L$
for all $\left( I,J\right) \in \bigcup_{k,i}\mathcal{R}_{L_{{out,out}}^{k,i}}^{\flat \mathcal{L}}$, we now see that the reduced
admissible collection $\mathcal{Q}^{{out,out}}$ \emph{substraddles} the
cube $L$. Hence the Substraddling Lemma \ref{substraddle ref} yields the
bound
\begin{equation}
\widehat{\mathfrak{N}}_{{stop},\bigtriangleup ^{\omega
}}^{A,\bigcup_{k,i}\mathcal{R}_{L_{{out,out}}^{k,i}}^{\flat 
\mathcal{L}}}=\widehat{\mathfrak{N}}_{{stop},\bigtriangleup ^{\omega
}}^{A,\mathcal{Q}^{{out,out}}}\leq C\mathcal{S}_{{aug}{
size}}^{\alpha ,A}\left( \mathcal{Q}^{{out,out}}\right) \leq C\mathcal{S
}_{{aug}{size}}^{\alpha ,A}\left( \mathcal{P}_{L,0}^{\flat 
\mathcal{H}-big}\right) .  \label{left bound}
\end{equation}
Combining the bounds (\ref{disjoint bound}), (\ref{right bound}) and (\ref
{left bound}), we obtain (\ref{big t 3}).

Finally, we observe that the collections $\mathcal{P}_{L,0}^{\flat \mathcal{H
}-big}$ themselves are \emph{mutually orthogonal}, namely 
\begin{eqnarray*}
\mathcal{P}_{L,0}^{\flat \mathcal{H}-big} 
&\subset &
\mathcal{C}_{L}^{
\mathcal{H}}\times \mathcal{C}_{L}^{\mathcal{H},\flat {shift}}\ ,\ \
\ \ \ L\in \mathcal{H}\ , \\
\sum\limits_{L\in \mathcal{H}}\mathbf{1}_{\mathcal{C}_{L}^{\mathcal{H}}}
&\leq &
\mathbf{1}\text{ and }\sum\limits_{L\in \mathcal{H}}\mathbf{1}_{
\mathcal{C}_{L}^{\mathcal{H},\flat {shift}}}\leq \mathbf{1}.
\end{eqnarray*}
Thus an application of the Orthogonality Lemma \ref{mut orth} shows that
\begin{equation*}
\widehat{\mathfrak{N}}_{{stop},\bigtriangleup ^{\omega }}^{A,
\mathcal{Q}_{0}^{\flat \mathcal{H}-big}}\leq \sup_{L\in \mathcal{L}}\widehat{
\mathfrak{N}}_{{stop},\bigtriangleup ^{\omega }}^{A,\mathcal{P}
_{L,0}^{\flat \mathcal{H}-big}}\leq C\mathcal{S}_{{aug}{size}
}^{\alpha ,A}\left( \mathcal{P}\right) .
\end{equation*}
Altogether, the proof of Proposition \ref{bottom up 3} is now complete.
\end{proof}
This finishes the proofs of the inequalities (\ref{First inequality}) and (\ref{B stop form 3}).

\section{Finishing the proof}\label{Sub wrapup}

At this point we have controlled, either directly or probabilistically, the
norms of all of the forms in our decompositions - namely the disjoint,
nearby, far below, paraproduct, neighbour, broken and stopping forms - in
terms of the Muckenhoupt, energy and \emph{functional energy} conditions,
along with an arbitrarily small multiple of the operator norm. Thus it only
remains to control the functional energy condition by the Muckenhoupt and
energy conditions, since then, using $\int \left( T_{\sigma }^{\alpha
}f\right) gd\omega =\Theta \left( f,g\right) +\Theta ^{\ast }\left(
f,g\right) $ with the further decompositions above, we will have shown that
for any fixed tangent line truncation of the operator $T_{\sigma }^{\alpha }$ we have 
\begin{eqnarray*}
\left\vert \int \left( T_{\sigma }^{\alpha }f\right) gd\omega \right\vert =
\boldsymbol{E}_{\Omega }^{\mathcal{D}}\boldsymbol{E}_{\Omega }^{\mathcal{G}%
}\left\vert \int \left( T_{\sigma }^{\alpha }f\right) gd\omega \right\vert
&\leq&\!\!\!\!\!
\boldsymbol{E}_{\Omega }^{\mathcal{D}}\boldsymbol{E}_{\Omega }^{%
\mathcal{G}}\sum_{i=1}^{3}\left( \left\vert \Theta _{i}\left( f,g\right)
\right\vert +\left\vert \Theta _{i}^{\ast }\left( f,g\right) \right\vert
\right) \\
&\leq&\!\!\!\!\!
\left( C_{\eta }\mathcal{NTV}_{\alpha }+\eta \mathfrak{N}%
_{T^{\alpha }}\right) \left\Vert f\right\Vert _{L^{2}\left( \sigma \right)
}\left\Vert g\right\Vert _{L^{2}\left( \omega \right) }
\end{eqnarray*}%
for $f\in L^{2}\left( \sigma \right) \text{ and }g\in L^{2}\left( \omega \right)$,
for an arbitarily small positive constant $\eta >0$, and a correspondingly
large finite constant $C_{\eta }$. Note that the testing constants $%
\mathfrak{T}_{T^{\alpha }}$ and $\mathfrak{T}_{T^{\alpha ,\ast }}$ in $%
\mathcal{NTV}_{\alpha }$ already include the supremum over all tangent line
truncations of $T^{\alpha }$, while the operator norm $\mathfrak{N}%
_{T^{\alpha }}$ on the left refers to a \emph{fixed} tangent line truncation
of $T^{\alpha }$. This gives%
\begin{equation*}
\mathfrak{N}_{T^{\alpha }}=\sup_{\left\Vert f\right\Vert _{L^{2}(\sigma
)}=1}\sup_{\left\Vert g\right\Vert _{L^{2}(\omega )}=1}\left\vert \int
\left( T_{\sigma }^{\alpha }f\right) gd\omega \right\vert \leq C_{\eta }%
\mathcal{NTV}_{\alpha }+\eta \mathfrak{N}_{T^{\alpha }},
\end{equation*}%
and since the truncated operators have finite operator norm $\mathfrak{N}%
_{T^{\alpha }}$, we can absorb the term $\eta \mathfrak{N}_{T^{\alpha }}$
into the left hand side for $\eta <1$ and obtain $\mathfrak{N}_{T^{\alpha
}}\leq C_{\eta }^{\prime }\mathcal{NTV}_{\alpha }$ for each tangent line
truncation of $T^{\alpha }$. Taking the supremum over all such truncations
of $T^{\alpha }$ finishes the proof of Theorem \ref{dim high}.

The task of controlling functional energy is taken up next in the Appendix.

\section{Appendix :\ Control of functional energy\label{equiv}}

Now we arrive at one of the main propositions used in the proof of our
theorem. This result is proved \emph{independently} of the main theorem. The
organization of the proof is almost identical to that of the corresponding
result in\ \cite[pages 128-151]{SaShUr7}, together with the modifications in 
\cite[pages 348-360]{SaShUr9} to accommodate common point masses, but we
repeat the organization here with modifications required for the use of two
independent grids, and the appearance of weak goodness entering through the
cubes $J^{\maltese }$. Recall that the functional energy constant $\mathfrak{%
F}_{\alpha }=\mathfrak{F}_{\alpha }^{\mathbf{b}^{\ast }}\left( \mathcal{D},%
\mathcal{G}\right) $ in (\ref{e.funcEnergy n}), $0\leq \alpha <n$, namely
the best constant in the inequality (see (\ref{def M_r-deep}) below for the
definition of $\mathcal{W}\left( F\right) $), 
\begin{equation}
\sum_{F\in \mathcal{F}}\sum_{M\in \mathcal{W}\left( F\right) }\left( \frac{%
\mathrm{P}^{\alpha }\left( M,h\sigma \right) }{\left\vert M\right\vert^\frac{1}{n} }%
\right) ^{2}\left\Vert \mathsf{Q}_{\mathcal{C}_{F}^{\mathcal{G},{%
shift}};M}^{\omega ,\mathbf{b}^{\ast }}x\right\Vert _{L^{2}\left( \omega
\right) }^{\spadesuit 2}\leq \mathfrak{F}_{\alpha }\lVert h\rVert
_{L^{2}\left( \sigma \right) }\,,  \label{fec}
\end{equation}%
depends on the grids $\mathcal{D}$ and $\mathcal{G}$, the goodness parameter 
$\varepsilon >0$ used in the definition of $J^{\maltese }$ through the
shifted corona $\mathcal{C}_{F}^{\mathcal{G},{shift}}$, and on the
family of martingale differences $\left\{ \bigtriangleup _{J}^{\omega ,%
\mathbf{b}^{\ast }}\right\} _{J\in \mathcal{G}}$ associated with $x\in
L_{loc}^{2}\left( \omega \right) $, but not on the family of dual martingale
differences $\left\{ \square _{I}^{\sigma ,\mathbf{b}}\right\} _{I\in 
\mathcal{D}}$, since the function $h\in L^{2}\left( \sigma \right) $
appearing in the definition of functional energy is not decomposed as a sum
of pseudoprojections $\square _{I}^{\sigma ,\mathbf{b}}h$. Finally, we
emphasize that the pseudoprojection 
\begin{equation}
\mathsf{Q}_{\mathcal{C}_{F}^{\mathcal{G},{shift}};M}^{\omega ,%
\mathbf{b}^{\ast }}\equiv \sum_{J\in \mathcal{C}_{F}^{\mathcal{G},{%
shift}}:\ J\subset M}\bigtriangleup _{J}^{\omega ,\mathbf{b}^{\ast }}
\label{def pseudo rest}
\end{equation}%
here uses the shifted restricted corona in%
\begin{eqnarray}
\mathcal{C}_{F}^{\mathcal{G},{shift}} &=&\left\{ J\in \mathcal{G}%
:J^{\maltese }\in \mathcal{C}_{F}^{\mathcal{D}}\right\} ,
\label{def shift cor rest} \\
\mathcal{C}_{F}^{\mathcal{G},{shift}};K &\equiv &\left\{ J\in 
\mathcal{C}_{F}^{\mathcal{G},{shift}}:J\subset K\right\} ,  \notag
\end{eqnarray}%
where $J^{\maltese }$ is defined using the ${body}$ of a cube
as in Definition \ref{def sharp cross}, and where we have defined here the
`restriction' $\mathcal{C}_{F}^{\mathcal{G},{shift}};K$ to the
cube $K$ of the corona $\mathcal{C}_{F}^{\mathcal{G},{shift}}$
(c.f. $\Pi _{2}^{K}\mathcal{P}$\ in Definition \ref{rest K}, which uses the
stronger requirement $J^{\maltese }\subset K$). Moreover, recall from
Notation in \ref{nonstandard norm} and the definition of $\bigtriangledown _{J}^{\omega }$
in (\ref{Carleson avg op}), that for any subset $\mathcal{H}$ of the grid $%
\mathcal{G}$,%
\begin{eqnarray*}
\left\Vert \mathsf{Q}_{\mathcal{H}}^{\omega ,\mathbf{b}^{\ast }}x\right\Vert
_{L^{2}\left( \omega \right) }^{\spadesuit 2}
&\equiv& 
\sum_{J\in \mathcal{H}%
}\left\Vert \bigtriangleup _{J}^{\omega ,\mathbf{b}^{\ast }}x\right\Vert
_{L^{2}\left( \omega \right) }^{\spadesuit 2}\\
&=&
\sum_{J\in \mathcal{H}}\left(
\left\Vert \bigtriangleup _{J}^{\omega ,\mathbf{b}^{\ast }}x\right\Vert
_{L^{2}\left( \omega \right) }^{2}+\inf_{z\in \mathbb{R}}\left\Vert \widehat{%
\bigtriangledown }_{J}^{\omega }\left( x-z\right) \right\Vert _{L^{2}\left( \omega
\right) }^{2}\right) ,
\end{eqnarray*}%
so that we never need to consider the norm squared $\left\Vert \mathsf{Q}_{%
\mathcal{C}_{F}^{\mathcal{G},{shift}};M}^{\omega ,\mathbf{b}^{\ast
}}x\right\Vert _{L^{2}\left( \omega \right) }^{2}$ of the pseudoprojection $%
\mathsf{Q}_{\mathcal{C}_{F}^{\mathcal{G},{shift}};M}^{\omega ,%
\mathbf{b}^{\ast }}x$, something for which we have no lower Riesz
inequality. Note moreover that for $J\in \mathcal{G}$ and an arbitrary
cube $K$, we have by the frame inequality in (\ref{FRAME})
, 
\begin{eqnarray}
\sum_{J\in \mathcal{G}:\ J\subset K}\left\Vert \bigtriangleup _{J}^{\omega ,%
\mathbf{b}^{\ast }}x\right\Vert _{L^{2}\left( \omega \right) }^{2} &\lesssim
&\left\Vert x-m_{K}^{\omega }\right\Vert _{L^{2}\left( \mathbf{1}_{K}\omega
\right) }^{2},  \label{note more} \\
\sum_{J\in \mathcal{G}:\ J\subset K}\inf_{z\in \mathbb{R}}\left\Vert 
\widehat{\bigtriangledown }_{J}^{\omega }\left( x-z\right) \right\Vert _{L^{2}\left(
\omega \right) }^{2} &\leq &\sum_{J\in \mathcal{G}:\ J\subset K}\left\Vert 
\widehat{\bigtriangledown }_{J}^{\omega }\left\{ \left( x-p\right) \mathbf{1}%
_{K}\left( x\right) \right\} \right\Vert _{L^{2}\left( \omega \right)
}^{2} \\ &\lesssim& \left\Vert \left( x-p\right) \right\Vert _{L^{2}\left( \mathbf{%
1}_{K}\omega \right) }^{2},\ \ \ p\in K,  \notag
\end{eqnarray}%
where the second line follows from (\ref{Car embed}).

\begin{description}
\item[Important note] If $J\in \mathcal{C}_{F}^{\mathcal{G},{shift}}$%
, then in particular $J\Subset _{\mathbf{\rho },\varepsilon }F$ with $%
\mathbf{\rho }=\left[ \frac{3}{\varepsilon }\right] $ as mentioned above notation \ref{Notation rho} , and so $J\cap M\neq \emptyset $ for a \emph{unique} $M\in \mathcal{W}%
\left( F\right) $.
\end{description}

We will show that, uniformly in pairs of grids $\mathcal{D}$ and $\mathcal{G}
$, the functional energy constants $\mathfrak{F}_{\alpha }\left( \mathcal{D},%
\mathcal{G}\right) $ in (\ref{e.funcEnergy n}) are controlled by $\mathcal{A}%
_{2}^{\alpha }$, $A_{2}^{\alpha ,{punct}}$ and the large energy
constant $\mathfrak{E}_{2}^{\alpha }$ - {actually the proof shows that we
have control by the Whitney plugged energy constant as defined in (\ref{def
deep plug}) below}. More precisely this is our control of functional energy
proposition.

\begin{prop}
\label{func ener control}For all grids $\mathcal{D}$ and $\mathcal{G}$, and $%
\varepsilon >0$ sufficiently small, we have%
\begin{eqnarray*}
\mathfrak{F}_{\alpha }^{\mathbf{b}^{\ast }}\left( \mathcal{D},\mathcal{G}%
\right) &\lesssim &\mathfrak{E}_{2}^{\alpha }+\sqrt{\mathcal{A}_{2}^{\alpha }%
}+\sqrt{\mathcal{A}_{2}^{\alpha ,\ast }}+\sqrt{A_{2}^{\alpha ,{punct}%
}}\ , \\
\mathfrak{F}_{\alpha }^{\mathbf{b},\ast }\left( \mathcal{G},\mathcal{D}%
\right) &\lesssim &\mathfrak{E}_{2}^{\alpha ,\ast }+\sqrt{\mathcal{A}%
_{2}^{\alpha }}+\sqrt{\mathcal{A}_{2}^{\alpha ,\ast }}+\sqrt{A_{2}^{\alpha
,\ast ,{punct}}}\ ,
\end{eqnarray*}%
with implied constants independent of the grids $\mathcal{D}$ and $\mathcal{G%
}$.
\end{prop}

In order to prove this proposition, we first turn to recalling these more
refined notions of energy constants.

\subsection{Various energy conditions\label{sub various energy cond}}

In this subsection we recall various refinements of the strong energy
conditions appearing in the main theorem above. Variants of this material
already appear in earlier papers, but we repeat it here both for convenience
and in order to introduce some arguments we will use repeatedly later on.
These refinements represent the `weakest' energy side conditions that
suffice for use in our proof, but despite this, we will usually use the
large energy constant $\mathfrak{E}_{2}^{\alpha }$ in estimates to avoid
having to pay too much attention to which of the energy conditions we need
to use - leaving the determination of the weakest conditions in such
situations to the interested reader. We begin with the notion of `deeply
embedded'. Recall that the goodness parameter $\mathbf{r}\in \mathbb{N}$ is
determined by $\varepsilon >0$ in (\ref{choice of r}), and that $%
0<\varepsilon <\frac{1}{n+1}<\frac{1}{n+1-\alpha }$.

For arbitrary cubes in $\,J,K\in \mathcal{P}$, we say that $J$ is $%
\left( \mathbf{\rho },\varepsilon \right) $-\emph{deeply embedded} in $K$,
which we write as $J\Subset _{\mathbf{\rho },\varepsilon }K$, when $J\subset
K$ and both 
\begin{eqnarray}
\ell \left( J\right) &\leq &2^{-\mathbf{\rho }}\ell \left( K\right) ,
\label{def deep embed} \\
d\left( J,\partial K\right) &\geq &2\ell \left( J\right) ^{\varepsilon }\ell
\left( K\right) ^{1-\varepsilon }.  \notag
\end{eqnarray}%
Note that we use the \emph{boundary} of $K$ for the definition of $J\Subset
_{\mathbf{\rho },\varepsilon }K$, rather than the \emph{skeleton} or \emph{%
body} of $K$, which would result in a more restrictive notion of $\left( 
\mathbf{\rho },\varepsilon \right) $-deeply embedded. We will use this
notion for the purpose of grouping $\varepsilon -{good}$ cubes
into the following collections. Fix grids $\mathcal{D}$ and $\mathcal{G}$.
For $K\in \mathcal{D}$, define the collections,%
\begin{eqnarray}
\mathcal{M}_{\left( \mathbf{\rho },\varepsilon \right) -{deep},%
\mathcal{G}}\left( K\right) &\equiv &\left\{ J\in \mathcal{G}:J\text{ is
maximal w.r.t }J\Subset _{\mathbf{\rho },\varepsilon }K\right\} ,
\label{def M_r-deep} \\
\mathcal{M}_{\left( \mathbf{\rho },\varepsilon \right) -{deep},%
\mathcal{D}}\left( K\right) &\equiv &\left\{ M\in \mathcal{D}:M\text{ is
maximal w.r.t }M\Subset _{\mathbf{\rho },\varepsilon }K\right\} ,  \notag \\
\mathcal{W}\left( K\right) &\equiv &\left\{ M\in \mathcal{D}:M\text{ is
maximal w.r.t }3M\subset K\right\}  \notag
\end{eqnarray}%
where the first two consist of \emph{maximal} $\left( \mathbf{\rho }%
,\varepsilon \right) $-deeply embedded dyadic $\mathcal{G}$-subcubes $J$%
, respectively $\mathcal{D}$-subcubes $M$, of a $\mathcal{D}$-cube $%
K $, and the third consists of the maximal $\mathcal{D}$-subcubes $M$
whose triples are contained in $K$.

Let $\gamma >1$. Then the following bounded overlap property holds where\\ $%
\mathcal{M}_{\left( \mathbf{\rho },\varepsilon \right) -{deep}%
}\left( K\right) $ can be taken to be either $\mathcal{M}_{\left( \mathbf{%
\rho },\varepsilon \right) -{deep},\mathcal{G}}\left( K\right) $ or $%
\mathcal{M}_{\left( \mathbf{\rho },\varepsilon \right) -{deep},%
\mathcal{D}}\left( K\right) $ or $\mathcal{W}\left( K\right) $ throughout.

\begin{lem}
Let $0<\varepsilon \leq 1<\gamma \leq 1+4\cdot 2^{\mathbf{\rho }\left(
1-\varepsilon \right) }$. Then%
\begin{equation}
\sum_{J\in \mathcal{M}_{\left( \mathbf{\rho },\varepsilon \right) -{%
deep}}\left( K\right) }\mathbf{1}_{\gamma J}\leq \beta \mathbf{1}_{\left[
\bigcup\limits_{J\in \mathcal{M}_{\left( \mathbf{\rho },\varepsilon \right)
-{deep}}\left( K\right) }\gamma J\right] }  \label{bounded overlap}
\end{equation}%
holds for some positive constant $\beta $ depending only on $n, \gamma ,\mathbf{%
\rho }$ and $\varepsilon $. In addition $\gamma J\subset K$ for all $J\in 
\mathcal{M}_{\left( \mathbf{\rho },\varepsilon \right) -{deep}%
}\left( K\right) $, and consequently%
\begin{equation}
\sum_{J\in \mathcal{M}_{\left( \mathbf{\rho },\varepsilon \right) -{%
deep}}\left( K\right) }\mathbf{1}_{\gamma J}\leq \beta \mathbf{1}_{K}\ .
\label{bounded overlap in K}
\end{equation}%
A similar result holds for $\mathcal{W}\left( K\right) $.
\end{lem}

\begin{proof}
We suppose $0<\varepsilon <1$ and leave the simpler case $\varepsilon =1$
for the reader. To prove (\ref{bounded overlap}), we first note that there
are at most $\frac{2^{n(\mathbf{\rho }+1)}-1}{2^n-1}$ cubes $J$ contained in $K$ for which $%
\ell \left( J\right) >2^{-\mathbf{\rho }}\ell \left( K\right) $. On the
other hand, the maximal $\left( \mathbf{\rho },\varepsilon \right) $-deeply
embedded subcubes $J$ of $K$ also satisfy the comparability condition%
\begin{eqnarray*}
2\ell \left( J\right) ^{\varepsilon }\ell \left( K\right) ^{1-\varepsilon
}\leq d\left( J,\partial K\right) \leq d\left( \pi J,\partial K\right) +\ell
\left( J\right) 
&\leq&
2\left( 2\ell \left( J\right) \right) ^{\varepsilon
}\ell \left( K\right) ^{1-\varepsilon }+\ell \left( J\right) \\
&\leq&
4\ell
\left( J\right) ^{\varepsilon }\ell \left( K\right) ^{1-\varepsilon }+\ell
\left( J\right) .
\end{eqnarray*}%
Now with $0<\varepsilon <1$ and $\gamma >1$ fixed, let $y\in K$. Then if $%
y\in \gamma J$, we have%
\begin{eqnarray*}
2\ell \left( J\right) ^{\varepsilon }\ell \left( K\right) ^{1-\varepsilon }
&\leq &d\left( J,\partial K\right) \leq \gamma \ell \left( J\right) +d\left(
\gamma J,\partial K\right) \\
&\leq &\gamma \ell \left( J\right) +d\left( y,\partial K\right) .
\end{eqnarray*}%
Now assume that $\frac{\ell \left( J\right) }{\ell \left( K\right) }\leq
\left( \frac{1}{\gamma }\right) ^{\frac{1}{1-\varepsilon }}$. Then we have $%
\gamma \ell \left( J\right) \leq \ell \left( J\right) ^{\varepsilon }\ell
\left( K\right) ^{1-\varepsilon }$ and so 
\begin{equation*}
\ell \left( J\right) ^{\varepsilon }\ell \left( K\right) ^{1-\varepsilon
}\leq d\left( y,\partial K\right) .
\end{equation*}%
But we also have 
\begin{equation*}
d\left( y,\partial K\right) \leq\gamma \ell \left( J\right) +d\left( J,\partial
K\right) \leq \gamma\ell \left( J\right) +4\ell \left( J\right) ^{\varepsilon
}\ell \left( K\right) ^{1-\varepsilon }+\ell \left( J\right) \leq 6\ell
\left( J\right) ^{\varepsilon }\ell \left( K\right) ^{1-\varepsilon },
\end{equation*}%
and so altogether, under the assumption that $\frac{\ell \left( J\right) }{%
\ell \left( K\right) }\leq \left( \frac{1}{\gamma }\right) ^{\frac{1}{%
1-\varepsilon }}$, we have%
\begin{eqnarray*}
\frac{1}{6}d\left( y,\partial K\right) &\leq &\ell \left( J\right)
^{\varepsilon }\ell \left( K\right) ^{1-\varepsilon }\leq d\left( y,\partial
K\right) , \\
\text{i.e. }\left( \frac{1}{6}\frac{d\left( y,\partial K\right) }{\ell
\left( K\right) ^{1-\varepsilon }}\right) ^{\frac{1}{\varepsilon }} 
&\leq&
\ell \left( J\right) \leq \left( \frac{d\left( y,\partial K\right) }{\ell
\left( K\right) ^{1-\varepsilon }}\right) ^{\frac{1}{\varepsilon }},
\end{eqnarray*}%
which shows that the number of $J$'s satisfying $y\in \gamma J$ and 
$\frac{\ell \left( J\right) }{\ell \left( K\right) }\leq \left( \frac{1}{%
\gamma }\right) ^{\frac{1}{1-\varepsilon }}$ is at most $C^{\prime }\frac{1}{%
\varepsilon }$. On the other hand, the number of $J$'s contained in 
$K$ satisfying $y\in \gamma J$ and $\frac{\ell \left( J\right) }{\ell \left(
K\right) }>\left( \frac{1}{\gamma }\right) ^{\frac{1}{1-\varepsilon }}$ is
at most $C^{\prime }\frac{1}{1-\varepsilon }\left( 1+\log _{2}\gamma \right) 
$. This proves (\ref{bounded overlap}) with 
\begin{equation*}
\beta =\frac{2^{n(\mathbf{\rho }+1)}-1}{2^n-1}+C^{\prime }\frac{1}{\varepsilon }+C^{\prime }%
\frac{1}{1-\varepsilon }\left( 1+\log _{2}\gamma \right) .
\end{equation*}

In order to prove (\ref{bounded overlap in K}) it suffices, by (\ref{bounded
overlap}), to prove $\gamma J\subset K$ for all $J\in \mathcal{M}_{\left( 
\mathbf{\rho },\varepsilon \right) -{deep}}\left( K\right) $. But $%
J\in \mathcal{M}_{\left( \mathbf{\rho },\varepsilon \right) -{deep}%
}\left( K\right) $ implies%
\begin{equation*}
2\ell \left( J\right) ^{\varepsilon }\ell \left( K\right) ^{1-\varepsilon
}\leq d\left( J,\partial K\right) =d\left( c_{J},\partial K\right) -\frac{1}{%
2}\ell \left( J\right) .
\end{equation*}%
We wish to show $\gamma J\subset K$, which is implied by 
\begin{equation*}
\gamma \frac{1}{2}\ell \left( J\right) \leq d\left( c_{J},K^{c}\right)
=d\left( J,\partial K\right) +\frac{1}{2}\ell \left( J\right) .
\end{equation*}%
But we have%
\begin{equation*}
d\left( J,\partial K\right) +\frac{1}{2}\ell \left( J\right) \geq 2\ell
\left( J\right) ^{\varepsilon }\ell \left( K\right) ^{1-\varepsilon }+\frac{1%
}{2}\ell \left( J\right) ,
\end{equation*}%
and so it suffices to show that%
\begin{equation*}
2\ell \left( J\right) ^{\varepsilon }\ell \left( K\right) ^{1-\varepsilon }+%
\frac{1}{2}\ell \left( J\right) \geq \gamma \frac{1}{2}\ell \left( J\right) ,
\end{equation*}%
which is equivalent to%
\begin{equation*}
\gamma -1\leq 4\ell \left( J\right) ^{\varepsilon -1}\ell \left( K\right)
^{1-\varepsilon }.
\end{equation*}%
But the smallest that $\ell \left( J\right) ^{\varepsilon -1}\ell \left(
K\right) ^{1-\varepsilon }$ can get for $J\in \mathcal{M}_{\left( \mathbf{%
\rho },\varepsilon \right) -{deep}}\left( K\right) $ is $2^{\mathbf{%
\rho }\left( 1-\varepsilon \right) }\geq 1$, and so $\gamma \leq 1+4\cdot 2^{%
\mathbf{\rho }\left( 1-\varepsilon \right) }$ implies $\gamma -1\leq 4\ell
\left( J\right) ^{\varepsilon -1}\ell \left( K\right) ^{1-\varepsilon }$,
which completes the proof.

The reader can easily verify the same argument works for the Whitney
collection $\mathcal{W}\left( K\right) $.
\end{proof}

Now we recall the notion of \emph{alternate} dyadic cubes from \cite%
{SaShUr7}, which we rename \emph{augmented} dyadic cubes here.

\begin{dfn}
\label{def dyadic}Given a dyadic grid $\mathcal{D}$, the \emph{augmented
dyadic grid} $\mathcal{AD}$ consists\ of those cubes $I$ whose dyadic
children $I^{\prime }$ belong to the grid $\mathcal{D}$.
\end{dfn}

Of course an augmented grid is not actually a grid because the nesting
property fails, but this terminology should cause no confusion. These
augmented grids will be needed in order to use the `prepare to puncture'
argument (introduced in \cite{SaShUr9}) at several places below.

Now we proceed to recall certain of the definitions of various energy
conditions from \cite{SaShUr5} and \cite{SaShUr7}. While these definitions
are not explicitly used in the proof of functional energy, some of the
arguments we give to control them will be appealed to later, and so we take
the time to develop these definitions in detail.

\subsubsection{Whitney energy conditions}

The following definition of Whitney energy condition uses the \emph{Whitney}
decomposition $\mathcal{M}_{\left( \mathbf{\rho },1\right) -{deep},%
\mathcal{D}}\left( I_{r}\right) $ into $\mathcal{D}$-dyadic cubes in
which $\varepsilon =1$, as well as the `large' pseudoprojections%
\begin{equation}
\mathsf{Q}_{K}^{\omega ,\mathbf{b}^{\ast }}\equiv \sum_{J\in \mathcal{G}:\
J\subset K}\bigtriangleup _{J}^{\omega ,\mathbf{b}^{\ast }}.
\label{large pseudo}
\end{equation}

\begin{dfn}
\label{energy condition}Suppose $\sigma $ and $\omega $ are locally finite
positive Borel measures on $\mathbb{R}^n$ and fix $\gamma >1$. Then the\
Whitney energy condition constant $\mathcal{E}_{2}^{\alpha ,{Whitney}}$
is given by%
\begin{equation*}
\left( \mathcal{E}_{2}^{\alpha ,{Whitney}}\right) ^{2}\equiv \sup_{%
\mathcal{D},\mathcal{G}}\sup_{I=\dot{\cup}I_{r}}\frac{1}{\left\vert
I\right\vert _{\sigma }}\sum_{r=1}^{\infty }\sum_{M\in \mathcal{W}\left(
I_{r}\right) }\left( \frac{\mathrm{P}^{\alpha }\left( M,\mathbf{1}%
_{I\setminus \gamma M}\sigma \right) }{\left\vert M\right\vert^\frac{1}{n} }\right)
^{2}\left\Vert \mathsf{Q}_{M}^{\omega ,\mathbf{b}^{\ast }}x\right\Vert
_{L^{2}\left( \omega \right) }^{\spadesuit 2},
\end{equation*}%
where $\sup_{\mathcal{D},\mathcal{G}}\sup_{I=\dot{\cup}I_{r}}$ is taken over

\begin{enumerate}
\item all dyadic grids $\mathcal{D}$ and $\mathcal{G}$,

\item all $\mathcal{D}$-dyadic cubes $I$,

\item and all partitions $\left\{ I_{r}\right\} _{r=1}^{N\text{ or }\infty }$
of the cube $I$ into $\mathcal{D}$-dyadic subcubes $I_{r}$.
\end{enumerate}
\end{dfn}

If the parameter $\gamma >1$ above is chosen sufficiently close to $1$, then
the collection of cubes $\left\{ \gamma M\right\} _{M\in \mathcal{W}%
\left( I_{r}\right) }$ has bounded overlap $\beta $ by (\ref{bounded overlap
in K}), and the Whitney energy constant $\mathcal{E}_{2}^{\alpha ,{%
Whitney}}$ is controlled by the strong energy constant $\mathcal{E}%
_{2}^{\alpha }$ in (\ref{strong b* energy}),%
\begin{equation}
\mathcal{E}_{2}^{\alpha ,{Whitney}}\lesssim \mathcal{E}_{2}^{\alpha }.
\label{en con}
\end{equation}%
Indeed, to see this, fix a decomposition of a cube 
\begin{equation}
I=\overset{\cdot }{\bigcup }_{1\leq r<\infty }\overset{\cdot }{\bigcup }%
_{M\in \mathcal{W}\left( I_{r}\right) }M  \label{decomp int}
\end{equation}%
as in Definition \ref{energy condition}. Then consider the \emph{sub}%
decomposition 
\begin{equation*}
I\supset \overset{\cdot }{\bigcup }_{1\leq r<\infty }\overset{\cdot }{%
\bigcup }_{M\in \mathcal{W}\left( I_{r}\right) }M
\end{equation*}%
of the cube $I$ given by the collection of cubes,%
\begin{equation*}
\mathcal{I}\equiv \overset{\cdot }{\bigcup }_{1\leq r<\infty }\mathcal{W}%
\left( I_{r}\right) .
\end{equation*}%
We then have%
\begin{equation*}
\left( \mathcal{E}_{2}^{\alpha }\right) ^{2}\geq \frac{1}{\left\vert
I\right\vert _{\sigma }}\sum_{r=1}^{\infty }\sum_{M\in \mathcal{W}\left(
I_{r}\right) }\left( \frac{\mathrm{P}^{\alpha }\left( M,\mathbf{1}_{I}\sigma
\right) }{\left\vert M\right\vert ^\frac{1}{n} }\right) ^{2}\left\Vert x-m_{M}^{\omega
}\right\Vert _{L^{2}\left( \mathbf{1}_{M}\omega \right) }^{2}\ .
\end{equation*}%
Now $\mathrm{P}^{\alpha }\left( M,\mathbf{1}_{I}\sigma \right) \geq \mathrm{P%
}^{\alpha }\left( M,\mathbf{1}_{I\setminus \gamma M}\sigma \right) $ and
from (\ref{note more}),%
\begin{equation*}
\left\Vert x-m_{M}^{\omega }\right\Vert _{L^{2}\left( \mathbf{1}_{M}\omega
\right) }^{2}\gtrsim \left\Vert \mathsf{Q}_{M}^{\omega ,\mathbf{b}^{\ast
}}x\right\Vert _{L^{2}\left( \omega \right) }^{\spadesuit 2},
\end{equation*}
and combining these two inequalities, we obtain that%
\begin{equation*}
\left( \mathcal{E}_{2}^{\alpha }\right) ^{2}\geq c\frac{1}{\left\vert
I\right\vert _{\sigma }}\sum_{r=1}^{\infty }\sum_{M\in \mathcal{W}\left(
I_{r}\right) }\left( \frac{\mathrm{P}^{\alpha }\left( M,\mathbf{1}%
_{I\setminus \gamma M}\sigma \right) }{\left\vert M\right\vert^\frac{1}{n}  }\right)
^{2}\left\Vert \mathsf{Q}_{M}^{\omega ,\mathbf{b}^{\ast }}x\right\Vert
_{L^{2}\left( \omega \right) }^{\spadesuit 2}\ .
\end{equation*}%
Thus we conclude that%
\begin{equation*}
\frac{1}{\left\vert I\right\vert _{\sigma }}\sum_{r=1}^{\infty }\sum_{M\in 
\mathcal{W}\left( I_{r}\right) }\left( \frac{\mathrm{P}^{\alpha }\left( M,%
\mathbf{1}_{I\setminus \gamma M}\sigma \right) }{\left\vert M\right\vert^\frac{1}{n}  }%
\right) ^{2}\left\Vert \mathsf{Q}_{M}^{\omega ,\mathbf{b}^{\ast
}}x\right\Vert _{L^{2}\left( \omega \right) }^{\spadesuit 2}\leq \frac{C}{c}%
\beta \left( \mathcal{E}_{2}^{\alpha }\right) ^{2},
\end{equation*}%
and taking the supremum over all decompositions (\ref{decomp int}) as in
Definition \ref{energy condition}, we obtain (\ref{en con}).

There is a similar definition for the dual (backward) Whitney energy
conditions that simply interchanges $\sigma $ and $\omega $ everywhere.
These definitions of\ the Whitney energy conditions depend on the choice of $%
\gamma >1$.

\begin{description}
\item[Commentary on proofs] We now introduce a number of results concerning
partial plugging of the hole for Whitney energy conditions.
\end{description}

Note that we can `partially' plug the $\gamma $-hole in the Poisson integral 
$\mathrm{P}^{\alpha }\left( J,\mathbf{1}_{I\setminus \gamma J}\sigma \right) 
$ for $\mathcal{E}_{2}^{\alpha ,{Whitney}}$ using the offset $%
A_{2}^{\alpha }$ condition and the bounded overlap property (\ref{bounded
overlap in K}). Indeed, define 
\begin{eqnarray}
&&  
\label{plug} \left( \mathcal{E}_{2}^{\alpha ,{Whitney}{partial}}\right)^{2}\\
&\equiv& 
\sup_{\mathcal{D},\mathcal{G}}\sup_{I=\dot{\cup}I_{r}}\frac{1}{%
\left\vert I\right\vert _{\sigma }}\sum_{r=1}^{\infty }\sum_{M\in \mathcal{W}%
\left( I_{r}\right) }\left( \frac{\mathrm{P}^{\alpha }\left( M,\mathbf{1}%
_{I\setminus M}\sigma \right) }{\left\vert M\right\vert ^\frac{1}{n} }\right)
^{2}\left\Vert \mathsf{Q}_{M}^{\omega ,\mathbf{b}^{\ast }}x\right\Vert
_{L^{2}\left( \omega \right) }^{\spadesuit 2}\ .  \notag
\end{eqnarray}%
Recall from (\ref{bounded overlap in K}) that%
\begin{equation*}
\gamma M\subset I_{r}\text{ for all }M\in \mathcal{W}\left( I_{r}\right) 
\text{ provided }\gamma \leq 5.
\end{equation*}%
At this point we need the following analogues of the `energy $A_{2}^{\alpha
} $ conditions' from \cite{SaShUr9}, which we denote by $A_{2}^{\alpha ,%
{energy}}$ and $A_{2}^{\alpha ,\ast ,{energy}}$, and define
by%
\begin{eqnarray}
A_{2}^{\alpha ,{energy}}\left( \sigma ,\omega \right) &\equiv
&\sup_{Q\in \mathcal{P}}\frac{\left\Vert \mathsf{Q}_{Q}^{\omega ,\mathbf{b}%
^{\ast }}\frac{x}{\ell \left( Q\right) }\right\Vert _{L^{2}\left( \omega
\right) }^{\spadesuit 2}}{\left\vert Q\right\vert ^{1-\frac{\alpha}{n} }}\frac{%
\left\vert Q\right\vert _{\sigma }}{\left\vert Q\right\vert ^{1-\frac{\alpha}{n} }},
\label{def energy A2} \\
A_{2}^{\alpha ,\ast ,{energy}}\left( \sigma ,\omega \right) 
&\equiv&
\sup_{Q\in \mathcal{P}}\frac{\left\vert Q\right\vert _{\omega }}{\left\vert
Q\right\vert ^{1-\frac{\alpha}{n} }}\frac{\left\Vert \mathsf{Q}_{Q}^{\sigma ,\mathbf{b}%
}\frac{x}{\ell \left( Q\right) }\right\Vert _{L^{2}\left( \sigma \right)
}^{\spadesuit 2}}{\left\vert Q\right\vert ^{1-\frac{\alpha}{n} }}.  \notag
\end{eqnarray}%
Then if $\gamma \leq 5$, we have%
\begin{eqnarray}
&&\left( \mathcal{E}_{2}^{\alpha ,{Whitney}{partial}}\right)
^{2}  \label{plug the hole deep} \notag \\
&\lesssim &\sup_{\mathcal{D},\mathcal{G}}\sup_{I=\dot{\cup}I_{r}}\frac{1}{%
\left\vert I\right\vert _{\sigma }}\sum_{r=1}^{\infty }\sum_{M\in \mathcal{W}%
\left( I_{r}\right) }\left( \frac{\mathrm{P}^{\alpha }\left( M,\mathbf{1}%
_{I\setminus \gamma M}\sigma \right) }{\left\vert M\right\vert^\frac{1}{n}  }\right)
^{2}\left\Vert \mathsf{Q}_{M}^{\omega ,\mathbf{b}^{\ast }}x\right\Vert
_{L^{2}\left( \omega \right) }^{\spadesuit 2}  \notag \\
&&
+\sup_{\mathcal{D},\mathcal{G}}\sup_{I=\dot{\cup}I_{r}}\frac{1}{\left\vert
I\right\vert _{\sigma }}\sum_{r=1}^{\infty }\sum_{M\in \mathcal{W}\left(
I_{r}\right) }\left( \frac{\mathrm{P}^{\alpha }\left( M,\mathbf{1}_{\gamma
M\setminus M}\sigma \right) }{\left\vert M\right\vert^\frac{1}{n}  }\right)
^{2}\left\Vert \mathsf{Q}_{M}^{\omega ,\mathbf{b}^{\ast }}x\right\Vert
_{L^{2}\left( \omega \right) }^{\spadesuit 2}  \notag \\
&\lesssim &
\left( \mathcal{E}_{2}^{\alpha ,{Whitney}}\right) ^{2}+\sup_{%
\mathcal{D},\mathcal{G}}\sup_{I=\dot{\cup}I_{r}}\frac{1}{\left\vert
I\right\vert _{\sigma }}\sum_{r=1}^{\infty }\sum_{M\in \mathcal{W}\left(
I_{r}\right) }A_{2}^{\alpha ,{energy}}\left\vert \gamma M\right\vert
_{\sigma }\\ &\lesssim& \left( \mathcal{E}_{2}^{\alpha ,{deep}}\right)
^{2}+\beta A_{2}^{\alpha ,{energy}}\ ,  \notag
\end{eqnarray}%
by (\ref{bounded overlap in K}).

\subsubsection{Plugged energy conditions}

We continue to recall some results from \cite{SaShUr9} and \cite{SaShUr10}
that we will use repeatedly here. For example, we will use the punctured
Muckenhoupt conditions $A_{2}^{\alpha ,{punct}}$ and $A_{2}^{\alpha
,\ast ,{punct}}$ to control
the \emph{plugged }energy conditions, where the hole in the argument of the
Poisson term $\mathrm{P}^{\alpha }\left( M,\mathbf{1}_{I\setminus M}\sigma
\right) $ in the partially plugged energy condition above, is replaced with
the `plugged' term $\mathrm{P}^{\alpha }\left( M,\mathbf{1}_{I}\sigma
\right) $, for example%
\begin{equation}
\left( \mathcal{E}_{2}^{\alpha ,{Whitney}{plug}}\right)
^{2}\equiv \sup_{\mathcal{D},\mathcal{G}}\sup_{I=\dot{\cup}I_{r}}\frac{1}{%
\left\vert I\right\vert _{\sigma }}\sum_{r=1}^{\infty }\sum_{M\in \mathcal{W}%
\left( I_{r}\right) }\left( \frac{\mathrm{P}^{\alpha }\left( M,\mathbf{1}%
_{I}\sigma \right) }{\left\vert M\right\vert^\frac{1}{n}  }\right) ^{2}\left\Vert \mathsf{%
Q}_{M}^{\omega ,\mathbf{b}^{\ast }}x\right\Vert _{L^{2}\left( \omega \right)
}^{\spadesuit 2}\ .  \label{def deep plug}
\end{equation}%
By an argument similar to that in (\ref{plug the hole deep}), we obtain%
\begin{equation}
\mathcal{E}_{2}^{\alpha ,{Whitney}{plug}}\lesssim \mathcal{E}%
_{2}^{\alpha ,{Whitney}{partial}}+A_{2}^{\alpha ,{energy}}.
\label{plug the hole deep'}
\end{equation}

We first show that the punctured Muckenhoupt conditions $A_{2}^{\alpha ,%
{punct}}$ and $A_{2}^{\alpha ,\ast ,{punct}}$ control
respectively the `energy $A_{2}^{\alpha }$ conditions' in (\ref{def energy
A2}). We will make reference to the proof of the next lemma (for the $T1$
theorem this is from \cite[Lemma 3.2 on page 328.]{SaShUr9}) several times
in the sequel. We repeat the proof from \cite[Lemma 3.2 on page 328.]%
{SaShUr9} but with modifications to accommodate the differences that arise
here in the setting of a local $Tb$ theorem. Recall that $\mathfrak{P}%
_{\left( \sigma ,\omega \right) }$ is defined below (\ref{puncture}%
) above.

\begin{lem}
\label{energy A2}For any positive locally finite Borel measures $\sigma
,\omega $ we have%
\begin{eqnarray*}
A_{2}^{\alpha ,{energy}}\left( \sigma ,\omega \right) &\lesssim
&A_{2}^{\alpha ,{punct}}\left( \sigma ,\omega \right) , \\
A_{2}^{\alpha ,\ast ,{energy}}\left( \sigma ,\omega \right)
&\lesssim &A_{2}^{\alpha ,\ast ,{punct}}\left( \sigma ,\omega
\right) .
\end{eqnarray*}
\end{lem}

\begin{proof}
Fix a cube $Q\in \mathcal{D}$. Recall the definition of $\omega \left(
Q,\mathfrak{P}_{\left( \sigma ,\omega \right) }\right) $ in (\ref{puncture}%
). If $\omega \left( Q,\mathfrak{P}_{\left( \sigma ,\omega \right) }\right)
\geq \frac{1}{2}\left\vert Q\right\vert _{\omega }$, then we trivially have%
\begin{eqnarray*}
\frac{\left\Vert \mathsf{Q}_{Q}^{\omega ,\mathbf{b}^{\ast }}\frac{x}{\ell
\left( Q\right) }\right\Vert _{L^{2}\left( \omega \right) }^{\spadesuit 2}}{%
\left\vert Q\right\vert ^{1-\frac{\alpha}{n} }}\frac{\left\vert Q\right\vert _{\sigma }%
}{\left\vert Q\right\vert ^{1-\frac{\alpha}{n} }} &\lesssim &\frac{\left\vert
Q\right\vert _{\omega }}{\left\vert Q\right\vert ^{1-\frac{\alpha}{n} }}\frac{%
\left\vert Q\right\vert _{\sigma }}{\left\vert Q\right\vert ^{1-\frac{\alpha}{n} }} \\
&\leq &2\frac{\omega \left( Q,\mathfrak{P}_{\left( \sigma ,\omega \right)
}\right) }{\left\vert Q\right\vert ^{1-\frac{\alpha}{n} }}\frac{\left\vert
Q\right\vert _{\sigma }}{\left\vert Q\right\vert ^{1-\frac{\alpha}{n} }}\leq
2A_{2}^{\alpha ,{punct}}\left( \sigma ,\omega \right) .
\end{eqnarray*}%
On the other hand, if $\omega \left( Q,\mathfrak{P}_{\left( \sigma ,\omega
\right) }\right) <\frac{1}{2}\left\vert Q\right\vert _{\omega }$ then there
is a point $p\in Q\cap \mathfrak{P}_{\left( \sigma ,\omega \right) }$ such
that%
\begin{equation*}
\omega \left( \left\{ p\right\} \right) >\frac{1}{2}\left\vert Q\right\vert
_{\omega }\ ,
\end{equation*}%
and consequently, $p$ is the largest $\omega $-point mass in $Q$. Thus if we
define $\widetilde{\omega }=\omega -\omega \left( \left\{ p\right\} \right)
\delta _{p}$, then we have%
\begin{equation*}
\omega \left( Q,\mathfrak{P}_{\left( \sigma ,\omega \right) }\right)
=\left\vert Q\right\vert _{\widetilde{\omega }}\ .
\end{equation*}%
Now we observe from the construction of martingale differences that%
\begin{equation*}
\bigtriangleup _{J}^{\widetilde{\omega },\mathbf{b}^{\ast }}=\bigtriangleup
_{J}^{\omega ,\mathbf{b}^{\ast }},\ \ \ \ \ \text{for all }J\in \mathcal{D}%
\text{ with }p\notin J.
\end{equation*}%
So for each $s\geq 0$ there is a unique cube $J_{s}\in \mathcal{D}$ with 
$\ell \left( J_{s}\right) =2^{-s}\ell \left( Q\right) $ that contains the
point $p$. Now observe that, just as for the Haar projection, the
one-dimensional projection $\bigtriangleup _{J_{s}}^{\omega ,\mathbf{b}%
^{\ast }}$ is given by $\bigtriangleup _{J_{s}}^{\omega ,\mathbf{b}^{\ast
}}f=\left\langle h_{J_{s}}^{\omega ,\mathbf{b}^{\ast }},f\right\rangle
_{\omega }h_{J_{s}}^{\omega ,\mathbf{b}^{\ast }}$ for a unique up to $\pm $
unit vector $h_{J_{s}}^{\omega ,\mathbf{b}^{\ast }}$. For this cube we
then have%
\begin{eqnarray*}
\left\Vert \bigtriangleup _{J_{s}}^{\omega ,\mathbf{b}^{\ast }}x\right\Vert
_{L^{2}\left( \omega \right) }^{2} &=&\left\vert \left\langle
h_{J_{s}}^{\omega ,\mathbf{b}^{\ast }},x\right\rangle _{\omega }\right\vert
^{2}=\left\vert \left\langle h_{J_{s}}^{\omega ,\mathbf{b}^{\ast
}},x-p\right\rangle _{\omega }\right\vert ^{2} \\
&=&\left\vert \int_{J_{s}}h_{J_{s}}^{\omega ,\mathbf{b}^{\ast }}\left(
x\right) \left( x-p\right) d\omega \left( x\right) \right\vert
^{2}=\left\vert \int_{J_{s}}h_{J_{s}}^{\omega ,\mathbf{b}^{\ast }}\left(
x\right) \left( x-p\right) d\widetilde{\omega }\left( x\right) \right\vert
^{2} \\
&\leq &\left\Vert h_{J_{s}}^{\omega ,\mathbf{b}^{\ast }}\right\Vert
_{L^{2}\left( \widetilde{\omega }\right) }^{2}\left\Vert \mathbf{1}%
_{J_{s}}\left( x-p\right) \right\Vert _{L^{2}\left( \widetilde{\omega }%
\right) }^{2}\leq \left\Vert h_{J_{s}}^{\omega ,\mathbf{b}^{\ast
}}\right\Vert _{L^{2}\left( \omega \right) }^{2}\left\Vert \mathbf{1}%
_{J_{s}}\left( x-p\right) \right\Vert _{L^{2}\left( \widetilde{\omega }%
\right) }^{2} \\
&\leq &\ell \left( J_{s}\right) ^{2}\left\vert J_{s}\right\vert _{\widetilde{%
\omega }}\leq 2^{-2s}\ell \left( Q\right) ^{2}\left\vert Q\right\vert _{%
\widetilde{\omega }},
\end{eqnarray*}%
as well as%
\begin{eqnarray*}
\inf_{z\in \mathbb{R}}\left\Vert \widehat{\nabla }_{J_{s}}^{\omega }\left(
x-z\right) \right\Vert _{L^{2}\left( \omega \right) }^{2} &\lesssim&
\left\Vert
\left( x-p\right) \right\Vert _{L^{2}\left( \mathbf{1}_{J_{s}}\omega \right)
}^{2}=\left\Vert \left( x-p\right) \right\Vert _{L^{2}\left( \mathbf{1}%
_{J_{s}}\widetilde{\omega }\right) }^{2}   \leq\ell \left( J_{s}\right)
^{2}\left\vert J_{s}\right\vert _{\widetilde{\omega }}\\
&\leq &
2^{-2s}\ell
\left( Q\right) ^{2}\left\vert Q\right\vert _{\widetilde{\omega }}\ ,
\end{eqnarray*}%
from (\ref{note more}). Thus we can estimate%
\begin{eqnarray}
&&\label{omega tilda}
\left\Vert \mathsf{Q}_{Q}^{\omega ,\mathbf{b}^{\ast }}\frac{x}{\ell \left(
Q\right) }\right\Vert _{L^{2}\left( \omega \right) }^{\spadesuit 2} \\
&\leq& 
\frac{1}{\ell \left( Q\right) ^{2}} \left( \sum_{J\in \mathcal{D}:\ J\subset
Q}\left\Vert \bigtriangleup _{J}^{\omega ,\mathbf{b}^{\ast }}x\right\Vert
_{L^{2}\left( \omega \right) }^{2} 
+
\inf_{z\in \mathbb{R}}\left\Vert \widehat{%
\nabla }_{J_{s}}^{\omega }\left( x-z\right) \right\Vert _{L^{2}\left( \omega
\right) }^{2} \right)\notag \\
&=&\notag
\frac{1}{\ell \left( Q\right) ^{2}} \Bigg( \sum_{J\in \mathcal{D}:\
p\notin J\subset Q}\left\Vert \bigtriangleup _{J}^{\widetilde{\omega },%
\mathbf{b}^{\ast }}x\right\Vert _{L^{2}\left( \widetilde{\omega }\right)
}^{2}+\sum_{s=0}^{\infty }\left\Vert \bigtriangleup _{J_{s}}^{\omega ,%
\mathbf{b}^{\ast }}x\right\Vert _{L^{2}\left( \omega \right)
}^{2} \\ 
&&
\hspace{3cm}+
\inf_{z\in \mathbb{R}}\left\Vert \widehat{\nabla }_{J_{s}}^{\omega
}\left( x-z\right) \right\Vert _{L^{2}\left( \omega \right) }^{2}  \Bigg) \notag \\
&\lesssim &\frac{1}{\ell \left( Q\right) ^{2}}\left( \left\Vert \mathsf{Q}%
_{Q}^{\widetilde{\omega },\mathbf{b}^{\ast }}x\right\Vert _{L^{2}\left( 
\widetilde{\omega }\right) }^{\spadesuit 2}+\sum_{s=0}^{\infty }2^{-2s}\ell
\left( Q\right) ^{2}\left\vert Q\right\vert _{\widetilde{\omega }}\right) 
\notag \\
&\lesssim &\frac{1}{\ell \left( Q\right) ^{2}}\left( \ell \left( Q\right)
^{2}\left\vert Q\right\vert _{\widetilde{\omega }}+\sum_{s=0}^{\infty
}2^{-2s}\ell \left( Q\right) ^{2}\left\vert Q\right\vert _{\widetilde{\omega 
}}\right)  \notag \\
&\leq &3\left\vert Q\right\vert _{\widetilde{\omega }}=3\omega \left( Q,%
\mathfrak{P}_{\left( \sigma ,\omega \right) }\right) ,  \notag
\end{eqnarray}%
and so 
\begin{equation*}
\frac{\left\Vert \mathsf{Q}_{Q}^{\omega ,\mathbf{b}^{\ast }}\frac{x}{\ell
\left( Q\right) }\right\Vert _{L^{2}\left( \omega \right) }^{\spadesuit 2}}{%
\left\vert Q\right\vert ^{1-\frac{\alpha}{n} }}\frac{\left\vert Q\right\vert _{\sigma }%
}{\left\vert Q\right\vert ^{1-\frac{\alpha}{n} }}\lesssim \frac{3\omega \left( Q,%
\mathfrak{P}_{\left( \sigma ,\omega \right) }\right) }{\left\vert
Q\right\vert ^{1-\frac{\alpha}{n} }}\frac{\left\vert Q\right\vert _{\sigma }}{%
\left\vert Q\right\vert ^{1-\frac{\alpha}{n} }}\leq 3A_{2}^{\alpha ,{punct}%
}\left( \sigma ,\omega \right) .
\end{equation*}%
Now take the supremum over $Q\in \mathcal{D}$ to obtain $A_{2}^{\alpha ,%
{energy}}\left( \sigma ,\omega \right) \lesssim A_{2}^{\alpha ,%
{punct}}\left( \sigma ,\omega \right) $. The dual inequality follows
upon interchanging the measures $\sigma $ and $\omega $.
\end{proof}

We isolate a simple but key fact that will be used repeatedly in what
follows:%
\begin{equation}
\sum_{Q\in \mathcal{D}:\ Q\subset P}\ell \left( Q\right) ^{2}\left\vert
Q\right\vert _{\mu }\lesssim \ell \left( P\right) ^{2}\left\vert
P\right\vert _{\mu }\ ,\ \ \ \ \ \text{for }P\in \mathcal{D}\text{ and }\mu 
\text{ a positive measure}.  \label{key fact}
\end{equation}%
Indeed, to see (\ref{key fact}), simply pigeonhole the length of $Q$
relative to that of $P$ and sum. The next corollary follows immediately from
Lemma \ref{energy A2}, (\ref{plug the hole deep}) and (\ref{plug the hole
deep'}).

\begin{cor}
\label{all plugged}Provided $1<\gamma \leq 5$,%
\begin{equation*}
\mathcal{E}_{2}^{\alpha ,{Whitney}{plug}}\lesssim \mathcal{E}%
_{2}^{\alpha ,{Whitney}{partial}}+A_{2}^{\alpha ,{punct}%
}\lesssim \mathcal{E}_{2}^{\alpha ,{Whitney}}+A_{2}^{\alpha ,{%
punct}}\ ,
\end{equation*}%
and similarly for the dual plugged energy condition.
\end{cor}

Using Lemma \ref{energy A2} we can control the `plugged' energy $\mathcal{A}%
_{2}^{\alpha }$ conditions:%
\begin{eqnarray*}
\mathcal{A}_{2}^{\alpha ,{energy}{plug}}\left( \sigma
,\omega \right) &\equiv &\sup_{Q\in \mathcal{P}}\frac{\left\Vert \mathsf{Q}%
_{Q}^{\omega ,\mathbf{b}^{\ast }}\frac{x}{\ell \left( Q\right) }\right\Vert
_{L^{2}\left( \omega \right) }^{\spadesuit 2}}{\left\vert Q\right\vert
^{1-\frac{\alpha}{n} }}\mathcal{P}^{\alpha }\left( Q,\sigma \right) , \\
\mathcal{A}_{2}^{\alpha ,\ast ,{energy}{plug}}\left( \sigma
,\omega \right) &\equiv &\sup_{Q\in \mathcal{P}}\mathcal{P}^{\alpha }\left(
Q,\omega \right) \frac{\left\Vert \mathsf{Q}_{Q}^{\sigma ,\mathbf{b}}\frac{x%
}{\ell \left( Q\right) }\right\Vert _{L^{2}\left( \sigma \right)
}^{\spadesuit 2}}{\left\vert Q\right\vert ^{1-\frac{\alpha}{n} }}.
\end{eqnarray*}

\begin{lem}
\label{energy A2 plugged}We have
\end{lem}

\begin{eqnarray*}
\mathcal{A}_{2}^{\alpha ,{energy}{plug}}\left( \sigma
,\omega \right) &\mathcal{\lesssim }&\mathcal{A}_{2}^{\alpha }\left( \sigma
,\omega \right) +A_{2}^{\alpha ,{energy}}\left( \sigma ,\omega
\right) , \\
\mathcal{A}_{2}^{\alpha ,\ast ,{energy}{plug}}\left( \sigma
,\omega \right) &\mathcal{\lesssim }&\mathcal{A}_{2}^{\alpha ,\ast }\left(
\sigma ,\omega \right) +A_{2}^{\alpha ,\ast ,{energy}}\left( \sigma
,\omega \right) .
\end{eqnarray*}

\begin{proof}
We have%
\begin{eqnarray*}
\frac{\left\Vert \mathsf{Q}_{Q}^{\omega ,\mathbf{b}^{\ast }}\frac{x}{\ell
\left( Q\right) }\right\Vert _{L^{2}\left( \omega \right) }^{\spadesuit 2}}{%
\left\vert Q\right\vert ^{1-\frac{\alpha}{n} }}\mathcal{P}^{\alpha }\left( Q,\sigma
\right) &=&\frac{\left\Vert \mathsf{Q}_{Q}^{\omega ,\mathbf{b}^{\ast }}\frac{%
x}{\ell \left( Q\right) }\right\Vert _{L^{2}\left( \omega \right)
}^{\spadesuit 2}}{\left\vert Q\right\vert ^{1-\frac{\alpha}{n} }}\mathcal{P}^{\alpha
}\left( Q,\mathbf{1}_{Q^{c}}\sigma \right) \\&+&\frac{\left\Vert \mathsf{Q}%
_{Q}^{\omega ,\mathbf{b}^{\ast }}\frac{x}{\ell \left( Q\right) }\right\Vert
_{L^{2}\left( \omega \right) }^{\spadesuit 2}}{\left\vert Q\right\vert
^{1-\frac{\alpha}{n} }}\mathcal{P}^{\alpha }\left( Q,\mathbf{1}_{Q}\sigma \right) \\
&\lesssim &\frac{\left\vert Q\right\vert _{\omega }}{\left\vert Q\right\vert
^{1-\frac{\alpha}{n} }}\mathcal{P}^{\alpha }\left( Q,\mathbf{1}_{Q^{c}}\sigma \right) +%
\frac{\left\Vert \mathsf{Q}_{Q}^{\omega ,\mathbf{b}^{\ast }}\frac{x}{\ell
\left( Q\right) }\right\Vert _{L^{2}\left( \omega \right) }^{\spadesuit 2}}{%
\left\vert Q\right\vert ^{1-\frac{\alpha}{n} }}\frac{\left\vert Q\right\vert _{\sigma }%
}{\left\vert Q\right\vert ^{1-\frac{\alpha}{n} }} \\
&\lesssim &\mathcal{A}_{2}^{\alpha }\left( \sigma ,\omega \right)
+A_{2}^{\alpha ,{energy}}\left( \sigma ,\omega \right) .
\end{eqnarray*}
\end{proof}

\subsection{The Poisson formulation}

Recall from Definitions \ref{def sharp cross} and \ref{shifted corona} that%
\begin{equation*}
\mathcal{C}_{F}^{\mathcal{G},{shift}}=\left\{ J\in \mathcal{G}%
:J^{\maltese }\in \mathcal{C}_{F}\right\} ,
\end{equation*}%
where $F\in \mathcal{F}$ is a stopping cube in the dyadic grid $\mathcal{%
D}$. For convenience we repeat here the main result of this section,
Proposition \ref{func ener control}.

\begin{prop}
\label{func ener control'}For all grids $\mathcal{D}$ and $\mathcal{G}$, and 
$\varepsilon >0$ sufficiently small, we have%
\begin{eqnarray*}
\mathfrak{F}_{\alpha }^{\mathbf{b}^{\ast }}\left( \mathcal{D},\mathcal{G}%
\right) &\lesssim &\mathfrak{E}_{2}^{\alpha }+\sqrt{\mathcal{A}_{2}^{\alpha }%
}+\sqrt{\mathcal{A}_{2}^{\alpha ,\ast }}+\sqrt{A_{2}^{\alpha ,{punct}%
}}\ , \\
\mathfrak{F}_{\alpha }^{\mathbf{b},\ast }\left( \mathcal{G},\mathcal{D}%
\right) &\lesssim &\mathfrak{E}_{2}^{\alpha ,\ast }+\sqrt{\mathcal{A}%
_{2}^{\alpha }}+\sqrt{\mathcal{A}_{2}^{\alpha ,\ast }}+\sqrt{A_{2}^{\alpha
,\ast ,{punct}}}\ ,
\end{eqnarray*}%
with implied constants independent of the grids $\mathcal{D}$ and $\mathcal{G%
}$.
\end{prop}

To prove Proposition \ref{func ener control'}, we fix grids $\mathcal{D}$
and $\mathcal{G}$ and a subgrid $\mathcal{F}$ of $\mathcal{D}$\ as in (\ref%
{e.funcEnergy n}), and set 
\begin{equation}
\mu \equiv \sum_{F\in \mathcal{F}}\sum_{M\in \mathcal{W}\left( F\right)
}\left\Vert \mathsf{Q}_{F,M}^{\omega ,\mathbf{b}^{\ast }}x\right\Vert
_{L^{2}\left( \omega \right) }^{\spadesuit 2}\cdot \delta _{\left(
c_{M},\ell \left( M\right) \right) }\text{ and }d\overline{\mu }\left(
x,t\right) \equiv \frac{1}{t^{2}}d\mu \left( x,t\right) \ ,  \label{def mu n}
\end{equation}%
where $\mathcal{W}\left( F\right) $ consists of the maximal $\mathcal{D}$%
-subcubes of $F$ whose triples are contained in $F$, and where $\delta
_{\left( c_{M},\ell \left( M\right) \right) }$ denotes the Dirac unit mass
at the point $\left( c_{M},\ell \left( M\right) \right) $ in the upper
half-space $\mathbb{R}_{+}^{n+1}$. Here $M\in \mathcal{D}$ is a dyadic
cube with center $c_{M}$ and side length $\ell \left( M\right) $, and
for any cube $K\in \mathcal{P}$, the shorthand notation $\mathsf{P}%
_{F,K}^{\omega ,\mathbf{b}^{\ast }}$ (resp. $\mathsf{Q}_{F,K}^{\omega ,%
\mathbf{b}^{\ast }}$) is used for the localized pseudoprojection $\mathsf{P}%
_{\mathcal{C}_{F}^{\mathcal{G},{shift}};K}^{\omega ,\mathbf{b}^{\ast
}}$ (resp. $\mathsf{Q}_{\mathcal{C}_{F}^{\mathcal{G},{shift}%
};K}^{\omega ,\mathbf{b}^{\ast }}$) given in (\ref{def localization}):%
\begin{equation}
\mathsf{P}_{F,K}^{\omega ,\mathbf{b}^{\ast }}\equiv \mathsf{P}_{\mathcal{C}%
_{F}^{\mathcal{G},{shift}};K}^{\omega ,\mathbf{b}^{\ast
}}=\sum_{J\subset K:\ J\in \mathcal{C}_{F}^{\mathcal{G},{shift}%
}}\square _{J}^{\omega ,\mathbf{b}^{\ast }}\text{ }
\end{equation}
\begin{equation}
\left( \text{resp. }%
\mathsf{Q}_{F,K}^{\omega ,\mathbf{b}^{\ast }}\equiv \mathsf{Q}_{\mathcal{C}%
_{F}^{\mathcal{G},{shift}};K}^{\omega ,\mathbf{b}^{\ast
}}=\sum_{J\subset K:\ J\in \mathcal{C}_{F}^{\mathcal{G},{shift}%
}}\bigtriangleup _{J}^{\omega ,\mathbf{b}^{\ast }}\right) .  \label{def F,K}
\end{equation}%
We emphasize that all the subcubes $J$ that arise in the projection $%
\mathsf{Q}_{F,M}^{\omega ,\mathbf{b}^{\ast }}$ are good inside the cubes 
$F$ and beyond since $J^{\maltese }\subset F$. Here $J^{\maltese }$ is
defined in Definition \ref{def sharp cross} using the body of a cube.
Thus every $J\in \mathsf{Q}_{F}^{\omega ,\mathbf{b}^{\ast }}$ is contained
in a unique $M\in \mathcal{W}\left( F\right) $, so that $\mathsf{Q}%
_{F}^{\omega ,\mathbf{b}^{\ast }}=\overset{\cdot }{\bigcup }_{M\in \mathcal{%
W}\left( F\right) }\mathsf{Q}_{F,M}^{\omega ,\mathbf{b}^{\ast }}$. We can
replace $x$ by $x-c$ inside the projection for any choice of $c$ we wish;
the projection is unchanged. More generally, $\delta _{q}$ denotes a Dirac
unit mass at a point $q$ in the upper half-space $\mathbb{R}_{+}^{n+1}$.

We will prove the two-weight inequality 
\begin{equation}
\left\Vert \mathbb{P}^{\alpha }\left( f\sigma \right) \right\Vert _{L^{2}(%
\mathbb{R}_{+}^{n+1},\overline{\mu })}\lesssim \left( \mathfrak{E}_{2}^{\alpha
}+\sqrt{\mathcal{A}_{2}^{\alpha }}+\sqrt{\mathcal{A}_{2}^{\alpha ,\ast }}+%
\sqrt{A_{2}^{\alpha ,{punct}}}\right) \lVert f\rVert _{L^{2}\left(
\sigma \right) }\,,  \label{two weight Poisson n}
\end{equation}%
for all nonnegative $f$ in $L^{2}\left( \sigma \right) $, noting that $%
\mathcal{F}$ and $f$ are \emph{not} related here. Above, $\mathbb{P}^{\alpha
}(\cdot )$ denotes the $\alpha $-fractional Poisson extension to the upper
half-space $\mathbb{R}_{+}^{n+1}$,

\begin{equation*}
\mathbb{P}^{\alpha }\rho \left( x,t\right) \equiv \int_{\mathbb{R}^n}\frac{t}{%
\left( t^{2}+\left\vert x-y\right\vert ^{2}\right) ^{\frac{n+1-\alpha }{2}}}%
d\rho \left( y\right) ,
\end{equation*}%
so that in particular 
\begin{equation*}
\left\Vert \mathbb{P}^{\alpha }(f\sigma )\right\Vert _{L^{2}(\mathbb{R}%
_{+}^{2},\overline{\mu })}^{2}=\sum_{F\in \mathcal{F}}\sum_{M\in \mathcal{W}%
\left( F\right) }\mathbb{P}^{\alpha }\left( f\sigma \right) (c(M),\ell
\left( M\right) )^{2}\left\Vert \mathsf{Q}_{F,M}^{\omega ,\mathbf{b}^{\ast }}%
\frac{x}{\left\vert M\right\vert^\frac{1}{n} }\right\Vert _{L^{2}\left( \omega \right)
}^{\spadesuit 2}\,,
\end{equation*}%
and so (\ref{two weight Poisson n}) proves the first line in Proposition \ref%
{func ener control} upon inspecting (\ref{e.funcEnergy n}). Note also that
we can equivalently write $\left\Vert \mathbb{P}^{\alpha }\left( f\sigma
\right) \right\Vert _{L^{2}(\mathbb{R}_{+}^{2},\overline{\mu })}=\left\Vert 
\widetilde{\mathbb{P}}^{\alpha }\left( f\sigma \right) \right\Vert _{L^{2}(%
\mathbb{R}_{+}^{2},\mu )}$ where $\widetilde{\mathbb{P}}^{\alpha }\nu \left(
x,t\right) \equiv \frac{1}{t}\mathbb{P}^{\alpha }\nu \left( x,t\right) $ is
the renormalized Poisson operator. Here we have simply shifted the factor $%
\frac{1}{t^{2}}$ in $\overline{\mu }$ to $\left\vert \widetilde{\mathbb{P}}%
^{\alpha }\left( f\sigma \right) \right\vert ^{2}$ instead, and we will do
this shifting often throughout the proof when it is convenient to do so.

One version of the characterization of the two-weight inequality for
fractional and Poisson integrals in \cite{Saw3} was stated in terms of a
fixed dyadic grid $\mathcal{D}$ of cubes in $\mathbb{R}$ with sides
parallel to the coordinate axes. Using this theorem for the two-weight
Poisson inequality, but adapted to the $\alpha $-fractional Poisson integral 
$\mathbb{P}^{\alpha }$,\footnote{%
The proof for $0\leq \alpha <1$ is essentially identical to that for $\alpha
=0$ given in \cite{Saw3}.}  we see that inequality (\ref{two weight Poisson
n}) requires checking these two inequalities for dyadic cubes $I\in 
\mathcal{D}$ and boxes $\widehat{I}=I\times \left[ 0,\ell \left( I\right)
\right) $ in the upper half-space $\mathbb{R}_{+}^{n+1}$: 
\begin{eqnarray}
\int_{\mathbb{R}_{+}^{2}}\mathbb{P}^{\alpha }\left( \mathbf{1}_{I}\sigma
\right) \left( x,t\right) ^{2}d\overline{\mu }\left( x,t\right) 
&\equiv&
\left\Vert \mathbb{P}^{\alpha }\left( \mathbf{1}_{I}\sigma \right)
\right\Vert _{L^{2}(\overline{\mu })}^{2} \notag\\
&\lesssim&
\left( \left( \mathfrak{E}%
_{2}^{\alpha }\right) ^{2}+\mathcal{A}_{2}^{\alpha }+\mathcal{A}_{2}^{\alpha
,\ast }+A_{2}^{\alpha ,{punct}}\right) \sigma (I)\,,  \label{e.t1 n}
\end{eqnarray}%
\begin{equation}
\int_{\mathbb{R}}[\mathbb{Q}^{\alpha }(t\mathbf{1}_{\widehat{I}}\overline{%
\mu })]^{2}d\sigma (x)\lesssim \left( \left( \mathfrak{E}_{2}^{\alpha
}\right) ^{2}+\mathcal{A}_{2}^{\alpha }+A_{2}^{\alpha ,{punct}%
}\right) \int_{\widehat{I}}t^{2}d\overline{\mu }(x,t),  \label{e.t2 n}
\end{equation}%
for all \emph{dyadic} cubes $I\in \mathcal{D}$, and where the dual
Poisson operator $\mathbb{Q}^{\alpha }$ is given by 
\begin{equation*}
\mathbb{Q}^{\alpha }(t\mathbf{1}_{\widehat{I}}\overline{\mu })\left(
x\right) =\int_{\widehat{I}}\frac{t^{2}}{\left( t^{2}+\lvert x-y\rvert
^{2}\right) ^{\frac{n+1-\alpha }{2}}}d\overline{\mu }\left( y,t\right) \,.
\end{equation*}%
It is important to note that we can choose for $\mathcal{D}$ any fixed
dyadic grid, the compensating point being that the integrations on the left
sides of (\ref{e.t1 n}) and (\ref{e.t2 n}) are taken over the entire spaces $%
\mathbb{R}_{+}^{2}$ and $\mathbb{R}$ respectively\footnote{%
There is a gap in the proof of the Poisson inequality at the top of page 542
in \cite{Saw3}. However, this gap can be fixed as in \cite{SaWh} or \cite%
{LaSaUr1}.}.

\subsection{Poisson testing}

We now turn to proving the Poisson testing conditions (\ref{e.t1 n}) and (%
\ref{e.t2 n}). Similar testing conditions have been considered in \cite%
{SaShUr5}, \cite{SaShUr7}, \cite{SaShUr9} and \cite{SaShUr10}, and the
proofs there essentially carry over to the situation here, but careful
attention must now be paid to the changed definition of functional energy
and the weaker notion of goodness. We continue to circumvent the difficulty
of permitting common point masses here by using the energy Muckenhoupt
constants $A_{2}^{\alpha ,{energy}}$ and $A_{2}^{\alpha ,\ast ,%
{energy}}$, which require control by the punctured Muckenhoupt
constants $A_{2}^{\alpha ,{punct}}$ and $A_{2}^{\alpha ,\ast ,%
{punct}}$. The following elementary Poisson inequalities (see e.g. 
\cite{Vol}) will be used extensively.

\begin{lem}
\label{Poisson inequalities}Suppose that $J,K,I$ are cubes in $\mathbb{R}^n$, and that $\mu $ is a positive measure supported in $\mathbb{R}^n\setminus I$%
. If $J\subset K\subset \beta K\subset I$ for some $\beta >1$, then%
\begin{equation*}
\frac{\mathrm{P}^{\alpha }\left( J,\mu \right) }{\left\vert J\right\vert^\frac{1}{n} }
\approx
\frac{\mathrm{P}^{\alpha }\left( K,\mu \right) }{\left\vert K\right\vert^\frac{1}{n}  },
\end{equation*}%
while if $J\subset \beta K$, then%
\begin{equation*}
\frac{\mathrm{P}^{\alpha }\left( K,\mu \right) }{\left\vert K\right\vert^\frac{1}{n}  }%
\lesssim \frac{\mathrm{P}^{\alpha }\left( J,\mu \right) }{\left\vert J\right\vert ^\frac{1}{n} }.
\end{equation*}
\end{lem}

\begin{proof}
We have%
\begin{equation*}
\frac{\mathrm{P}^{\alpha }\left( J,\mu \right) }{\left\vert J\right\vert^\frac{1}{n}  }=%
\frac{1}{\left\vert J\right\vert^\frac{1}{n}  }\int \frac{\left\vert J\right\vert^\frac{1}{n}  }{%
\left( \left\vert J\right\vert^\frac{1}{n}  +\left\vert x-c_{J}\right\vert \right)
^{n+1-\alpha }}d\mu \left( x\right) ,
\end{equation*}%
where $J\subset K\subset \beta K\subset I$ implies that%
\begin{equation*}
\left\vert J\right\vert^\frac{1}{n}  +\left\vert x-c_{J}\right\vert \approx \left\vert
K\right\vert^\frac{1}{n}  +\left\vert x-c_{K}\right\vert ,\ \ \ \ \ x\in \mathbb{R}^n%
\setminus I,
\end{equation*}%
and where $J\subset \beta K$ implies that%
\begin{equation*}
\left\vert J\right\vert^\frac{1}{n}  +\left\vert x-c_{J}\right\vert \lesssim \left\vert
J\right\vert^\frac{1}{n}  +\left\vert c_{K}-c_{J}\right\vert +\left\vert
x-c_{K}\right\vert \lesssim \left\vert K\right\vert^\frac{1}{n}  +\left\vert
x-c_{K}\right\vert ,\ \ \ \ \ x\in \mathbb{R}^n.
\end{equation*}
\end{proof}

Recall that in the case\ of the $T1$ theorem in \cite{SaShUr7}, where we
assumed \emph{traditional} goodness in a single family of grids $\mathcal{D}$%
, we had a \emph{strong} bounded overlap property associated with the
projections $\mathsf{P}_{F,J}^{\omega ,\mathbf{b}^{\ast }}$ defined there;
namely, that for each cube $I_{0}\in \mathcal{D}$, there were a bounded
number of cubes $F\in \mathcal{F}$ with the property that $F\supsetneqq
I_{0}\supset J$ for some $J\in \mathcal{M}_{\left( \mathbf{\rho }%
,\varepsilon \right) -{deep}}\left( F\right) $ with $\mathsf{P}%
_{F,J}^{\omega ,\mathbf{b}^{\ast }}\neq 0$ (see the first part of Lemma 10.4
in \cite{SaShUr7}). However, we no longer have this strong bounded overlap
property when ordinary goodness is replaced with the \emph{weak} goodness of
Hyt\"{o}nen and Martikainen. Indeed, there may now be an \emph{unbounded}
number of cubes $F\in \mathcal{F}$ with $F\supsetneqq I_{0}\supset J$
and $\mathsf{P}_{F,J}^{\omega ,\mathbf{b}^{\ast }}\neq 0$, simply because
there can be $J^{\prime }\in \mathcal{G}$ with both $J^{\prime }\subset
I_{0} $ and $\left( J^{\prime }\right) ^{\maltese }$ \emph{arbitrarily}
large.

What will save us in obtaining the following lemma is that the Whitney
cubes $M$ in $\mathcal{W}\left( F\right) $ that happen to lie in some $%
I\in \mathcal{D}$ with $I\subset F$ have one of just two different forms: if 
$I$ shares an endpoint with $F$ then the cubes $M$ near that endpoint
are the same as those in $\mathcal{W}\left( I\right) $ - note that $F$ has
been replaced with $I$ here - while otherwise there are a bounded number of
Whitney cubes $M$ in $I$, and each such $M$ has side length comparable
to $\ell \left( I\right) $.

The next lemma will be used in bounding both of the local Poisson testing
conditions. Recall from Definition \ref{def dyadic}\ that $\mathcal{AD}$
consists of all augmented $\mathcal{D}$-dyadic cubes where $K$ is an
augmented dyadic cube if it is a union of $2$ $\mathcal{D}$-dyadic
cubes $K^{\prime }$ with $\ell \left( K^{\prime }\right) =\frac{1}{2}%
\ell \left( K\right) $.

\begin{lem}
\label{refined lemma}Let $\mathcal{D}$ and $\mathcal{G}$ and $\mathcal{%
F\subset D}$ be grids and let $\left\{ \mathsf{Q}_{F,M}^{\omega ,\mathbf{b}%
^{\ast }}\right\} _{\substack{ F\in \mathcal{F}  \\ M\in \mathcal{W}\left(
F\right) }}$ be as in (\ref{def F,K}) above. For any augmented cube $%
I\in \mathcal{AD}$ define%
\begin{equation}
B\left( I\right) \equiv \sum_{F\in \mathcal{F}:\ F\supsetneqq I^{\prime }%
\text{ for some }I^{\prime }\in \mathfrak{C}\left( I\right) }\sum_{M\in 
\mathcal{W}\left( F\right) :\ M\subset I}\left( \frac{\mathrm{P}^{\alpha
}\left( M,\mathbf{1}_{I}\sigma \right) }{\left\vert M\right\vert^\frac{1}{n}  }\right)
^{2}\left\Vert \mathsf{Q}_{F,M}^{\omega ,\mathbf{b}^{\ast }}x\right\Vert
_{L^{2}\left( \omega \right) }^{\spadesuit 2}\ .  \label{term B}
\end{equation}%
Then%
\begin{equation}
B\left( I\right) \lesssim \left( \left( \mathfrak{E}_{2}^{\alpha }\right)
^{2}+A_{2}^{\alpha ,{energy}}\right) \left\vert I\right\vert
_{\sigma }\ .  \label{B bound}
\end{equation}
\end{lem}

\begin{proof}
We first prove the bound (\ref{B bound}) for $B\left( I\right) $ ignoring
for the moment the possible case when $M=I$ in the sum defining $B\left(
I\right) $. So suppose that $I\in \mathcal{AD}$ is an augmented $\mathcal{D}$%
-dyadic cube. Define%
\begin{equation*}
\Lambda ^{\ast }\left( I\right) \equiv \left\{ M\subsetneqq I:M\in \mathcal{W%
}\left( F\right) \text{ for some }F\supsetneqq I^{\prime }\text{, }I^{\prime
}\in \mathfrak{C}\left( I\right) \text{ with }\mathsf{Q}_{F,M}^{\omega ,%
\mathbf{b}^{\ast }}x\neq 0\right\} ,
\end{equation*}%
and pigeonhole this collection as $\Lambda ^{\ast }\left( I\right)
=\bigcup\limits_{I^{\prime }\in \mathfrak{C}\left( I\right) }\Lambda \left(
I^{\prime }\right) $, where for each $I^{\prime }\in \mathfrak{C}\left(
I\right) $ we define 
\begin{equation*}
\Lambda \left( I^{\prime }\right) \equiv \left\{ M\subset I^{\prime }:M\in 
\mathcal{W}\left( F\right) \text{ for some }F\supsetneqq I^{\prime }\text{
with }\mathsf{Q}_{F,M}^{\omega ,\mathbf{b}^{\ast }}x\neq 0\right\} .
\end{equation*}%
Consider first the case when $3I^{\prime }\subset F$, so that $d\left(
I^{\prime },\partial F\right) \geq \ell \left( I^{\prime }\right) $. Then if 
$M\in \mathcal{W}\left( F\right) $ for some $F\supsetneqq I^{\prime }$ we
have $\ell \left( M\right) =d\left( M,\partial F\right) $, and if in
addition $M\subset I^{\prime }$, then $M=I^{\prime }$. Consider the sum over
all $F\supsetneqq I^{\prime }=M$:%
\begin{eqnarray*}
B_{M}\left( I\right) &\equiv &\sum_{F\in \mathcal{F}:\ F\supsetneqq M\text{
for some }M\in \mathfrak{C}\left( I\right) \cap \mathcal{W}\left( F\right)
}\left( \frac{\mathrm{P}^{\alpha }\left( M,\mathbf{1}_{I}\sigma \right) }{%
\left\vert M\right\vert^\frac{1}{n}  }\right) ^{2}\left\Vert \mathsf{Q}_{F,M}^{\omega ,%
\mathbf{b}^{\ast }}x\right\Vert _{L^{2}\left( \omega \right) }^{\spadesuit 2}
\\
&\leq &
\left( \frac{\mathrm{P}^{\alpha }\left( M,\mathbf{1}_{I}\sigma
\right) }{\left\vert M\right\vert^\frac{1}{n}  }\right) ^{2}\left\Vert \mathsf{Q}%
_{M}^{\omega ,\mathbf{b}^{\ast }}x\right\Vert _{L^{2}\left( \omega \right)
}^{\spadesuit 2}\lesssim \left( \frac{\mathrm{P}^{\alpha }\left( I,\mathbf{1}%
_{I}\sigma \right) }{\left\vert I\right\vert^\frac{1}{n}  }\right) ^{2}\left\Vert \mathsf{%
Q}_{I}^{\omega ,\mathbf{b}^{\ast }}x\right\Vert _{L^{2}\left( \omega \right)
}^{\spadesuit 2}\\&\lesssim& A_{2}^{\alpha ,{energy}}\left\vert
I\right\vert _{\sigma }\ ,
\end{eqnarray*}%
where we have used the definitions (\ref{def F,K}) and (\ref{large pseudo}).
Thus we have obtained the bound%
\begin{equation*}
\sum_{F\in \mathcal{F}:\ F\supsetneqq M\text{ for some }M\in \mathfrak{C}%
\left( I\right) \cap \mathcal{W}\left( F\right) }\left( \frac{\mathrm{P}%
^{\alpha }\left( M,\mathbf{1}_{I}\sigma \right) }{\left\vert M\right\vert^\frac{1}{n}  }%
\right) ^{2}\left\Vert \mathsf{Q}_{F,M}^{\omega ,\mathbf{b}^{\ast
}}x\right\Vert _{L^{2}\left( \omega \right) }^{\spadesuit 2}\lesssim
A_{2}^{\alpha ,{energy}}\left\vert I\right\vert _{\sigma }\ .
\end{equation*}

Now we turn to the case $3I^{\prime }\not\subset F$, i.e. when $\partial
I^{\prime }\cap \partial F$ consists of exactly one boundary point. In this
case, if both $M\subset I^{\prime }$ and $M\in \mathcal{W}\left( F\right) $
for some $F\supsetneqq I^{\prime }$, then we must have either $M\in \mathcal{%
W}\left( I^{\prime }\right) $ or $M\in \mathfrak{C}\left( I^{\prime }\right) 
$, since both $M$ and $I^{\prime }$ are then close to the same boundary
point in $\partial F$. Note that it is here that we use the Whitney
decompositions to full advantage. So again we can estimate%
\begin{eqnarray*}
&&\sum_{\substack{ F\in \mathcal{F}:\ F\supsetneqq I^{\prime }\text{ for
some }I^{\prime }\in \mathfrak{C}\left( I\right)  \\ 3I^{\prime }\not\subset
F}}\sum_{M\in \mathcal{W}\left( F\right) :\ M\subset I^{\prime }}\left( 
\frac{\mathrm{P}^{\alpha }\left( M,\mathbf{1}_{I}\sigma \right) }{\left\vert
M\right\vert ^\frac{1}{n} }\right) ^{2}\left\Vert \mathsf{Q}_{F,M}^{\omega ,\mathbf{b}%
^{\ast }}x\right\Vert _{L^{2}\left( \omega \right) }^{\spadesuit 2} \\
&\leq &
\sum_{M\in \left\{ \mathcal{W}\left( I^{\prime }\right) \cup 
\mathfrak{C}\left( I^{\prime }\right) \right\} \cap \mathcal{W}\left(
F\right) }\left( \frac{\mathrm{P}^{\alpha }\left( M,\mathbf{1}_{I}\sigma
\right) }{\left\vert M\right\vert^\frac{1}{n}  }\right) ^{2}\left\Vert \mathsf{Q}%
_{M}^{\omega ,\mathbf{b}^{\ast }}x\right\Vert _{L^{2}\left( \omega \right)
}^{\spadesuit 2}\lesssim \left( \mathfrak{E}_{2}^{\alpha }\right)
^{2}\left\vert I\right\vert _{\sigma }\ .
\end{eqnarray*}

Finally, we consider the case $M=I$. In this case $I\in \mathcal{D}$ and so $%
F\supsetneqq I^{\prime }$ implies $F\supset I$ and we can estimate%
\begin{eqnarray*}
\sum_{F\in \mathcal{F}:\ F\supset I}\left( \frac{\mathrm{P}^{\alpha }\left(
I,\mathbf{1}_{I}\sigma \right) }{\left\vert I\right\vert^\frac{1}{n}  }\right)
^{2}\left\Vert \mathsf{Q}_{F,I}^{\omega ,\mathbf{b}^{\ast }}x\right\Vert
_{L^{2}\left( \omega \right) }^{\spadesuit 2}
&\leq&
\left( \frac{\mathrm{P}%
^{\alpha }\left( I,\mathbf{1}_{I}\sigma \right) }{\left\vert I\right\vert^\frac{1}{n}  }%
\right) ^{2}\left\Vert \mathsf{Q}_{I}^{\omega ,\mathbf{b}^{\ast
}}x\right\Vert _{L^{2}\left( \omega \right) }^{\spadesuit 2}\\
&\lesssim&
A_{2}^{\alpha ,{energy}}\left\vert I\right\vert _{\sigma }\ .
\end{eqnarray*}%
This completes the proof of Lemma \ref{refined lemma}.
\end{proof}

\subsection{The forward Poisson testing inequality}

Fix $I\in \mathcal{D}$. We split the integration on the left side of (\ref%
{e.t1 n}) into a local and global piece:%
\begin{equation*}
\int_{\mathbb{R}_{+}^{n+1}}\mathbb{P}^{\alpha }\left( \mathbf{1}_{I}\sigma
\right) ^{2}d\overline{\mu }=\!\!\!\int_{\widehat{I}}\mathbb{P}^{\alpha }\left( 
\mathbf{1}_{I}\sigma \right) ^{2}d\overline{\mu }+\!\!\int_{\mathbb{R}%
_{+}^{n+1}\setminus \widehat{I}}\!\!\!\!\!\!\!\mathbb{P}^{\alpha }\left( \mathbf{1}%
_{I}\sigma \right) ^{2}d\overline{\mu }\equiv \mathbf{Local}\left( I\right) +%
\mathbf{Global}\left( I\right)
\end{equation*}%
where more explicitly,%
\begin{eqnarray}
&&\mathbf{Local}\left( I\right) \equiv \int_{\widehat{I}}\left[ \mathbb{P}%
^{\alpha }\left( \mathbf{1}_{I}\sigma \right) \left( x,t\right) \right] ^{2}d%
\overline{\mu }\left( x,t\right) ;\ \ \ \ \ \overline{\mu }\equiv \frac{1}{%
t^{2}}\mu ,  \label{def local forward} \\
\text{i.e. }\overline{\mu } &\equiv &\ \sum_{F\in \mathcal{F}}\sum_{M\in 
\mathcal{W}\left( F\right) }\left\Vert \mathsf{Q}_{F,M}^{\omega ,\mathbf{b}%
^{\ast }}\frac{x}{\ell \left( M\right) }\right\Vert _{L^{2}\left( \omega
\right) }^{\spadesuit 2}\cdot \delta _{\left( c_{M},\ell \left( M\right)
\right) },  \notag
\end{eqnarray}%
where we recall $\mathsf{Q}_{F,M}^{\omega ,\mathbf{b}^{\ast }}$ is defined
in (\ref{def F,K}) above. Here is a brief schematic diagram of the
decompositions, with bounds in $\fbox{}$, used in this subsection:%
\begin{equation*}
\fbox{$%
\begin{array}{ccc}
\mathbf{Local}\left( I\right) &  &  \\ 
\downarrow &  &  \\ 
\mathbf{Local}^{{plug}}\left( I\right) & + & \mathbf{Local}^{%
{hole}}\left( I\right) \\ 
\downarrow &  & \fbox{$\left( \mathfrak{E}_{2}^{\alpha }\right) ^{2}$} \\ 
\downarrow &  &  \\ 
A & + & B \\ 
\fbox{$\left( \mathfrak{E}_{2}^{\alpha }\right) ^{2}+A_{2}^{\alpha ,{%
energy}}$} &  & \fbox{$\left( \mathfrak{E}_{2}^{\alpha }\right)
^{2}+A_{2}^{\alpha ,{energy}}$}%
\end{array}%
$}
\end{equation*}%
and%
\begin{equation*}
\fbox{$%
\begin{array}{ccccccc}
\mathbf{Global}\left( I\right) &  &  &  &  &  &  \\ 
\downarrow &  &  &  &  &  &  \\ 
A & + & B & + & C & + & D \\ 
\fbox{$A_{2}^{\alpha }$} &  & \fbox{$\left( \mathfrak{E}_{2}^{\alpha
}\right) ^{2}+A_{2}^{\alpha }+A_{2}^{\alpha ,{energy}}$} &  & \fbox{$%
\mathcal{A}_{2}^{\alpha ,\ast }$} &  & \fbox{$\mathcal{A}_{2}^{\alpha ,\ast
}+A_{2}^{\alpha ,{punct}}$}%
\end{array}%
$}
\end{equation*}

As in our earlier papers \cite{SaShUr2}-\cite{SaShUr10} that used a single
family of random grids, we have the useful equivalence that%
\begin{equation}
\left( c\left( M\right) ,\ell \left( M\right) \right) \in \widehat{I}\text{ 
\textbf{if and only if} }M\subset I,  \label{tent consequence}
\end{equation}%
since $M$ and $I$ live in the common grid $\mathcal{D}$. We thus have

\begin{eqnarray*}
&&\mathbf{Local}\left( I\right) =\int_{\widehat{I}}\mathbb{P}^{\alpha
}\left( \mathbf{1}_{I}\sigma \right) \left( x,t\right) ^{2}d\overline{\mu }%
\left( x,t\right) \\
&=&\sum_{F\in \mathcal{F}}\sum_{M\in \mathcal{W}\left( F\right) :\ M\subset
I}\mathbb{P}^{\alpha }\left( \mathbf{1}_{I}\sigma \right) \left( c_{M},\ell
\left( M\right) \right) ^{2}\left\Vert \mathsf{Q}_{F,M}^{\omega ,\mathbf{b}%
^{\ast }}\frac{x}{\left\vert M\right\vert^\frac{1}{n}  }\right\Vert _{L^{2}\left( \omega
\right) }^{\spadesuit 2} \\
&\approx &
\sum_{F\in \mathcal{F}}\sum_{M\in \mathcal{W}\left( F\right) :\
M\subset I}\mathrm{P}^{\alpha }\left( M,\mathbf{1}_{I}\sigma \right)
^{2}\lVert \mathsf{Q}_{F,M}^{\omega ,\mathbf{b}^{\ast }}\frac{x}{\left\vert
M\right\vert^\frac{1}{n}  }\rVert _{L^{2}\left( \omega \right) }^{\spadesuit 2} \\
&\approx &
\mathbf{Local}^{{plug}}\left( I\right) +\mathbf{Local}^{%
{hole}}\left( I\right) ,
\end{eqnarray*}%
where%
\begin{eqnarray*}
\mathbf{Local}^{{plug}}\left( I\right) &\equiv &\sum_{F\in \mathcal{F%
}}\sum_{M\in \mathcal{W}\left( F\right) ):\ M\subset I}\left( \frac{\mathrm{P%
}^{\alpha }\left( M,\mathbf{1}_{F\cap I}\sigma \right) }{\left\vert
M\right\vert^\frac{1}{n}  }\right) ^{2}\left\Vert \mathsf{Q}_{F,M}^{\omega ,\mathbf{b}%
^{\ast }}x\right\Vert _{L^{2}\left( \omega \right) }^{2}, \\
\mathbf{Local}^{{hole}}\left( I\right) &\equiv &\sum_{F\in \mathcal{F}%
}\sum_{M\in \mathcal{W}\left( F\right) :\ M\subset I}\left( \frac{\mathrm{P}%
^{\alpha }\left( M,\mathbf{1}_{I\setminus F}\sigma \right) }{\left\vert
M\right\vert^\frac{1}{n}  }\right) ^{2}\left\Vert \mathsf{Q}_{F,M}^{\omega ,\mathbf{b}%
^{\ast }}x\right\Vert _{L^{2}\left( \omega \right) }^{\spadesuit 2}.
\end{eqnarray*}%
The `plugged' local sum $\mathbf{Local}^{{plug}}\left( I\right) $
can be further decomposed into 
\begin{eqnarray*}
\mathbf{Local}^{{plug}}\left( I\right) &=&\left\{ \sum_{\substack{F\in 
\mathcal{F}\\ F\subset I}}\!\!+\!\!\sum_{\substack{F\in \mathcal{F}\\ F\supsetneqq I}}\right\} 
\sum_{\substack{M\in \mathcal{W}\left( F\right) \\ M\subset I}}\left( \frac{\mathrm{P}%
^{\alpha }\left( M,\mathbf{1}_{F\cap I}\sigma \right) }{\left\vert
M\right\vert^\frac{1}{n}  }\right) ^{2}\left\Vert \mathsf{Q}_{F,M}^{\omega ,\mathbf{b}%
^{\ast }}x\right\Vert _{L^{2}\left( \omega \right) }^{\spadesuit 2} \\
& =& A+B.
\end{eqnarray*}%
Then an application of the Whitney plugged energy condition gives 
\begin{eqnarray*}
A &=&\sum_{F\in \mathcal{F}:\ F\subset I}\sum_{M\in \mathcal{W}\left(
F\right) }\left( \frac{\mathrm{P}^{\alpha }\left( M,\mathbf{1}_{F\cap
I}\sigma \right) }{\left\vert M\right\vert^\frac{1}{n}  }\right) ^{2}\left\Vert \mathsf{Q}%
_{F,M}^{\omega ,\mathbf{b}^{\ast }}x\right\Vert _{L^{2}\left( \omega \right)
}^{\spadesuit 2} \\
&\leq &\sum_{F\in \mathcal{F}:\ F\subset I}\left( \mathfrak{E}_{2}^{\alpha }+%
\sqrt{A_{2}^{\alpha ,{energy}}}\right) ^{2}\left\vert F\right\vert
_{\sigma }\lesssim \left( \mathfrak{E}_{2}^{\alpha }+\sqrt{A_{2}^{\alpha ,%
{energy}}}\right) ^{2}\left\vert I\right\vert _{\sigma }\,,
\end{eqnarray*}%
since $\left\Vert \mathsf{Q}_{F,M}^{\omega ,\mathbf{b}^{\ast }}x\right\Vert
_{L^{2}\left( \omega \right) }^{\spadesuit 2}\leq \left\Vert \mathsf{Q}%
_{M}^{\omega ,\mathbf{b}^{\ast }}x\right\Vert _{L^{2}\left( \omega \right)
}^{\spadesuit 2}$. We also used here that the stopping cubes $\mathcal{F}
$ satisfy a $\sigma $-Carleson measure estimate, 
\begin{equation*}
\sum_{F\in \mathcal{F}:\ F\subset F_{0}}\left\vert F\right\vert _{\sigma
}\lesssim \left\vert F_{0}\right\vert _{\sigma }.
\end{equation*}%
Lemma \ref{refined lemma} applies to the remaining term $B$ to obtain the
bound%
\begin{equation*}
B\lesssim \left( \left( \mathfrak{E}_{2}^{\alpha }\right) ^{2}+A_{2}^{\alpha
,{energy}}\right) \left\vert I\right\vert _{\sigma }\ .
\end{equation*}

Next we show the inequality with `holes', where the support of $\sigma $ is
restricted to the complement of the cube $F$.

\begin{lem}
\label{local hole}We have 
\begin{equation}
\mathbf{Local}^{{hole}}\left( I\right) \lesssim \left( \mathfrak{E}%
_{2}^{\alpha }\right) ^{2}\left\vert I\right\vert _{\sigma }\,.
\label{RTS n}
\end{equation}
\end{lem}

\begin{proof}
Fix $I\in \mathcal{D}$ and define%
\begin{equation*}
\mathcal{F}_{I}\equiv \left\{ F\in \mathcal{F}:F\subset I\right\} \cup
\left\{ I\right\} ,
\end{equation*}%
and denote by $\pi F$, for this proof only, the parent of $F$ in the tree $%
\mathcal{F}_{I}$. Also denote by $d\left( F,F^{\prime }\right) \equiv d_{%
\mathcal{F}_{I}}\left( F,F^{\prime }\right) $ the distance from $F$ to $%
F^{\prime }$ in the tree $\mathcal{F}_{I}$, and denote by $d\left( F\right)
\equiv d_{\mathcal{F}_{I}}\left( F,I\right) $ the distance of $F$ from the
root $I$. Since $I\setminus F$ appears in the argument of the Poisson
integral, those $F\in \mathcal{F}\setminus \mathcal{F}_{I}$ do not
contribute to the sum and so we estimate%
\begin{equation*}
S\equiv \mathbf{Local}^{{hole}}\left( I\right) =\sum_{F\in \mathcal{F}%
_{I}}\sum_{M\in \mathcal{W}\left( F\right) :\ M\subset I}\left( \frac{%
\mathrm{P}^{\alpha }\left( M,\mathbf{1}_{I\setminus F}\sigma \right) }{%
\left\vert M\right\vert^\frac{1}{n}  }\right) ^{2}\left\Vert \mathsf{Q}_{F,M}^{\omega ,%
\mathbf{b}^{\ast }}x\right\Vert _{L^{2}\left( \omega \right) }^{\spadesuit 2}
\end{equation*}%
by using $\sum_{F^{\prime }\in \mathcal{F}:\ F\subset F^{\prime }\subsetneqq
I}\frac{1}{d\left( F^{\prime }\right) ^{2}}\leq C$ to obtain\footnote{%
In \cite{SaShUr7} and \cite{SaShUr6} the first line of this display
incorrectly avoided the use of the Cauchy-Schwarz inequality. In the earlier
versions \cite{SaShUr5} and version \#2 of \cite{SaShUr6}, the argument was
correctly given by duality. The fix used here is taken from pages 94-95 of
version \#4 of \cite{SaShUr5}.} 

\begin{eqnarray*}
S &=&\sum_{F\in \mathcal{F}_{I}}\sum_{\substack{M\in \mathcal{W}\left( F\right) \\
M\subset I}}\!\!\left( \sum_{F^{\prime }\in \mathcal{F}:\ F\subset F^{\prime
}\subsetneqq I}\frac{d\left( F^{\prime }\right) }{d\left( F^{\prime }\right) 
}\frac{\mathrm{P}^{\alpha }\left( M,\mathbf{1}_{\pi F^{\prime }\setminus
F^{\prime }}\sigma \right) }{\left\vert M\right\vert^\frac{1}{n}  }\right) ^{2}\left\Vert 
\mathsf{Q}_{F,M}^{\omega ,\mathbf{b}^{\ast }}x\right\Vert _{L^{2}\left(
\omega \right) }^{\spadesuit 2} \\
&\leq &
\sum_{F\in \mathcal{F}_{I}}\sum_{\substack{M\in \mathcal{W}\left( F\right) \\
M\subset I}}\!\!\left( \sum_{\substack{F^{\prime }\in \mathcal{F}\\ F\subset F^{\prime
}\subsetneqq I}}\!\!\!\!\frac{1}{d\left( F^{\prime }\right) ^{2}}\right) 
\left( \sum_{\substack{F^{\prime }\in 
\mathcal{F}\\ F\subset F^{\prime }\subsetneqq I}}\!\!\!\!\!d\left( F^{\prime }\right)
^{2}\left( \frac{\mathrm{P}^{\alpha }\left( M,\mathbf{1}_{\pi F^{\prime
}\setminus F^{\prime }}\sigma \right) }{\left\vert M\right\vert^\frac{1}{n}  }\right)
^{2}\right) \left\Vert \mathsf{Q}_{F,M}^{\omega ,\mathbf{b}^{\ast
}}x\right\Vert _{L^{2}\left( \omega \right) }^{\spadesuit 2} \\
&\leq &
C\sum_{F^{\prime }\in \mathcal{F}_{I}}d\left( F^{\prime }\right)
^{2}\sum_{\substack{F\in \mathcal{F}\\ F\subset F^{\prime }}}\sum_{M\in \mathcal{W}%
\left( F\right) :\ M\subset I}\left( \frac{\mathrm{P}^{\alpha }\left( M,%
\mathbf{1}_{\pi F^{\prime }\setminus F^{\prime }}\sigma \right) }{\left\vert
M\right\vert^\frac{1}{n}  }\right) ^{2}\left\Vert \mathsf{Q}_{F,M}^{\omega ,\mathbf{b}%
^{\ast }}x\right\Vert _{L^{2}\left( \omega \right) }^{\spadesuit 2} \\
&=&
C\sum_{F^{\prime }\in \mathcal{F}_{I}}d\left( F^{\prime }\right)
^{2}\!\!\!\!\!\sum_{K\in \mathcal{W}\left( F^{\prime }\right) }\sum_{\substack{F\in \mathcal{F}%
\\ F\subset F^{\prime }}}\sum_{\substack{M\in \mathcal{W}\left( F\right) \\ M\subset
I}}
\!\!\!\left( \frac{\mathrm{P}^{\alpha }\left( M,\mathbf{1}_{\pi F^{\prime
}\setminus F^{\prime }}\sigma \right) }{\left\vert M\right\vert^\frac{1}{n}  }\right)
^{2}\left\Vert \mathsf{Q}_{F,M\cap K}^{\omega ,\mathbf{b}^{\ast
}}x\right\Vert _{L^{2}\left( \omega \right) }^{\spadesuit 2} \\
&\lesssim &
\sum_{F^{\prime }\in \mathcal{F}_{I}}d\left( F^{\prime }\right)
^{2}\!\!\!\!\!\sum_{K\in \mathcal{W}\left( F^{\prime }\right) }\left( \frac{\mathrm{P}%
^{\alpha }\left( K,\mathbf{1}_{\pi F^{\prime }\setminus F^{\prime }}\sigma
\right) }{\left\vert K\right\vert^\frac{1}{n}  }\right) ^{2}
\sum_{\substack{F\in \mathcal{F}\\
F\subset F^{\prime }}}\sum_{\substack{M\in \mathcal{W}\left( F\right) \\  M\subset
I}}\left\Vert \mathsf{Q}_{F,M\cap K}^{\omega ,\mathbf{b}^{\ast }}x\right\Vert
_{L^{2}\left( \omega \right) }^{\spadesuit 2}
\end{eqnarray*}%
where in the fifth line we have used that each $J^{\prime }$ appearing in $%
\mathsf{Q}_{F,M}^{\omega ,\mathbf{b}^{\ast }}$ occurs in one of the $\mathsf{%
Q}_{F,M\cap K}^{\omega ,\mathbf{b}^{\ast }}$ since each $M$ is contained in
a unique $K$. We have also used there the Poisson inequalities in Lemma \ref%
{Poisson inequalities}.

We now use the lower frame inequality applied to the
function $\mathbf{1}_{K}\left( x-m_{K}^{\omega }\right) $ to obtain%
\begin{equation*}
\sum_{F\in \mathcal{F}:\ F\subset F^{\prime }}\sum_{M\in \mathcal{W}\left(
F\right) :\ M\subset I}\left\Vert \mathsf{Q}_{F,M\cap K}^{\omega ,\mathbf{b}%
^{\ast }}x\right\Vert _{L^{2}\left( \omega \right) }^{\spadesuit 2}\lesssim
\left\Vert \mathbf{1}_{K}\left( x-m_{K}^{\omega }\right) \right\Vert
_{L^{2}\left( \omega \right) }^{\spadesuit 2}\ .
\end{equation*}

Since the collection $\mathcal{F}_{I}$ satisfies a Carleson condition,
namely $\sum_{F\in \mathcal{F}_{I}}\left\vert F\cap I^{\prime }\right\vert
_{\sigma }\leq C\left\vert I^{\prime }\right\vert _{\sigma }$ for all
cubes $I^{\prime }$, we have geometric decay in generations:%
\begin{equation}
\sum_{F\in \mathcal{F}_{I}:\ d\left( F\right) =k}\left\vert F\right\vert
_{\sigma }\lesssim 2^{-\delta k}\left\vert I\right\vert _{\sigma }\ ,\ \ \ \
\ k\geq 0.  \label{geometric decay}
\end{equation}%
Indeed, with $m>2C$ we have for each $F^{\prime }\in \mathcal{F}_{I}$,%
\begin{equation}
\sum_{F\in \mathcal{F}_{I}:\ F\subset F^{\prime }\text{ and }d\left(
F,F^{\prime }\right) =m}\left\vert F\cap F^{\prime }\right\vert _{\sigma }<%
\frac{1}{2}\left\vert F^{\prime }\right\vert _{\sigma }\ ,  \label{half}
\end{equation}%
since otherwise%
\begin{equation*}
\sum_{F\in \mathcal{F}_{I}:\ F\subset F^{\prime }\text{ and }d\left(
F,F^{\prime }\right) \leq m}\left\vert F\cap F^{\prime }\right\vert _{\sigma
}\geq m\frac{1}{2}\left\vert F^{\prime }\right\vert _{\sigma }\ ,
\end{equation*}%
a contradiction. Now iterate (\ref{half}) to obtain (\ref{geometric decay}).

Thus we can write%
\begin{eqnarray*}
S &\lesssim &\sum_{F^{\prime }\in \mathcal{F}_{I}}d\left( F^{\prime }\right)
^{2}\sum_{K\in \mathcal{W}\left( F^{\prime }\right) }\left( \frac{\mathrm{P}%
^{\alpha }\left( K,\mathbf{1}_{\pi F^{\prime }\setminus F^{\prime }}\sigma
\right) }{\left\vert K\right\vert^\frac{1}{n}  }\right) ^{2}\left\Vert \mathbf{1}%
_{K}\left( x-m_{K}^{\omega }\right) \right\Vert _{L^{2}\left( \omega \right)
}^{\spadesuit 2} \\
&=&\sum_{k=1}^{\infty }k^{2}\sum_{F^{\prime }\in \mathcal{F}_{I}:\ d\left(
F^{\prime }\right) =k}\sum_{K\in \mathcal{W}\left( F^{\prime }\right)
}\left( \frac{\mathrm{P}^{\alpha }\left( K,\mathbf{1}_{\pi F^{\prime
}\setminus F^{\prime }}\sigma \right) }{\left\vert K\right\vert^\frac{1}{n}  }\right)
^{2}\left\Vert \mathbf{1}_{K}\left( x-m_{K}^{\omega }\right) \right\Vert
_{L^{2}\left( \omega \right) }^{\spadesuit 2}\\&\equiv& \sum_{k=1}^{\infty
}A_{k}\ ,
\end{eqnarray*}%
where $A_{k}$ is defined at the end of the above display. Hence using the
strong energy condition,%
\begin{eqnarray*}
A_{k} &=&k^{2}\sum_{F^{\prime }\in \mathcal{F}_{I}:\ d\left( F^{\prime
}\right) =k}\sum_{K\in \mathcal{W}\left( F^{\prime }\right) }\left( \frac{%
\mathrm{P}^{\alpha }\left( K,\mathbf{1}_{\pi F^{\prime }\setminus F^{\prime
}}\sigma \right) }{\left\vert K\right\vert ^\frac{1}{n} }\right) ^{2}\left\Vert \mathbf{1}%
_{K}\left( x-m_{K}^{\omega }\right) \right\Vert _{L^{2}\left( \omega \right)
}^{\spadesuit 2} \\
&\lesssim &k^{2}\left( \mathfrak{E}_{2}^{\alpha }\right) ^{2}\sum_{F^{\prime
\prime }\in \mathcal{F}_{I}:\ d\left( F^{\prime \prime }\right)
=k-1}\left\vert F^{\prime \prime }\right\vert _{\sigma }\lesssim \left( 
\mathfrak{E}_{2}^{\alpha }\right) ^{2}k^{2}2^{-\delta k}\left\vert
I\right\vert _{\sigma }\ ,
\end{eqnarray*}%
where we have applied the strong energy condition for each $F^{\prime \prime
}\in \mathcal{F}_{I}$ with $d\left( F^{\prime \prime }\right) =k-1$ to obtain%
\begin{equation}
\sum_{F^{\prime }\in \mathcal{F}_{I}:\ \pi F^{\prime }=F^{\prime \prime
}}\sum_{K\in \mathcal{W}\left( F^{\prime }\right) }\left( \frac{\mathrm{P}%
^{\alpha }\left( K,\mathbf{1}_{F^{\prime \prime }\setminus F^{\prime
}}\sigma \right) }{\left\vert K\right\vert^\frac{1}{n}  }\right) ^{2}\left\Vert \mathbf{1}%
_{K}\left( x-m_{K}^{\omega }\right) \right\Vert _{L^{2}\left( \omega \right)
}^{\spadesuit 2}\leq \left( \mathfrak{E}_{2}^{\alpha }\right) ^{2}\left\vert
F^{\prime \prime }\right\vert _{\sigma }\ .  \label{to obtain}
\end{equation}%
Finally then we obtain%
\begin{equation*}
S\lesssim \sum_{k=1}^{\infty }\left( \mathfrak{E}_{2}^{\alpha }\right)
^{2}k^{2}2^{-\delta k}\left\vert I\right\vert _{\sigma }\lesssim \left( 
\mathfrak{E}_{2}^{\alpha }\right) ^{2}\left\vert I\right\vert _{\sigma }\ ,
\end{equation*}%
which is (\ref{RTS n}).
\end{proof}

Altogether we have now proved the estimate $\mathbf{Local}\left( I\right)
\!\lesssim\! \left( \left( \mathfrak{E}_{2}^{\alpha }\right) ^{2}\!+\!A_{2}^{\alpha ,%
{energy}}\right) \left\vert I\right\vert _{\sigma }$ when $I\in 
\mathcal{D}$, i.e. for every dyadic cube $I\in \mathcal{D}$,%
\begin{eqnarray}
&&  \label{local} \\
\mathbf{Local}\left( I\right) &\approx &\sum_{F\in \mathcal{F}}\sum_{M\in 
\mathcal{W}\left( F\right) :\ M\subset I}\left( \frac{\mathrm{P}^{\alpha
}\left( M,\mathbf{1}_{I}\sigma \right) }{\left\vert M\right\vert ^\frac{1}{n} }\right)
^{2}\left\Vert \mathsf{Q}_{F,M}^{\omega ,\mathbf{b}^{\ast }}x\right\Vert
_{L^{2}\left( \omega \right) }^{\spadesuit 2}  \notag \\
&\lesssim &\left( \left( \mathfrak{E}_{2}^{\alpha }\right)
^{2}+A_{2}^{\alpha ,{energy}}\right) \left\vert I\right\vert
_{\sigma },\ \ \ I\in \mathcal{D}.  \notag
\end{eqnarray}

\subsubsection{The augmented local estimate}

For future use in the `prepare to puncture' arguments below, we prove a
strengthening of the local estimate $\mathbf{Local}\left( I\right) $ to 
\emph{augmented} cubes $L\in \mathcal{AD}$.

\begin{lem}
\label{shifted}With notation as above and $L\in \mathcal{AD}$ an augmented
cube, we have 
\begin{eqnarray}
&&  \label{shifted local} \\
\mathbf{Local}\left( L\right) &\equiv &\sum_{F\in \mathcal{F}}\sum_{M\in 
\mathcal{W}\left( F\right) :\ M\subset L}\left( \frac{\mathrm{P}^{\alpha
}\left( M,\mathbf{1}_{L}\sigma \right) }{\left\vert M\right\vert ^\frac{1}{n} }\right)
^{2}\left\Vert \mathsf{Q}_{F,M}^{\omega ,\mathbf{b}^{\ast }}x\right\Vert
_{L^{2}\left( \omega \right) }^{\spadesuit 2}  \notag \\
&\lesssim &\left( \left( \mathfrak{E}_{2}^{\alpha }\right)
^{2}+A_{2}^{\alpha ,{energy}}\right) \left\vert L\right\vert
_{\sigma },\ \ \ L\in \mathcal{AD}.  \notag
\end{eqnarray}
\end{lem}

\begin{proof}
We prove (\ref{shifted local}) by repeating the above proof of (\ref{local})
and noting the points requiring change. First we decompose 
\begin{equation*}
\mathbf{Local}\left( L\right) \lesssim \mathbf{Local}^{{plug}}\left(
L\right) +\mathbf{Local}^{{hole}}\left( L\right) +\mathbf{Local}^{%
{offset}}\left( L\right)
\end{equation*}%
where $\mathbf{Local}^{{plug}}\left( L\right) $,  $\mathbf{Local}^{%
{hole}}\left( L\right) $ are analogous to $\mathbf{Local}^{{plug%
}}\left( I\right) $ and $\mathbf{Local}^{{hole}}\left( I\right) $
above, and where $\mathbf{Local}^{{offset}}\left( L\right) $ is an
additional term arising because $L\setminus F$ need not be empty when $L\cap
F\neq \emptyset $ and $F$ is not contained in $L$:%
\begin{eqnarray*}
\mathbf{Local}^{{plug}}\left( L\right) &\equiv &\sum_{F\in \mathcal{F%
}}\sum_{M\in \mathcal{W}\left( F\right) :\ M\subset L}\left( \frac{\mathrm{P}%
^{\alpha }\left( M,\mathbf{1}_{L\cap F}\sigma \right) }{\left\vert
M\right\vert^\frac{1}{n}  }\right) ^{2}\left\Vert \mathsf{Q}_{F,M}^{\omega ,\mathbf{b}%
^{\ast }}x\right\Vert _{L^{2}\left( \omega \right) }^{\spadesuit 2}\ , \\
\mathbf{Local}^{{hole}}\left( L\right) &\equiv &\sum_{F\in \mathcal{F}%
:\ F\subset L}\sum_{M\in \mathcal{W}\left( F\right) :\ M\subset L}\left( 
\frac{\mathrm{P}^{\alpha }\left( M,\mathbf{1}_{L\setminus F}\sigma \right) }{%
\left\vert M\right\vert^\frac{1}{n}  }\right) ^{2}\left\Vert \mathsf{Q}_{F,M}^{\omega ,%
\mathbf{b}^{\ast }}x\right\Vert _{L^{2}\left( \omega \right) }^{\spadesuit
2}\ \\
\mathbf{Local}^{{offset}}\left( L\right) &\equiv &\sum_{F\in 
\mathcal{F}:\ F\not\subset L}\sum_{M\in \mathcal{W}\left( F\right) :\
M\subset L}\left( \frac{\mathrm{P}^{\alpha }\left( M,\mathbf{1}_{L\setminus
F}\sigma \right) }{\left\vert M\right\vert^\frac{1}{n}  }\right) ^{2}\left\Vert \mathsf{Q}%
_{F,M}^{\omega ,\mathbf{b}^{\ast }}x\right\Vert _{L^{2}\left( \omega \right)
}^{\spadesuit 2}
\end{eqnarray*}%
We have%
\begin{eqnarray*}
 \mathbf{Local}^{{plug}}\left( L\right) &=& \!\!\!\!
 \left\{ \sum_{F\in 
\mathcal{F}:\ F\subset \text{ some }L^{\prime }\in \mathfrak{C}\left(
L\right) }+\sum_{F\in \mathcal{F}:\ F\supsetneqq \text{ some }L^{\prime }\in 
\mathfrak{C}_{\mathcal{D}}\left( L\right) }\right\} \sum_{M\in \mathcal{W}%
\left( F\right) :\ M\subset L} \\
& &\ \ \ \ \ \ \ \ \ \ \ \ \ \ \ \ \ \ \ \ \ \ \ \ \ \ \ \ \ \ \times \left( 
\frac{\mathrm{P}^{\alpha }\left( M,\mathbf{1}_{F\cap L}\sigma \right) }{%
\left\vert M\right\vert^\frac{1}{n}  }\right) ^{2}\left\Vert \mathsf{Q}_{F,M}^{\omega ,%
\mathbf{b}^{\ast }}x\right\Vert _{L^{2}\left( \omega \right) }^{\spadesuit 2}\\
&=& A+B.
\end{eqnarray*}
Term $A$ satisfies%
\begin{equation*}
A\lesssim \left( \mathfrak{E}_{2}^{\alpha }+\sqrt{A_{2}^{\alpha ,{%
energy}}}\right) ^{2}\left\vert L\right\vert _{\sigma }\ ,
\end{equation*}%
just as above using $\left\Vert \mathsf{Q}_{F,M}^{\omega }x\right\Vert
_{L^{2}\left( \omega \right) }^{2}\leq \left\Vert \mathsf{Q}_{M}^{\omega
}x\right\Vert _{L^{2}\left( \omega \right) }^{2}$, and the fact that the
stopping cubes $\mathcal{F}$ satisfy a $\sigma $-Carleson measure
estimate, 
\begin{equation*}
\sum_{F\in \mathcal{F}:\ F\subset L}\left\vert F\right\vert _{\sigma
}\lesssim \left\vert L\right\vert _{\sigma }.
\end{equation*}

Term $B$ is handled directly by Lemma \ref{refined lemma} with the augmented
cube $I=L$ to obtain%
\begin{equation*}
B\lesssim \left( \left( \mathfrak{E}_{2}^{\alpha }\right) ^{2}+A_{2}^{\alpha
,{energy}}\right) \left\vert L\right\vert _{\sigma }\ .
\end{equation*}

To handle $\mathbf{Local}^{{hole}}\left( L\right) $, we define%
\begin{equation*}
\mathcal{F}_{L}\equiv \left\{ F\in \mathcal{F}:F\subset L\right\} \cup
\left\{ L\right\} ,
\end{equation*}%
and follow along the proof there with only trivial changes. The analogue of (%
\ref{to obtain}) is now%
\begin{equation*}
\sum_{F^{\prime }\in \mathcal{F}_{L}:\ \pi F^{\prime }=F^{\prime \prime
}}\sum_{K\in \mathcal{W}\left( F^{\prime }\right) }\left( \frac{\mathrm{P}%
^{\alpha }\left( K,\mathbf{1}_{F^{\prime \prime }\setminus F^{\prime
}}\sigma \right) }{\left\vert K\right\vert^\frac{1}{n}  }\right) ^{2}\left\Vert \mathbf{1}%
_{K}\left( x-m_{K}^{\omega }\right) \right\Vert _{L^{2}\left( \omega \right)
}^{\spadesuit 2}\leq \left( \mathfrak{E}_{2}^{\alpha }\right) ^{2}\left\vert
F^{\prime \prime }\right\vert _{\sigma }\ ,
\end{equation*}%
the only change being that $\mathcal{F}_{L}$ now appears in place of $%
\mathcal{F}_{I}$, so that the energy condition still applies. We conclude
that 
\begin{equation*}
\mathbf{Local}^{{hole}}\left( L\right) \lesssim \left( \mathfrak{E}%
_{2}^{\alpha }\right) ^{2}\left\vert L\right\vert _{\sigma }\ .
\end{equation*}

Finally, the additional term $\mathbf{Local}^{{offset}}\left(
L\right) $ is handled directly by Lemma \ref{refined lemma}, and this
completes the proof of the estimate (\ref{shifted local}) in Lemma \ref%
{shifted}.
\end{proof}

\subsubsection{The global estimate}

Now we turn to proving the following estimate for the global part of the
first testing condition \eqref{e.t1 n}:%
\begin{equation*}
\mathbf{Global}\left( I\right) =\int_{\mathbb{R}_{+}^{n+1}\setminus \widehat{I}%
}\mathbb{P}^{\alpha }\left( \mathbf{1}_{I}\sigma \right) ^{2}d\overline{\mu }%
\lesssim \left( \left( \mathfrak{E}_{2}^{\alpha }\right) ^{2}+\mathcal{A}%
_{2}^{\alpha ,\ast }+A_{2}^{\alpha ,{punct}}\right) \left\vert
I\right\vert _{\sigma }.
\end{equation*}%
We begin by decomposing the integral above into four pieces. We have from (%
\ref{tent consequence}):%
\begin{eqnarray*}
&&\int_{\mathbb{R}_{+}^{n+1}\setminus \widehat{I}}\mathbb{P}^{\alpha }\left( 
\mathbf{1}_{I}\sigma \right) ^{2}d\overline{\mu } \\ &=&\sum_{M:\ \left(
c_{M},\ell \left( M\right) \right) \in \mathbb{R}_{+}^{n+1}\setminus \widehat{I%
}}\mathbb{P}^{\alpha }\left( \mathbf{1}_{I}\sigma \right) \left( c_{M},\ell
\left( M\right) \right) ^{2}\sum_{\substack{ F\in \mathcal{F}:  \\ M\in 
\mathcal{W}\left( F\right) }}\left\Vert \mathsf{Q}_{F,M}^{\omega ,\mathbf{b}%
^{\ast }}\frac{x}{\left\vert M\right\vert^\frac{1}{n}  }\right\Vert _{L^{2}\left( \omega
\right) }^{\spadesuit 2} \\
&=&\left\{ \sum_{\substack{ M\cap 3I=\emptyset  \\ \ell \left( M\right) \leq
\ell \left( I\right) }}+\sum_{M\subset 3I\setminus I}+\sum_{\substack{ M\cap
I=\emptyset  \\ \ell \left( M\right) >\ell \left( I\right) }}%
+\sum_{M\supsetneqq I}\right\} \mathbb{P}^{\alpha }\left( \mathbf{1}%
_{I}\sigma \right) \left( c_{M},\ell \left( M\right) \right) ^{2} \cdot \\ 
&&\hspace{6.25cm}\cdot\sum 
_{\substack{ F\in \mathcal{F}:  \\ M\in \mathcal{W}\left( F\right) }}%
\left\Vert \mathsf{Q}_{F,M}^{\omega ,\mathbf{b}^{\ast }}\frac{x}{\left\vert
M\right\vert^\frac{1}{n}  }\right\Vert _{L^{2}\left( \omega \right) }^{\spadesuit 2} \\
&=&A+B+C+D.
\end{eqnarray*}

We further decompose term $A$ according to the length of $M$ and its
distance from $I$, and then use the pairwise disjointedness of the
projections $\mathsf{Q}_{F,M}^{\omega ,\mathbf{b}^{\ast }}$ in $F$ (see the
definition in (\ref{def F,K})) to obtain:%
\begin{eqnarray*}
A &\lesssim &\sum_{m=0}^{\infty }\sum_{k=1}^{\infty }\sum_{\substack{ %
M\subset 3^{k+1}I\setminus 3^{k}I  \\ \ell \left( M\right) =2^{-m}\ell
\left( I\right) }}\left( \frac{2^{-m}\left\vert I\right\vert }{d\left(
M,I\right) ^{n+1-\alpha }}\left\vert I\right\vert _{\sigma }\right)
^{2}\left\vert M\right\vert _{\omega } \\
&\lesssim &\sum_{m=0}^{\infty }2^{-2m}\sum_{k=1}^{\infty }\frac{\left\vert
I\right\vert ^{2}\left\vert I\right\vert _{\sigma }\left\vert
3^{k+1}I\setminus 3^{k}I\right\vert _{\omega }}{\left\vert 3^{k}I\right\vert
^{2\left( n+1-\alpha \right) }}\left\vert I\right\vert _{\sigma } \\
&\lesssim &\sum_{m=0}^{\infty }2^{-2m}\sum_{k=1}^{\infty }3^{-2k}\left\{ 
\frac{\left\vert 3^{k+1}I\setminus 3^{k}I\right\vert _{\omega }\left\vert
3^{k}I\right\vert _{\sigma }}{\left\vert 3^{k}I\right\vert ^{2\left(
1-\alpha \right) }}\right\} \left\vert I\right\vert _{\sigma }\lesssim
A_{2}^{\alpha }\left\vert I\right\vert _{\sigma },
\end{eqnarray*}%
where the offset Muckenhoupt constant $A_{2}^{\alpha }$ applies because $%
3^{k+1}I$ has only three times the side length of $3^{k}I$.

\medskip

For term $B$ we first dispose of the nearby sum $B_{{nearby}}$ that
consists of the sum over those $M$ which satisfy in addition $2^{-\mathbf{%
\rho }}\ell \left( I\right) \leq \ell \left( M\right) \leq \ell \left(
I\right) $. But it is a straightforward task to bound $B_{{nearby}}$
by $CA_{2}^{\alpha ,{energy}}\left\vert I\right\vert _{\sigma }$ as
there are at most $2^{\mathbf{\rho }+1}$ such cubes $M$. To bound $B_{%
{away}}\equiv B-B_{{nearby}}$, we further decompose the sum
over $F\in \mathcal{F}$ according to whether or not $F\subset 3I\setminus I$:%
\begin{eqnarray*}
B_{{away}} &\approx &\sum_{M\subset 3I\setminus I\text{ and }\ell
\left( M\right) <2^{-\mathbf{\rho }}\ell \left( I\right) }\left( \frac{%
\mathrm{P}^{\alpha }\left( M,\mathbf{1}_{I}\sigma \right) }{\left\vert
M\right\vert ^\frac{1}{n} }\right) ^{2}\sum_{\substack{ F\in \mathcal{F}:\ F\subset
3I\setminus I  \\ M\in \mathcal{W}\left( F\right) }}\left\Vert \mathsf{Q}%
_{F,M}^{\omega ,\mathbf{b}^{\ast }}x\right\Vert _{L^{2}\left( \omega \right)
}^{\spadesuit 2} \\
&&+\sum_{M\subset 3I\setminus I\text{ and }\ell \left( M\right) <2^{-\mathbf{%
\rho }}\ell \left( I\right) }\left( \frac{\mathrm{P}^{\alpha }\left( M,%
\mathbf{1}_{I}\sigma \right) }{\left\vert M\right\vert^\frac{1}{n}  }\right) ^{2}\sum 
_{\substack{ F\in \mathcal{F}:\ F\not\subset 3I\setminus I  \\ M\in \mathcal{%
W}\left( F\right) }}\left\Vert \mathsf{Q}_{F,M}^{\omega ,\mathbf{b}^{\ast
}}x\right\Vert _{L^{2}\left( \omega \right) }^{\spadesuit 2} \\
&\equiv &B_{{away}}^{1}+B_{{away}}^{2}\ .
\end{eqnarray*}%
\ \ 

To estimate $B_{{away}}^{1}$, let 
\begin{equation}
\mathcal{J}^{\ast }\equiv \bigcup\limits_{\substack{ F\in \mathcal{F}  \\ %
F\subset 3I\setminus I}}\bigcup\limits_{\substack{ M\in \mathcal{W}\left(
F\right)  \\ M\subset 3I\setminus I\text{ and }\ell \left( M\right) <2^{-%
\mathbf{\rho }}\ell \left( I\right) }}\left\{ J\in \mathcal{C}_{F}^{\mathcal{%
G},{shift}}:J\subset M\right\}  \label{def J*}
\end{equation}%
consist of all cubes $J\in \mathcal{G}$ for which the projection $%
\triangle _{J}^{\omega ,\mathbf{b}^{\ast }}$ occurs in one of the
projections $\mathsf{Q}_{F,M}^{\omega ,\mathbf{b}^{\ast }}$ in term $B_{%
{away}}^{1}$. In order to use $\mathcal{J}^{\ast }$ in the estimate for 
$B_{{away}}^{1}$ we need the following inequality. For any cube $%
M\in \mathcal{W}\left( F\right) $ we have%
\begin{eqnarray}
\left( \frac{\mathrm{P}^{\alpha }\left( M,\mathbf{1}_{I}\sigma \right) }{%
\left\vert M\right\vert }\right) ^{2}\left\Vert \mathsf{Q}_{F;M}^{\omega ,%
\mathbf{b}^{\ast }}x\right\Vert _{L^{2}\left( \omega \right) }^{\spadesuit
2}
\!\!\!\!&=&\!\!\!\!
\left( \frac{\mathrm{P}^{\alpha }\left( M,\mathbf{1}_{I}\sigma \right) 
}{\left\vert M\right\vert }\right) ^{2}\sum_{J\in \mathcal{C}_{F}^{\mathcal{G%
},{shift}}:\ J\subset M}\!\!\!\left\Vert \triangle _{J}^{\omega ,\mathbf{b}%
^{\ast }}x\right\Vert _{L^{2}\left( \omega \right) }^{\spadesuit 2}\notag \\
&\lesssim &\!\!\!\!
\label{accomplished}\sum_{J\in \mathcal{C}_{F}^{\mathcal{G},{shift}}:\ J\subset
M}\left( \frac{\mathrm{P}^{\alpha }\left( J,\mathbf{1}_{I}\sigma \right) }{%
\left\vert J\right\vert }\right) ^{2}\left\Vert \triangle _{J}^{\omega ,%
\mathbf{b}^{\ast }}x\right\Vert _{L^{2}\left( \omega \right) }^{\spadesuit
2}
\end{eqnarray}%
since%
\begin{eqnarray*}
\frac{\mathrm{P}^{\alpha }\left( M,\mathbf{1}_{I}\sigma \right) }{\left\vert
M\right\vert^\frac{1}{n}  } &=&\int_{I}\frac{1}{\left( \ell \left( M\right) +\left\vert
x-c_{M}\right\vert \right) ^{n+1-\alpha }}d\sigma \left( x\right) \\
&\lesssim &\int_{I}\frac{1}{\left( \ell \left( J\right) +\left\vert
x-c_{J}\right\vert \right) ^{n+1-\alpha }}d\sigma \left( x\right) =\frac{%
\mathrm{P}^{\alpha }\left( J,\mathbf{1}_{I}\sigma \right) }{\left\vert
J\right\vert^\frac{1}{n}  }
\end{eqnarray*}%
for $J\subset M$ because%
\begin{equation*}
\ell \left( J\right) +\left\vert x-c_{J}\right\vert \lesssim \ell \left(
M\right) +\left\vert x-c_{M}\right\vert ,\ \ \ \ \ J\subset M\text{ and }%
x\in \mathbb{R}^n.
\end{equation*}%
We now use (\ref{accomplished}) to replace the sum over $M\in \mathcal{W}%
\left( F\right) $ in $B_{{away}}^{1}$, with a sum over $J\in \mathcal{J}%
^{\ast }$:%
\begin{eqnarray*}
B_{{away}}^{1} &=&\sum_{M\subset 3I\setminus I\text{ and }\ell \left(
M\right) <2^{-\mathbf{\rho }}\ell \left( I\right) }\left( \frac{\mathrm{P}%
^{\alpha }\left( M,\mathbf{1}_{I}\sigma \right) }{\left\vert M\right\vert^\frac{1}{n}  }%
\right) ^{2}\sum_{\substack{ F\in \mathcal{F}:\ F\subset 3I\setminus I  \\ %
M\in \mathcal{W}\left( F\right) }}\left\Vert \mathsf{Q}_{F,M}^{\omega ,%
\mathbf{b}^{\ast }}x\right\Vert _{L^{2}\left( \omega \right) }^{\spadesuit 2}
\\
&\lesssim &\!\!\!\!
\sum_{M\subset 3I\setminus I \& \ell \left( M\right) <2^{-%
\mathbf{\rho }}\ell \left( I\right) }\sum_{\substack{ F\in \mathcal{F}:\
F\subset 3I\setminus I  \\ M\in \mathcal{W}\left( F\right) }}\sum_{\substack{J\in 
\mathcal{C}_{F}^{\mathcal{G},{shift}}\\ J\subset M}} \!\!\!\left( \frac{\mathrm{%
P}^{\alpha }\left( J,\mathbf{1}_{I}\sigma \right) }{\left\vert J\right\vert^\frac{1}{n}  }%
\right) ^{\!2}\!\!\left\Vert \triangle _{J}^{\omega ,\mathbf{b}^{\ast
}}x\right\Vert _{L^{2}\left( \omega \right) }^{\spadesuit 2} \\
&\lesssim &\sum_{J\in \mathcal{J}^{\ast }}\left( \frac{\mathrm{P}^{\alpha
}\left( J,\mathbf{1}_{I}\sigma \right) }{\left\vert J\right\vert }\right)
^{2}\mathbf{\ }\left\Vert \bigtriangleup _{J}^{\omega ,\mathbf{b}^{\ast
}}x\right\Vert _{L^{2}\left( \omega \right) }^{\spadesuit 2}\ ,
\end{eqnarray*}%
where the final line follows since for each $J\in \mathcal{J}^{\ast }$ there
is a unique pair $\left( F,M\right) $ satisfying the conditions in the
second line.\bigskip

We will now exploit the smallness of $\varepsilon >0$ in the weak goodness
condition by decomposing the sum over $J\in \mathcal{J}^{\ast }$ according
to the length of $J$, and then using the fractional Poisson inequality (\ref%
{e.Jsimeq}) in Lemma \ref{Poisson inequality} on the neighbour $I^{\prime }$%
\ of $I$ containing $J$. Indeed, for $J\subset I^{\prime }\subset \mathbb{R}$
and $I\subset \mathbb{R}\setminus I^{\prime }$, we have 
\begin{equation}
\mathrm{P}^{\alpha }\left( J,\mathbf{1}_{I}\sigma \right) ^{2}\lesssim
\left( \frac{\ell \left( J\right) }{\ell \left( I\right) }\right)
^{2-2\left( n+1-\alpha \right) \varepsilon }\mathrm{P}^{\alpha }\left( I,%
\mathbf{1}_{I}\sigma \right) ^{2},\ \ \ \ \ J\in \mathcal{J}^{\ast },
\label{Poisson inequalities 2}
\end{equation}%
where we have used that $\ell \left( I^{\prime }\right) =\ell \left(
I\right) $ and $\mathrm{P}^{\alpha }\left( I^{\prime },\mathbf{1}_{I}\sigma
\right) \approx \mathrm{P}^{\alpha }\left( I,\mathbf{1}_{I}\sigma \right) $,
and that the cubes $J\in \mathcal{J}^{\ast }$ are good in $I^{\prime }$
and beyond, and have side length at most $2^{-\mathbf{\rho }}\ell \left(
I\right) $, all because $J^{\maltese }\subset F\subset 3I\setminus I$ and we
have already dealt with the term $B_{{nearby}}$. Moreover, we may
also assume here that the exponent $2-2\left( n+1-\alpha \right) \varepsilon $
is positive, i.e.$\ \varepsilon <\frac{1}{n+1-\alpha }$, which is of course
implied by $0<\varepsilon <\frac{1}{2}$. We then obtain from (\ref{Poisson
inequalities 2}), the inequality $\left\Vert \bigtriangleup _{J}^{\omega ,%
\mathbf{b}^{\ast }}x\right\Vert _{L^{2}\left( \omega \right) }^{\spadesuit
2}\lesssim \left\vert J\right\vert ^{2}\left\vert J\right\vert _{\omega }$,
the pairwise disjointedness of the $M\in \mathcal{W}\left( F\right) $, the
uniqueness of $F$ with $J\in \mathcal{C}_{F}^{\mathcal{G},{shift}}$,
and since $F\subset 3I\setminus I$ in the sum over $J\in \mathcal{J}^{\ast }$%
, that%
\begin{eqnarray*}
B_{{away}}^{1} &\lesssim &\sum_{J\in \mathcal{J}^{\ast }}\left( \frac{%
\mathrm{P}^{\alpha }\left( J,\mathbf{1}_{I}\sigma \right) }{\left\vert
J\right\vert^\frac{1}{n}  }\right) ^{2}\mathbf{\ }\left\Vert \bigtriangleup _{J}^{\omega ,%
\mathbf{b}^{\ast }}x\right\Vert _{L^{2}\left( \omega \right) }^{\spadesuit
2} \\ 
&\lesssim&
\sum_{m=\mathbf{\rho }}^{\infty }\sum_{\substack{ J\in \mathcal{J}%
^{\ast }  \\ \ell \left( J\right) =2^{-m}\ell \left( I\right) }}\left(
2^{-m}\right) ^{2-2\left( n+1-\alpha \right) \varepsilon }\mathrm{P}^{\alpha
}\left( I,\mathbf{1}_{I}\sigma \right) ^{2}\left\vert J\right\vert _{\omega }
\\
&\lesssim &
\sum_{m=\mathbf{\rho }}^{\infty }\left( 2^{-m}\right) ^{2-2\left(
n+1-\alpha \right) \varepsilon }\left( \frac{\left\vert I\right\vert _{\sigma }%
}{\left\vert I\right\vert ^{1-\frac{\alpha}{n} }}\right) ^{2}\sum_{\substack{ J\subset
3I\setminus I  \\ \ell \left( J\right) =2^{-m}\ell \left( I\right) }}%
\left\vert J\right\vert _{\omega } \\ 
&\lesssim&
\sum_{m=\mathbf{\rho }}^{\infty
}\left( 2^{-m}\right) ^{2-2\left( n+1-\alpha \right) \varepsilon }\frac{%
\left\vert I\right\vert _{\sigma }\left\vert 3I\setminus I\right\vert
_{\omega }}{\left\vert 3I\right\vert ^{2\left( 1-\frac{\alpha}{n} \right) }}\left\vert
I\right\vert _{\sigma }\lesssim A_{2}^{\alpha }\left\vert I\right\vert
_{\sigma }\ ,
\end{eqnarray*}%
since $2-2\left( n+1-\alpha \right) \varepsilon >0$.

To complete the bound for term $B=B_{{nearby}}+B_{{away}%
}^{1}+B_{{away}}^{2}$, it remains to estimate term $B_{{away}}^{2}$
in which we sum over $F\not\subset 3I\setminus I$. In this case $%
F\varsupsetneqq I^{\prime }$ for one of the two neighbours $I^{\prime }$ of $%
I$, and so we can apply Lemma \ref{refined lemma}, with $I$ there replaced
by the augmented cubes $I^{\prime }\cup I$, to obtain the estimate%
\begin{equation*}
B_{{away}}^{2}\lesssim \left( \left( \mathfrak{E}_{2}^{\alpha }\right)
^{2}+A_{2}^{\alpha ,{energy}}\right) \left\vert I\right\vert
_{\sigma }\ .
\end{equation*}

Next we turn to term $D$. The cubes $M$ occurring here are included in
the set of ancestors $A_{k}\equiv \pi _{\mathcal{D}}^{\left( k\right) }I$ of 
$I$, $1\leq k<\infty $. Then $D$ is equal to
\begin{eqnarray*}
&&
\sum_{k=1}^{\infty }\mathbb{P}^{\alpha }\left( \mathbf{1}_{I}\sigma
\right) \left( c\left( A_{k}\right) ,\left\vert A_{k}\right\vert \right)
^{2}\sum_{\substack{ F\in \mathcal{F}:  \\ A_{k}\in \mathcal{W}\left(
F\right) }}\left\Vert \mathsf{Q}_{F,A_{k}}^{\omega ,\mathbf{b}^{\ast }}\frac{x}{\lvert A_{k}\rvert^\frac{1}{n} }\right\Vert _{L^{2}\left( \omega \right)
}^{\spadesuit 2} \\
&=&
\sum_{k=1}^{\infty }\mathbb{P}^{\alpha }\left( \mathbf{1}_{I}\sigma
\right) \left( c\left( A_{k}\right) ,\left\vert A_{k}\right\vert \right)
^{2}\!\!\!\sum_{\substack{ F\in \mathcal{F}:  \\ A_{k}\in \mathcal{W}\left(
F\right) }}\sum_{J^{\prime }\in \mathcal{C}_{F}^{\mathcal{G},{shift}%
}:\ J^{\prime }\subset A_{k}\setminus I}\left\Vert \bigtriangleup
_{J^{\prime }}^{\omega ,\mathbf{b}^{\ast }}\frac{x}{\lvert A_{k}\rvert^\frac{1}{n}  }%
\right\Vert _{L^{2}\left( \omega \right) }^{\spadesuit 2} \\
&&
+\sum_{k=1}^{\infty }\mathbb{P}^{\alpha }\left( \mathbf{1}_{I}\sigma
\right) \left( c\left( A_{k}\right) ,\left\vert A_{k}\right\vert \right)
^{2}\!\!\!\sum_{\substack{ F\in \mathcal{F}:  \\ A_{k}\in \mathcal{W}\left(
F\right) }}\sum_{\substack{ J^{\prime }\in \mathcal{C}_{F}^{\mathcal{G},%
{shift}}:\ J^{\prime }\subset A_{k}  \\ J^{\prime }\cap I\neq
\emptyset \text{ and }\ell \left( J^{\prime }\right) \leq \ell \left(
I\right) }}\left\Vert \bigtriangleup _{J^{\prime }}^{\omega ,\mathbf{b}%
^{\ast }}\frac{x}{\lvert A_{k}\rvert^\frac{1}{n}  }\right\Vert _{L^{2}\left( \omega
\right) }^{\spadesuit 2} \\
&&+\sum_{k=1}^{\infty }\mathbb{P}^{\alpha }\left( \mathbf{1}_{I}\sigma
\right) \left( c\left( A_{k}\right) ,\left\vert A_{k}\right\vert \right)
^{2}\!\!\!\sum_{\substack{ F\in \mathcal{F}:  \\ A_{k}\in \mathcal{W}\left(
F\right) }}\sum_{\substack{ J^{\prime }\in \mathcal{C}_{F}^{\mathcal{G},%
{shift}}:\ J^{\prime }\subset A_{k}  \\ J^{\prime }\cap I\neq
\emptyset \text{ and }\ell \left( J^{\prime }\right) >\ell \left( I\right) }}%
\left\Vert \bigtriangleup _{J^{\prime }}^{\omega ,\mathbf{b}^{\ast }}\frac{x%
}{\lvert A_{k}\rvert^\frac{1}{n}  }\right\Vert _{L^{2}\left( \omega \right) }^{\spadesuit
2} \\
&\equiv &D_{{disjoint}}+D_{{descendent}}+D_{{ancestor%
}}\ .
\end{eqnarray*}%
We thus have from the pairwise disjointedness of the projections $\mathsf{Q}%
_{F,A_{k}}^{\omega ,\mathbf{b}^{\ast }}$ in $F$ once again, $D_{{disjoint}}$ equals
\begin{eqnarray*}
&&\sum_{k=1}^{\infty }\mathbb{P}^{\alpha }\left( 
\mathbf{1}_{I}\sigma \right) \left( c\left( A_{k}\right) ,\left\vert
A_{k}\right\vert \right) ^{2}\!\!\!\!\sum_{\substack{ F\in \mathcal{F}:  \\ A_{k}\in 
\mathcal{W}\left( F\right) }}\sum_{\substack{J^{\prime }\in \mathcal{C}_{F}^{\mathcal{G%
},{shift}}:\\ J^{\prime }\subset A_{k}\setminus I}}\left\Vert
\bigtriangleup _{J^{\prime }}^{\omega ,\mathbf{b}^{\ast }}\frac{x}{\lvert
A_{k}\rvert^\frac{1}{n}  }\right\Vert _{L^{2}\left( \omega \right) }^{\spadesuit 2} \\
&\lesssim &
\sum_{k=1}^{\infty }\left( \frac{\left\vert I\right\vert _{\sigma
}\left\vert A_{k}\right\vert }{\left\vert A_{k}\right\vert ^{n+1-\alpha }}%
\right) ^{2}\mathbf{\;}\left\vert A_{k}\setminus I\right\vert _{\omega
}
=
\left\{ \frac{\left\vert I\right\vert _{\sigma }}{\left\vert I\right\vert
^{1-\frac{\alpha}{n} }}\sum_{k=1}^{\infty }\frac{\left\vert I\right\vert ^{1-\frac{\alpha}{n} }}{%
\left\vert A_{k}\right\vert ^{2\left( n-\alpha \right) }}\left\vert
A_{k}\setminus I\right\vert _{\omega }\right\} \left\vert I\right\vert
_{\sigma } \\
&\lesssim &
\left\{ \frac{\left\vert I\right\vert _{\sigma }}{\left\vert
I\right\vert ^{1-\frac{\alpha}{n} }}\mathcal{P}^{\alpha }\left( I,\mathbf{1}%
_{I^{c}}\omega \right) \right\} \left\vert I\right\vert _{\sigma }\lesssim 
\mathcal{A}_{2}^{\alpha ,\ast }\left\vert I\right\vert _{\sigma },
\end{eqnarray*}%
since%
\begin{eqnarray*}
\sum_{k=1}^{\infty }\frac{\left\vert I\right\vert ^{1-\frac{\alpha}{n} }}{\left\vert
A_{k}\right\vert ^{2\left( n-\alpha \right) }}\left\vert A_{k}\setminus
I\right\vert _{\omega } 
&=&
\int \sum_{k=1}^{\infty }\frac{\left\vert
I\right\vert ^{1-\frac{\alpha}{n} }}{\left\vert A_{k}\right\vert ^{2\left( n-\alpha
\right) }}\mathbf{1}_{A_{k}\setminus I}\left( x\right) d\omega \left(
x\right) \\
&=&
\int \sum_{k=1}^{\infty }\frac{1}{2^{2\left( n-\alpha \right) k}}\frac{%
\left\vert I\right\vert ^{1-\frac{\alpha}{n} }}{\left\vert I\right\vert ^{2\left(
n-\alpha \right) }}\mathbf{1}_{A_{k}\setminus I}\left( x\right) d\omega
\left( x\right) \\
&\lesssim &
\int_{I^{c}}\left( \frac{\left\vert I\right\vert^\frac{1}{n}  }{\left[
\left\vert I\right\vert^\frac{1}{n}  +d\left( x,I\right) \right] ^{2}}\right) ^{n-\alpha
}d\omega \left( x\right) =\mathcal{P}^{\alpha }\left( I,\mathbf{1}%
_{I^{c}}\omega \right) ,
\end{eqnarray*}%
upon summing a geometric series with $2\left( n-\alpha \right) >0$.

The next term $D_{{descendent}}$ satisfies%
\begin{eqnarray*}
D_{{descendent}} &\lesssim &\sum_{k=1}^{\infty }\left( \frac{%
\left\vert I\right\vert _{\sigma }\left\vert A_{k}\right\vert }{\left\vert
A_{k}\right\vert ^{n+1-\alpha }}\right) ^{2}\mathbf{\;}\left\Vert \mathsf{Q}%
_{3I}^{\omega ,\mathbf{b}^{\ast }}\frac{x}{2^{k}\lvert I\rvert^\frac{1}{n}  }\right\Vert
_{L^{2}\left( \omega \right) }^{\spadesuit 2} \\
&=&
\sum_{k=1}^{\infty }2^{-2k\left( n+1-\alpha \right) }\left( \frac{%
\left\vert I\right\vert _{\sigma }}{\left\vert I\right\vert ^{n-\alpha }}%
\right) ^{2}\left\Vert \mathsf{Q}_{3I}^{\omega ,\mathbf{b}^{\ast }}\frac{x}{%
\lvert I\rvert^\frac{1}{n}  }\right\Vert _{L^{2}\left( \omega \right) }^{\spadesuit 2} \\
&\lesssim &
\left\{ \frac{\left\vert I\right\vert _{\sigma }\left\Vert 
\mathsf{Q}_{3I}^{\omega ,\mathbf{b}^{\ast }}\frac{x}{\lvert I\rvert^\frac{1}{n}  }%
\right\Vert _{L^{2}\left( \omega \right) }^{\spadesuit 2}}{\left\vert
I\right\vert ^{2\left(n-\alpha \right) }}\right\} \left\vert I\right\vert
_{\sigma }\lesssim A_{2}^{\alpha ,{energy}}\left\vert I\right\vert
_{\sigma }\ .
\end{eqnarray*}

Lastly, for $D_{{ancestor}}$ we note that there are at most two
cubes $K_{1}$ and $K_{2}$ in $\mathcal{G}$ having side length $\ell
\left( I\right) $ and such that $K_{i}\cap I\neq \emptyset $. Then each $%
J^{\prime }$ occurring in the sum in $D_{{ancestor}}$ is of the form 
$J^{\prime }=A_{i}^{\ell }\equiv \pi _{\mathcal{G}}^{\left( \ell \right)
}K_{i}$ with $J^{\prime }\subset A_{k}$ for some $1\leq \ell \leq k$ and $%
i\in \left\{ 1,2\right\} $. Now we write%
\begin{eqnarray*}
D_{{ancestor}} \!\!\!\!&=&\!\!\!\!\sum_{k=1}^{\infty }\mathbb{P}^{\alpha }\left( 
\mathbf{1}_{I}\sigma \right) \left( c\left( A_{k}\right) ,\left\vert
A_{k}\right\vert \right) ^{2}\!\!\!\!\!\!\sum_{\substack{ F\in \mathcal{F}:  \\ A_{k}\in 
\mathcal{W}\left( F\right) }}\!\!\!\sum_{\substack{ J^{\prime }\in \mathcal{C}%
_{F}^{\mathcal{G},{shift}}:\ J^{\prime }\subset A_{k}  \\ J^{\prime
}\cap I\neq \emptyset \text{ and }\ell \left( J^{\prime }\right) >\ell
\left( I\right) }}\!\!\!\!\!\! \!\left\Vert \bigtriangleup _{J^{\prime }}^{\omega ,\mathbf{b%
}^{\ast }}\frac{x}{\lvert A_{k}\rvert^\frac{1}{n}  }\right\Vert _{L^{2}\left( \omega
\right) }^{\spadesuit 2} \\
&\lesssim &\!\!\!\!
\sum_{k=1}^{\infty }\left( \frac{\left\vert I\right\vert _{\sigma
}\left\vert A_{k}\right\vert }{\left\vert A_{k}\right\vert ^{n+1-\alpha }}%
\right) ^{2}\sum_{i=1}^{2}\sum_{\ell =1}^{k}\left\Vert \bigtriangleup
_{A_{i}^{\ell }}^{\omega ,\mathbf{b}^{\ast }}\frac{x}{\lvert A_{k}\rvert^\frac{1}{n}  }%
\right\Vert _{L^{2}\left( \omega \right) }^{\spadesuit 2} \\
&\leq &2\sum_{k=1}^{\infty }\left( \frac{\left\vert I\right\vert _{\sigma
}\left\vert A_{k}\right\vert }{\left\vert A_{k}\right\vert ^{n+1-\alpha }}%
\right) ^{2}\left\Vert \mathsf{Q}_{A_{k}}^{\omega ,\mathbf{b}^{\ast }}\frac{x%
}{\lvert A_{k}\rvert^\frac{1}{n}  }\right\Vert _{L^{2}\left( \omega \right) }^{\spadesuit
2}.
\end{eqnarray*}%
At this point we need a \emph{`prepare to puncture'} argument, as we will
want to derive geometric decay from $\left\Vert \mathsf{Q}_{J^{\prime
}}^{\omega ,\mathbf{b}^{\ast }}x\right\Vert _{L^{2}\left( \omega \right)
}^{\spadesuit 2}$ by dominating it by the `nonenergy' term $\left\vert
J^{\prime }\right\vert ^{2}\left\vert J^{\prime }\right\vert _{\omega }$, as
well as using the Muckenhoupt energy constant. For this we define $%
\widetilde{\omega }=\omega -\omega \left( \left\{ p\right\} \right) \delta
_{p}$ where $p$ is an atomic point in $I$ for which\\ 
$\displaystyle
\omega \left( \left\{ p\right\} \right) =\sup_{q\in \mathfrak{P}_{\left(
\sigma ,\omega \right) }:\ q\in I}\omega \left( \left\{ q\right\} \right) .
$
(If $\omega $ has no atomic point in common with $\sigma $ in $I$ set $%
\widetilde{\omega }=\omega $.) Then we have $\left\vert I\right\vert _{%
\widetilde{\omega }}=\omega \left( I,\mathfrak{P}_{\left( \sigma ,\omega
\right) }\right) $ and%
\begin{equation*}
\frac{\left\vert I\right\vert _{\widetilde{\omega }}}{\left\vert
I\right\vert ^{1-\frac{\alpha}{n} }}\frac{\left\vert I\right\vert _{\sigma }}{%
\left\vert I\right\vert ^{1-\frac{\alpha}{n} }}=\frac{\omega \left( I,\mathfrak{P}%
_{\left( \sigma ,\omega \right) }\right) }{\left\vert I\right\vert
^{1-\frac{\alpha}{n} }}\frac{\left\vert I\right\vert _{\sigma }}{\left\vert
I\right\vert ^{1-\frac{\alpha}{n} }}\leq A_{2}^{\alpha ,{punct}}.
\end{equation*}%
A key observation, already noted in the proof of Lemma \ref{energy A2}
above, is that%
\begin{equation}
\left\Vert \bigtriangleup _{K}^{\omega ,\mathbf{b}^{\ast }}x\right\Vert
_{L^{2}\left( \omega \right) }^{2}=\left\{ 
\begin{array}{ccc}
\left\Vert \bigtriangleup _{K}^{\omega ,\mathbf{b}^{\ast }}\left( x-p\right)
\right\Vert _{L^{2}\left( \omega \right) }^{2} & \text{ if } & p\in K \\ 
\left\Vert \bigtriangleup _{K}^{\omega ,\mathbf{b}^{\ast }}x\right\Vert
_{L^{2}\left( \widetilde{\omega }\right) }^{2} & \text{ if } & p\notin K%
\end{array}%
\right. \leq \ell \left( K\right) ^{2}\left\vert K\right\vert _{\widetilde{%
\omega }}, \label{key obs}
\end{equation}%
and so, as in the proof of (\ref{omega tilda}) in Lemma \ref{energy A2},%
\begin{equation*}
\left\Vert \mathsf{Q}_{A_{k}}^{\omega ,\mathbf{b}^{\ast }}\frac{x}{%
\left\vert A_{k}\right\vert^\frac{1}{n}  }\right\Vert _{L^{2}\left( \omega \right)
}^{\spadesuit 2}\lesssim \left\vert A_{k}\right\vert _{\widetilde{\omega }}\
.
\end{equation*}%
Then we continue with%
\begin{eqnarray*}
&&\sum_{k=1}^{\infty }\left( \frac{\left\vert I\right\vert _{\sigma
}\left\vert A_{k}\right\vert }{\left\vert A_{k}\right\vert ^{n+1-\alpha }}%
\right) ^{2}\left\Vert \mathsf{Q}_{A_{k}}^{\omega ,\mathbf{b}^{\ast }}\frac{x%
}{\lvert A_{k}\rvert^\frac{1}{n}  }\right\Vert _{L^{2}\left( \omega \right) }^{\spadesuit
2} \\
&\lesssim &\sum_{k=1}^{\infty }\left( \frac{\left\vert I\right\vert _{\sigma
}\left\vert A_{k}\right\vert }{\left\vert A_{k}\right\vert ^{n+1-\alpha }}%
\right) ^{2}\left\vert A_{k}\right\vert _{\widetilde{\omega }} \\
&=&
\sum_{k=1}^{\infty }\left( \frac{\left\vert I\right\vert _{\sigma }}{%
\left\vert A_{k}\right\vert ^{n-\alpha }}\right) ^{2}\left\vert
A_{k}\setminus I\right\vert _{\omega }+\sum_{k=1}^{\infty }\left( \frac{%
\left\vert I\right\vert _{\sigma }}{2^{k\left( n-\alpha \right) }\left\vert
I\right\vert ^{n-\alpha }}\right) ^{2}\left\vert I\right\vert _{\widetilde{%
\omega }} \\
&\lesssim &\left( \mathcal{A}_{2}^{\alpha ,\ast }+A_{2}^{\alpha ,{%
punct}}\right) \left\vert I\right\vert _{\sigma },
\end{eqnarray*}%
where the inequality $\sum_{k=1}^{\infty }\left( \frac{\left\vert
I\right\vert _{\sigma }}{\left\vert A_{k}\right\vert ^{n-\alpha }}\right)
^{2}\left\vert A_{k}\setminus I\right\vert _{\omega }\lesssim \mathcal{A}%
_{2}^{\alpha ,\ast }\left\vert I\right\vert _{\sigma }$ is already proved
above in the display estimating $D_{{disjoint}}$.

\medskip

Finally, for term $C$ we will have to group the cubes $M$ into blocks $%
B_{i}$. We first split the sum according to whether or not $I$ intersects
the triple of $M$:%
\begin{eqnarray*}
C &\approx &\left\{ \sum_{\substack{ M:\ I\cap 3M=\emptyset  \\ \ell \left(
M\right) >\ell \left( I\right) }}+\sum_{\substack{ M:\ I\subset 3M\setminus
M  \\ \ell \left( M\right) >\ell \left( I\right) }}\right\} \left( \frac{%
\left\vert M\right\vert^\frac{1}{n}  }{\left( \left\vert M\right\vert^\frac{1}{n}  +d\left( M,I\right)
\right) ^{n+1-\alpha }}\left\vert I\right\vert _{\sigma }\right) ^{2} \cdot \\ && 
\hspace{3cm}\cdot\sum 
_{\substack{ F\in \mathcal{F}:  \\ M\in \mathcal{W}\left( F\right) }}%
\left\Vert \mathsf{Q}_{F,M}^{\omega \mathbf{b}^{\ast }}\frac{x}{\left\vert
M\right\vert^\frac{1}{n}  }\right\Vert _{L^{2}\left( \omega \right) }^{\spadesuit 2} \\
&=&C_{1}+C_{2}.
\end{eqnarray*}%
We first consider $C_{1}$. Let $\mathcal{M}$ consist of the maximal dyadic
cubes in the collection $\left\{ Q:3Q\cap I=\emptyset \right\} $, and
then let $\left\{ B_{i}\right\} _{i=1}^{\infty }$ be an enumeration of those 
$Q\in \mathcal{M}$ whose side length is at least $\ell \left( I\right) $.
Note in particular that $3B_{i}\cap I=\emptyset $. Now we further decompose
the sum in $C_{1}$ by grouping the cubes $M$ into the `Whitney'
cubes $B_{i}$, and then using the pairwise disjointedness of the
martingale supports of the pseudoprojections $\mathsf{Q}_{F,M}^{\omega ,%
\mathbf{b}^{\ast }}$ in $F$: 
\begin{eqnarray*}
C_{1}
&\leq &
\sum_{i=1}^{\infty }\sum_{M:\ M\subset B_{i}}\left( \frac{1}{%
\left( \left\vert M\right\vert^\frac{1}{n}  +d\left( M,I\right) \right) ^{n+1-\alpha }}%
\left\vert I\right\vert _{\sigma }\right) ^{2}\sum_{\substack{ F\in \mathcal{%
F}:  \\ M\in \mathcal{W}\left( F\right) }}\left\Vert \mathsf{Q}%
_{F,M}^{\omega ,\mathbf{b}^{\ast }}x\right\Vert _{L^{2}\left( \omega \right)
}^{\spadesuit 2} \\
&\lesssim &
\sum_{i=1}^{\infty }\left( \frac{1}{\left( \left\vert
B_{i}\right\vert^\frac{1}{n}  +d\left( B_{i},I\right) \right) ^{n+1-\alpha }}\left\vert
I\right\vert _{\sigma }\right) ^{2}\sum_{M:\ M\subset B_{i}}\sum_{\substack{ %
F\in \mathcal{F}:  \\ M\in \mathcal{W}\left( F\right) }}\left\Vert \mathsf{Q}%
_{F,M}^{\omega ,\mathbf{b}^{\ast }}x\right\Vert _{L^{2}\left( \omega \right)
}^{\spadesuit 2} \\
&\lesssim &
\sum_{i=1}^{\infty }\left( \frac{1}{\left( \left\vert
B_{i}\right\vert^\frac{1}{n}  +d\left( B_{i},I\right) \right) ^{n+1-\alpha }}\left\vert
I\right\vert _{\sigma }\right) ^{2}\sum_{M:\ M\subset B_{i}}\left\vert
M\right\vert ^{^\frac{2}{n} }\left\vert M\right\vert _{\omega } \\
&\lesssim &
\sum_{i=1}^{\infty }\left( \frac{1}{\left( \left\vert
B_{i}\right\vert^\frac{1}{n}  +d\left( B_{i},I\right) \right) ^{n+1-\alpha }}\left\vert
I\right\vert _{\sigma }\right) ^{2}\mathbf{\ }\left\vert B_{i}\right\vert
^{^\frac{2}{n}} \left\vert B_{i}\right\vert _{\omega } \\
&\lesssim &
\left\{ \sum_{i=1}^{\infty }\frac{\left\vert B_{i}\right\vert
_{\omega }\left\vert I\right\vert _{\sigma }}{\left\vert B_{i}\right\vert
^{2\left( n-\alpha \right) }}\right\} \left\vert I\right\vert _{\sigma }\ ,
\end{eqnarray*}%
Now since $\left\vert B_{i}\right\vert \approx d\left( x,I\right) $ for $%
x\in B_{i}$, 

\begin{eqnarray*}
\sum_{i=1}^{\infty }\frac{\left\vert B_{i}\right\vert _{\omega }\left\vert
I\right\vert _{\sigma }}{\left\vert B_{i}\right\vert ^{2\left( n-\alpha
\right) }} 
\!\!\!\!\!&=&\!\!\!\!
\frac{\left\vert I\right\vert _{\sigma }}{\left\vert
I\right\vert ^{1-\alpha }}\sum_{i=1}^{\infty }\frac{\left\vert I\right\vert
^{1-\alpha }}{\left\vert B_{i}\right\vert ^{2\left( n-\alpha \right) }}%
\left\vert B_{i}\right\vert _{\omega } 
\approx 
\frac{\left\vert I\right\vert _{\sigma }}{\left\vert I\right\vert
^{1-\frac{\alpha}{n} }}\sum_{i=1}^{\infty }\int_{B_{i}}\frac{\left\vert I\right\vert
^{1-\frac{\alpha}{n} }}{d\left( x,I\right) ^{2\left( n-\alpha \right) }}d\omega \left(
x\right) \\
&\approx &
\frac{\left\vert I\right\vert _{\sigma }}{\left\vert I\right\vert
^{1-\frac{\alpha}{n} }}\sum_{i=1}^{\infty }\int_{B_{i}}\left( \frac{\left\vert
I\right\vert^\frac{1}{n}  }{\left[ \left\vert I\right\vert^\frac{1}{n}  +d\left( x,I\right) \right]
^{2}}\right) ^{n-\alpha }d\omega \left( x\right) \\
&\leq &
\frac{\left\vert I\right\vert _{\sigma }}{\left\vert I\right\vert
^{1-\frac{\alpha}{n} }}\mathcal{P}^{\alpha }\left( I,\mathbf{1}_{I^{c}}\omega \right)
\leq \mathcal{A}_{2}^{\alpha ,\ast },
\end{eqnarray*}%
we obtain%
\begin{equation*}
C_{1}\lesssim \mathcal{A}_{2}^{\alpha ,\ast }\left\vert I\right\vert
_{\sigma }\ .
\end{equation*}

Next we turn to estimating term $C_{2}$ where the triple of $M$ contains $I$
but $M$ itself does not. Note that there are at most two such cubes $M$
of a given side length. So with this in mind, we sum over the cubes $M$
according to their lengths to obtain%
\begin{eqnarray*}
C_{2}
\!\!\!\!&=&\!\!\!\!
\sum_{m=1}^{\infty }\!\!\!\!\!\sum_{\substack{ M:\\ I\subset 3M\setminus M  \\ %
\ell \left( M\right) =2^{m}\ell \left( I\right) }}\!\!\!\!\!\left( \frac{\left\vert
M\right\vert^\frac{1}{n}  }{\left( \left\vert M\right\vert^\frac{1}{n}  +\dist\left(
M,I\right) \right) ^{n+1-\alpha }}\left\vert I\right\vert _{\sigma }\right)
^{2}\!\!\!\!\!\sum_{\substack{ F\in \mathcal{F}:  \\ M\in \mathcal{W}\left( F\right) }}%
\left\Vert \mathsf{Q}_{F,M}^{\omega ,\mathbf{b}^{\ast }}\frac{x}{\left\vert
M\right\vert^\frac{1}{n}  }\right\Vert _{L^{2}\left( \omega \right) }^{\spadesuit 2} \\
&\lesssim &
\!\!\!\!\sum_{m=1}^{\infty }\left( \frac{\left\vert I\right\vert _{\sigma
}}{\left\vert 2^{m}I\right\vert ^{n-\alpha }}\right) ^{\!\!\!2}
\!\left\vert \left( 5\!\cdot\! 2^{m}I\right)\!\backslash I\right\vert _{\omega
}=\left\{ \frac{\left\vert I\right\vert _{\sigma }}{\left\vert I\right\vert
^{1-\frac{\alpha}{n} }}\sum_{m=1}^{\infty }\frac{\left\vert I\right\vert ^{1-\frac{\alpha}{n}
}\left\vert \left( 5\cdot 2^{m}I\right) \setminus I\right\vert _{\omega }}{%
\left\vert 2^{m}I\right\vert ^{2\left( n-\alpha \right) }}\right\}
\left\vert I\right\vert _{\sigma } \\
&\lesssim &
\!\!\!\!
\left\{ \frac{\left\vert I\right\vert _{\sigma }}{\left\vert
I\right\vert ^{1-\frac{\alpha}{n} }}\mathcal{P}^{\alpha }\left( I,\mathbf{1}%
_{I^{c}}\omega \right) \right\} \left\vert I\right\vert _{\sigma }
\leq 
\mathcal{A}_{2}^{\alpha ,\ast }\left\vert I\right\vert _{\sigma },
\end{eqnarray*}%
since in analogy with the corresponding estimate above,%
\begin{equation*}
\sum_{m=1}^{\infty }\frac{\left\vert I\right\vert ^{1-\frac{\alpha}{n} }\left\vert
\left( 5\!\cdot\! 2^{m}I\right)\! \setminus\! I\right\vert _{\omega }}{\left\vert
2^{m}I\right\vert ^{2\left( n-\alpha \right) }}=\int \sum_{m=1}^{\infty }%
\frac{\left\vert I\right\vert ^{1-\frac{\alpha}{n} }}{\left\vert 2^{m}I\right\vert
^{2\left( n-\alpha \right) }}\mathbf{1}_{\left( 5\cdot 2^{m}I\right)
\setminus I}\!\!\left( x\right) \ d\omega\!\! \left( x\right) \lesssim \mathcal{P}%
^{\alpha }\left( I,\mathbf{1}_{I^{c}}\omega \right) .
\end{equation*}

\subsection{The backward Poisson testing inequality\label{Subsec back test}}

The argument here follows the broad outline of the analogous argument in 
\cite{SaShUr7}, but using modifications from \cite{SaShUr9} that involve
`prepare to puncture arguments', using decompositions $\mathcal{W}\left(
F\right) $ in place of $\left( \mathbf{\rho },\varepsilon \right) $%
-decompositions, and using pseudoprojections $\mathsf{Q}_{F,M}^{\omega ,%
\mathbf{b}^{\ast }}x$ (see (\ref{def F,K}) for the definition). The final
change here is that there is no splitting into local and global parts as in 
\cite{SaShUr7} - instead, we follow the treatment in \cite{SaShUr6} in this
regard.

Fix $I\in \mathcal{D}$. It suffices to prove%
\begin{eqnarray}
\mathbf{Back}\left( \widehat{I}\right) 
&\equiv& \int_{\mathbb{R}^n}\left[ 
\mathbb{Q}^{\alpha }\left( t\mathbf{1}_{\widehat{I}}\overline{\mu }\right)
\left( y\right) \right] ^{2}d\sigma (y)
\notag\\
&\lesssim& \left\{ \mathcal{A}%
_{2}^{\alpha }+\left( \mathfrak{E}_{2}^{\alpha }+\sqrt{A_{2}^{\alpha ,%
{energy}}}\right) \sqrt{A_{2}^{\alpha ,{punct}}}\right\}
\int_{\widehat{I}}t^{2}d\overline{\mu }(x,t).  \label{e.t2 n'}
\end{eqnarray}%
Note that for a `Poisson integral with holes' and a measure $\mu $ built
with Haar projections, Hyt\"{o}nen obtained in \cite{Hyt2} the simpler bound 
$A_{2}^{\alpha }$ for a term analogous to, but significantly smaller than, (%
\ref{e.t2 n'}). 
Using (\ref{tent consequence}) we see that the integral on the right hand
side of (\ref{e.t2 n'}) is 
\begin{equation}
\int_{\widehat{I}}t^{2}d\overline{\mu }=\sum_{F\in \mathcal{F}}\sum_{M\in 
\mathcal{W}\left( F\right) :\ M\subset I}\lVert \mathsf{Q}_{F,M}^{\omega ,%
\mathbf{b}^{\ast }}x\rVert _{L^{2}\left( \omega \right) }^{\spadesuit 2}\,.
\label{mu I hat}
\end{equation}%
where $\mathsf{Q}_{F,M}^{\omega ,\mathbf{b}^{\ast }}$ was defined in (\ref%
{def F,K}).

We now compute using (\ref{tent consequence}) again that 
\begin{eqnarray}
\mathbb{Q}^{\alpha }\left( t\mathbf{1}_{\widehat{I}}\overline{\mu }\right)
\left( y\right) &=&\int_{\widehat{I}}\frac{t^{2}}{\left( t^{2}+\left\vert
x-y\right\vert ^{2}\right) ^{\frac{n+1-\alpha }{2}}}d\overline{\mu }\left(
x,t\right)  \label{PI hat} \\
&\approx &\sum_{F\in \mathcal{F}}\sum_{M\in \mathcal{W}\left( F\right) :\
M\subset I}\frac{\lVert \mathsf{Q}_{F,M}^{\omega ,\mathbf{b}^{\ast }}x\rVert
_{L^{2}\left( \omega \right) }^{\spadesuit 2}}{\left( \left\vert
M\right\vert +\left\vert y-c_{M}\right\vert \right) ^{n+1-\alpha }},  \notag
\end{eqnarray}%
and then expand the square and integrate to obtain that the term $\mathbf{%
Back}\left( \widehat{I}\right) $ is 
\begin{equation*}
\sum_{\substack{ F\in \mathcal{F}  \\ M\in \mathcal{W}\left( F\right)  \\ %
M\subset I}}\sum_{\substack{ F^{\prime }\in \mathcal{F}:  \\ M^{\prime }\in 
\mathcal{W}\left( F^{\prime }\right)  \\ M^{\prime }\subset I}}\int_{\mathbb{%
R}}\frac{\left\Vert \mathsf{Q}_{F,M}^{\omega ,\mathbf{b}^{\ast
}}x\right\Vert _{L^{2}\left( \omega \right) }^{\spadesuit 2}}{\left(
\left\vert M\right\vert +\left\vert y-c_{M}\right\vert \right) ^{n+1-\alpha }}%
\frac{\left\Vert \mathsf{Q}_{F^{\prime },M^{\prime }}^{\omega ,\mathbf{b}%
^{\ast }}x\right\Vert _{L^{2}\left( \omega \right) }^{\spadesuit 2}}{\left(
\left\vert M^{\prime }\right\vert +\left\vert y-c_{M^{\prime }}\right\vert
\right) ^{n+1-\alpha }}d\sigma \left( y\right) .
\end{equation*}

By symmetry we may assume that $\ell \left( M^{\prime }\right) \leq \ell
\left( M\right) $. We fix a nonnegative integer $s$, and consider those
cubes $M$ and $M^{\prime }$ with $\ell \left( M^{\prime }\right)
=2^{-s}\ell \left( M\right) $. For fixed $s$ we will control the expression 
\begin{eqnarray}
&& \label{def Us}\\
U_{s} \!\!\!\!
&\equiv &\!\!\!\!\!\!\!\!\!
\sum_{\substack{ F,F^{\prime }\in \mathcal{F}}}\!\!\!\!\!\!\!\!\!\!\!\!\!\!\sum 
_{\substack{ M\in \mathcal{W}\left( F\right) \\ M^{\prime }\in \mathcal{W}%
\left( F^{\prime }\right)  \\ M,M^{\prime }\subset I,\ \ell \left( M^{\prime
}\right) =2^{-s}\ell \left( M\right) }} \!\!\!\!\!\!\!\!\!\!\!\int_{\mathbb{R}^n}\frac{\left\Vert \mathsf{Q}_{F,M}^{\omega ,\mathbf{%
b}^{\ast }}x\right\Vert _{L^{2}\left( \omega \right) }^{\spadesuit 2}}{%
\left( \left\vert M\right\vert +\left\vert y-c_{M}\right\vert \right)
^{n+1-\alpha }}\frac{\left\Vert \mathsf{Q}_{F^{\prime },M^{\prime }}^{\omega ,
\mathbf{b}^{\ast }}x\right\Vert _{L^{2}\left( \omega \right) }^{\spadesuit 2}%
}{\left( \left\vert M^{\prime }\right\vert +\left\vert y-c_{M^{\prime
}}\right\vert \right) ^{n+1-\alpha }}d\sigma \left( y\right)  \notag
\end{eqnarray}%
by proving that%
\begin{equation}
U_{s}\lesssim 2^{-\delta s}\left\{ \mathcal{A}_{2}^{\alpha }+\left( 
\mathfrak{E}_{2}^{\alpha }+\sqrt{A_{2}^{\alpha ,{energy}}}\right) 
\sqrt{A_{2}^{\alpha ,{punct}}}\right\} \int_{\widehat{I}}t^{2}d%
\overline{\mu },\ \ \ \ \ \text{where }\delta =\frac{1}{2}.  \label{Us bound}
\end{equation}%
With this accomplished, we can sum in $s\geq 0$ to control the term $\mathbf{%
Back}\left( \widehat{I}\right) $. We now decompose $U_{s}=T_{s}^{{%
proximal}}+T_{s}^{{difference}}+T_{s}^{{intersection}}$ into
three pieces.

Our first decomposition is to write%
\begin{equation}
U_{s}=T_{s}^{{proximal}}+V_{s}^{{remote}}\ ,
\label{initial decomp}
\end{equation}%
where in the `proximal' term $T_{s}^{{proximal}}$ we restrict the
summation over pairs of cubes $M,M^{\prime }$ to those satisfying $%
d\left( c_{M},c_{M^{\prime }}\right) <2^{s\delta }\ell \left( M\right) $;
while in the `remote' term $V_{s}^{{remote}}$ we restrict the
summation over pairs of cubes $M,M^{\prime }$ to those satisfying the
opposite inequality $d\left( c_{M},c_{M^{\prime }}\right) \geq 2^{s\delta
}\ell \left( M\right) $. Then we further decompose 
\begin{equation*}
V_{s}^{{remote}}=T_{s}^{{difference}}+T_{s}^{{%
intersection}},
\end{equation*}%
where in the `difference' term $T_{s}^{{difference}}$ we restrict
integration in $y$ to the difference $\mathbb{R}\setminus B\left(
M,M^{\prime }\right) $ of $\mathbb{R}$ and 
\begin{equation}
B\left( M,M^{\prime }\right) \equiv B\left( c_{M},\frac{1}{2}d\left(
c_{M},c_{M^{\prime }}\right) \right) ,  \label{def BMM'}
\end{equation}%
the ball centered at $c_{M}$ with radius $\frac{1}{2}d\left(
c_{M},c_{M^{\prime }}\right) $; while in the `intersection' term $T_{s}^{%
{intersection}}$ we restrict integration in $y$ to the intersection $%
\mathbb{R}^n\cap B\left( M,M^{\prime }\right) $ of $\mathbb{R}^n$ with the ball $%
B\left( M,M^{\prime }\right) $; i.e. 
\begin{eqnarray}
&&\label{def Tints}\\
T_{s}^{{intersection}} &\equiv &\sum_{\substack{ F,F^{\prime }\in 
\mathcal{F}}}\sum_{\substack{ M\in \mathcal{W}\left( F\right) ,\ M^{\prime
}\in \mathcal{W}\left( F^{\prime }\right)  \\ M,M^{\prime }\subset I,\ \ell
\left( M^{\prime }\right) =2^{-s}\ell \left( M\right)  \\ d\left(
c_{M},c_{M^{\prime }}\right) \geq 2^{s\left( 1+\delta \right) }\ell \left(
M^{\prime }\right) }} \notag \\
&&\int_{B\left( M,M^{\prime }\right) }\frac{\left\Vert \mathsf{Q}%
_{F,M}^{\omega ,\mathbf{b}^{\ast }}x\right\Vert _{L^{2}\left( \omega \right)
}^{\spadesuit 2}}{\left( \left\vert M\right\vert +\left\vert
y-c_{M}\right\vert \right) ^{n+1-\alpha }}\frac{\left\Vert \mathsf{Q}%
_{F^{\prime },M^{\prime }}^{\omega ,\mathbf{b}^{\ast }}x\right\Vert
_{L^{2}\left( \omega \right) }^{\spadesuit 2}}{\left( \left\vert M^{\prime
}\right\vert +\left\vert y-c_{M^{\prime }}\right\vert \right) ^{n+1-\alpha }}%
d\sigma \left( y\right) .  \notag
\end{eqnarray}%
Here is a schematic reminder of the these decompositions with the
distinguishing points of the definitions boxed:{}

\begin{equation*}
\fbox{$%
\begin{array}{ccccc}
U_{s} &  &  &  &  \\ 
\downarrow &  &  &  &  \\ 
T_{s}^{{proximal}} & + & V_{s}^{{remote}} &  &  \\ 
\fbox{$d\left( c_{M},c_{M^{\prime }}\right) <2^{s\delta }\ell \left(
M\right) $} &  & \fbox{$d\left( c_{M},c_{M^{\prime }}\right) \geq 2^{s\delta
}\ell \left( M\right) $} &  &  \\ 
&  & \downarrow &  &  \\ 
&  & T_{s}^{{difference}} & + & T_{s}^{{intersection}} \\ 
&  & \fbox{$\int_{\mathbb{R}\setminus B\left( M,M^{\prime }\right) }$} &  & 
\fbox{$\fbox{$\int_{B\left( M,M^{\prime }\right) }$}$}%
\end{array}%
$}
\end{equation*}

We will exploit the restriction of integration to $B\left( M,M^{\prime
}\right) $, together with the condition 
\begin{equation*}
d\left( c_{M},c_{M^{\prime }}\right) \geq 2^{s\left( 1+\delta \right) }\ell
\left( M^{\prime }\right) =2^{s\delta }\ell \left( M\right) ,
\end{equation*}%
which will then give an estimate for the term $T_{s}^{{intersection}%
} $ using an argument dual to that used for the other terms $T_{s}^{{%
proximal}}$ and $T_{s}^{{difference}}$, to which we now turn.

\subsubsection{The proximal and difference terms}

We have%
\begin{align}
T_{s}^{{proximal}}& \equiv \sum_{\substack{ F,F^{\prime }\in 
\mathcal{F}}}\sum_{\substack{ M\in \mathcal{W}\left( F\right) ,\ M^{\prime
}\in \mathcal{W}\left( F^{\prime }\right)  \\ M,M^{\prime }\subset I,\ \ell
\left( M^{\prime }\right) =2^{-s}\ell \left( M\right) \text{ and }d\left(
c_{M},c_{M^{\prime }}\right) <2^{s\delta }\ell \left( M\right) }}
\label{def Tproxs} \\
& \times \int_{\mathbb{R}^n}\frac{\left\Vert \mathsf{Q}_{F,M}^{\omega ,\mathbf{%
b}^{\ast }}x\right\Vert _{L^{2}\left( \omega \right) }^{\spadesuit 2}}{%
\left( \left\vert M\right\vert +\left\vert y-c_{M}\right\vert \right)
^{n+1-\alpha }}\frac{\left\Vert \mathsf{Q}_{F^{\prime },M^{\prime }}^{\omega ,%
\mathbf{b}^{\ast }}x\right\Vert _{L^{2}\left( \omega \right) }^{\spadesuit 2}%
}{\left( \left\vert M^{\prime }\right\vert +\left\vert y-c_{M^{\prime
}}\right\vert \right) ^{n+1-\alpha }}d\sigma \left( y\right)  \notag \\
& \leq \mathcal{M}_{s}^{{proximal}}\sum_{F\in \mathcal{F}}\sum 
_{\substack{ M\in \mathcal{W}\left( F\right)  \\ M\subset I}}\lVert \mathsf{Q%
}_{F,M}^{\omega ,\mathbf{b}^{\ast }}z\rVert _{\omega }^{\spadesuit 2}=%
\mathcal{M}_{s}^{{proximal}}\int_{\widehat{I}}t^{2}d\overline{\mu },
\notag
\end{align}%
where%
\begin{align*}
\mathcal{M}_{s}^{{proximal}}& \equiv \sup_{F\in \mathcal{F}}\sup 
_{\substack{ M\in \mathcal{W}\left( F\right)  \\ M\subset I}}\mathcal{A}%
_{s}^{{proximal}}\left( M\right) ; \\
\mathcal{A}_{s}^{{proximal}}\left( M\right) & \equiv \sum_{F^{\prime
}\in \mathcal{F}}\sum_{\substack{ M^{\prime }\in \mathcal{W}\left( F^{\prime
}\right)  \\ M^{\prime }\subset I,\ \ell \left( M^{\prime }\right)
=2^{-s}\ell \left( M\right) \text{ and }\\ d\left( c_{M},c_{M^{\prime }}\right)
<2^{s\delta }\ell \left( M\right) }}\int_{\mathbb{R}^n}S_{\left( M^{\prime
},M\right) }^{F^{\prime }}\left( y\right) d\sigma \left( y\right) ; \\
S_{\left( M^{\prime },M\right) }^{F^{\prime }}\left( x\right)
& \equiv \frac{%
1}{\left( \left\vert M\right\vert +\left\vert y-c_{M}\right\vert \right)
^{n+1-\alpha }}\frac{\left\Vert \mathsf{Q}_{F^{\prime },M^{\prime }}^{\omega
}x\right\Vert _{L^{2}\left( \omega \right) }^{\spadesuit 2}}{\left(
\left\vert M^{\prime }\right\vert +\left\vert y-c_{M^{\prime }}\right\vert
\right) ^{n+1-\alpha }},
\end{align*}%
and similarly%
\begin{align}
T_{s}^{{difference}}& \equiv \sum_{\substack{ F,F^{\prime }\in 
\mathcal{F}}}\sum_{\substack{ M\in \mathcal{W}\left( F\right) ,\ M^{\prime
}\in \mathcal{W}\left( F^{\prime }\right)  \\ M,M^{\prime }\subset I,\ \ell
\left( M^{\prime }\right) =2^{-s}\ell \left( M\right) \text{ and }d\left(
c_{M},c_{M^{\prime }}\right) \geq 2^{s\delta }\ell \left( M\right) }}
\label{def Tdiffs} \\
& \times \int_{\mathbb{R}^n\setminus B\left( M,M^{\prime }\right) }\frac{%
\left\Vert \mathsf{Q}_{F,M}^{\omega ,\mathbf{b}^{\ast }}x\right\Vert
_{L^{2}\left( \omega \right) }^{\spadesuit 2}}{\left( \left\vert
M\right\vert +\left\vert y-c_{M}\right\vert \right) ^{n+1-\alpha }}\frac{%
\left\Vert \mathsf{Q}_{F^{\prime },M^{\prime }}^{\omega ,\mathbf{b}^{\ast
}}x\right\Vert _{L^{2}\left( \omega \right) }^{\spadesuit 2}}{\left(
\left\vert M^{\prime }\right\vert +\left\vert y-c_{M^{\prime }}\right\vert
\right) ^{n+1-\alpha }}d\sigma \left( y\right)  \notag \\
& \leq \mathcal{M}_{s}^{{difference}}\sum_{F\in \mathcal{F}}\sum 
_{\substack{ M\in \mathcal{W}\left( F\right)  \\ M\subset I}}\lVert \mathsf{Q%
}_{F,M}^{\omega ,\mathbf{b}^{\ast }}z\rVert _{\omega }^{\spadesuit 2}=%
\mathcal{M}_{s}^{{difference}}\int_{\widehat{I}}t^{2}d\overline{\mu }%
;  \notag
\end{align}%
where%
\begin{eqnarray*}
\mathcal{M}_{s}^{{difference}} &\equiv &\sup_{F\in \mathcal{F}}\sup 
_{\substack{ _{\substack{ M\in \mathcal{W}\left( F\right) }}  \\ M\subset I}}%
\mathcal{A}_{s}^{{difference}}\left( M\right) ; \\
\mathcal{A}_{s}^{{difference}}\left( M\right) &\equiv
&\sum_{F^{\prime }\in \mathcal{F}}\sum_{\substack{ M^{\prime }\in \mathcal{W}%
\left( F^{\prime }\right)  \\ M^{\prime }\subset I,\ \ell \left( M^{\prime
}\right) =2^{-s}\ell \left( M\right) \text{ and }\\ d\left( c_{M},c_{M^{\prime
}}\right) \geq 2^{s\delta }\ell \left( M\right) }}\!\!\!\!\! \int_{\mathbb{R}^n\setminus
B\left( M,M^{\prime }\right) }S_{\left( M^{\prime },M\right) }^{F^{\prime
}}\left( y\right) d\sigma \left( y\right) .
\end{eqnarray*}%
The restriction of integration in $\mathcal{A}_{s}^{{difference}}$
to $\mathbb{R}^n\setminus B\left( M,M^{\prime }\right) $ will be used to
establish (\ref{vanishing close}) below.

\begin{notation}
\label{Sum *}Since the cubes $F,M,F^{\prime },M^{\prime }$ that arise in
all of the sums here satisfy 
\begin{equation*}
M\in \mathcal{W}\left( F\right) ,\ M^{\prime }\in \mathcal{W}\left(
F^{\prime }\right) \text{ and }\ell \left( M^{\prime }\right) =2^{-s}\ell
\left( M\right) \text{ and }M,M^{\prime }\subset I,
\end{equation*}%
we will often employ the notation $\overset{\ast }{\sum }$ to remind the
reader that, as applicable, these four conditions are in force even when
they are\ not explictly mentioned.
\end{notation}

Now fix $M$ as in $\mathcal{M}_{s}^{{proximal}}$ respectively $%
\mathcal{M}_{s}^{{difference}}$, and decompose the sum over $%
M^{\prime }$ in $\mathcal{A}_{s}^{{proximal}}\left( M\right) $
respectively $\mathcal{A}_{s}^{{difference}}\left( M\right) $ by%
\begin{eqnarray*}
&&\mathcal{A}_{s}^{{proximal}}\left( M\right) =\sum_{F^{\prime }\in 
\mathcal{F}}\sum_{\substack{ M^{\prime }\in \mathcal{W}\left( F^{\prime
}\right)  \\ M^{\prime }\subset I,\ \ell \left( M^{\prime }\right)
=2^{-s}\ell \left( M\right) \text{ and }\\d\left( c_{M},c_{M^{\prime }}\right)
<2^{s\delta }\ell \left( M\right) }}  \int_{\mathbb{R}^n}S_{\left( M^{\prime
},M\right) }^{F^{\prime }}\left( y\right) d\sigma \left( y\right) \\
&=&\sum_{F^{\prime }\in \mathcal{F}}\overset{\ast }{\sum_{\substack{ %
c_{M^{\prime }}\in 2M  \\ d\left( c_{M},c_{M^{\prime }}\right) <2^{s\delta
}\ell \left( M\right) }}}\int_{\mathbb{R}^n}S_{\left( M^{\prime },M\right)
}^{F^{\prime }}\left( y\right) d\sigma \left( y\right) \\
&&
\hspace{3cm}+\sum_{F^{\prime }\in 
\mathcal{F}}\sum_{\ell =1}^{\infty }\overset{\ast }{\sum_{\substack{ %
c_{M^{\prime }}\in 2^{\ell +1}M\setminus 2^{\ell }M  \\ d\left(
c_{M},c_{M^{\prime }}\right) <2^{s\delta }\ell \left( M\right) }}}\int_{%
\mathbb{R}^n}S_{\left( M^{\prime },M\right) }^{F^{\prime }}\left( y\right)
d\sigma \left( y\right) \\
&\equiv &
\sum_{\ell =0}^{\infty }\mathcal{A}_{s}^{{proximal},\ell
}\left( M\right) ,
\end{eqnarray*}%
respectively%
\begin{eqnarray*}
&&\mathcal{A}_{s}^{{difference}}\left( M\right) =\sum_{F^{\prime
}\in \mathcal{F}}\sum_{\substack{ M^{\prime }\in \mathcal{W}\left( F^{\prime
}\right)  \\ M^{\prime }\subset I,\ \ell \left( M^{\prime }\right)
=2^{-s}\ell \left( M\right) \text{ and }\\ d\left( c_{M},c_{M^{\prime }}\right)
\geq 2^{s\delta }\ell \left( M\right) }} \!\!\!\!\!\!\!\! \int_{\mathbb{R}^n\setminus B\left(
M,M^{\prime }\right) }S_{\left( M^{\prime },M\right) }^{F^{\prime }}\left(
y\right) d\sigma \left( y\right) \\
&=&\sum_{F^{\prime }\in \mathcal{F}}\overset{\ast }{\sum_{\substack{ %
c_{M^{\prime }}\in 2M  \\ d\left( c_{M},c_{M^{\prime }}\right) \geq
2^{s\delta }\ell \left( M\right) }}}\int_{\mathbb{R}^n\setminus B\left(
M,M^{\prime }\right) }S_{\left( M^{\prime },M\right) }^{F^{\prime }}\left(
y\right) d\sigma \left( y\right) \\
&&\hspace{2.5cm}+\sum_{\ell =1}^{\infty }\sum_{F^{\prime }\in \mathcal{F}}\overset{\ast }{%
\sum_{\substack{ c_{M^{\prime }}\in 2^{\ell +1}M\setminus 2^{\ell }M  \\ %
d\left( c_{M},c_{M^{\prime }}\right) \geq 2^{s\delta }\ell \left( M\right) }}%
}\int_{\mathbb{R}^n\setminus B\left( M,M^{\prime }\right) }S_{\left( M^{\prime
},M\right) }^{F^{\prime }}\left( y\right) d\sigma \left( y\right) \\
&\equiv &\sum_{\ell =0}^{\infty }\mathcal{A}_{s}^{{difference},\ell
}\left( M\right) .
\end{eqnarray*}%
Let $m=2$ so that 
\begin{equation}
2^{-m}\leq \frac{1}{3}.  \label{smallest m}
\end{equation}%
Now decompose the integrals over $\mathbb{R}^n$ in $\mathcal{A}_{s}^{{%
proximal},\ell }\left( M\right) $ by%
\begin{eqnarray*}
\mathcal{A}_{s}^{{proximal},0}\left( M\right) &=&\sum_{F^{\prime
}\in \mathcal{F}}\overset{\ast }{\sum_{\substack{ c_{M^{\prime }}\in 2M  \\ %
d\left( c_{M},c_{M^{\prime }}\right) <2^{s\delta }\ell \left( M\right) }}}%
\int_{\mathbb{R}^n\setminus 4M}S_{\left( M^{\prime },M\right) }^{F^{\prime
}}\left( y\right) d\sigma \left( y\right) \\
&&+\sum_{F^{\prime }\in \mathcal{F}}\overset{\ast }{\sum_{\substack{ %
c_{M^{\prime }}\in 2M  \\ d\left( c_{M},c_{M^{\prime }}\right) <2^{s\delta
}\ell \left( M\right) }}}\int_{4M}S_{\left( M^{\prime },M\right)
}^{F^{\prime }}\left( y\right) d\sigma \left( y\right) \\
&\equiv &\mathcal{A}_{s,far}^{{proximal},0}\left( M\right) +\mathcal{%
A}_{s,near}^{{proximal},0}\left( M\right) ,
\end{eqnarray*}%
and for $\ell\geq 1$
\begin{eqnarray*}
\mathcal{A}_{s}^{{proximal},\ell }\left( M\right) &=&\sum_{F^{\prime
}\in \mathcal{F}}\overset{\ast }{\sum_{\substack{ c_{M^{\prime }}\in 2^{\ell
+1}M\setminus 2^{\ell }M  \\ d\left( c_{M},c_{M^{\prime }}\right)
<2^{s\delta }\ell \left( M\right) }}}\int_{\mathbb{R}^n\setminus 2^{\ell
+2}M}S_{\left( M^{\prime },M\right) }^{F^{\prime }}\left( y\right) d\sigma
\left( y\right) \\
&&+\sum_{F^{\prime }\in \mathcal{F}}\overset{\ast }{\sum_{\substack{ %
c_{M^{\prime }}\in 2^{\ell +1}M\setminus 2^{\ell }M  \\ d\left(
c_{M},c_{M^{\prime }}\right) <2^{s\delta }\ell \left( M\right) }}}%
\int_{2^{\ell +2}M\setminus 2^{\ell -m}M}S_{\left( M^{\prime },M\right)
}^{F^{\prime }}\left( y\right) d\sigma \left( y\right) \\
&&+\sum_{F^{\prime }\in \mathcal{F}}\overset{\ast }{\sum_{\substack{ %
c_{M^{\prime }}\in 2^{\ell +1}M\setminus 2^{\ell }M  \\ d\left(
c_{M},c_{M^{\prime }}\right) <2^{s\delta }\ell \left( M\right) }}}%
\int_{2^{\ell -m}M}S_{\left( M^{\prime },M\right) }^{F^{\prime }}\left(
y\right) d\sigma \left( y\right) \\
&\equiv &\mathcal{A}_{s,far}^{{proximal},\ell }\left( M\right) +%
\mathcal{A}_{s,near}^{{proximal},\ell }\left( M\right) +\mathcal{A}%
_{s,close}^{{proximal},\ell }\left( M\right) 
\end{eqnarray*}%
Similarly we decompose the integrals over the difference 
\begin{equation*}
B^{\ast }\equiv \mathbb{R}^n\setminus B\left( M,M^{\prime }\right)
\end{equation*}%
in $\mathcal{A}_{s}^{{difference},\ell }\left( M\right) $ by%
\begin{eqnarray*}
\mathcal{A}_{s}^{{difference},0}\left( M\right) &=&\sum_{F^{\prime
}\in \mathcal{F}}\overset{\ast }{\sum_{\substack{ c_{M^{\prime }}\in 2M  \\ %
d\left( c_{M},c_{M^{\prime }}\right) \geq 2^{s\delta }\ell \left( M\right) }}%
}\int_{B^{\ast }\setminus 4M}S_{\left( M^{\prime },M\right) }^{F^{\prime
}}\left( y\right) d\sigma \left( y\right) \\
&&+\sum_{F^{\prime }\in \mathcal{F}}\overset{\ast }{\sum_{\substack{ %
c_{M^{\prime }}\in 2M  \\ d\left( c_{M},c_{M^{\prime }}\right) \geq
2^{s\delta }\ell \left( M\right) }}}\int_{B^{\ast }\cap 4M}S_{\left(
M^{\prime },M\right) }^{F^{\prime }}\left( y\right) d\sigma \left( y\right)
\\
&\equiv &\mathcal{A}_{s,far}^{{difference},0}\left( M\right) +%
\mathcal{A}_{s,near}^{{difference},0}\left( M\right) ,
\end{eqnarray*}%
and%
\begin{eqnarray*}
&&\mathcal{A}_{s}^{{difference},\ell }\left( M\right)
=\sum_{F^{\prime }\in \mathcal{F}}\overset{\ast }{\sum_{\substack{ %
c_{M^{\prime }}\in 2^{\ell +1}M\setminus 2^{\ell }M  \\ d\left(
c_{M},c_{M^{\prime }}\right) \geq 2^{s\delta }\ell \left( M\right) }}}%
\int_{B^{\ast }\setminus 2^{\ell +2}M}S_{\left( M^{\prime },M\right)
}^{F^{\prime }}\left( y\right) d\sigma \left( y\right) \\
&&+\sum_{F^{\prime }\in \mathcal{F}}\overset{\ast }{\sum_{\substack{ %
c_{M^{\prime }}\in 2^{\ell +1}M\setminus 2^{\ell }M  \\ d\left(
c_{M},c_{M^{\prime }}\right) \geq 2^{s\delta }\ell \left( M\right) }}}%
\int_{B^{\ast }\cap \left( 2^{\ell +2}M\setminus 2^{\ell -m}M\right)
}S_{\left( M^{\prime },M\right) }^{F^{\prime }}\left( y\right) d\sigma
\left( y\right) \\
&&+\sum_{F^{\prime }\in \mathcal{F}}\overset{\ast }{\sum_{\substack{ %
c_{M^{\prime }}\in 2^{\ell +1}M\setminus 2^{\ell }M  \\ d\left(
c_{M},c_{M^{\prime }}\right) \geq 2^{s\delta }\ell \left( M\right) }}}%
\int_{B^{\ast }\cap 2^{\ell -m}M}S_{\left( M^{\prime },M\right) }^{F^{\prime
}}\left( y\right) d\sigma \left( y\right) \\
&\equiv &\mathcal{A}_{s,far}^{{difference},\ell }\left( M\right) +%
\mathcal{A}_{s,near}^{{difference},\ell }\left( M\right) +\mathcal{A}%
_{s,close}^{{difference},\ell }\left( M\right) ,\ \ \ \ \ \ell \geq
1.
\end{eqnarray*}

We now note the important point that the close terms $\mathcal{A}_{s,close}^{%
{proximal},\ell }\left( M\right) $ and $\mathcal{A}_{s,close}^{%
{difference},\ell }\left( M\right) $ both $\emph{vanish}$ for $\ell
>\delta s$ because of the decomposition (\ref{initial decomp}):%
\begin{equation}
\mathcal{A}_{s,close}^{{proximal},\ell }\left( M\right) =\mathcal{A}%
_{s,close}^{{difference},\ell }\left( M\right) =0,\ \ \ \ \ \ell
\geq 1+\delta s.  \label{vanishing close}
\end{equation}%
Indeed, if $c_{M^{\prime }}\in 2^{\ell +1}M\setminus 2^{\ell }M$, then we
have%
\begin{equation}
\frac{1}{2}2^{\ell }\ell \left( M\right) \leq d\left( c_{M},c_{M^{\prime
}}\right) .  \label{distJJ'}
\end{equation}%
Now the summands in $\mathcal{A}_{s,close}^{{proximal},\ell }\left(
M\right) $ satisfy $d\left( c_{M},c_{M^{\prime }}\right) <2^{\delta s}\ell
\left( M\right) $, which by (\ref{distJJ'}) is impossible if $\ell \geq
1+\delta s$ - indeed, if $\ell \geq 1+\delta s$, we get the contradiction%
\begin{equation*}
2^{\delta s}\ell \left( M\right) =\frac{1}{2}2^{1+\delta s}\ell \left(
M\right) \leq \frac{1}{2}2^{\ell }\ell \left( M\right) \leq d\left(
c_{M},c_{M^{\prime }}\right) <2^{\delta s}\ell \left( M\right) .
\end{equation*}%
It now follows that $\mathcal{A}_{s,close}^{{proximal},\ell }\left(
M\right) =0$. Thus we are left to consider the term $\mathcal{A}_{s,close}^{%
{difference},\ell }\left( M\right) $, where the integration is taken
over the set $\mathbb{R}^n\setminus B\left( M,M^{\prime }\right) $. But we are
also restricted in $\mathcal{A}_{s,close}^{{difference},\ell }\left(
M\right) $ to integrating over the cube $2^{\ell -m}M$, which is
contained in $B\left( M,M^{\prime }\right) $ by (\ref{distJJ'}). Indeed, the
smallest\ ball centered at $c_{M}$ that contains $2^{\ell -m}M$ has radius $%
\frac{1}{2}2^{\ell -m}\ell \left( M\right) $, which by (\ref{smallest m})
and (\ref{distJJ'}) is at most $\frac{1}{4}2^{\ell }\ell \left( M\right)
\leq \frac{1}{2}d\left( c_{M},c_{M^{\prime }}\right) $, the radius of $%
B\left( M,M^{\prime }\right) $. Thus the range of integration in the term $%
\mathcal{A}_{s,close}^{{difference},\ell }\left( M\right) $ is the
empty set, and so $\mathcal{A}_{s,close}^{{difference},\ell }\left(
M\right) =0$ as well as $\mathcal{A}_{s,close}^{{proximal},\ell
}\left( M\right) =0$. This proves (\ref{vanishing close}).

From now on we treat $T_{s}^{{proximal}}$ and $T_{s}^{{%
difference}}$ in the same way since the terms $\mathcal{A}_{s,close}^{%
{proximal},\ell }\left( M\right) $ and $\mathcal{A}_{s,close}^{%
{difference},\ell }\left( M\right) $ both vanish for $\ell \geq
1+\delta s$. Thus we will suppress the superscripts ${proximal}$ and 
${difference}$ in the $far$, $near$ and $close$ decomposition of $%
\mathcal{A}_{s}^{{proximal},\ell }\left( M\right) $ and $\mathcal{A}%
_{s}^{{difference},\ell }\left( M\right) $, and we will also
suppress the conditions $d\left( c_{M},c_{M^{\prime }}\right) <2^{s\delta
}\ell \left( M\right) $ and $d\left( c_{M},c_{M^{\prime }}\right) \geq
2^{s\delta }\ell \left( M\right) $ in the proximal and difference terms
since they no longer play a role. Using the pairwise disjointedness of the
shifted coronas $\mathcal{C}_{F}^{\mathcal{G},{shift}}$, we have 
\begin{equation*}
\sum_{F^{\prime }\in \mathcal{F}}\left\Vert \mathsf{Q}_{F^{\prime
},A}^{\omega ,\mathbf{b}^{\ast }}x\right\Vert _{L^{2}\left( \omega \right)
}^{\spadesuit 2}\lesssim \left\vert A\right\vert ^{2}\left\vert A\right\vert
_{\omega }\ ,\ \ \ \ \ \text{for any cube }A.
\end{equation*}%
Note that if $c_{M^{\prime }}\in 2M$, then $M^{\prime }\subset 3M$. Then
with 
\begin{equation}
\mathcal{W}_{M}^{s}\equiv \bigcup\limits_{F^{\prime }\in \mathcal{F}%
}\left\{ M^{\prime }\in \mathcal{W}\left( F^{\prime }\right) :M^{\prime
}\subset 3M\text{ and }\ell \left( M^{\prime }\right) =2^{-s}\ell \left(
M\right) \right\} ,  \label{def WMs}
\end{equation}%
we have%
\begin{eqnarray*}
\mathcal{A}_{s,far}^{0}\left( M\right) &\leq &\sum_{F^{\prime }\in \mathcal{F%
}}\overset{\ast }{\sum_{c_{M^{\prime }}\in 2M}}\int_{\mathbb{R}^n\setminus
4M}S_{\left( M^{\prime },M\right) }^{F^{\prime }}\left( y\right) d\sigma
\left( y\right) \\
&\lesssim &\sum_{A\in \mathcal{W}_{M}^{s}}\sum_{F^{\prime }\in \mathcal{F}:\
A\in \mathcal{W}\left( F^{\prime }\right) }\int_{\mathbb{R}^n\setminus 4M}%
\frac{\left\Vert \mathsf{Q}_{F^{\prime },M^{\prime }}^{\omega ,\mathbf{b}%
^{\ast }}x\right\Vert _{L^{2}\left( \omega \right) }^{\spadesuit 2}}{\left(
\left\vert M\right\vert +\left\vert y-c_{M}\right\vert \right) ^{2\left(
n+1-\alpha \right) }}d\sigma \left( y\right) \\
&\lesssim &\sum_{A\in \mathcal{W}_{M}^{s}}\int_{\mathbb{R}^n\setminus 4M}\frac{%
\left\vert A\right\vert ^{2}\left\vert A\right\vert _{\omega }}{\left(
\left\vert M\right\vert +\left\vert y-c_{M}\right\vert \right) ^{2\left(
n+1-\alpha \right) }}d\sigma \left( y\right) \\
&=&\left( \sum_{A\in \mathcal{W}_{M}^{s}}\left\vert A\right\vert
^{2}\left\vert A\right\vert _{\omega }\right) \int_{\mathbb{R}^n\setminus 4M}%
\frac{1}{\left( \left\vert M\right\vert +\left\vert y-c_{M}\right\vert
\right) ^{2\left( n+1-\alpha \right) }}d\sigma \left( y\right) .
\end{eqnarray*}%
Now we use the standard pigeonholing of side length of $A$ to conclude that 
\begin{eqnarray}
\sum_{A\in \mathcal{W}_{M}^{s}}\left\vert A\right\vert ^{2}\left\vert
A\right\vert _{\omega } &=&\sum_{k=s}^{\infty }\sum_{A\in \mathcal{W}%
_{M}^{s}:\ \ell \left( A\right) =2^{-k}\ell \left( M\right) }\left\vert
A\right\vert ^{2}\left\vert A\right\vert _{\omega } \notag \\ &\leq& \sum_{k=s}^{\infty
}2^{-2k}\left\vert M\right\vert ^{2}\sum_{A\in \mathcal{W}_{M}^{s}:\ \ell
\left( A\right) =2^{-k}\ell \left( M\right) }\left\vert A\right\vert
_{\omega }  \label{stan pig} \\
&\leq &\sum_{k=s}^{\infty }2^{-2k}\left\vert M\right\vert ^{2}\left\vert
3M\right\vert _{\omega }\lesssim 2^{-2s}\left\vert M\right\vert
^{2}\left\vert 3M\right\vert _{\omega },  \notag
\end{eqnarray}%
so that combining the previous two displays we have%
\begin{eqnarray*}
\mathcal{A}_{s,far}^{0}\left( M\right) 
&\lesssim &
2^{-2s}\left\vert
M\right\vert ^{2}\left\vert 3M\right\vert _{\omega }\int_{\mathbb{R}^n
\setminus 4M}\frac{1}{\left( \left\vert M\right\vert +\left\vert
y-c_{M}\right\vert \right) ^{2\left( n+1-\alpha \right) }}d\sigma \left(
y\right) \\
&\leq &
2^{-2s}\left\vert 4M\right\vert _{\omega }\int_{\mathbb{R}^n \setminus
4M}\frac{1}{\left( \left\vert M\right\vert +\left\vert y-c_{M}\right\vert
\right) ^{2\left( 1-\alpha \right) }}d\sigma \left( y\right) \\
&\approx &
2^{-2s}\frac{\left\vert 4M\right\vert _{\omega }}{\left\vert
4M\right\vert ^{1-\alpha }}\int_{\mathbb{R}^n \setminus 4M}\left( \frac{%
\left\vert M\right\vert }{\left( \left\vert M\right\vert +\left\vert
y-c_{M}\right\vert \right) ^{2}}\right) ^{1-\alpha }d\sigma \left( y\right)
\\
&\lesssim &
2^{-2s}\frac{\left\vert 4M\right\vert _{\omega }}{\left\vert
4M\right\vert ^{1-\alpha }}\mathcal{P}^{\alpha }\left( 4M,\mathbf{1}_{%
\mathbb{R}^n\setminus 4M}\sigma \right) \lesssim 2^{-2s}\mathcal{A}%
_{2}^{\alpha }\ .
\end{eqnarray*}
To estimate the near term $\mathcal{A}_{s,near}^{0}\left( M\right) $, we
initially keep the energy $\left\Vert \mathsf{Q}_{F^{\prime },M^{\prime
}}^{\omega ,\mathbf{b}^{\ast }}z\right\Vert _{L^{2}\left( \omega \right)
}^{2}$ and write 
\begin{eqnarray}
&& \mathcal{A}_{s,near}^{0}\left( M\right) \leq \sum_{F^{\prime }\in \mathcal{%
F}}\overset{\ast }{\sum_{c_{M^{\prime }}\in 2M}}\int_{4M}S_{\left( M^{\prime
},M\right) }^{F^{\prime }}\left( y\right) d\sigma \left( y\right)
\label{A0snear} \notag \\
&\approx&
\sum_{F^{\prime }\in \mathcal{F}}\overset{\ast }{%
\sum_{c_{M^{\prime }}\in 2M}}\int_{4M}\frac{1}{\left\vert M\right\vert
^{n+1-\alpha }}\frac{\left\Vert \mathsf{Q}_{F^{\prime },M^{\prime }}^{\omega ,%
\mathbf{b}^{\ast }}x\right\Vert _{L^{2}\left( \omega \right) }^{\spadesuit 2}%
}{\left( \left\vert M^{\prime }\right\vert^\frac{1}{n} +\left\vert y-c_{M^{\prime
}}\right\vert \right) ^{n+1-\alpha }}d\sigma \left( y\right)  \notag \\
&=&
\sum_{F^{\prime }\in \mathcal{F}}\frac{1}{\left\vert M\right\vert
^{n+1-\alpha }}\!\!\overset{\ast }{\sum_{c_{M^{\prime }}\in 2M}}\left\Vert \mathsf{%
Q}_{F^{\prime },M^{\prime }}^{\omega ,\mathbf{b}^{\ast }}x\right\Vert
_{L^{2}\left( \omega \right) }^{\spadesuit 2}\int_{4M}\frac{d\sigma \left( y\right) }{\left(
\left\vert M^{\prime }\right\vert^\frac{1}{n} +\left\vert y-c_{M^{\prime }}\right\vert
\right) ^{n+1-\alpha }} \notag \\
&=&
\sum_{F^{\prime }\in \mathcal{F}}\frac{1}{\left\vert M\right\vert
^{n+1-\alpha }}\overset{\ast }{\sum_{c_{M^{\prime }}\in 2M}}\left\Vert \mathsf{%
Q}_{F^{\prime },M^{\prime }}^{\omega ,\mathbf{b}^{\ast }}x\right\Vert
_{L^{2}\left( \omega \right) }^{\spadesuit 2}\frac{\mathrm{P}^{\alpha
}\left( M^{\prime },\mathbf{1}_{4M}\sigma \right) }{\left\vert M^{\prime
}\right\vert^\frac{1}{n} }.  \notag
\end{eqnarray}%
In order to estimate the final sum above, we must invoke the `prepare to
puncture' argument above, as we will want to derive geometric decay from $%
\left\Vert \mathsf{Q}_{M^{\prime }}^{\omega ,\mathbf{b}^{\ast }}x\right\Vert
_{L^{2}\left( \omega \right) }^{\spadesuit 2}$ by dominating it by the
`nonenergy' term $\left\vert M^{\prime }\right\vert ^{2}\left\vert M^{\prime
}\right\vert _{\omega }$, as well as using the Muckenhoupt energy constant.
Choose an augmented cube $\widetilde{M}\!\in\! \mathcal{AD}$ satisfying $%
\bigcup\limits_{c_{M^{\prime }}\in 2M}\!\! M^{\prime }\subset 4M\subset 
\widetilde{M}$ and $\ell \left( \widetilde{M}\right) \leq C\ell \left(
M\right) $. Define $\widetilde{\omega }=\omega -\omega \left( \left\{
p\right\} \right) \delta _{p}$ where $p$ is an atomic point in $\widetilde{M}
$ for which 
\begin{equation*}
\omega \left( \left\{ p\right\} \right) =\sup_{q\in \mathfrak{P}_{\left(
\sigma ,\omega \right) }:\ q\in \widetilde{M}}\omega \left( \left\{
q\right\} \right) .
\end{equation*}%
(If $\omega $ has no atomic point in common with $\sigma $ in $\widetilde{M}$%
, set $\widetilde{\omega }=\omega $). Then we have $\left\vert \widetilde{M}%
\right\vert _{\widetilde{\omega }}=\omega \left( \widetilde{M},\mathfrak{P}%
_{\left( \sigma ,\omega \right) }\right) $ and%

\begin{equation*}
\frac{\left\vert \widetilde{M}\right\vert _{\widetilde{\omega }}}{\left\vert 
\widetilde{M}\right\vert ^{1-\frac{\alpha}{n} }}\frac{\left\vert \widetilde{M}%
\right\vert _{\sigma }}{\left\vert \widetilde{M}\right\vert ^{1-\frac{\alpha}{n} }}=%
\frac{\omega \left( \widetilde{M},\mathfrak{P}_{\left( \sigma ,\omega
\right) }\right) }{\left\vert \widetilde{M}\right\vert ^{1-\frac{\alpha}{n} }}\frac{%
\left\vert \widetilde{M}\right\vert _{\sigma }}{\left\vert \widetilde{M}%
\right\vert ^{1-\frac{\alpha}{n} }}\leq A_{2}^{\alpha ,{punct}}.
\end{equation*}%
From (\ref{key obs}) and (\ref{key fact}) we also have%
\begin{equation*}
\sum_{F^{\prime }\in \mathcal{F}}\left\Vert \mathsf{Q}_{F^{\prime
},A}^{\omega ,\mathbf{b}^{\ast }}x\right\Vert _{L^{2}\left( \omega \right)
}^{\spadesuit 2}\lesssim \ell \left( A\right) ^{2}\left\vert A\right\vert _{%
\widetilde{\omega }}\ ,\ \ \ \ \ \text{for any cube }A.
\end{equation*}
Now by Cauchy-Schwarz and the augmented local estimate (\ref{shifted local})
in Lemma \ref{shifted} with $M=\widetilde{M}$ applied to the second line
below, and with $\mathcal{W}_{M}^{s}$ as in (\ref{def WMs}), and noting (\ref%
{stan pig}), the last sum in (\ref{A0snear}) is dominated by%
\begin{eqnarray}
&&\frac{1}{\left\vert M\right\vert ^{n+1-\alpha }}\left( \sum_{F^{\prime }\in 
\mathcal{F}}\overset{\ast }{\sum_{c\left( M^{\prime }\right) \in 2M}}%
\left\Vert \mathsf{Q}_{F^{\prime },M^{\prime }}^{\omega ,\mathbf{b}^{\ast
}}x\right\Vert _{L^{2}\left( \omega \right) }^{\spadesuit 2}\right) ^{\frac{1%
}{2}}  \label{dom by} \\
&&\ \ \ \ \ \ \ \ \ \ \times \left( \sum_{F^{\prime }\in \mathcal{F}}\overset%
{\ast }{\sum_{c_{M^{\prime }}\in 2M}}\left\Vert \mathsf{Q}_{F^{\prime
},M^{\prime }}^{\omega ,\mathbf{b}^{\ast }}x\right\Vert _{L^{2}\left( \omega
\right) }^{\spadesuit 2}\left( \frac{\mathrm{P}^{\alpha }\left( M^{\prime },%
\mathbf{1}_{4M}\sigma \right) }{\left\vert M^{\prime }\right\vert^\frac{1}{n} }\right)
^{2}\right) ^{\frac{1}{2}}  \notag \\
&\lesssim &
\frac{1}{\left\vert M\right\vert ^{n+1-\alpha }}\left( \sum_{A\in 
\mathcal{W}_{M}^{s}}\left\vert A\right\vert ^{2}\left\vert A\right\vert _{%
\widetilde{\omega }}\right) ^{\frac{1}{2}}\sqrt{\left( \mathfrak{E}%
_{2}^{\alpha }\right) ^{2}+A_{2}^{\alpha ,{energy}}}\sqrt{\left\vert 
\widetilde{M}\right\vert _{\sigma }}  \notag \\
&\lesssim &
\frac{2^{-s}\left\vert M\right\vert }{\left\vert M\right\vert
^{n+1-\alpha }}\sqrt{\left\vert 4M\right\vert _{\widetilde{\omega }}}\sqrt{%
\left( \mathfrak{E}_{2}^{\alpha }\right) ^{2}+A_{2}^{\alpha ,{energy}%
}}\sqrt{\left\vert \widetilde{M}\right\vert _{\sigma }}  \notag \\
&\lesssim &
2^{-s}\sqrt{\left( \mathfrak{E}_{\alpha }\right)
^{2}+A_{2}^{\alpha ,{energy}}}\sqrt{\frac{\left\vert \widetilde{M}%
\right\vert _{\widetilde{\omega }}}{\left\vert \widetilde{M}\right\vert
^{n+1-\alpha }}\frac{\left\vert \widetilde{M}\right\vert _{\sigma }}{%
\left\vert \widetilde{M}\right\vert ^{n+1-\alpha }}}\notag \\
&\lesssim&
2^{-s}\sqrt{%
\left( \mathfrak{E}_{2}^{\alpha }\right) ^{2}+A_{2}^{\alpha ,{energy}%
}}\sqrt{A_{2}^{\alpha ,{punct}}}\ .  \notag
\end{eqnarray}

Similarly, for $\ell \geq 1$, we can estimate the far term $\mathcal{A}%
_{s,far}^{\ell }\left( M\right) $ by the argument used for $\mathcal{A}%
_{s,far}^{0}\left( M\right) $ but applied to $2^{\ell }M$ in place of $M$.
For this need the following variant of $\mathcal{W}_{M}^{s}$ in (\ref{def
WMs}) given by%
\begin{equation}
\mathcal{W}_{M}^{s,\ell }\equiv \bigcup\limits_{F^{\prime }\in \mathcal{F}%
}\left\{ M^{\prime }\in \mathcal{W}\left( F^{\prime }\right) :M^{\prime
}\subset 3\left( 2^{\ell }M\right) \text{ and }\ell \left( M^{\prime
}\right) =2^{-s-\ell }\ell \left( 2^{\ell }M\right) \right\} .
\label{def WMsell}
\end{equation}%
Then we have

\begin{eqnarray*}
\mathcal{A}_{s,far}^{\ell }\left( M\right) 
&\leq &
\sum_{F^{\prime }\in 
\mathcal{F}}\overset{\ast }{\sum_{c_{M^{\prime }}\in \left( 2^{\ell
+1}M\right) \setminus \left( 2^{\ell }M\right) }}\int_{\mathbb{R}^n\setminus
2^{\ell +2}M}S_{\left( M^{\prime },M\right) }^{F^{\prime }}\left( y\right)
d\sigma \left( y\right) \\
&\lesssim &
\sum_{A\in \mathcal{W}_{M}^{s,\ell }}\sum_{F^{\prime }\in 
\mathcal{F}:\ A\in \mathcal{W}\left( F^{\prime }\right) }\int_{\mathbb{R}^n%
\setminus 4\left( 2^{\ell }M\right) }\frac{\left\Vert \mathsf{Q}_{F^{\prime
},A}^{\omega ,\mathbf{b}^{\ast }}x\right\Vert _{L^{2}\left( \omega \right)
}^{\spadesuit 2}}{\left( \left\vert M\right\vert +\left\vert
y-c_{M}\right\vert \right) ^{2\left( n+1-\alpha \right) }}d\sigma \left(
y\right) \\
&\lesssim &
\sum_{A\in \mathcal{W}_{M}^{s,\ell }}\int_{\mathbb{R}^n\setminus
4\left( 2^{\ell }M\right) }\frac{\left\vert A\right\vert ^{2}\left\vert
A\right\vert _{\omega }}{\left( \left\vert M\right\vert +\left\vert
y-c_{M}\right\vert \right) ^{2\left( n+1-\alpha \right) }}d\sigma \left(
y\right) \\
&=&
\left( \sum_{A\in \mathcal{W}_{M}^{s,\ell }}\left\vert A\right\vert
^{2}\left\vert A\right\vert _{\omega }\right) \int_{\mathbb{R}^n\setminus
4\left( 2^{\ell }M\right) }\frac{1}{\left( \left\vert M\right\vert
+\left\vert y-c_{M}\right\vert \right) ^{2\left( n+1-\alpha \right) }}d\sigma
\left( y\right) ,
\end{eqnarray*}%
where, just as for the sum over $A\in \mathcal{W}_{M}^{s,0}$, we have%
\begin{eqnarray}
\sum_{A\in \mathcal{W}_{M}^{s,\ell }}\left\vert A\right\vert
^{2}\left\vert A\right\vert _{\omega }\notag  
&=&
\sum_{k=s}^{\infty }\sum_{A\in \mathcal{W}_{M}^{s,\ell }:\ \ell \left(
A\right) =2^{-k-\ell }\ell \left( 2^{\ell }M\right) }\left\vert A\right\vert
^{2}\left\vert A\right\vert _{\omega } \notag \\ 
&\leq& 
\sum_{k=s}^{\infty }2^{-2k-2\ell
}\left\vert 2^{\ell }M\right\vert ^{2}\sum_{A\in \mathcal{W}_{M}^{s,\ell }:\
\ell \left( A\right) =2^{-k-\ell }\ell \left( 2^{\ell }M\right) }\left\vert
A\right\vert _{\omega }  \label{stan pig'}  \\
&\leq &
\sum_{k=s}^{\infty }2^{-2k-2\ell }\left\vert 2^{\ell }M\right\vert
^{2}\left\vert 3\left( 2^{\ell }M\right) \right\vert _{\omega }\lesssim
2^{-2s-2\ell }\left\vert 2^{\ell }M\right\vert ^{2}\left\vert 3\left(
2^{\ell }M\right) \right\vert _{\omega }\ .  \notag
\end{eqnarray}%
Now using $\frac{\left\vert 2^{\ell }M\right\vert ^{2}}{\left( \left\vert
M\right\vert +\left\vert y-c_{M}\right\vert \right) ^{2\left( n+1-\alpha
\right) }}\leq \frac{1}{\left( \left\vert 2^{\ell }M\right\vert +\left\vert
y-c_{2^{\ell }M}\right\vert \right) ^{2\left( 1-\alpha \right) }}$ for $%
y\notin 2^{\ell +2}M$, we can continue with%
\begin{eqnarray*}
&&
\mathcal{A}_{s,far}^{\ell }\left( M\right) \lesssim 2^{-2s}2^{-2\ell
}\left\vert 2^{\ell +2}M\right\vert _{\omega }\int_{\mathbb{R}^n\setminus
2^{\ell +2}M}\frac{d\sigma
\left( y\right)}{\left( \left\vert 2^{\ell }M\right\vert +\left\vert
y-c_{2^{\ell }M}\right\vert \right) ^{2\left( 1-\alpha \right) }} \\
&\approx &
2^{-2s}2^{-2\ell }\frac{\left\vert 2^{\ell +2}M\right\vert
_{\omega }}{\left\vert 2^{\ell }M\right\vert ^{1-\alpha }}\int_{\mathbb{R}^n
\setminus 2^{\ell +2}M}\left( \frac{\left\vert 2^{\ell }M\right\vert }{\left( \left\vert 2^{\ell }M\right\vert +\left\vert y-c_{2^{\ell
}M}\right\vert \right) ^{2}}\right) ^{1-\alpha }d\sigma \left( y\right) \\
&\lesssim &
2^{-2s}2^{-2\ell }\left\{ \frac{\left\vert 2^{\ell
+2}M\right\vert _{\omega }}{\left\vert 2^{\ell }M\right\vert ^{1-\alpha }}%
\mathcal{P}^{\alpha }\left( 2^{\ell +2}M,1_{\mathbb{R}^n\setminus 2^{\ell
+2}M}\sigma \right) \right\} \lesssim 2^{-2s}2^{-2\ell }\mathcal{A}%
_{2}^{\alpha }\ .
\end{eqnarray*}

To estimate the near term $\mathcal{A}_{s,near}^{\ell }\left( M\right) $ we
must again invoke the \emph{`prepare to puncture'} argument. Choose an
augmented cube $\widetilde{M}\in \mathcal{AD}$ such that $\ell \left( 
\widetilde{M}\right) \leq C2^{\ell }\ell \left( M\right) $ and $
\bigcup\limits_{c_{M^{\prime }}\in 2^{\ell +1}M\setminus 2^{\ell
}M}M^{\prime }\subset 2^{\ell +2}M\subset \widetilde{M}$. Define $\widetilde{\omega }=\omega -\omega \left( \left\{ p\right\} \right) \delta
_{p}$ where $p$ is an atomic point in $\widetilde{M}$ for which 
$\displaystyle
\omega \left( \left\{ p\right\} \right) =\sup_{q\in \mathfrak{P}_{\left(
\sigma ,\omega \right) }:\ q\in \widetilde{M}}\omega \left( \left\{
q\right\} \right) .
$\\
(If $\omega $ has no atomic point in common with $\sigma $ in $\widetilde{M}$
set $\widetilde{\omega }=\omega $.) Then we have $\left\vert \widetilde{M}%
\right\vert _{\widetilde{\omega }}=\omega \left( \widetilde{M},\mathfrak{P}%
_{\left( \sigma ,\omega \right) }\right) $, and just as in the argument
above following (\ref{A0snear}), we have from (\ref{key obs}) and (\ref{key
fact}) that both%
\begin{equation*}
\frac{\left\vert \widetilde{M}\right\vert _{\widetilde{\omega }}}{\left\vert 
\widetilde{M}\right\vert ^{1-\frac{\alpha}{n} }}\frac{\left\vert \widetilde{M}%
\right\vert _{\sigma }}{\left\vert \widetilde{M}\right\vert ^{1-\frac{\alpha}{n} }}%
\leq A_{2}^{\alpha ,{punct}}\text{ and }\sum_{F^{\prime }\in 
\mathcal{F}}\left\Vert \mathsf{Q}_{F^{\prime },M^{\prime }}^{\omega ,\mathbf{%
b}^{\ast }}x\right\Vert _{L^{2}\left( \omega \right) }^{\spadesuit
2}\lesssim \ell \left( M^{\prime }\right) ^{2}\left\vert M^{\prime
}\right\vert _{\widetilde{\omega }}\ .
\end{equation*}%
Thus using that $m=2$ in the definition of $A_{s,near}^{\ell }\left(M\right) $, we see that
\begin{eqnarray*}
&&
\mathcal{A}_{s,near}^{\ell }\left( M\right) \leq \sum_{F^{\prime }\in 
\mathcal{F}}\overset{\ast }{\sum_{c_{M^{\prime }}\in 2^{\ell +1}M\setminus
2^{\ell }M}}\int_{2^{\ell +2}M\setminus 2^{\ell -m}M}S_{\left( M^{\prime
},M\right) }^{F^{\prime }}\left( y\right) d\sigma \left( y\right) \\
&\approx &\!\!\!\!\!\!
\sum_{F^{\prime }\in \mathcal{F}}\overset{\ast }{%
\sum_{c_{M^{\prime }}\in 2^{\ell +1}M\setminus 2^{\ell }M}}\int_{2^{\ell
+2}M\setminus 2^{\ell -m}M}\frac{1}{\left\vert 2^{\ell }M\right\vert
^{n+1-\alpha }}\frac{\left\Vert \mathsf{Q}_{F^{\prime },M^{\prime }}^{\omega ,%
\mathbf{b}^{\ast }}x\right\Vert _{L^{2}\left( \omega \right) }^{\spadesuit 2}%
}{\left(\!
\left\vert M^{\prime }\right\vert^\frac{1}{n}\! +\! \left\vert y-c_{M^{\prime }}\right\vert\!
\right)^{\!n\!+\!1\!-\!\alpha }}d\sigma \left( y\right) \\
&\lesssim &\!\!\!\!\!\!
\frac{1}{\left\vert 2^{\ell }M\right\vert ^{n+1-\alpha }}%
\!\!\!\sum_{F^{\prime }\in \mathcal{F}}\overset{\ast }{\sum_{c_{M^{\prime }}\in
2^{\ell +1}M\setminus 2^{\ell }M}}\!\!\!\!\left\Vert \mathsf{Q}_{F^{\prime
},M^{\prime }}^{\omega ,\mathbf{b}^{\ast }}x\right\Vert _{L^{2}\left( \omega
\right) }^{\spadesuit 2}
\int_{2^{\ell +2}M}\!\frac{d\sigma \left( y\right)}{\left(\!
\left\vert M^{\prime }\right\vert^\frac{1}{n}\! +\! \left\vert y-c_{M^{\prime }}\right\vert\!
\right)^{\!n\!+\!1\!-\!\alpha }}
\end{eqnarray*}%
is dominated by%
\begin{eqnarray*}
&&\frac{1}{\left\vert 2^{\ell }M\right\vert ^{n+1-\alpha }}\sum_{F^{\prime
}\in \mathcal{F}}\overset{\ast }{\sum_{c_{M^{\prime }}\in 2^{\ell
+1}M\setminus 2^{\ell }M}}\left\Vert \mathsf{Q}_{F^{\prime },M^{\prime
}}^{\omega ,\mathbf{b}^{\ast }}x\right\Vert _{L^{2}\left( \omega \right)
}^{\spadesuit 2}\frac{\mathrm{P}^{\alpha }\left( M^{\prime },\mathbf{1}%
_{2^{\ell +2}M}\sigma \right) }{\left\vert M^{\prime }\right\vert^\frac{1}{n} } \\
&\leq &
\frac{1}{\left\vert 2^{\ell }M\right\vert ^{n+1-\alpha }}\left(
\sum_{F^{\prime }\in \mathcal{F}}\overset{\ast }{\sum_{c_{M^{\prime }}\in
2^{\ell +1}M\setminus 2^{\ell }M}}\left\Vert \mathsf{Q}_{F^{\prime
},M^{\prime }}^{\omega ,\mathbf{b}^{\ast }}x\right\Vert _{L^{2}\left( \omega
\right) }^{\spadesuit 2}\right) ^{\frac{1}{2}} \\
&&\times \left( \sum_{F^{\prime }\in \mathcal{F}}\overset{\ast }{%
\sum_{c_{M^{\prime }}\in 2^{\ell +1}M\setminus 2^{\ell }M}}\left\Vert 
\mathsf{Q}_{F^{\prime },M^{\prime }}^{\omega ,\mathbf{b}^{\ast
}}x\right\Vert _{L^{2}\left( \omega \right) }^{\spadesuit 2}\left( \frac{%
\mathrm{P}^{\alpha }\left( M^{\prime },\mathbf{1}_{2^{\ell +2}M}\sigma
\right) }{\left\vert M^{\prime }\right\vert^\frac{1}{n} }\right) ^{2}\right) ^{\frac{1}{2%
}}.
\end{eqnarray*}

This can now be estimated as for the term $\mathcal{A}_{s,near}^{0}\left(
M\right) $, along with the augmented local estimate (\ref{shifted local}) in
Lemma \ref{shifted} with $M=\widetilde{M}$ applied to the final line above,
to get 
\begin{eqnarray*}
\mathcal{A}_{s,near}^{\ell }\left( M\right) &\lesssim &2^{-s}2^{-\ell }\frac{%
\left\vert 2^{\ell }M\right\vert }{\left\vert 2^{\ell }M\right\vert
^{n+1-\alpha }}\sqrt{\left\vert \widetilde{M}\right\vert _{\widetilde{\omega }}%
}\sqrt{\left( \mathfrak{E}_{2}^{\alpha }\right) ^{2}+A_{2}^{\alpha ,{%
energy}}}\sqrt{\left\vert \widetilde{M}\right\vert _{\sigma }} \\
&\lesssim &2^{-s}2^{-\ell }\sqrt{\left( \mathfrak{E}_{2}^{\alpha }\right)
^{2}+A_{2}^{\alpha ,{energy}}}\sqrt{\frac{\left\vert \widetilde{M}%
\right\vert _{\widetilde{\omega }}}{\left\vert \widetilde{M}\right\vert
^{1-\frac{\alpha}{n} }}\frac{\left\vert \widetilde{M}\right\vert _{\sigma }}{%
\left\vert \widetilde{M}\right\vert ^{1-\frac{\alpha}{n} }}} \\
&\lesssim &2^{-s}2^{-\ell }\sqrt{\left( \mathfrak{E}_{2}^{\alpha }\right)
^{2}+A_{2}^{\alpha ,{energy}}}\sqrt{A_{2}^{\alpha ,{punct}}}%
\ .
\end{eqnarray*}%
Each of the estimates for $\mathcal{A}_{s,far}^{\ell }\left( M\right) $ and $%
\mathcal{A}_{s,near}^{\ell }\left( M\right) $ is summable in both $s$ and $%
\ell $. 

Now we turn to the terms $\mathcal{A}_{s,close}^{\ell }\left( M\right) $,
and recall from (\ref{vanishing close}) that $\mathcal{A}_{s,close}^{\ell
}\left( M\right) =0$ if $\ell \geq 1+\delta s$. So we now suppose that $\ell
\leq \delta s$. We have, with $m=2$ as in (\ref{smallest m}),
\begin{eqnarray*}
&&\mathcal{A}_{s,close}^{\ell }\left( M\right) \leq \sum_{F^{\prime }\in 
\mathcal{F}}\overset{\ast }{\sum_{c_{M^{\prime }}\in 2^{\ell +1}M\setminus
2^{\ell }M}}\int_{2^{\ell -m}M}S_{\left( M^{\prime },M\right) }^{F^{\prime
}}\left( y\right) d\sigma \left( y\right) \\
&\approx &
\sum_{F^{\prime }\in \mathcal{F}}\overset{\ast }{%
\sum_{c_{M^{\prime }}\in 2^{\ell +1}M\setminus 2^{\ell }M}}\int_{2^{\ell
-m}M}\frac{1}{\left( \left\vert M\right\vert^\frac{1}{n} +\left\vert y-c_{M}\right\vert
\right) ^{n+1-\alpha }}\frac{\left\Vert \mathsf{Q}_{F^{\prime },M^{\prime
}}^{\omega ,\mathbf{b}^{\ast }}x\right\Vert _{L^{2}\left( \omega \right)
}^{\spadesuit 2}}{\left\vert 2^{\ell }M\right\vert ^{n+1-\alpha }}d\sigma
\left( y\right) \\
&=&
\left( \sum_{F^{\prime }\in \mathcal{F}}\overset{\ast }{%
\sum_{c_{M^{\prime }}\in 2^{\ell +1}M\setminus 2^{\ell }M}}\left\Vert 
\mathsf{Q}_{F^{\prime },M^{\prime }}^{\omega ,\mathbf{b}^{\ast
}}x\right\Vert _{L^{2}\left( \omega \right) }^{\spadesuit 2}\right) \frac{1}{%
\left\vert 2^{\ell }M\right\vert ^{n+1-\alpha }} \\
&&\ \ \ \ \ \ \ \ \ \ \ \ \ \ \ \ \ \ \ \ \times \int_{2^{\ell -m}M}\frac{1}{%
\left( \left\vert M\right\vert^\frac{1}{n} +\left\vert y-c_{M}\right\vert \right)
^{n+1-\alpha }}d\sigma \left( y\right) .
\end{eqnarray*}%
The argument used to prove (\ref{stan pig'}) gives the\ analogous inequality
with a hole $2^{\ell -1}M$, 
\begin{equation*}
\sum_{F^{\prime }\in \mathcal{F}}\overset{\ast }{\sum_{c_{M^{\prime }}\in
2^{\ell +1}M\setminus 2^{\ell }M}}\left\Vert \mathsf{Q}_{F^{\prime
},M^{\prime }}^{\omega ,\mathbf{b}^{\ast }}x\right\Vert _{L^{2}\left( \omega
\right) }^{\spadesuit 2}\lesssim 2^{-2s}\left\vert 2^{\ell }M\right\vert
^{\frac{2}{n}}\left\vert 2^{\ell +2}M\setminus 2^{\ell -1}M\right\vert _{\omega }\ .
\end{equation*}%
Thus we get that $\mathcal{A}_{s,close}^{\ell }\left( M\right)$ is bounded by
\begin{eqnarray*}
&\lesssim &
2^{-2s}\left\vert 2^{\ell }M\right\vert ^{\frac{2}{n}}\!\!\left\vert 2^{\ell
+2}M\!\setminus\! 2^{\ell -1}M\right\vert _{\omega }\frac{1}{\left\vert 2^{\ell
}M\right\vert ^{\frac{1}{n}({n+1-\alpha })}}\!\!\int_{2^{\ell -m}M}\frac{d\sigma
\left( y\right)}{\left( \!\left\vert
M\right\vert^\frac{1}{n} \!+\!\left\vert y-c_{M}\right\vert \!\right) ^{\!\!n\!+\!1\!-\!\alpha }} \\
&\lesssim &
2^{-2s}\left\vert 2^{\ell }M\right\vert ^{\frac{2}{n}}\frac{\left\vert
2^{\ell +2}M\setminus 2^{\ell -1}M\right\vert _{\omega }}{\left\vert 2^{\ell
}M\right\vert ^{\frac{1}{n}(n+1-\alpha) }}\frac{\left\vert 2^{\ell -m}M\right\vert _{\sigma
}}{\left\vert M\right\vert ^{\frac{1}{n}(n+1-\alpha) }} \\
&\lesssim &
2^{-2s}2^{\left( n+1-\alpha \right) \ell }\frac{\left\vert 2^{\ell
+2}M\setminus 2^{\ell -1}M\right\vert _{\omega }}{\left\vert 2^{\ell
+2}M\right\vert ^{1-\frac{\alpha}{n} }}\frac{\left\vert 2^{\ell -m}M\right\vert
_{\sigma }}{\left\vert 2^{\ell -m}M\right\vert ^{1-\frac{\alpha}{n} }}\lesssim
2^{-2s}2^{\left( n+1-\alpha \right) \ell }A_{2}^{\alpha },
\end{eqnarray*}%
provided that $m=2>1$. Note that we can use the offset Muckenhoupt constant $%
A_{2}^{\alpha }$ here since $2^{\ell +2}M\setminus 2^{\ell -1}M$ and $%
2^{\ell -m}M$ are disjoint. If $\ell \leq s$, then we have the relatively
crude estimate $\mathcal{A}_{s,close}^{\ell }\left( M\right) \lesssim
2^{-s}A_{2}^{\alpha }\ $without decay in $\ell $. But we are assuming $\ell
\leq \delta s$ here, and so we obtain a suitable estimate for $\mathcal{A}%
_{s,close}^{\ell }\left( M\right) $ provided we choose $0<\delta \leq \frac{1%
}{n+1-\alpha }$. Indeed, we then have%
\begin{equation*}
\sum_{l=1}^{\delta s}2^{-2s}2^{\left( n+1-\alpha \right) \ell }A_{2}^{\alpha
}=2^{-2s}\left( \sum_{l=1}^{\delta s}2^{\left( n+1-\alpha \right) \ell
}\right) A_{2}^{\alpha }\lesssim 2^{-2s}2^{\left( n+1-\alpha \right) \delta
s}A_{2}^{\alpha }\leq 2^{-s}A_{2}^{\alpha }\ ,
\end{equation*}%
provided $\delta \leq \frac{1}{n+1-\alpha }$, and in particular we may take $%
\delta =\frac{1}{2}$. Altogether, the above estimates prove%
\begin{equation*}
T_{s}^{{proximal}}+T_{s}^{{difference}}\lesssim 2^{-s}\left( 
\mathcal{A}_{2}^{\alpha }+\sqrt{\left( \mathfrak{E}_{2}^{\alpha }\right)
^{2}+A_{2}^{\alpha ,{energy}}}\sqrt{A_{2}^{\alpha ,{punct}}}%
\right)\int_{\widehat{I}}t^2d\mu ,
\end{equation*}%
which is summable in $s$.

\subsubsection{The intersection term}

Now we return to the term $T_{s}^{{intersection}}\equiv$
\begin{eqnarray*}
\sum_{\substack{ F,F^{\prime }\in 
\mathcal{F}}}\!\!\!\!\!\!\!\!\!\!\!\!\!\! \sum_{\substack{ M\in \mathcal{W}\left( F\right) \\ M^{\prime
}\in \mathcal{W}\left( F^{\prime }\right)  \\ M,M^{\prime }\subset I,\ \ell
\left( M^{\prime }\right) =2^{-s}\ell \left( M\right)  \\ d\left(
c_{M},c_{M^{\prime }}\right) \geq 2^{s\left( 1+\delta \right) }\ell \left(
M^{\prime }\right) }} 
\!\!\!\!\!\!\!\!\!\!\!\!\!\! \int_{B\left( M,M^{\prime }\right) }\frac{\left\Vert \mathsf{Q}%
_{F,M}^{\omega ,\mathbf{b}^{\ast }}x\right\Vert _{L^{2}\left( \omega \right)
}^{\spadesuit 2}}{\left( \left\vert M\right\vert +\left\vert
y-c_{M}\right\vert \right) ^{n\!+\!1\!-\!\alpha }}\frac{\left\Vert \mathsf{Q}%
_{F^{\prime },M^{\prime }}^{\omega ,\mathbf{b}^{\ast }}x\right\Vert
_{L^{2}\left( \omega \right) }^{\spadesuit 2}}{\left( \left\vert M^{\prime
}\right\vert +\left\vert y-c_{M^{\prime }}\right\vert \right) ^{n\!+\!1\!-\!\alpha }}%
d\sigma \left( y\right) .
\end{eqnarray*}%
It will suffice to show that $T_{s}^{{intersection}}$ satisfies the
estimate,%
\begin{eqnarray*}
T_{s}^{{intersection}} 
\!\!\!\!&\lesssim &\!\!\!\!
2^{-s\delta }\sqrt{\left( 
\mathfrak{E}_{2}^{\alpha }\right) ^{2}\!+\!A_{2}^{\alpha ,{energy}}}%
\sqrt{A_{2}^{\alpha ,{punct}}} 
\!\!\!\sum_{F^{\prime }\in \mathcal{F}%
^{\prime }}\!\sum_{\substack{ M^{\prime }\in \mathcal{M}_{\left( \mathbf{\rho }%
,\varepsilon \right) -{deep}}\left( F^{\prime }\right)  \\ M^{\prime
}\subset I}}\!\!\!\!\!\!\!\!\!\!\!\!\!\! \lVert \mathsf{Q}_{F^{\prime },M^{\prime }}^{\omega ,\mathbf{b}%
^{\ast }}x\rVert _{L^{2}\left( \omega \right) }^{\spadesuit 2} \\
&=&
\!\!\!\! 2^{-s\delta }\sqrt{\left( \mathfrak{E}_{2}^{\alpha }\right)
^{2}+A_{2}^{\alpha ,{energy}}}\sqrt{A_{2}^{\alpha ,{punct}}}%
\int_{\widehat{I}}t^{2}\overline{\mu }\ .
\end{eqnarray*}%
Recalling $B\left( M,M^{\prime }\right) =B\left( c_{M},\frac{1}{2}d\left(
c_{M},c_{M^{\prime }}\right) \right) $, we can write (suppressing some
notation for clarity) $T_{s}^{{intersection}}$ as
\begin{eqnarray*}
\!\!&=&\!\!\!\!\!
\sum_{F,F^{\prime }}\sum_{\substack{ M,M^{\prime }}}\int_{B\left(
M,M^{\prime }\right) }\frac{\left\Vert \mathsf{Q}_{F,M}^{\omega ,\mathbf{b}%
^{\ast }}x\right\Vert _{L^{2}\left( \omega \right) }^{\spadesuit 2}}{\left(
\left\vert M\right\vert^\frac{1}{n} +\left\vert y-c_{M}\right\vert \right) ^{\!n\!+\!1\!-\!\alpha }}%
\frac{\left\Vert \mathsf{Q}_{F^{\prime },M^{\prime }}^{\omega ,\mathbf{b}%
^{\ast }}x\right\Vert _{L^{2}\left( \omega \right) }^{\spadesuit 2}}{\left(
\left\vert M^{\prime }\right\vert^\frac{1}{n} +\left\vert y-c_{M^{\prime }}\right\vert
\right) ^{\!n\!+\!1\!-\!\alpha }}d\sigma \left( y\right) \\
&\approx &\!\!\!\!\!
\sum_{F,F^{\prime }}\sum_{\substack{ M,M^{\prime }}}
\frac{\left\Vert 
\mathsf{Q}_{F,M}^{\omega ,\mathbf{b}^{\ast }}x\right\Vert _{L^{2}\left(
\omega \right) }^{\spadesuit 2}\left\Vert \mathsf{Q}_{F^{\prime },M^{\prime
}}^{\omega ,\mathbf{b}^{\ast }}x\right\Vert _{L^{2}\left( \omega \right)
}^{\spadesuit 2}}{\left\vert c_{M}-c_{M^{\prime }}\right\vert
^{n+1-\alpha }}\!\!\int_{B\left( M,M^{\prime }\right) }\frac{d\sigma \left(
y\right) }{\left( \left\vert M\right\vert^\frac{1}{n} +\left\vert y-c_{M}\right\vert
\right) ^{n+1-\alpha }} \\
&\leq &\!\!\!\!\!
\sum_{F^{\prime }}\sum_{\substack{ M^{\prime }}}\left\Vert \mathsf{Q}%
_{F^{\prime },M^{\prime }}^{\omega ,\mathbf{b}^{\ast }}x\right\Vert
_{L^{2}\left( \omega \right) }^{\spadesuit 2}\!\!\sum_{F}\sum_{\substack{ M}}%
\frac{\left\Vert 
\mathsf{Q}_{F,M}^{\omega ,\mathbf{b}^{\ast }}x\right\Vert _{L^{2}\left(
\omega \right) }^{\spadesuit 2}}{\left\vert c_{M}-c_{M^{\prime }}\right\vert ^{n\!+\!1\!-\!\alpha }}
\!\!\int_{B\left( M,M^{\prime }\right) }\!\!\frac{
d\sigma \left( y\right) }{\left(\! \left\vert M\right\vert^\frac{1}{n} \!+\! \left\vert
y-c_{M}\right\vert \!\right) ^{\!n\!+\!1\!-\!\alpha }} \\
&\equiv &\!\!\!\!\!
\sum_{F^{\prime }}\sum_{\substack{ M^{\prime }}}\left\Vert \mathsf{Q%
}_{F^{\prime },M^{\prime }}^{\omega ,\mathbf{b}^{\ast }}x\right\Vert
_{L^{2}\left( \omega \right) }^{\spadesuit 2}S_{s}\left( M^{\prime }\right) ,
\end{eqnarray*}%
and since $\int_{B\left( M,M^{\prime }\right) }\frac{d\sigma \left( y\right) 
}{\left( \left\vert M\right\vert^\frac{1}{n} +\left\vert y-c_{M}\right\vert \right)
^{n+1-\alpha }}\approx \frac{\mathrm{P}^{\alpha }\left( M,\mathbf{1}_{B\left(
M,M^{\prime }\right) }\sigma \right) }{\left\vert M\right\vert^\frac{1}{n} }$, it
remains to show that for each fixed $M^{\prime }$,%
\begin{eqnarray*}
S_{s}\left( M^{\prime }\right) 
&\approx &
\sum_{F}\overset{\ast }{\sum 
_{\substack{ M:\ d\left( c_{M},c_{M^{\prime }}\right) \geq 2^{s\left(
1+\delta \right) }\ell \left( M^{\prime }\right) }}}\frac{\left\Vert \mathsf{%
Q}_{F,M}^{\omega ,\mathbf{b}^{\ast }}x\right\Vert _{L^{2}\left( \omega
\right) }^{\spadesuit 2}}{\left\vert c_{M}-c_{M^{\prime }}\right\vert
^{n+1-\alpha }}\frac{\mathrm{P}^{\alpha }\left( M,\mathbf{1}_{B\left(
M,M^{\prime }\right) }\sigma \right) }{\left\vert M\right\vert ^\frac{1}{n}} \\
&\lesssim &
2^{-\delta s}\sqrt{\left( \mathfrak{E}_{2}^{\alpha }\right)
^{2}+A_{2}^{\alpha ,{energy}}}\sqrt{A_{2}^{\alpha }}\ .
\end{eqnarray*}

We write%
\begin{eqnarray}
S_{s}\left( M^{\prime }\right)     \notag\label{def Sks} 
&\approx&
\sum_{k\geq s\left( 1+\delta \right) }\frac{1}{\left( 2^{k}\left\vert
M^{\prime }\right\vert \right) ^{n+1-\alpha }}S_{s}^{k}\left( M^{\prime
}\right) \ ;  \notag \\
S_{s}^{k}\left( M^{\prime }\right) 
&\equiv &
\sum_{F}\overset{\ast }{%
\sum_{M:\ d\left( c_{M},c_{M^{\prime }}\right) \approx 2^{k}\ell \left(
M^{\prime }\right) }}\left\Vert \mathsf{Q}_{F,M}^{\omega ,\mathbf{b}^{\ast
}}x\right\Vert _{L^{2}\left( \omega \right) }^{\spadesuit 2}\frac{\mathrm{P}%
^{\alpha }\left( M,\mathbf{1}_{B\left( M,M^{\prime }\right) }\sigma \right) 
}{\left\vert M\right\vert^\frac{1}{n} },  \notag
\end{eqnarray}%
where by $d\left( c_{M},c_{M^{\prime }}\right) \approx 2^{k}\ell \left(
M^{\prime }\right) $ we mean $2^{k}\ell \left( M^{\prime }\right) \leq
d\left( c_{M},c_{M^{\prime }}\right) \leq 2^{k+1}\ell \left( M^{\prime
}\right) $. Moreover, if $d\left( c_{M},c_{M^{\prime }}\right) \approx
2^{k}\ell \left( M^{\prime }\right) $, then from the fact that the radius of 
$B\left( M,M^{\prime }\right) $ is $\frac{1}{2}d\left( c_{M},c_{M^{\prime
}}\right) $, we obtain 
\begin{equation*}
B\left( M,M^{\prime }\right) \subset C_{0}2^{k}M^{\prime },
\end{equation*}%
where $C_{0}$ is a positive constant ($C_{0}=6$ works).

For fixed $k\geq s\left( 1+\delta \right) $, we invoke yet again the \emph{%
`prepare to puncture'} argument. Choose an augmented cube $\widetilde{%
M^{\prime }}\in \mathcal{AD}$ such that $C_{0}2^{k}M\subset \widetilde{%
M^{\prime }}$ and $\ell \left( \widetilde{M^{\prime }}\right) \leq
C2^{k}\ell \left( M^{\prime }\right) $. Define $\widetilde{\omega }=\omega
-\omega \left( \left\{ p\right\} \right) \delta _{p}$ where $p$ is an atomic
point in $\widetilde{M^{\prime }}$ for which 
\begin{equation*}
\omega \left( \left\{ p\right\} \right) =\sup_{q\in \mathfrak{P}_{\left(
\sigma ,\omega \right) }:\ q\in \widetilde{M^{\prime }}}\omega \left(
\left\{ q\right\} \right) .
\end{equation*}%
(If $\omega $ has no atomic point in common with $\sigma $ in $\widetilde{%
M^{\prime }}$, set $\widetilde{\omega }=\omega $.) Then we have $\left\vert 
\widetilde{M^{\prime }}\right\vert _{\widetilde{\omega }}=\omega \left( 
\widetilde{M^{\prime }},\mathfrak{P}_{\left( \sigma ,\omega \right) }\right) 
$ and so from (\ref{key obs}) and (\ref{key fact}), for any cube $A$,
\begin{equation*}
\frac{\left\vert \widetilde{M^{\prime }}\right\vert _{\widetilde{\omega }}}{%
\left\vert \widetilde{M^{\prime }}\right\vert ^{1-\frac{\alpha}{n} }}\frac{\left\vert 
\widetilde{M^{\prime }}\right\vert _{\sigma }}{\left\vert \widetilde{%
M^{\prime }}\right\vert ^{1-\frac{\alpha}{n} }}\leq A_{2}^{\alpha ,{punct}}%
\text{ and }\sum_{F\in \mathcal{F}}\left\Vert \mathsf{Q}_{F,A}^{\omega ,%
\mathbf{b}^{\ast }}x\right\Vert _{L^{2}\left( \omega \right) }^{\spadesuit
2}\lesssim \ell \left( A\right) ^{2}\left\vert A\right\vert _{\widetilde{%
\omega }}
\end{equation*}

Now we are ready to apply Cauchy-Schwarz and the augmented local estimate (%
\ref{shifted local}) in Lemma \ref{shifted} with $M=\widetilde{M^{\prime }}\ 
$to the second line below, and to apply the argument in (\ref{stan pig'}) to
the first line below, to get the following estimate for $S_{s}^{k}\left(
M^{\prime }\right) $ defined in (\ref{def Sks}) above:%
\begin{eqnarray*}
S_{s}^{k}\left( M^{\prime }\right) &\leq &\left( \sum_{F}\sum_{M:\ d\left(
c_{M},c_{M^{\prime }}\right) \approx 2^{k}\ell \left( M^{\prime }\right)
} \left\Vert \mathsf{Q}_{F,M}^{\omega ,\mathbf{b}^{\ast }}x\right\Vert
_{L^{2}\left( \omega \right) }^{\spadesuit 2}\right) ^{\frac{1}{2}} \\
&&
\times \left( \sum_{F}\sum_{M:\ d\left( c_{M},c_{M^{\prime }}\right)
\approx 2^{k}\ell \left( M^{\prime }\right) }\!\!\!\!\!\!\!\!\!\!\!\!\left\Vert \mathsf{Q}%
_{F,M}^{\omega ,\mathbf{b}^{\ast }}x\right\Vert _{L^{2}\left( \omega \right)
}^{\spadesuit 2}\left( \frac{\mathrm{P}^{\alpha }\left( M,\mathbf{1}%
_{B\left( M,M^{\prime }\right) }\sigma \right) }{\left\vert M\right\vert^\frac{1}{n} }%
\right) ^{2}\right) ^{\frac{1}{2}} \\
&\lesssim &
\left( 2^{2s}\left\vert \widetilde{M^{\prime }}\right\vert
^{2}\left\vert \widetilde{M^{\prime }}\right\vert _{\widetilde{\omega }%
}\right) ^{\frac{1}{2}}\left( \left[ \left( \mathfrak{E}_{2}^{\alpha
}\right) ^{2}+A_{2}^{\alpha ,{energy}}\right] \left\vert \widetilde{%
M^{\prime }}\right\vert _{\sigma }\right) ^{\frac{1}{2}} \\
&\lesssim &
\sqrt{\left( \mathfrak{E}_{2}^{\alpha }\right) ^{2}+A_{2}^{\alpha
,{energy}}}2^{s}\left\vert \widetilde{M^{\prime }}\right\vert \sqrt{%
\left\vert \widetilde{M^{\prime }}\right\vert _{\widetilde{\omega }}}\sqrt{%
\left\vert \widetilde{M^{\prime }}\right\vert _{\sigma }} \end{eqnarray*}
\begin{eqnarray*}
&\lesssim &
\sqrt{\left( \mathfrak{E}_{2}^{\alpha }\right) ^{2}+A_{2}^{\alpha
,{energy}}}\sqrt{A_{2}^{\alpha ,{punct}}}2^{s}\left\vert 
\widetilde{M^{\prime }}\right\vert \left\vert \widetilde{M^{\prime }}%
\right\vert ^{1-\alpha } \\
&\approx &
\sqrt{\left( \mathfrak{E}_{2}^{\alpha }\right) ^{2}+A_{2}^{\alpha ,%
{energy}}}\sqrt{A_{2}^{\alpha ,{punct}}}2^{s}2^{k\left(
1-\alpha \right) }\left\vert M^{\prime }\right\vert ^{n+1-\alpha }, 
\end{eqnarray*}
because
$ \ell \left( \widetilde{M^{\prime }}\right) \approx 2^{k}\ell
\left( M^{\prime }\right).$

Altogether then we have%
\begin{eqnarray*}
S_{s}\left( M^{\prime }\right) &=&\sum_{k\geq \left( 1+\delta \right) s}%
\frac{1}{\left( 2^{k}\left\vert M^{\prime }\right\vert \right) ^{n+1-\alpha }}%
S_{s}^{k}\left( M^{\prime }\right) \\
&\lesssim &
\sqrt{\left( \mathfrak{E}_{2}^{\alpha }\right) ^{2}+A_{2}^{\alpha
,{energy}}}\sqrt{A_{2}^{\alpha ,{punct}}}\!\!\!\!\sum_{k\geq \left(
1+\delta \right) s}\frac{2^{s}2^{k\left( 1-\alpha \right) }}{\left( 2^{k}\left\vert M^{\prime }\right\vert
\right) ^{n+1-\alpha }}\left\vert M^{\prime
}\right\vert ^{n+1-\alpha } \\
&=&
\sqrt{\left( \mathfrak{E}_{2}^{\alpha }\right) ^{2}+A_{2}^{\alpha ,%
{energy}}}\sqrt{A_{2}^{\alpha ,{punct}}}\sum_{k\geq \left(
1+\delta \right) s}2^{s-k} \\ &\lesssim& 2^{-\delta s}\sqrt{\left( \mathfrak{E}%
_{2}^{\alpha }\right) ^{2}+A_{2}^{\alpha ,{energy}}}\sqrt{%
A_{2}^{\alpha ,{punct}}},
\end{eqnarray*}%
which is summable in $s$. This completes the proof of (\ref{Us bound}), and
hence of the estimate for $\mathbf{Back}\left( \widehat{I}\right) $ in (\ref%
{e.t2 n'}).

The proof of Proposition \ref{func ener control} is now complete.

\section{Glossary}
\subsubsection{Section 1}

\begin{enumerate}
\item $C_{CZ}$; (\ref{sizeandsmoothness'})

\item $\mathfrak{N}_{T_{\sigma }^{\alpha }}$; (\ref{two weight'})

\item $T_{\sigma ,\delta ,R}^{\alpha }f\left( x\right) $; (\ref{def
truncation})

\item $p$\emph{-weakly }$\mu $\emph{-accretive} family; (\ref{acc})

\item $\mathfrak{T}_{T^{\alpha }}^{\mathbf{b}},\mathfrak{T}_{T^{\alpha ,\ast
}}^{\mathbf{b}^{\ast }}$; (\ref{b testing cond})

\item $\mathrm{P}^{\alpha }\left( Q,\mu \right) ,\mathcal{P}^{\alpha }\left(
Q,\mu \right) $; Subsection \ref{def Poisson}

\item $\mathcal{A}_{2}^{\alpha },\mathcal{A}_{2}^{\alpha ,\ast }$; Definition \ref{def
call A2}

\item $\mathfrak{P}_{\left( \sigma ,\omega \right) }$; Subsection \ref{def punct}

\item $A_{2}^{\alpha ,{punct}},A_{2}^{\alpha ,\ast ,{punct}}$%
; Subsection \ref{def punct}

\item $\mathfrak{A}_{2}^{\alpha }$; (\ref{def A2})

\item $\mathcal{E}_{2}^{\alpha },\mathcal{E}_{2}^{\alpha ,\ast },\mathfrak{E}%
_{2}^{\alpha },\mathcal{NTV}_{\alpha }$;$\ $(\ref{strong b* energy}), (\ref%
{strong b energy}), (\ref{def NTV})
\end{enumerate}

\subsubsection{Section 2}

\begin{enumerate}
\item reverse H\"{o}lder control of children (\ref{rev Hol con})

\item Calder\'{o}n-Zygmund stopping intervals; Definition \ref{CZ stopping
times}

\item $\mathbf{b}$-accretive/weak testing stopping intervals; Definition \ref%
{accretive stopping times gen}

\item $\sigma $-energy stopping intervals; Definition \ref{def energy corona
3}

\item $\mathbf{X}_{\alpha }\left( \mathcal{C}_{S}\right) $ energy stopping
times; (\ref{def stopping energy 3})

\item $\mathcal{D}_{\beta }$: (\ref{def dyadic grid})

\item $\mathbb{E}_{Q}^{\mu ,\mathbf{b}}f\left( x\right) $; (\ref{def
expectation})

\item $\widehat{\mathbb{F}}_{Q}^{\mu ,\mathbf{b}}f\left( x\right) $; (\ref{F
hat})

\item $\bigtriangleup _{Q}^{\mu ,\mathbf{b}}f\left( x\right) $, $\square
_{Q}^{\mu ,\mathbf{b}}f\left( x\right) $; (\ref{def diff})

\item $\bigtriangledown _{Q}^{\mu }$, $\widehat{\bigtriangledown }_{Q}^{\mu
} $; (\ref{Carleson avg op})

\item $\mathbb{E}_{Q}^{\mu ,\pi ,\mathbf{b}}f$, $\mathbb{F}_{Q}^{\mu ,\pi ,%
\mathbf{b}}f$; (\ref{def pi exp})

\item $\bigtriangleup _{Q}^{\mu ,\pi ,\mathbf{b}}f$, $\square _{Q}^{\mu ,\pi
,\mathbf{b}}f$; (\ref{def pi box})

\item $\square _{I}^{\sigma ,\flat ,\mathbf{b}}f$; (\ref{flat box})

\item $\widehat{\square }_{I}^{\sigma ,\flat ,\mathbf{b}}f$; (\ref{flat box
hat})

\item $\bigtriangleup _{I,{broken}}^{\mu ,\flat ,\mathbf{b}}f$, $%
\square _{I,{broken}}^{\mu ,\flat ,\mathbf{b}}f$; (\ref{def flat
broken})

\item $\Psi _{\mathcal{B}}^{\mu ,\mathbf{b}}f$:(
\ref{Psi op})
\end{enumerate}

\subsubsection{Section 3}

\begin{enumerate}
\item ${body}K$, body of an interval; (\ref{body})

\item $\varepsilon -{good}$ \emph{with respect to} an interval; Definition \ref{good arb}

\item $\mathcal{G}_{\left( k,\varepsilon \right) -{good}}^{\mathcal{D%
}}$; Definition \ref{good two grids}

\item $k$-${bad}$ in a grid; Definition \ref{bad in grid}

\item $R^{\maltese },\kappa \left( R\right) $: Definition \ref{def sharp cross}

\item $\Theta _{2}^{{bad}\natural }\left( f,g\right) $; (\ref%
{Theta_2^bad sharp})

\item $\mathcal{G}_{k-{bad}}^{\mathcal{D}}=\mathcal{G}_{\left(
k,\varepsilon \right) -{bad}}^{\mathcal{D}}$; Definition \ref{def
Gbad}

\item $\Theta _{2}^{{good}}\left( f,g\right) $; (\ref{def Theta 2
good})

\end{enumerate}

\subsubsection{Section 5}

\begin{enumerate}

\item $\left\{ E,F\right\} $; (\ref{def E,F})

\item $\mathrm{P}_{\delta }^{\alpha }\mathsf{Q}^{\omega }\left( J,\upsilon
\right) $; (\ref{def compact})

\end{enumerate}

\subsubsection{Section 6}

\begin{enumerate}
\item $\mathcal{C}_{B}^{\mathcal{G},{shift}}$; Definition \ref%
{shifted corona}

\item $\left\langle T_{\sigma }^{\alpha }\left( \mathsf{P}_{\mathcal{C}_{A}^{%
\mathcal{D}}}^{\sigma ,\mathbf{b}}f\right) ,\mathsf{P}_{\mathcal{C}_{B}^{%
\mathcal{G},{shift}}}^{\omega ,\mathbf{b}^{\ast }}g\right\rangle
_{\omega }^{\Subset _{\mathbf{r},\varepsilon }}$; (\ref{def shorthand})

\item $\mathsf{B}_{\Subset _{\mathbf{r},\varepsilon }}^{A}\left( f,g\right) $%
; (\ref{def local})

\item $\mathfrak{F}_{\alpha }=\mathfrak{F}_{\alpha }^{\mathbf{b}^{\ast
}}\left( \mathcal{D},\mathcal{G}\right) $; (\ref{e.funcEnergy n})

\item $\mathsf{B}_{{stop}}^{A}\left( f,g\right) $; (\ref{bounded
stopping form}), (\ref{def stop})

\item $\mathcal{C}_{F}^{\mathcal{G},{shift}}$; Definition \ref{def
shift}

\item $\mathsf{Q}_{\mathcal{H}}^{\omega ,\mathbf{b}^{\ast }}$; (\ref{def
localization})

\item $\mathsf{B}_{{broken}}^{A}\left( f,g\right) $; (\ref{broken
vanish})

\item $\mathsf{B}_{{neighbour}}^{A}\left( f,g\right) $; (\ref{def
neighbour})

\end{enumerate}

\subsubsection{Section 7}

\begin{enumerate}

\item $\mathsf{B}_{{stop}}^{A,\mathcal{P}}\left( f,g\right) $; (\ref%
{def stop P})

\item $\varphi _{J}^{\mathcal{P}}$; above (\ref{phi bound})

\item $\widehat{\mathfrak{N}}_{{stop},\bigtriangleup ^{\omega }}^{A,%
\mathcal{P}}$; (\ref{best hat})

\item $\Pi _{2}^{K}\mathcal{P}$; (\ref{rest K})

\item $\Pi _{1}^{{below}}\mathcal{P}$; (\ref{Pi below})

\item $\mathcal{S}_{{init}{size}}^{\alpha ,A}\left( \mathcal{%
P}\right) ^{2}$; (\ref{def ext size})

\item $\Pi _{2}^{K,{aug}}\mathcal{P}$; Definition \ref{augs}

\item $\mathcal{S}_{{aug}{size}}^{\alpha ,A}\left( \mathcal{P%
}\right) $; Definition \ref{augs}, (\ref{def P stop energy' 3})

\item $\mathcal{P}_{{cor}}^{A}$; (\ref{def cor})

\item $\mathcal{Q}$ $\flat $\textbf{straddles} $\mathcal{S}$; Definition \ref%
{flat straddles}

\item $\mathcal{W}^{\ast }\left( S\right) $; (\ref{def Whit})

\item $\mathcal{S}_{{loc}{size}}^{\alpha ,A;S}\left( 
\mathcal{Q}\right) ^{2}$; (\ref{localized size ref})

\item $\mathcal{Q}$ \textbf{substraddles} $L$; (\ref{def substraddles})

\item $\omega _{\mathcal{P}}$ and $\omega _{\flat \mathcal{P}}$; (\ref{def
atomic})

\item $\mathcal{G}_{d}$; (\ref{geom depth})

\item $\mathcal{P}^{\flat big}$, $\mathcal{P}^{\flat small}$; (\ref{def big
small flat})
\end{enumerate}

\subsubsection{Section 9}

\begin{enumerate}
\item $\mathsf{Q}_{\mathcal{C}_{F}^{\mathcal{G},{shift}};M}^{\omega ,%
\mathbf{b}^{\ast }}$; (\ref{def pseudo rest})

\item $\mathcal{C}_{F}^{\mathcal{G},{shift}}$, $\mathcal{C}_{F}^{%
\mathcal{G},{shift}};K$; (\ref{def shift cor rest})

\item $J\Subset _{\mathbf{\rho },\varepsilon }K$; (\ref{def deep embed})

\item $\mathcal{M}_{\left( \mathbf{\rho },\varepsilon \right) -{deep}%
,\mathcal{G}}\left( K\right) $, $\mathcal{M}_{\left( \mathbf{\rho }%
,\varepsilon \right) -{deep},\mathcal{D}}\left( K\right) $, $%
\mathcal{W}\left( K\right) $; (\ref{def M_r-deep})

\item augmented dyadic grid $\mathcal{AD}$; Definition \ref{def dyadic}

\item $\mathsf{Q}_{K}^{\omega ,\mathbf{b}^{\ast }}$; (\ref{large pseudo})

\item $\mathcal{E}_{2}^{\alpha ,{Whitney}{partial}}$; (\ref%
{plug})

\item $A_{2}^{\alpha ,{energy}}$, $A_{2}^{\alpha ,\ast ,{%
energy}}$; (\ref{def energy A2})

\item $\mathcal{E}_{2}^{\alpha ,{Whitney}{plug}}$; (\ref{def
deep plug})

\item $\mu $; (\ref{def mu n})

\item $\mathsf{P}_{F,K}^{\omega ,\mathbf{b}^{\ast }}\equiv \mathsf{P}_{%
\mathcal{C}_{F}^{\mathcal{G},{shift}};K}^{\omega ,\mathbf{b}^{\ast
}} $; (\ref{def F,K})

\item $\mathbf{Local}\left( I\right) $, $\overline{\mu }$; (\ref{def local
forward})

\item $\mathcal{J}^{\ast }$; (\ref{def J*})

\item $\mathbf{Back}\left( \widehat{I}\right) $; (\ref{e.t2 n'})

\item $U_{s}$; (\ref{def Us})

\item $B\left( M,M^{\prime }\right) $; (\ref{def BMM'})

\item $T_{s}^{{intersection}}$; (\ref{def Tints})

\item $T_{s}^{{proximal}}$; (\ref{def Tproxs})

\item $T_{s}^{{difference}}$; (\ref{def Tdiffs})

\item $\overset{\ast }{\sum }$; Notation \ref{Sum *}

\item $\mathcal{W}_{M}^{s}$; (\ref{def WMs})

\item $\mathcal{W}_{M}^{s,\ell }$; (\ref{def WMsell})
\end{enumerate}

\end{document}